\documentclass[10pt]{article}
\usepackage{amsmath,amssymb,amsthm}
\usepackage{color}
\usepackage[nohead,margin=1.0in]{geometry}
\usepackage{booktabs}
\usepackage{graphicx,caption,subcaption}
\usepackage{grffile}
\usepackage{algorithm,algorithmic}
\usepackage{verbatim}
\usepackage{url}

\theoremstyle{plain}

\numberwithin{equation}{section}

\newcommand{\N}{\mathbb{N}}
\newcommand{\R}{\mathbb{R}}

\title{L1 data fitting for robust reconstruction in magnetic particle imaging: quantitative evaluation on \texttt{Open MPI dataset}}
\author{Tobias Kluth\thanks{Center for Industrial Mathematics, University of Bremen, Bibliothekstr. 5, 28357 Bremen, Germany (\texttt{tkluth@math.uni-bremen.de})} \and
Bangti Jin\thanks{Department of Computer Science, University College London, Gower Street, London WC1E 6BT, UK (\texttt{b.jin@ucl.ac.uk, bangti.jin@gmail.com})}}

\begin{document}
\maketitle

\begin{abstract}
Magnetic particle imaging is an emerging quantitative imaging modality, exploiting the unique nonlinear magnetization
phenomenon of superparamagnetic iron oxide nanoparticles for recovering the concentration. Traditionally the reconstruction
is formulated into a penalized least-squares problem with nonnegativity constraint, and then solved using a variant of
Kaczmarz method which is often stopped early after a small number of iterations. Besides the phantom signal, measurements additionally include a background signal and a noise signal. In order to obtain good reconstructions, a preprocessing step of frequency selection to remove the
deleterious influences of the noise is often adopted. In this work, we propose a complementary pure variational approach
to noise treatment, by viewing highly noisy measurements as outliers, and employing the l1 data fitting, one popular
approach from robust statistics. When compared with the standard approach, it is easy to implement with a comparable computational
complexity. Experiments with a public domain dataset, i.e., \texttt{Open MPI dataset}
\cite{KnoppOpenMPI}, show that it can give accurate reconstructions, and is less prone to noisy measurements, which is illustrated by quantitative (PSNR / SSIM) and qualitative comparisons with the Kaczmarz method. We also investigate the performance of the Kaczmarz method for small iteration numbers quantitatively.
\end{abstract}

\section{Introduction}

Magnetic particle imaging (MPI), invented by Gleich and Weizenecker in 2005 \cite{Gleich2005}, is a relatively new medical  imaging modality.
It exploits the unique nonlinear magnetization behavior of super-paramagnetic iron oxide nanoparticles in an applied magnetic field. In the
experiment, a static magnetic field (selection field), given by a gradient field, generates a field free point (FFP) (or a field free line (FFL)
\cite{Weizenecker2008ffl}), and its superposition with a spatially homogeneous but time-dependent field (drive field) moves the field free
region along a predefined trajectory defining the field-of-view. The change of the applied field causes a change of the nanoparticle magnetization,
which can be measured for recovering the spatially dependent concentration of nanoparticles. See the surveys
\cite{KnoppBuzug:2012,Knopp2017,Kluth:2018} for relevant physics, instrumentation and mathematical modeling.

In comparison with more traditional imaging modalities, e.g., ultrasound, MRI and PET, MPI has a number of distinct features: high
temporal / spatial resolution, high sensitivity and free from the need of harmful radiation. Thus it is especially attractive
for {\it in-vivo} applications, and the list of potential medical applications is long and fast growing, including imaging blood
flow \cite{weizenecker2009three}, long-term tracer monitoring \cite{khandhar2017evaluation}, estimating potential flow \cite{Franke:2017},
tracking medical instruments \cite{haegele2012magnetic}, tracking and guiding instruments for angioplasty \cite{Salamon:2016cz}, cancer
detection \cite{Yu2017} and cancer treatment by hyperthermia \cite{murase2015usefulness}.

Hence, the MPI reconstruction problem is of great importance, and has received much attention \cite{weizenecker2009three,RahmerWeizenecker:2009,
Knopp2010e,GoodwillConolly:2011,GoodwillConolly:2012,RahmerHalkola:2015,Knopp_Online_2016,Storath:2017,KluthMaass:2017,BrandtSeppanen:2018,KluthJin:2019,
DittmerKluth:2018} (see \cite{Knopp2017} for an overview). Roughly, existing approaches can be categorized into two groups, i.e.,
data-based v.s. model-based, dependent of the description of the forward map. The former uses a measured forward map,
whereas the latter employs a mathematical model to describe the forward map, where the equilibrium model based on Langevin
theory is popular \cite{Kluth:2018,KluthJinLi:2017}. The data-based approach is predominant in practice, since it can deliver better quality reconstructions.

\begin{figure}[hbt!]
\centering
\scalebox{1}{
\begin{minipage}{0.23\textwidth}
\centering
$x$-coil \\
\includegraphics[width=\textwidth]{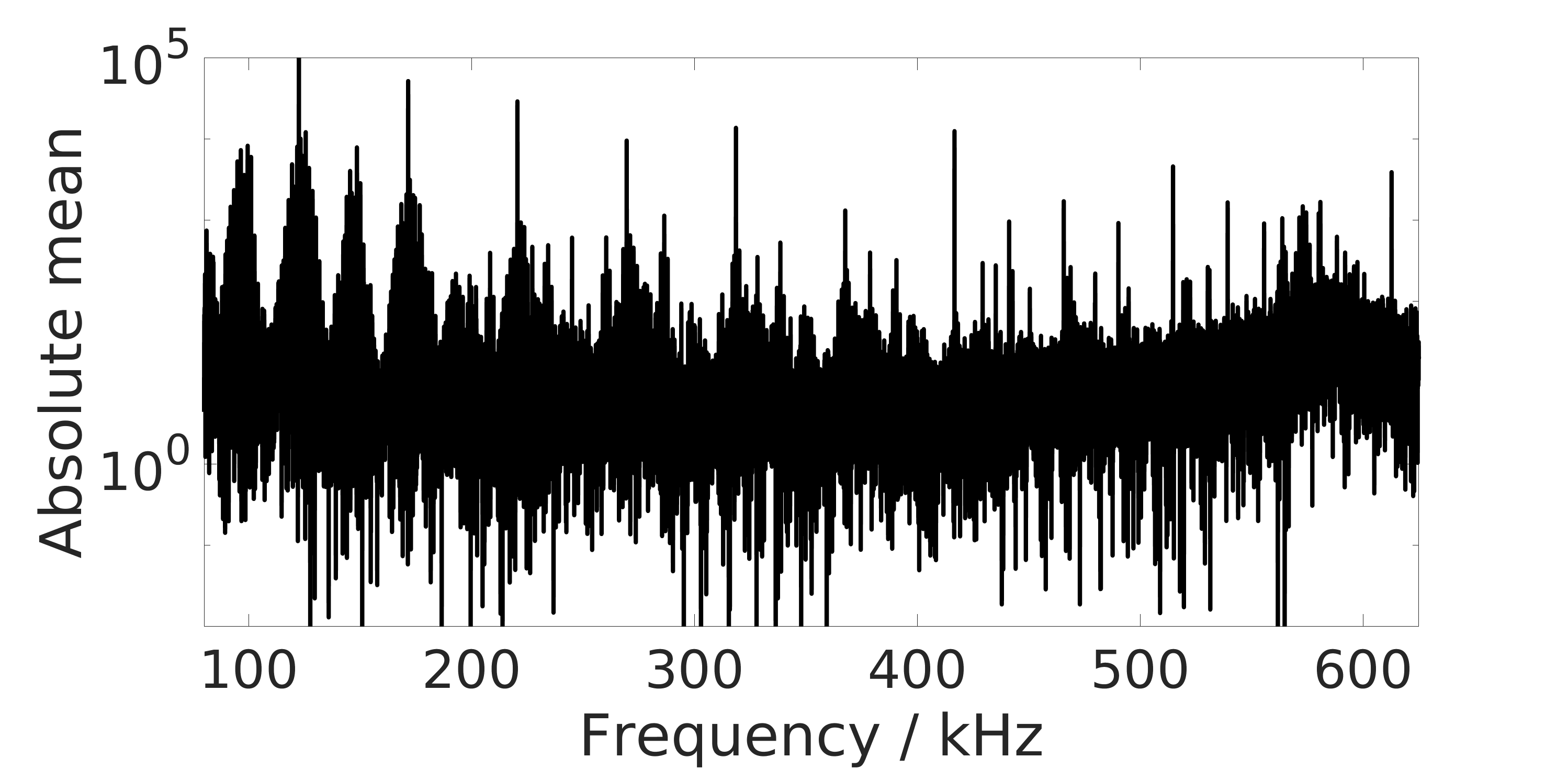} \\
\includegraphics[width=\textwidth]{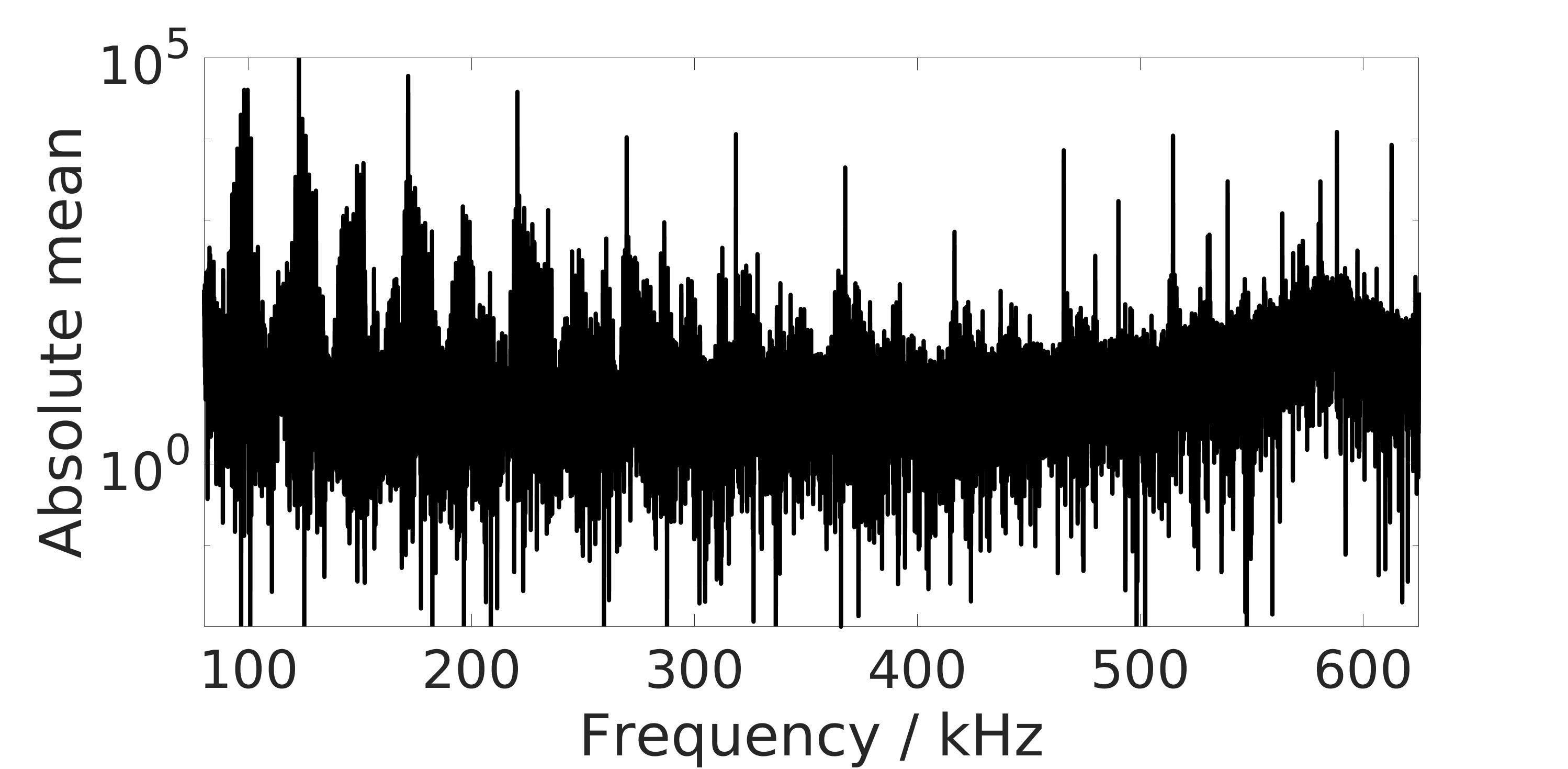} \\
 \end{minipage}
 \begin{minipage}{0.23\textwidth}
\centering
\ \\
\includegraphics[width=\textwidth]{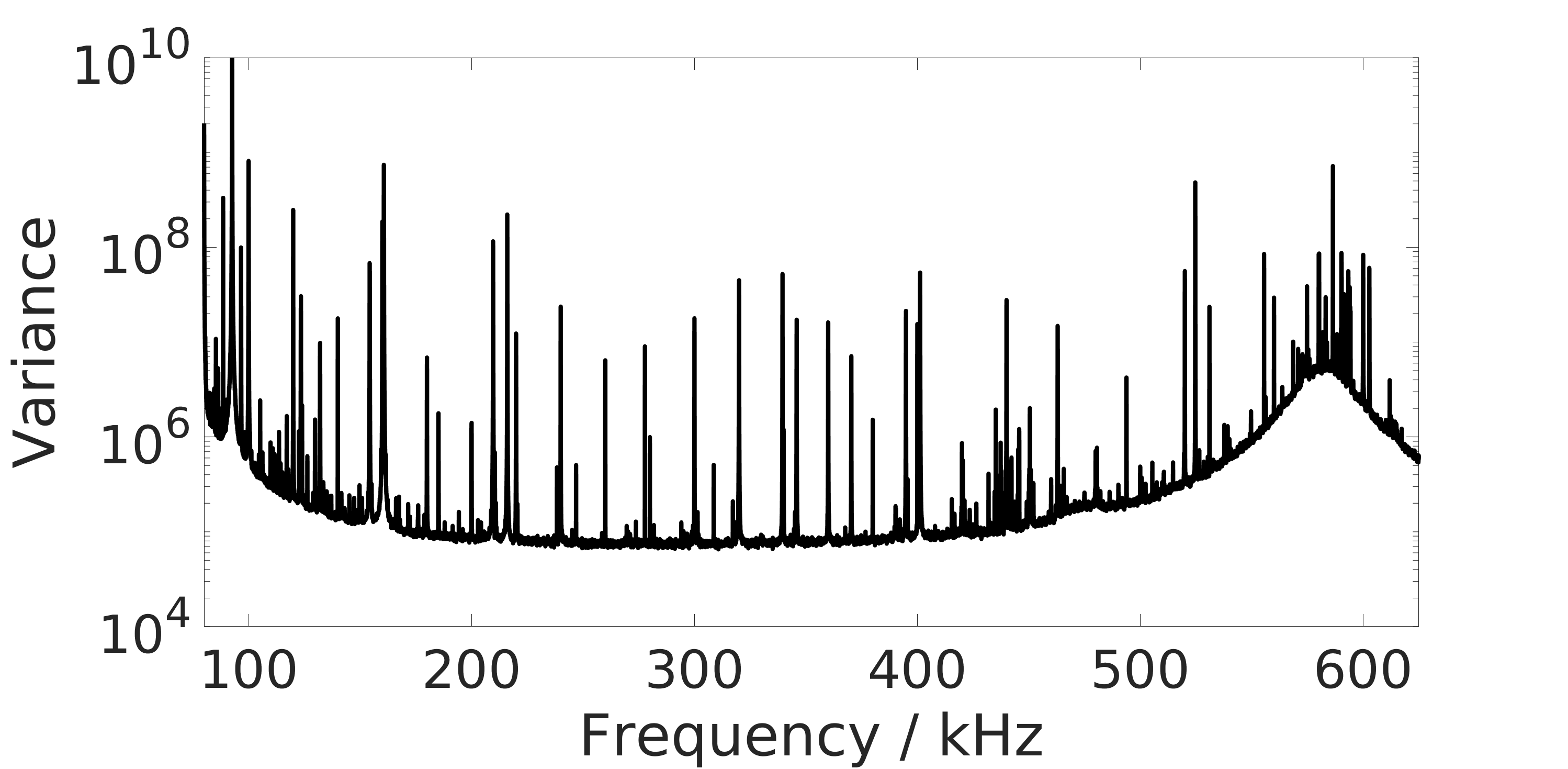} \\
\includegraphics[width=\textwidth]{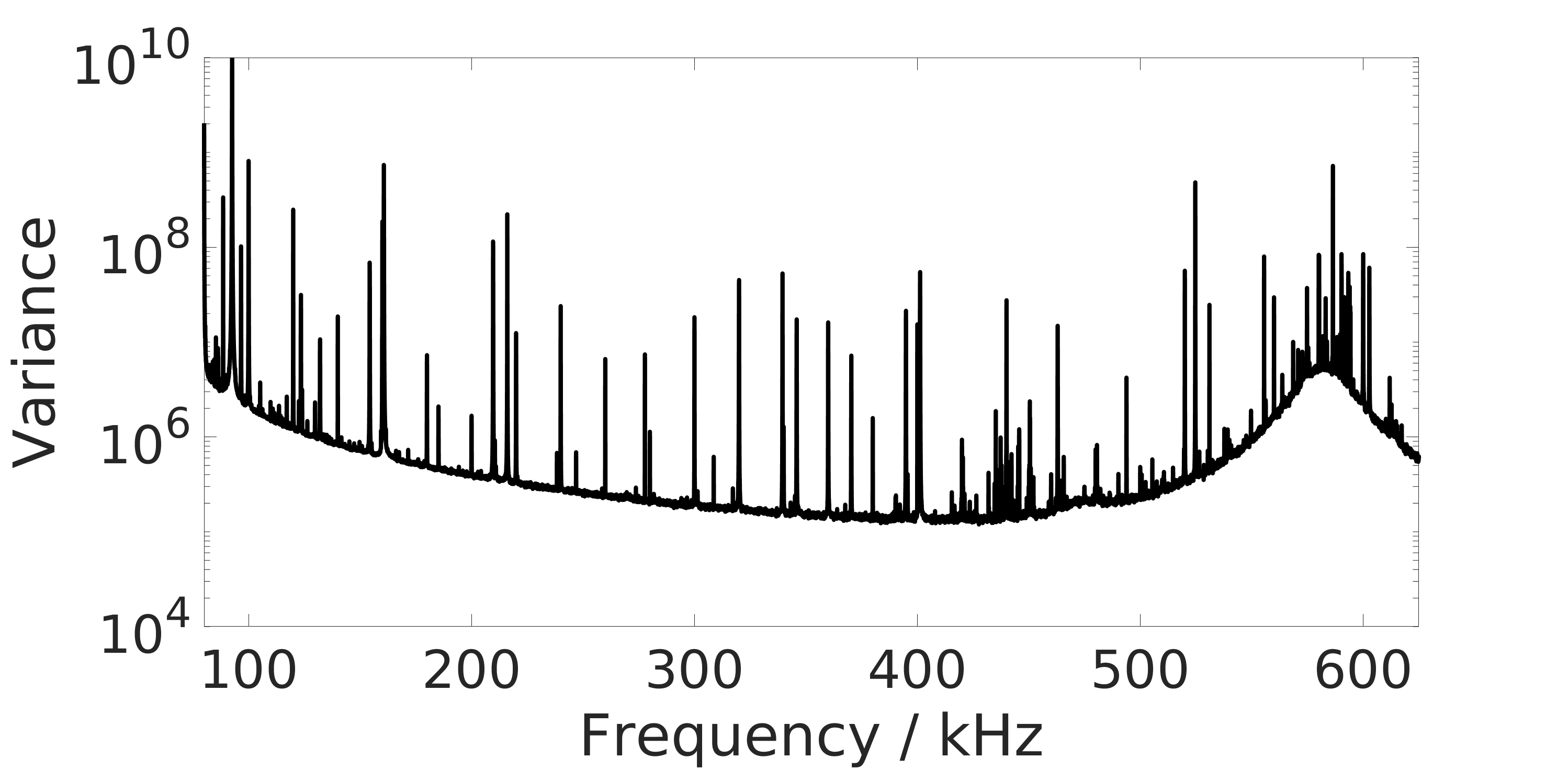} \\
 \end{minipage}
\begin{minipage}{0.23\textwidth}
\centering
$y$-coil \\
\includegraphics[width=\textwidth]{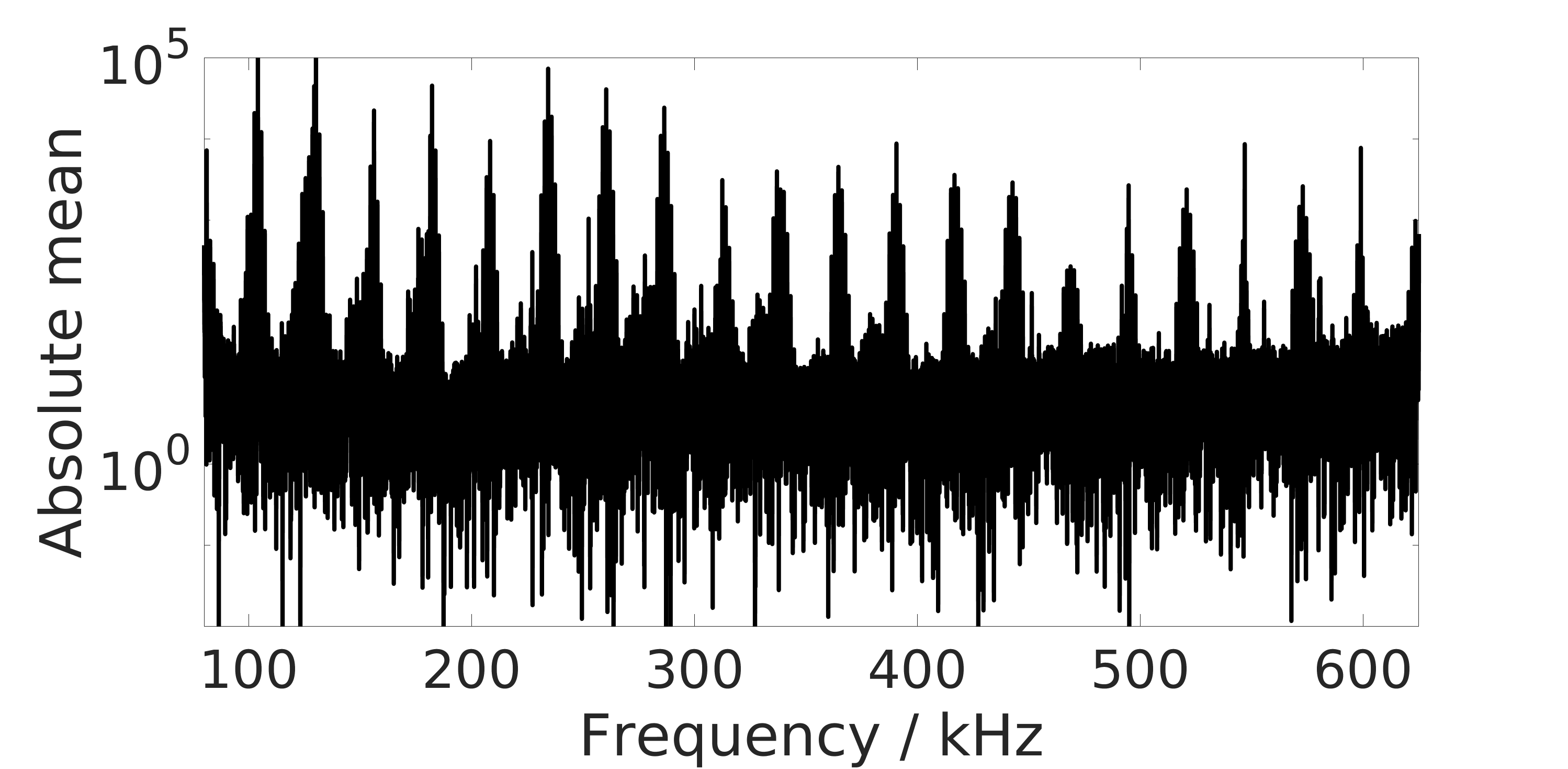} \\
\includegraphics[width=\textwidth]{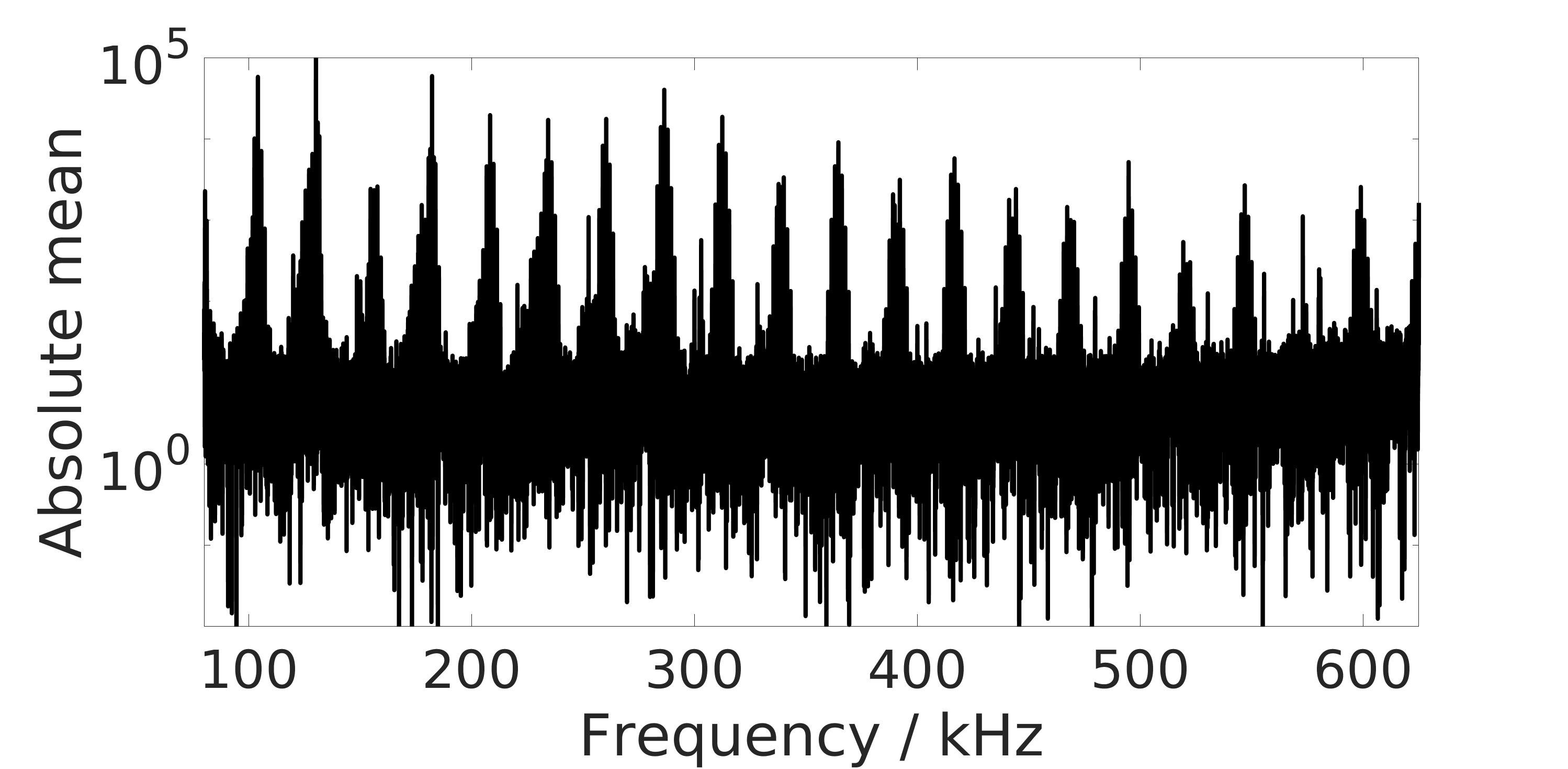} \\
 \end{minipage}
 \begin{minipage}{0.23\textwidth}
\centering
\ \\
\includegraphics[width=\textwidth]{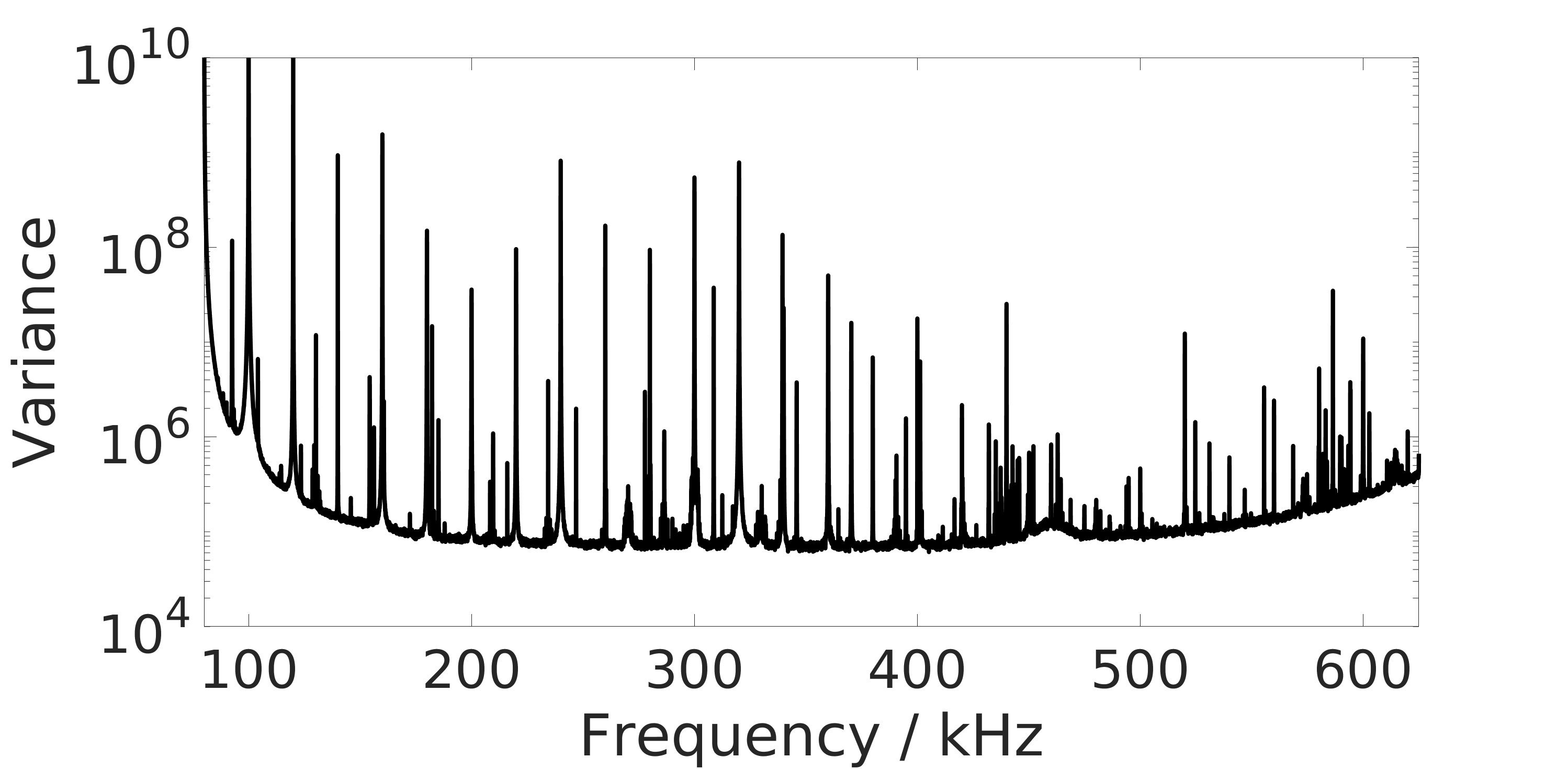} \\
\includegraphics[width=\textwidth]{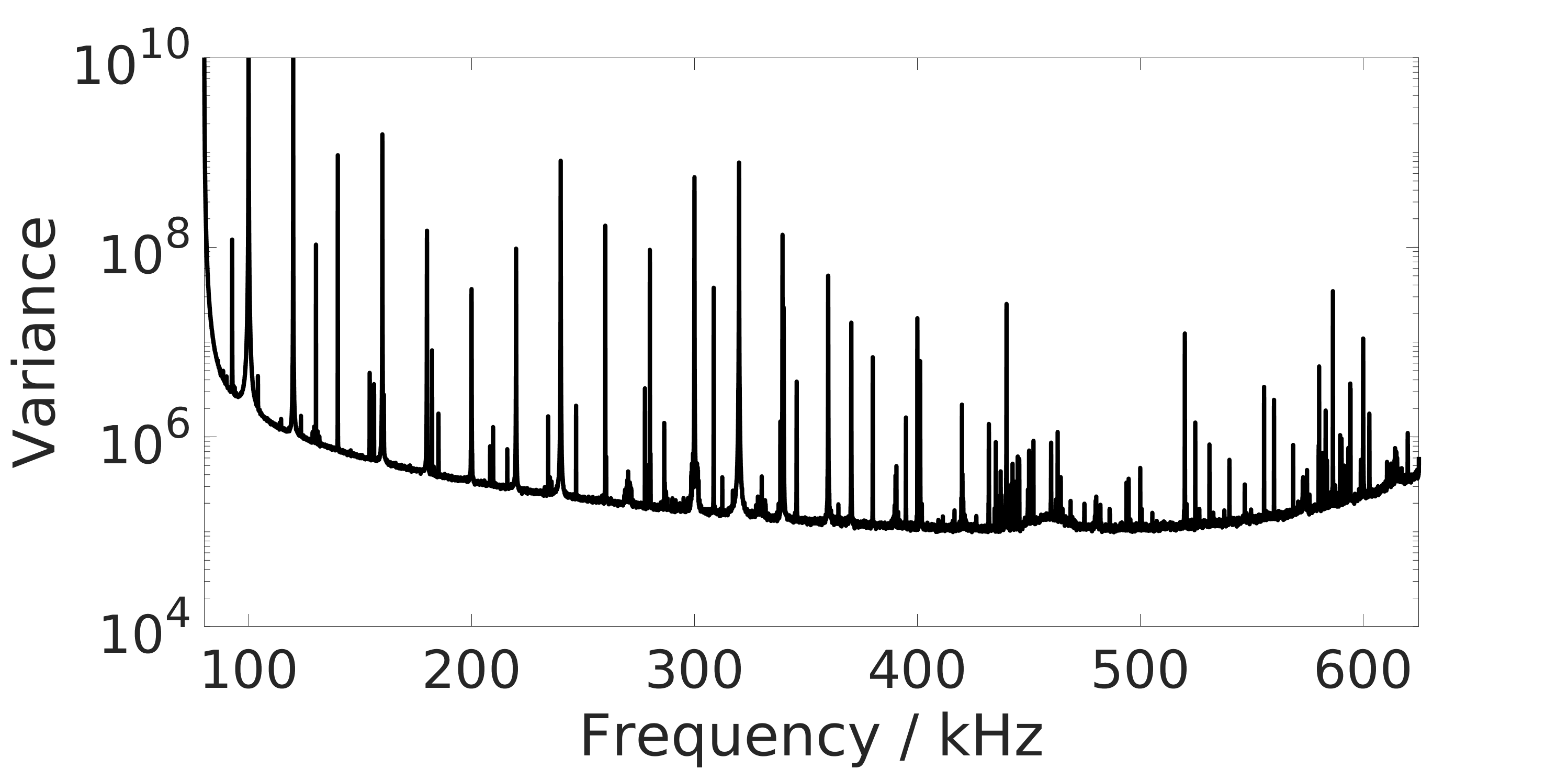} \\
 \end{minipage}

}
\caption{Mean and variance of the background measurement provided with the ``shape'' phantom from \texttt{Open MPI dataset},
computed from 1000 empty scans when using a 3D excitation in the preclinical Bruker MPI system. Visualized
individually for $x$ and $y$ receive coils with respect to the frequency; real part (top),
imaginary part (bottom).
}\label{fig:cov_diag}
\end{figure}

The starting point of this study is the following empirical observation: the recorded
signal often contains a significant amount of frequency-dependent anomalies. In Fig. \ref{fig:cov_diag}, we show mean and variance of
background measurements provided with the ``shape'' phantom from \texttt{Open MPI dataset} \cite{KnoppOpenMPI} (available
at \url{www.tuhh.de/ibi/research/open-mpi-data.html}, accessed on September 20, 2019) acquired at $x$- and $y$-receive coils.
It is observed that the variance of some data points is huge and deviate enormously from the bulk (noting the log scale for $y$-axis) in both
real and imaginary parts and do not follow an independent and identically distributed (i.i.d.) Gaussian model. In practice,
they are often deemed data ``noise'', and an important question in image reconstruction is how to deal with such noise.

In data-based approaches, this noise contributes to two important components in the MPI reconstruction setup,
i.e., calibrated system matrix and phantom measurement. In both cases one needs to distinguish three components in the
individual measurements: the phantom/system matrix signal, the background signal, and a noise signal. Note that the variance
in Fig. \ref{fig:cov_diag} illustrates the statistical characteristics of the noise signal, while the presented absolute mean
illustrates the mean estimate of the background signal. The mean is characterized by some dominant peaks at higher harmonics
of the corresponding excitation frequency (i.e., each excitation frequency corresponds to a certain space dimension like the
respective receive coil, particularly for a transmit-receive configuration). The interested reader is referred to Fig.
\ref{fig:cov_diag_harmonics_xcoil} in the supplementary material, where this is highlighted for the $x$ receive coil from
Fig. \ref{fig:cov_diag}. Surprisingly, the relationship between large peaks in the variance and the higher harmonics of the
excitation frequency remains unclear. On the system matrix side, there exist different approaches to address the background
removal issue, e.g., subtraction methods \cite{ThemKaul:2016,KnoppGdaniec:2019} prior the reconstruction or joint reconstruction
methods \cite{KluthMaass:2017,StraubSchulz:2018}. Particularly, in \cite{StraubSchulz:2018}, the authors performed a joint
reconstruction of the same background signal in the system matrix and the phantom measurement. However, while the works
\cite{ThemKaul:2016,KnoppGdaniec:2019} include signal drift assumptions on the background signal and a Gaussian assumption
on the phantom measurement noise the other works \cite{KluthMaass:2017,StraubSchulz:2018} incorporate a Gaussian assumption
on the noise distribution for the phantom measurement as well as for the background measurement (in terms of a l2 norm in a
Tikhonov-type functional). Even if one includes drift assumptions on the background signal, it is not guaranteed that the
noise characteristic is taken into account properly, when using a Gaussian assumption (often implicitly by using an l2 data
fitting term, see also further specifications below). Thus, it is of interest to develop statistically consistent MPI
reconstruction under the given circumstances.

The data-based image reconstruction in MPI (given a background-corrected system matrix) often proceeds as follows.
First, one selects a number of frequencies in a background-corrected phantom measurement based on a
suitable ``signal-to-noise ratio (SNR)'' criterion \cite{Franke:2016,Knopp_Online_2016}, and discards the remaining frequencies; see Section
\ref{ssec:freqsel} for details. This step partly removes the nongaussian component of the noise. The reconstruction is
performed using only the selected frequencies, often formulated into a penalized least-squares problem with nonnegativity
constraint and solved by a variant of Kaczmarz method \cite{Kaczmarz:1937}. This strategy has achieved great empirical
success, and is now a standard MPI reconstruction algorithm.

Nonetheless, there are still several issues on the overall reconstruction procedure. First, frequency selection as a noise treatment procedure is \textit{ad hoc} in nature,
since the threshold $\tau$ is nontrivial to set due to dependence on noise, while exerts big effect on imaging
quality. A too large $\tau$ may throw away informative data points, whereas a too small $\tau$ may invalidate i.i.d. Gaussian assumption of
the l2 fitting. Second, the performance relies on hybrid regularizing effects of both penalty and early stopping of Kaczmarz iteration (e.g., a few sweeps over the data),
and the delicate interplay has not been fully exposed. The hybridization greatly complicates the choice of the regularization
parameter and the stopping index. Third and last, the background signal and its noise distribution in the measured forward map is often not fully accounted for,
even though there are important efforts in that direction \cite{KluthMaass:2017,StraubSchulz:2018,BrandtSeppanen:2018}. This issue is also related to the
proper mathematical modeling in MPI and is fundamental towards an accurate model-based reconstruction algorithm as one commonly needs to solve a parameter identification problem to determine a model-based system function (see
\cite{KluthKnopp:2019} for recent progress).

Common methods for image reconstruction in MPI aim at minimizing a Tikhonov-type functional. 
Given an ill-posed operator equation $Ax=y$, the standard Tikhonov functional reads
\begin{equation*}
 \frac12 \|Ax-y\|_2^2 + \frac{\alpha}{2} \|x\|^2_2
\end{equation*}
where the first term is the \textit{data fitting} or \textit{discrepancy} term incorporating a certain noise model and where the second term is the regularization or penalty term incorporating a priori knowledge on the solution.
While several methods applied to MPI focus on various choices of a priori knowledge on the solution \cite{weizenecker2009three,Knopp2010e,Konkle2015a,Storath:2017,KluthMaass:2017,Bathke2018IWMPI} (see also the review \cite{Knopp2017}), i.e., formulations for the penalty term, in the present work the focus is on the data fitting term.
In this work, we present a complementary approach to noise treatment by SNR type frequency selection which partly addresses the first challenge raised above.
It is based on the standard l1 data fitting (i.e., using the l1 norm in the data fitting term), or equivalently a Laplace model on the noise (see Section \ref{ssec:Laplace} for details), which has been popular in several areas, e.g., signal processing \cite{AllineyRuzinsky:1994} and image processing \cite{ClasonJinKunisch:2010},
but it has not been applied to MPI reconstruction yet, to the best of our knowledge. The rational is that the nonselected
frequencies deviate largely from the bulk of the signals and thus can be viewed as outliers, and l1 fitting is known to be more robust
with respect to outliers  than the l2 fitting \cite{GelmanCarlinSternRubin:2004}. The approach allows adaptive use of the data and thus
full exploitation of the given data for better reconstructions. Numerically, l1 fitting leads to a convex but nonsmooth optimization
problem, which can be solved efficiently by many modern stand-alone optimization solvers. We employ a popular variant of the limited
memory BFGS algorithm, i.e., L-BFGS-B \cite{ByrdLu:1995,ZhuByrd:1997}. We carry out extensive numerical experiments with \texttt{Open
MPI dataset}. Our findings include that the l1 approach can indeed yield excellent reconstructions both quantitatively in terms of
PSNR and SSIM and qualitatively in terms of background and sharpness. The l1 approach is able to compete with the standard iterative 
Kaczmarz-type approach which yields high quality MPI reconstructions when using small numbers of iterations, while the variational l2 
approach fails. Thus, these techniques may facilitate fast and accurate MPI
reconstruction using variational regularization techniques. To the best of our knowledge, this is the first work
presenting quantitative results in terms of standard image quality measures (PSNR/SSIM) for phantom MPI data in the \texttt{Open MPI dataset}.

Note that one should not mix l1 (data) fitting with l1 norm penalty
that has been widely used in compressed sensing \cite{Donoho:2006} and recently also in MPI
reconstruction (see, e.g., \cite{Storath:2017,KluthMaass:2017,IlbeyTop:2019}). The latter assumes the sparsity of the solution instead of
noise, and thus is drastically different from l1 fitting of this work.

The rest of the paper is organized as follows. In Section \ref{sec:method}, we motivate and develop the robust formulation,
and describe the limited-memory BFGS algorithm. In Section \ref{sec:numer}, we present extensive quantitative and qualitative
numerical results to showcase the performance of the proposed approach and to investigate the interplay of regularization methods in a standard method. In Section \ref{sec:conc}, we give concluding remarks and further discussions. In the supplements, we provide additional numerical results.

\section{Methodologies}\label{sec:method}

In this part we describe the standard approach and develop the l1 approach.

\subsection{Standard approach}\label{ssec:freqsel}

\begin{figure*}[hbt!]
\centering
\scalebox{0.9}{
\begin{minipage}{0.3\textwidth}
\centering
$\tau=1$ \\
\includegraphics[width=\textwidth]{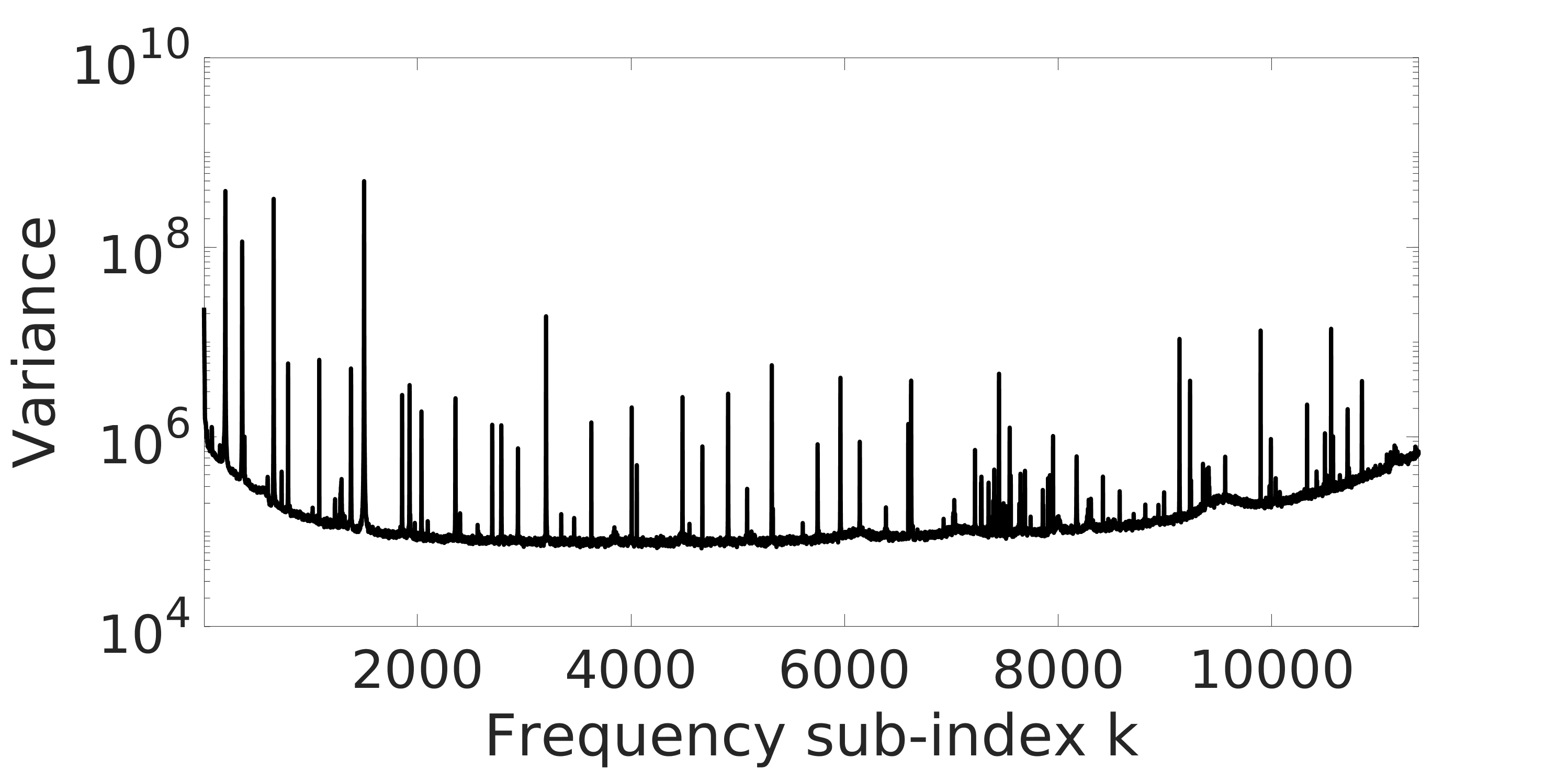} \\
\includegraphics[width=\textwidth]{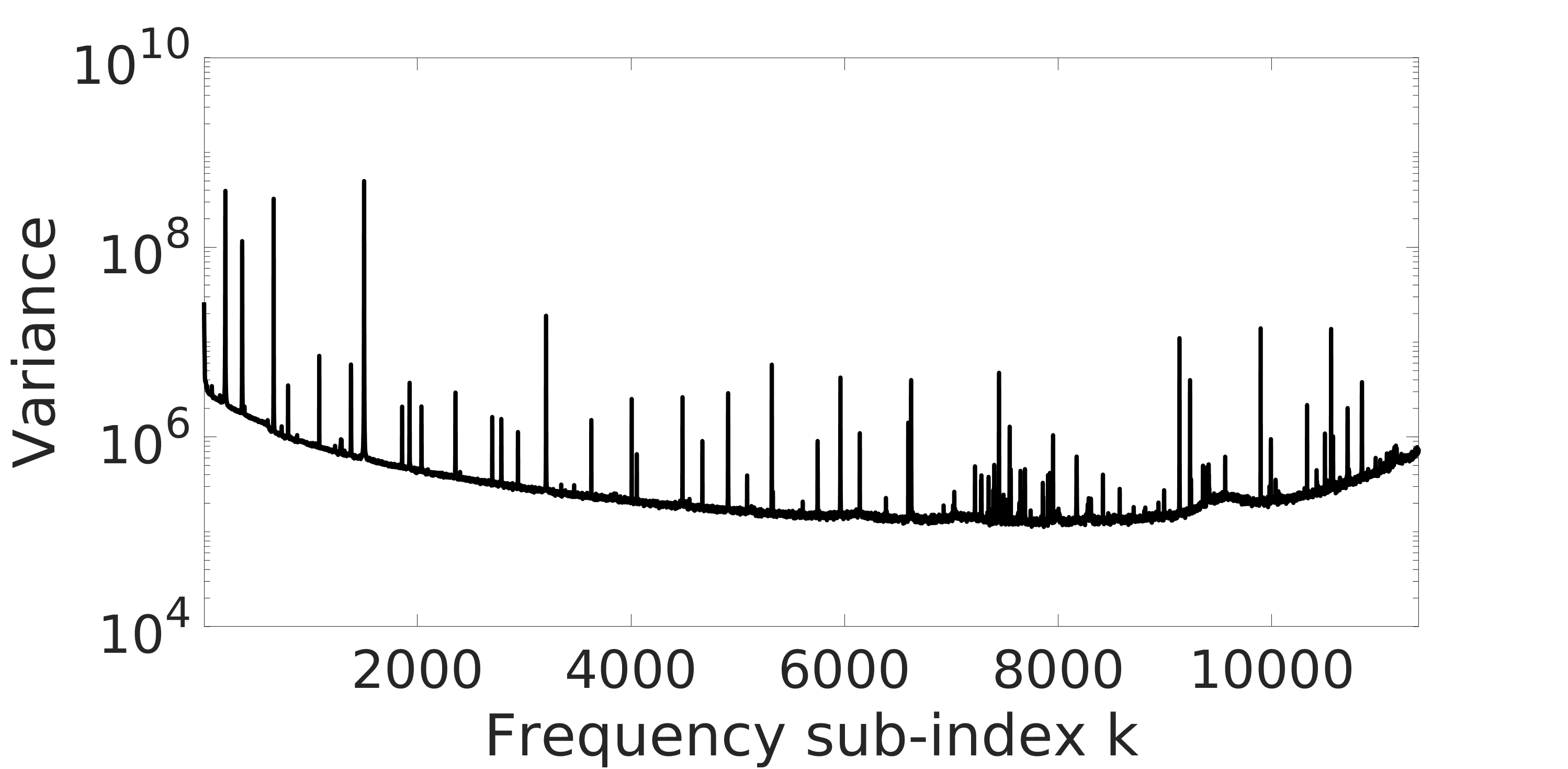} \\
 \end{minipage}
\begin{minipage}{0.3\textwidth}
\centering
$\tau=3$ \\
\includegraphics[width=\textwidth]{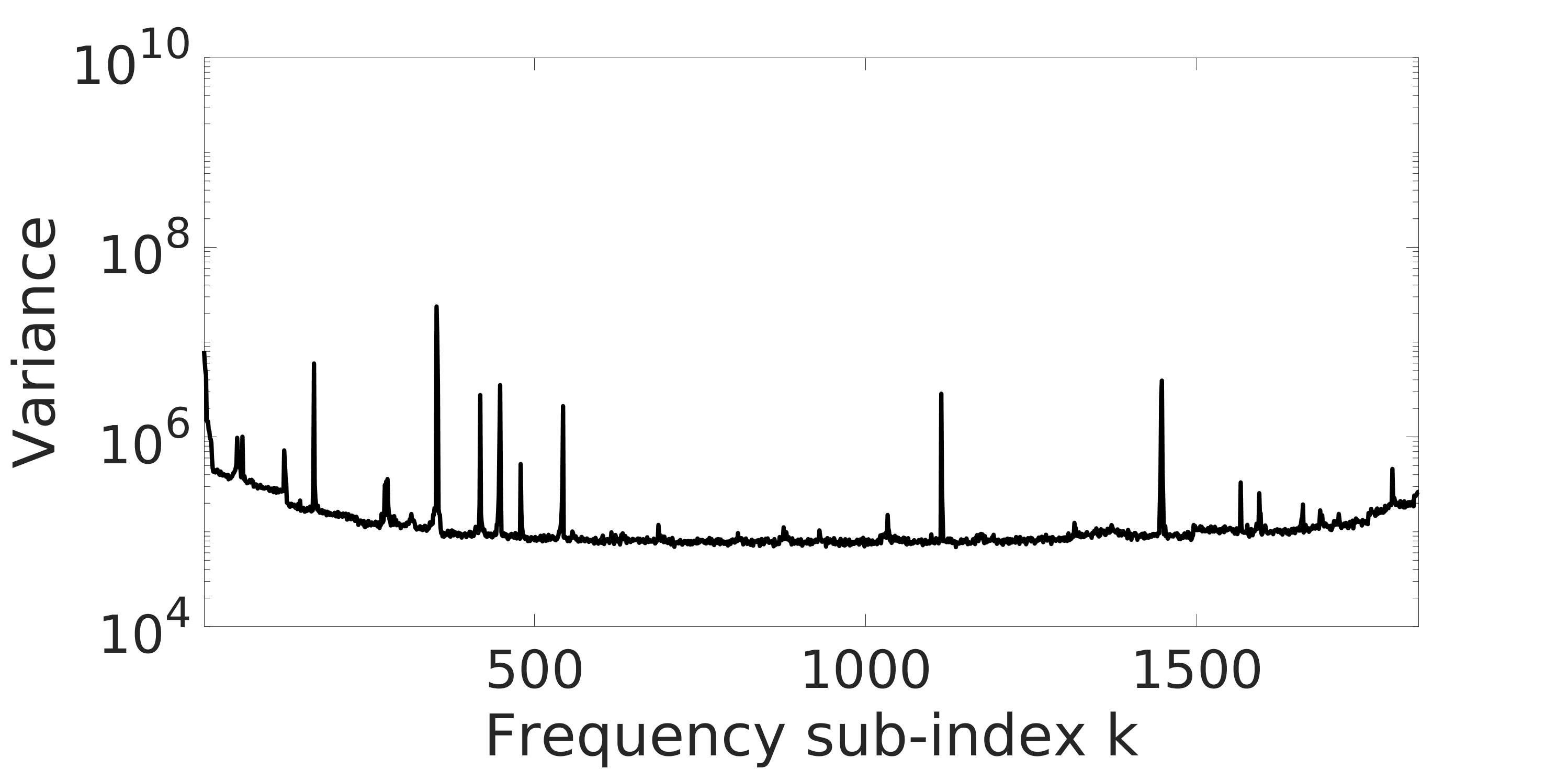} \\
\includegraphics[width=\textwidth]{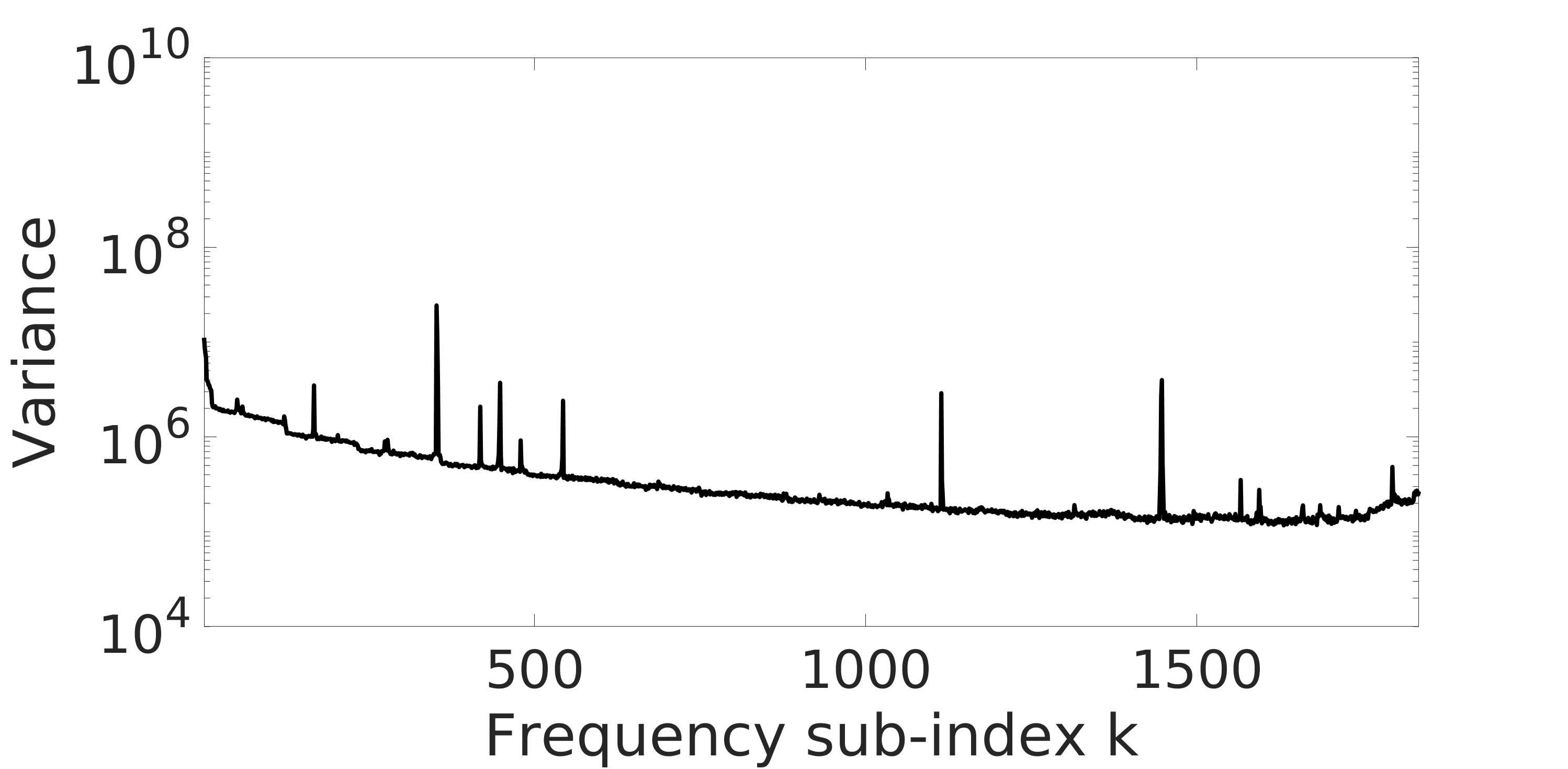} \\
 \end{minipage}
\begin{minipage}{0.3\textwidth}
\centering
$\tau=5$ \\
\includegraphics[width=\textwidth]{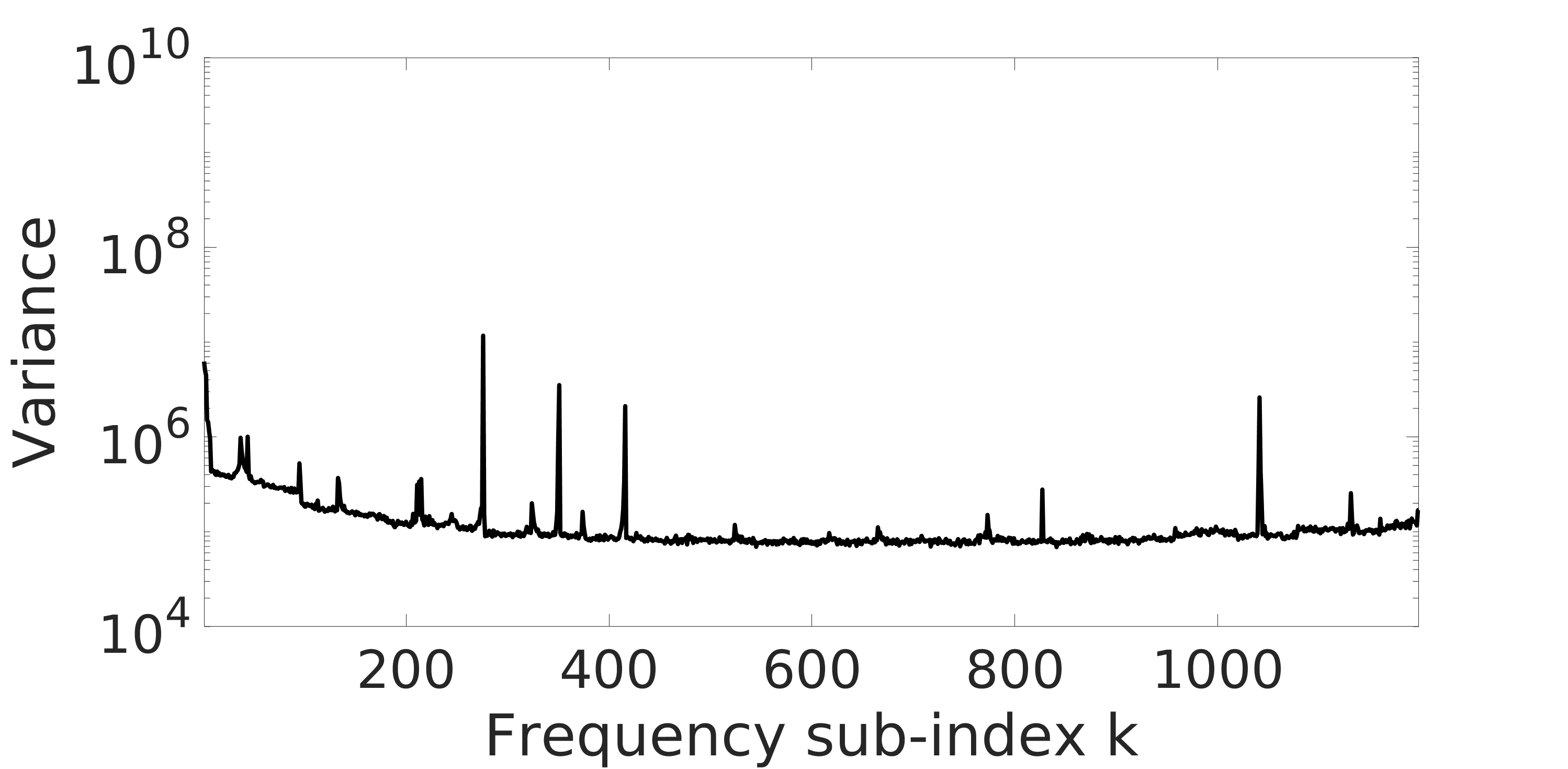} \\
\includegraphics[width=\textwidth]{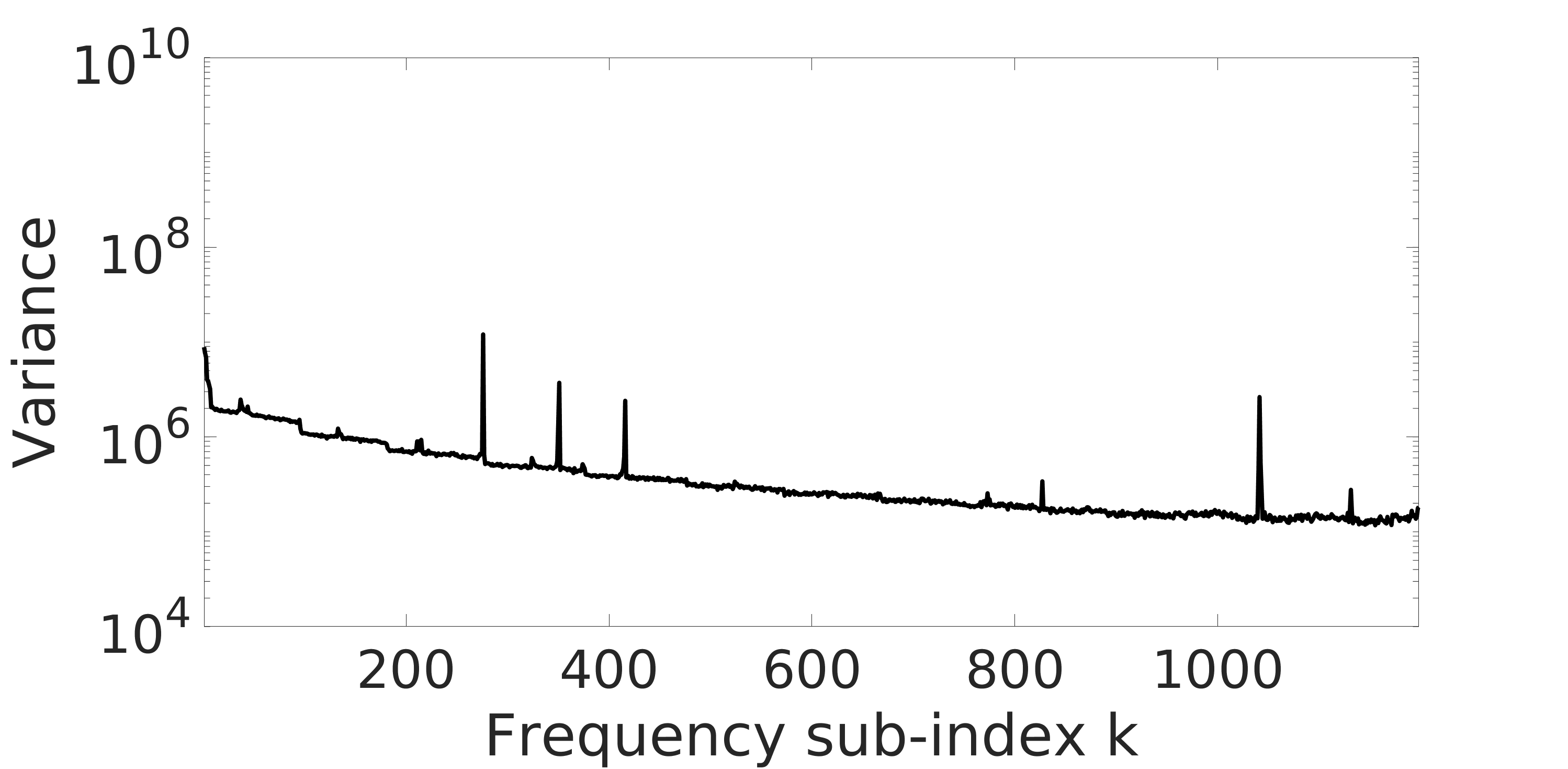} \\
 \end{minipage}
 }
\caption{Variance of the background measurement provided with the ``shape'' phantom from \texttt{Open MPI dataset},
computed from 1000 empty scans when using a 3D excitation in the preclinical Bruker MPI system. Visualized
for the receive coil in $x$-direction for different thresholds $\tau$ versus
frequency indices $j_k\in J_1$; real part (top), imaginary part (bottom).}\label{fig:cov_diag_tau_xcoil}
\end{figure*}

The now standard preprocessing approaches to treat the noise is frequency selection, including band pass approach and SNR-type thresholding.
The description of these approaches below largely follows \cite{KluthJin:2019}. Let $J_\mathrm{BP}=\{j\in \mathbb{Z}|\ b_1 \leq |j|/T\leq
b_2 \}$ be the indices for frequency band limits $0\leq b_1 < b_2 \leq \infty$ and measurement time $T$. This step is to further filter out
remaining signal contributions of the analogously filtered direct feedthrough induced by the applied magnetic field. For SNR-type thresholding,
one standard quality measure is the ratio of mean absolute values from individual measurements $v_\ell^{(i)}$ (for the $i$-th calibration
scan at the $\ell$-th receive coil) and empty scanner measurements $\{ v_{\ell,0}^{(k)} \}_{k=1}^K$ \cite{Franke:2016}. Specifically, let
$I_\mathrm{SNR} \subset \{1,\hdots,N\}$ be the index set of individual measurements. Let $\{\psi_j\}_{j\in\mathbb{N}}$ be an orthonormal basis,
e.g., discrete Fourier basis, for $L^2(I)$, where $I$ is the time interval for measurement. Then we define
\begin{equation}
 d_{\ell,j}= \frac{\frac{1}{{|I_\mathrm{SNR} |}}\sum_{{ i\in I_\mathrm{SNR} }} |\langle v_\ell^{(i)} -\mu_\ell^{(i)}, \psi_j \rangle
 |}{\frac{1}{K} \sum_{k=1}^K |\langle v_{\ell,0}^{(k)} - \mu_\ell, \psi_j \rangle |},
\end{equation}
where $\mu_\ell=\frac{1}{K} \sum_{k=1}^K v_{\ell,0}^{(k)}$ is the mean background measurement, and
$\mu_\ell^{(i)}=\kappa_i v_{\ell,0}^{(k_i)} + (1-\kappa_i) v_{\ell,0}^{(k_i+1)}$ is a convex
combination of the $k_i$-th and $k_i+1$-th empty scanner measurements
for the $i$-th calibration scan. The parameters $\kappa_i \in [0,1]$ are chosen to be
equidistant for all calibration scans between two consecutive empty scanner measurements.
That is, if there are $Q$ calibration measurements between the $k_i$-th and the $k_i+1$-th
empty scanner measurement, then $\kappa_i\in \{0, \frac{1}{Q-1}, \frac{2}{Q-1}, \hdots,1\}$.
For a given threshold $\tau \geq 0$, we define
\begin{equation}
 J_\ell=\{ j \in J_\mathrm{BP} |  d_{\ell,j}\geq \tau \},\quad \ \ell =1,\hdots,L,
\end{equation}
which comprises all frequency indices within a certain frequency band and fulfilling an SNR-type measure for the $\ell$-th
receive coil. The threshold $\tau$ determines the size of the reduced system and its accuracy: with a large $\tau$, the
procedure is more conservative but may erroneously remove informative data, whereas with a small
$\tau$, it may risk including highly corrupted data points. In Fig. \ref{fig:cov_diag_tau_xcoil},
we present the SNR-type frequency selection with three thresholds.
With a proper $\tau$, the number of outliers is reduced, but not completely removed,
even for $\tau=5$, which may still greatly influence the reconstruction.
Note the logarithmic scale on the vertical axis.

After applying band passing, thresholding, and splitting real and imaginary part of the Fourier-transformed and background-subtracted measured system matrix, we obtain a (reduced) linear system
\begin{equation}\label{eqn:mpi-inv}
  A x =y,
\end{equation}
where $A$ is the (processed) calibrated system matrix (and thus contains noise) and $y$ is the noisy background-subtracted phantom data.
For a detailed description of the processing chain we refer to \cite[Sec. 2.1]{KluthJin:2019}
The standard reconstruction approach in the MPI literature employs an l2 data fitting, which leads to the following constrained approach
\begin{equation}\label{eqn:standard}
  \min_{x\geq 0} \tfrac{1}{2}\|Ax-y\|^2 + \tfrac{\alpha}{2}\|x\|^2,
\end{equation}
where $\|\cdot\|$ denotes the Euclidean norm, and $\alpha>0$ is the penalty parameter, controlling the tradeoff
between data fitting and penalty \cite{ItoJin:2015}. The constraint $x\geq0$ is interpreted componentwise. Problem \eqref{eqn:standard} is often
minimized by a variant of Kaczmarz method \cite{Knopp2017} which is often used with a small number of iterations to obtain reasonable reconstructions.
Indeed, a fixed small iteration number does not guarantee convergence to a minimizer of the functional, i.e., it is rather a hybrid regularization method including an iterative mechanism with early stopping than a pure variational Tikhonov regularization.
Nevertheless, in practice, dimension reduction techniques (via SNR type criterion or randomized SVD) and
proper weighting \cite{Knopp_Online_2016,KluthJin:2019} may also be incorporated to accelerate and enhance the
reconstructions.

Note that the discussion so far assumes that background subtraction has been carried out so that the noise has a zero
mean, which is also assumed below. This condition is implicit in the standard formulation \eqref{eqn:standard}. If
the mean of the noise is nonzero, then the l2 fitting should incorporate the mean as a drift term. The influence
of background subtraction differs from calibrating the noise statistics, e.g. whitening and heavy-tailed modeling.
In practice, background subtraction is not always direct; see the works \cite{ThemKaul:2016,KluthMaass:2017,StraubSchulz:2018,KnoppGdaniec:2019}
for in-depth study, including joint estimation of the background and foreground.

\subsection{L1 fitting}\label{ssec:Laplace}

Statistically, the formulation \eqref{eqn:standard} assumes an i.i.d. Gaussian noise with zero mean. This is often justified
by appealing to a version of central limit theorems, i.e., the Gaussian is suitable for data that are formed
from the sum of a large number of independent components. A well known limitation of the Gaussian model
formulation is its lack of robustness against outliers, i.e., data points that lie far away from the bulk of the data: A
single aberrant data point can greatly influence all the parameters in the model, even for these with little substantive
connection to the outlying observations \cite[p. 443]{GelmanCarlinSternRubin:2004}.

Thus, the validity of the approach \eqref{eqn:standard} resides on validity of the i.i.d. Gaussian assumption on the noise.
However, not all data in MPI can be adequately described by a Gaussian model. From Fig. \ref{fig:cov_diag}, Gaussianity at
best holds true only for some frequencies, whereas for the others, the data contains a significant amount of error, with
outlier like noise. The precise mechanism for the noise remains largely elusive, and there are multiple sources, related to unmodeled
physics of the experimental process, e.g., imperfect analog filter, direct feedthrough, and unexpected magnetization of
scanner components. In practice, the electronic noise is often assumed to be Gaussian \cite{Weizenecker2007}, which is then weighted with
a transfer function, and also there exist (generally nongaussian) noise artifacts; see \cite{SchmaleGleich:2010}
for a study on the noise sources in the receive chain of an MPI scanner. It is known in image
processing \cite{Bovik:2010}, signal processing \cite{AllineyRuzinsky:1994} and statistics \cite{Huber:1981} that noise
with outliers is more adequately described by heavy-tailed distributions. In the presence of outliers, an inadvertent
adoption of the Gaussian model can seriously compromise the reconstruction accuracy \cite{GelmanCarlinSternRubin:2004},
and often does not allow full extraction of the information provided by the data. This calls for methods that are robust to outliers.

There are several ways to derive robust estimators. One classical approach is to first identify outliers with
noise detectors, e.g., by adaptive median filter and nonlocal mean filter \cite{HwangHaddad:1995,BuadesMorel:2008},
and then to perform inversion on the dataset with outliers excluded \cite{GelmanCarlinSternRubin:2004}.
Frequency selection in Section \ref{ssec:freqsel} is a special noise detector (with an SNR type criterion). This
approach depends on the accuracy of the noise detector. It can be highly nontrivial to accurately
identify all outliers, and misidentification can adversely affect the reconstruction quality. See Fig.
\ref{fig:cov_diag_tau_xcoil} for an illustration.
These observations necessitate developing more systematic strategies for
handling outliers, which can be achieved by modeling them explicitly with a heavy-tailed
distribution, e.g., Laplace, Student $t$  and Cauchy  \cite{FossKorshunov:2013}. Laplace
distribution is one of most popular choices, with its density $p(\xi)$ in one-dimension given by
\begin{equation*}
   p(\xi) = \tfrac{\lambda}{2} e^{-\lambda|\xi-\mu|},
\end{equation*}
where $\mu$ and $\lambda>0$ denote the mean (location) and (inverse) scale, respectively.

The proposed approach is based on an i.i.d. Laplace distribution with zero mean assumption on the noise, so
as to allow outliers in the data. Assuming a Gaussian prior (for the solution $x$) in the Bayesian formalism as in \eqref{eqn:standard}
and then considering the maximum a posteriori estimator lead to
\begin{equation}\label{eqn:l1}
  \min_{x\geq0}\|Ax-y\|_1 + \tfrac\alpha2 \|x\|^2,
\end{equation}
where the notation $\|\cdot\|_1$ denotes the $\ell^1$ norm, i.e.,
\begin{equation*}
  \|z\|_1 = \sum_i|z_i|,
\end{equation*}
and the scalar $\alpha>0$ is the corresponding regularization parameter. In the absence of nonnegativity constraint,
this model was analyzed in \cite{ClasonJinKunisch:2010}.

The difference of \eqref{eqn:l1} from \eqref{eqn:standard} is that it employs the l1 fitting, which is more robust
to outliers, i.e., the outliers influence less the reconstructions, instead of the usual l2 fitting. It partly
avoids the frequency selection step in the two-step procedure, and allows using more systematically the
given data. In passing, one may also employ alternatives, e.g., student $t$ likelihood or Huber's robust
statistics, but they will not be explored below.

\subsection{Numerical algorithm}

The formulation \eqref{eqn:l1} involves solving a convex but nonsmooth constrained optimization problem, and it can
be solved efficiently in several different ways, e.g., iteratively
reweighted least-squares \cite{RodriguezWohlberg:2009}, alternating direction method of multipliers \cite{YangZhang:2011},
semismooth Newton method \cite{ClasonJinKunisch:2010} and limited-memory BFGS. These algorithms are easy
to implement and converge steadily, if relevant tuning parameters are properly chosen. We employ a version of
limited-memory BFGS, i.e., L-BFGS-B \cite{ByrdLu:1995,ZhuByrd:1997}. It can ensure that problem \eqref{eqn:l1} is
solved accurately in the sense of optimization, i.e., finding a near global minimizer, so as to avoid extra
regularizing effect from the optimizer due to early stopping.

Limited-memory BFGS-B is a popular quasi-Newton type method using a limited amount of computer memory for
a differentiable objective function, approximating the inverse Hessian matrix using the BFGS approximation, and
handling the simple box constraint (i.e., upper and lower bounds) by an active set type strategy \cite{ByrdLu:1995}.
It also includes a line search step to safeguard the progress, and speeds up the computation using a compact
representation of the BFGS Hessian approximation. It is well suited for large-scale optimization problems
with simple constraint, and there are several well tested implementations \cite{ZhuByrd:1997} (see
\url{https://github.com/stephenbeckr/L-BFGS-B-C} for a MATLAB wrapper).

Since the l1 fitting is nondifferentiable, we approximate \eqref{eqn:l1} by
\begin{equation*}
 \min_{x\geq 0} \|Ax-y\|_{1,\epsilon} + \tfrac{\alpha}{2}\|x\|^2,
\end{equation*}
where $\epsilon>0$ is small, and $\|\cdot\|_{1,\epsilon}$ is defined by
\begin{equation*}
  \|v\|_{1,\epsilon} = \sum_i\sqrt{v_i^2+\epsilon^2}.
\end{equation*}
Upon smoothing, the objective function is differentiable, and thus the limited memory BFGS-B can be
applied directly. This smoothing is simple and easy to implement.

\section{Numerical experiments}\label{sec:numer}

Now we present numerical results to illustrate the potential and performance of the proposed l1 fitting on real data.
The experimental setup is as follows. We employ a measured system matrix, where a band pass filter is
applied (with $b_1=80$~kHz and $b_2=625$~kHz) and frequency selection (with discrete Fourier basis $\{\psi_j\}_{j\in\N}$) with a SNR threshold $\tau$
is optionally applied, which yields a system matrix $A_\tau \in \R^{n\times m}$ for the $L=3$ receive
channels (see \cite[Sec. 2.1]{KluthJin:2019} for the description). Optionally, $A_\tau$ can also be whitened \cite[Sec. 2.3]{KluthJin:2019}, where background measurements
are used to obtain a diagonal whitening matrix $W_\tau\in \R^{n\times n}$.
System matrices and measurements are concatenated and background-subtracted \cite[Sec. 2.1]{KluthJin:2019}.
For frequency selection, we consider four thresholds, i.e., $\tau=0,1,3,5$, and the corresponding number $n$ of
rows of $A_\tau$ is 70446, 68566, 9564 and 6146.
All forward maps are scaled to have a unit operator norm and phantom measurements $y$ are obtained analogously.

Below we compare results obtained from the following reconstruction methods.
\begin{itemize}
   \item \textbf{[l1-L]}: The l1 fitted reconstructions $x_{\rm l1}$ and $x_{W;\rm l1}$ are respectively obtained by
   \begin{align*}
      x_\mathrm{l1}&=\arg\min_{x\geq 0}  \|A_\tau x-y\|_{1,\epsilon} + \tfrac{\alpha}{2}\|x\|^2 \quad \mbox{and}\quad \notag \\ x_{W;\mathrm{l1}}
   &=\arg\min_{x\geq 0}   \|W_\tau A_{\tau}x-W_\tau y\|_{1,\epsilon} + \tfrac{\alpha}{2}\|x\|^2,
   \end{align*}
    where the minimization is performed with L-BFGS-B.

  \item \textbf{[l2-K, l2-L]}: The reconstructions $x_{\rm l2}$ and $x_{W;\rm l2}$ are respectively obtained by
    \begin{align*}
      x_\mathrm{l2}&=\arg\min_{x\geq 0}\tfrac{1}{2}\|A_\tau x-y\|^2 + \tfrac{\alpha}{2}\|x\|^2 \quad \mbox{and}\quad \notag \\ x_{W;\mathrm{l2}}
   &=\arg\min_{x\geq 0}\tfrac{1}{2}\|W_\tau A_{\tau}x-W_\tau y\|^2 + \tfrac{\alpha}{2}\|x\|^2,
   \end{align*}
   where l2-K denotes minimization by Kaczmarz method (see, e.g., \cite[Algorithm 1]{KluthJin:2019}) with $N$ iterations (i.e., one loop over the entire matrix); respectively l2-L denotes minimization by L-BFGS-B.
\end{itemize}

These methods are evaluated on a public 3D dataset \texttt{open MPI dataset} (downloaded
from \url{https://www.tuhh.de/ibi/research/open-mpi-data.html}, accessed on September 20, 2019) provided
in the MPI Data Format (MDF) \cite{knopp2018mdf}. The system matrix data $\{v_\ell^{(i)}\}_{i=1}^m$,
$\ell=1,2,3$, is obtained using a cuboid sample of size 2 mm $\times$ 2 mm $\times$ 1 mm and a 3D Lissajous-type FFP excitation. The calibration is
carried out with Perimag tracer with a concentration 100 mmol/l.
The field-of-view has a size of 38 mm $\times$ 38 mm $\times$ 19 mm and the sample positions have a distance
of 2 mm in $x$- and $y$-direction and 1 mm in $z$-direction, resulting in $19\times19\times19=6859$ voxels, which
gives the number $m$ of columns in the full matrix $A$. The entries of $A$ are averaged over 1000 repetitions
and empty scanner measurements are performed and averaged every 19 calibration scans. The measurements
are averaged over 1000 repetitions of the excitation sequence, and with each phantom, an empty measurement
with 1000 repetitions is provided, which are used for the background correction of the measurement and
$A$ \cite[Sec. 2.1]{KluthJin:2019} and also for approximating the diagonal covariance $C$
respectively the whitening matrix $W$ \cite[Sec. 2.3]{KluthJin:2019}. For the comparison below,
the Kaczmarz method \cite[Algorithm 1]{KluthJin:2019} is run up to 200 iterations (one iteration means
one loop over the entire matrix). The L-BFGS-B algorithm in l1-L and l2-L is used with 20
limited-memory vectors, 1e-10 for pgtol (tolerance for the $\ell^\infty$ norm of the
projected gradient), and 10000 for maximum number of iterations, and for l1-L, $\epsilon=10^{-12}$ is chosen.

\begin{figure}[hbt!]
\centering
\scalebox{0.7}{
\begin{tabular}{ccc}
\includegraphics[height=5cm]{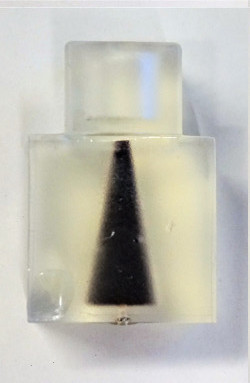}&
\includegraphics[height=5cm]{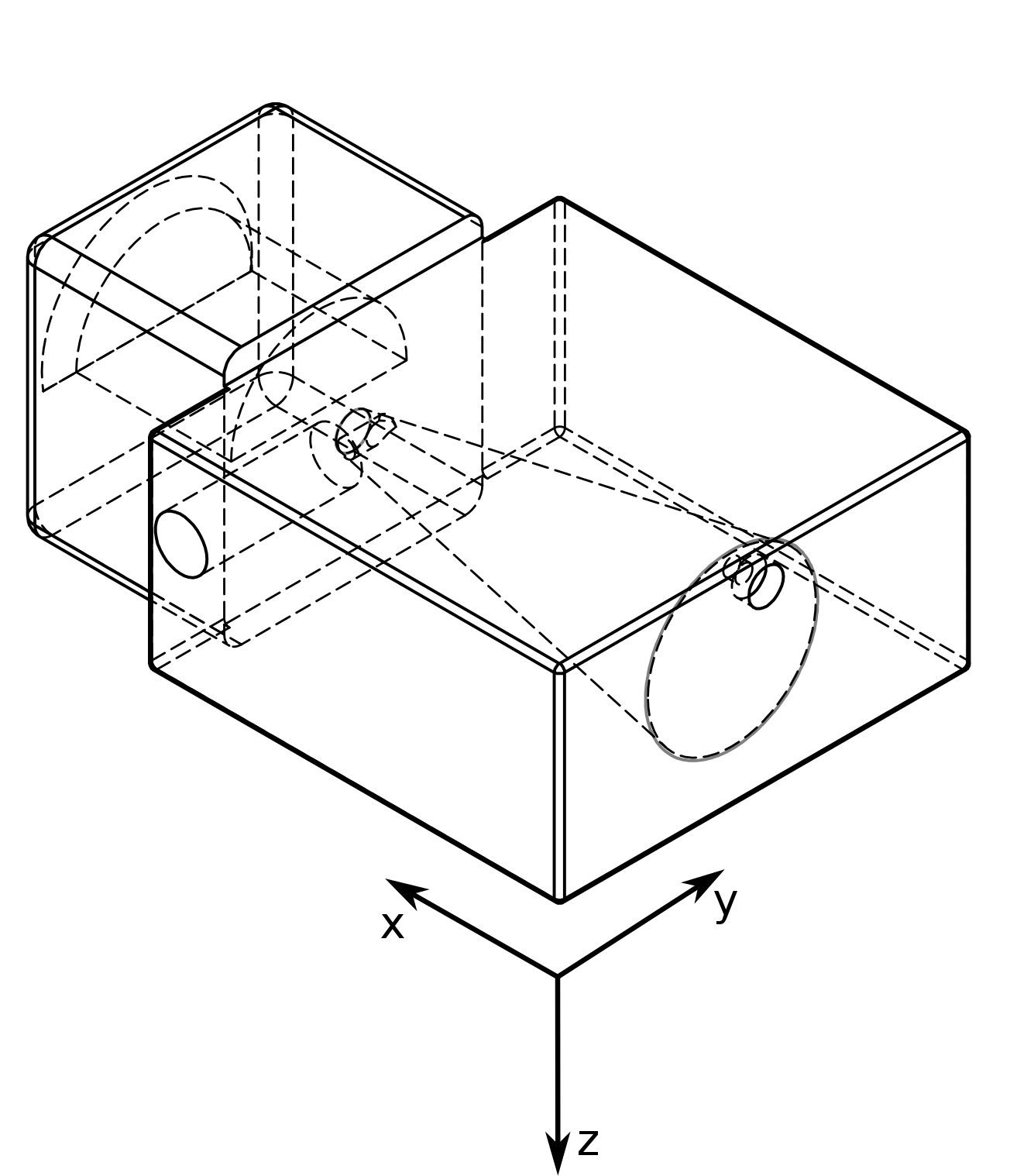}&
\includegraphics[height=5cm]{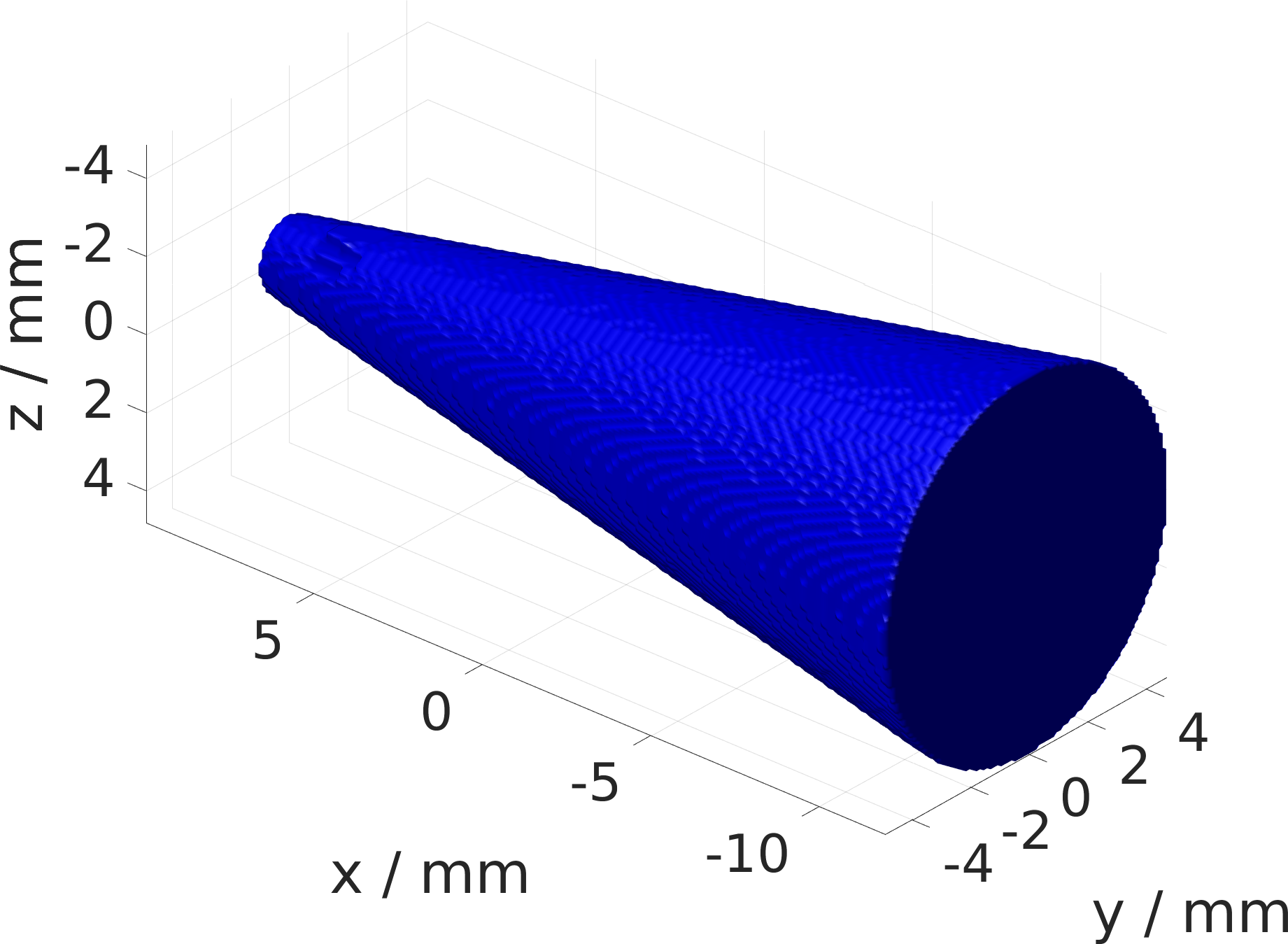}\\
(a) photo & (b) CAD drawing &(c) 3D voxel image\\
\end{tabular}
}
\caption{``Shape'' phantom from the \texttt{open MPI dataset}.}\label{fig:phantom_shape}
\end{figure}

\begin{figure}%
\centering
\scalebox{0.7}{
\begin{tabular}{ccc}
\includegraphics[height=4cm]{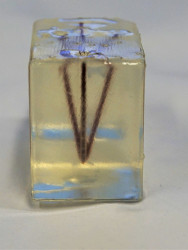} &
\includegraphics[height=4cm]{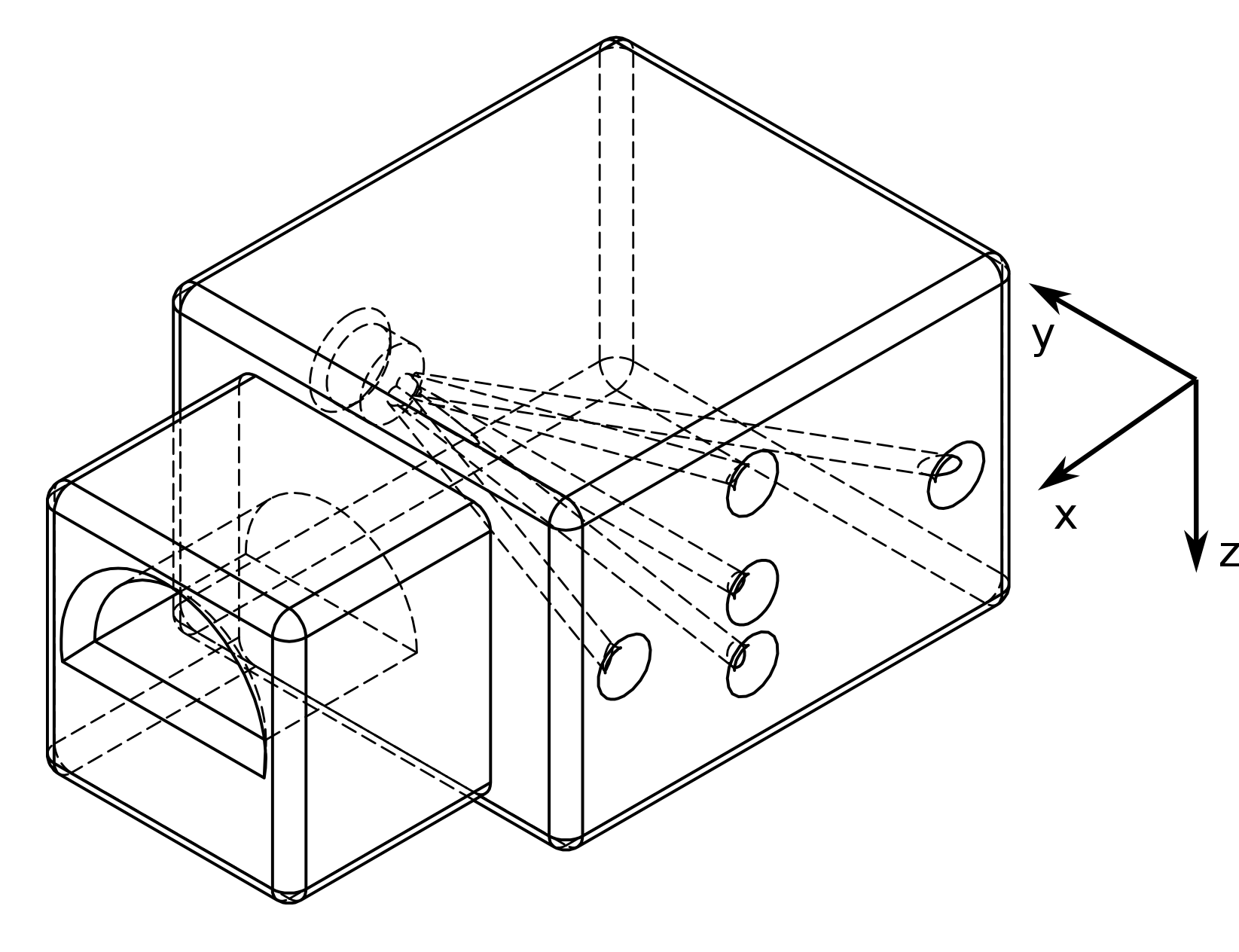}
& \includegraphics[height=4cm]{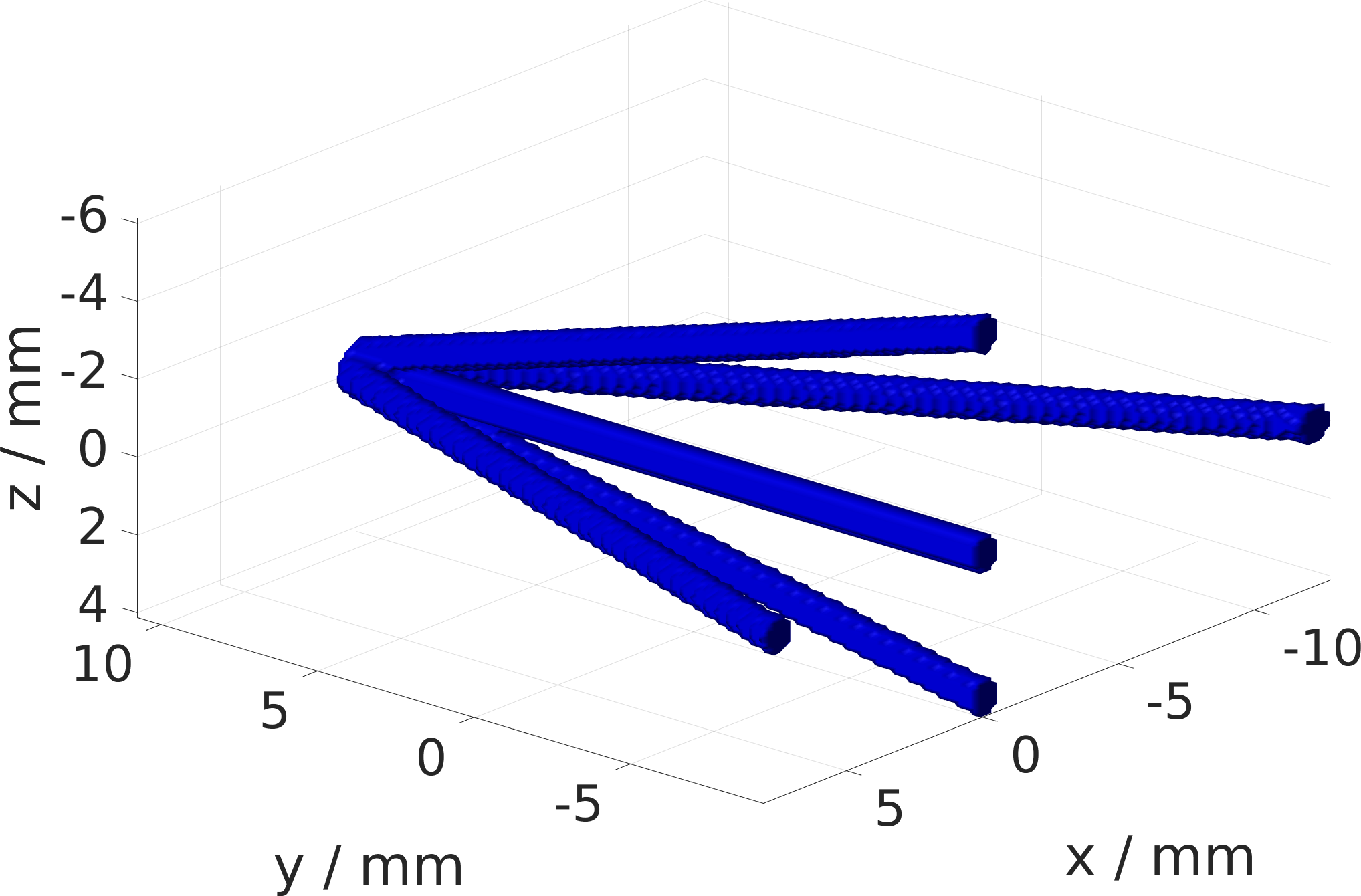}\\
(a) photo & (b) CAD drawing & (c) 3D voxel image\\
\end{tabular}
}
\caption{``Resolution'' phantom from the \texttt{open MPI dataset}.  }\label{fig:phantom_resolution}
\end{figure}

\begin{figure}[hbt!]
\centering
\scalebox{0.75}{
\begin{tabular}{ccc}
\includegraphics[height=5cm]{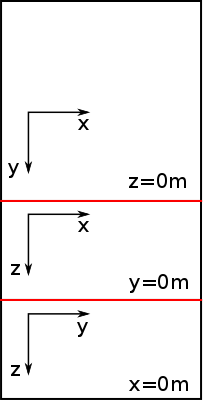}&
\includegraphics[height=5cm]{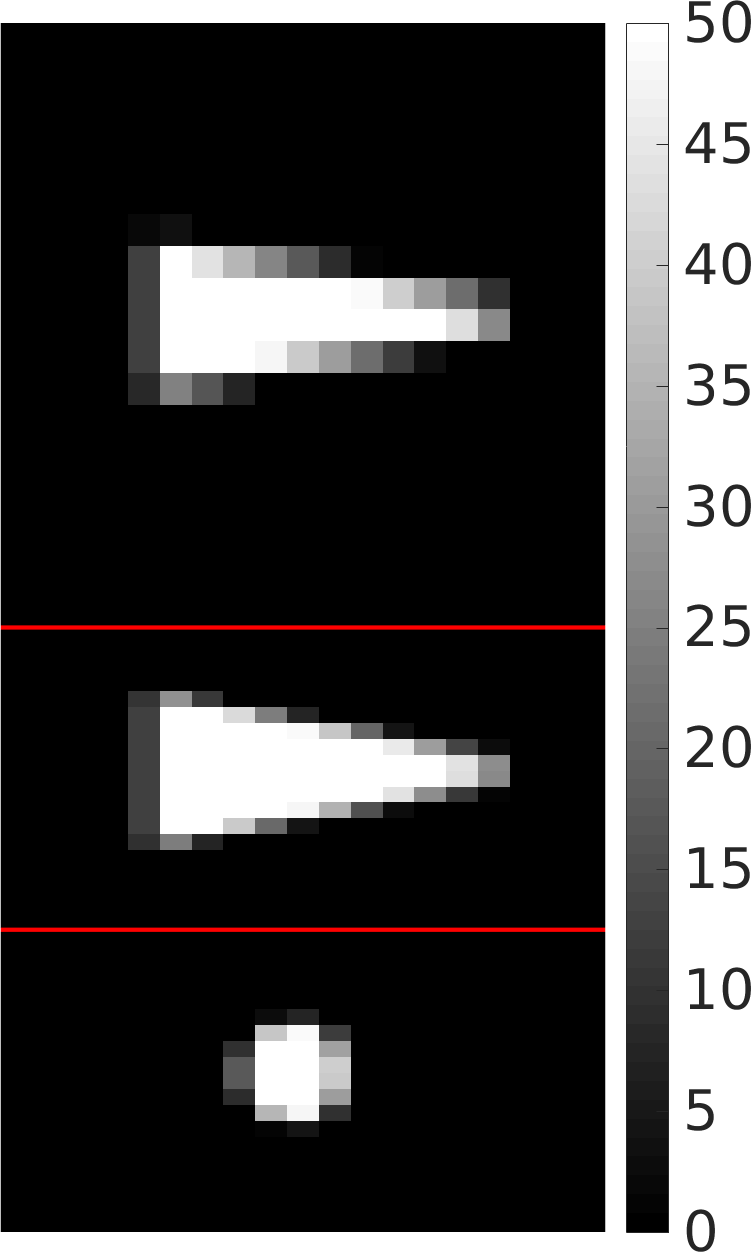}&
\includegraphics[height=5cm]{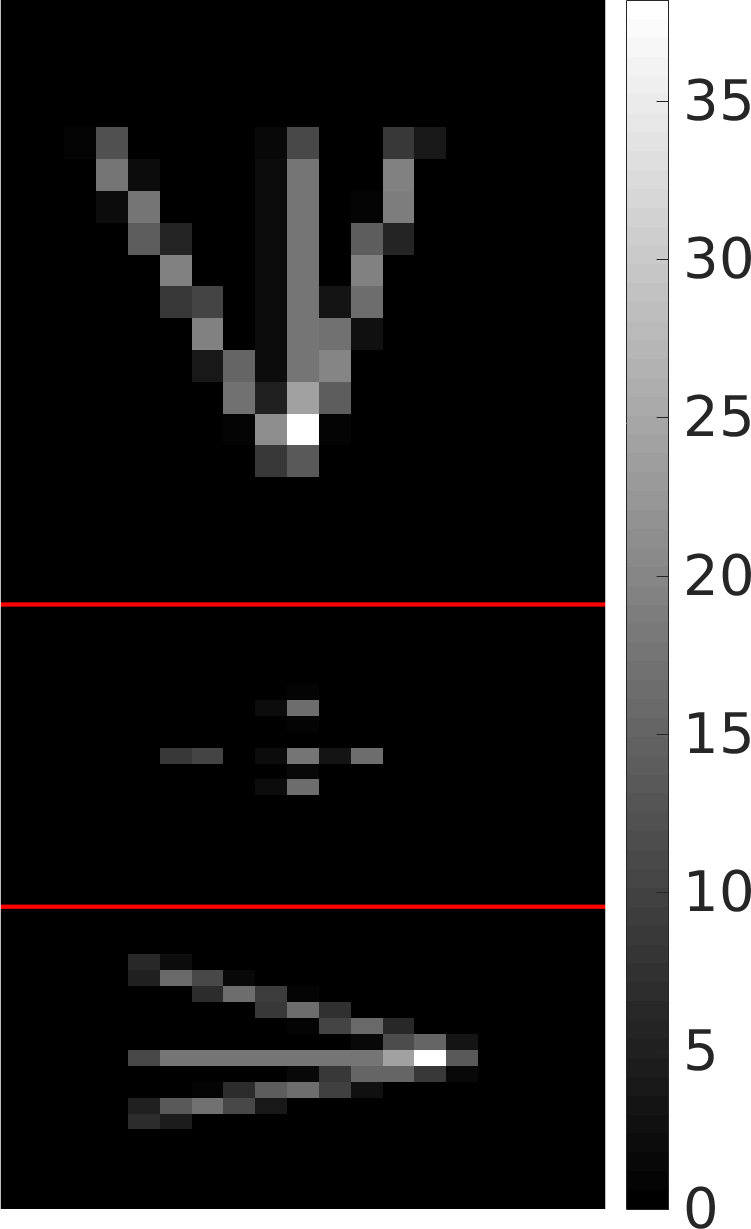}\\
\end{tabular}
}
\caption{Visualization structure
for the 3D reconstructions (left) and the ground truth ``shape'' (middle) and ``resolution'' (right) phantoms on the spatial grid. %
}\label{fig:phantom_ground_truth}
\end{figure}

We validate the proposed method on the ``shape'' and ``resolution'' phantoms in the dataset. The ``shape''
phantom is a cone defined by a 1 mm radius tip, an apex angle of 10 degree, and a height of 22 mm. The total
volume is 683.9 $\mu$l. Perimag tracer with a concentration of 50 mmol/l is used. See Fig. \ref{fig:phantom_shape}
for a schematic illustration, where the plots are adapted from \texttt{Open MPI dataset}.
The ``resolution'' phantom consists of 5 tubes filled with Perimag
tracer with a concentration of 50 mmol/l. The 5 tubes have a common origin on one side of the phantom,
and extend in different angles from the origin within the $x$-$y$- and $y$-$z$-planes.
In the $z$-direction, the angles in the $y$-$z$-plane are chosen smaller (10 deg and 15 deg) than in
$x$-$y$-plane (20 deg and 30 deg); see Fig. \ref{fig:phantom_resolution} for the illustration. In all the
reconstructions below, the concentration unit is mmol/l. See Fig. \ref{fig:phantom_ground_truth} for the visualization
structure of the 3D reconstructions below.

In the \texttt{Open MPI dataset}, CAD drawings of the phantoms are provided; see Figs. \ref{fig:phantom_shape} and
\ref{fig:phantom_resolution}.
From these drawings we extracted the support $\Gamma\subset\mathbb{R}^3$ of the the respective phantom.
Together with the known tracer concentration $\tilde{c}_0 > 0$, it allows extracting voxel images as a ground truth reference. However, there is
uncertainty with the actual phantom position, and also the robot arm moving the phantom in the bore has an unknown
standard deviation. Thus, we formulate an image quality measure as follows. First we estimate the
phantom position from the reconstructed image to define a reasonable estimated reference image $c_0: \mathbb{R}^3
\to \mathbb{R}_+$. Then we account for position uncertainty by shifts $\Delta r \in
\mathcal{R}$, where $\mathcal{R}$ is the set of all possible position shifts in the neighborhood $[-3\text{mm},3\text{mm}]^3$
with a step size $0.5$mm in each direction ($|\mathcal{R}|=2197$). The $L^2$ scalar products of
$c_0(r+\Delta r)=\tilde{c}_0 \chi_{\Gamma+\Delta r}(r)$ and the piecewise constant basis functions $\{\psi_j\}_j$ with respect to the
voxel grid yields the reference image $x_{\mathrm{ref},\Delta r}$, which
is then used to define the following uncertainty-aware image quality measures
\begin{align}
 \epsilon_\mathrm{PSNR}(x)&= \max_{\Delta r \in \mathcal{R}} \mathrm{PSNR}(x,x_{\mathrm{ref},\Delta r})\\
 \epsilon_\mathrm{SSIM}(x)&= \max_{\Delta r \in \mathcal{R}} \mathrm{SSIM}(x,x_{\mathrm{ref},\Delta r})
\end{align}
exploiting the standard image quality measures, i.e., peak-signal-to-noise-ratio (PSNR) and structural similarity
measure (SSIM) \cite{HoreZiou:2010}. Two example reference images are illustrated in Fig. \ref{fig:phantom_ground_truth}.
These two metrics are used for quantitative comparisons below.

\begin{table*}[hbt!]
\setlength{\tabcolsep}{1pt}
\centering
\caption{The $\epsilon_\mathrm{PSNR}$ values for l1-L, l2-L, and l2-K.
The numbers in brackets refer to $\alpha$, respectively $\alpha$ and the iteration number $N$ for l2-K. %
}\label{tab:psnr_nonwhitened_vs_whitened}
\scalebox{0.95}{
\begin{tabular}{c|c|c|c||c|c|c}
&\multicolumn{6}{c}{``Shape'' phantom} \\
\hline
& \multicolumn{3}{c||}{non-whitened} & \multicolumn{3}{c}{whitened} \\
\hline
  $\tau$ & l1-L & l2-L & l2-K & l1-L & l2-L & l2-K \\
\hline
0 &  $ 19.687 \ (2^{-5}) $  &  $ 19.848 \ (2^{-7}) $  &  $ \mathbf{28.615} \ (2^{-17}, 2) $  &  $ 19.379 \ (2^{-5}) $  &  $ 20.139 \ (2^{-7}) $  &  $ \mathbf{29.430} \ (2^{-16}, 2) $  \\
1 &  $ 23.997 \ (2^{-4}) $  &  $ 21.327 \ (2^{-6}) $  &  $ \mathbf{29.075} \ (2^{-15}, 2) $  &  $ 24.240 \ (2^{-2}) $  &  $ 24.475 \ (2^{-6}) $  &  $ \mathbf{29.866} \ (2^{-13}, 2) $  \\
3 &  $ 27.738 \ (2^{-2}) $  &  $ 25.305 \ (2^{-8}) $  &  $ \mathbf{29.233} \ (2^{-14}, 2) $  &  $ 27.888 \ (2^{-1}) $  &  $ 26.152 \ (2^{-7}) $  &  $ \mathbf{29.702} \ (2^{-13}, 2) $  \\
5 &  $ 27.669 \ (2^{-3}) $  &  $ 25.680 \ (2^{-8}) $  &  $ \mathbf{28.907} \ (2^{-14}, 2) $  &  $ 27.616 \ (2^{-2}) $  &  $ 26.443 \ (2^{-7}) $  &  $ \mathbf{29.393} \ (2^{-13}, 2) $  \\
\hline

\end{tabular}}

\vspace{0.5cm}

\scalebox{0.95}{
\begin{tabular}{c|c|c|c||c|c|c}
&\multicolumn{6}{c}{``Resolution'' phantom} \\
\hline
& \multicolumn{3}{c||}{non-whitened} & \multicolumn{3}{c}{whitened} \\
\hline
  $\tau$ & l1-L & l2-L & l2-K & l1-L & l2-L & l2-K \\
\hline
0 &  $ 29.713 \ (2^{-3}) $  &  $ 29.488 \ (2^{-11}) $  &  $ \mathbf{31.673} \ (2^{-18}, 1) $  &  $ 29.812 \ (2^{-3}) $  &  $ 29.512 \ (2^{-10}) $  &  $ \mathbf{31.880} \ (2^{-17}, 1) $  \\
1 &  $ 30.765 \ (2^{-2}) $  &  $ 29.534 \ (2^{-10}) $  &  $ \mathbf{31.812} \ (2^{-18}, 37) $  &  $ 30.990 \ (2^{-2}) $  &  $ 30.248 \ (2^{-9}) $  &  $ \mathbf{32.419} \ (2^{-17}, 84) $  \\
3 &  $ 31.634 \ (2^{-2}) $  &  $ 30.426 \ (2^{-11}) $  &  $ \mathbf{31.908} \ (2^{-17}, 27) $  &  $ 31.707 \ (2^{-1}) $  &  $ 30.875 \ (2^{-9}) $  &  $ \mathbf{32.160} \ (2^{-16}, 26) $  \\
5 &  $ 31.510 \ (2^{-3}) $  &  $ 30.432 \ (2^{-11}) $  &  $ \mathbf{32.152} \ (2^{-18}, 56) $  &  $ 31.544 \ (2^{-2}) $  &  $ 31.243 \ (2^{-12}) $  &  $ \mathbf{32.134} \ (2^{-16}, 53) $  \\
\hline
\end{tabular}}
\end{table*}

\begin{figure*}[hbt!]
\centering
\scalebox{0.85}{
\begin{tabular}{ccc|ccc}
\multicolumn{3}{c|}{non-whitened} & \multicolumn{3}{c}{whitened} \\
\hline
l1-L & l2-L & l2-K & l1-L & l2-L & l2-K \\
\hline
\multicolumn{6}{l}{$\tau=0$} \\
 \includegraphics[height=3.4cm]{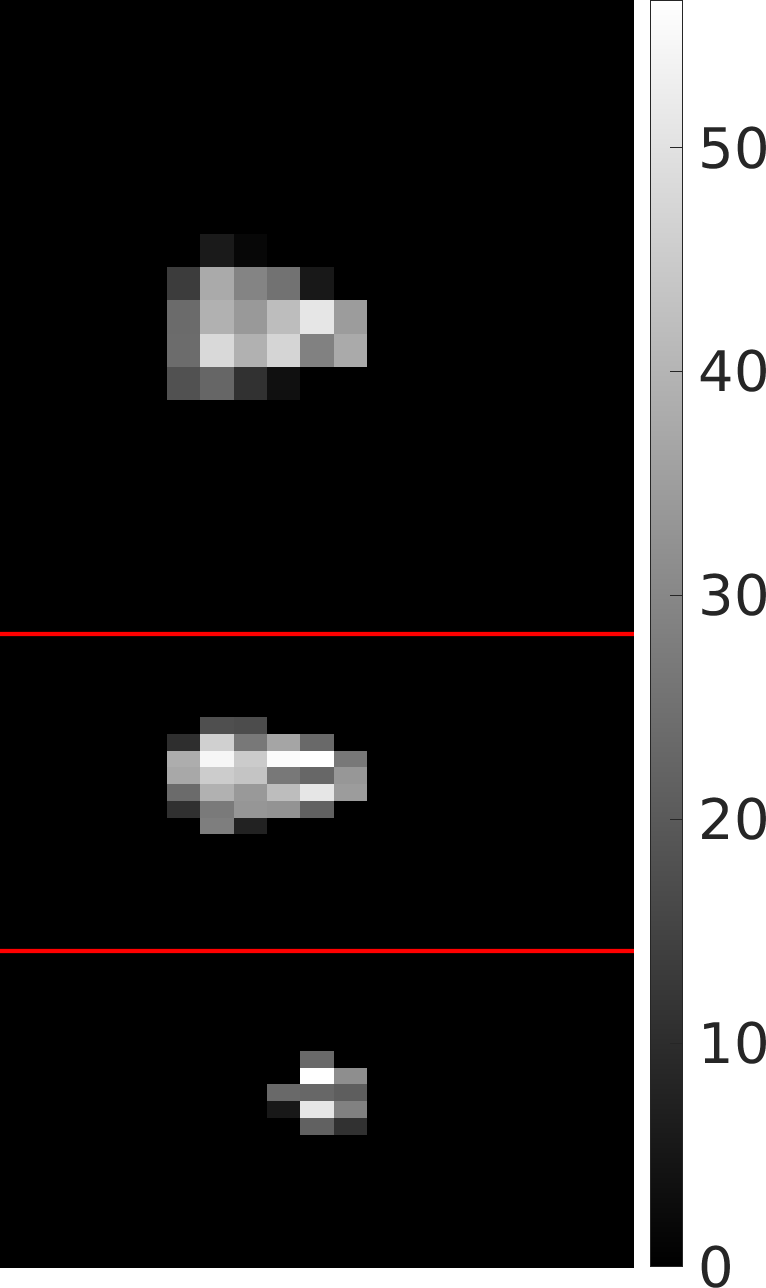}&
 \includegraphics[height=3.4cm]{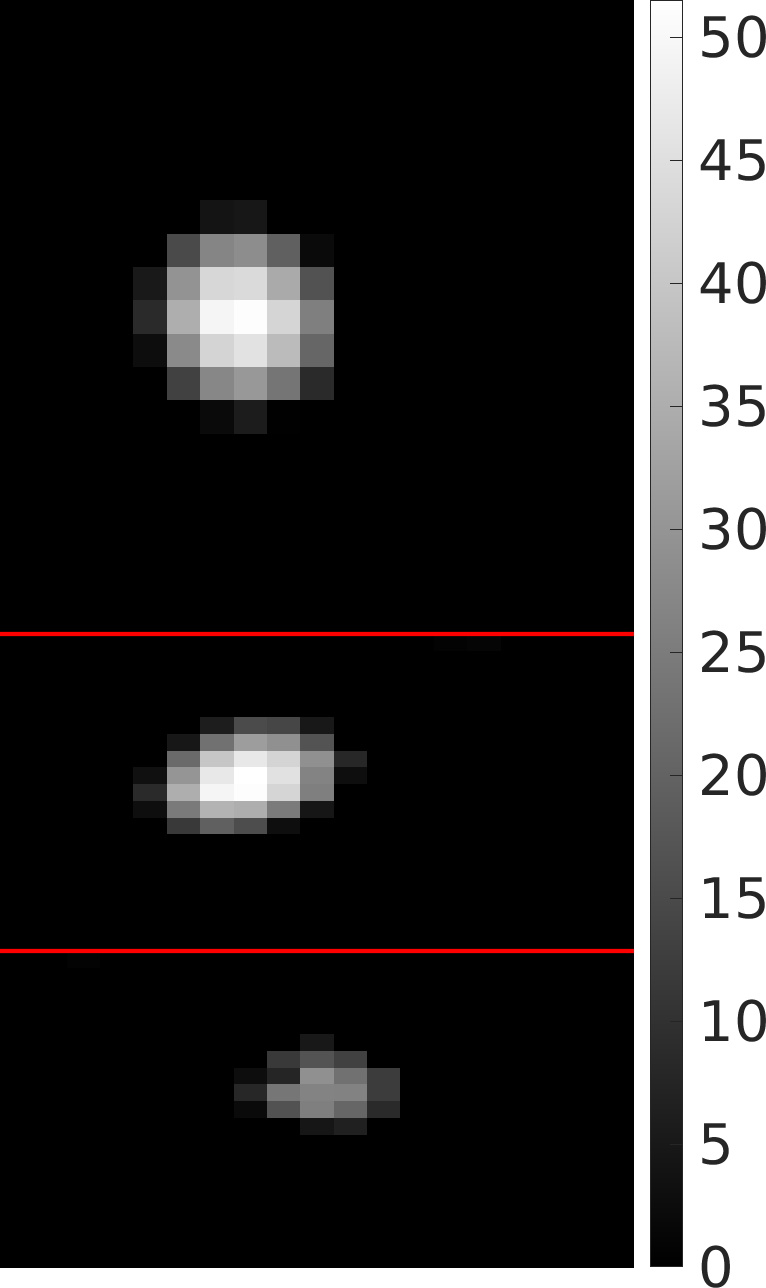}&
 \includegraphics[height=3.4cm]{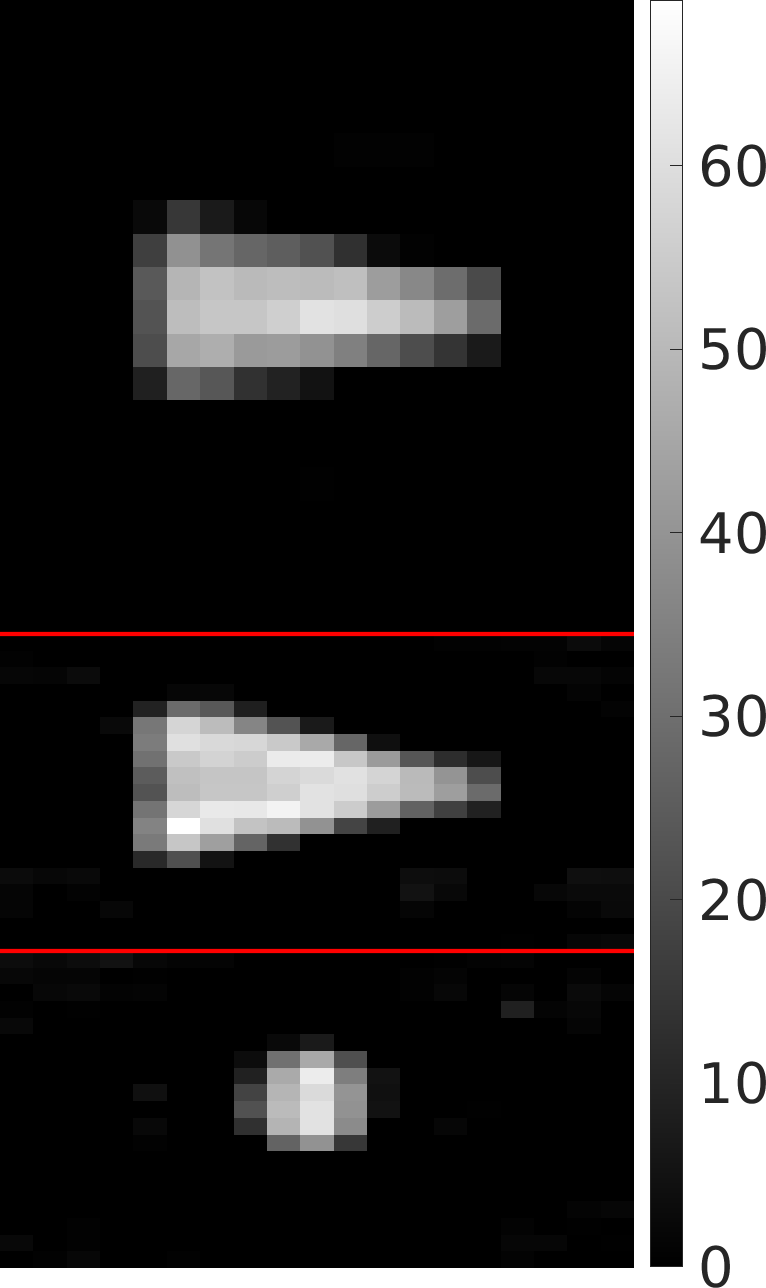}&
 \includegraphics[height=3.4cm]{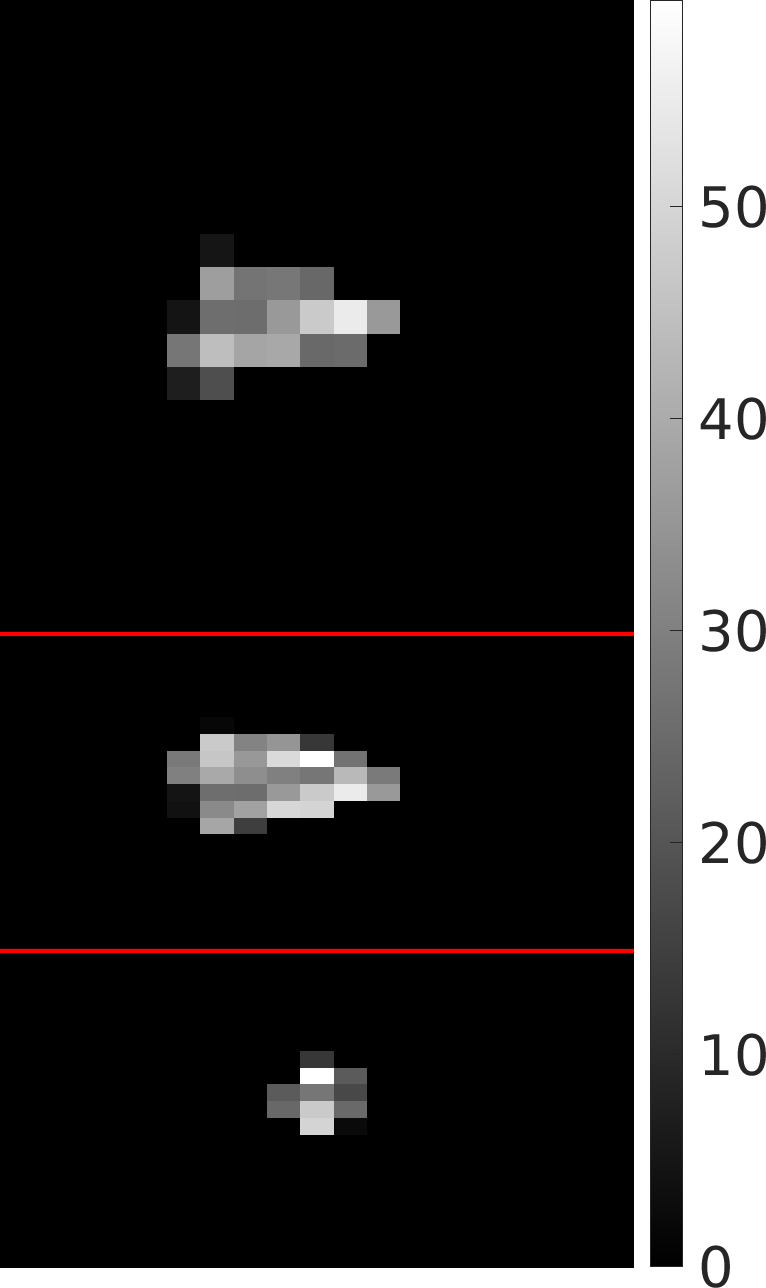}&
 \includegraphics[height=3.4cm]{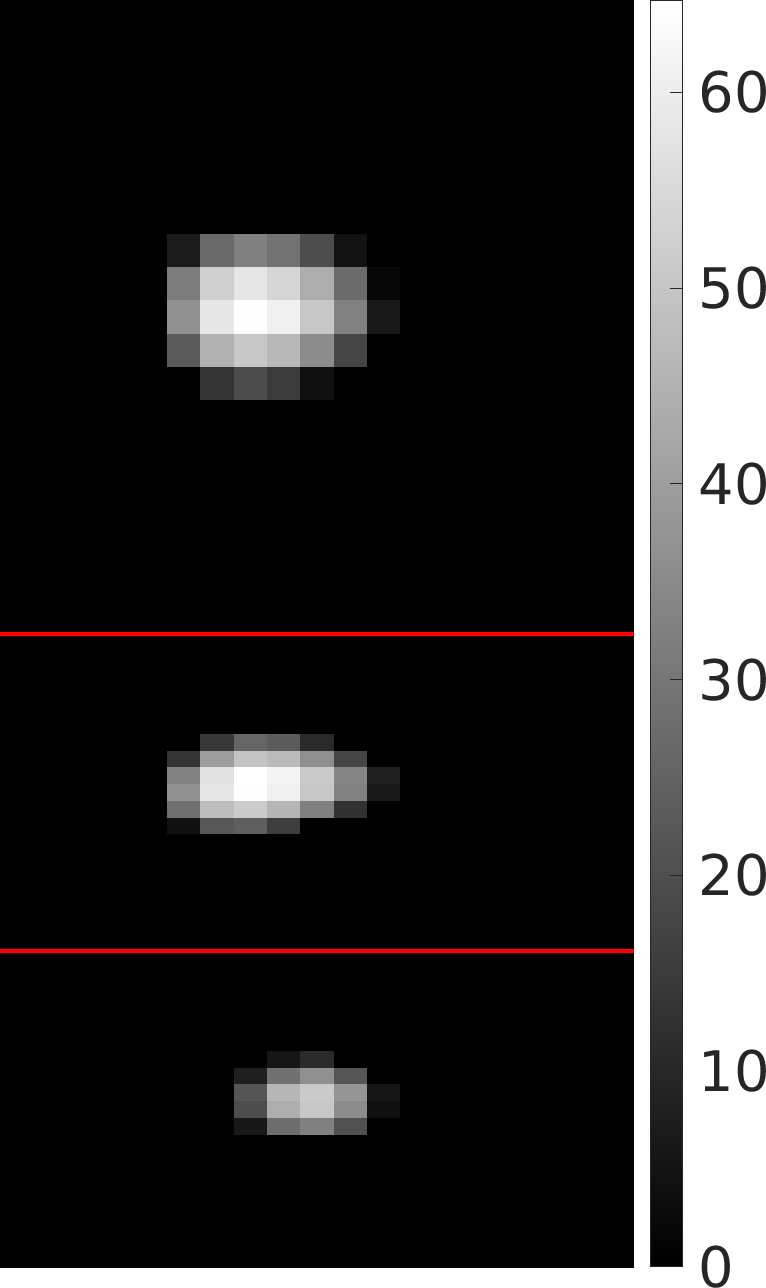}&
 \includegraphics[height=3.4cm]{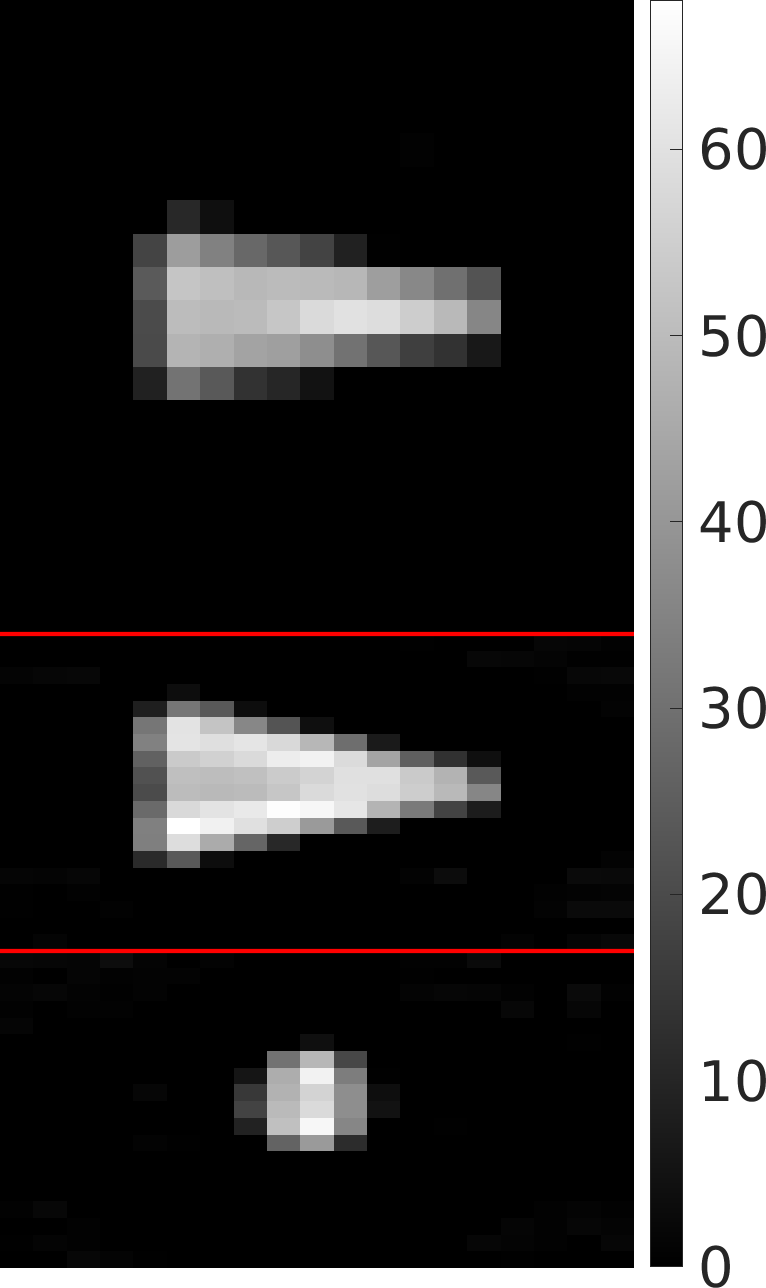}\\
\hline
\multicolumn{6}{l}{$\tau=1$} \\
 \includegraphics[height=3.4cm]{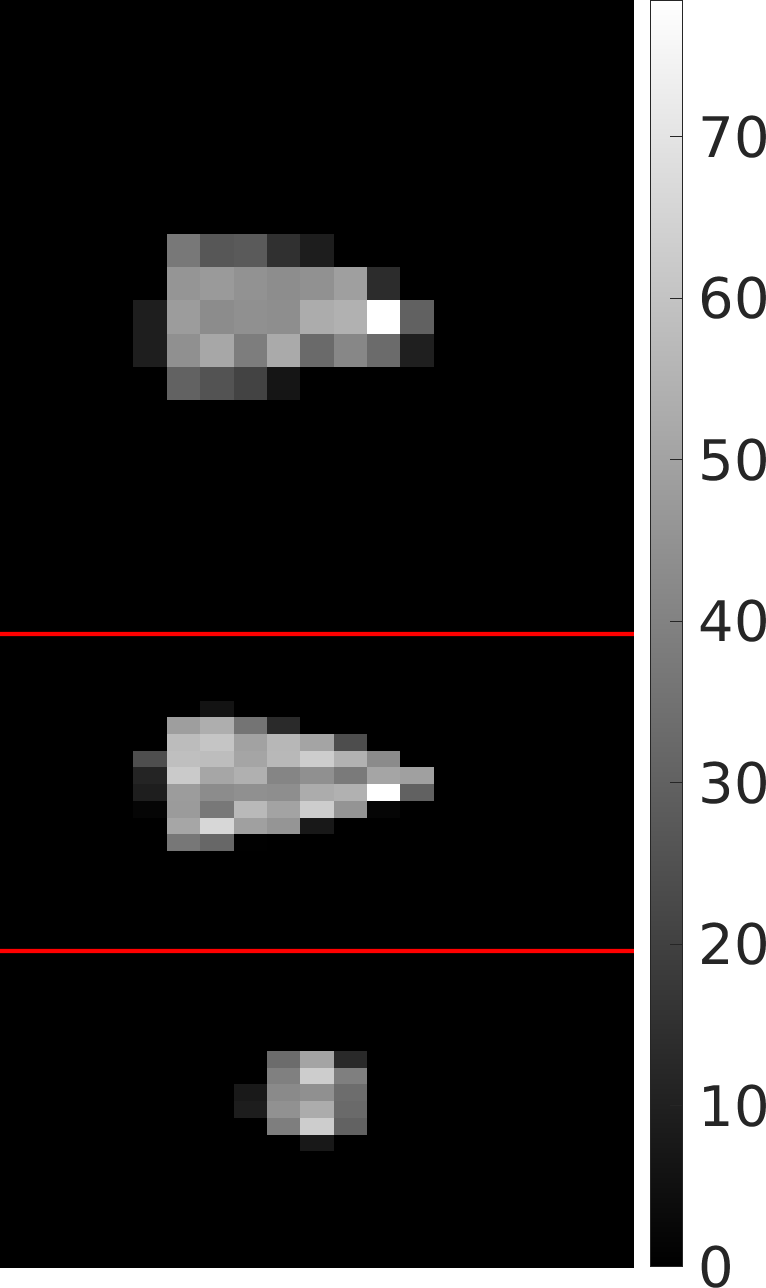}&
 \includegraphics[height=3.4cm]{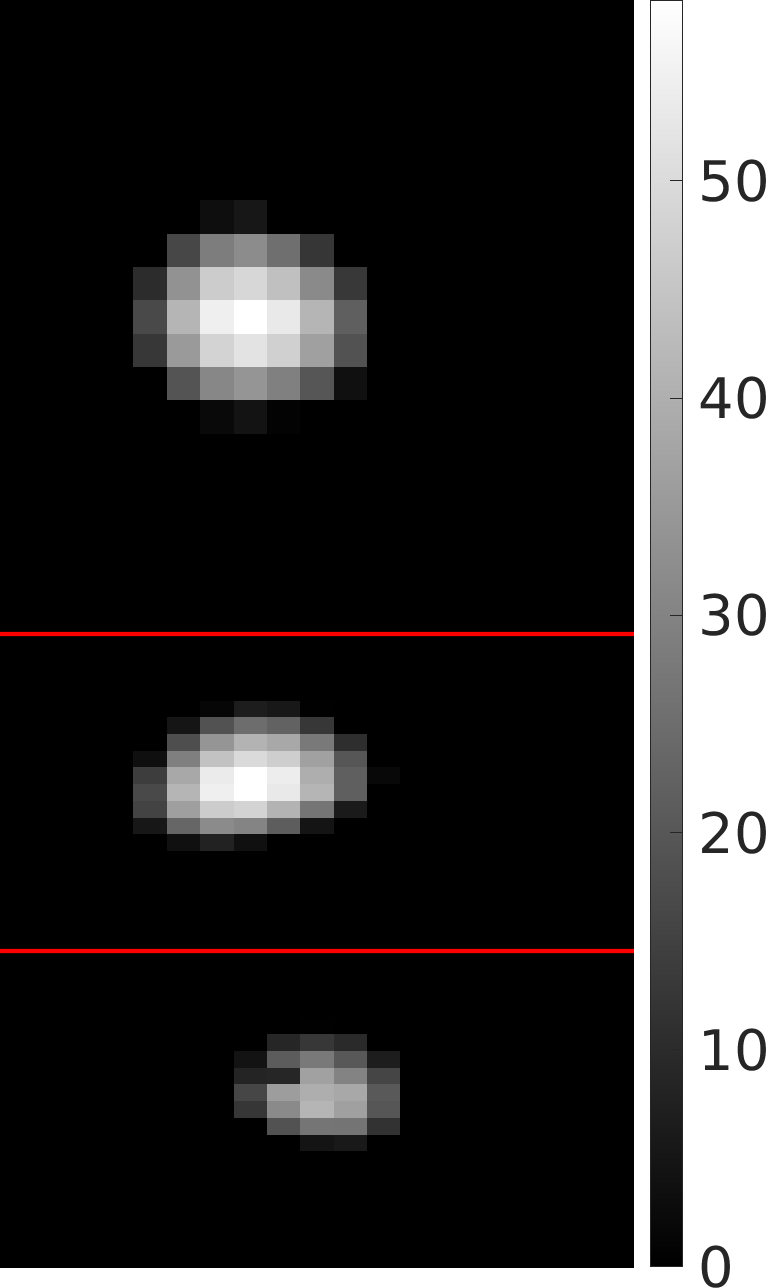}&
 \includegraphics[height=3.4cm]{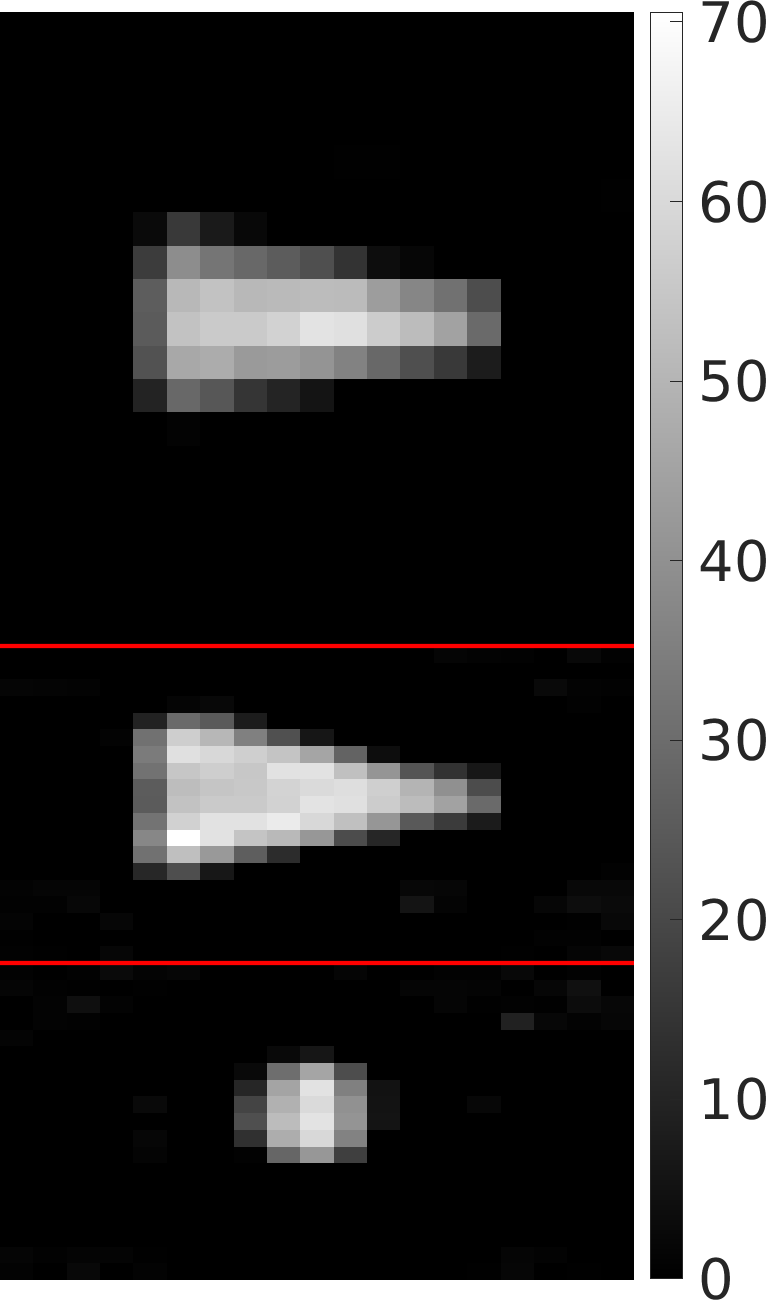}&
 \includegraphics[height=3.4cm]{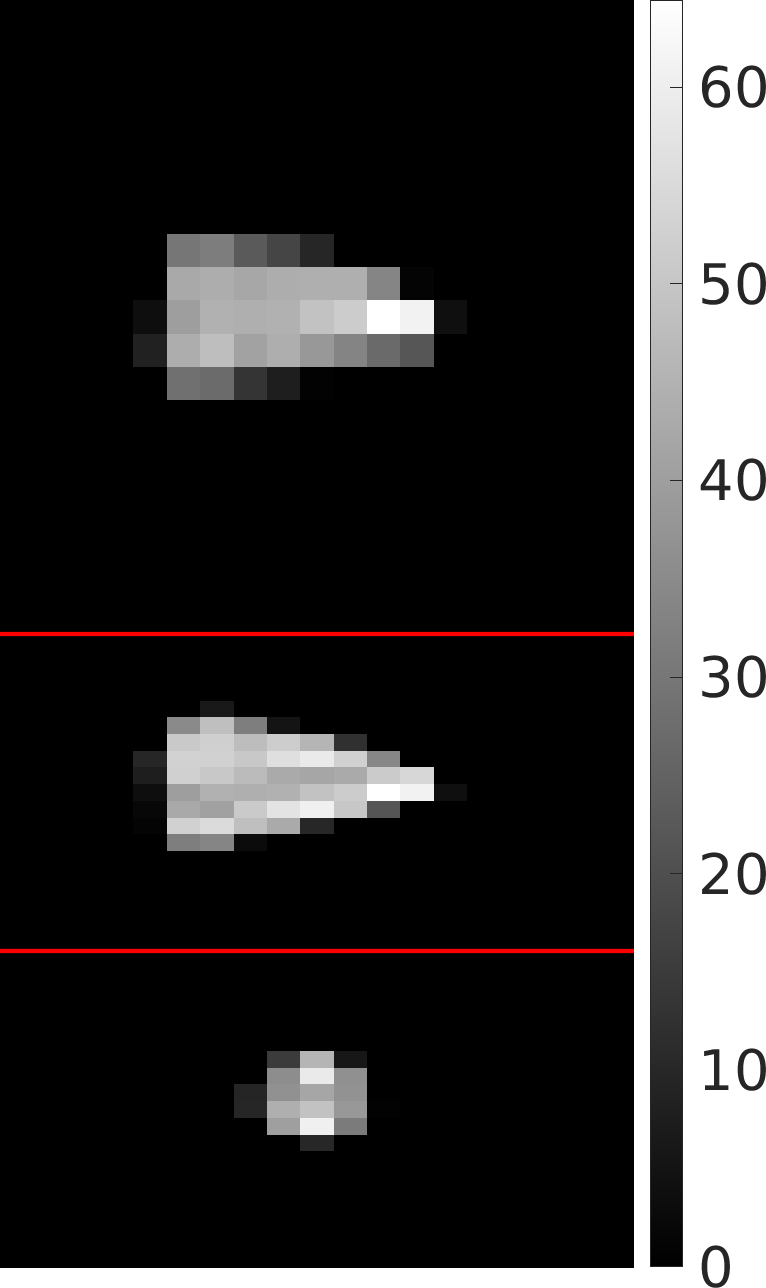}&
 \includegraphics[height=3.4cm]{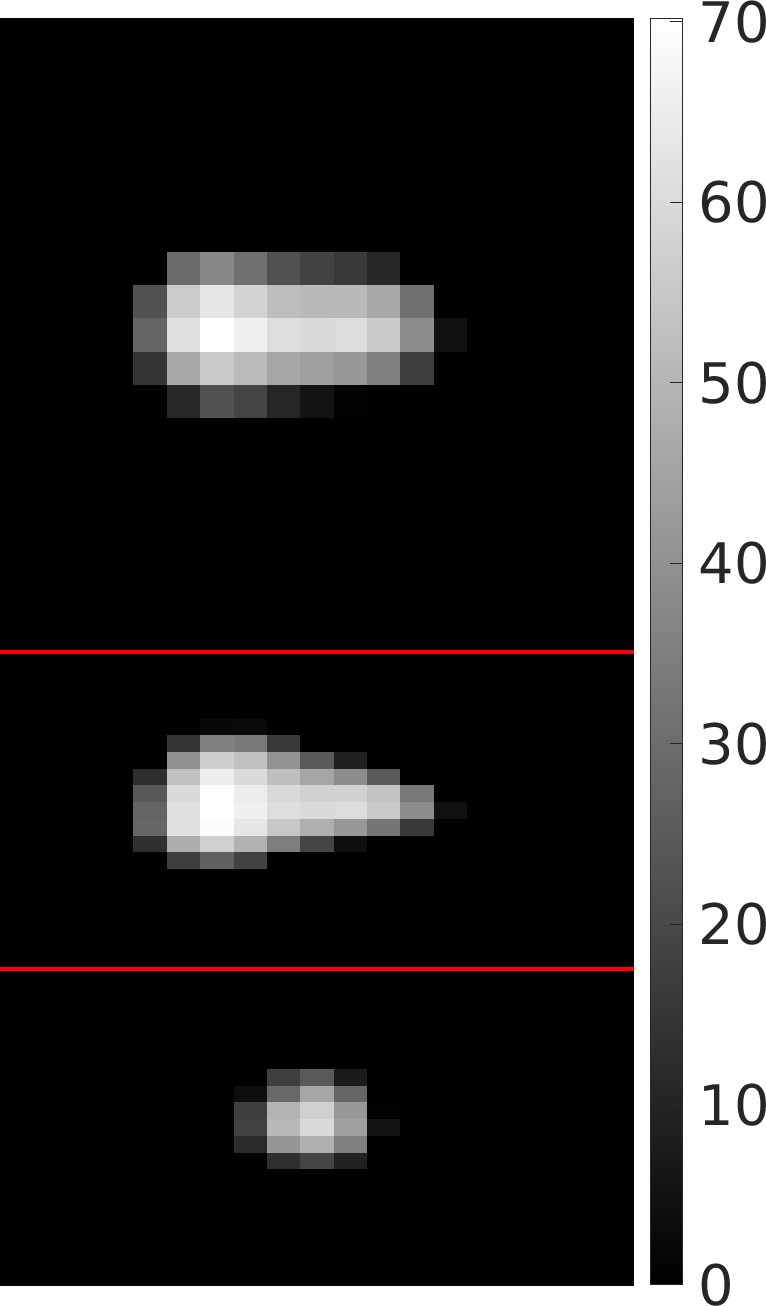}&
 \includegraphics[height=3.4cm]{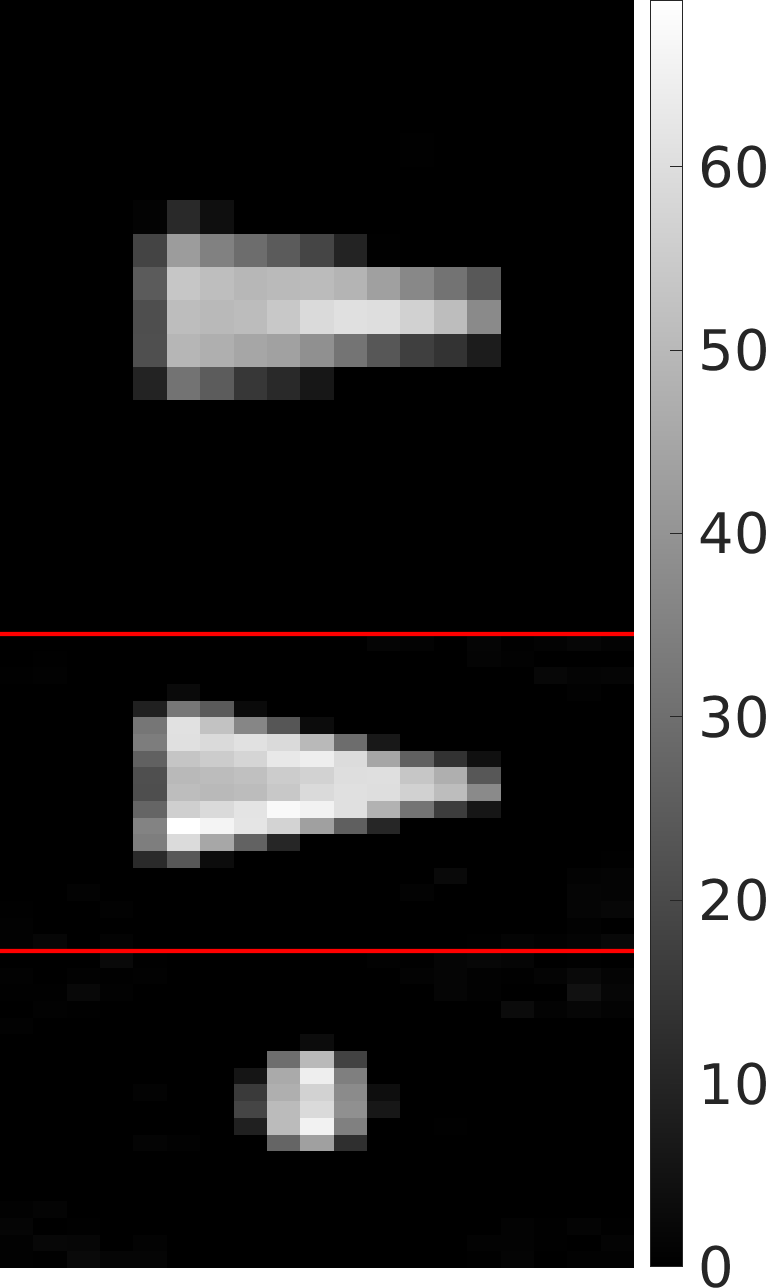}\\
\hline
\multicolumn{6}{l}{$\tau=3$} \\
 \includegraphics[height=3.4cm]{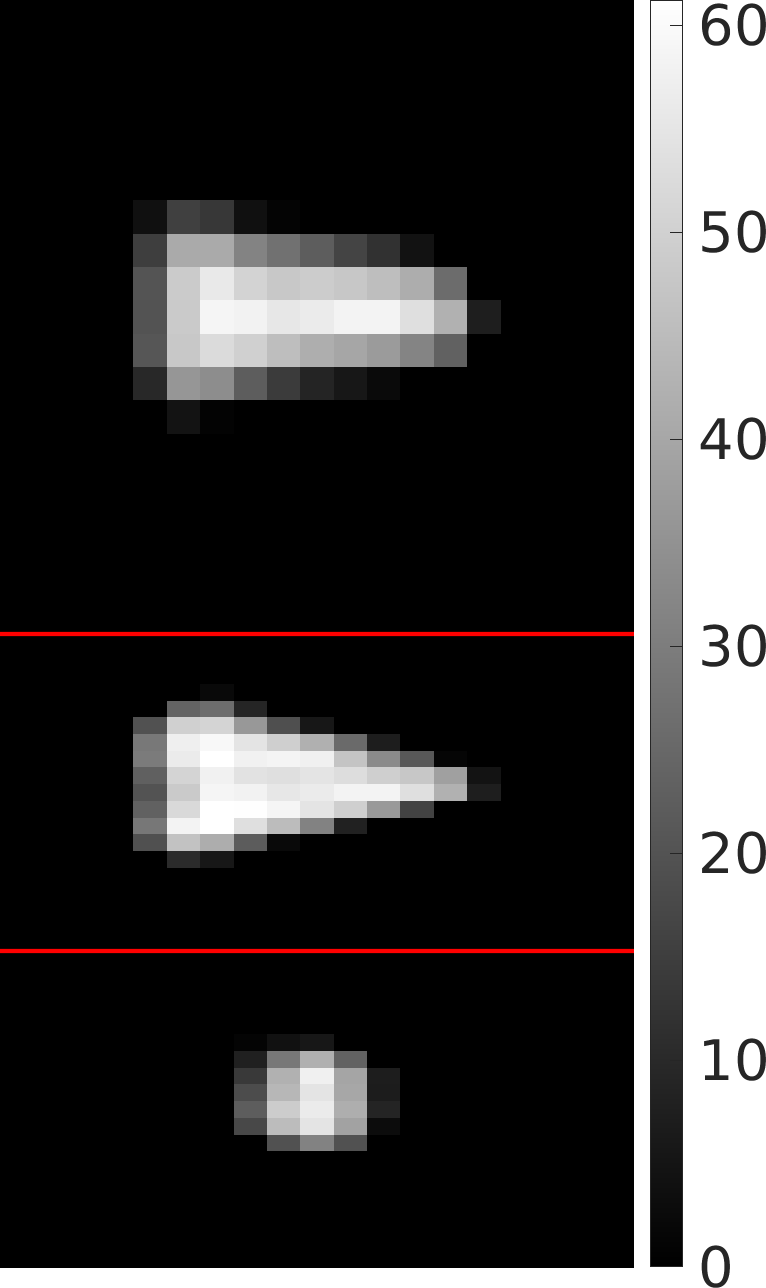}&
 \includegraphics[height=3.4cm]{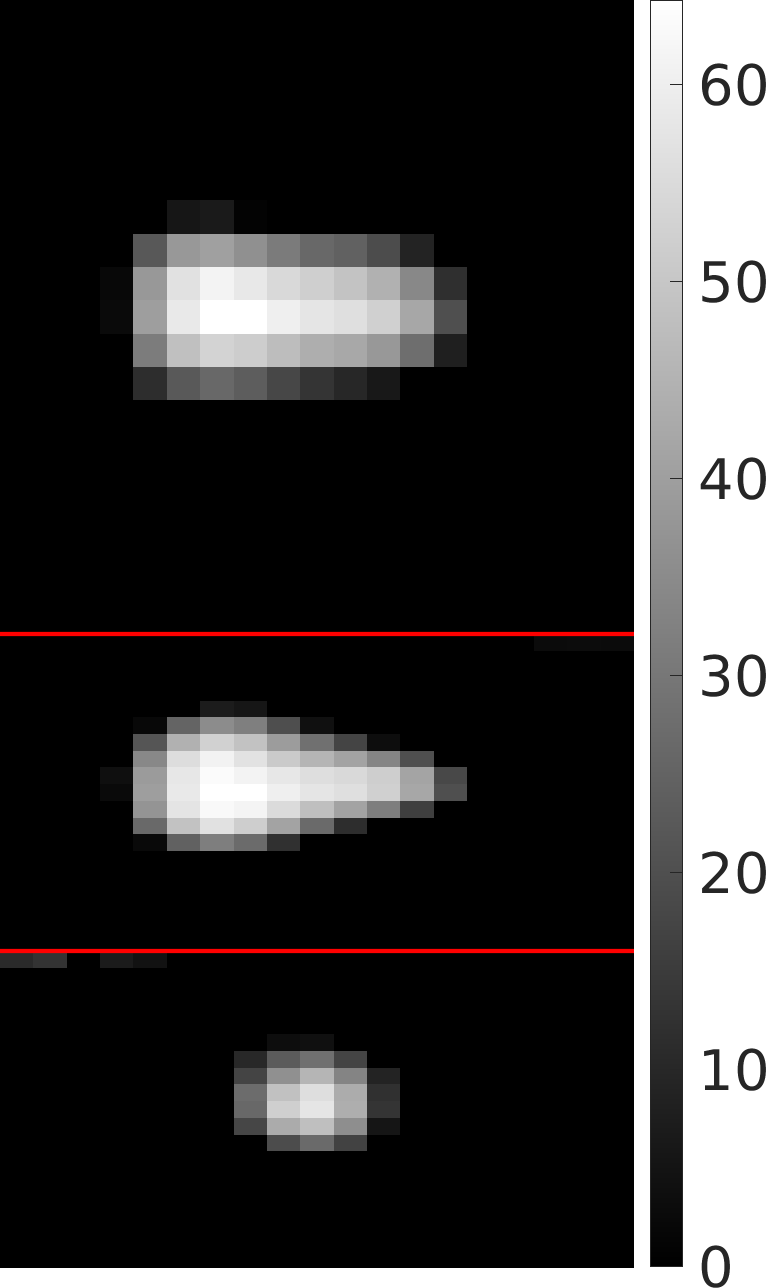}&
 \includegraphics[height=3.4cm]{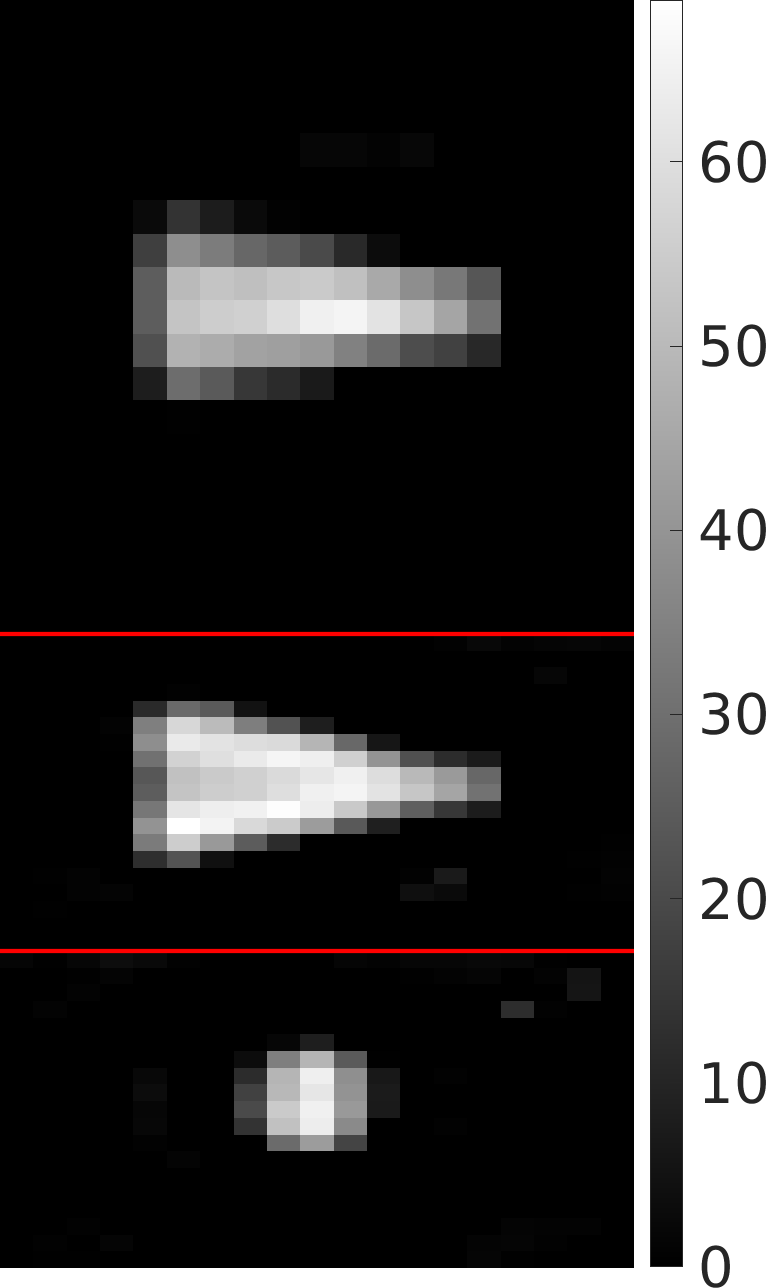}&
 \includegraphics[height=3.4cm]{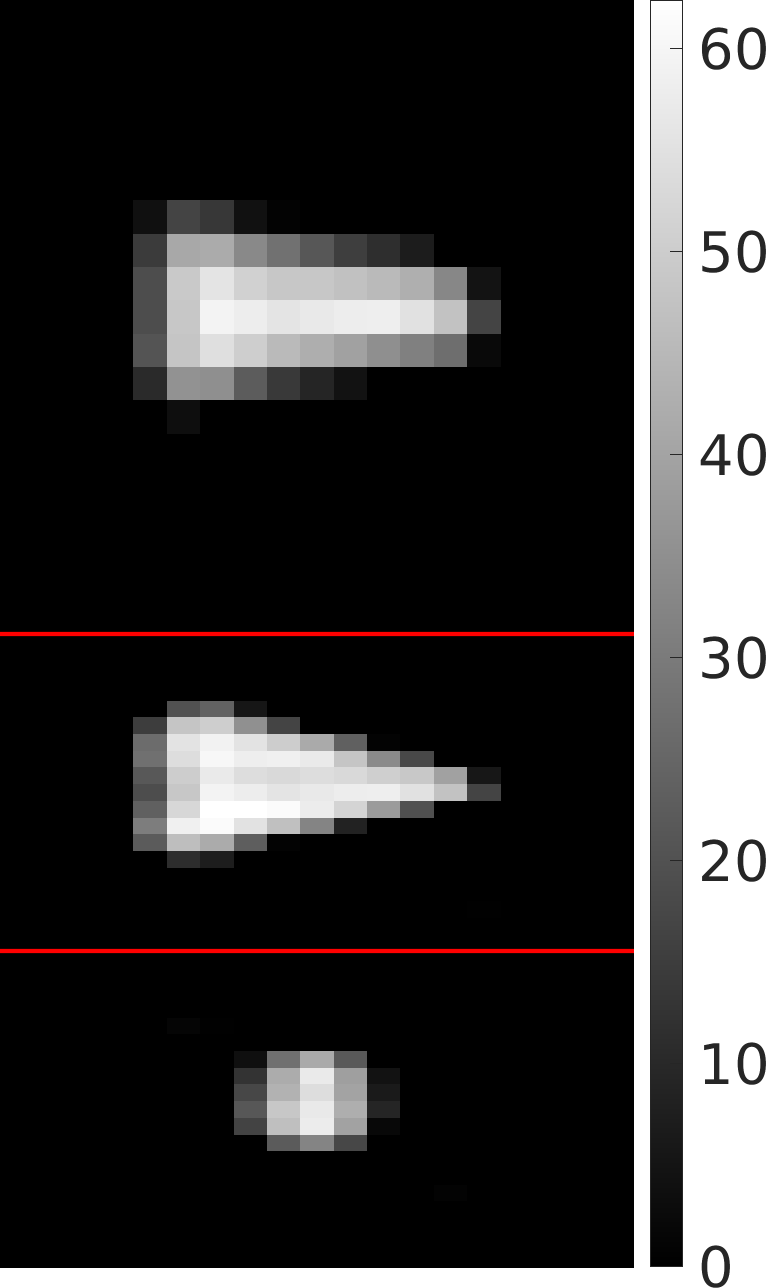}&
 \includegraphics[height=3.4cm]{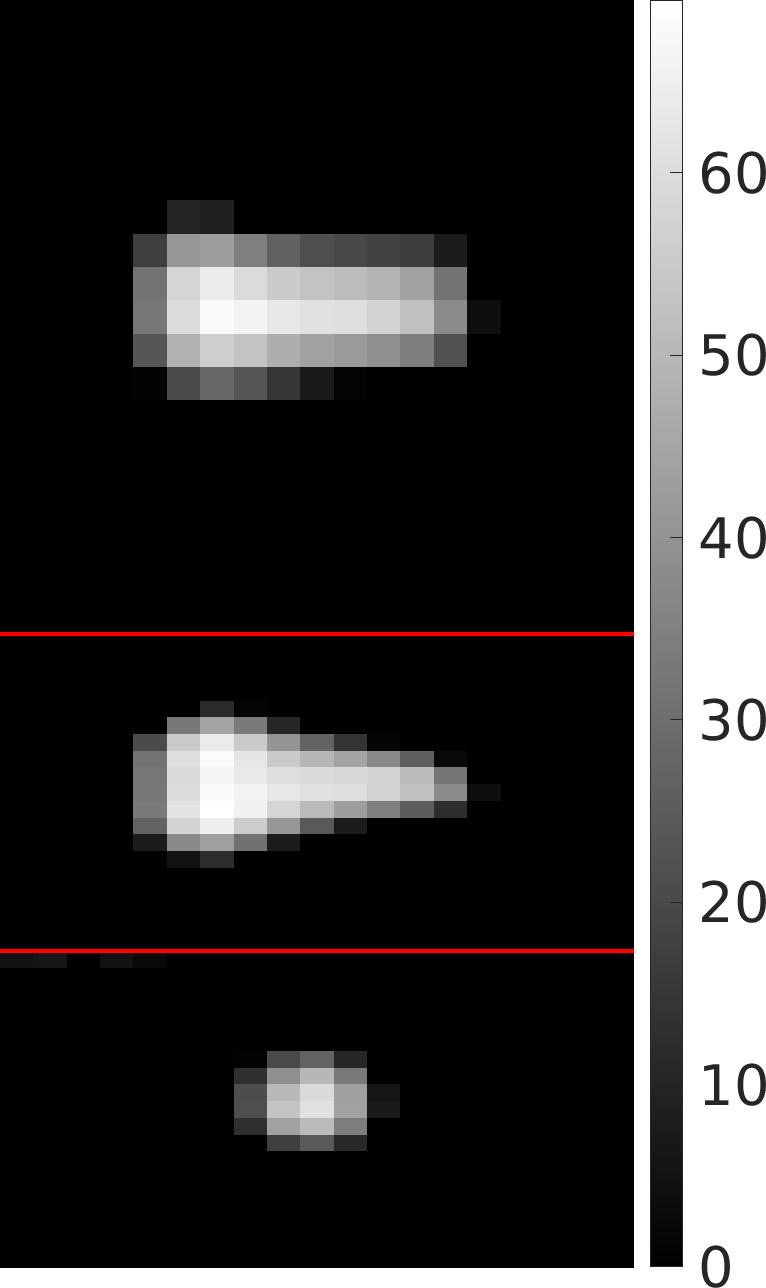}&
 \includegraphics[height=3.4cm]{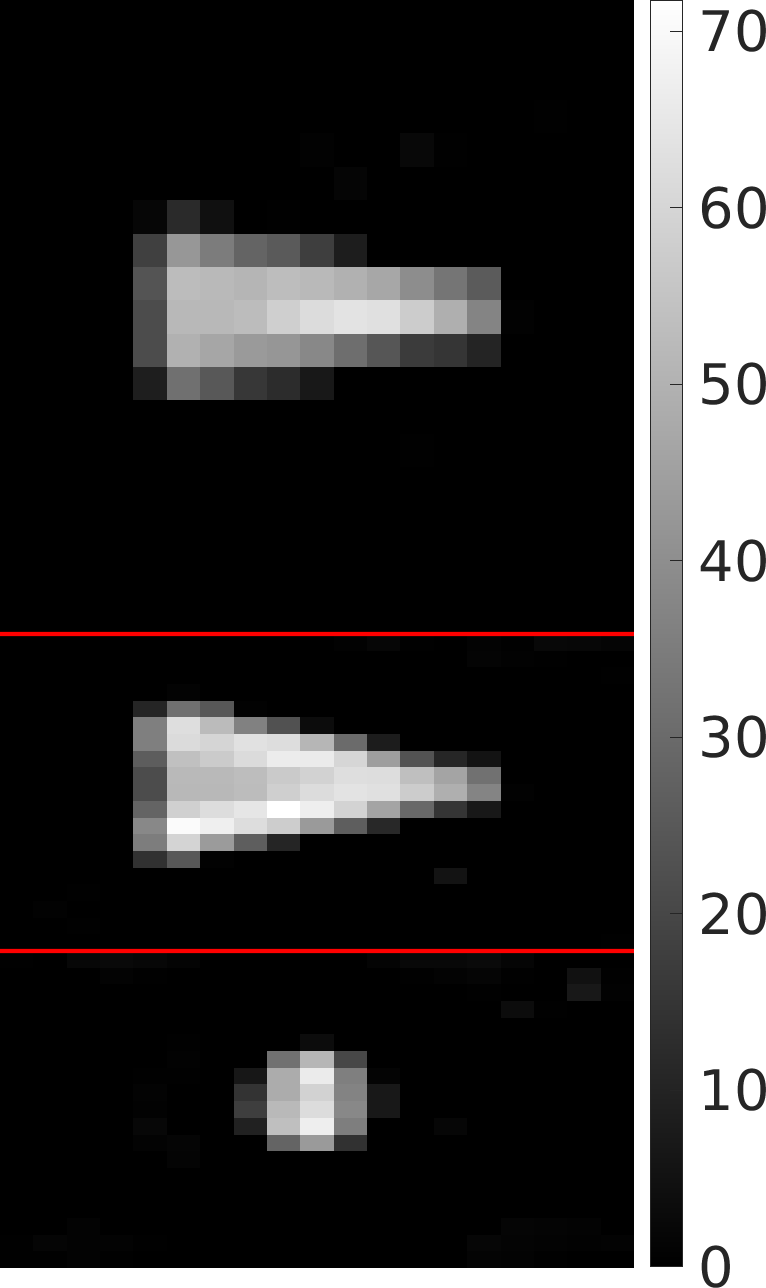}\\
 \hline
\multicolumn{6}{l}{$\tau=5$} \\
 \includegraphics[height=3.4cm]{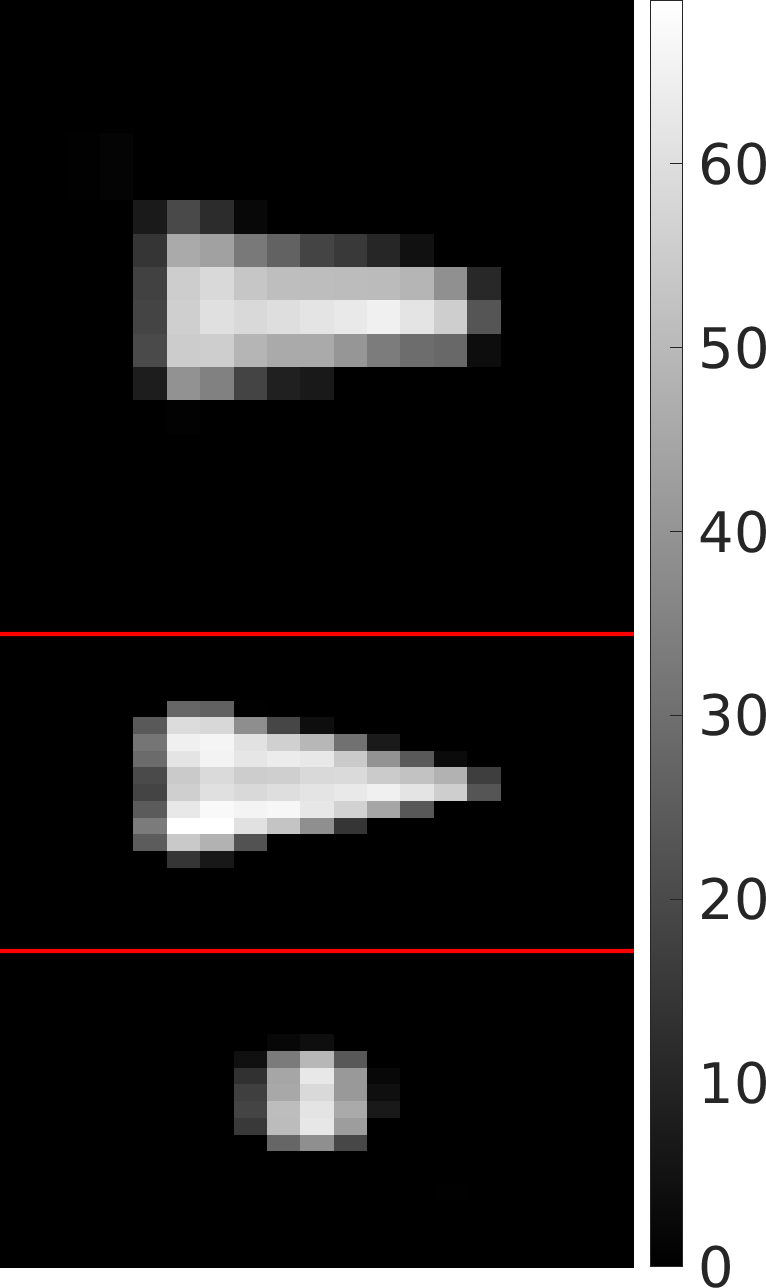}&
 \includegraphics[height=3.4cm]{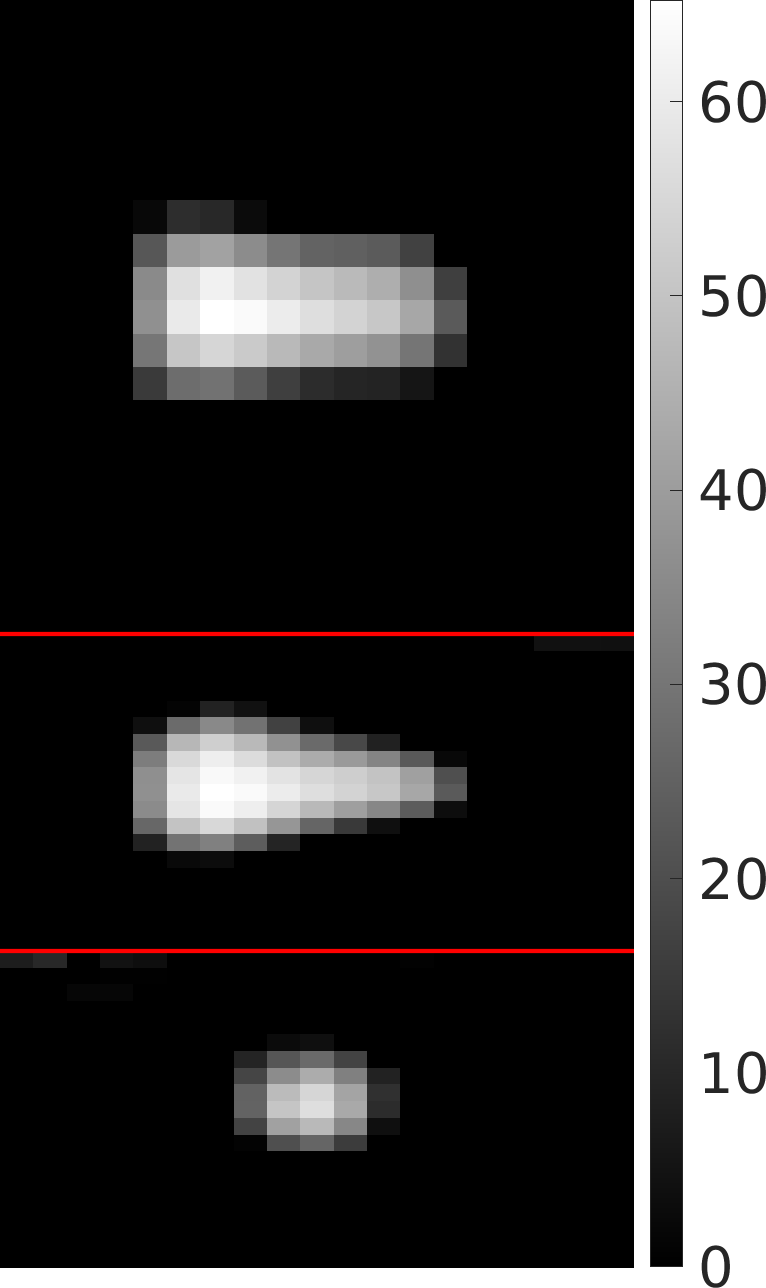}&
 \includegraphics[height=3.4cm]{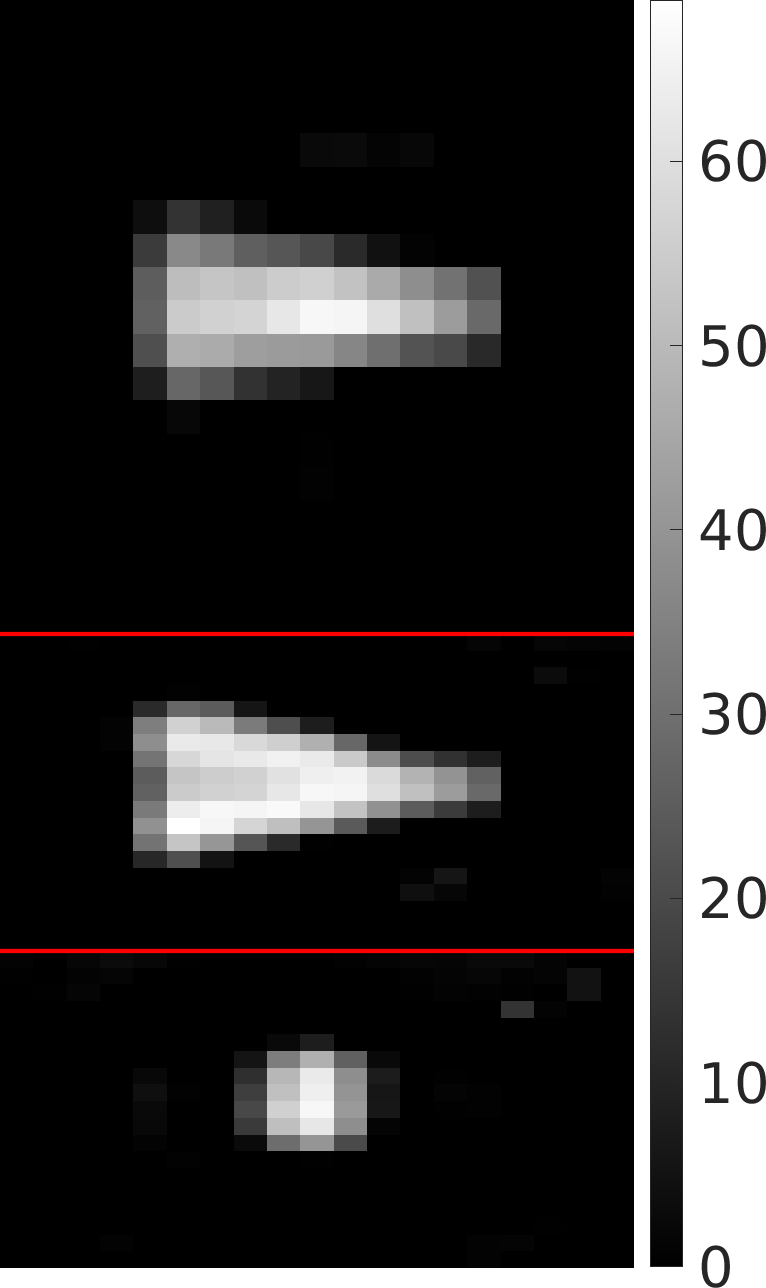}&
 \includegraphics[height=3.4cm]{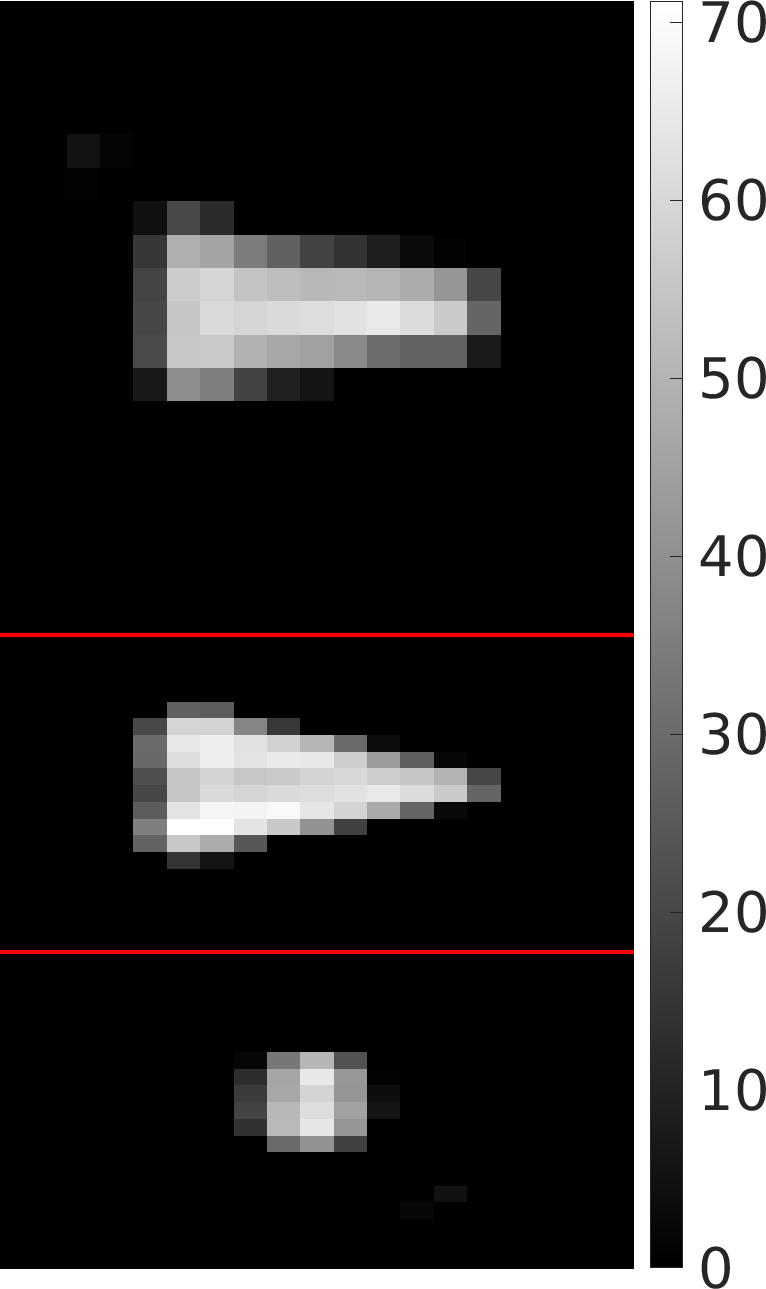}&
 \includegraphics[height=3.4cm]{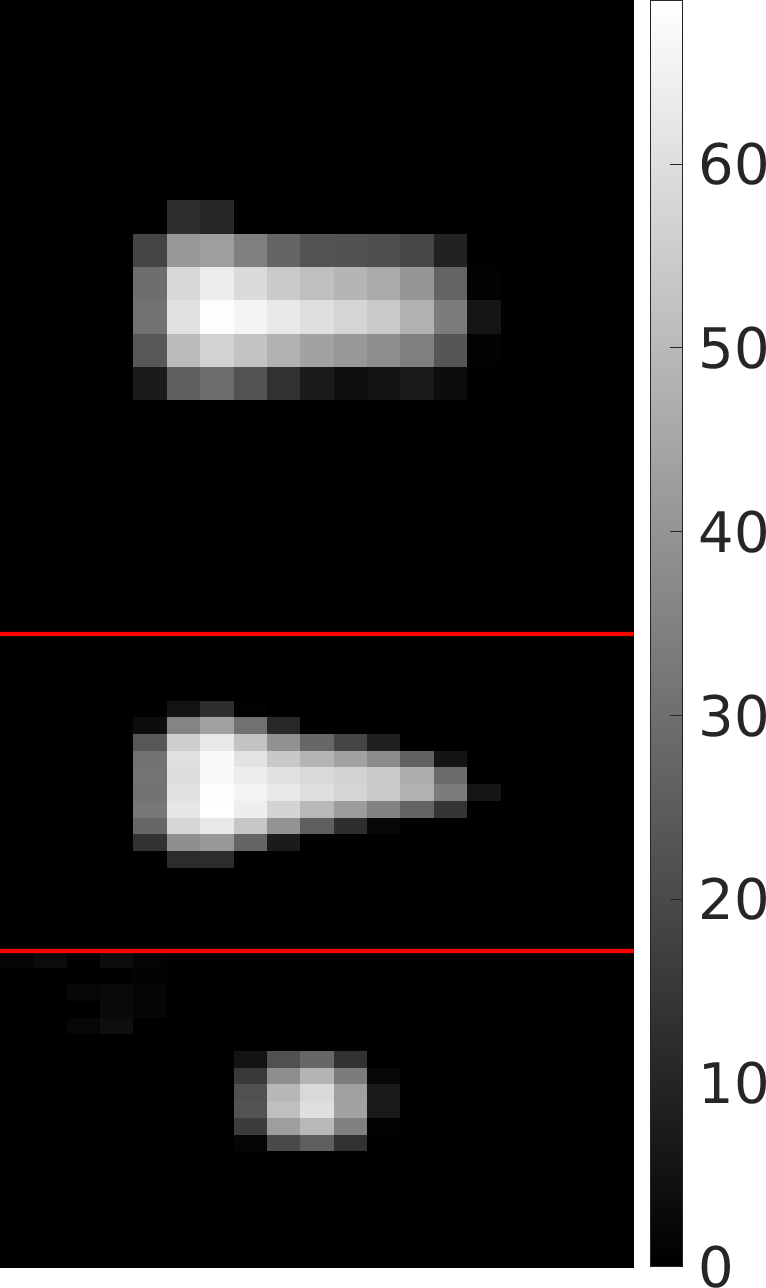}&
 \includegraphics[height=3.4cm]{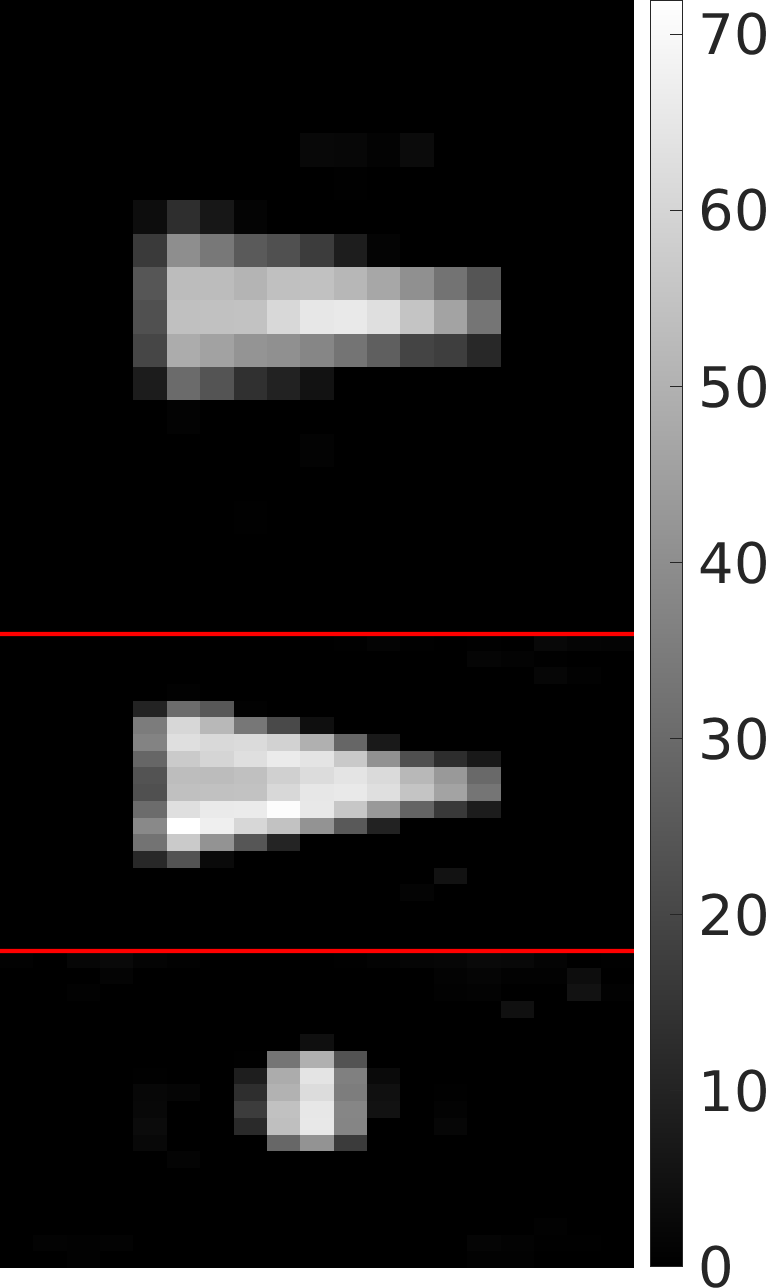}\\
\end{tabular}
}
\caption{``Shape'' phantom reconstructions, PSNR-optimized $\alpha$ and iteration number $N$ (for l2-K only) according to Table \ref{tab:psnr_nonwhitened_vs_whitened}.}
\label{fig:methods_nonwhitened_vs_whitened_shape_psnr}
\end{figure*}

\begin{figure*}[hbt!]
\centering
\scalebox{0.85}{
\begin{tabular}{ccc|ccc}
\multicolumn{3}{c|}{non-whitened} & \multicolumn{3}{c}{whitened} \\
\hline
l1-L & l2-L & l2-K & l1-L & l2-L & l2-K \\
\hline
\multicolumn{6}{l}{$\tau=0$} \\
 \includegraphics[height=3.4cm]{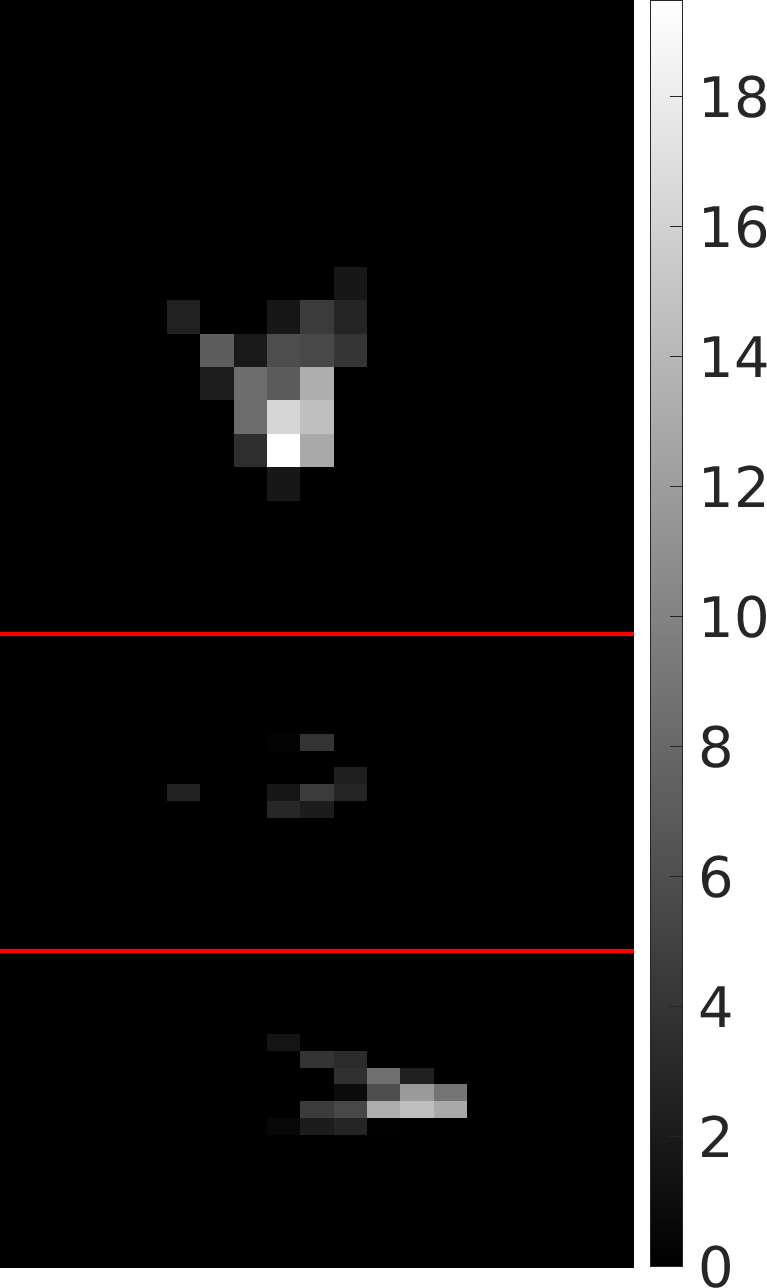}&
 \includegraphics[height=3.4cm]{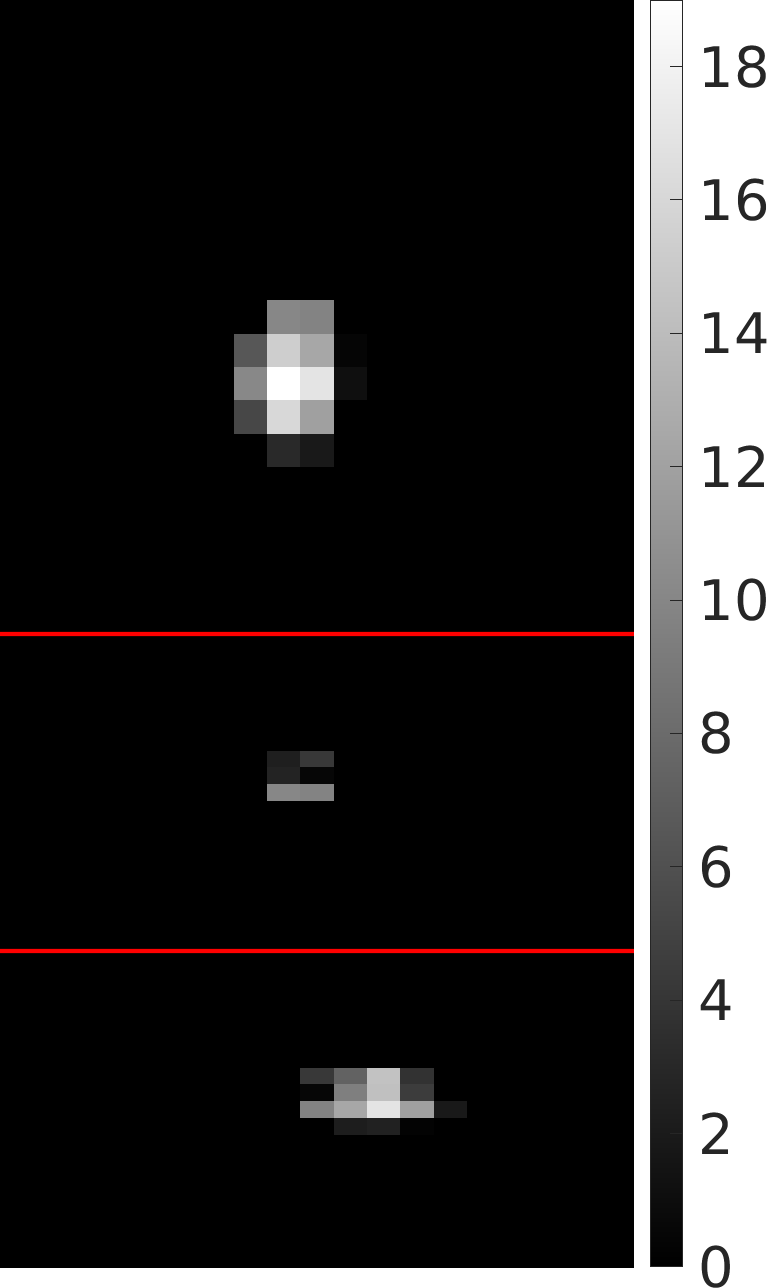}&
 \includegraphics[height=3.4cm]{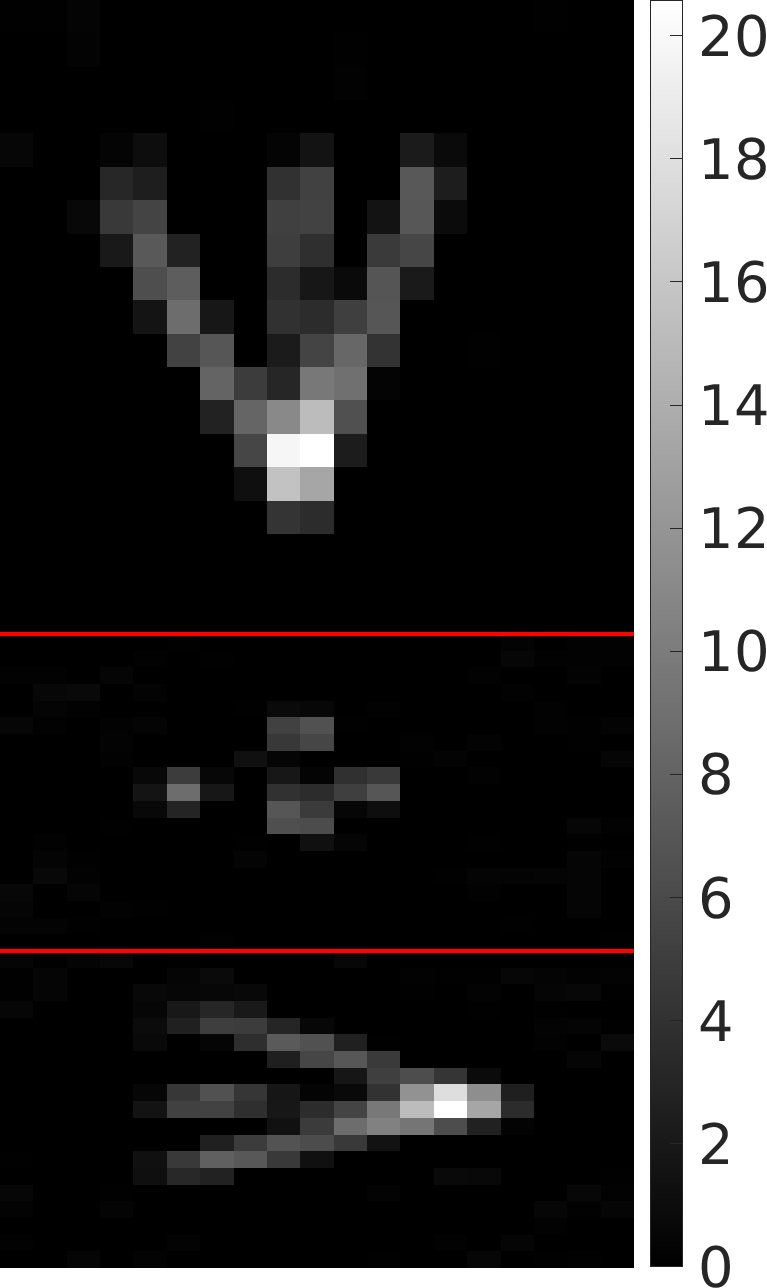}&
 \includegraphics[height=3.4cm]{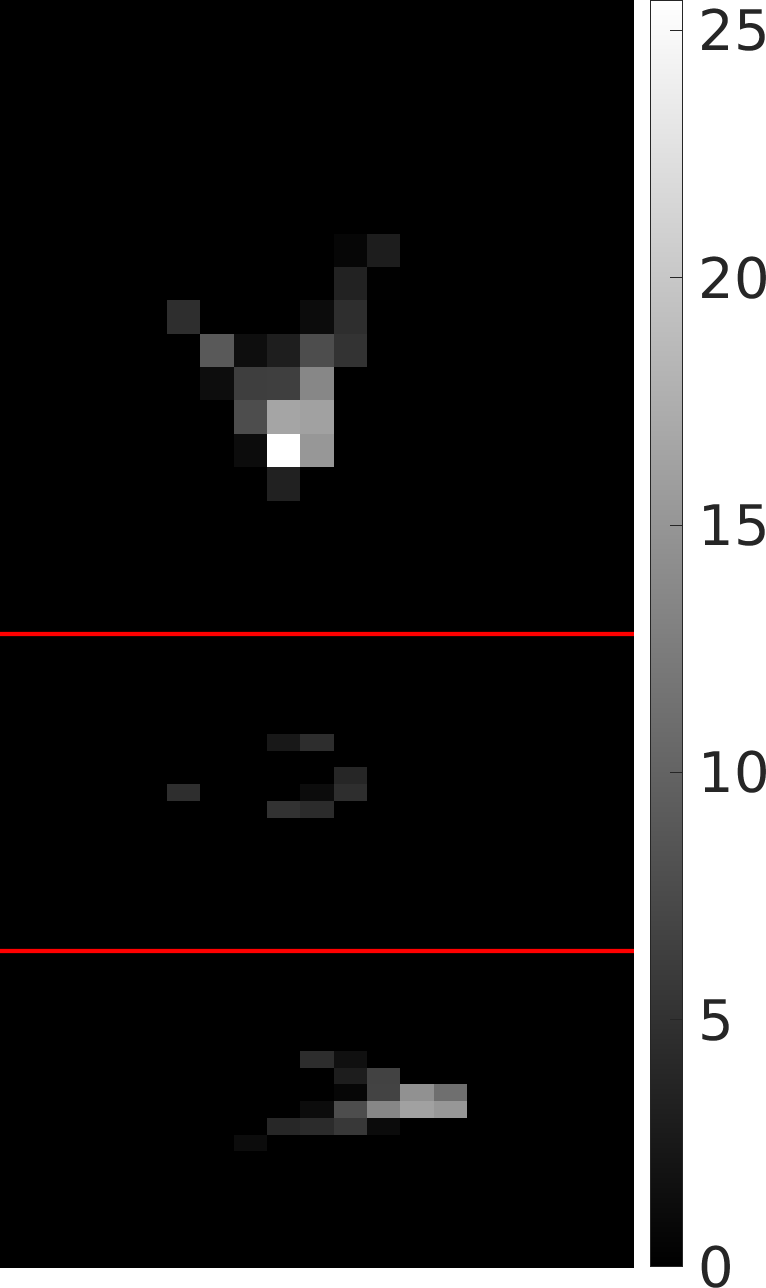}&
 \includegraphics[height=3.4cm]{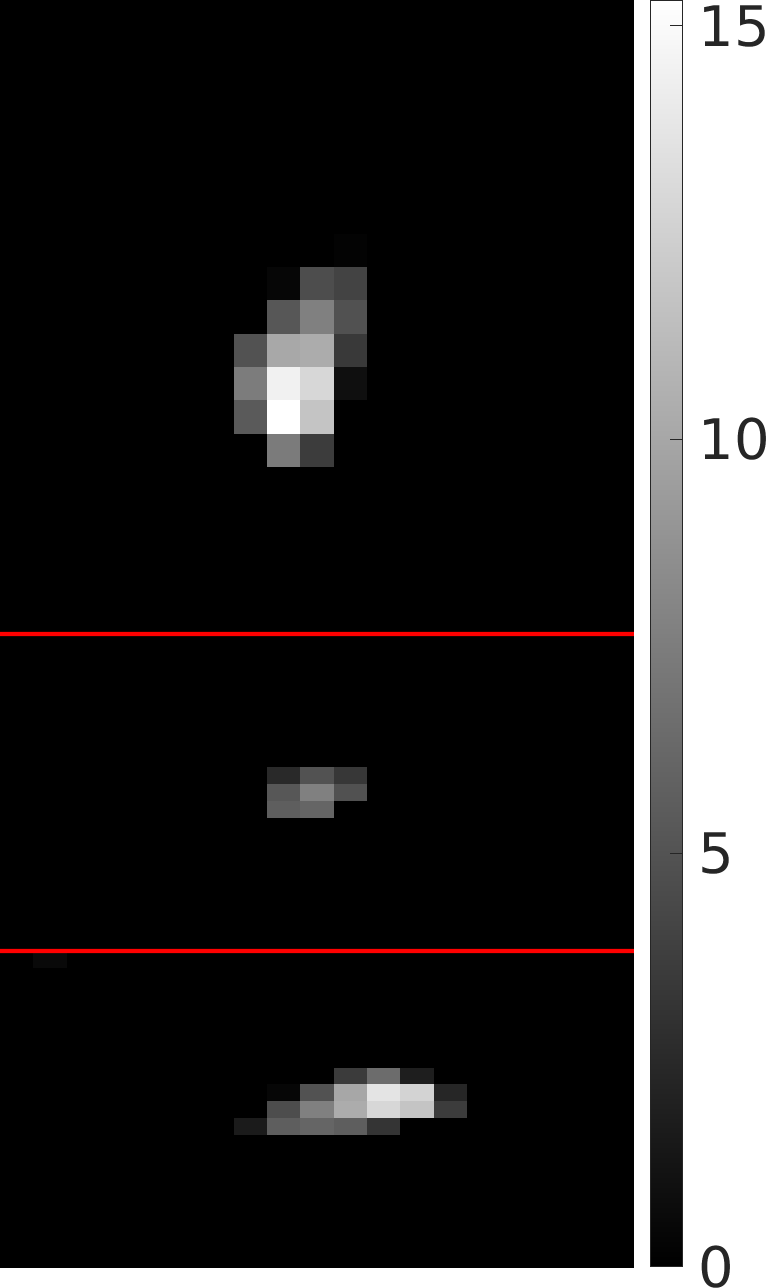}&
 \includegraphics[height=3.4cm]{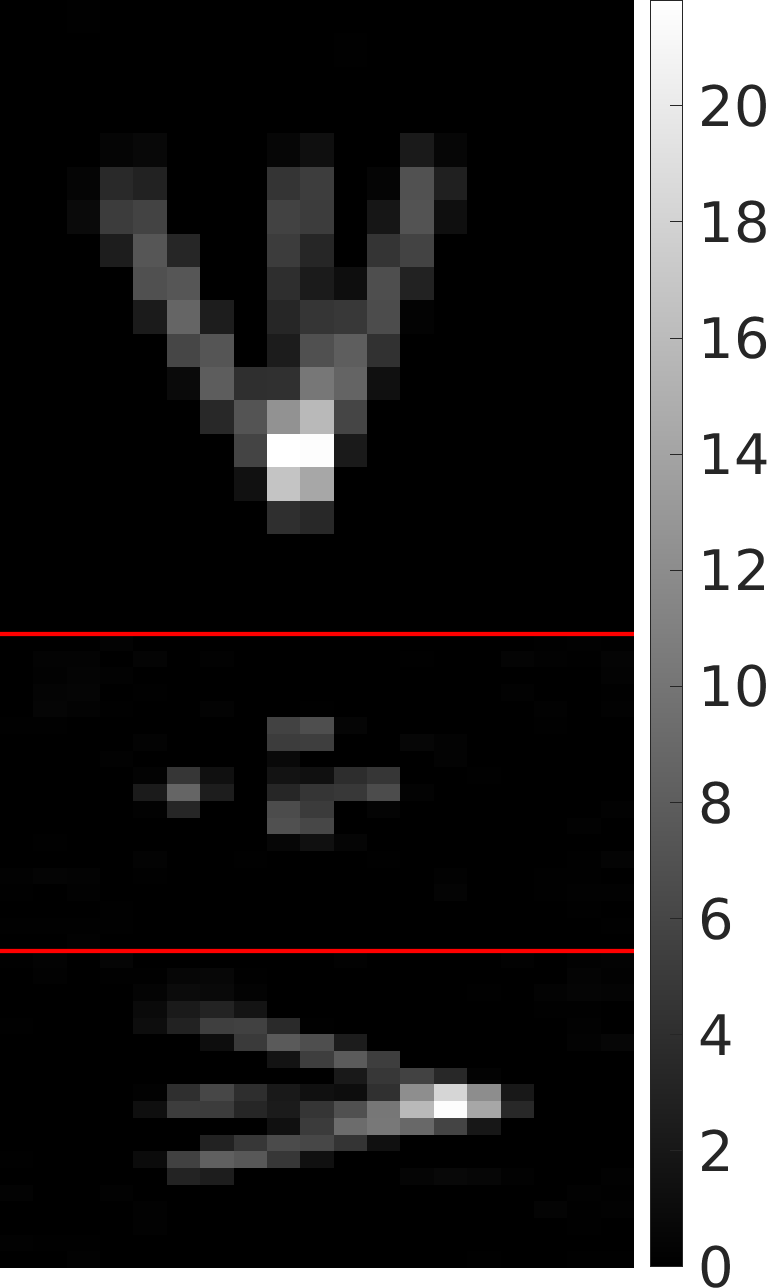}\\
\hline
\multicolumn{6}{l}{$\tau=1$} \\
 \includegraphics[height=3.4cm]{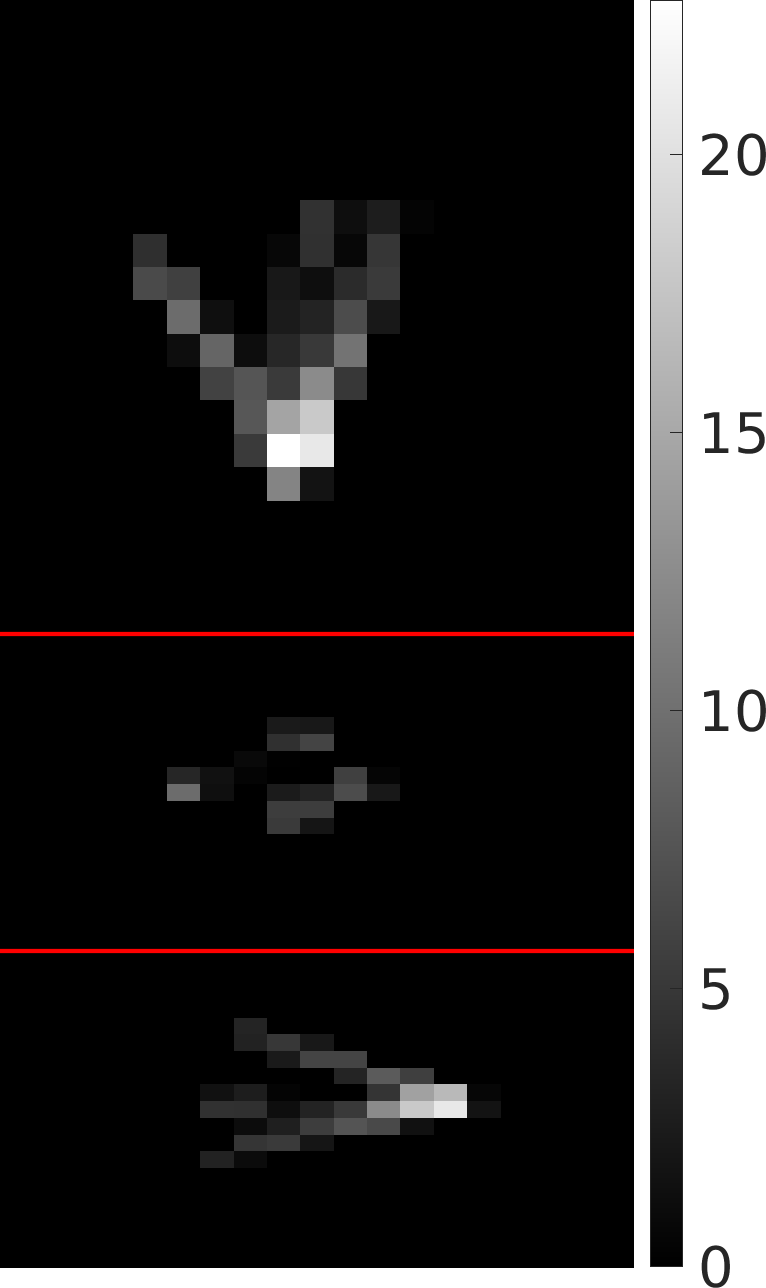}&
 \includegraphics[height=3.4cm]{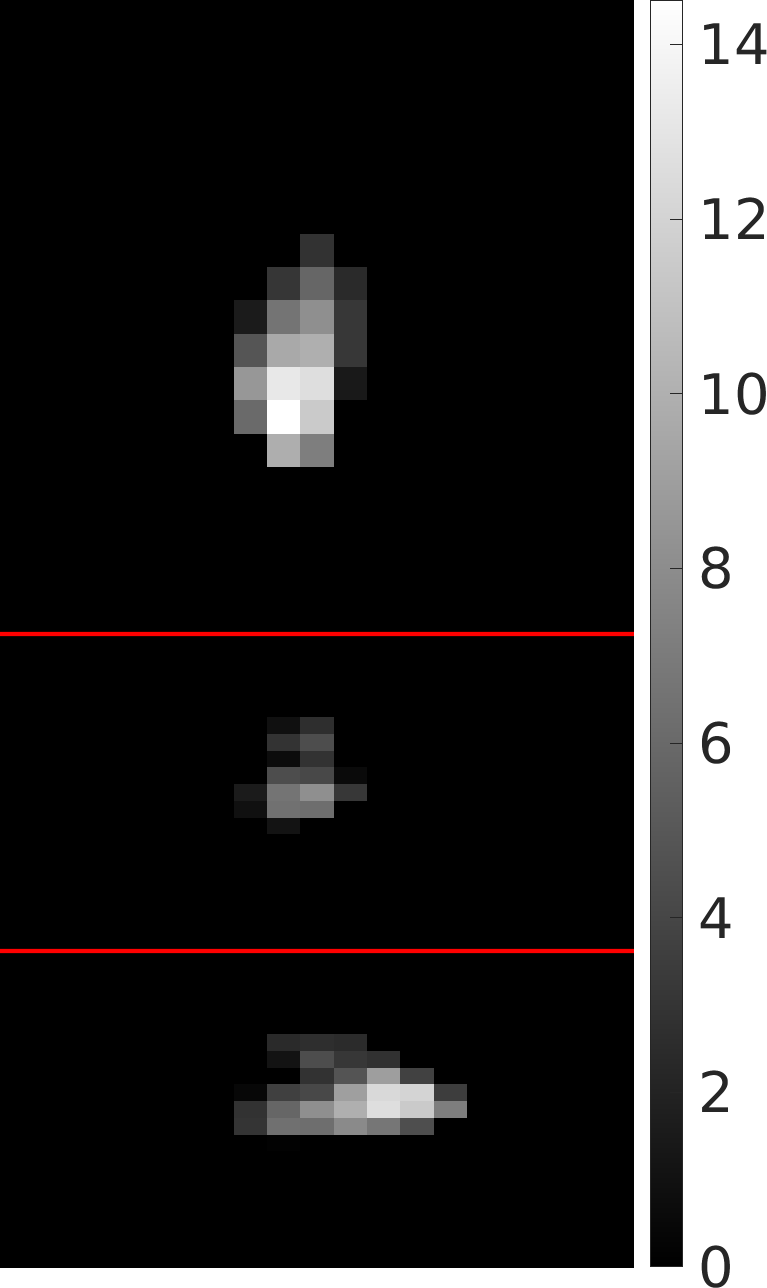}&
 \includegraphics[height=3.4cm]{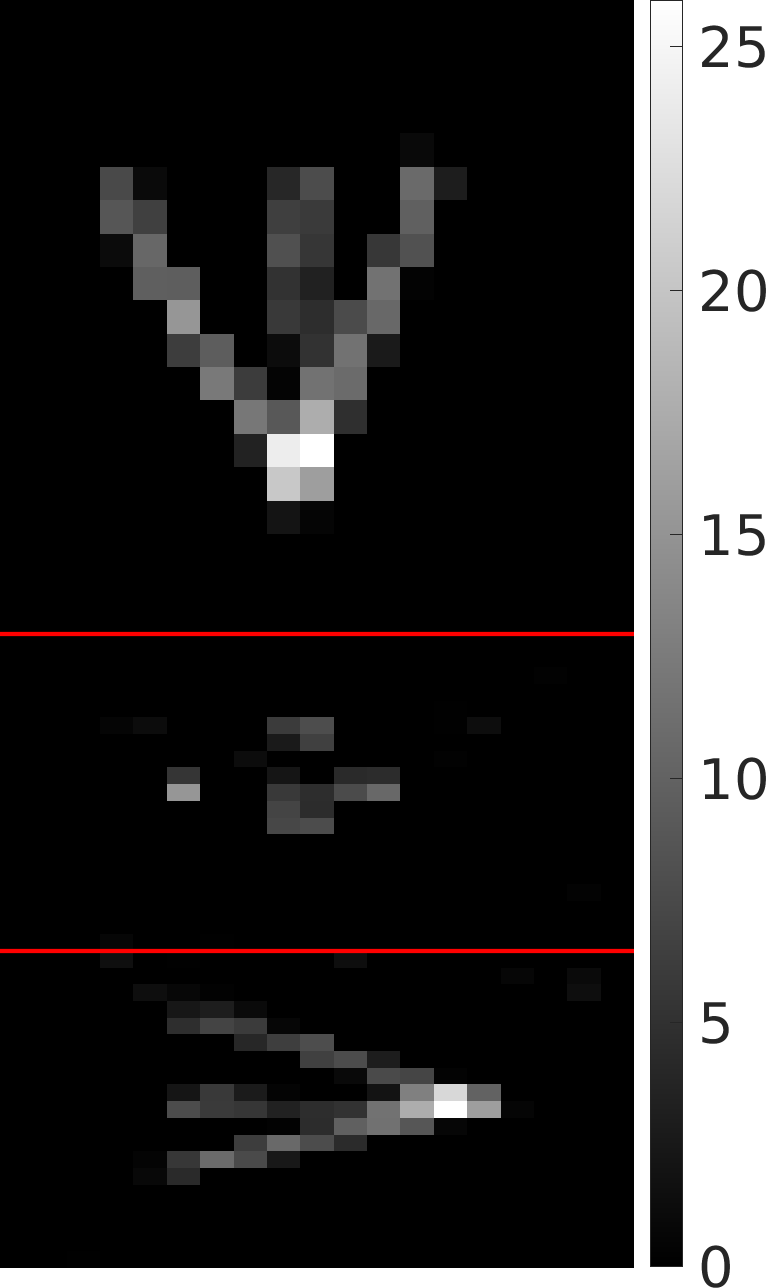}&
 \includegraphics[height=3.4cm]{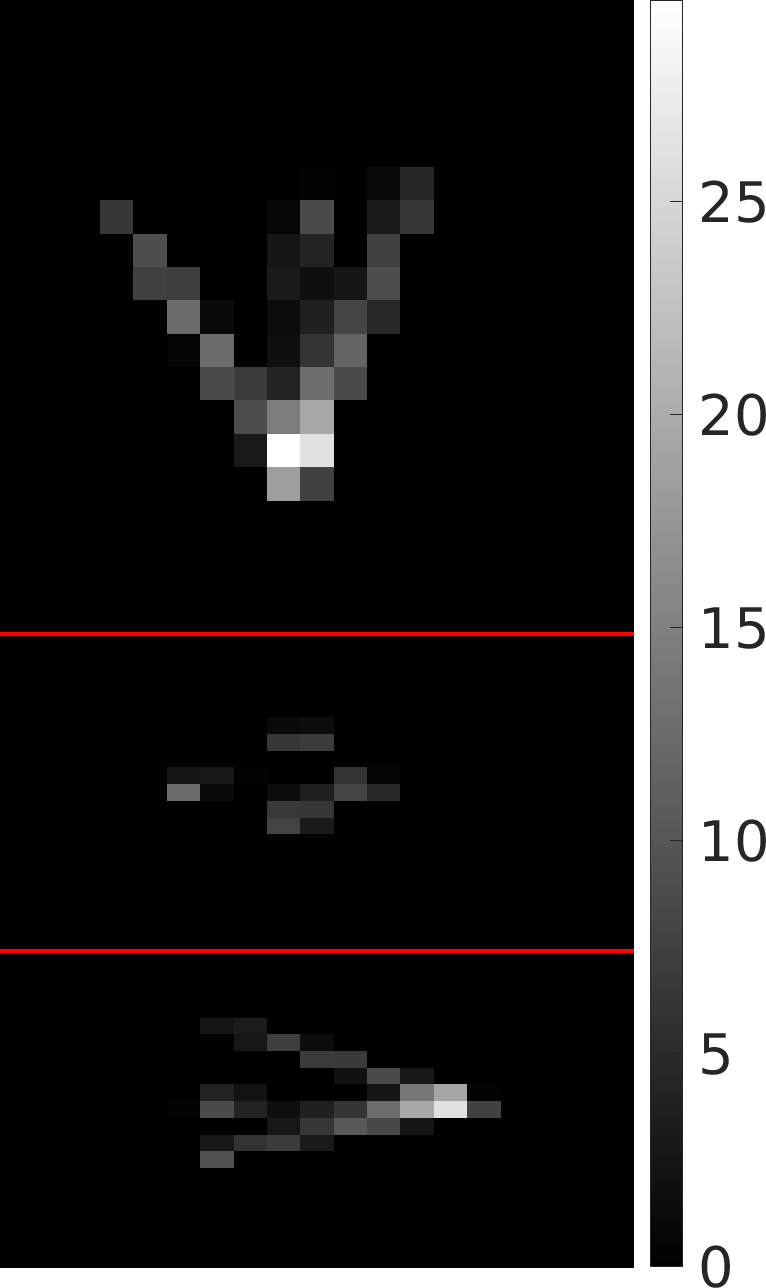}&
 \includegraphics[height=3.4cm]{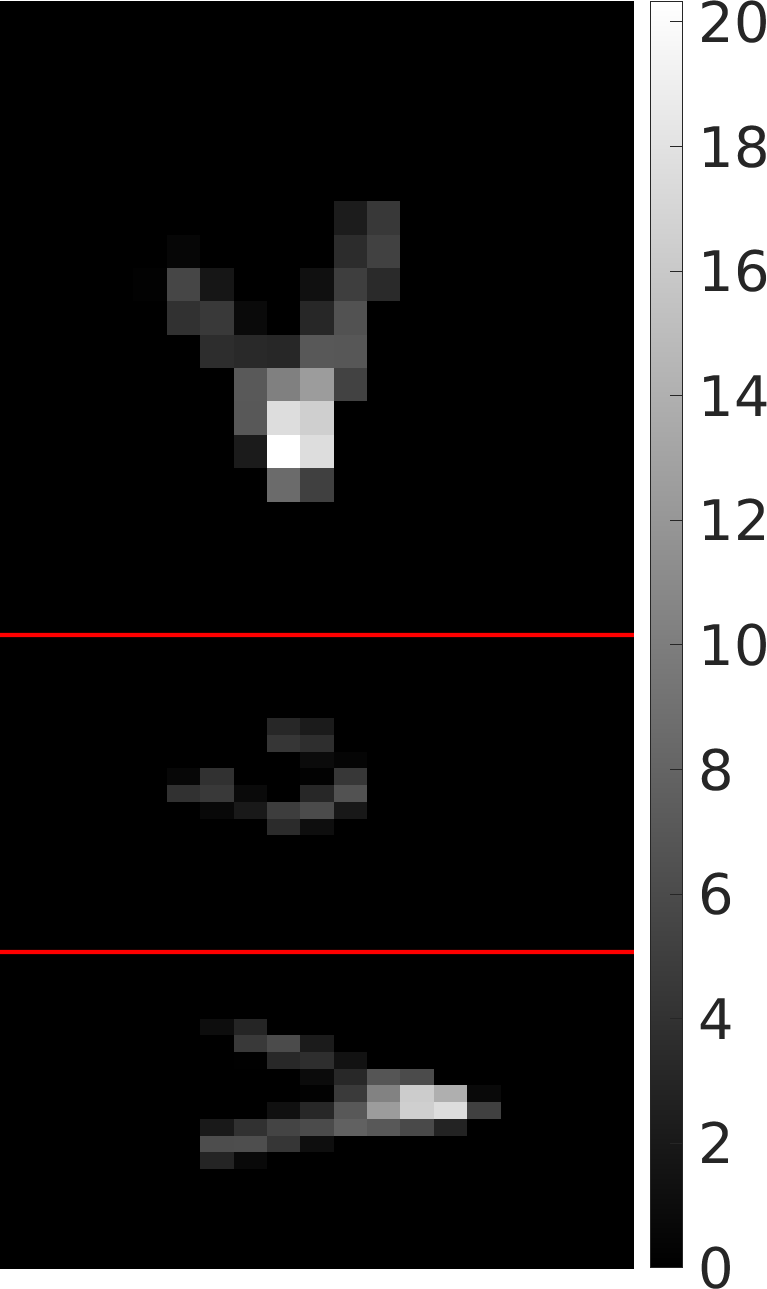}&
 \includegraphics[height=3.4cm]{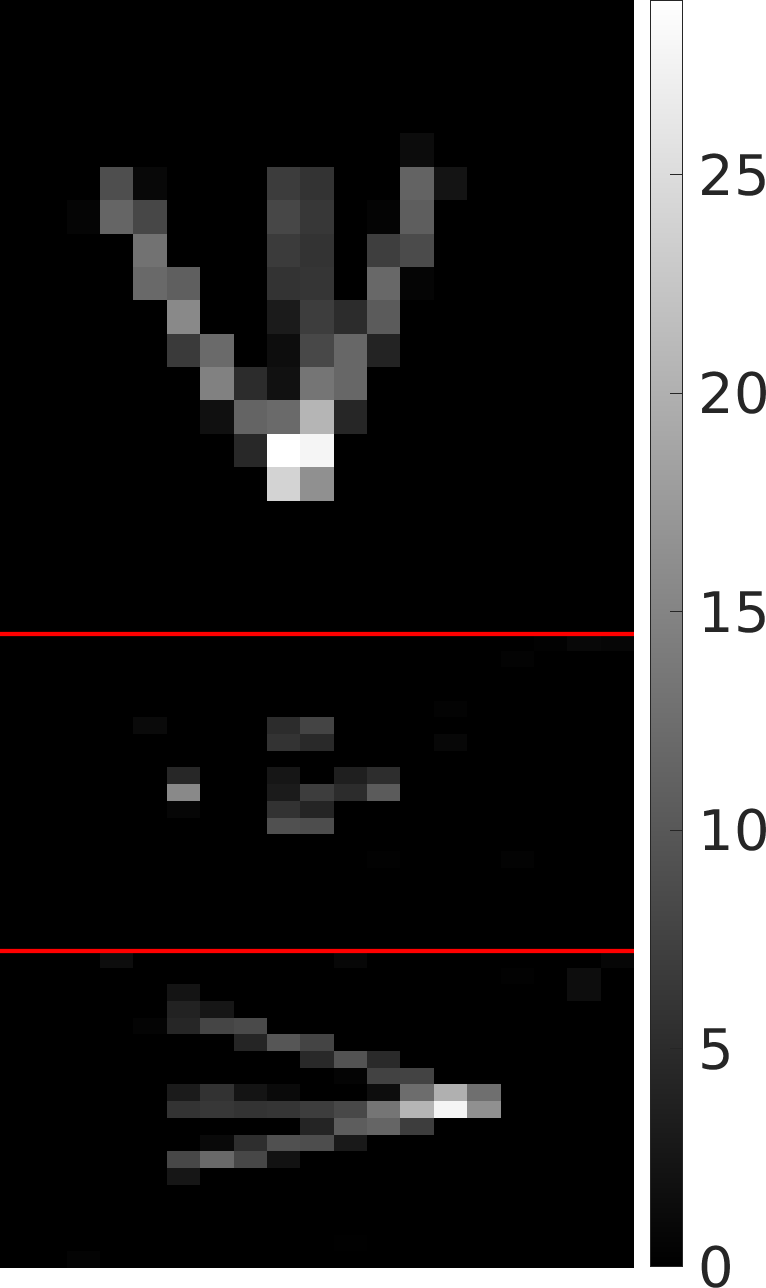}\\
\hline
\multicolumn{6}{l}{$\tau=3$} \\
 \includegraphics[height=3.4cm]{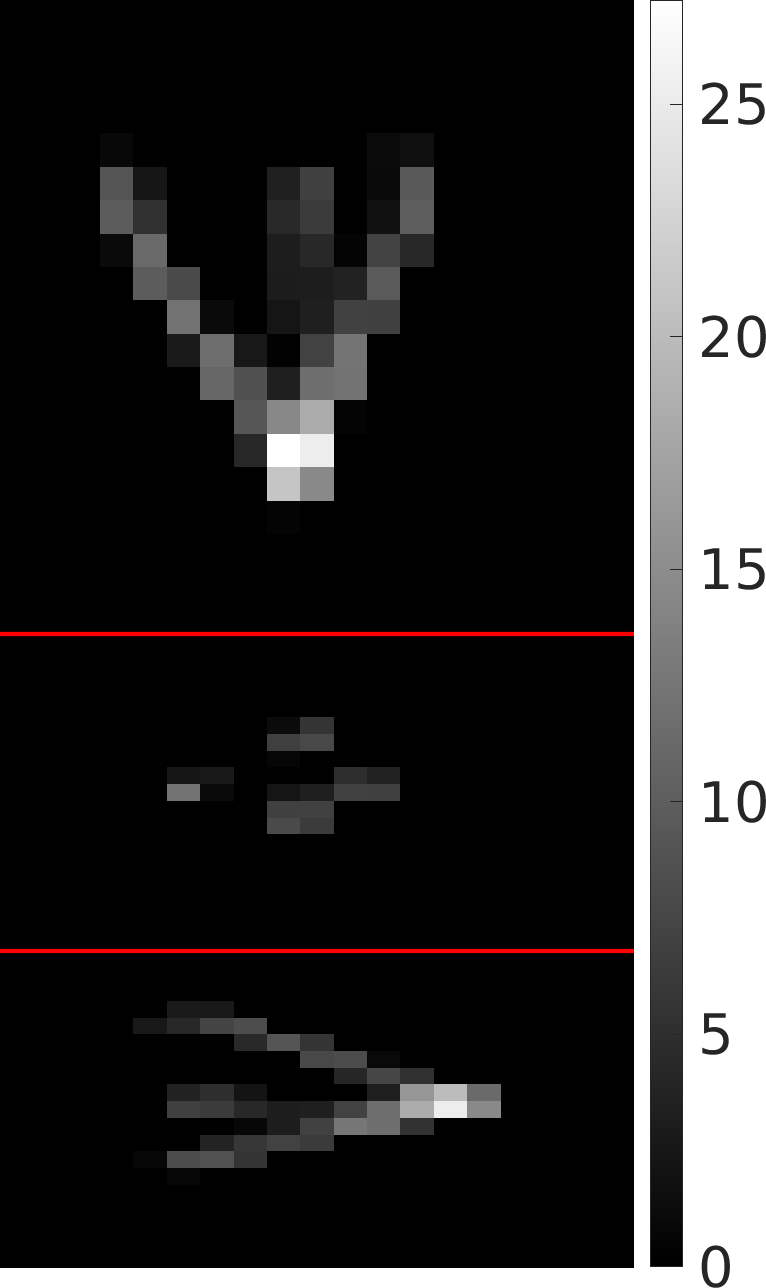}&
 \includegraphics[height=3.4cm]{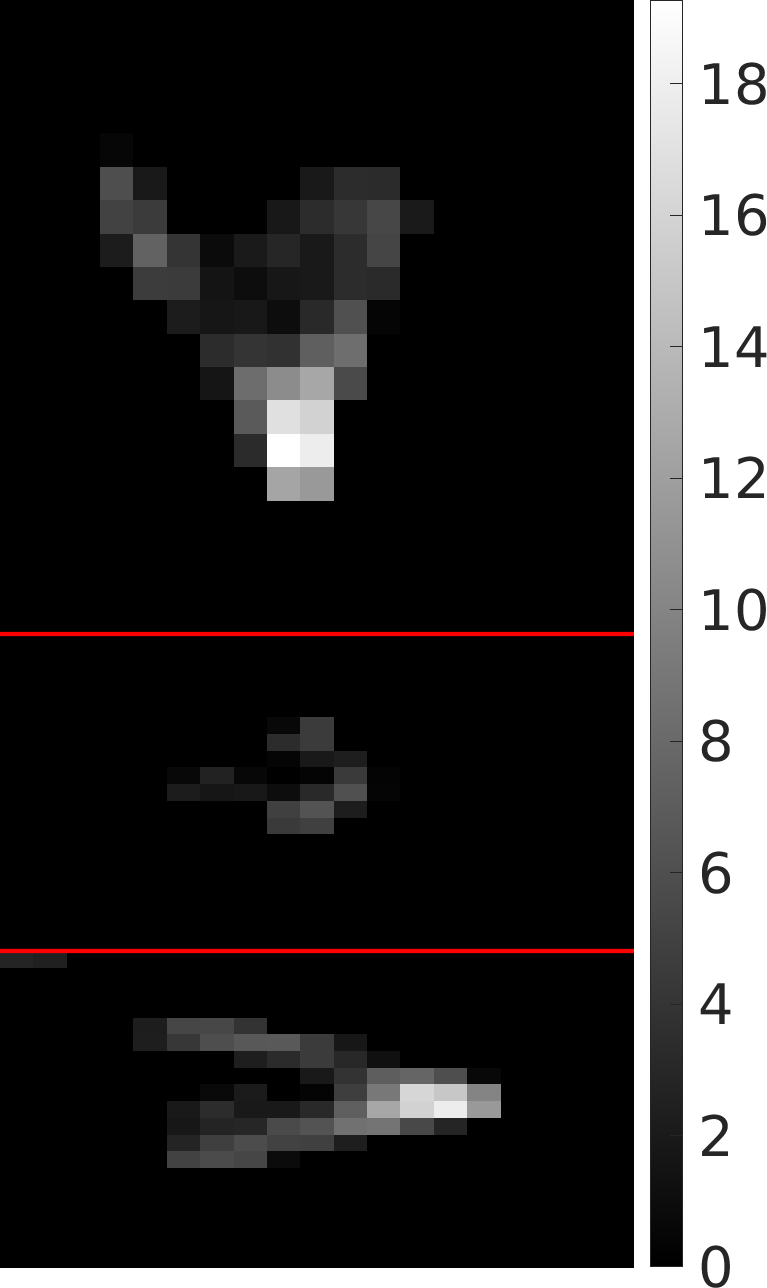}&
 \includegraphics[height=3.4cm]{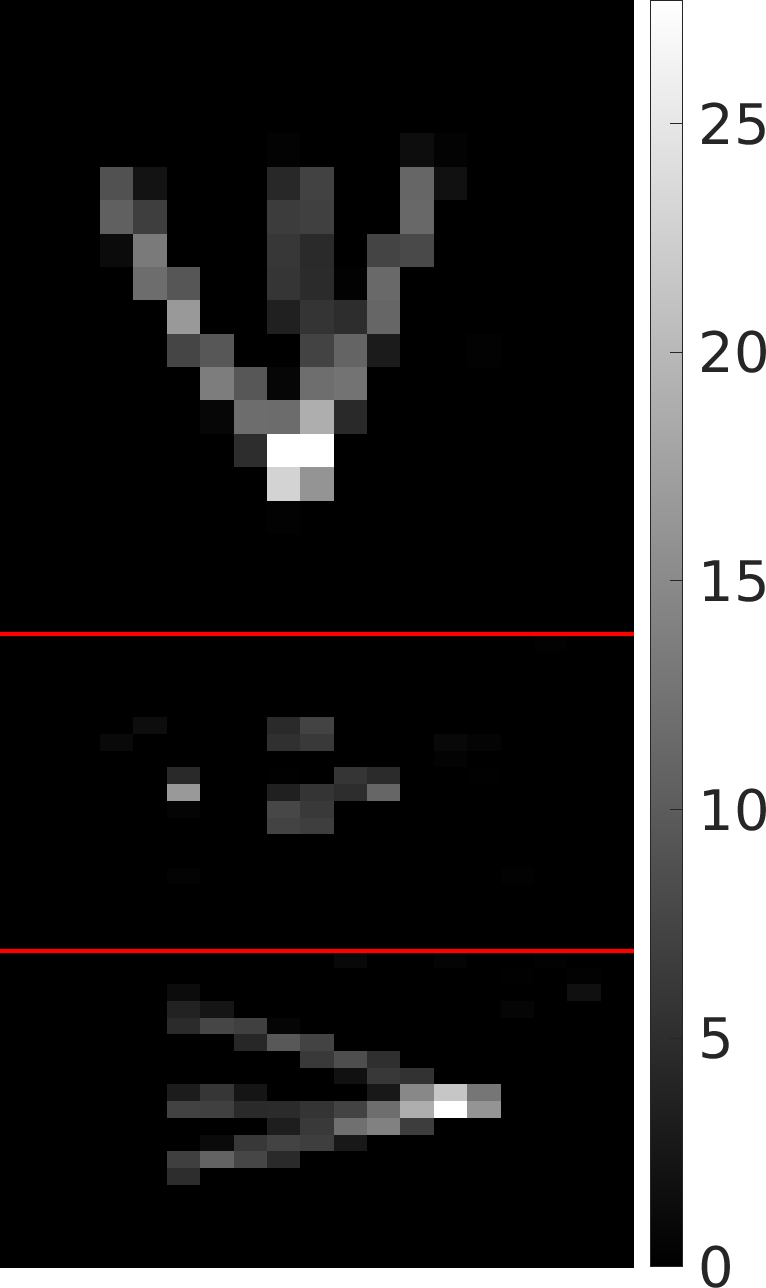}&
 \includegraphics[height=3.4cm]{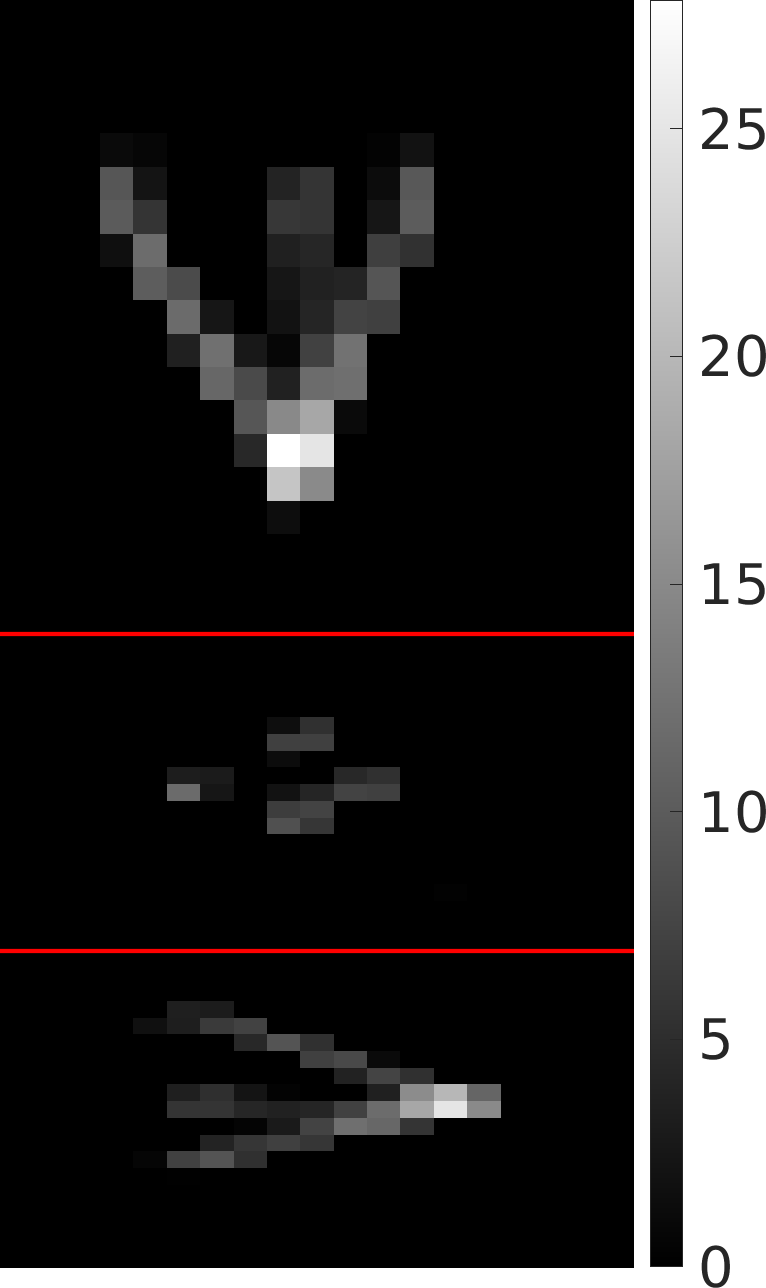}&
 \includegraphics[height=3.4cm]{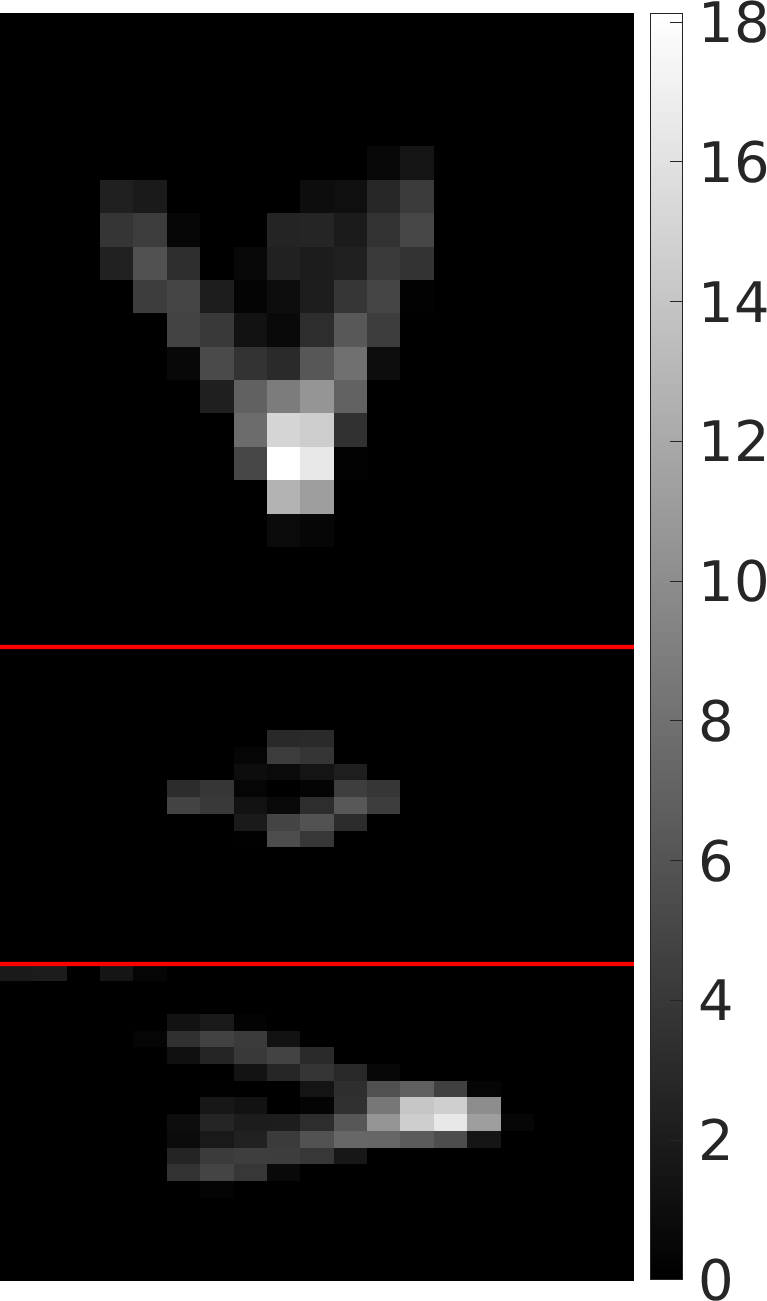}&
 \includegraphics[height=3.4cm]{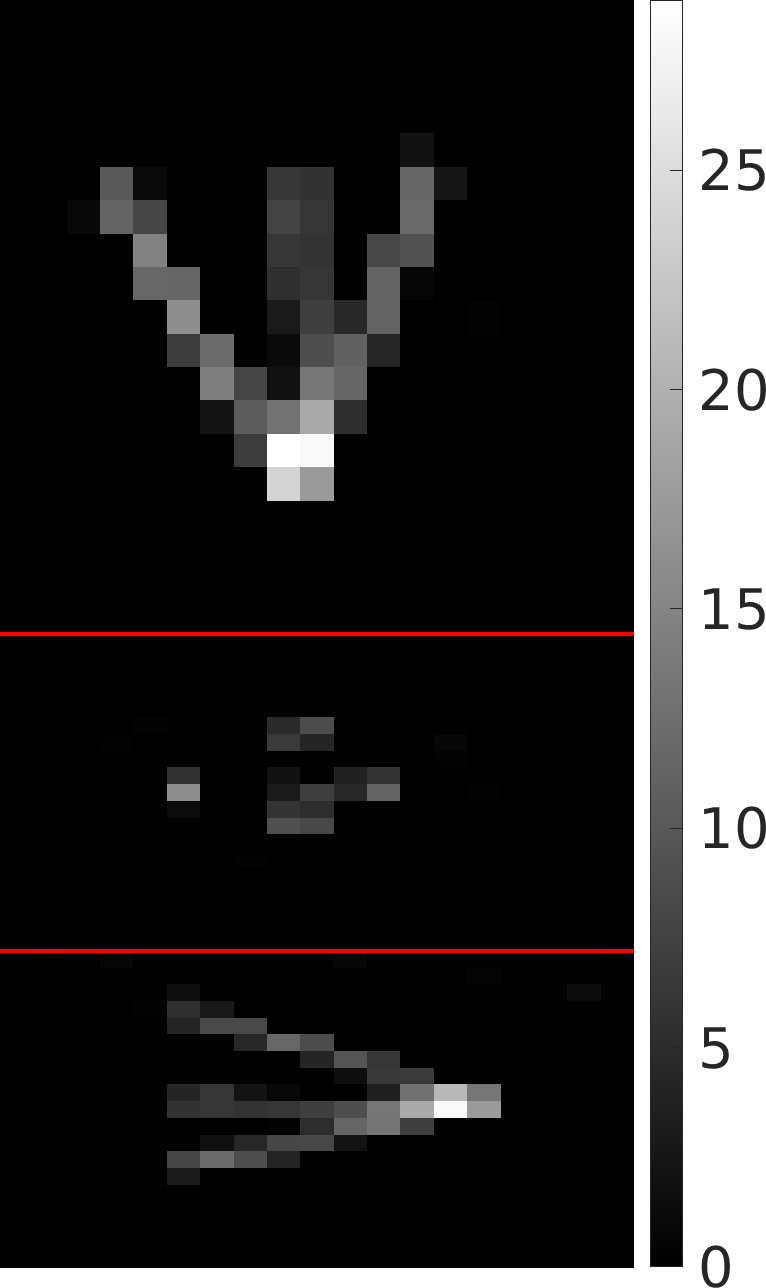}\\
 \hline
\multicolumn{6}{l}{$\tau=5$} \\
 \includegraphics[height=3.4cm]{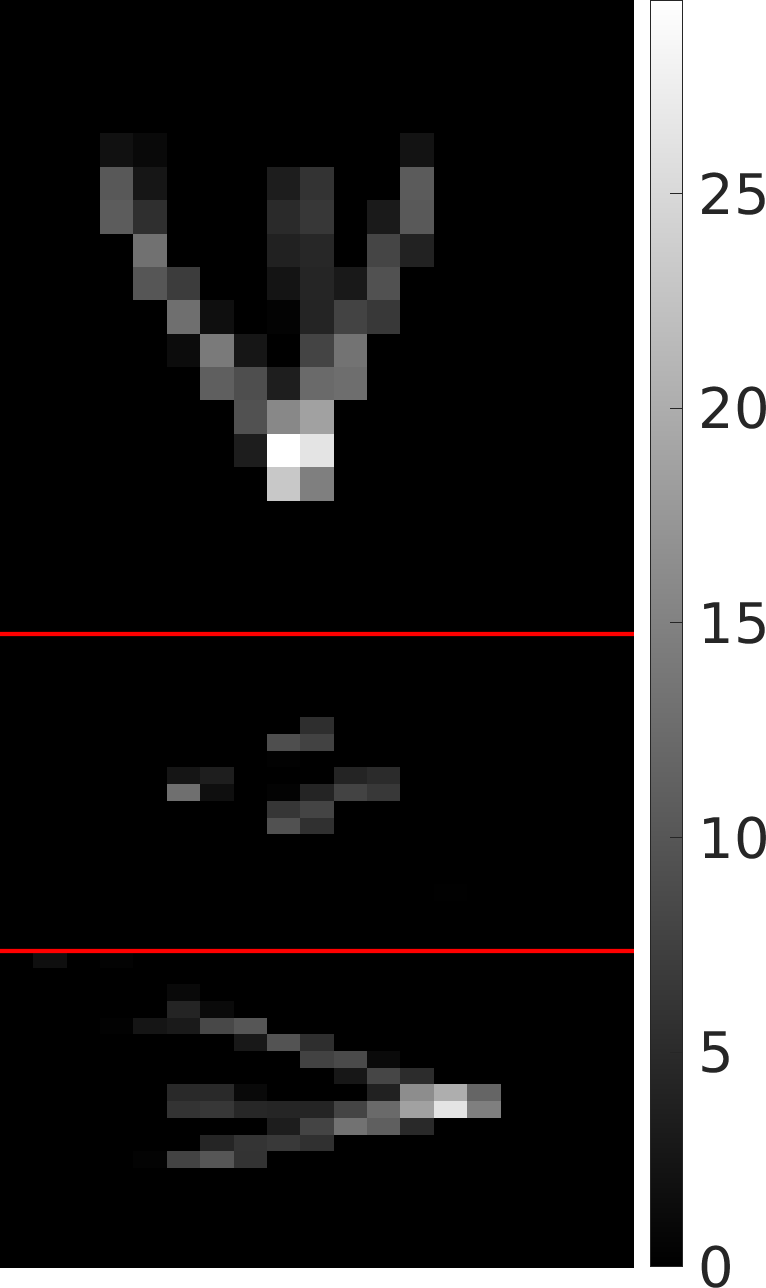}&
 \includegraphics[height=3.4cm]{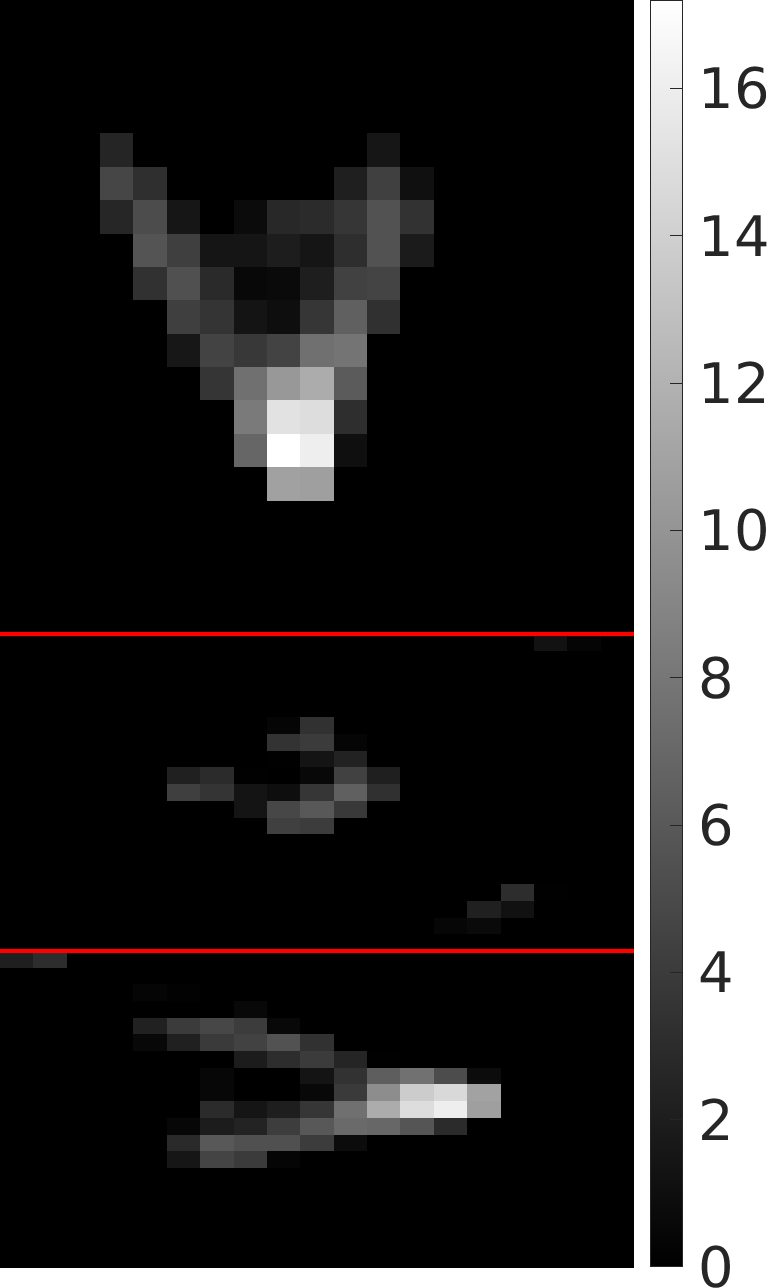}&
 \includegraphics[height=3.4cm]{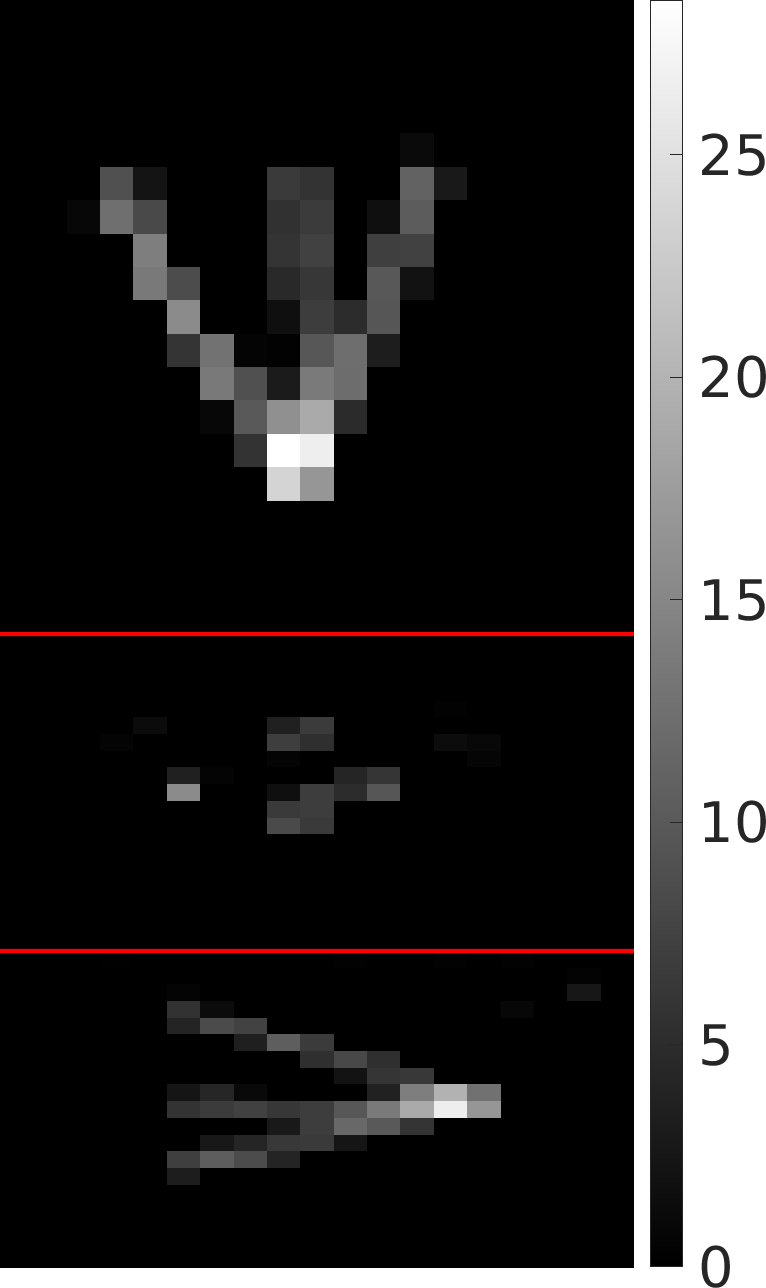}&
 \includegraphics[height=3.4cm]{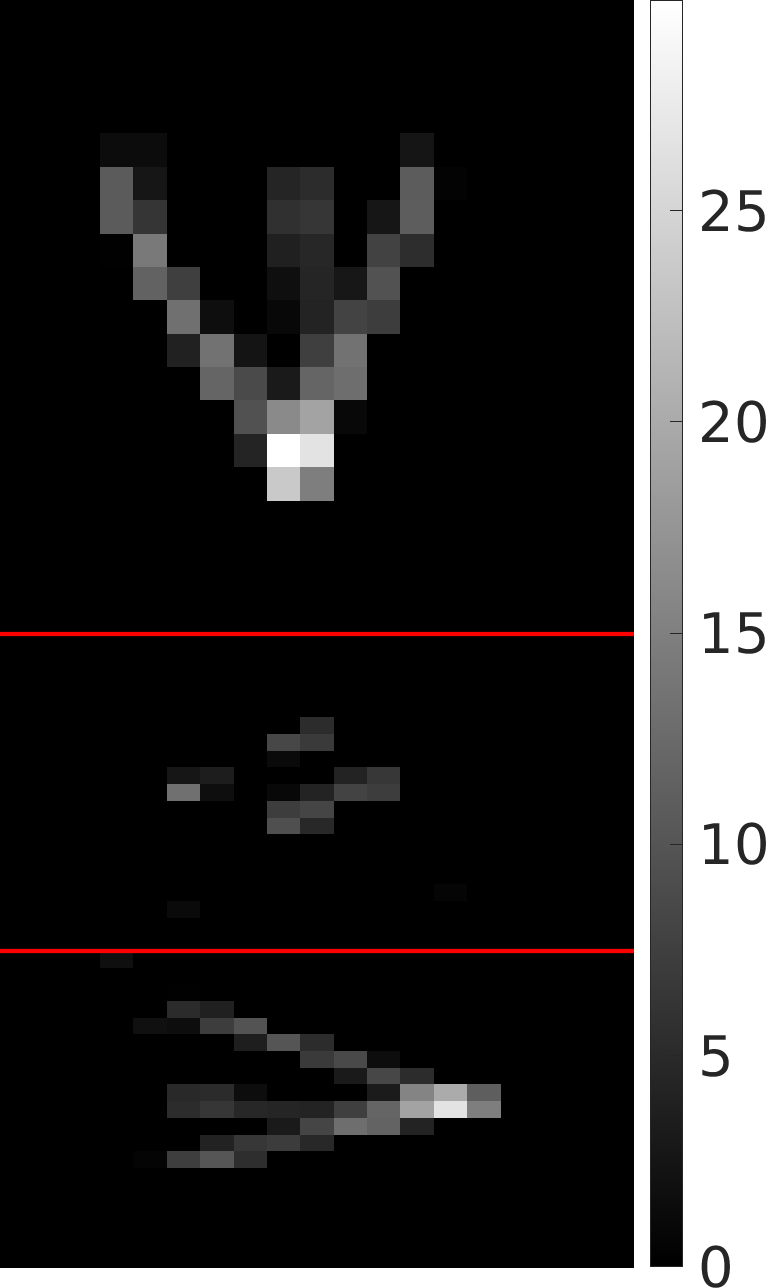}&
 \includegraphics[height=3.4cm]{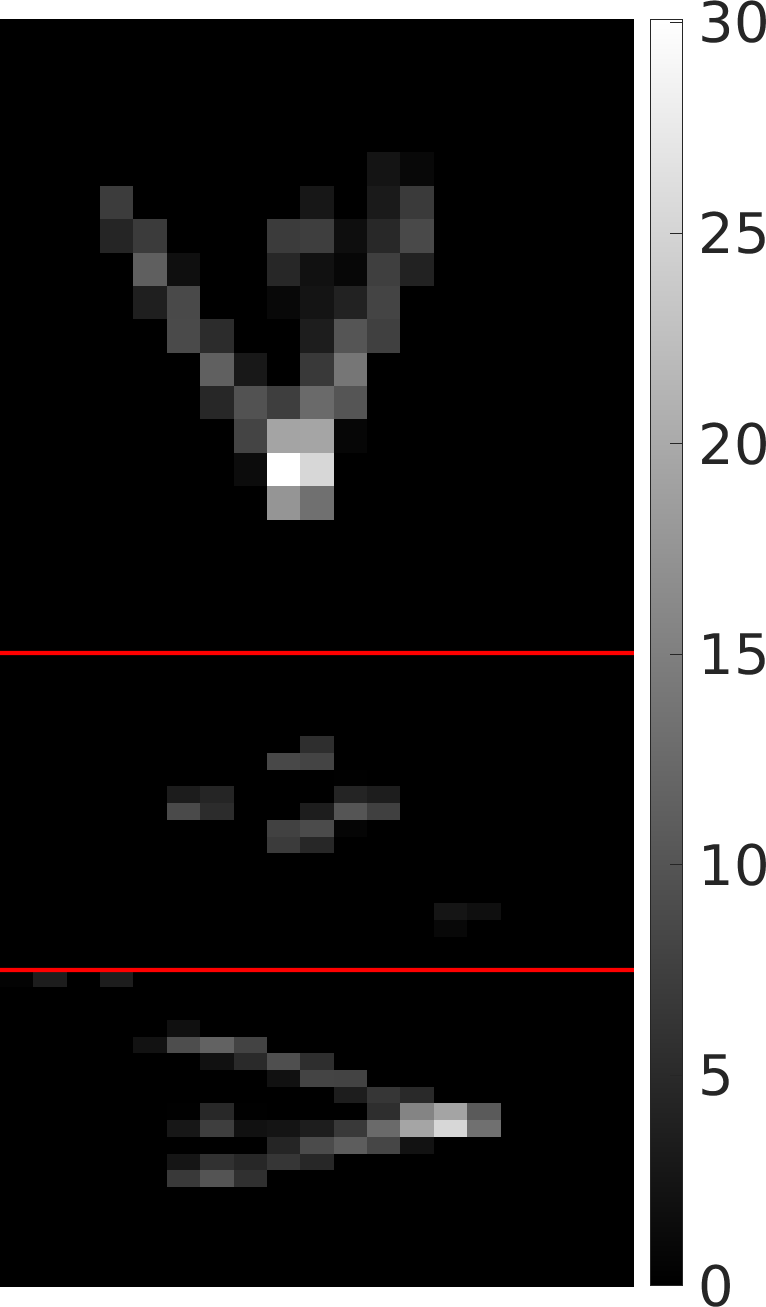}&
 \includegraphics[height=3.4cm]{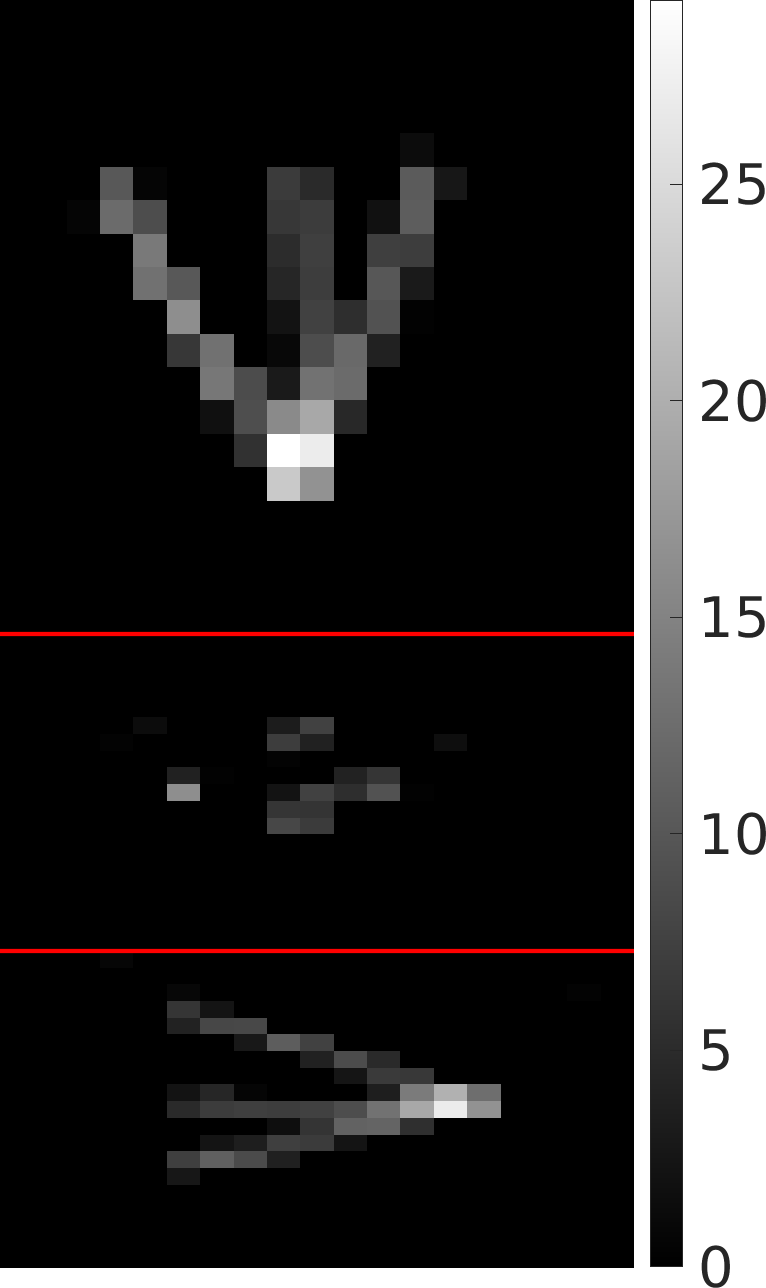}\\
\end{tabular}
}
\caption{``Resolution'' phantom reconstructions, PSNR-optimized $\alpha$ and iteration number $N$ (for l2-K only) according to Table \ref{tab:psnr_nonwhitened_vs_whitened}.}
\label{fig:methods_nonwhitened_vs_whitened_resolution_psnr}
\end{figure*}

\begin{table*}[hbt!]
\setlength{\tabcolsep}{1pt}
\centering
\caption{The $\epsilon_\mathrm{SSIM}$ values (dynamic range 100 mmol/l) for l1-L, l2-L, and l2-K.
 The numbers in brackets refer to $\alpha$, respectively $\alpha$ and the iteration number $N$ for l2-K. %
}\label{tab:ssim_nonwhitened_vs_whitened}
\scalebox{0.95}{
\begin{tabular}{c|c|c|c||c|c|c}
&\multicolumn{6}{c}{``Shape'' phantom} \\
\hline
& \multicolumn{3}{c||}{non-whitened} & \multicolumn{3}{c}{whitened} \\
\hline
  $\tau$ & l1-L & l2-L & l2-K & l1-L & l2-L & l2-K \\
\hline
0 &  $ 0.844 \ (2^{-6}) $  &  $ 0.830 \ (2^{-8}) $  &  $ \mathbf{0.933} \ (2^{-12}, 38) $  &  $ 0.846 \ (2^{-8}) $  &  $ 0.845 \ (2^{-8}) $  &  $ \mathbf{0.949} \ (2^{-11}, 16) $  \\
1 &  $ 0.915 \ (2^{-4}) $  &  $ 0.868 \ (2^{-6}) $  &  $ \mathbf{0.947} \ (2^{-13}, 11) $  &  $ 0.918 \ (2^{-3}) $  &  $ 0.925 \ (2^{-6}) $  &  $ \mathbf{0.963} \ (2^{-10}, 5) $  \\
3 &  $ \mathbf{0.970} \ (2^{-3}) $  &  $ 0.905 \ (2^{-9}) $  &  $ 0.929 \ (2^{-13}, 3) $  &  $ \mathbf{0.972} \ (2^{-2}) $  &  $ 0.943 \ (2^{-6}) $  &  $ 0.946 \ (2^{-9}, 3) $  \\
5 &  $ \mathbf{0.972} \ (2^{-4}) $  &  $ 0.897 \ (2^{-8}) $  &  $ 0.944 \ (2^{-12}, 4) $  &  $ \mathbf{0.973} \ (2^{-3}) $  &  $ 0.942 \ (2^{-7}) $  &  $ 0.964 \ (2^{-10}, 3) $  \\
\hline
\end{tabular}}

\vspace{0.5cm}

\scalebox{0.95}{
\begin{tabular}{c|c|c|c||c|c|c}
&\multicolumn{6}{c}{``Resolution'' phantom} \\
\hline
& \multicolumn{3}{c||}{non-whitened} & \multicolumn{3}{c}{whitened} \\
\hline
  $\tau$ & l1-L & l2-L & l2-K & l1-L & l2-L & l2-K \\
\hline
0 &  $ 0.939 \ (2^{-3}) $  &  $ 0.931 \ (2^{-10}) $  &  $ \mathbf{0.977} \ (2^{-17}, 22) $  &  $ 0.943 \ (2^{-2}) $  &  $ 0.939 \ (2^{-10}) $  &  $ \mathbf{0.980} \ (2^{-16}, 16) $  \\
1 &  $ 0.964 \ (2^{-3}) $  &  $ 0.942 \ (2^{-10}) $  &  $ \mathbf{0.980} \ (2^{-16}, 20) $  &  $ 0.968 \ (2^{-1}) $  &  $ 0.957 \ (2^{-9}) $  &  $ \mathbf{0.981} \ (2^{-13}, 12) $  \\
3 &  $ \mathbf{0.981} \ (2^{-1}) $  &  $ 0.957 \ (2^{-9}) $  &  $ 0.980 \ (2^{-15}, 12) $  &  $ \mathbf{0.981} \ (2^{0}) $  &  $ 0.968 \ (2^{-9}) $  &  $ 0.980 \ (2^{-13}, 10) $  \\
5 &  $ \mathbf{0.980} \ (2^{-1}) $  &  $ 0.957 \ (2^{-9}) $  &  $ \mathbf{0.980} \ (2^{-15}, 18) $  &  $ \mathbf{0.981} \ (2^{0}) $  &  $ 0.970 \ (2^{-10}) $  &  $ \mathbf{0.981} \ (2^{-13}, 15) $  \\
\hline
\end{tabular}}
\end{table*}

\begin{figure*}[hbt!]
\centering
\scalebox{0.85}{
\begin{tabular}{ccc|ccc}
\multicolumn{3}{c|}{non-whitened} & \multicolumn{3}{c}{whitened} \\
\hline
l1-L & l2-L & l2-K & l1-L & l2-L & l2-K \\
\hline
\multicolumn{6}{l}{$\tau=0$} \\
 \includegraphics[height=3.4cm]{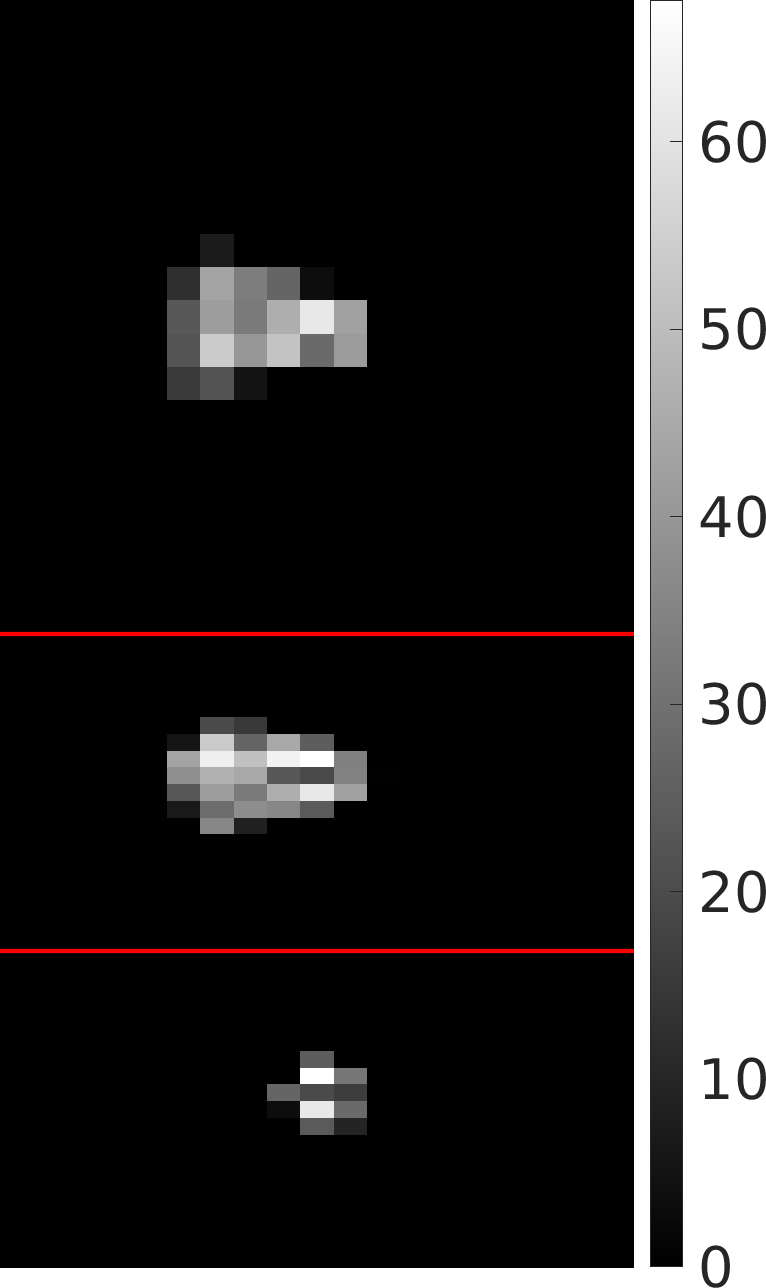}&
 \includegraphics[height=3.4cm]{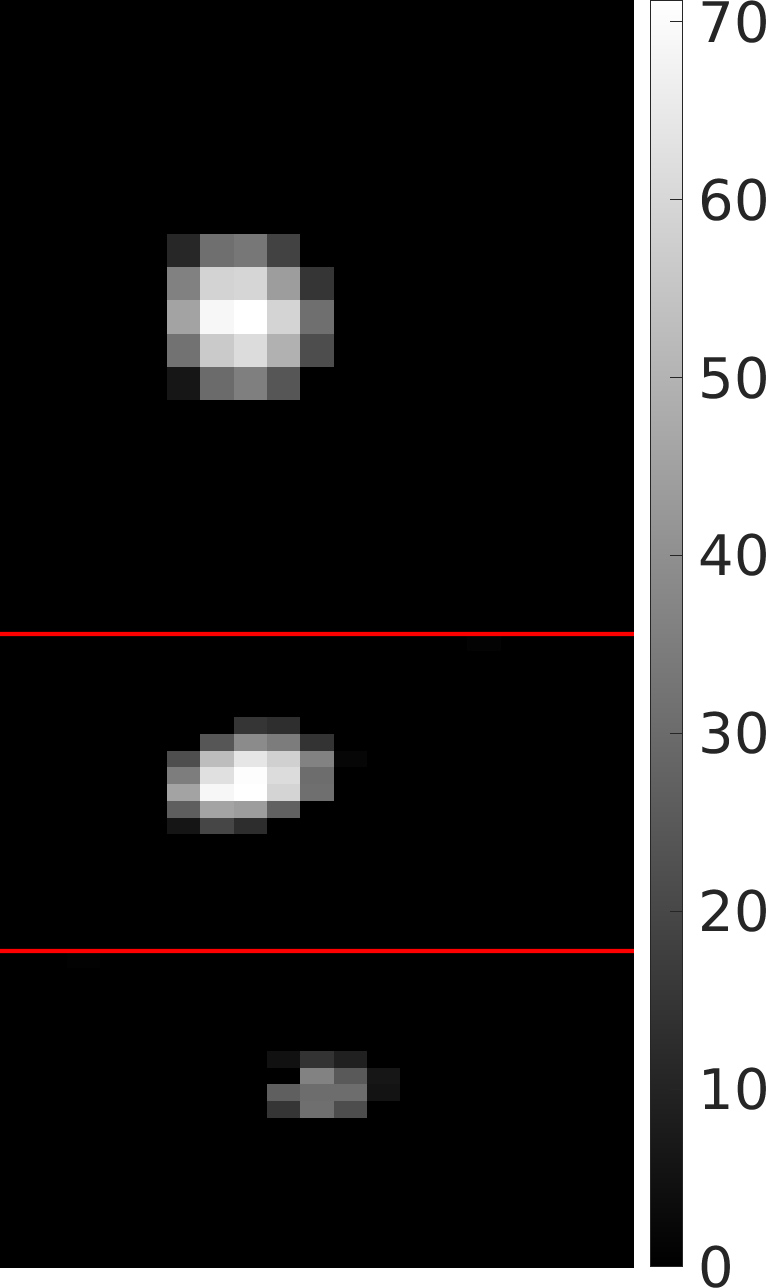}&
 \includegraphics[height=3.4cm]{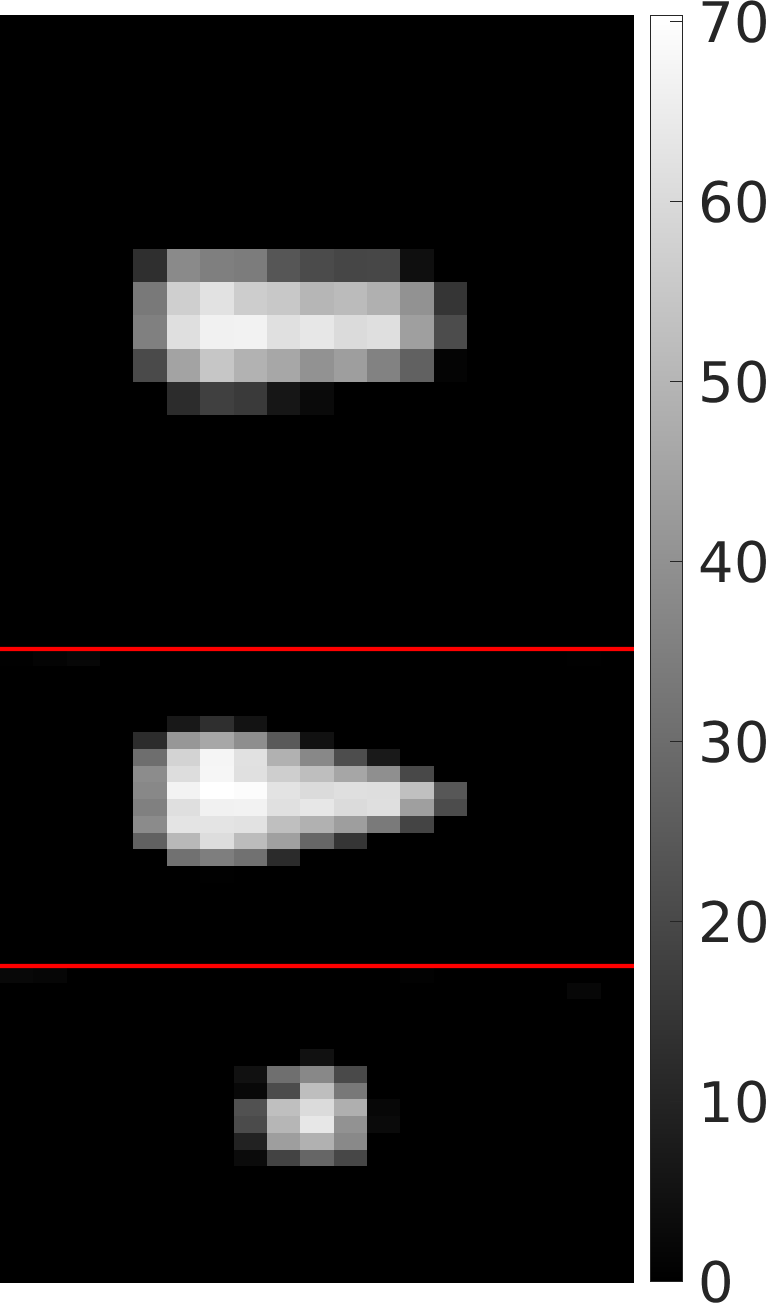}&
 \includegraphics[height=3.4cm]{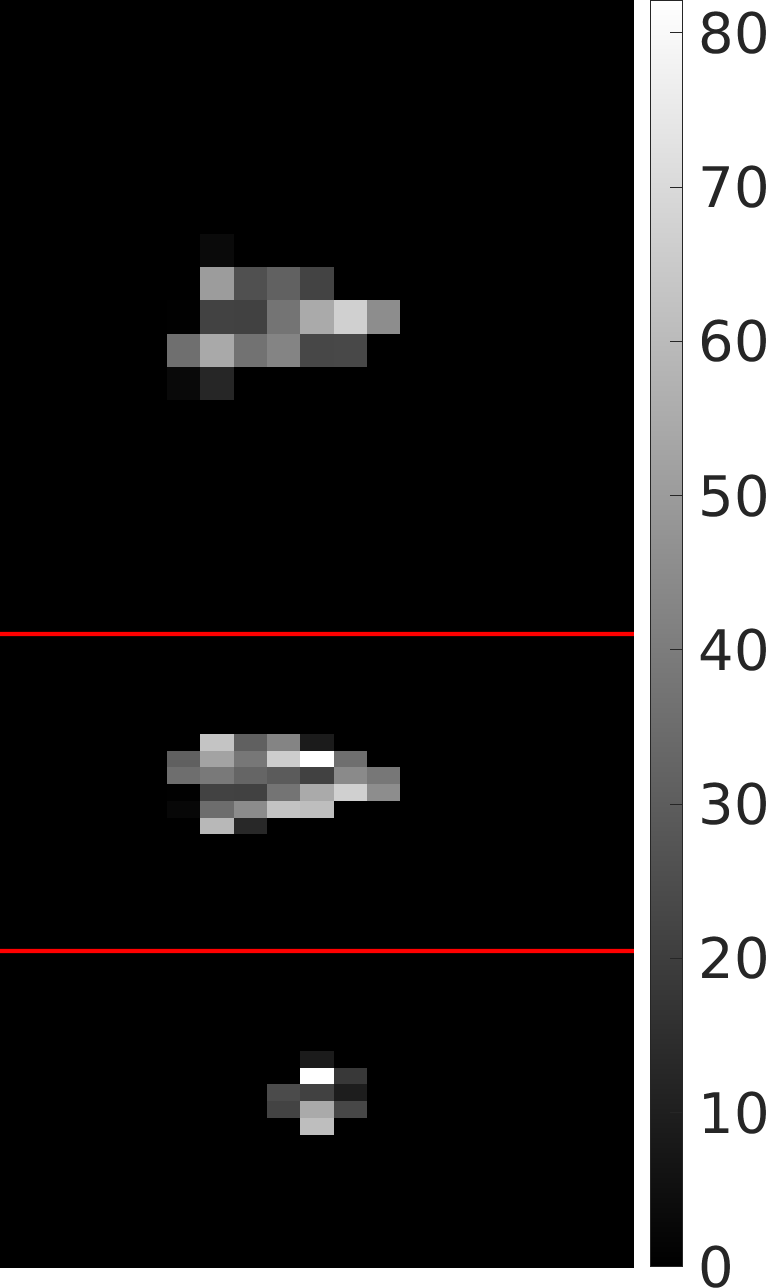}&
 \includegraphics[height=3.4cm]{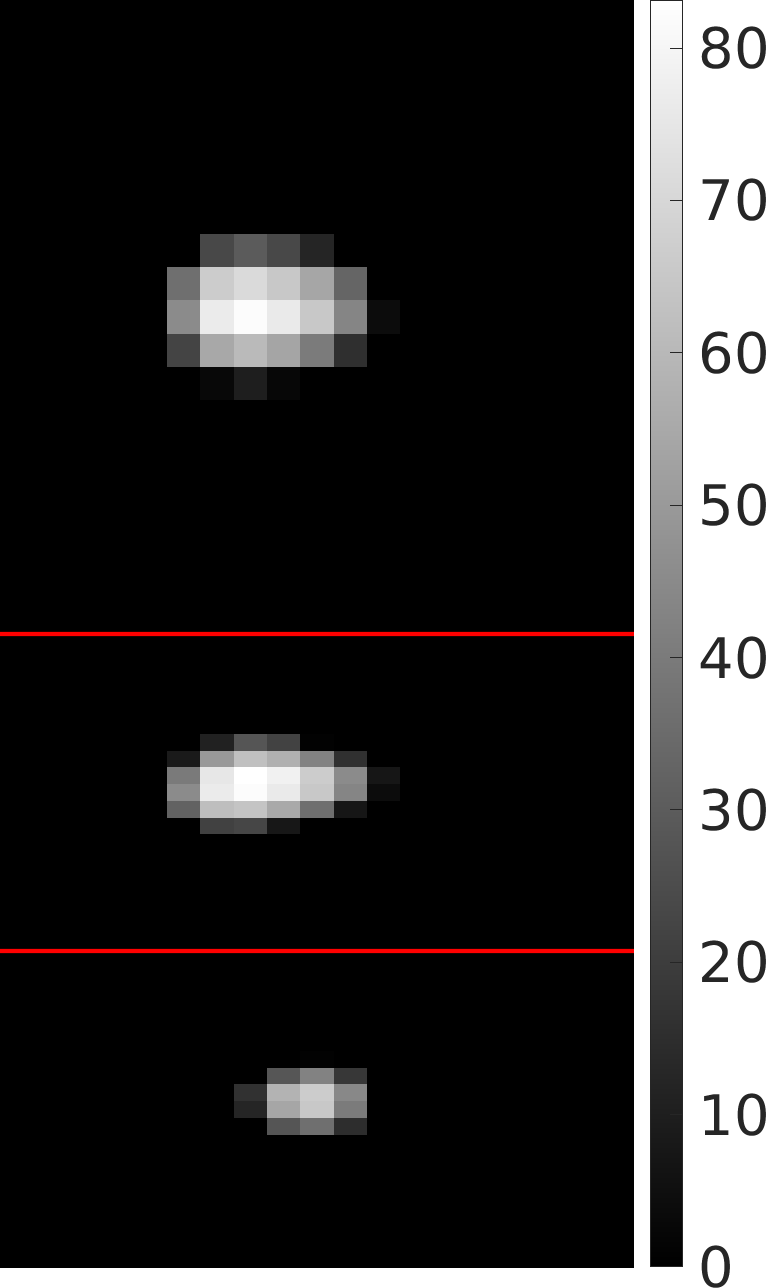}&
 \includegraphics[height=3.4cm]{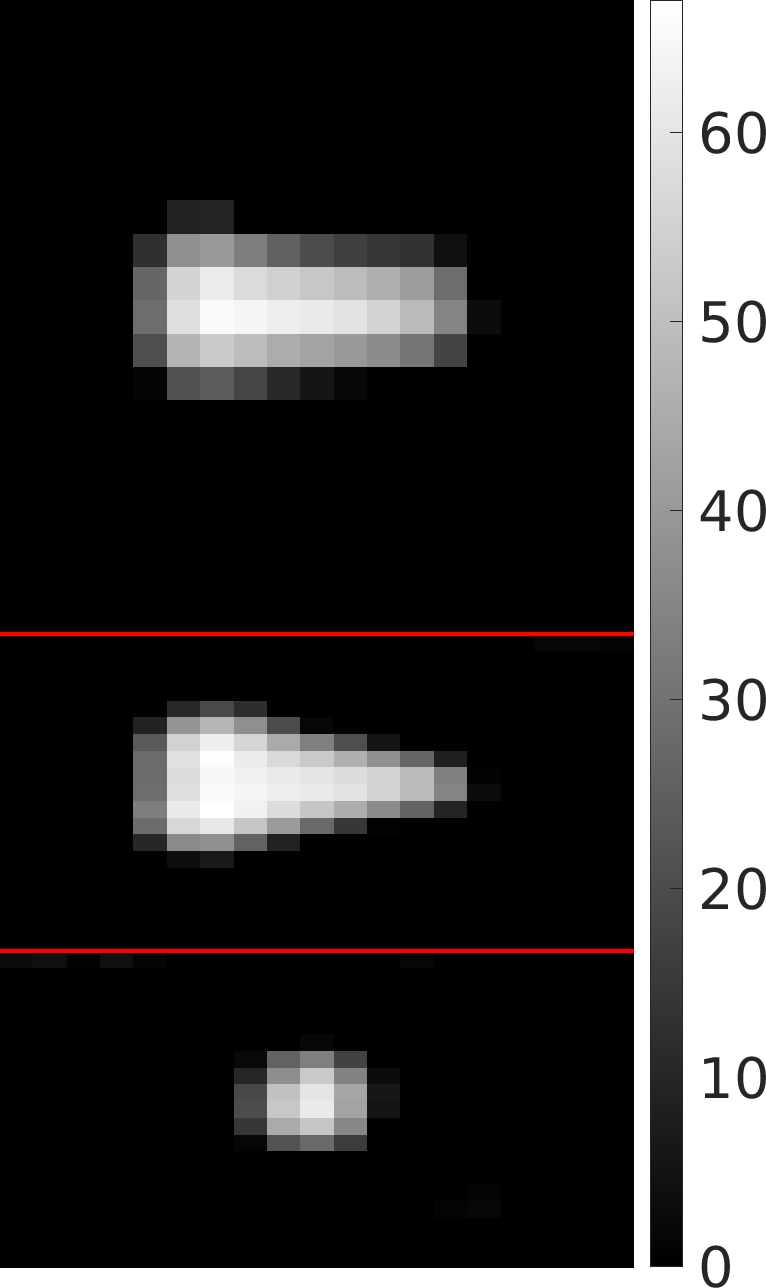}\\
\hline
\multicolumn{6}{l}{$\tau=1$} \\
 \includegraphics[height=3.4cm]{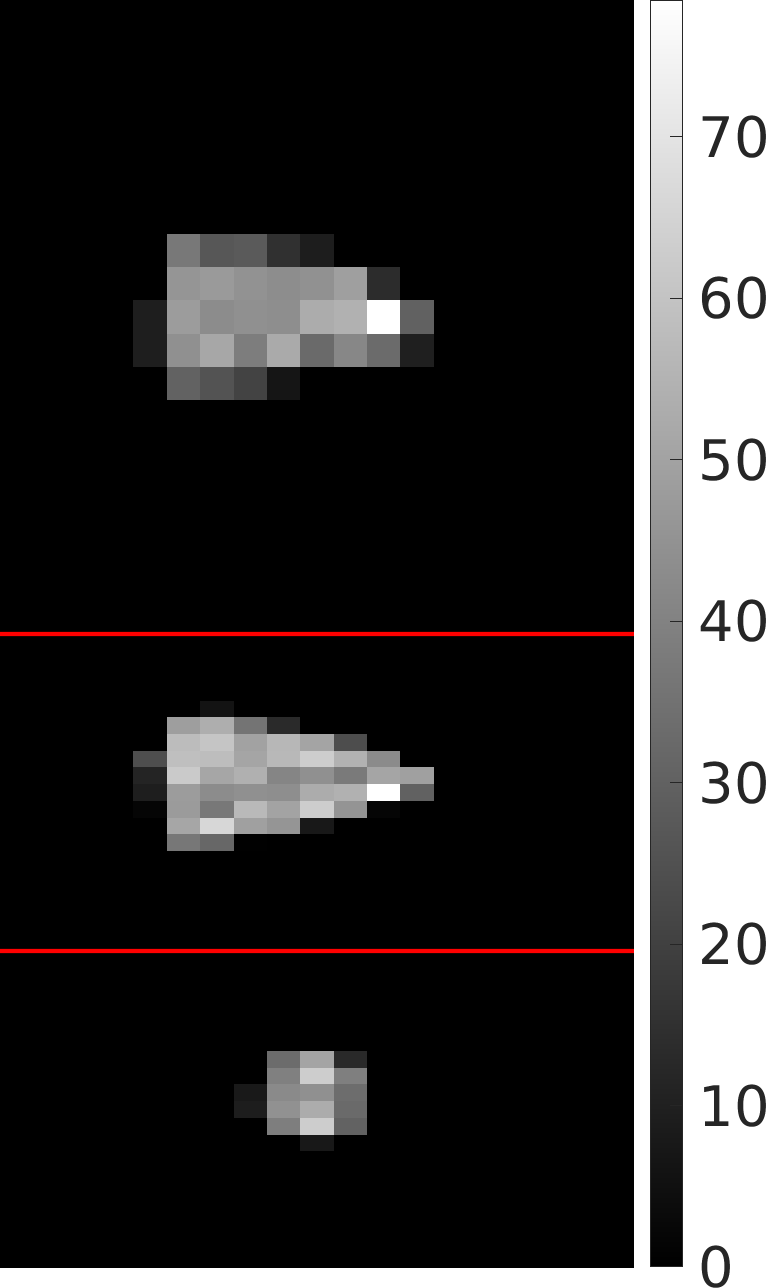}&
 \includegraphics[height=3.4cm]{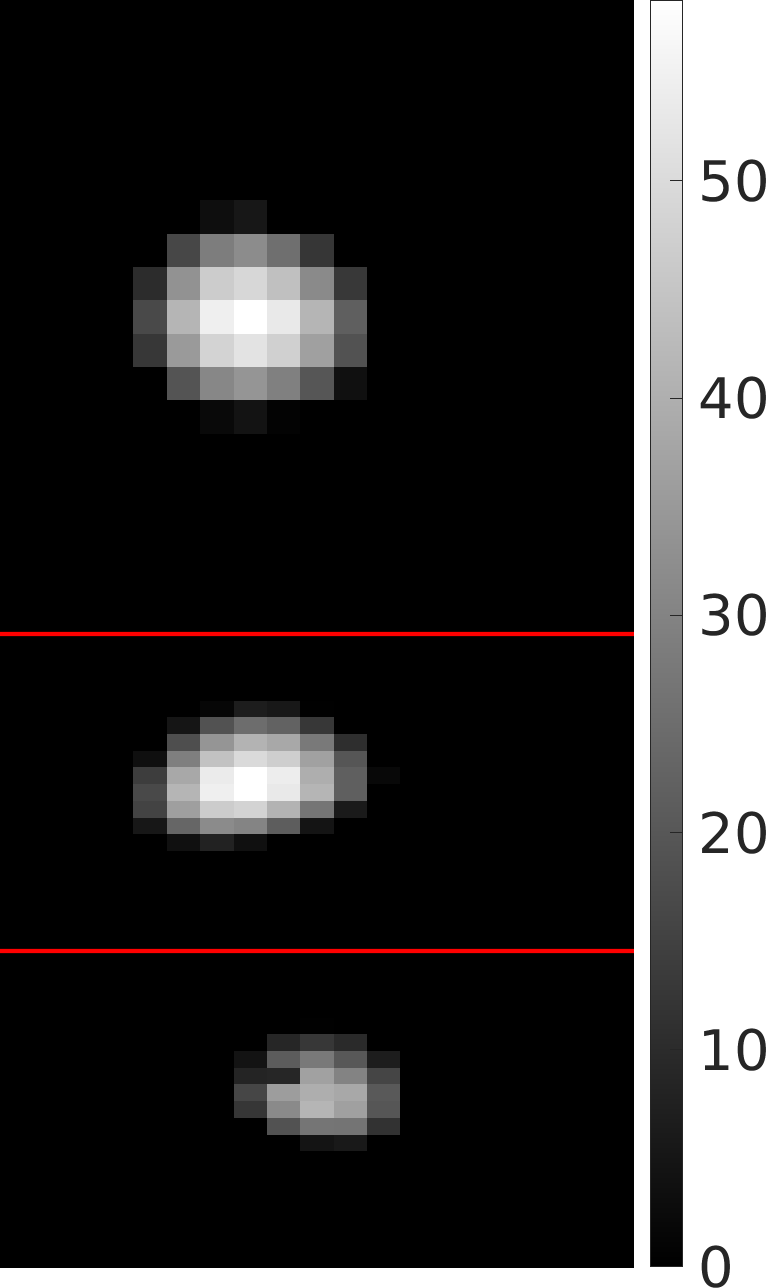}&
 \includegraphics[height=3.4cm]{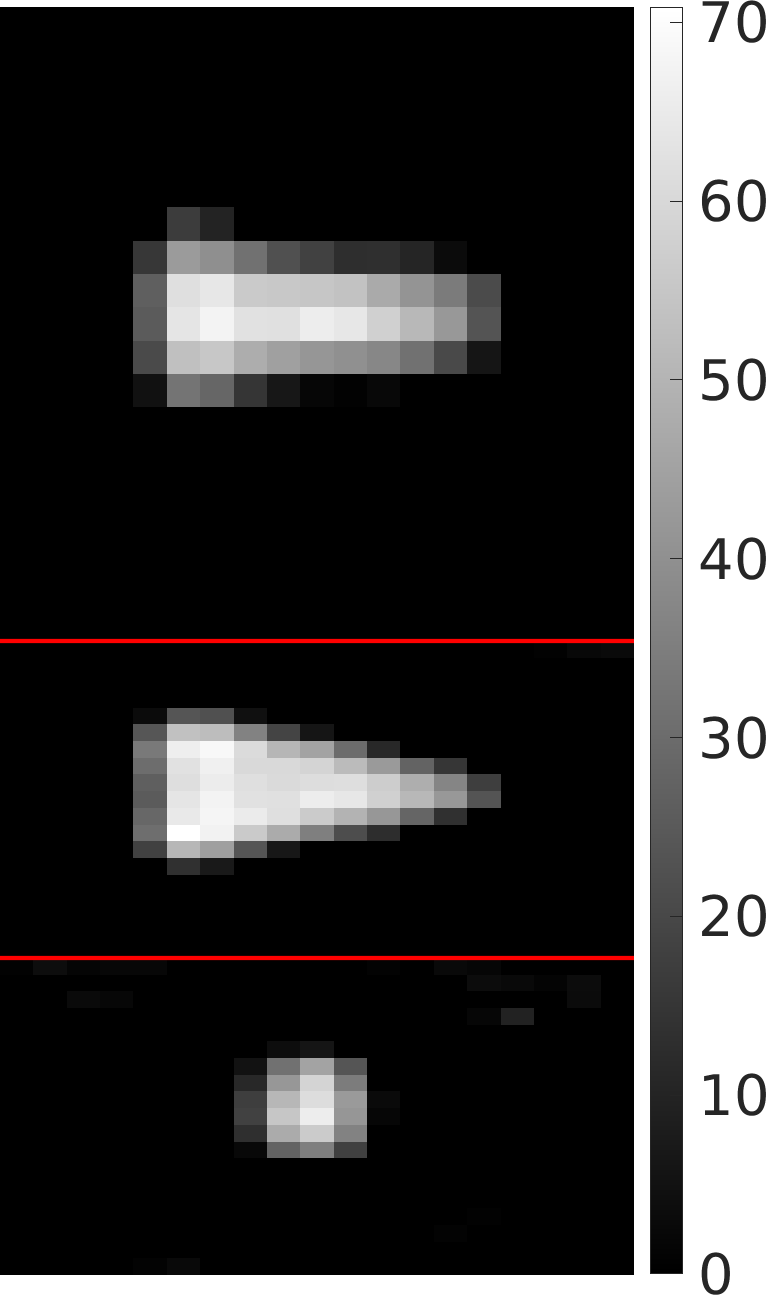}&
 \includegraphics[height=3.4cm]{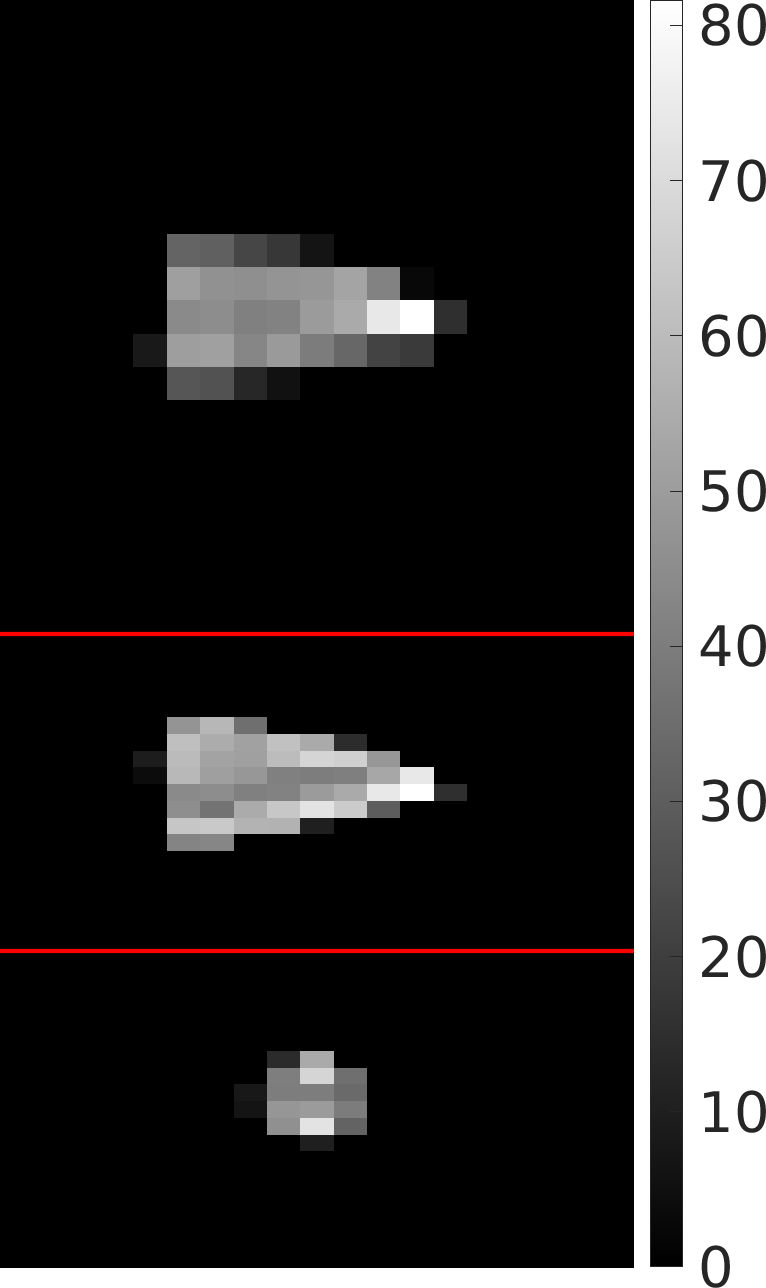}&
 \includegraphics[height=3.4cm]{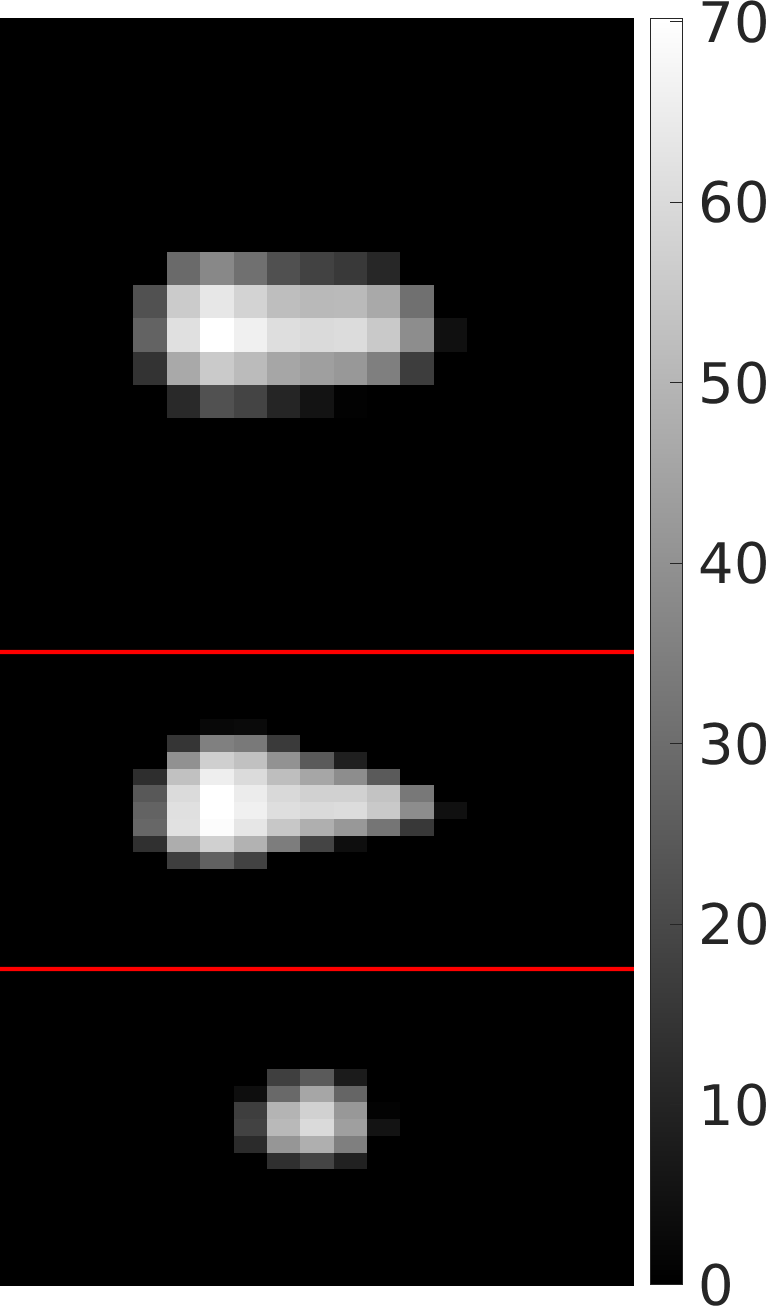}&
 \includegraphics[height=3.4cm]{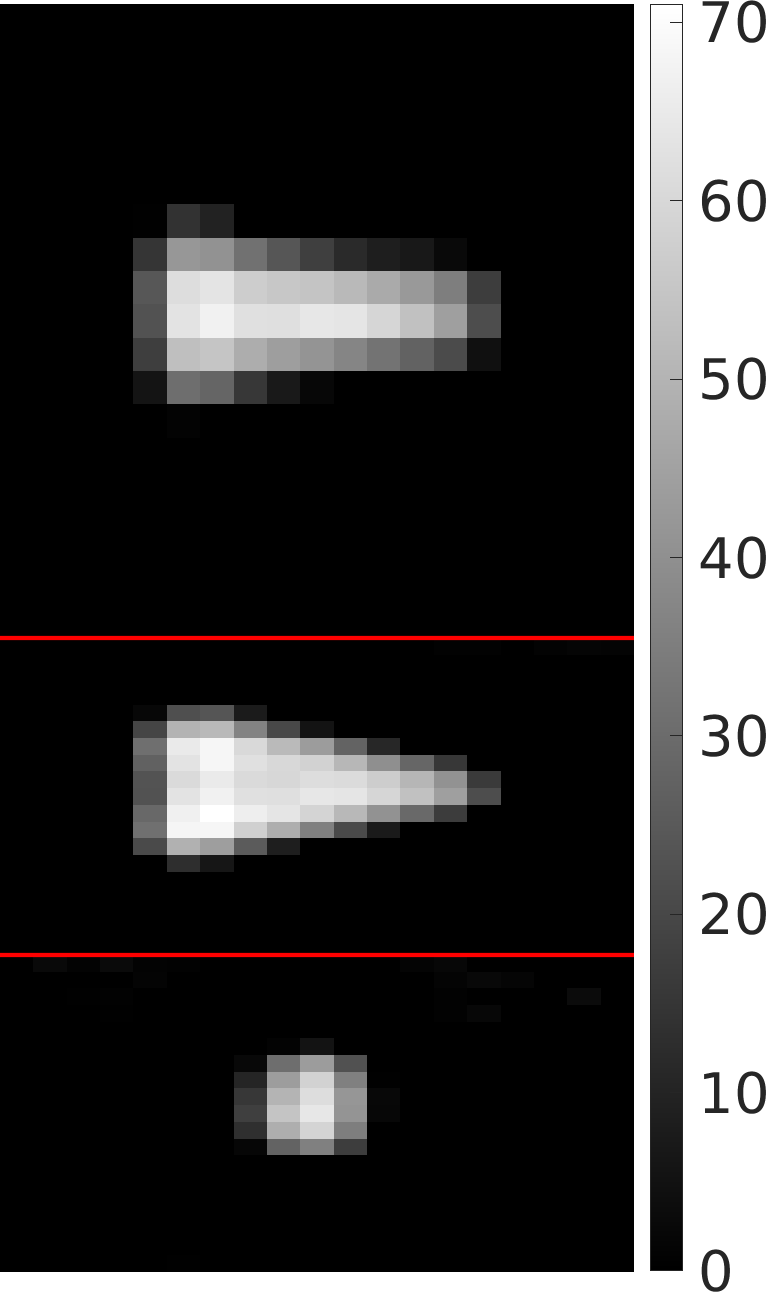}\\
\hline
\multicolumn{6}{l}{$\tau=3$} \\
 \includegraphics[height=3.4cm]{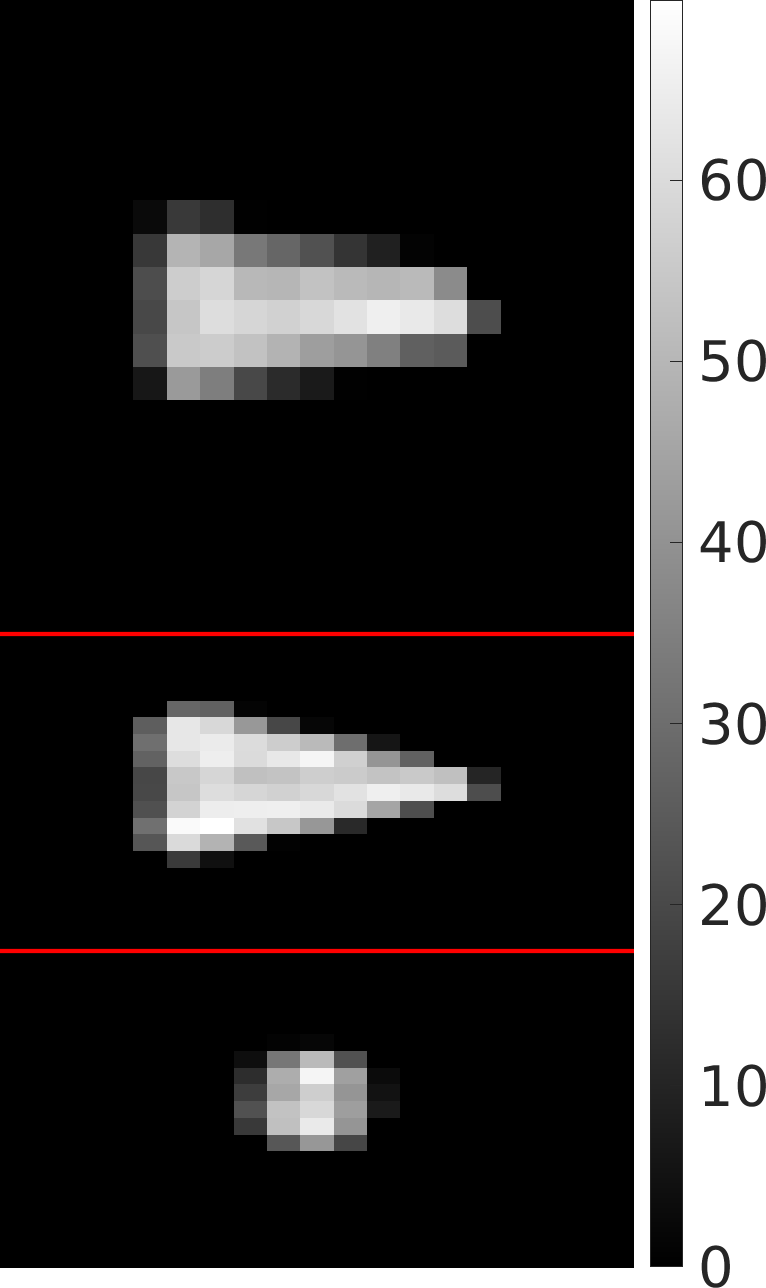}&
 \includegraphics[height=3.4cm]{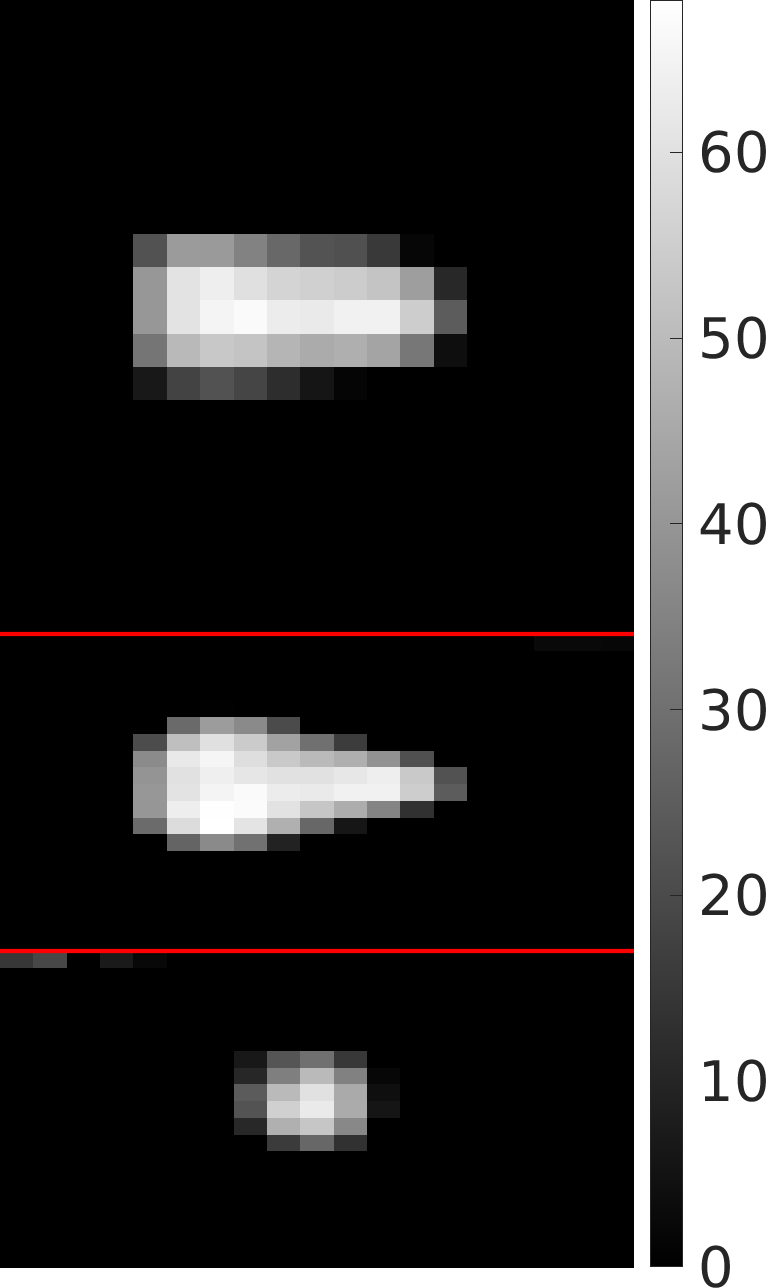}&
 \includegraphics[height=3.4cm]{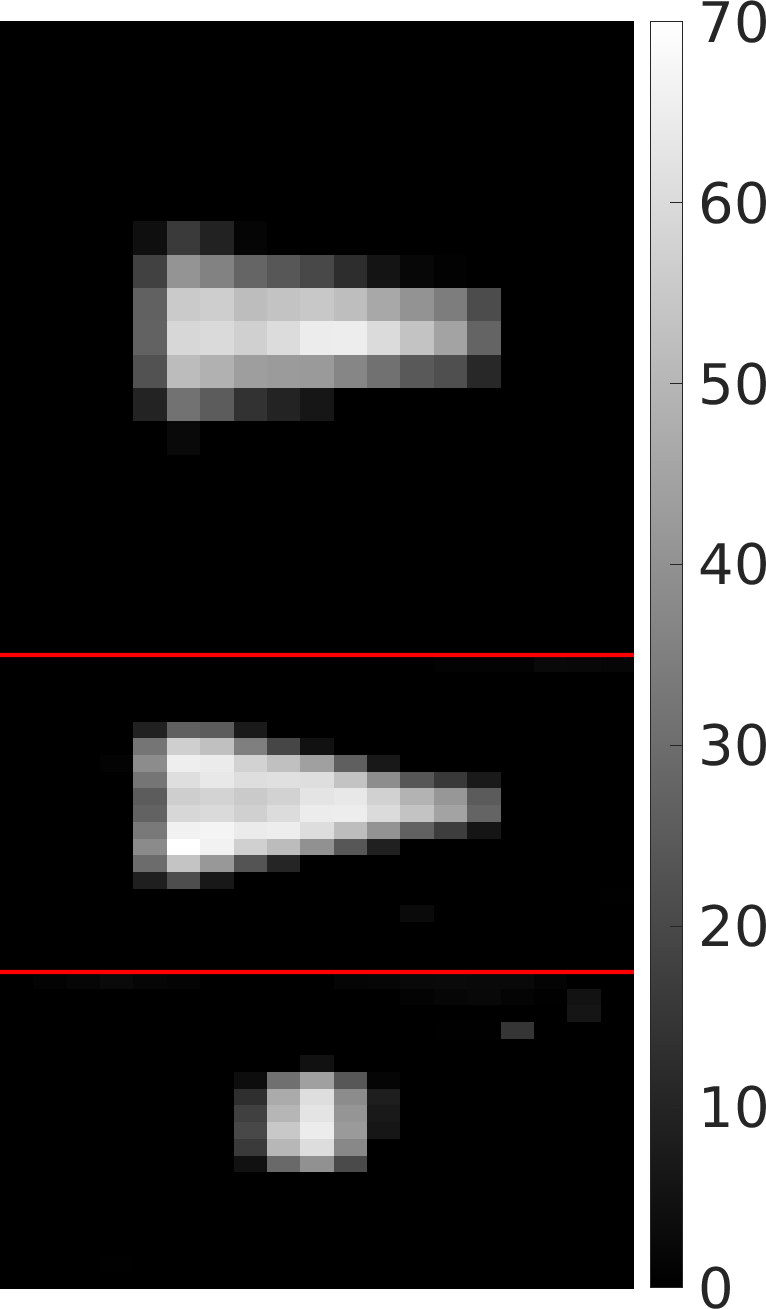}&
 \includegraphics[height=3.4cm]{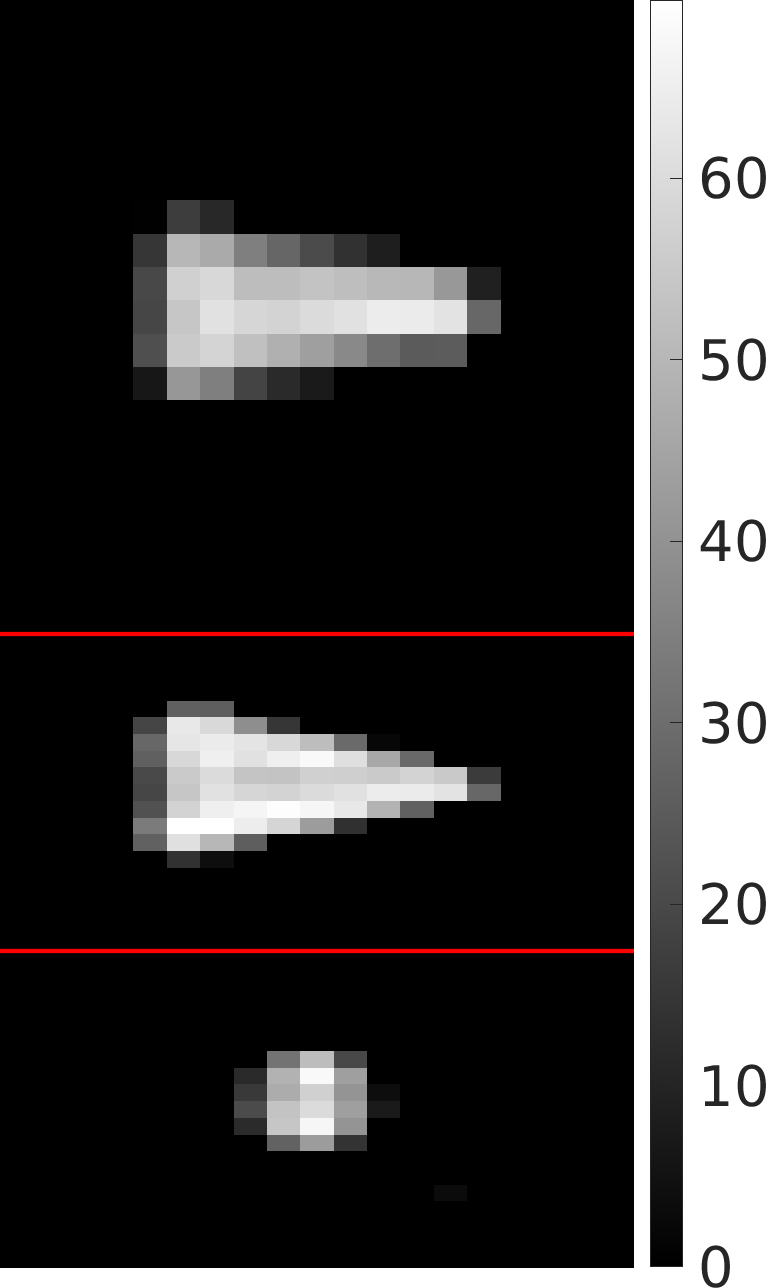}&
 \includegraphics[height=3.4cm]{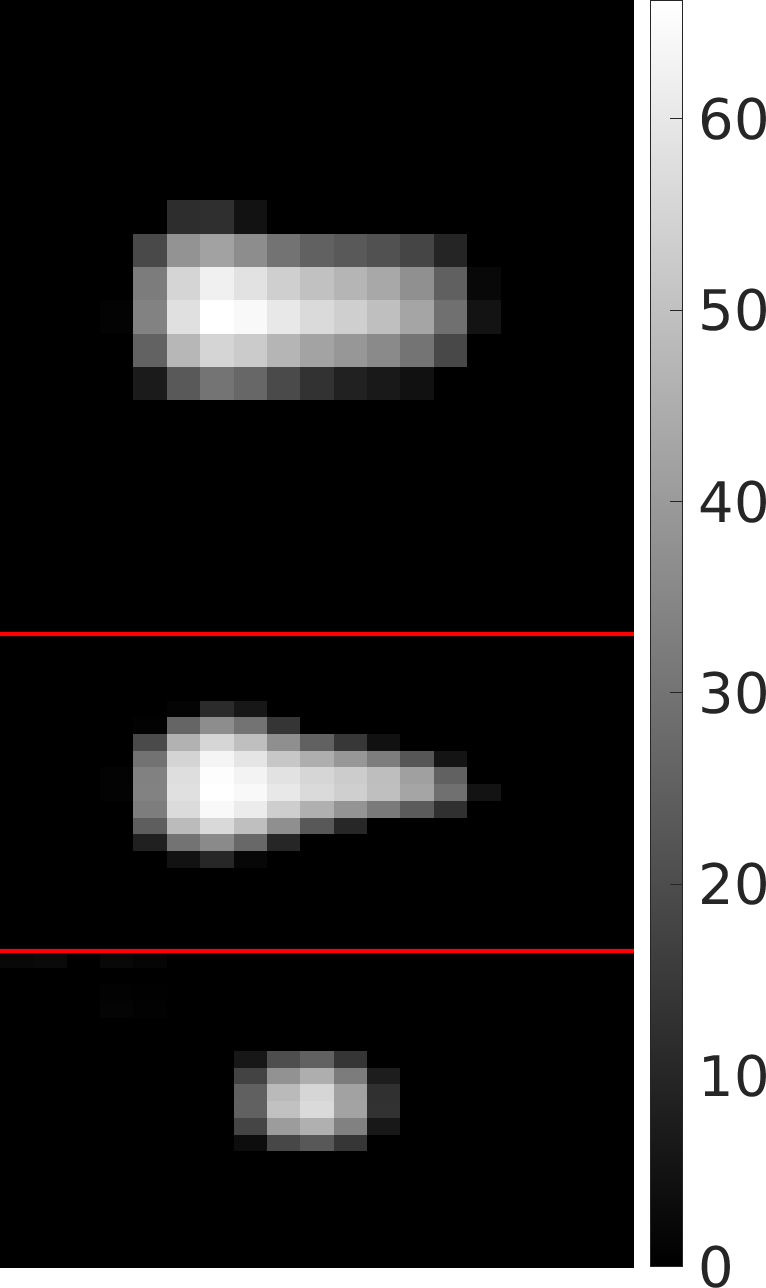}&
 \includegraphics[height=3.4cm]{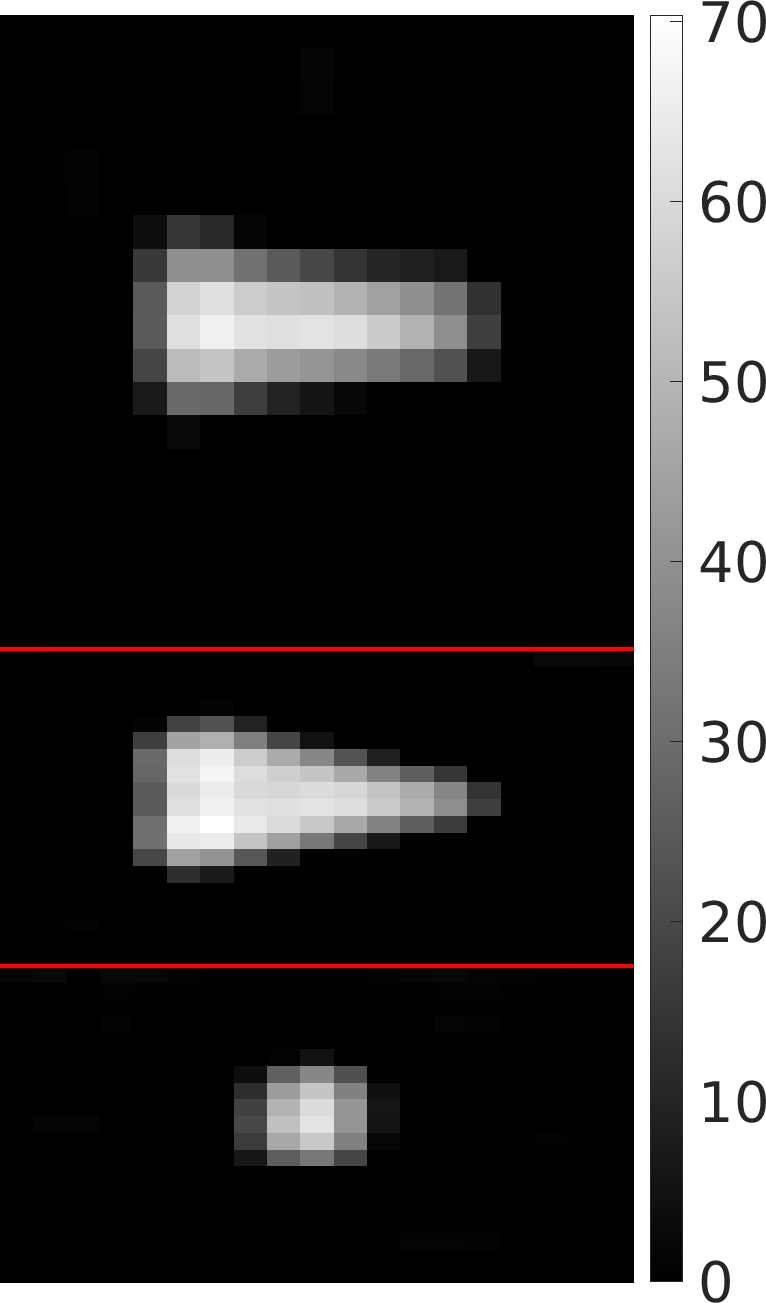}\\
 \hline
\multicolumn{6}{l}{$\tau=5$} \\
 \includegraphics[height=3.4cm]{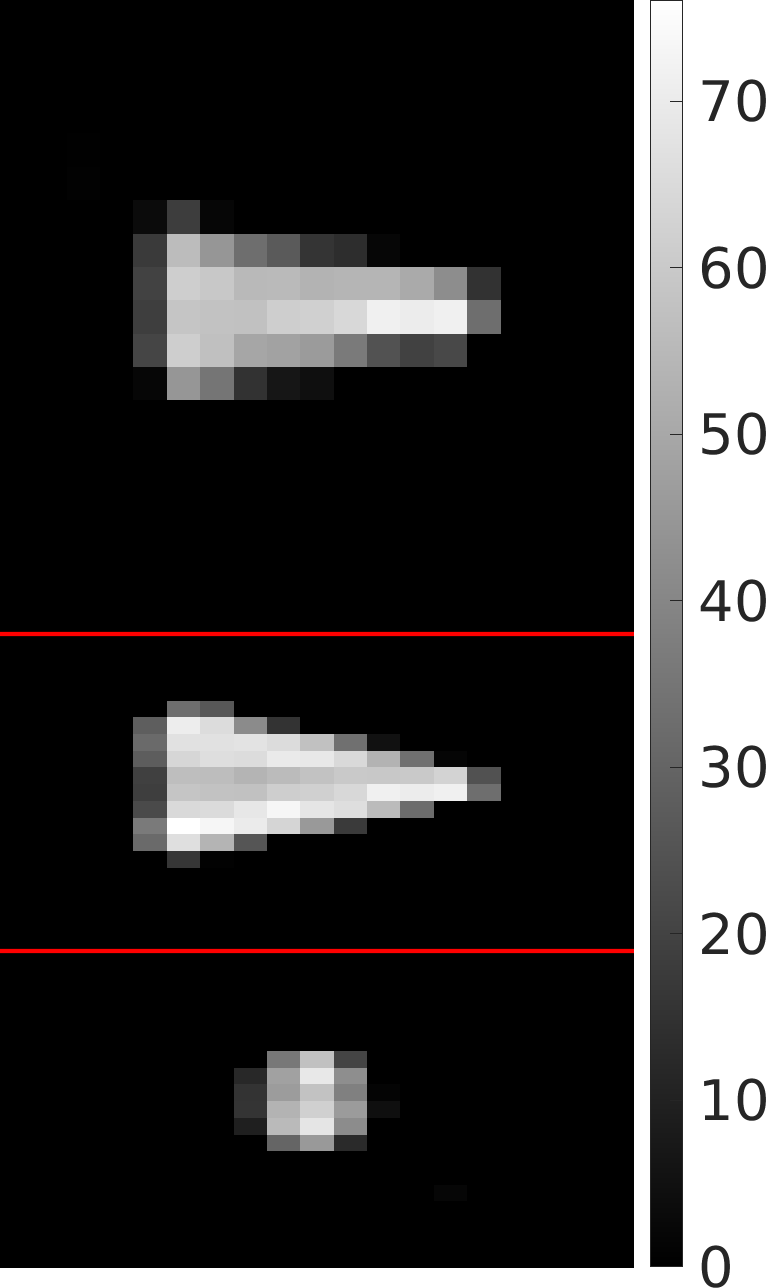}&
 \includegraphics[height=3.4cm]{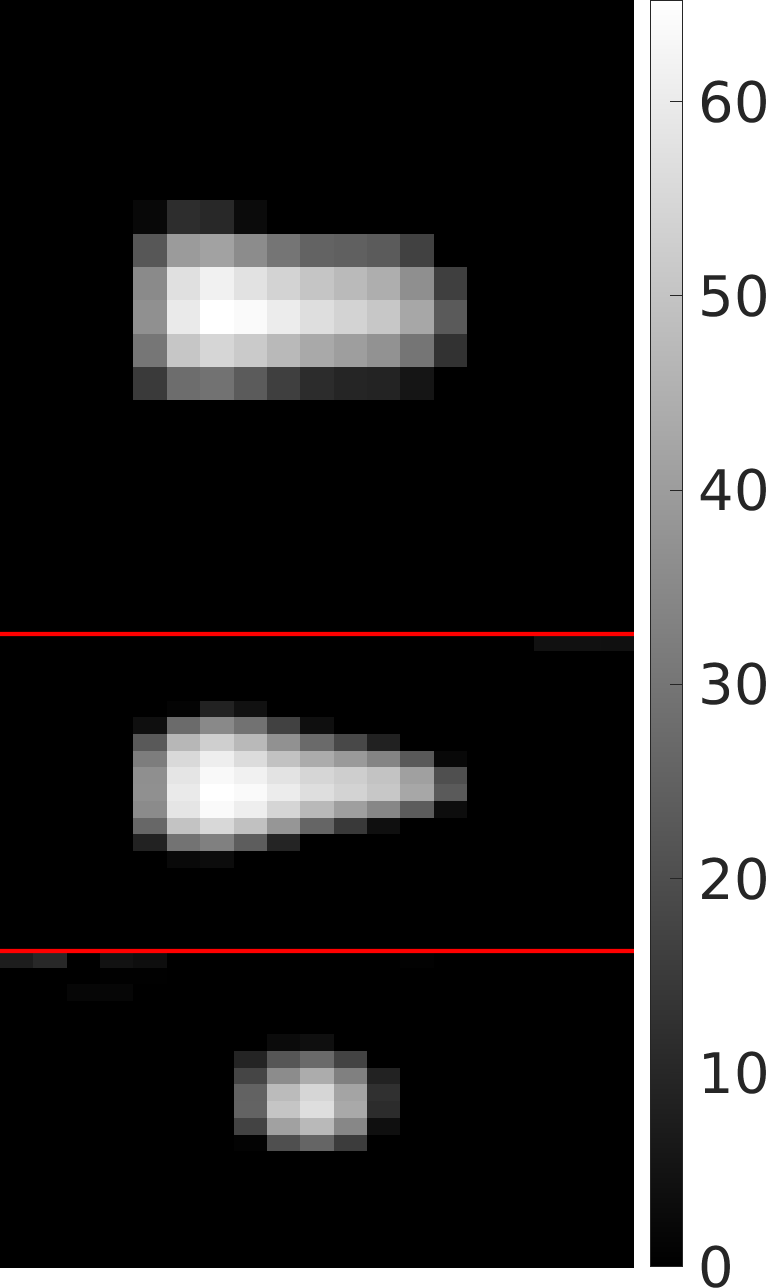}&
 \includegraphics[height=3.4cm]{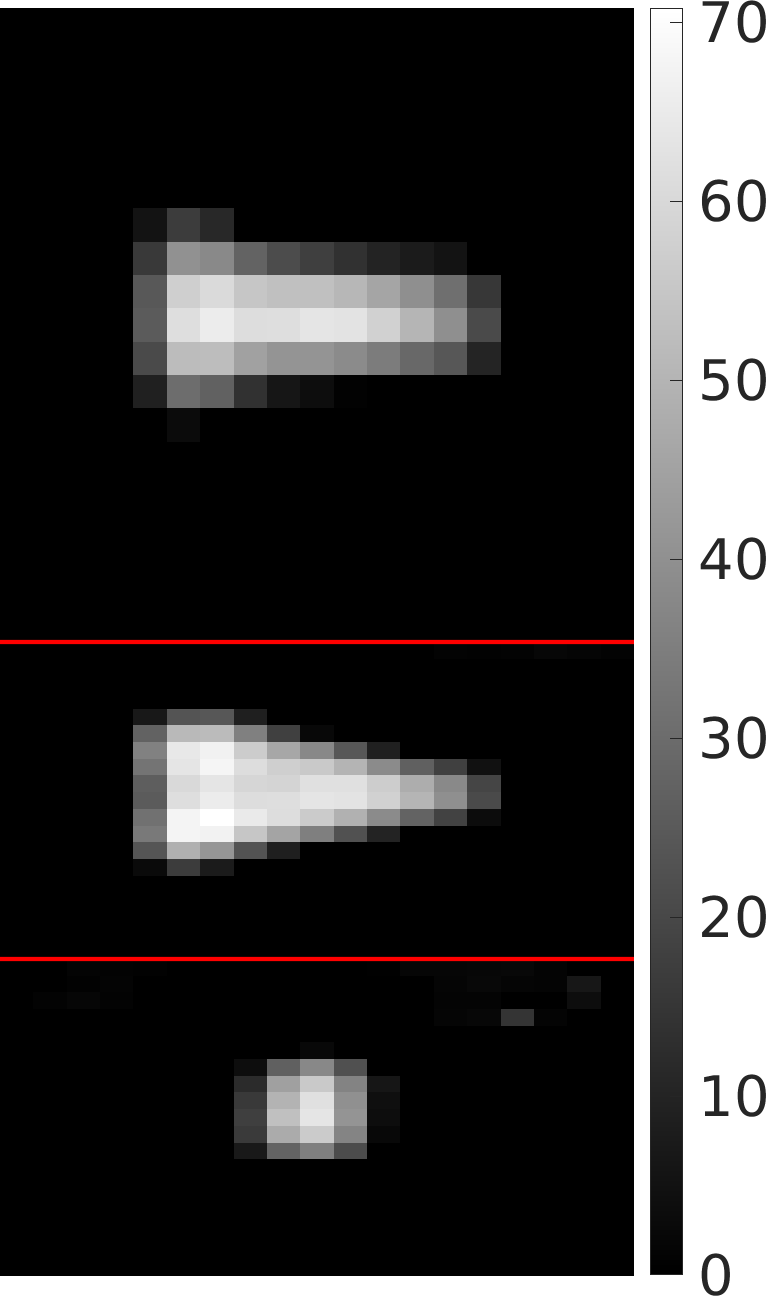}&
 \includegraphics[height=3.4cm]{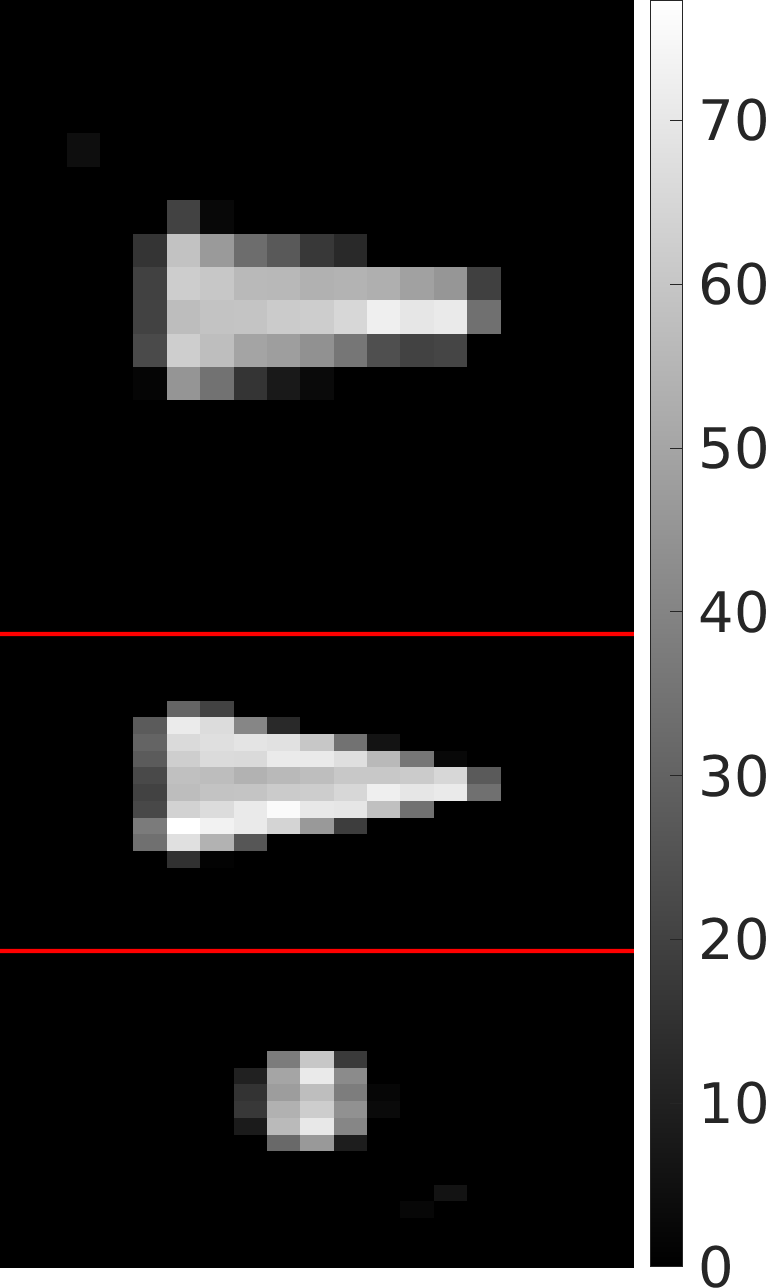}&
 \includegraphics[height=3.4cm]{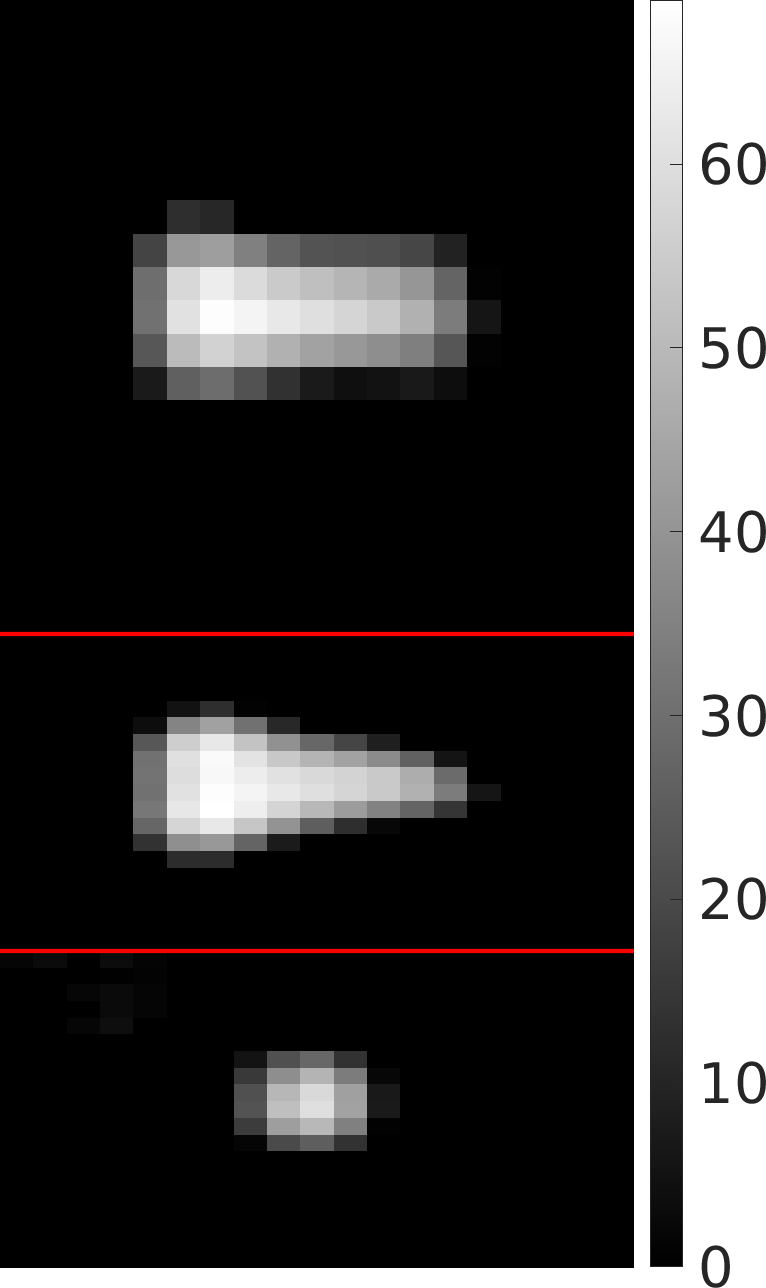}&
 \includegraphics[height=3.4cm]{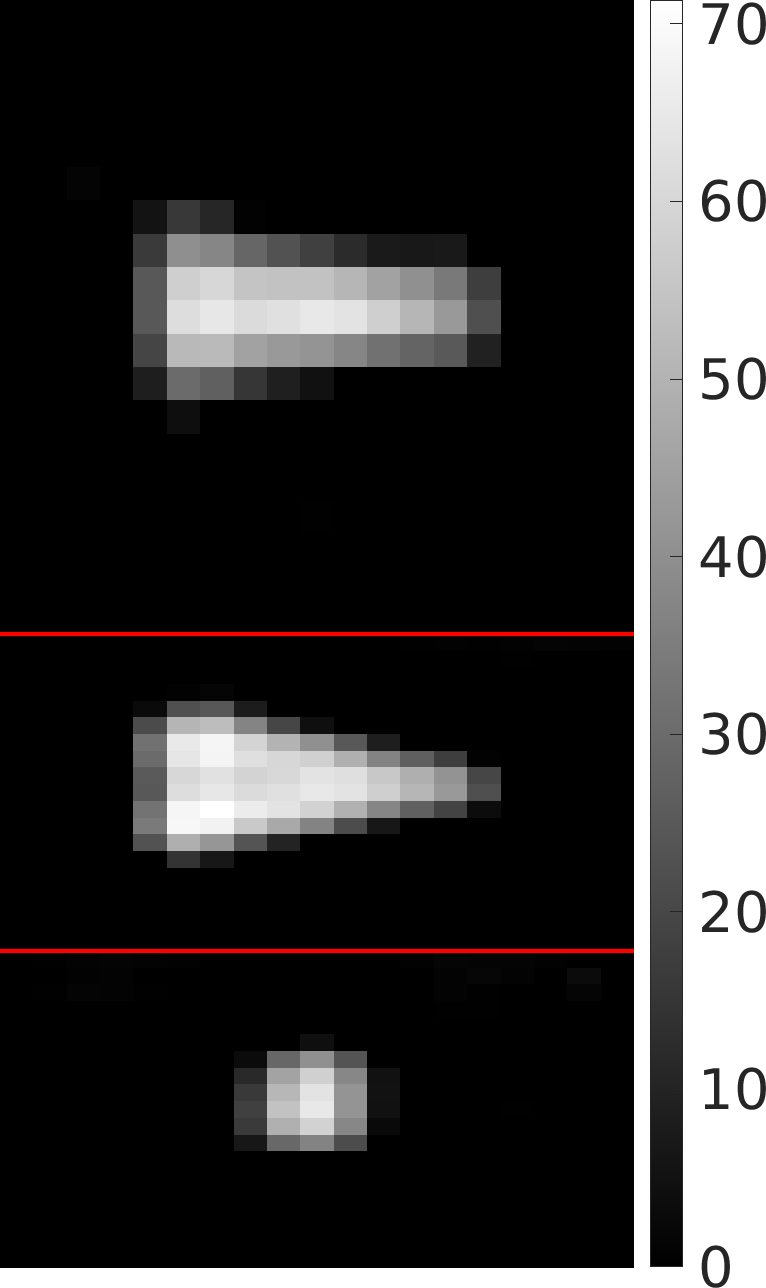}\\
\end{tabular}
}
\caption{``Shape'' phantom reconstructions, SSIM-optimized $\alpha$ and iteration number $N$ (for l2-K only) according to Table \ref{tab:ssim_nonwhitened_vs_whitened}. }
\label{fig:methods_nonwhitened_vs_whitened_shape_ssim}
\end{figure*}

\begin{figure*}[hbt!]%
\centering
\scalebox{0.85}{
\begin{tabular}{ccc|ccc}
\multicolumn{3}{c|}{non-whitened} & \multicolumn{3}{c}{whitened} \\
\hline
l1-L & l2-L & l2-K & l1-L & l2-L & l2-K \\
\hline
\multicolumn{6}{l}{$\tau=0$} \\
 \includegraphics[height=3.4cm]{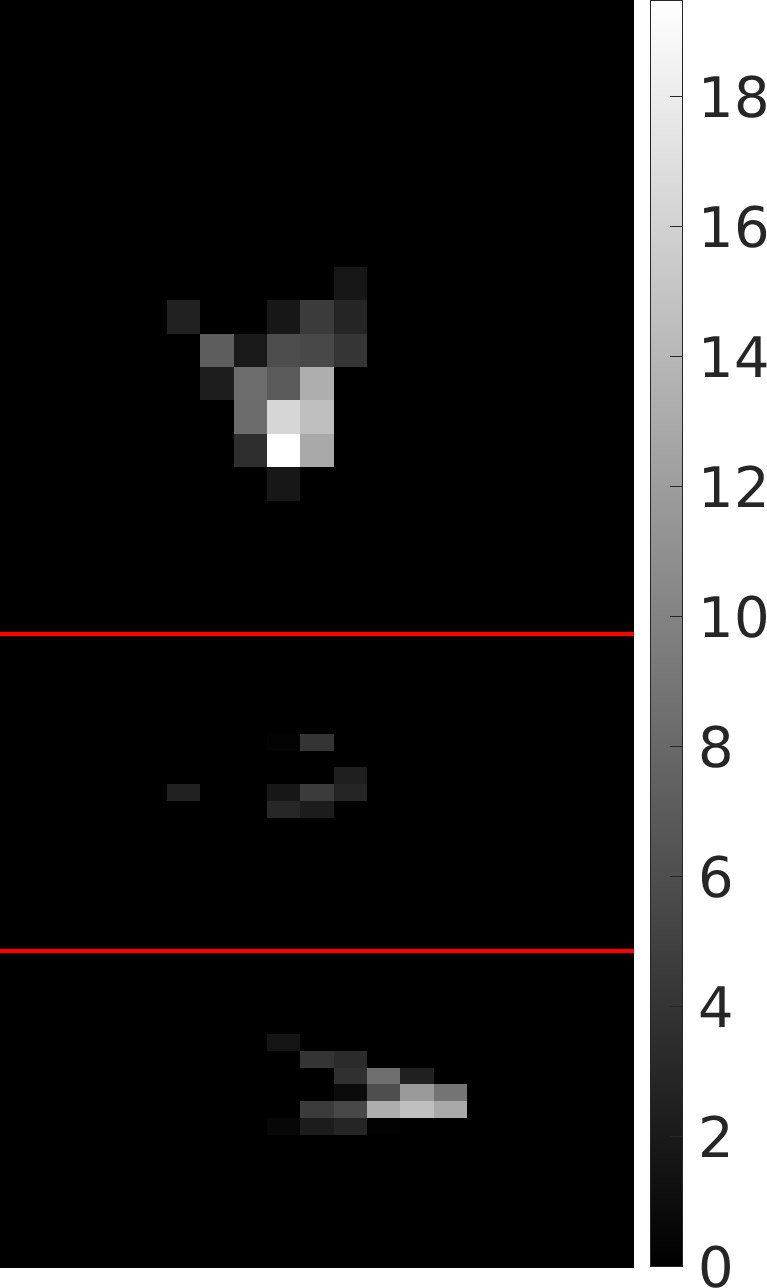}&
 \includegraphics[height=3.4cm]{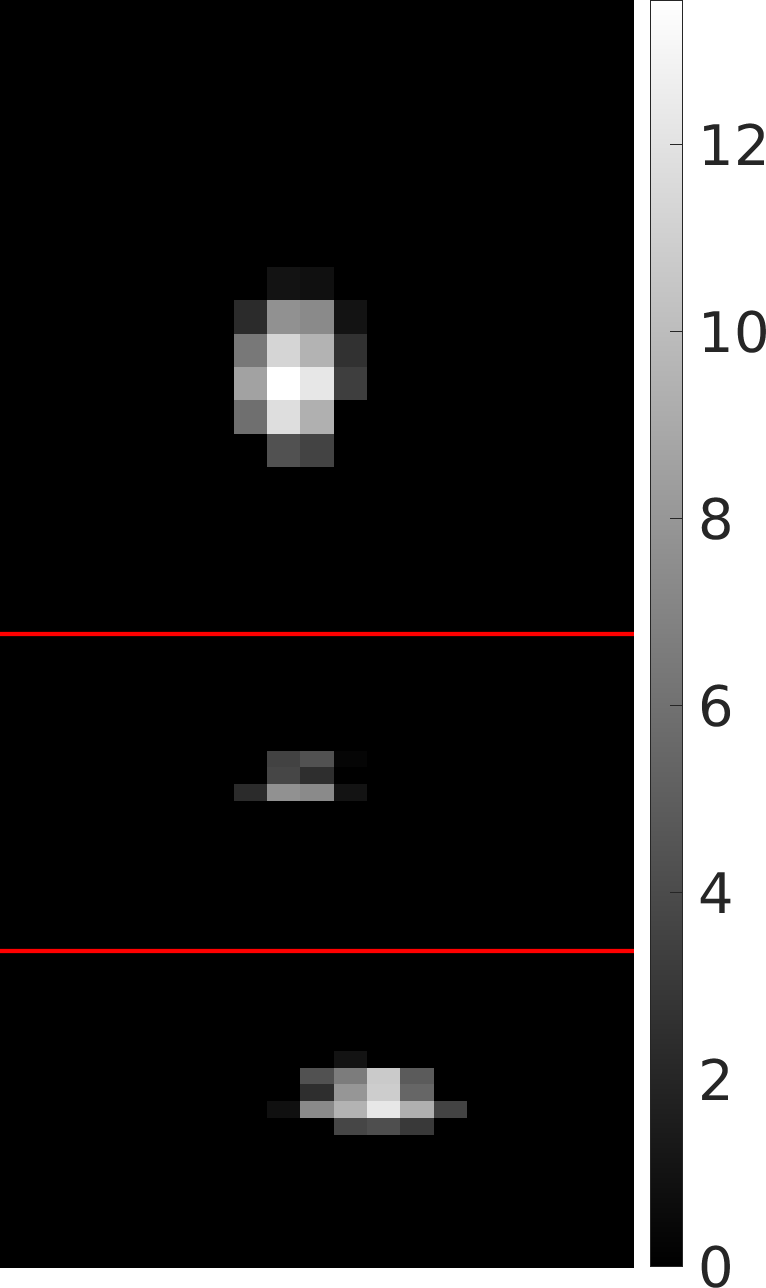}&
 \includegraphics[height=3.4cm]{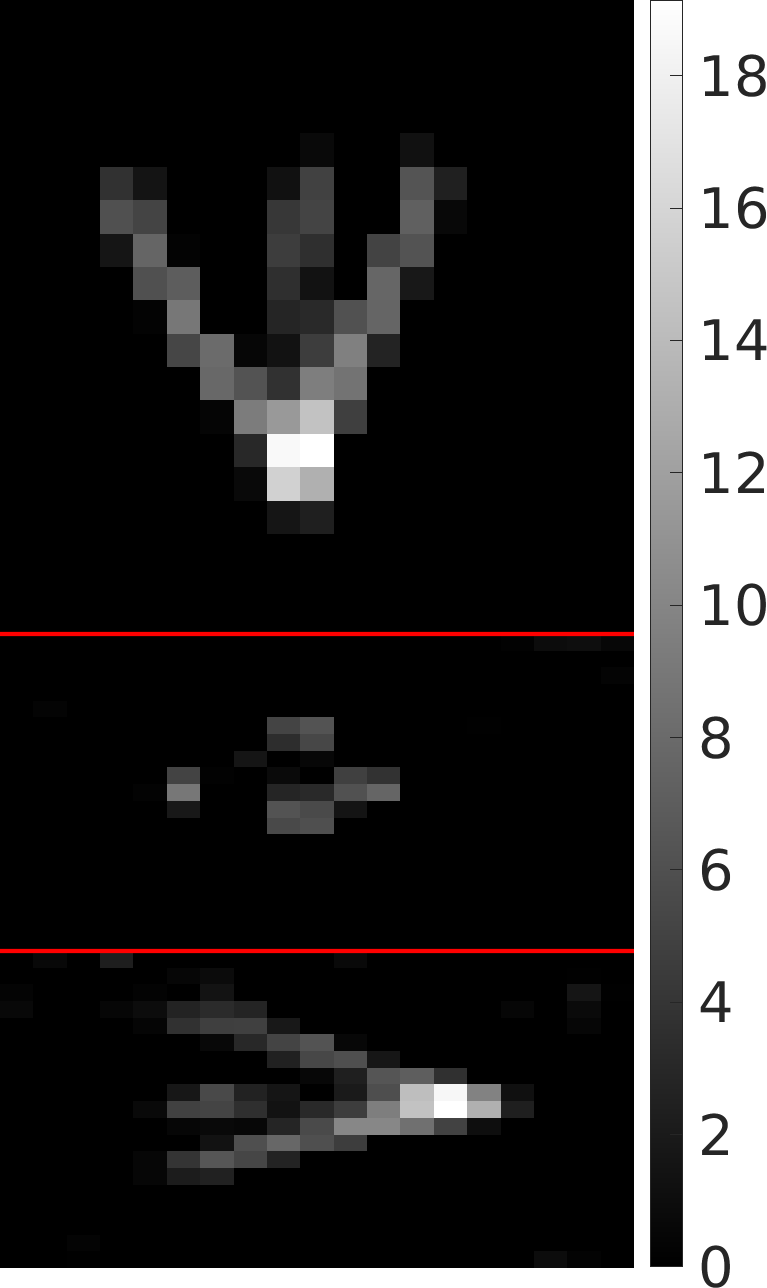}&
 \includegraphics[height=3.4cm]{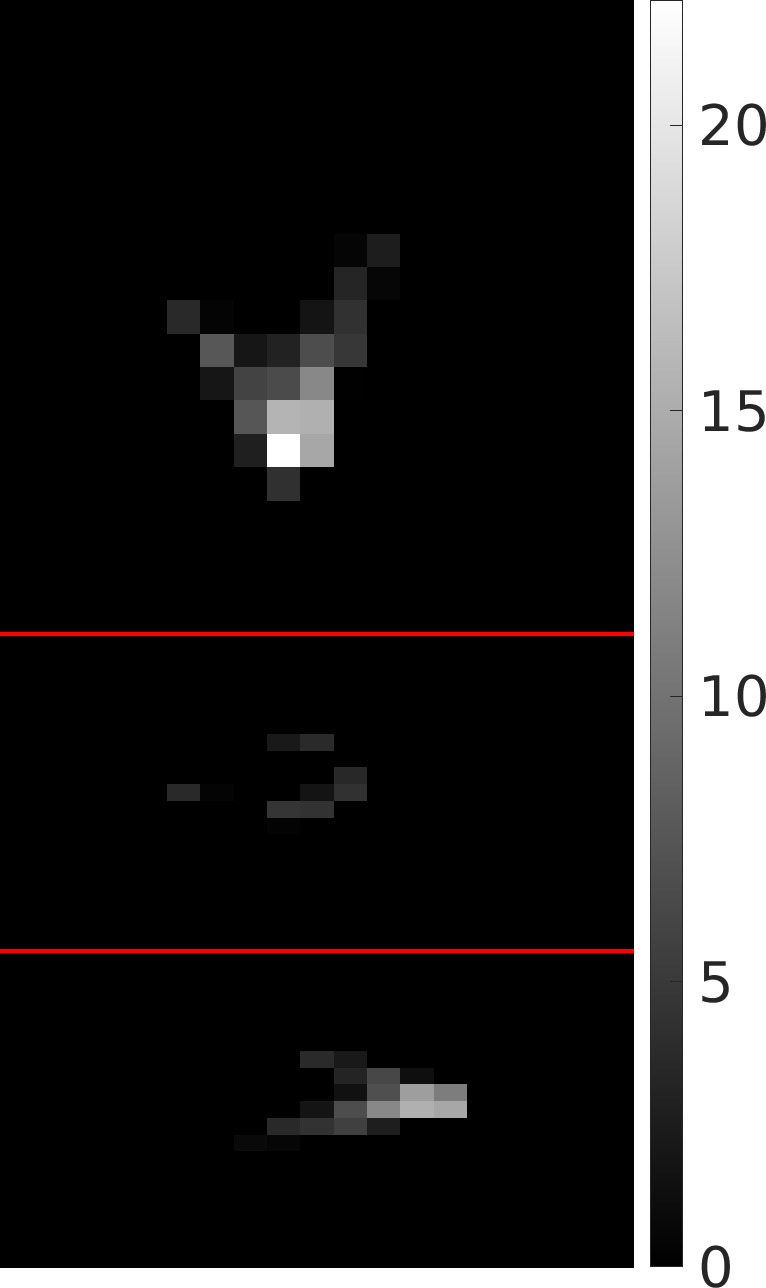}&
 \includegraphics[height=3.4cm]{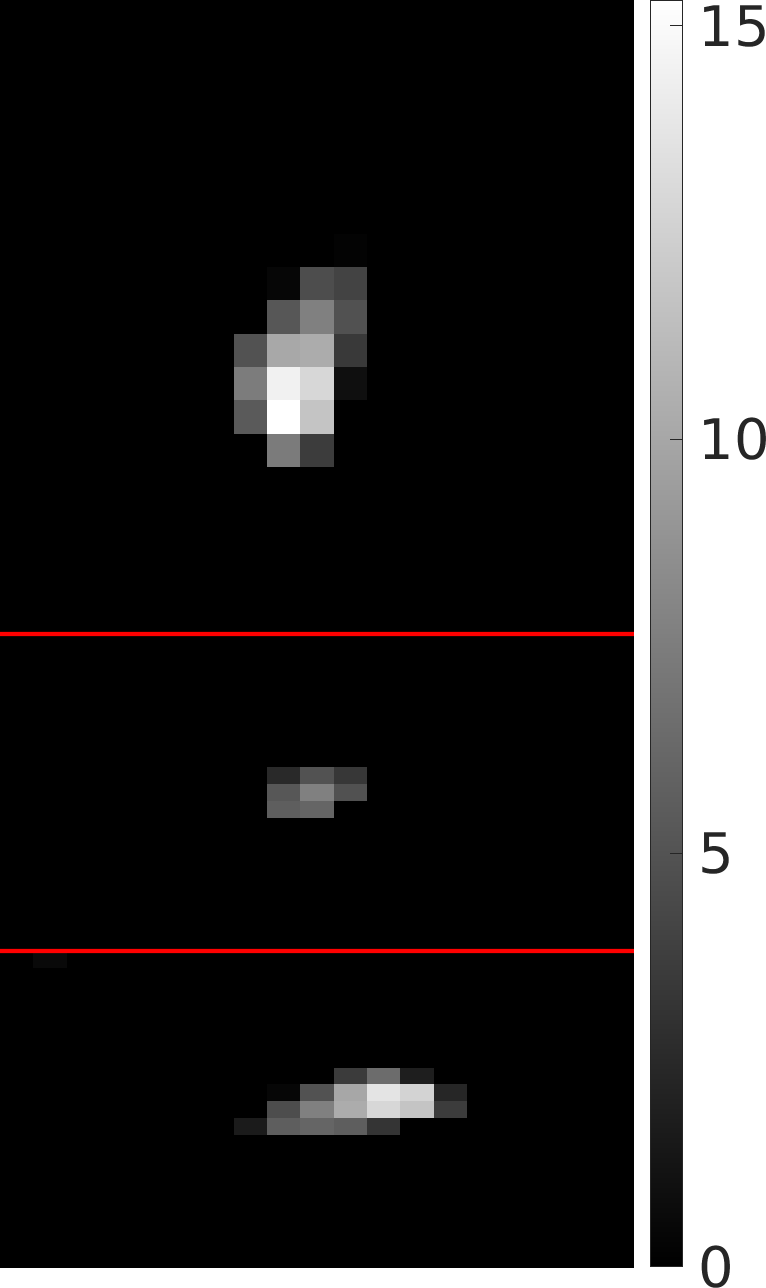}&
 \includegraphics[height=3.4cm]{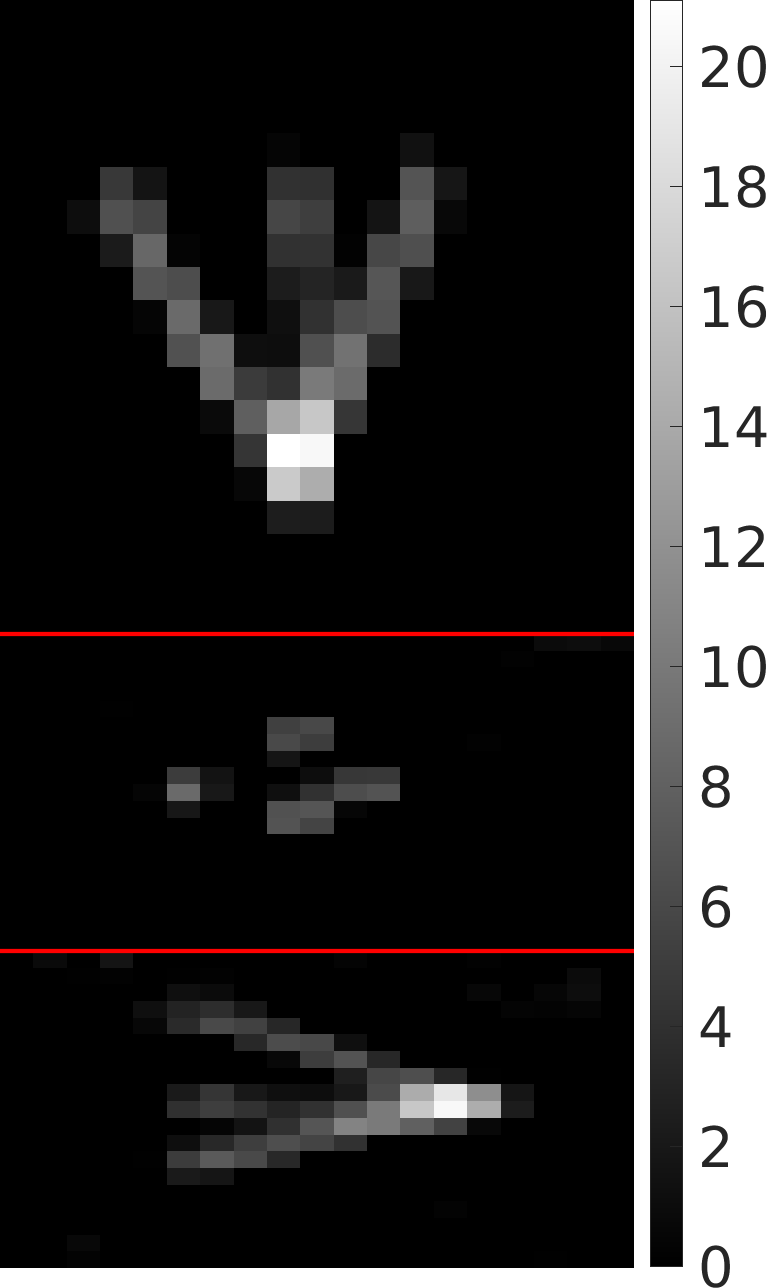}\\
\hline
\multicolumn{6}{l}{$\tau=1$} \\
 \includegraphics[height=3.4cm]{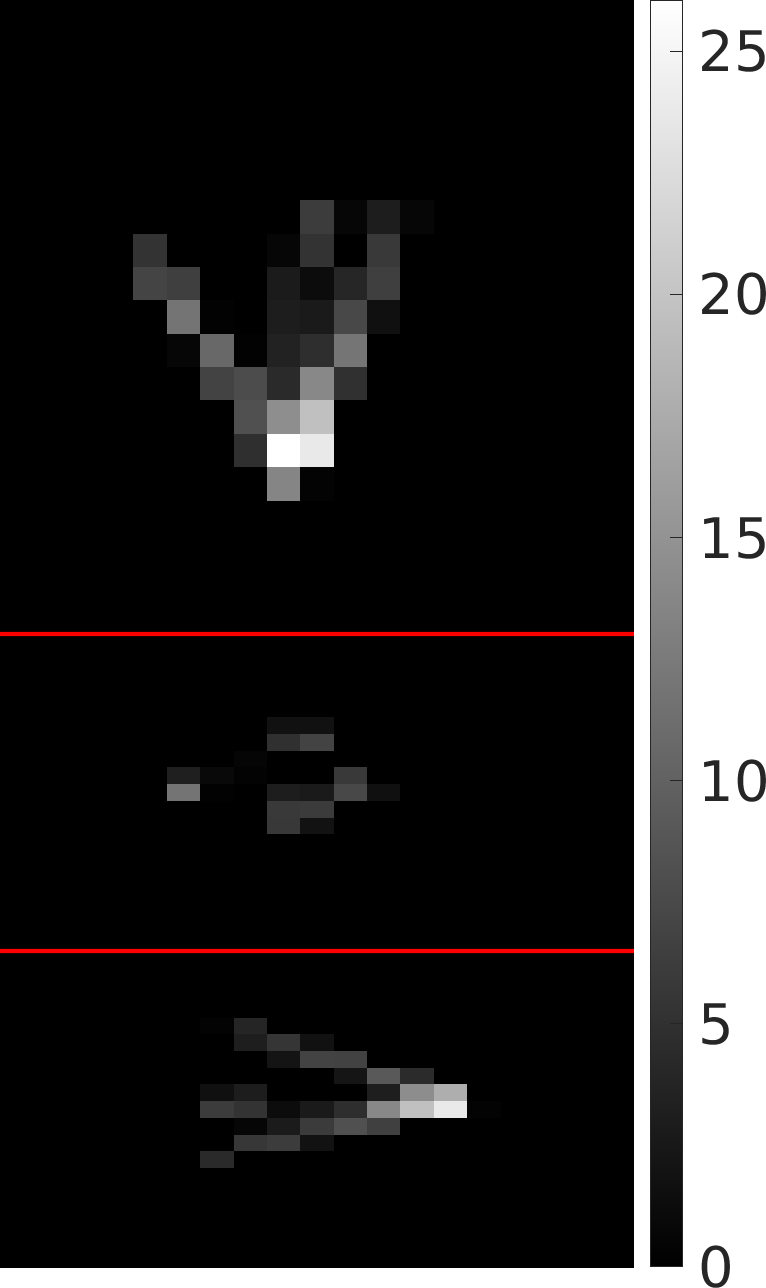}&
 \includegraphics[height=3.4cm]{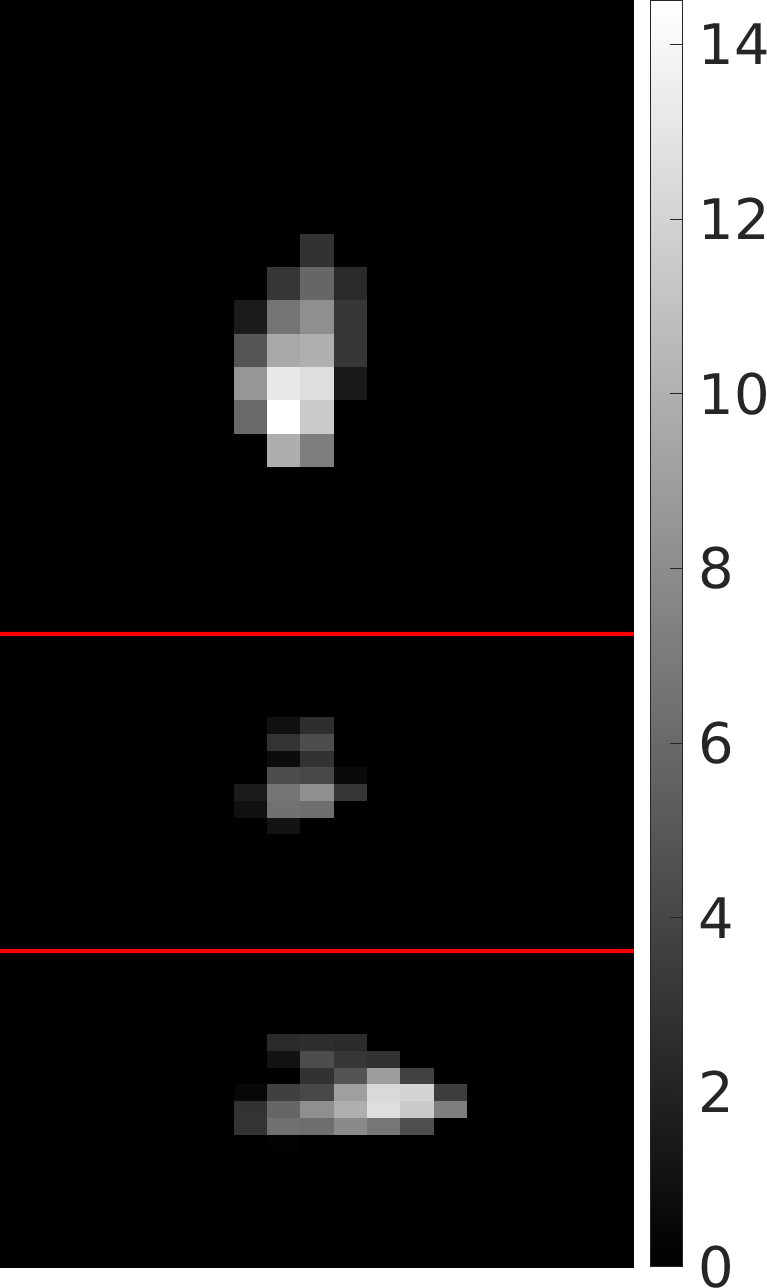}&
 \includegraphics[height=3.4cm]{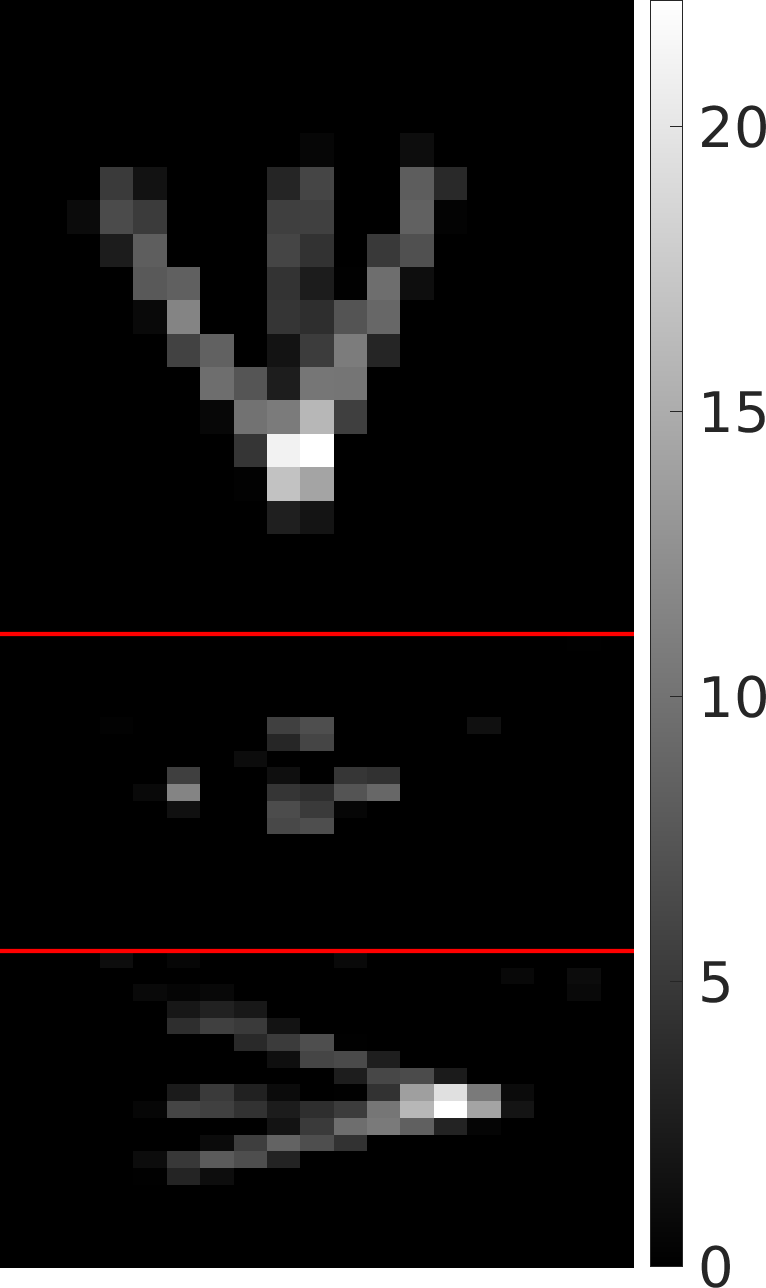}&
 \includegraphics[height=3.4cm]{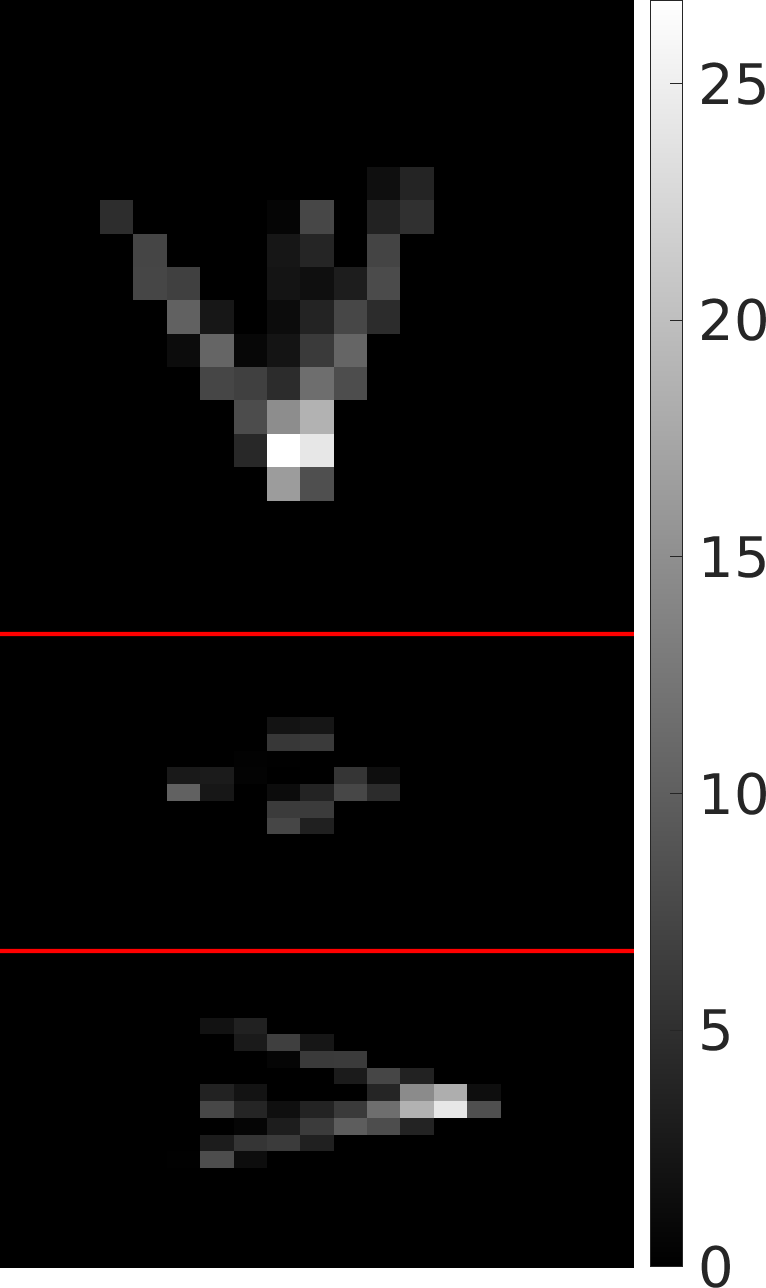}&
 \includegraphics[height=3.4cm]{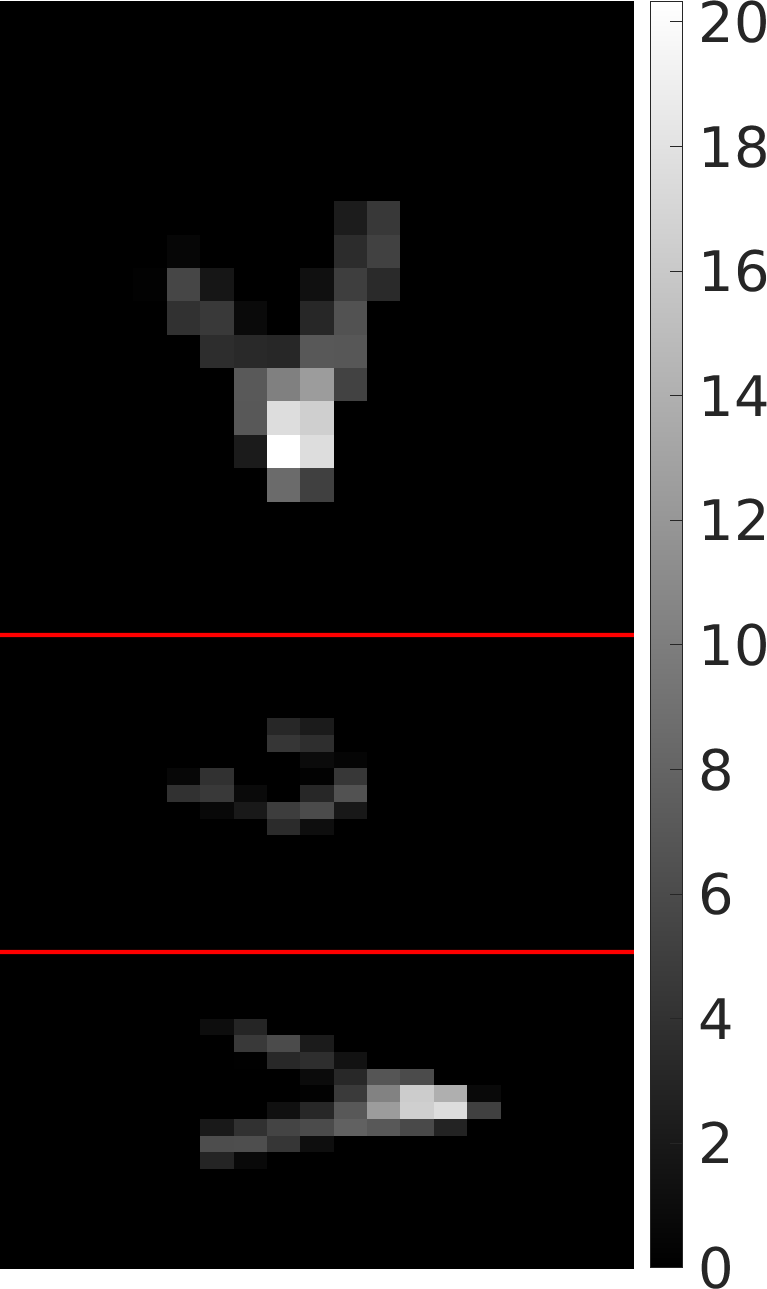}&
 \includegraphics[height=3.4cm]{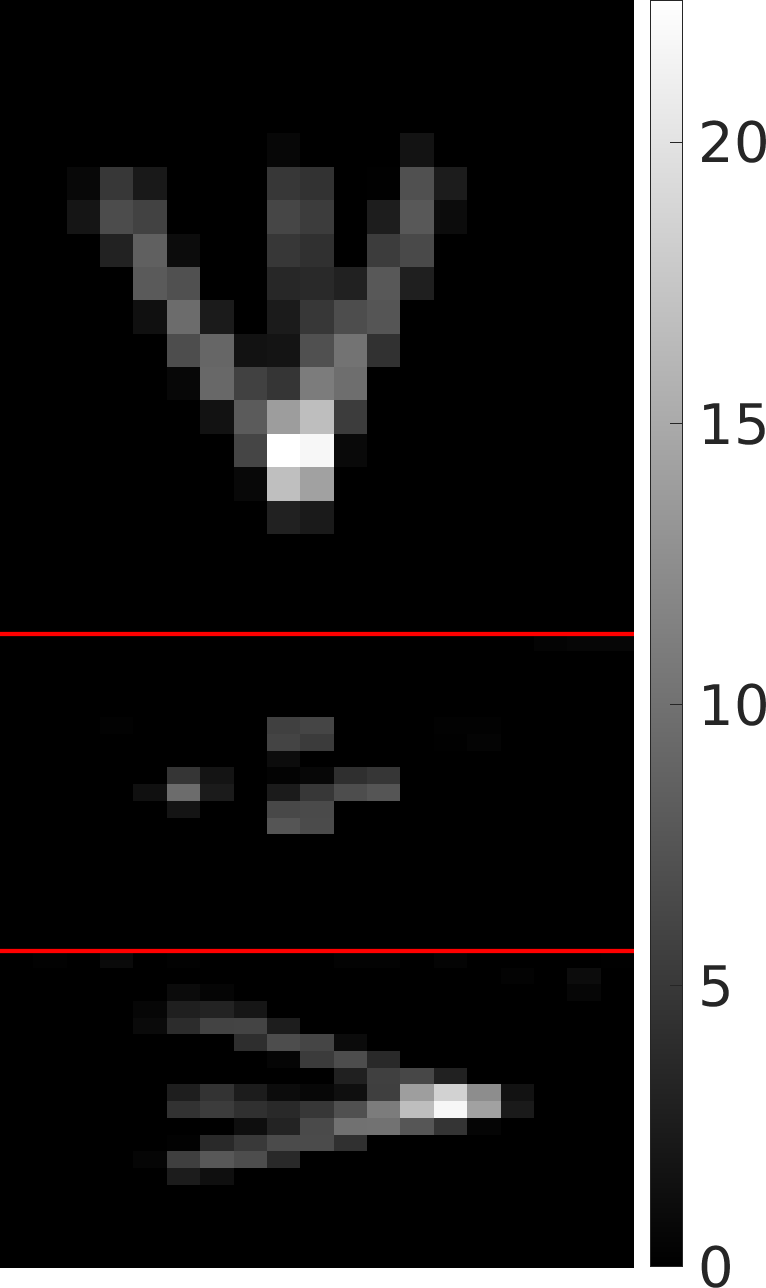}\\
\hline
\multicolumn{6}{l}{$\tau=3$} \\
 \includegraphics[height=3.4cm]{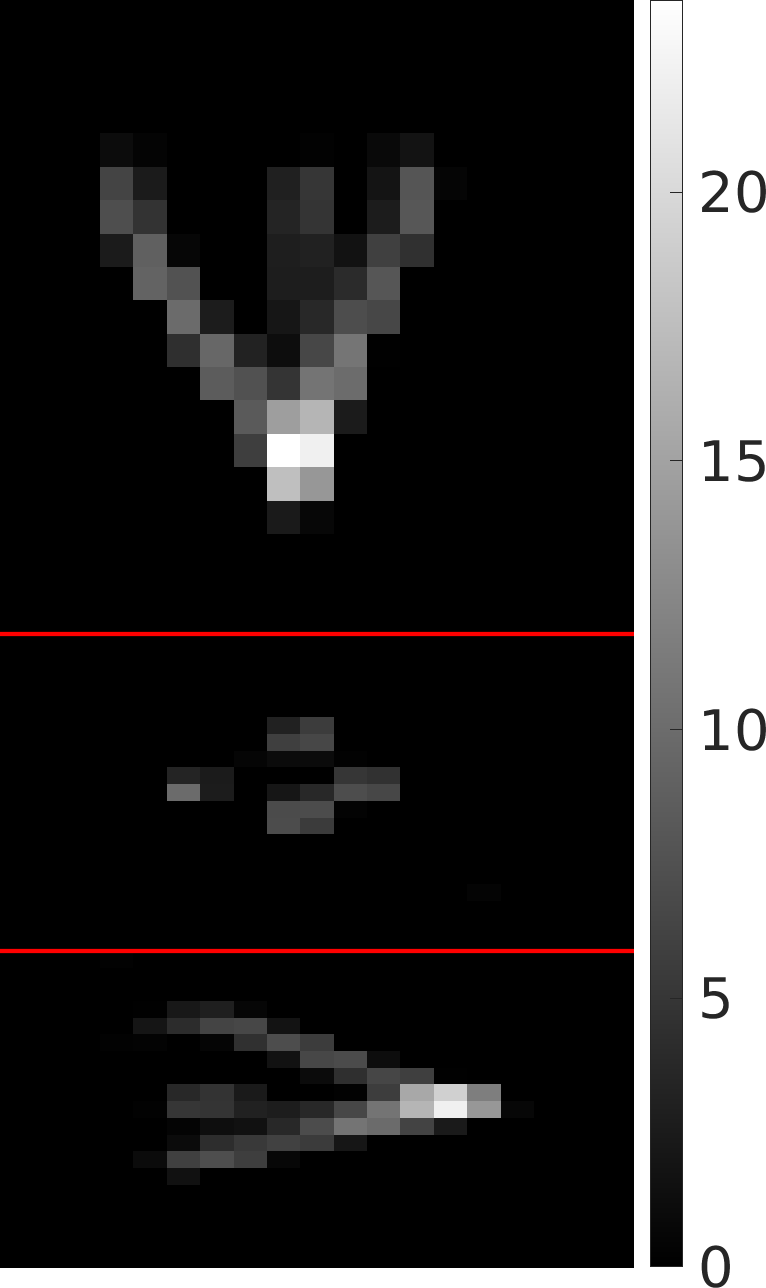}&
 \includegraphics[height=3.4cm]{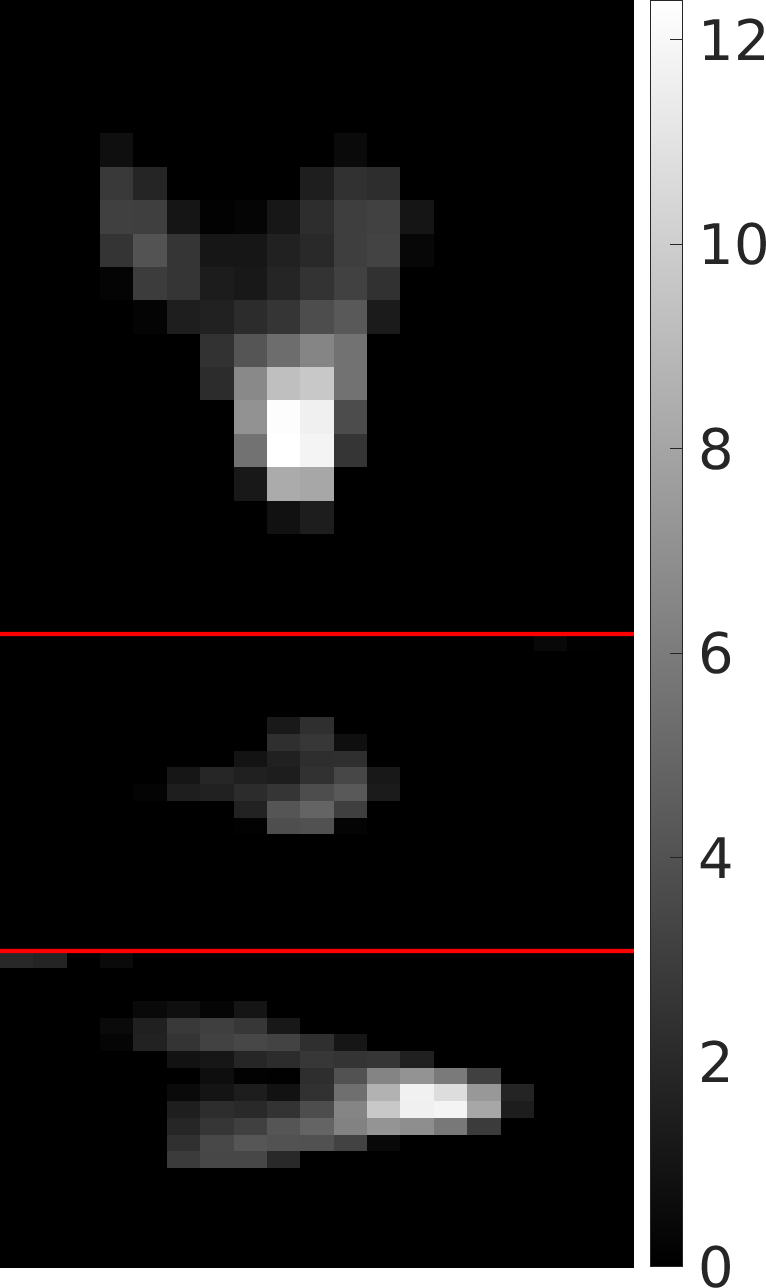}&
 \includegraphics[height=3.4cm]{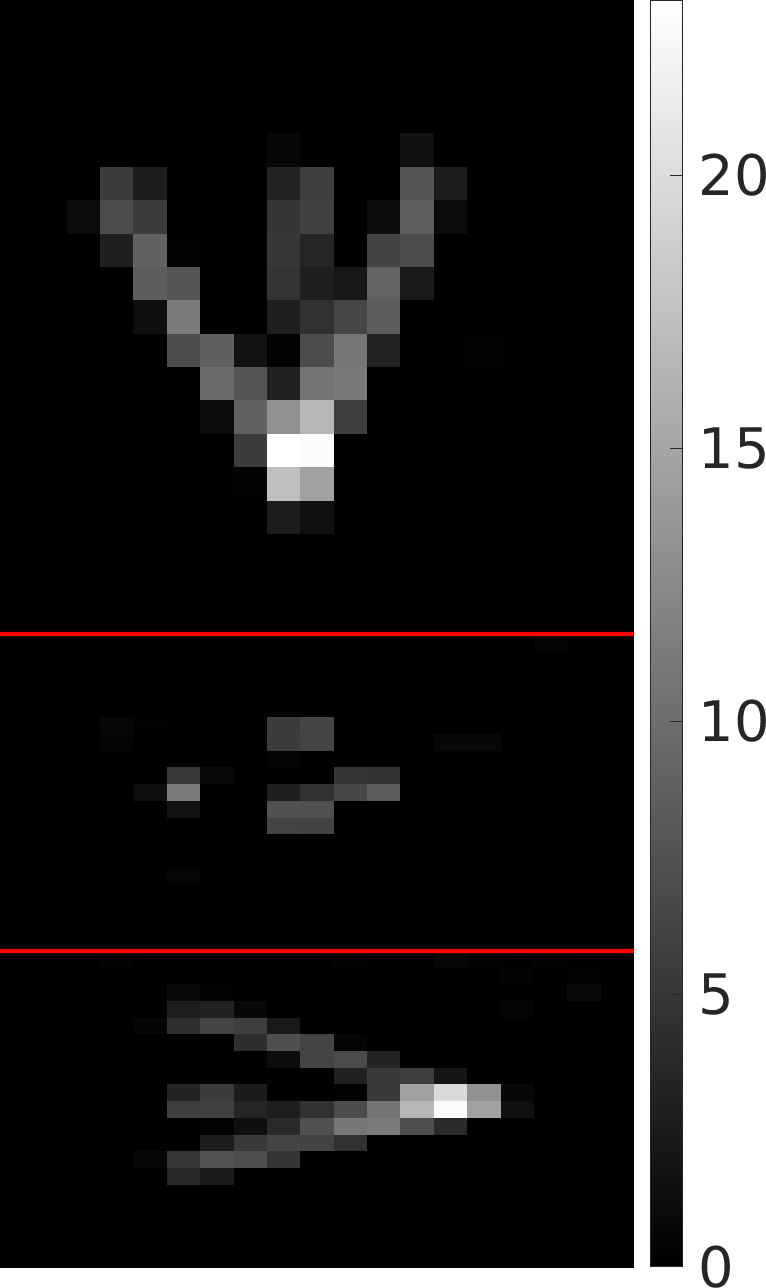}&
 \includegraphics[height=3.4cm]{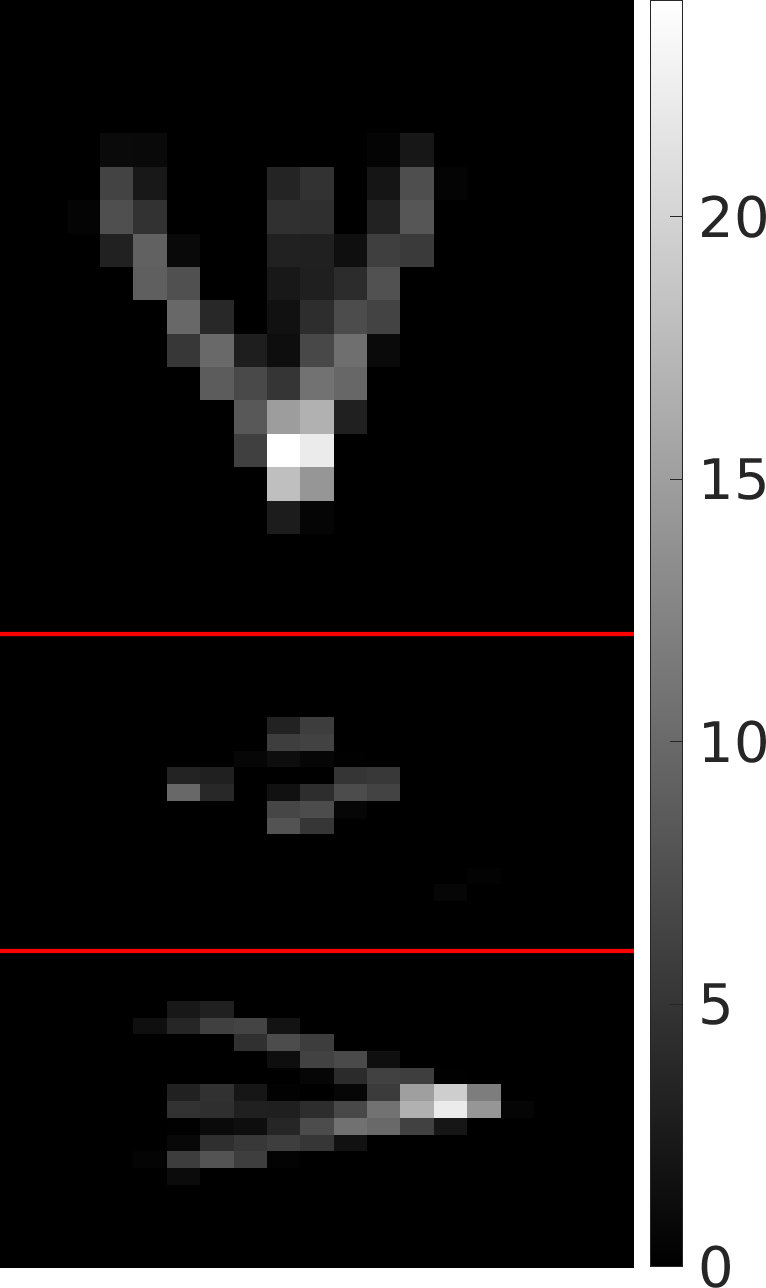}&
 \includegraphics[height=3.4cm]{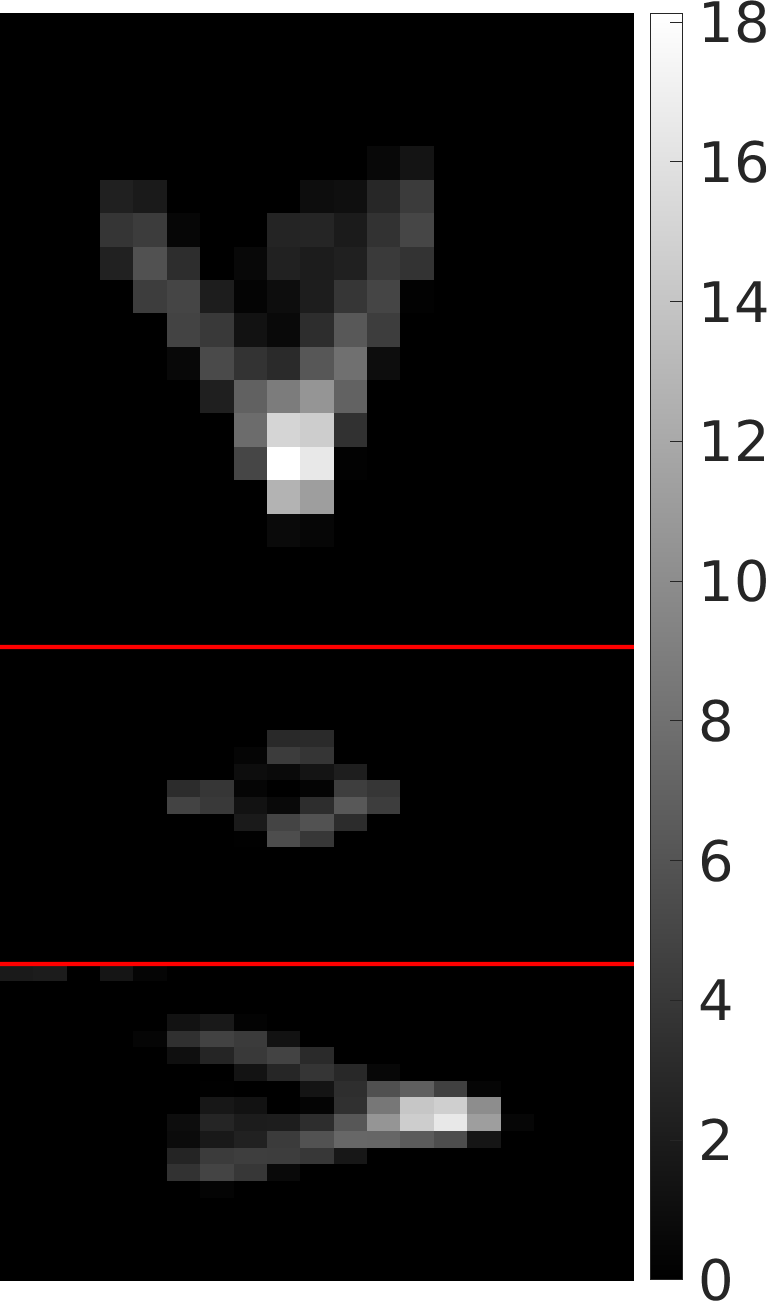}&
 \includegraphics[height=3.4cm]{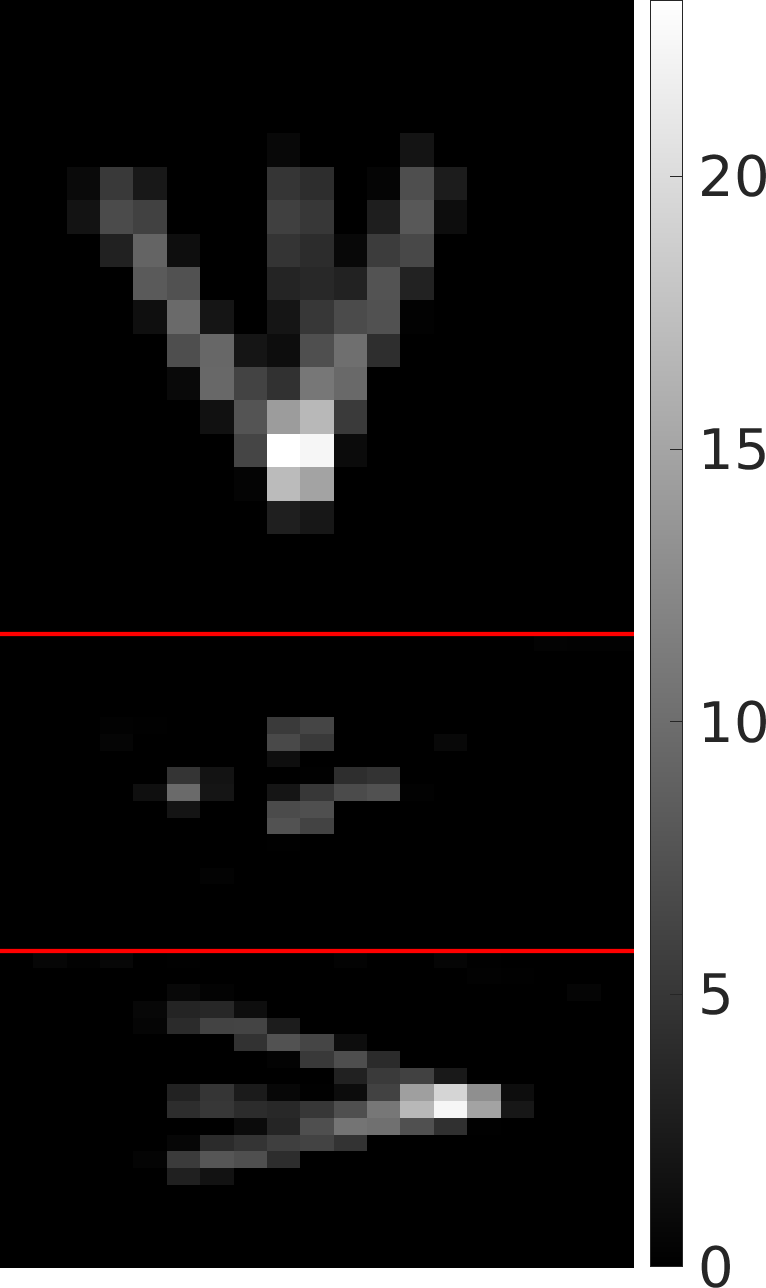}\\
 \hline
\multicolumn{6}{l}{$\tau=5$} \\
 \includegraphics[height=3.4cm]{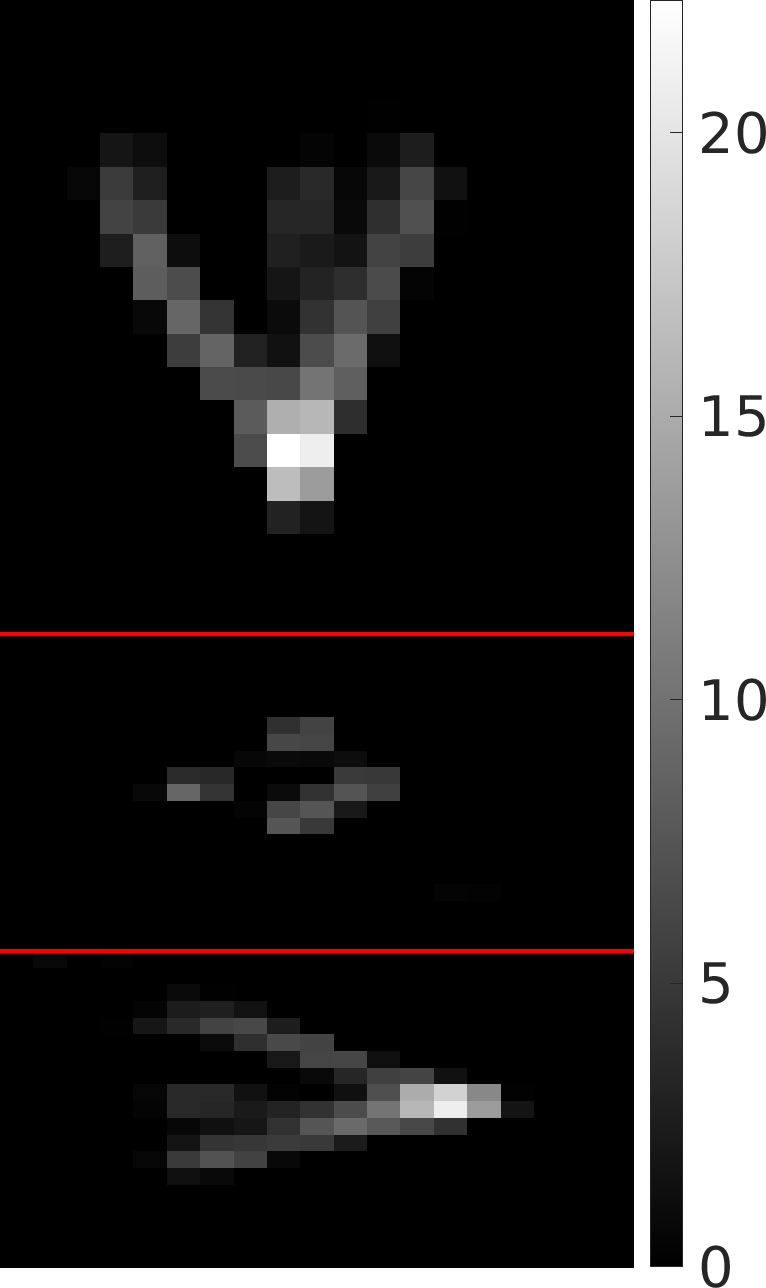}&
 \includegraphics[height=3.4cm]{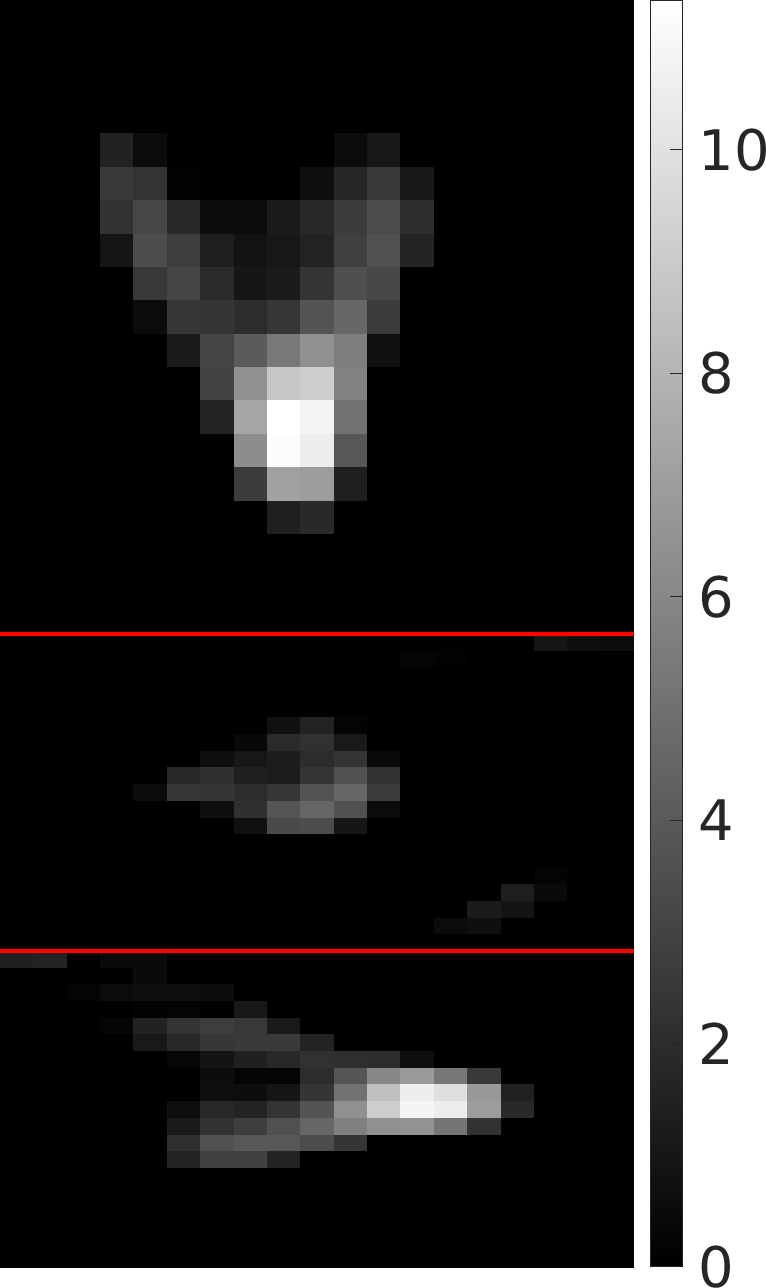}&
 \includegraphics[height=3.4cm]{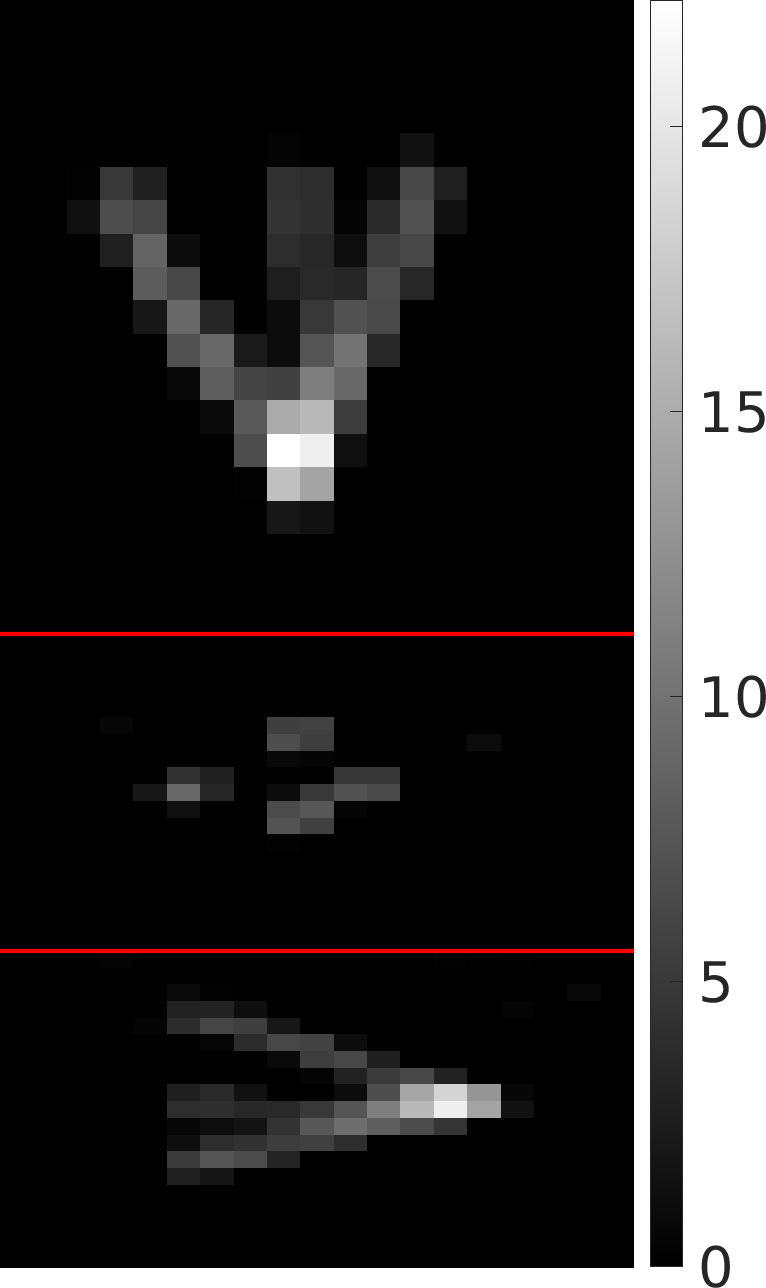}&
 \includegraphics[height=3.4cm]{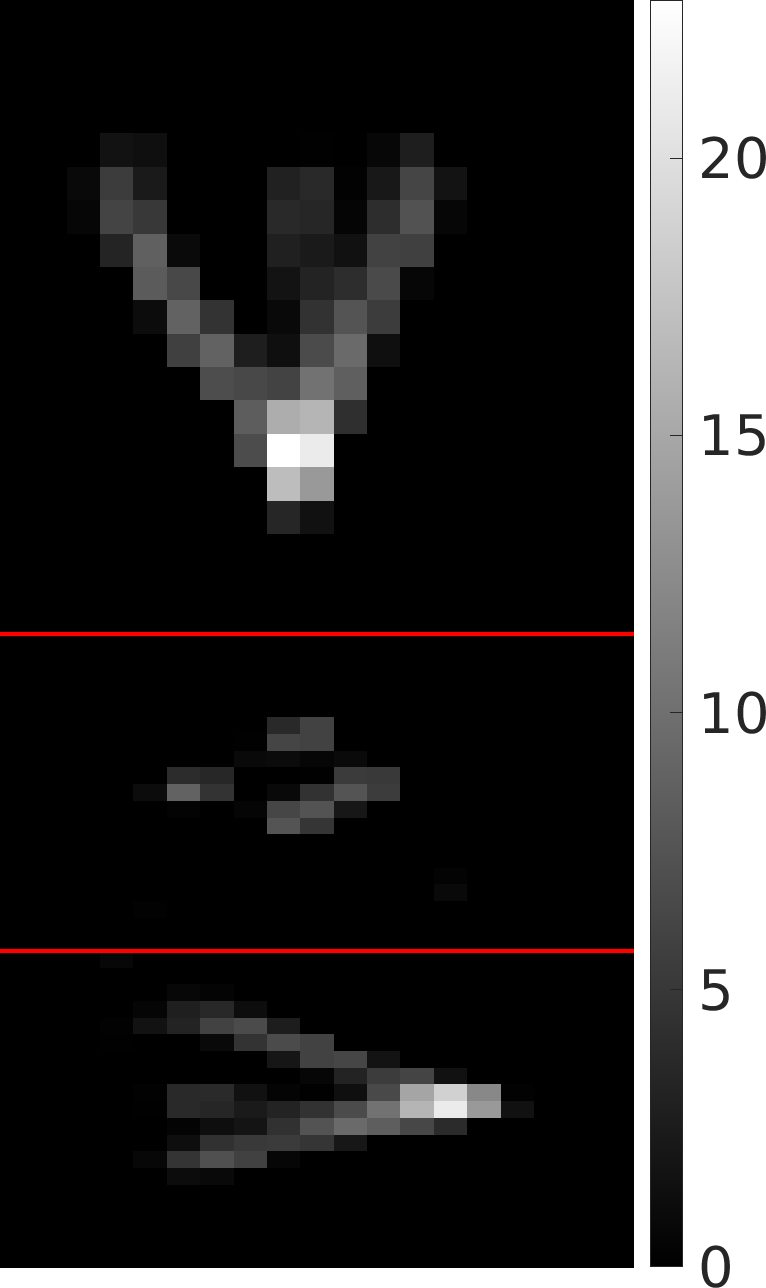}&
 \includegraphics[height=3.4cm]{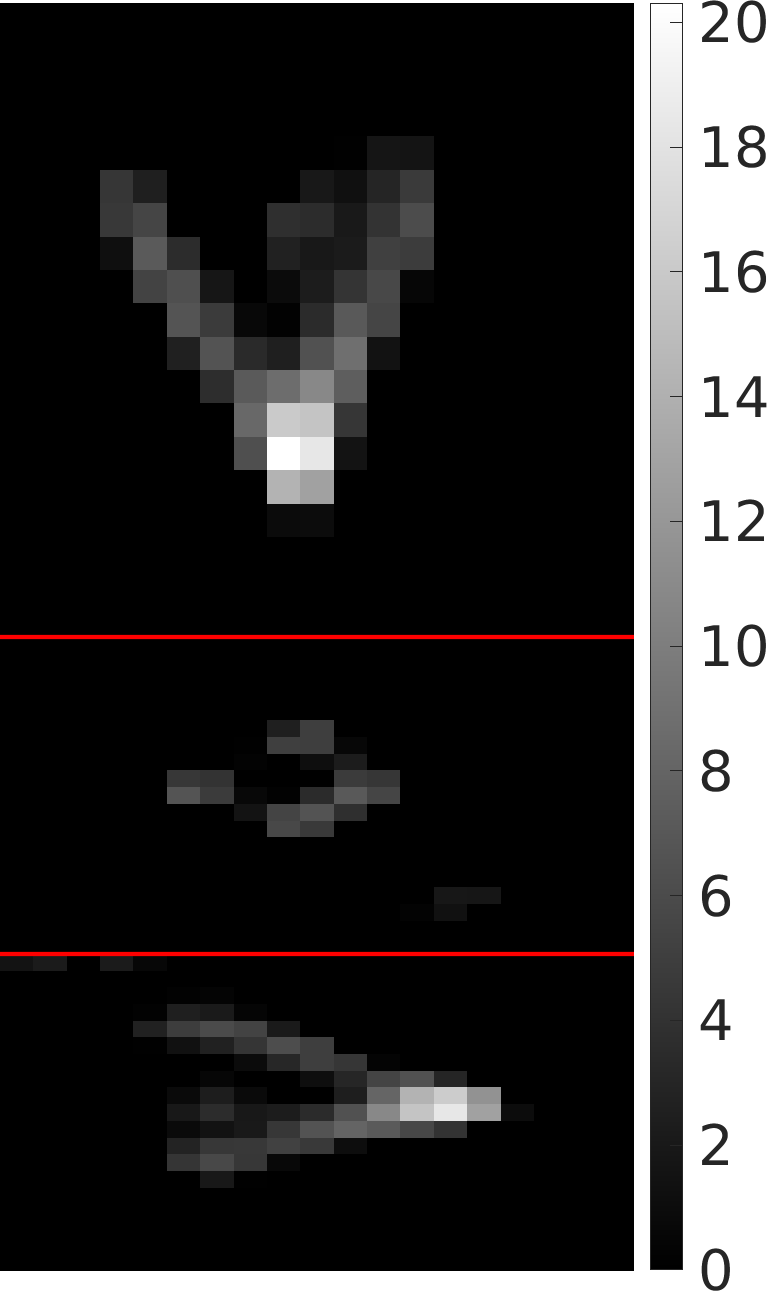}&
 \includegraphics[height=3.4cm]{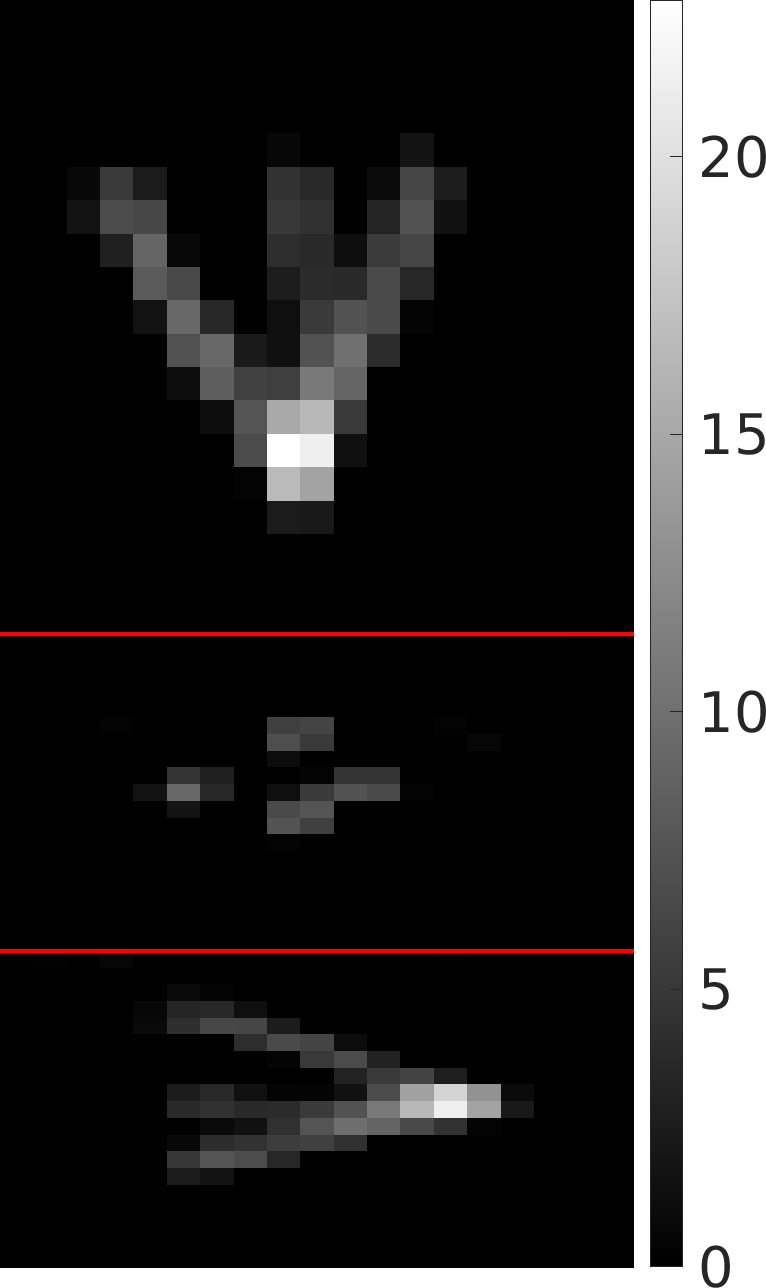}\\
\end{tabular}
}
\caption{``Resolution'' phantom reconstructions, SSIM-optimized $\alpha$ and iteration number $N$ (for l2-K only) according to Table \ref{tab:ssim_nonwhitened_vs_whitened}.}
\label{fig:methods_nonwhitened_vs_whitened_resolution_ssim}
\end{figure*}

\subsection{Algorithmic performance, quantitative comparison and the influence of SNR-type thresholding}

First, we compare the performance of the methods quantitatively and qualitatively. To this end, we employ the best image
quality measure for each method and SNR-type threshold $\tau$, using the corresponding ``optimal'' regularization parameter
$\alpha$, respectively $\alpha$ and iteration number $N$ for l2-K, and then visualize the respective reconstructions for a qualitative comparison. We analyze the
measures PSNR and SSIM separately.

\subsubsection*{PSNR}
The PSNR results are given in Table \ref{tab:psnr_nonwhitened_vs_whitened};
see Figs. \ref{fig:methods_nonwhitened_vs_whitened_shape_psnr} and \ref{fig:methods_nonwhitened_vs_whitened_resolution_psnr} for the reconstructions.
First, we compare the two pure variational regularization approaches, i.e., l1-L and l2-L. Clearly, l1-L greatly outperforms
l2-L in all cases, showing robustness of the l1 fitting with respect to outliers. Second, we compare l1-L to l2-K.
The comparison with l2-K is more difficult due to the presence of two different sources of regularizing effect, i.e., variational
and iterative, controlled respectively by the regularization parameter $\alpha$ (in all methods) and iteration number $N$. Using a fixed and small $N$,
as often done in practice, l2-K does not reach convergence in the sense of optimization (i.e., finding a global minimizer
to \eqref{eqn:standard}), and instead is actually early stopping in the spirit of
iterative regularization \cite{EnglHankeNeubauer:1996}. Nonetheless, for a sufficiently large $N$, the result
by l2-K is similar to that by l2-L. The behavior of l2-K with regard to the early stopping issue will be examined more closely below.
Here, we choose the optimum with respect to the tuple $(\alpha,N)$.
Now we examine the PSNR results more closely.
For both phantoms l2-K yields the best reconstruction for all SNR thresholds.
The overall optimum is found in the whitened case for $\tau=1$ for both phantoms.
In particular, the ``resolution'' phantom results obtained with l1-L yields comparable PSNR values when using larger $\tau$.

Qualitatively, the reconstructions in Figs. \ref{fig:methods_nonwhitened_vs_whitened_shape_psnr} and
\ref{fig:methods_nonwhitened_vs_whitened_resolution_psnr} (corresponding to the PSNR-optimal $\alpha$, resp. $N$)
exhibit severe background artifacts for l2-K; The artifacts are more visible in the inverted map shown in
Appendix \ref{app:supplements_inverted_colormap}. In particular, l2-L fails to give reasonable results,
and l1-L gives far more reasonable reconstructions for $\tau\geq 1$. For both phantoms, l1-L and l2-K give similar
results for $\tau \geq 3$, and both outperform l2-L.

\subsubsection*{SSIM}

The SSIM results are given in Table \ref{tab:ssim_nonwhitened_vs_whitened} and the reconstructions
in Figs. \ref{fig:methods_nonwhitened_vs_whitened_shape_ssim} and \ref{fig:methods_nonwhitened_vs_whitened_resolution_ssim}.
Similar to the PSNR case, l1-L and l2-K outperform l2-L in all cases. In the non-whitened case, l1-L performs best for both phantoms phantom
 for large $\tau\geq 3$. For both phantoms, the overall best possible SSIM is obtained with l1-L
in the non-whitened. In the whitened case, l1-L and l2-K yield the best possible SSIM.

In contrast to the PSNR results, the reconstructions in Figs. \ref{fig:methods_nonwhitened_vs_whitened_shape_ssim} and
\ref{fig:methods_nonwhitened_vs_whitened_resolution_ssim} indicate less severe but still pronounced background artifacts for l2-K.
For small thresholds, l2-L fails to give reasonable results, while l1-L can give good reconstructions for $\tau \geq 1$.
For the ``resolution'' phantom, l1-L and l2-K give similar results for $\tau \geq 3$, and both perform superior to l2-L.

\subsubsection*{Influence of SNR-threshold and general observations}

Now we study the behavior of the quality measures with respect to the SNR threshold $\tau$ by examining
Tables \ref{tab:psnr_nonwhitened_vs_whitened} and \ref{tab:ssim_nonwhitened_vs_whitened} columnwise more closely. For PSNR,
the observations vary across the methods. For l2-L, PSNR increases monotonically with $\tau$ in all cases,
for l2-K, it does not show a steady trend, and for l1-L, it reaches a maximum before decreases again for $\tau=5$. In contrast,
SSIM in Table \ref{tab:ssim_nonwhitened_vs_whitened} first increases and then decreases for l1-L and l2-L in most cases, while l2-K does not show a steady trend like in the PSNR case.
Thus, frequency selection with a proper $\tau$ benefits all variational methods and is recommended for MPI reconstruction, but
a too large $\tau$ may compromise imaging quality, as observed earlier \cite{KluthJin:2019}.

With whitening, for both phantoms, the performance of l2-L and l2-K can be improved in almost all cases, but l1-L benefits less from whitening
in most cases. Statistically, whitening ensures that the i.i.d. assumption in the least-squares formulation is more
adequately fulfilled (if the variance estimate is accurate), and thus it is more beneficial to the standard approach
\eqref{eqn:standard}. l1-L is more resilient to noise type, and thus whitening plays a less important role.

The worse performance of l1-L for $\tau = 0$ indicate that solely l1 fitting without SNR thresholding is not able to compensate all data outliers.
Visually inspecting Fig. \ref{fig:cov_diag} allows identifying severe outliers when no SNR-type thresholding
is applied. Thresholding still results in a background signal with a variance structure and large outliers, but
the maximum variance is several orders smaller in magnitudes than the no-thresholding
case, cf. Fig. \ref{fig:cov_diag_tau_xcoil}. In sum, frequency selection is beneficial for all variational methods when done
carefully: A minimum $\tau$ is necessary but a too large $\tau$ can compromise image quality.

\subsection{The influence of multiple regularization techniques in the standard approach}\label{ssec:Kaczmarz}

\begin{figure*}[hbt!]%
\centering
\scalebox{1}{
\begin{tabular}{c|c}
PSNR & SSIM \\
 \hline
\multicolumn{2}{l}{$\tau=0$} \\
 \includegraphics[width=0.4\textwidth]{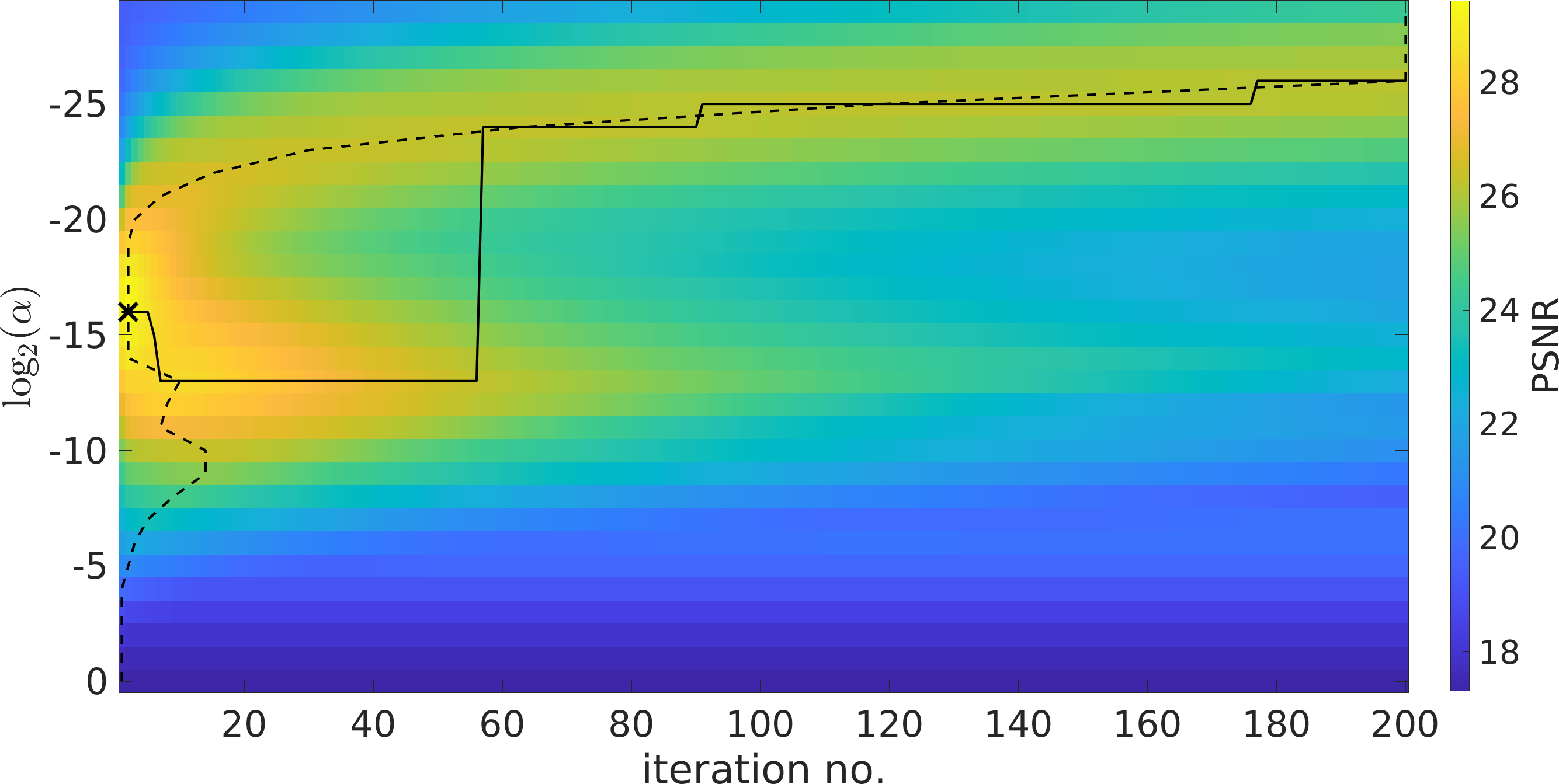} &
 \includegraphics[width=0.4\textwidth]{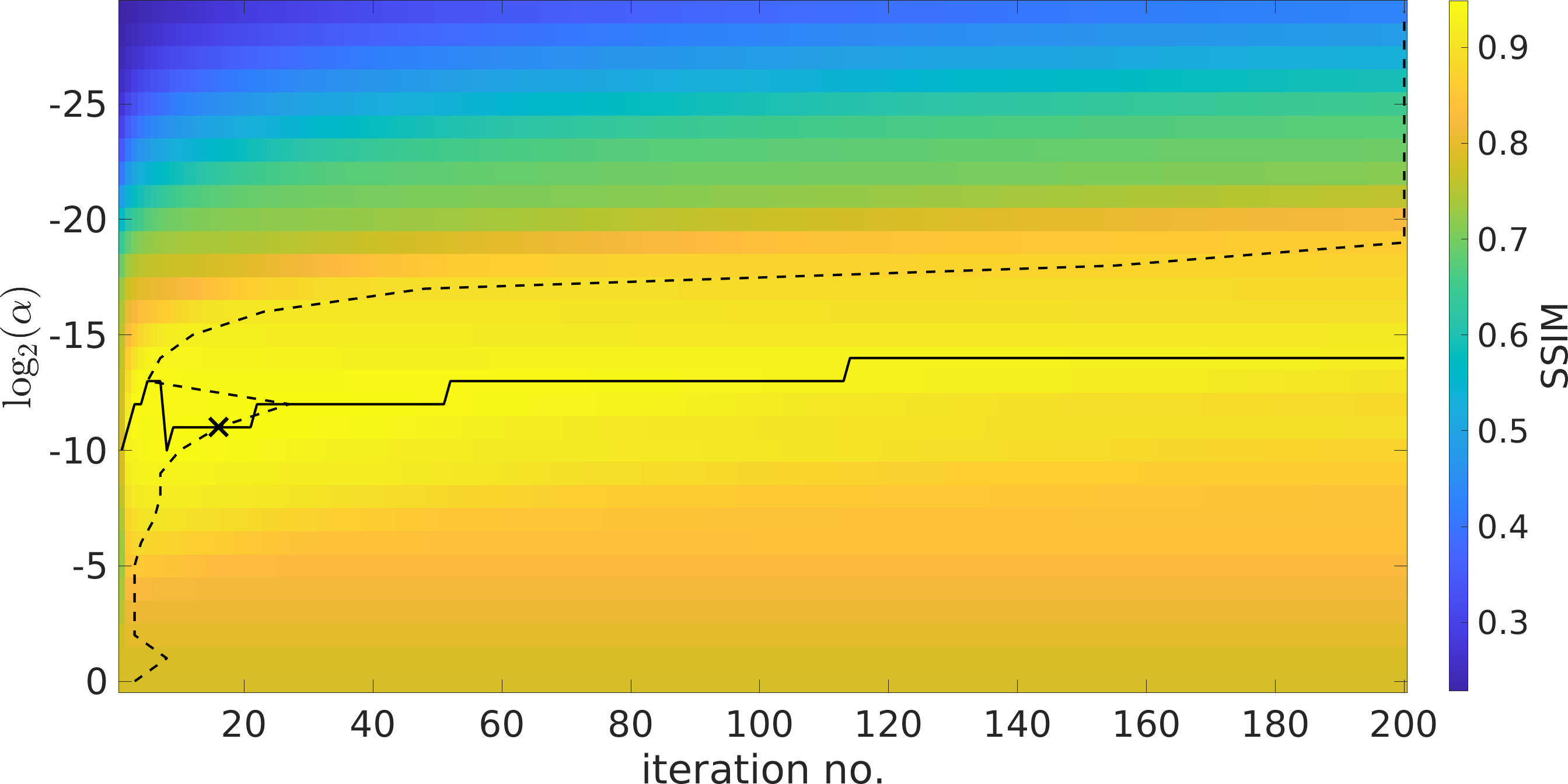} \\
 \hline
\multicolumn{2}{l}{$\tau=1$} \\
 \includegraphics[width=0.4\textwidth]{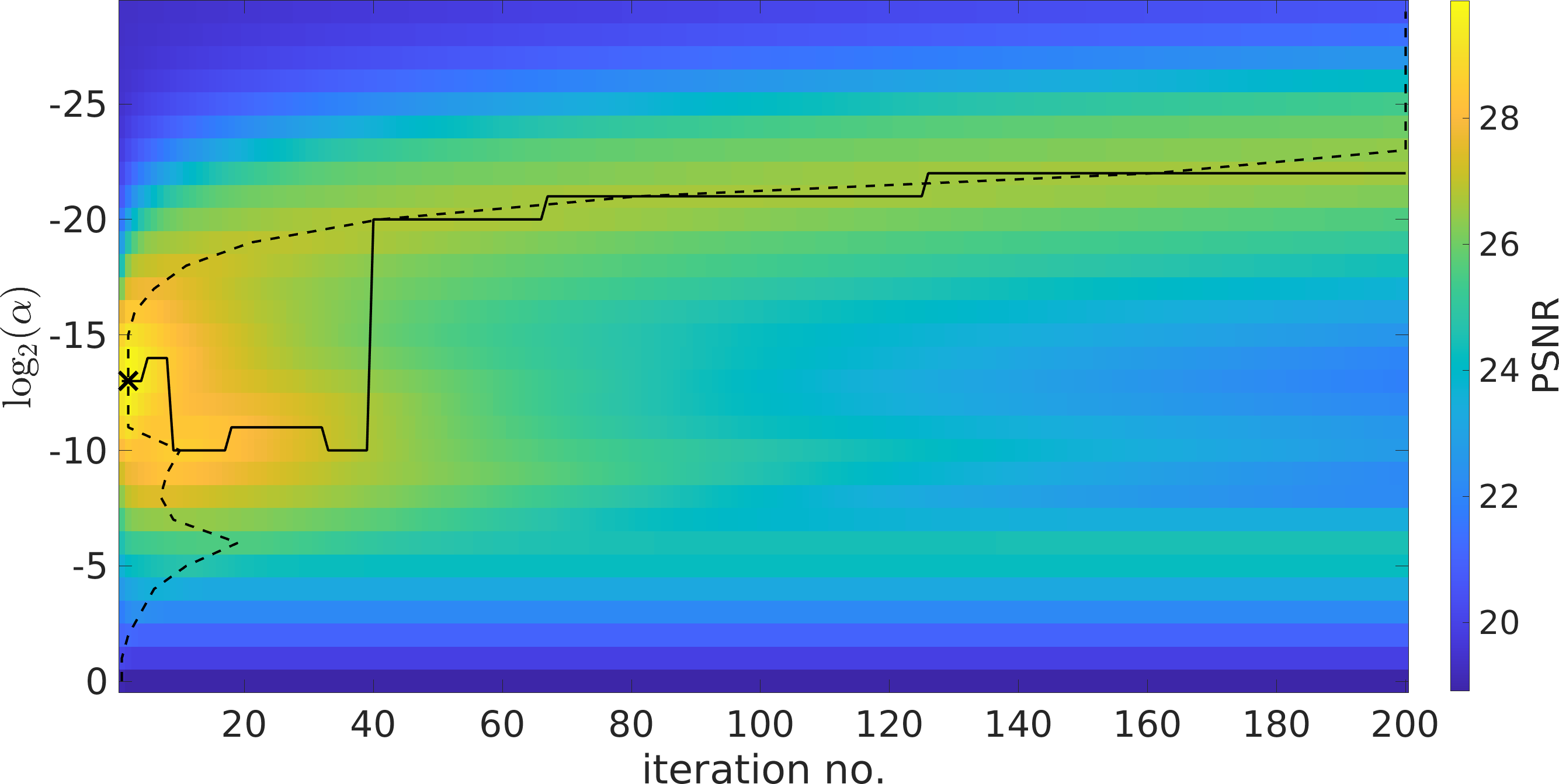} &
 \includegraphics[width=0.4\textwidth]{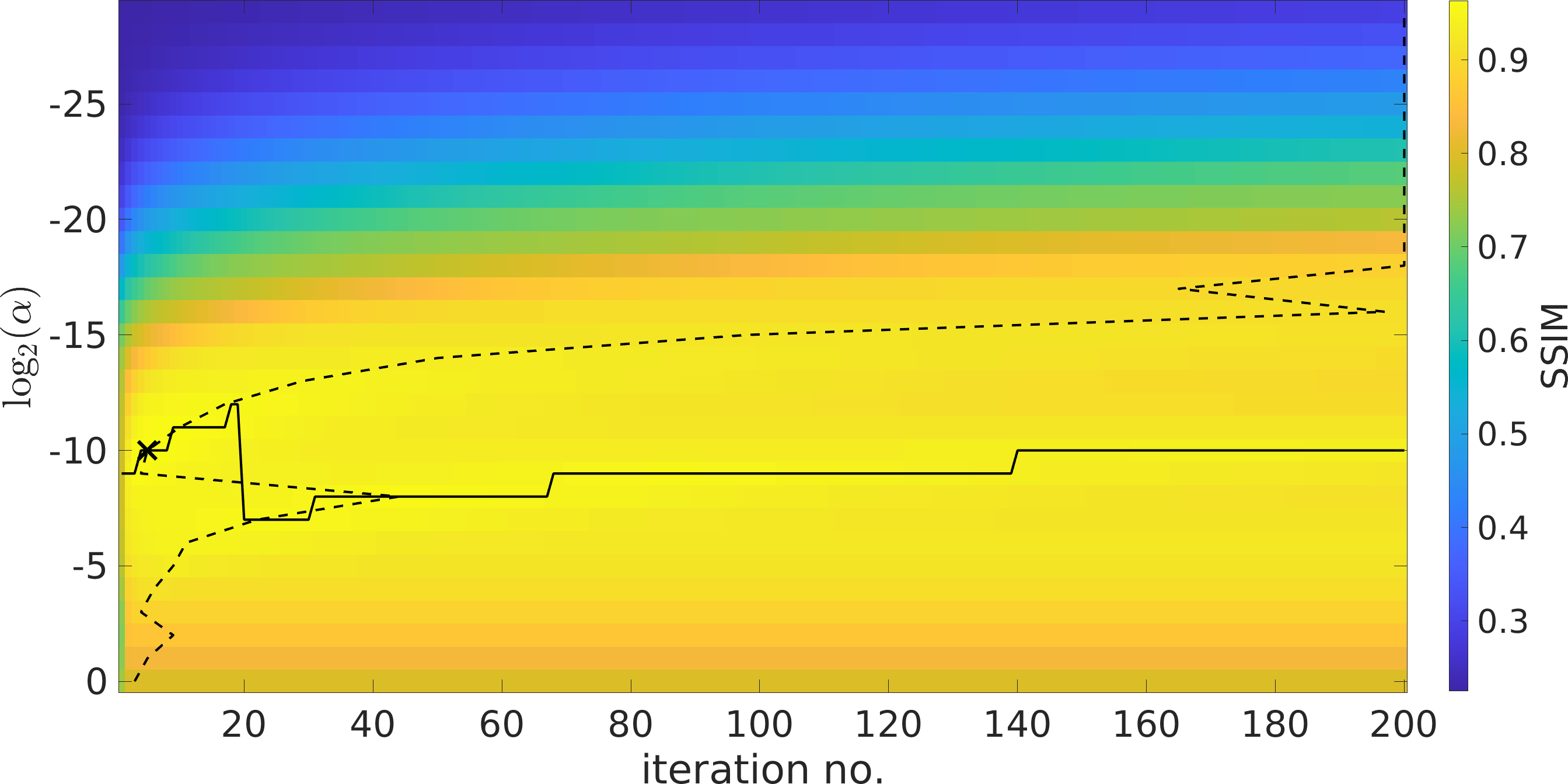} \\
  \hline
\multicolumn{2}{l}{$\tau=3$} \\
  \includegraphics[width=0.4\textwidth]{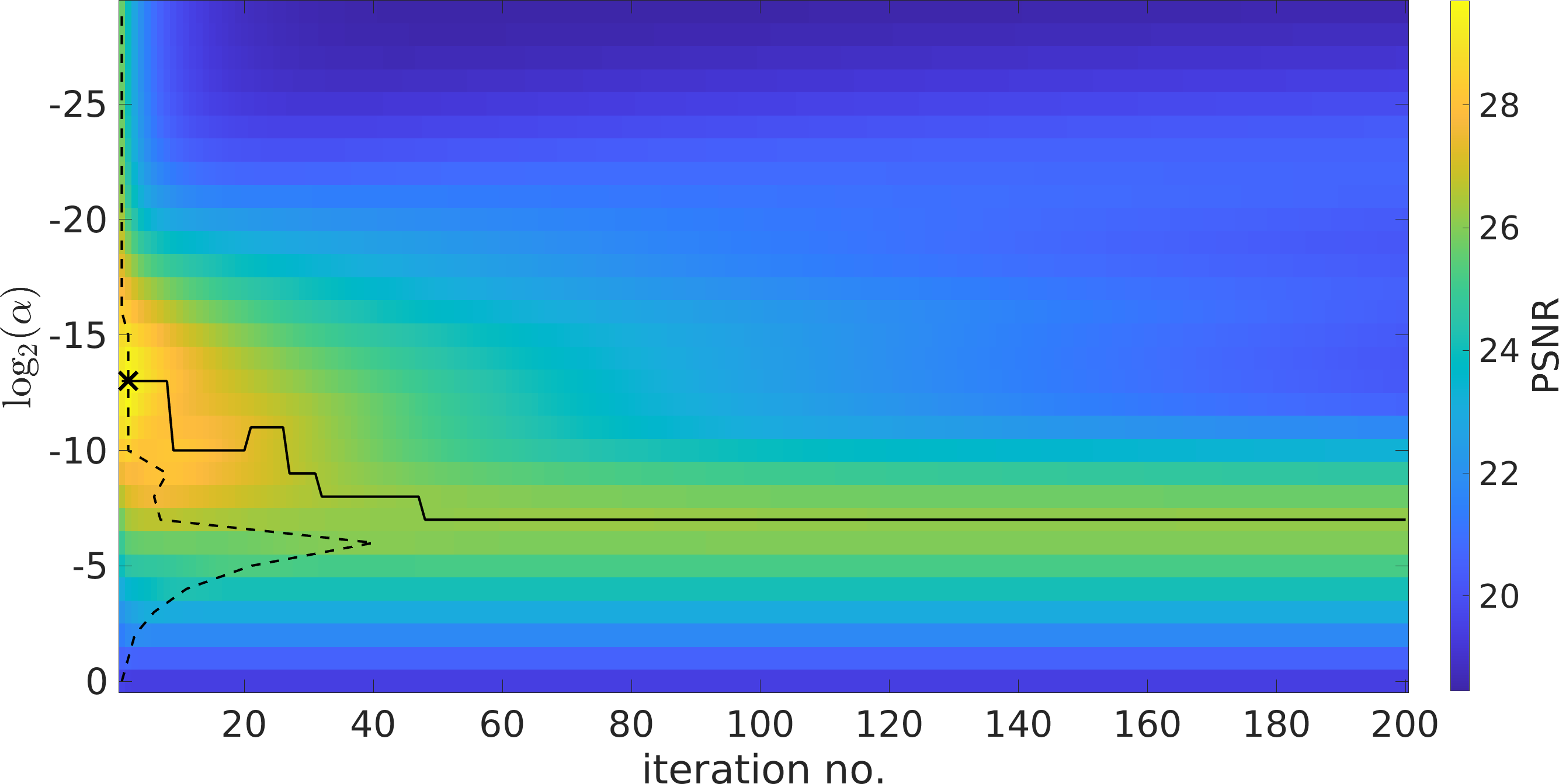} &
 \includegraphics[width=0.4\textwidth]{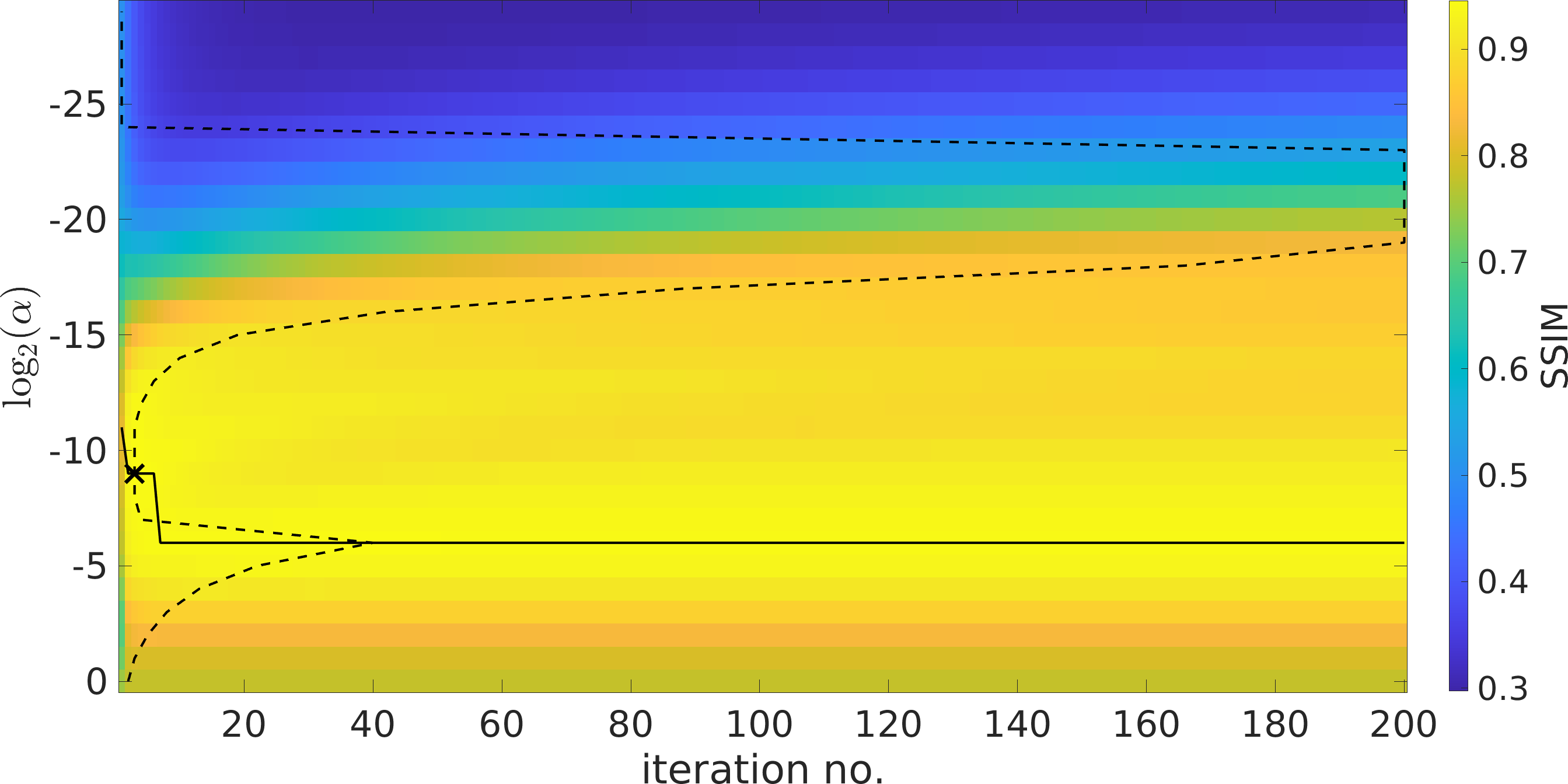} \\
  \hline
\multicolumn{2}{l}{$\tau=5$} \\
 \includegraphics[width=0.4\textwidth]{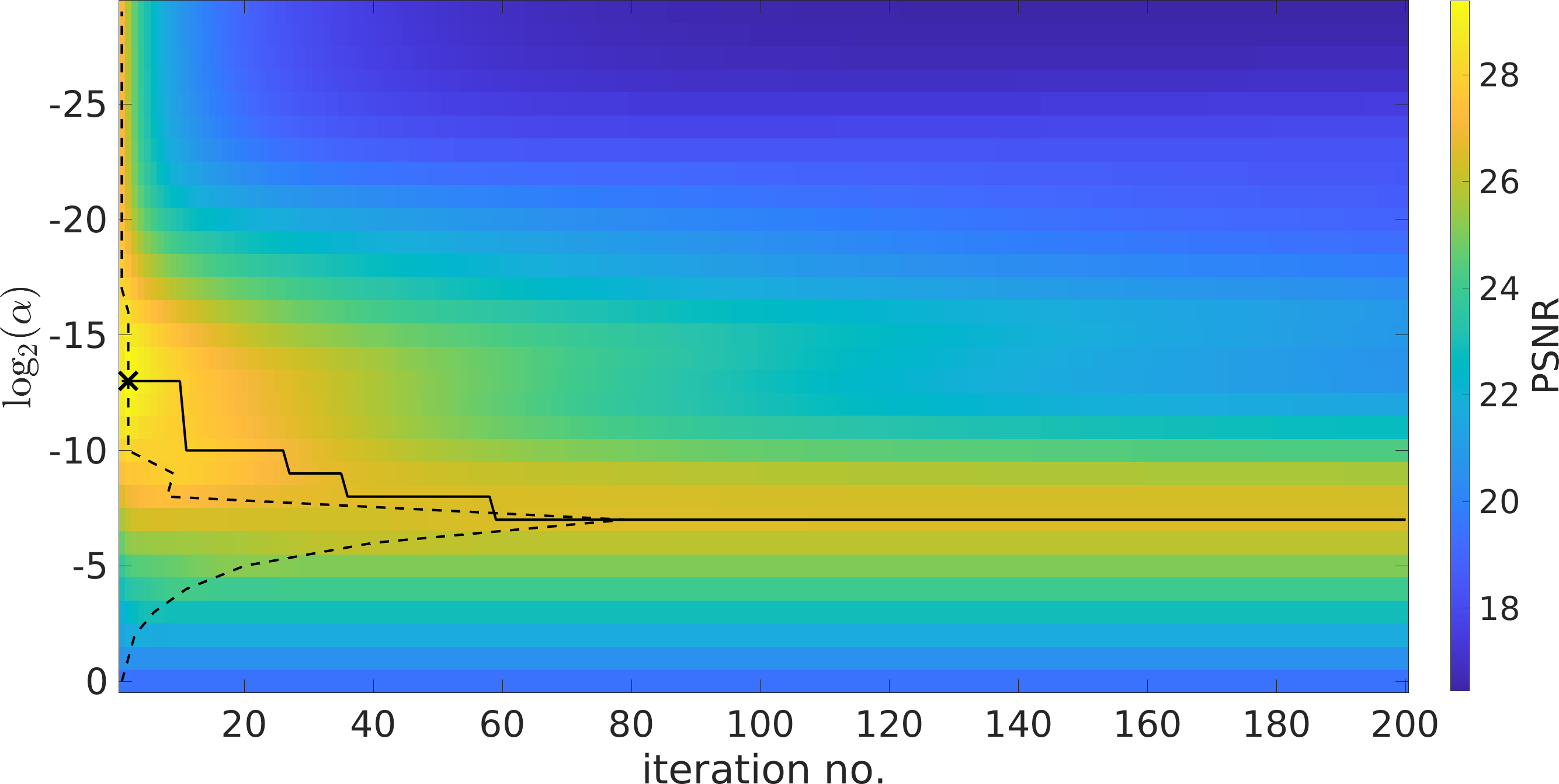} &
 \includegraphics[width=0.4\textwidth]{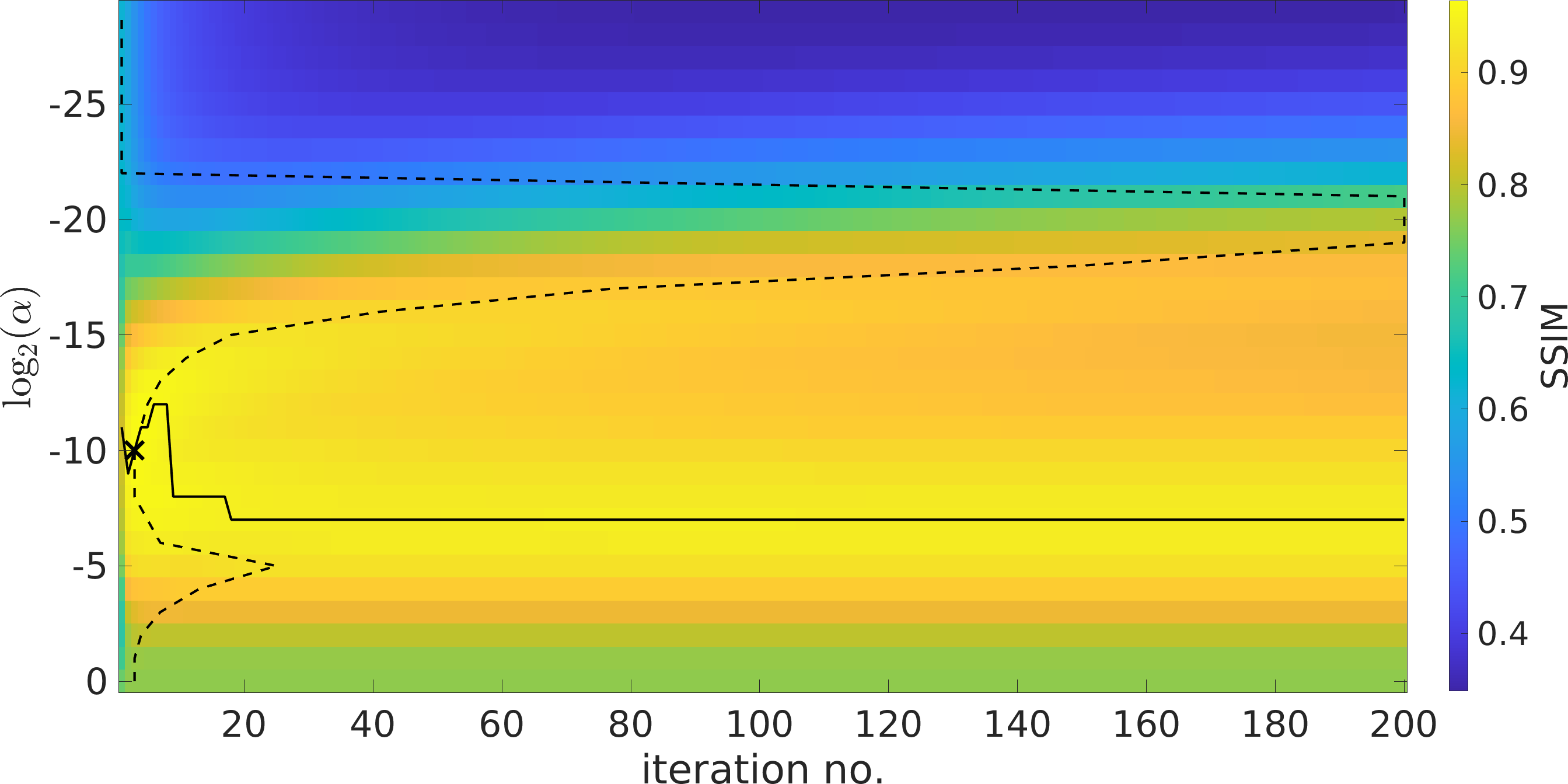} \\
\end{tabular}
}

\caption{Image quality measures with respect to $\alpha$ values and iteration number of l2-K  visualized for the ``shape'' phantom in the whitened case. The ``x'' marks the optimum. The solid line highlights the contour of the maximum image quality over $\alpha$ values for fixed iteration number $N$ (highlights the maximum of each column). The dashed line highlights the contour of the maximum image quality over iteration numbers $N$ for fixed $\alpha$ values (highlights the maximum of each row).    }
\label{fig:shape_Kaczmarz_image_quality_whitened}
\end{figure*}

\begin{figure*}[hbt!]%
\centering
\scalebox{1}{
\begin{tabular}{c|c}
PSNR & SSIM \\
 \hline
\multicolumn{2}{l}{$\tau=0$} \\
 \includegraphics[width=0.4\textwidth]{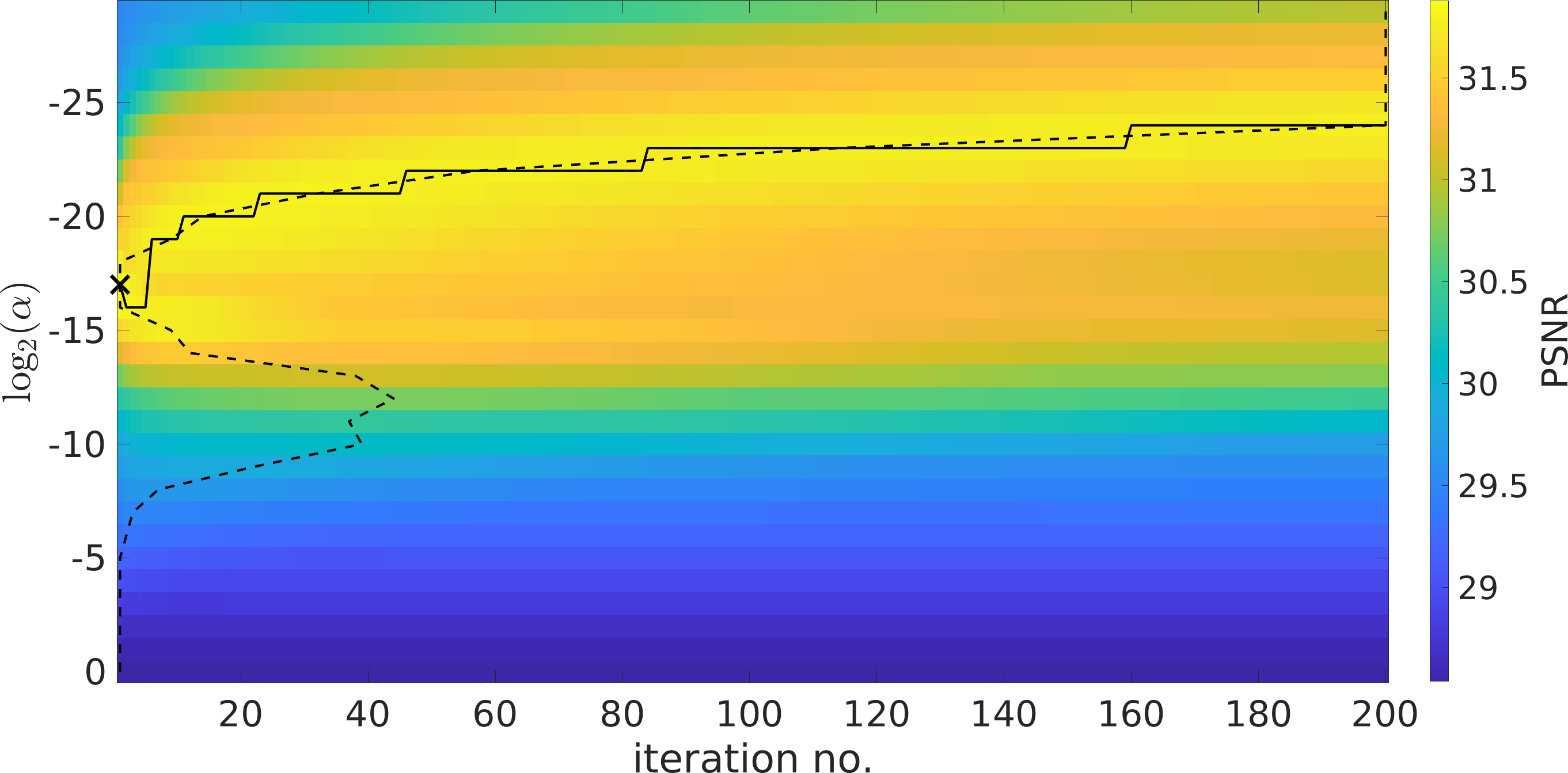} &
 \includegraphics[width=0.4\textwidth]{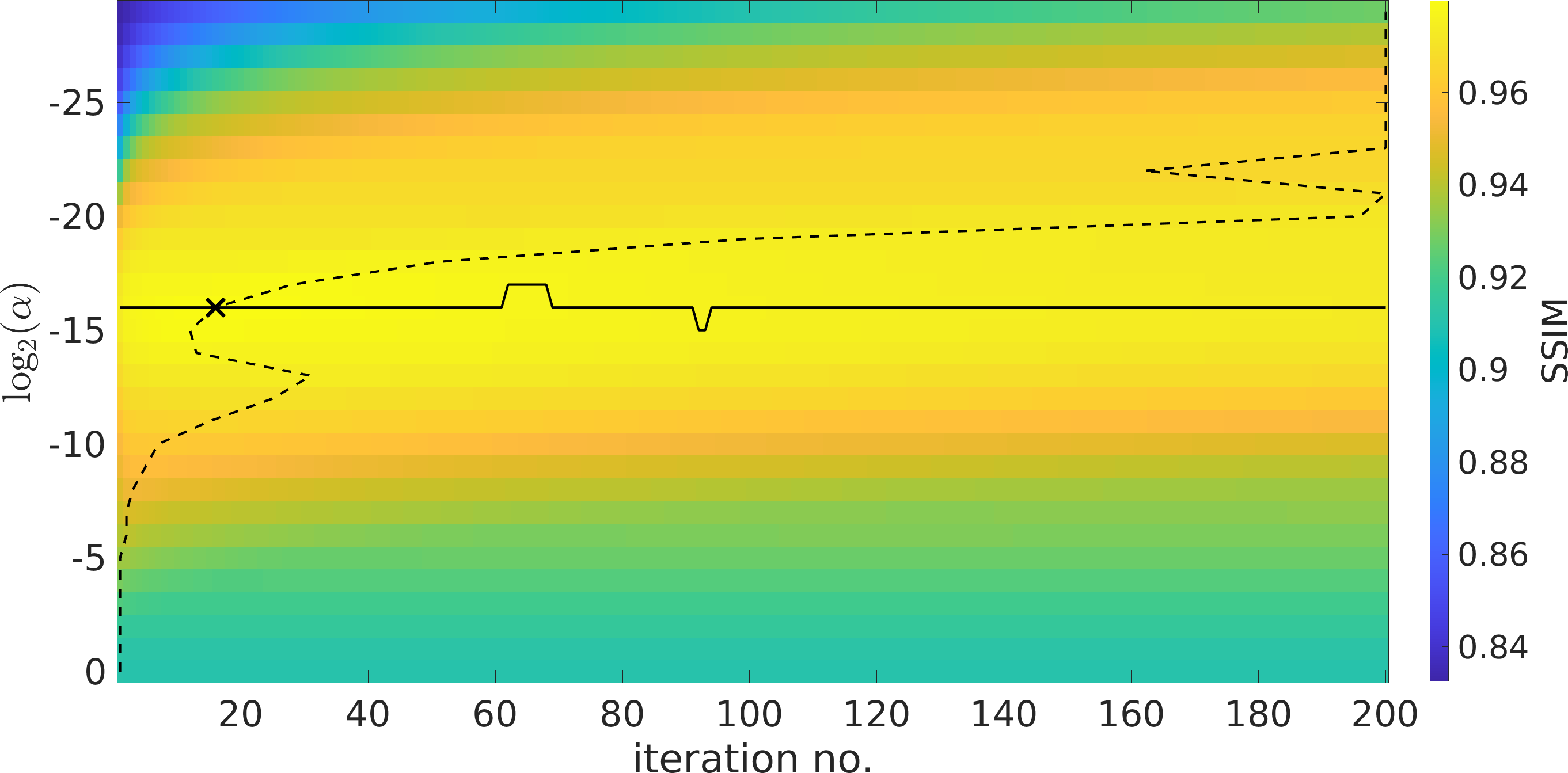} \\
 \hline
\multicolumn{2}{l}{$\tau=1$} \\
 \includegraphics[width=0.4\textwidth]{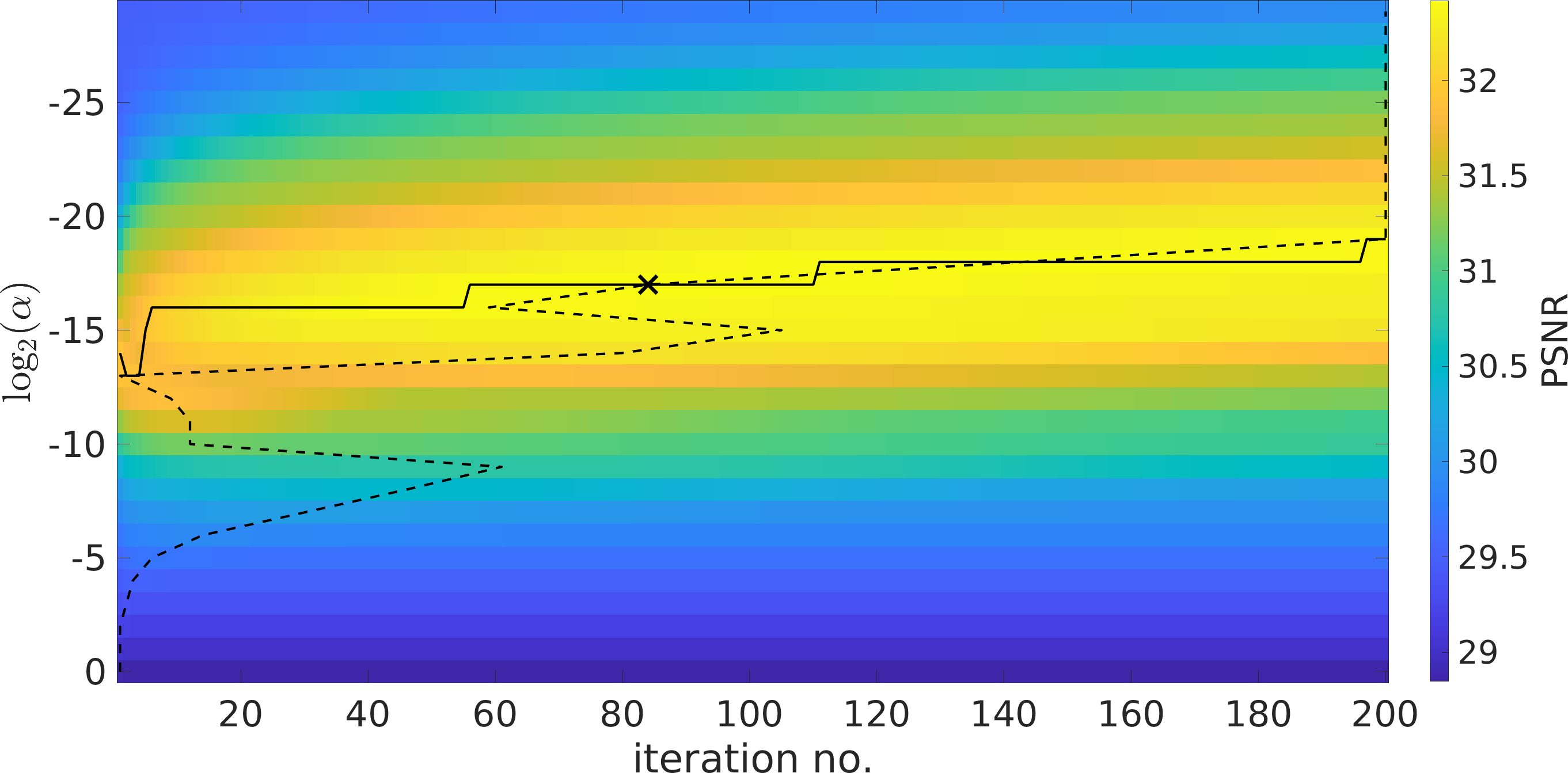} &
 \includegraphics[width=0.4\textwidth]{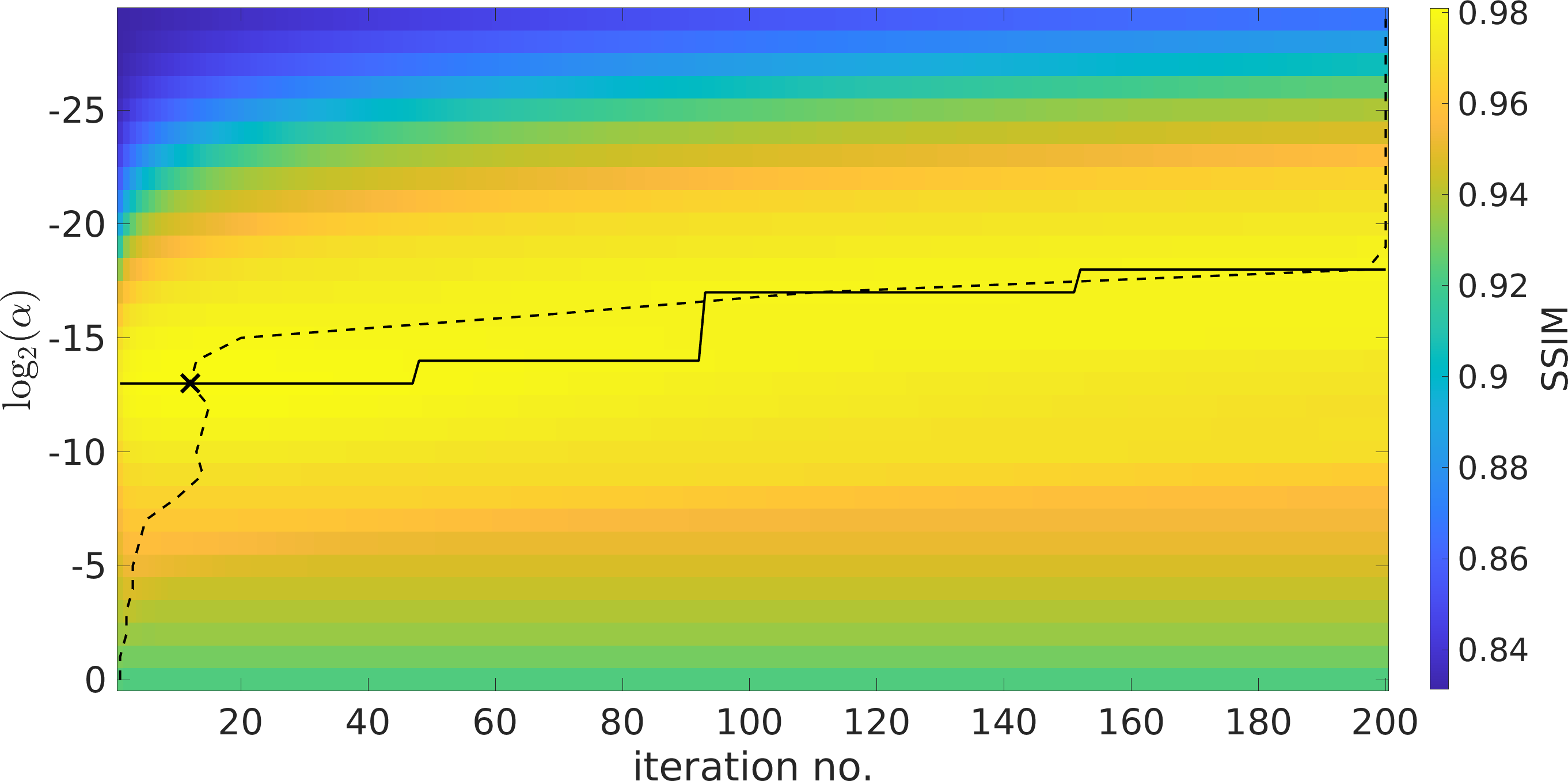} \\
  \hline
\multicolumn{2}{l}{$\tau=3$} \\
  \includegraphics[width=0.4\textwidth]{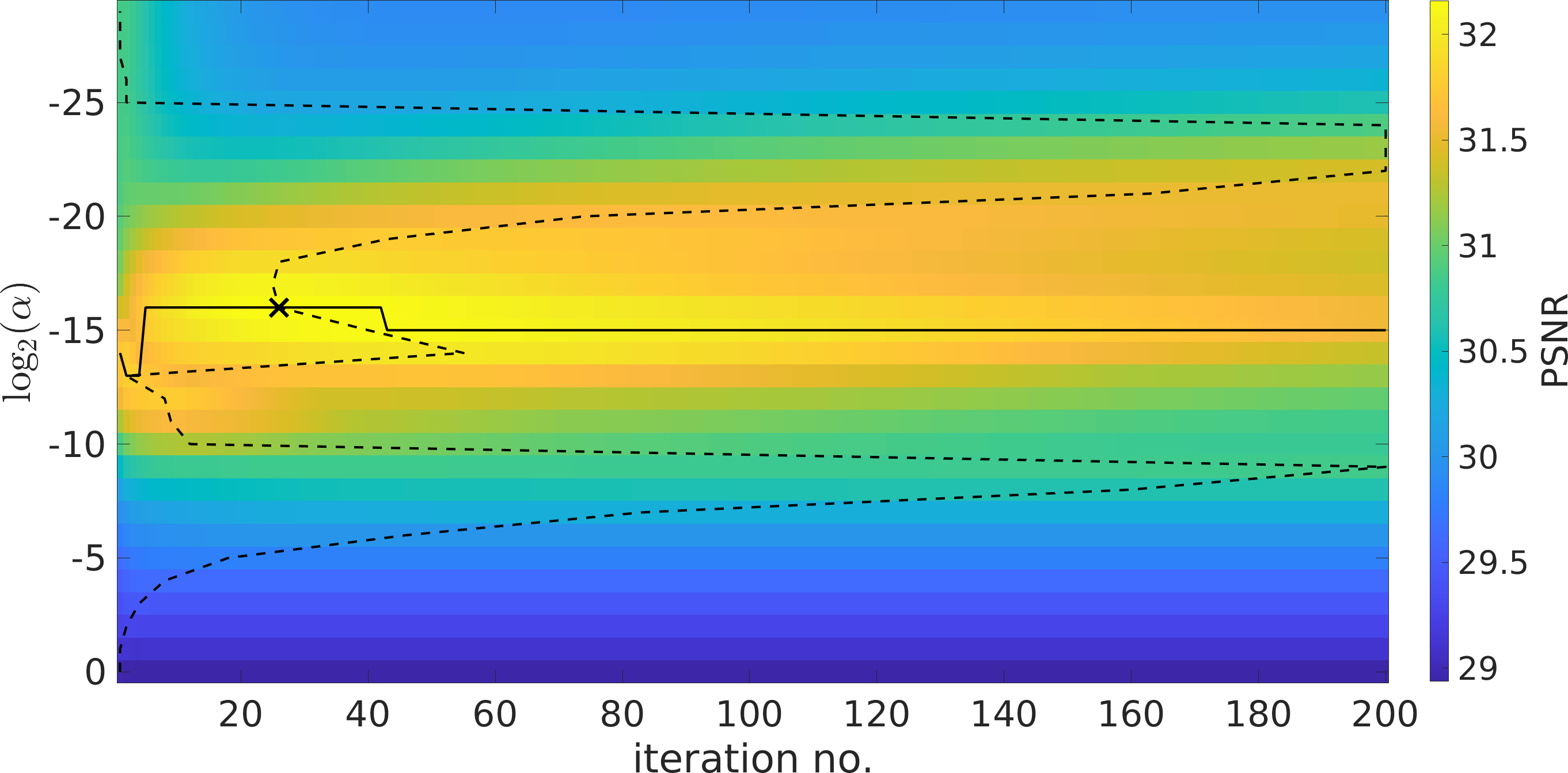} &
 \includegraphics[width=0.4\textwidth]{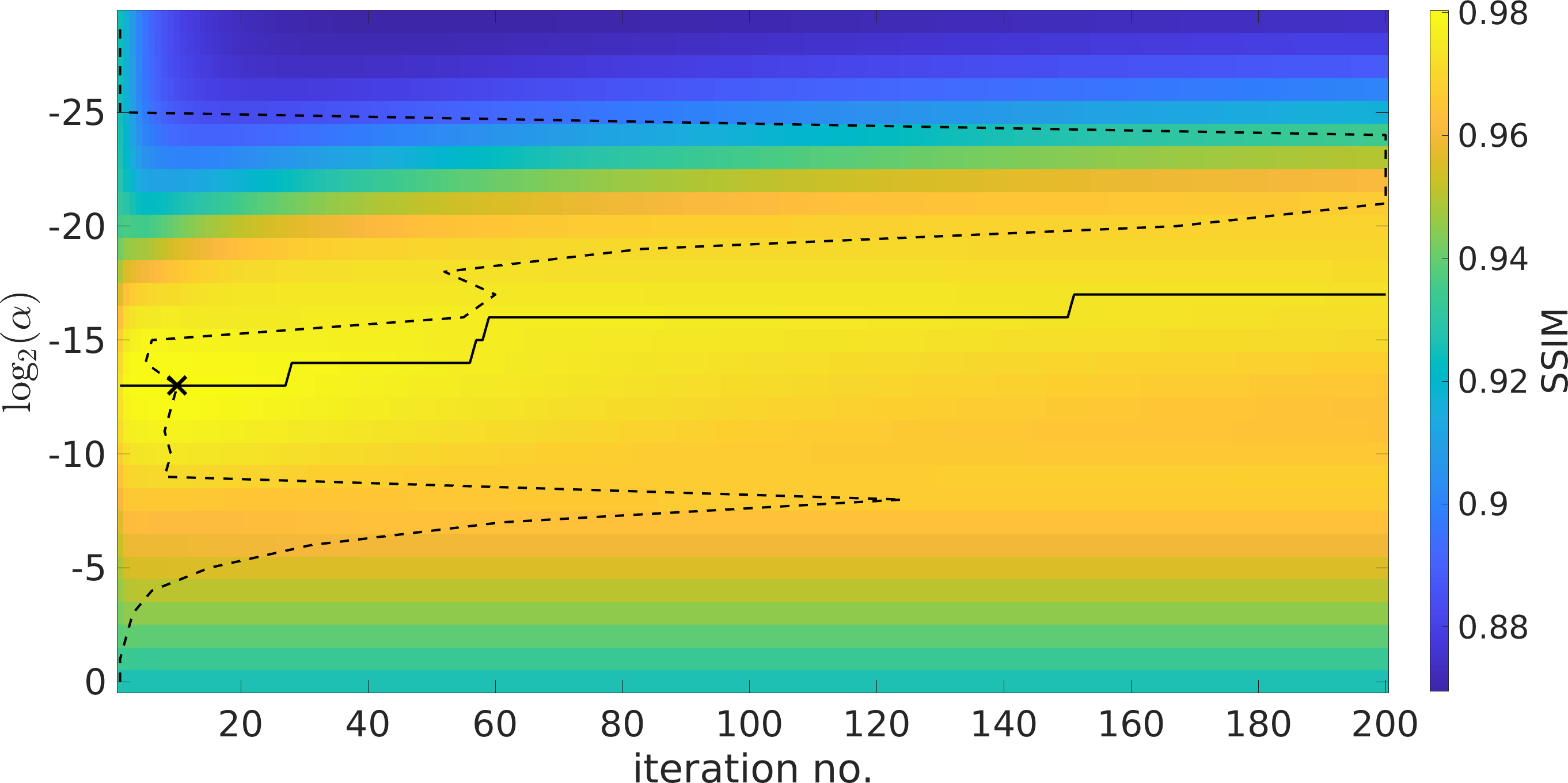} \\
  \hline
\multicolumn{2}{l}{$\tau=5$} \\
 \includegraphics[width=0.4\textwidth]{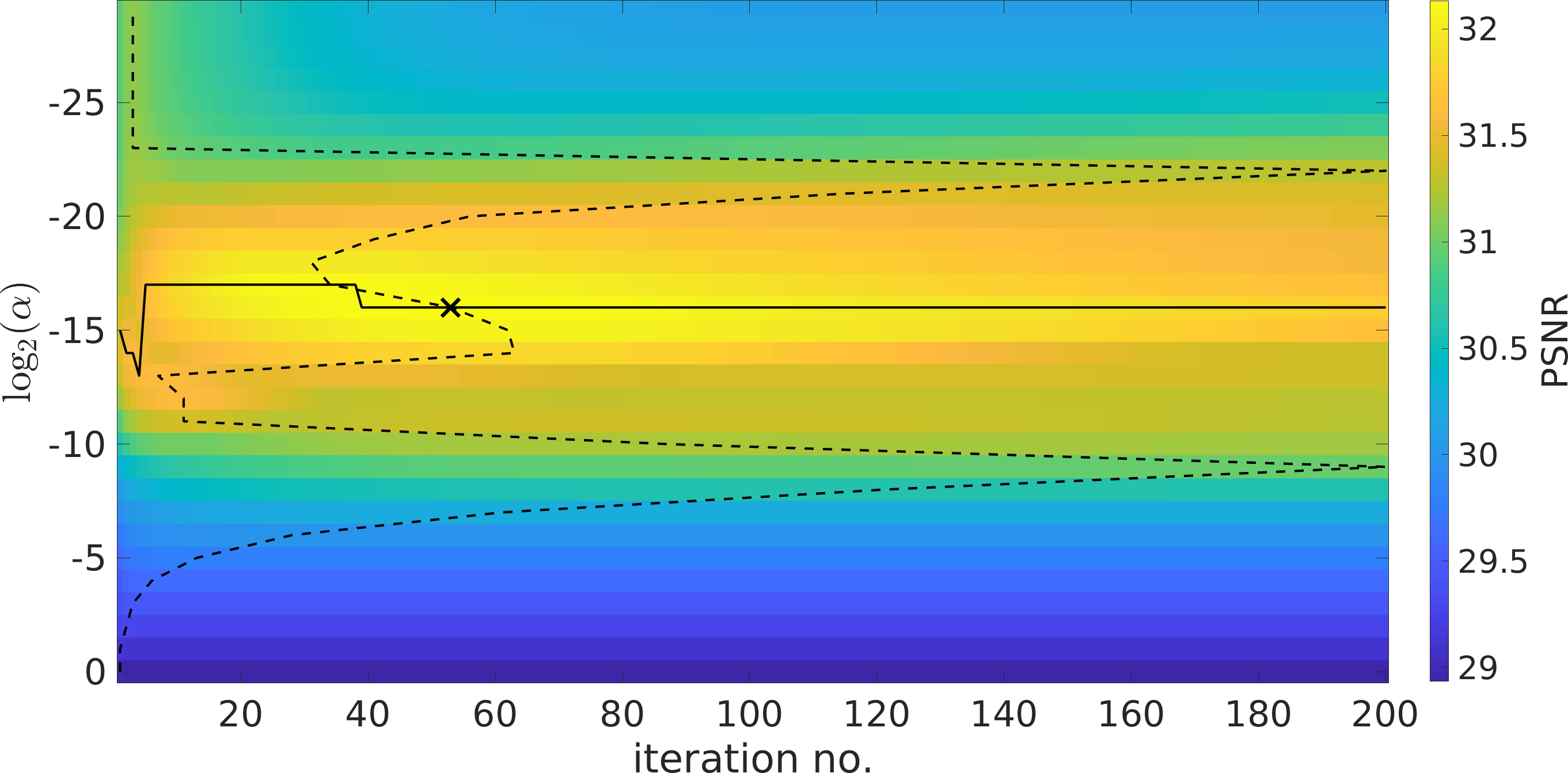} &
 \includegraphics[width=0.4\textwidth]{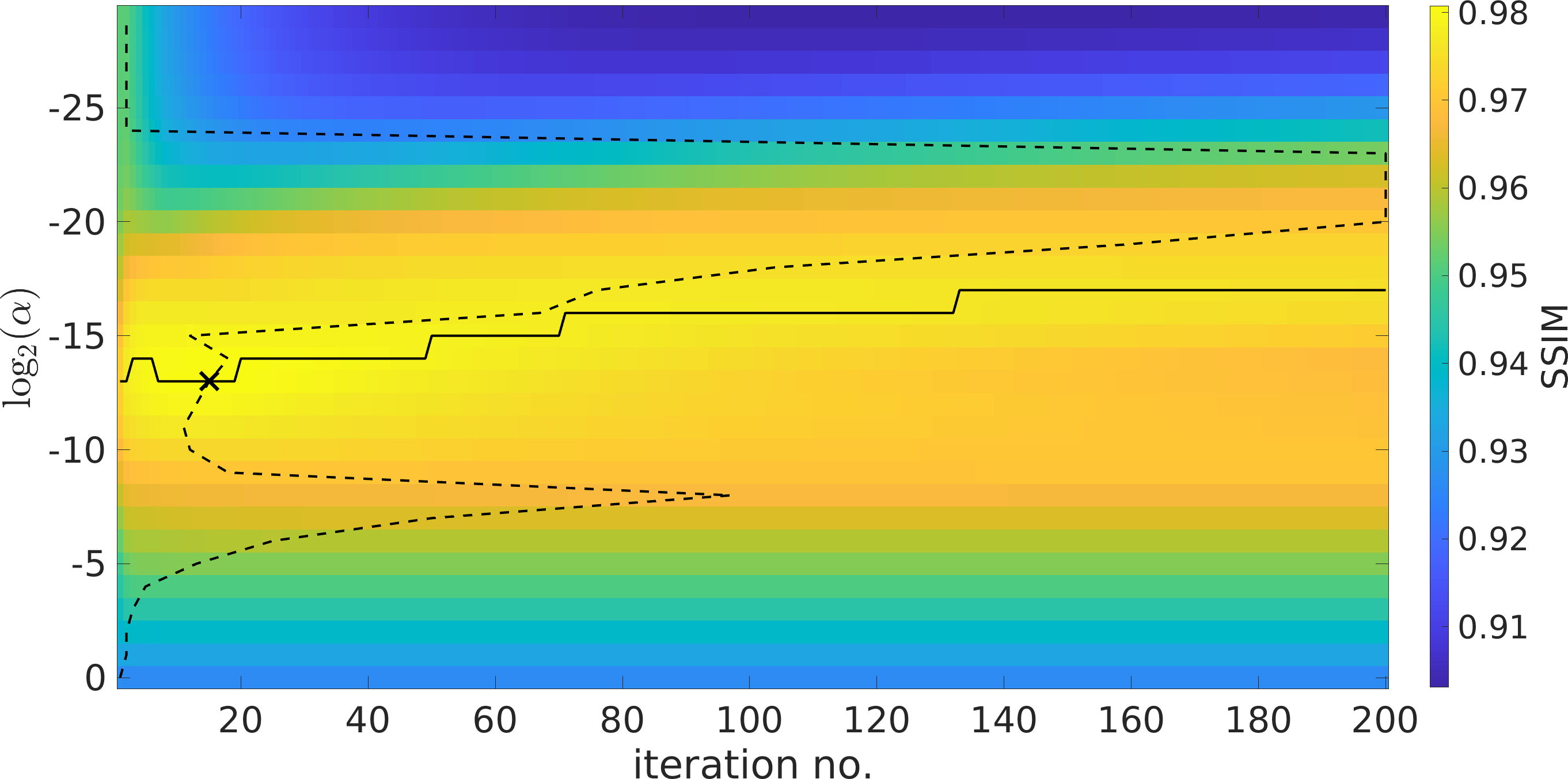} \\
\end{tabular}
}

\caption{Image quality measures with respect to $\alpha$ values and iteration number of l2-K  visualized for the ``resolution'' phantom in the whitened case. The ``x'' marks the optimum. The solid line highlights the contour of the maximum image quality over $\alpha$ values for fixed iteration number $N$ (highlights the maximum of each column). The dashed line highlights the contour of the maximum image quality over iteration numbers $N$ for fixed $\alpha$ values (highlights the maximum of each row).}
\label{fig:resolution_Kaczmarz_image_quality_whitened}
\end{figure*}

Now we present a more detailed comparative study of l2-K, which is one of the most commonly used techniques
in MPI. Since within l2-K, there are two different sources of regularizing effect (i.e., variational and iterative), we test
different iteration numbers $N$ in l2-K to shed further insights, and present the PSNR and SSIM results for ``shape'' and ``resolution'' phantom in Figs. \ref{fig:shape_Kaczmarz_image_quality_whitened} and
\ref{fig:resolution_Kaczmarz_image_quality_whitened} for the whitened case.
For the non-whitened case one can make analogous observations (see Figs. \ref{fig:shape_Kaczmarz_image_quality_nonwhitened} and \ref{fig:resolution_Kaczmarz_image_quality_nonwhitened} in the appendix).
This shows the delicate interplay of the different regularization techniques:
The smaller $\alpha$ (vertical direction) is, the less is the reconstruction influenced
by the variational regularization. The influence of an iterative regularization mechanism decreases with
larger iteration numbers (horizontal direction).

For the ``shape'' phantom, PSNR and SSIM results give similar pictures. In terms of PSNR, one already obtains an optimal
value after two iterations (i.e., $N=2$) and also for the SSIM case the necessary iteration numbers to reach the optimum are all below 20.
For the ``resolution'' phantom we can make analogous observations for SSIM but we observe larger iteration numbers up to 100 to reach the optimal PSNR.
Nevertheless, 100 iterations are still insufficient to reach convergence of the Kaczmarz methods as the l2-L results (if l2-K does reach convergence) still deviate largely. In particular, for $\tau\geq3$ the trend of the solid lines already indicates the optimal $\alpha$ values found for l2-L,
showing the beneficial regularizing effect built into l2-K due to early stopping.
However, Figs. \ref{fig:shape_Kaczmarz_image_quality_whitened} and
\ref{fig:resolution_Kaczmarz_image_quality_whitened} indicate that sole iterative regularization is not enough
as l2-K fails for small $\alpha$ tending to zero (see trend from bottom to top in each plot).
Interestingly, sole variational regularization also fails, since the image quality decreases for larger
iteration numbers in l2-K (approaching the worse l2-L results for sufficiently large $N$).
This indicates the need of tuning multiple parameters, which is in practice undesirable.
Provided that the optimal stopping index (and $\alpha$) is found, l2-K yields
superior PSNR values and at least superior SSIM values for small $\tau$. For $\tau\geq 3$ l1-L is a pure variational alternative depending on one single parameter only.

\begin{table*}[hbt!]
\centering
\caption{ The mean and standard deviation of computing times are computed over 30 $\alpha$ values and 2 phantoms.
The number in the bracket is the iteration number $N$ for l2-K.
Computations are carried out on a server with 2$\times$Intel\textsuperscript{\textregistered} Xeon\textsuperscript{\textregistered} Broadwell-EP
Series Processor E5-2687W v4, 3.00 GHz, 12-Core, and 1.5 TB DDR4 PC2666 main memory.}\label{tab:cmp_time}
 \scalebox{0.95}{
\begin{tabular}{c|c|cccc}
& \multicolumn{5}{c}{non-whitened} \\
\hline
  $\tau$ & l1-L  & l2-K(200) & (50) & (20) & (1) \\
\hline
0 &  $ 71.09  \pm  91.46 $  &  $ 556.64 \pm 12.21 $  &  $ 140.78 \pm 3.66$  &  $ 58.96 \pm 1.84$  &  $ 7.08\pm 0.25 $  \\
1 &  $ 125.60 \pm 151.77 $  &  $ 547.17 \pm 11.74 $  &  $ 137.53 \pm 3.54$  &  $ 57.25 \pm 1.88$  &  $ 6.92\pm 0.26 $  \\
3 &  $ 43.30  \pm  26.91 $  &  $  78.60 \pm  2.87 $  &  $  19.74 \pm 0.72 $  &  $ 8.19 \pm 0.26 $  & $ 0.96 \pm 0.02 $  \\
5 &  $ 33.51  \pm  17.32 $  &  $  48.62 \pm  2.57 $  &  $  12.64 \pm 0.46 $  &  $ 5.34 \pm 0.22 $  & $ 0.66 \pm 0.02 $  \\
\hline

\end{tabular}
 }
\vspace{0.2cm}
 \scalebox{0.95}{
\begin{tabular}{c|c|cccc}
& \multicolumn{5}{c}{whitened} \\
\hline
  $\tau$ & l1-L  & l2-K(200) & (50) & (20) & (1) \\
\hline
0 &  $  23.63 \pm  21.44$  &   $ 552.91  \pm 9.69 $  &  $ 139.68\pm 3.52$  &  $ 58.82 \pm 1.93$  &  $ 7.27\pm 0.24$  \\
1 &  $ 113.24 \pm 130.26 $  &  $ 544.22 \pm 10.16 $  &  $ 137.10\pm 3.14$  &  $ 57.43 \pm 1.81$  &  $ 6.97\pm 0.32$  \\
3 &  $ 48.61 \pm   25.53 $  &   $ 78.37  \pm 2.44 $  &  $  19.59\pm 0.78 $  &  $ 8.14 \pm 0.28$  &  $ 0.97 \pm 0.02$  \\
5 &  $ 34.11 \pm   17.21 $  &   $ 48.61  \pm 1.96 $  &  $  12.60\pm 0.47 $  &  $ 5.30 \pm 0.22 $  &  $0.65 \pm 0.03$  \\
\hline
\end{tabular}
}
\end{table*}

Finally, the computing time for the methods is summarized in Table \ref{tab:cmp_time}, which shows that computationally l1-L is comparable
with l2-K for $N=50$. Note that in the literature the Kaczmarz method is also sometimes exploited for online reconstruction \cite{Knopp_Online_2016,KluthJin:2019}.
This can be realized using a dimension reduction technique together with a sufficiently small number of iterations.

\section{Concluding remarks and comments}\label{sec:conc}

In this work we have investigated the potential of the l1 data fitting for MPI reconstruction and compared it with the
standard method l2-K in MPI. After applying a bandpass filter only, one can observe severe outliers
in the MPI signal which are of orders of magnitudes larger and prevent obtaining reasonable reconstructions. A data fidelity
term based on the l1 norm has been successfully applied to various applications, where the noise is characterized by
severe outliers. However, within MPI, the l1 fitting can only solve the problem to a certain extent,
and the recommendation is to combine it with the established SNR-type frequency selection (which still results in an MPI signal with large outliers).
Then the l1 fitting can compete with the standard method's reconstruction performance for both non-whitened and whitened cases in terms of
the popular PSNR and SSIM measures and visual quality. However, the l1 fitting only relies on one tunable parameter compared to two
in the standard method, which is advantageous for the development of one-click solutions for applicants. To the best of our knowledge, we have presented a
first quantitative study on phantom MPI data with respect to image quality measures for \texttt{Open MPI dataset},
which shows the benefits of using l1 fitting for pure variational methods, and we quantitatively explored the beneficial interplay iterative and variational regularization in the standard method. There is strong implicit regularization built into the popular l2-K, whose precise mechanism is to be ascertained.

In the context of l2 fitting, whitening is known to be beneficial for the reconstruction \cite{KluthJin:2019}. This suggests
that whitening might adjust the noise characteristic such that it is closer to the i.i.d. Gaussian case, for which the l2
fitting is most suitable. Numerically, this clearly allows enhancing the reconstruction quality.
The influence of whitening on the l1 fitting is less dramatic due to its robustness with respect to outliers.
Variational regularization methods need to respect the MPI noise characteristic in order to compete with
l2-K, which is confirmed by the failure of l2-L and the success of l1-L for larger SNR-thresholds in terms of SSIM.
Besides MPI reconstructions, these findings also can have implications for the calibration procedure in the model-based approach. Whenever model parameters
have to be identified for the purpose of system calibration, the noise characteristic should
be properly accounted for. Furthermore, simultaneous background-removal and image reconstruction approaches \cite{KluthMaass:2017,StraubSchulz:2018} have only exploited Gaussian assumptions on the noise so far and might also benefit from an extension taking into account the noise characteristic by l1 fitting terms in the respective variational approach. The identification of the noise in the system and its proper modeling are also related but they are different directions of research and deserve further research.

This study has only focused on the influence of the data fidelity, and does not touch the important issue of penalty for
best possible image reconstruction. Advanced variational penalties \cite{KluthMaass:2017,IlbeyTop:2019}, e.g., l1, total variation and their variants
and more recent learning based approaches \cite{DittmerKluth:2018}, promise highly desirable features, e.g., edge
preservation, at the expense of increased computational efforts, but largely remain to be systematically explored,
naturally also with the l1 fitting. In addition the development of a larger phantom dataset or extension of the existing one as a benchmark for reconstruction methods is highly desirable.  We leave these important issues to future works.

\section*{Acknowledgements}

T. Kluth acknowledges funding by the Deutsche Forschungsgemeinschaft
(DFG, German Research Foundation) - project number 281474342/GRK2224/1 ``Pi$^3$ : Parameter Identification
- Analysis, Algorithms, Applications'' and support by the project ``MPI$^2$'' funded by the Federal Ministry of Education and Research (BMBF, project no.
05M16LBA). The work of B Jin is partly supported by UK EPSRC EP/T000864/1.

\bibliographystyle{IEEEtran}
\bibliography{literature}

\newpage
\appendix
\section{Supplementary material: Higher harmonics in the background signal}

\begin{figure}[hbt!]
\centering
\scalebox{0.9}{
\begin{minipage}{0.3\textwidth}
\centering
Real part ($x$ receive coil) \\
$f_x$ \\
\includegraphics[width=\textwidth]{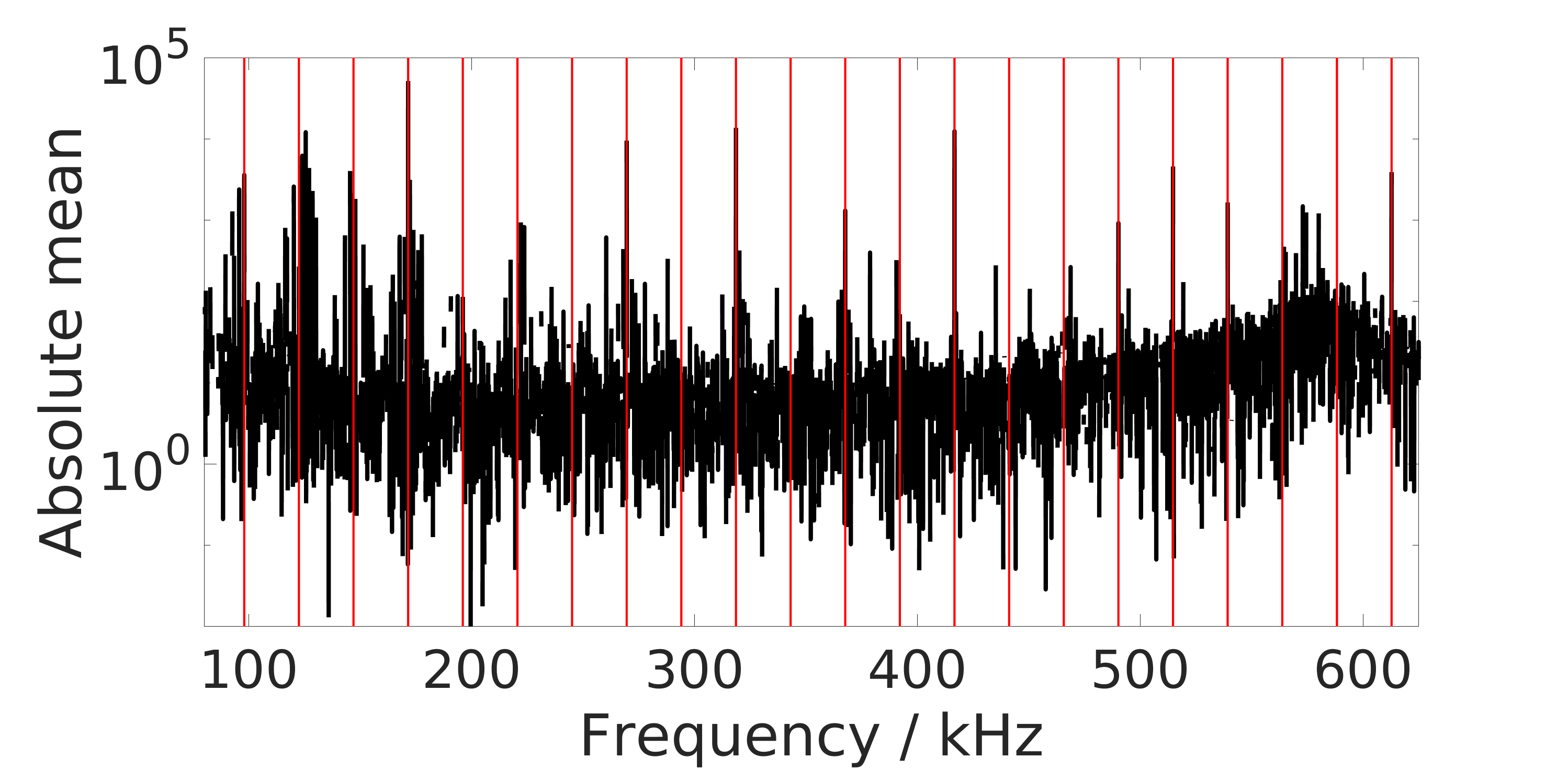} \\
\includegraphics[width=\textwidth]{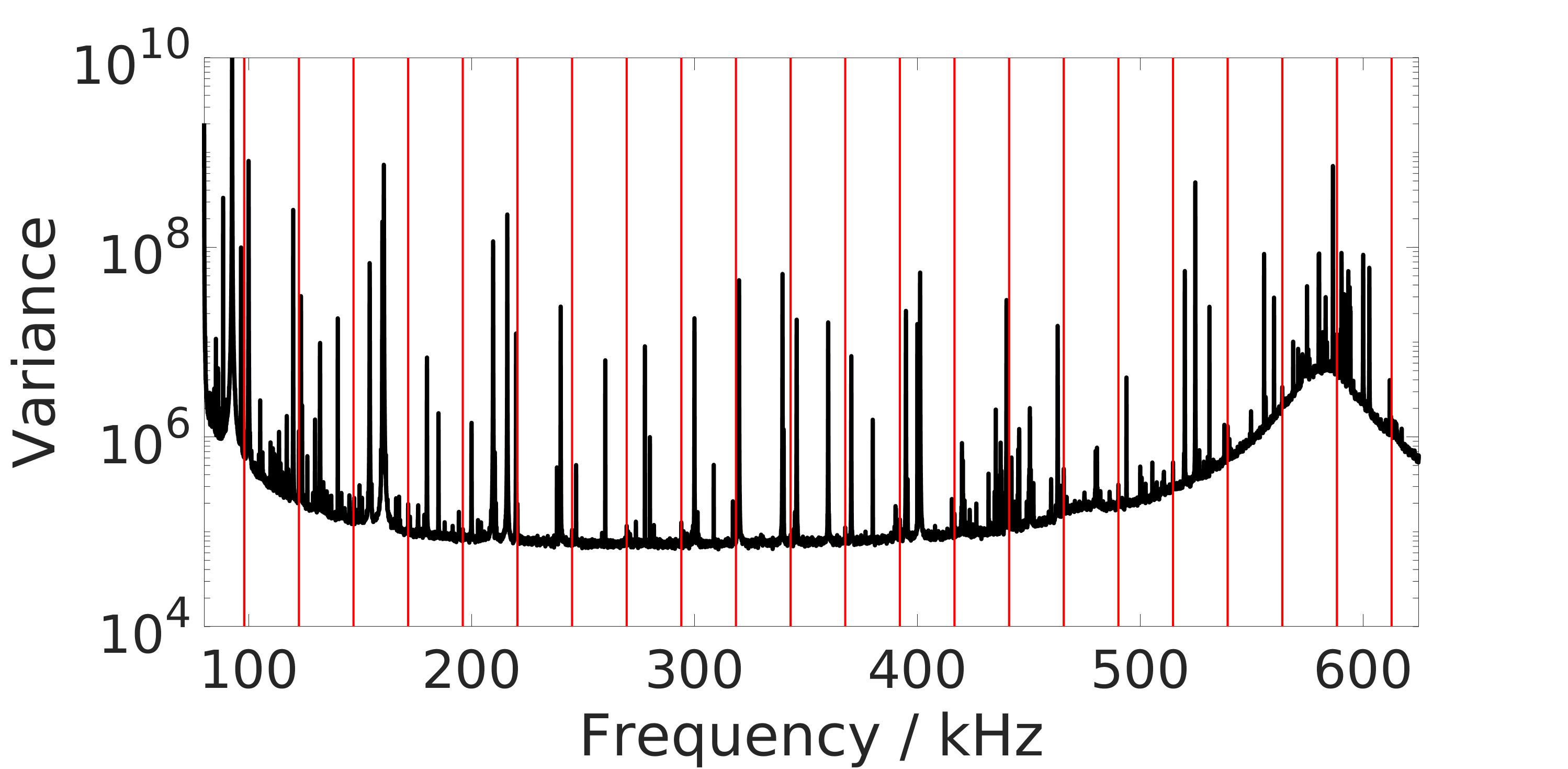} \\
 \end{minipage}
 \begin{minipage}{0.3\textwidth}
\centering
\ \\
$f_y$ \\
\includegraphics[width=\textwidth]{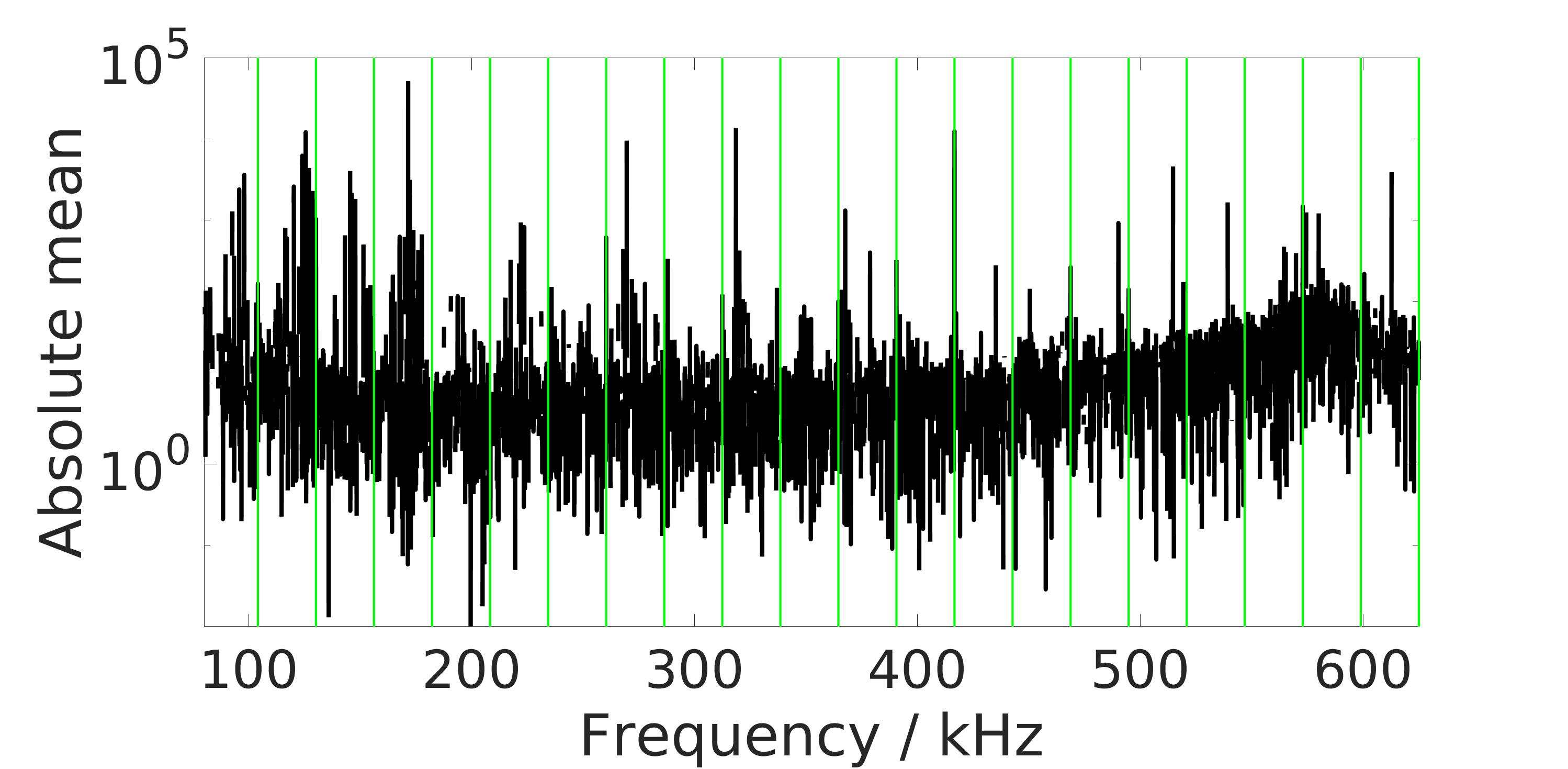} \\
\includegraphics[width=\textwidth]{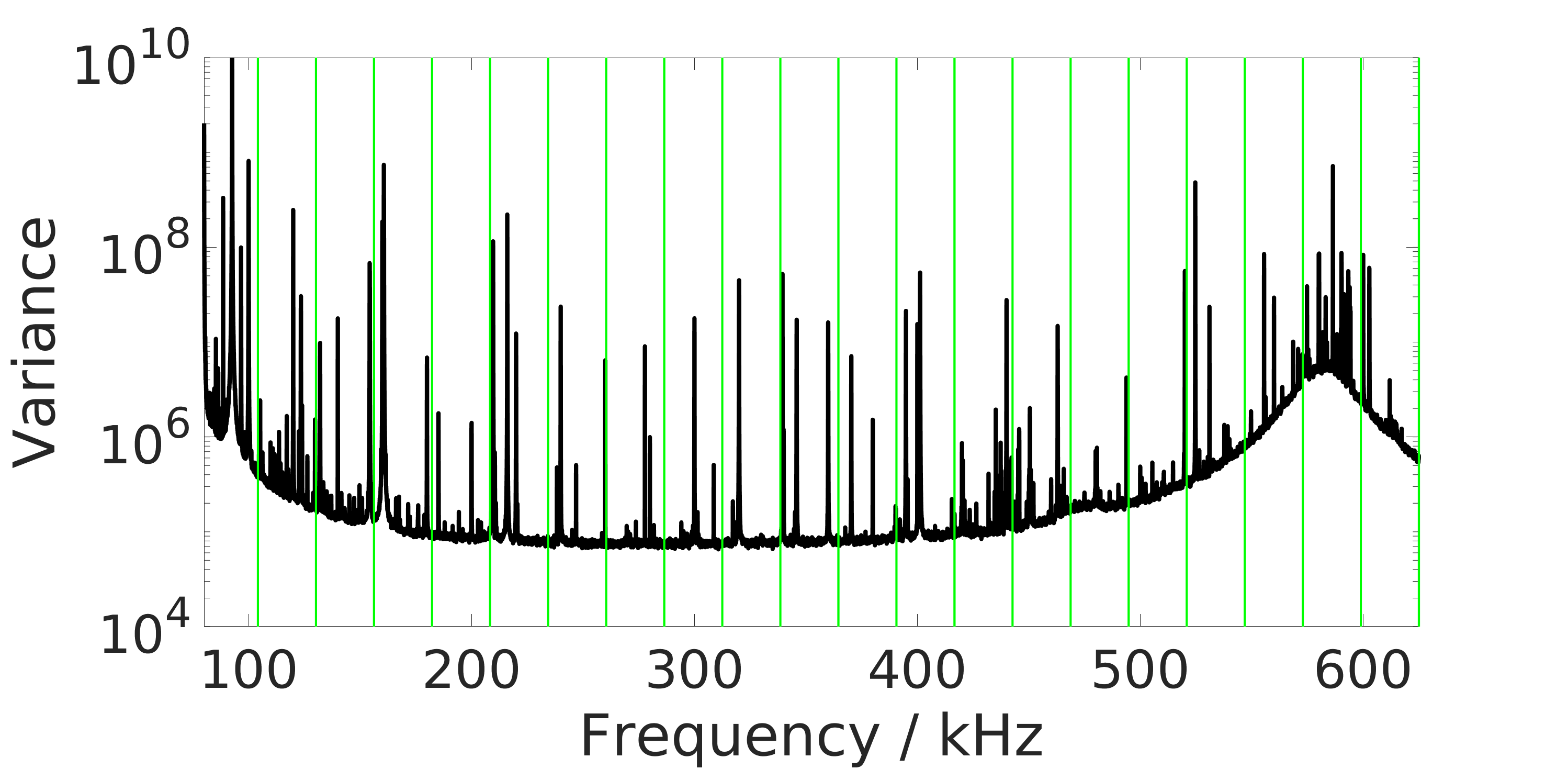} \\
 \end{minipage}
\begin{minipage}{0.3\textwidth}
\centering
\ \\
$f_z$ \\
\includegraphics[width=\textwidth]{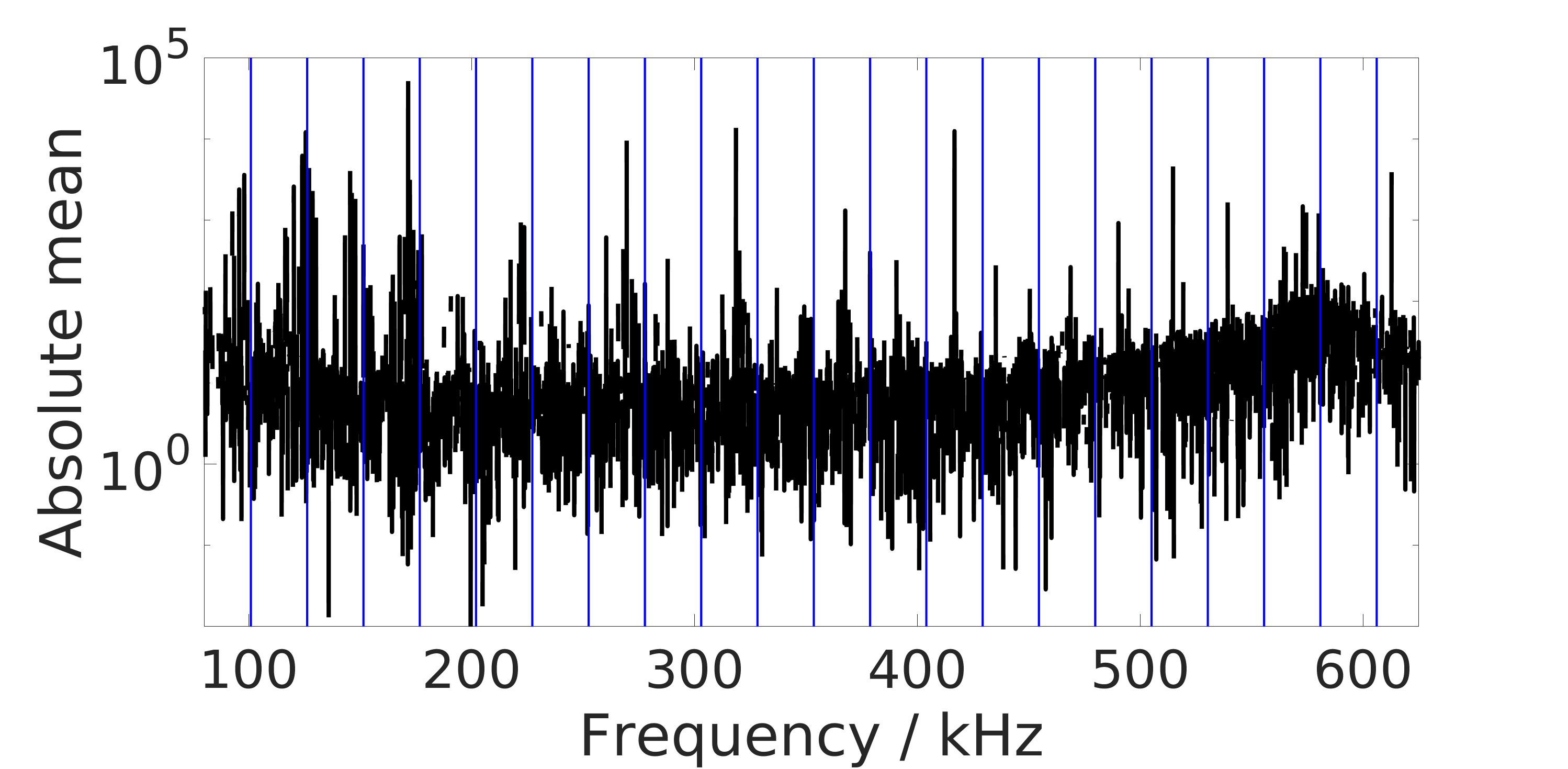} \\
\includegraphics[width=\textwidth]{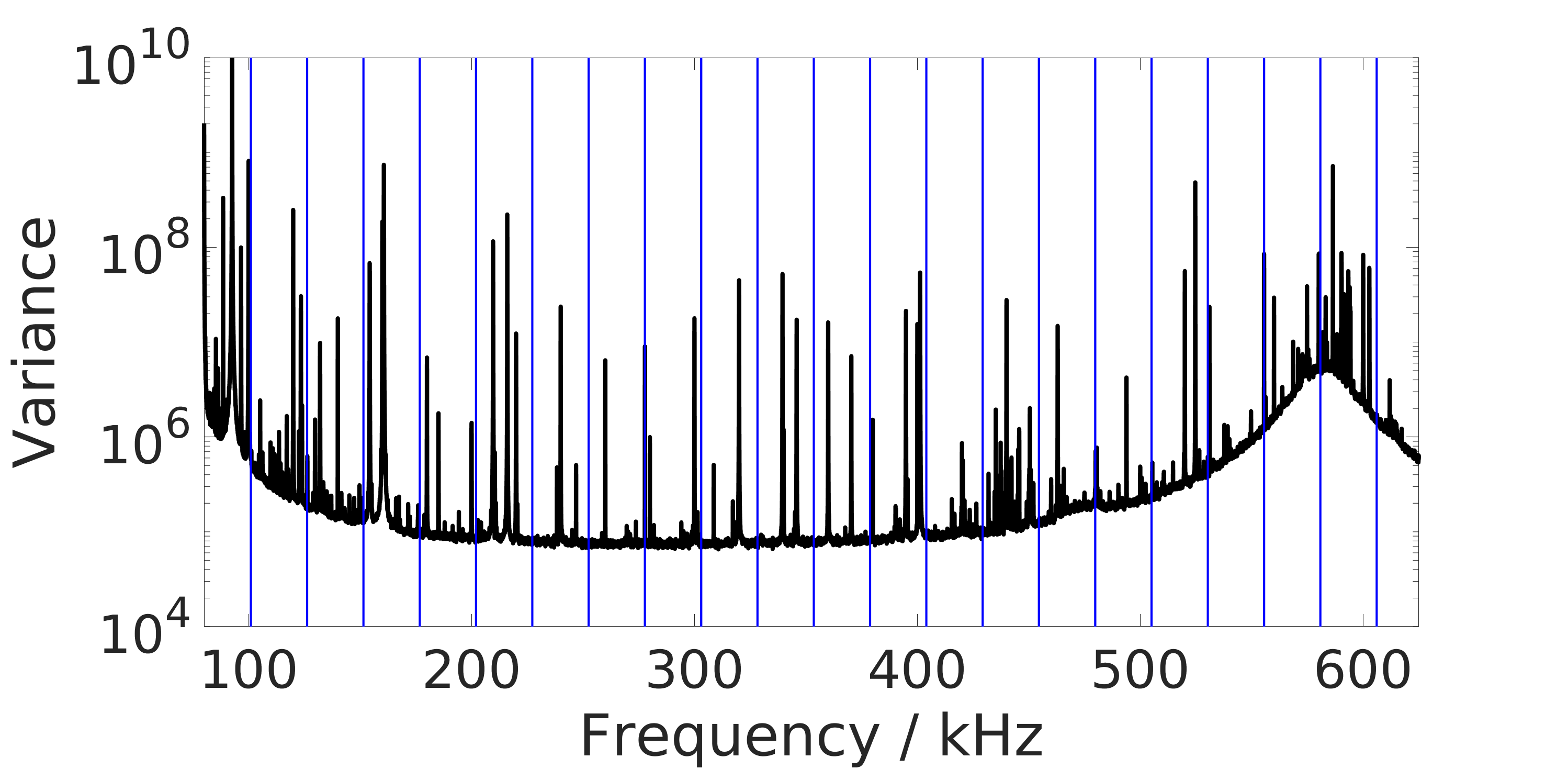} \\
 \end{minipage}
}

 \vspace{0.8cm}
 \scalebox{0.9}{
 \begin{minipage}{0.3\textwidth}
\centering
Imaginary part ($x$ receive coil) \\
$f_x$ \\
\includegraphics[width=\textwidth]{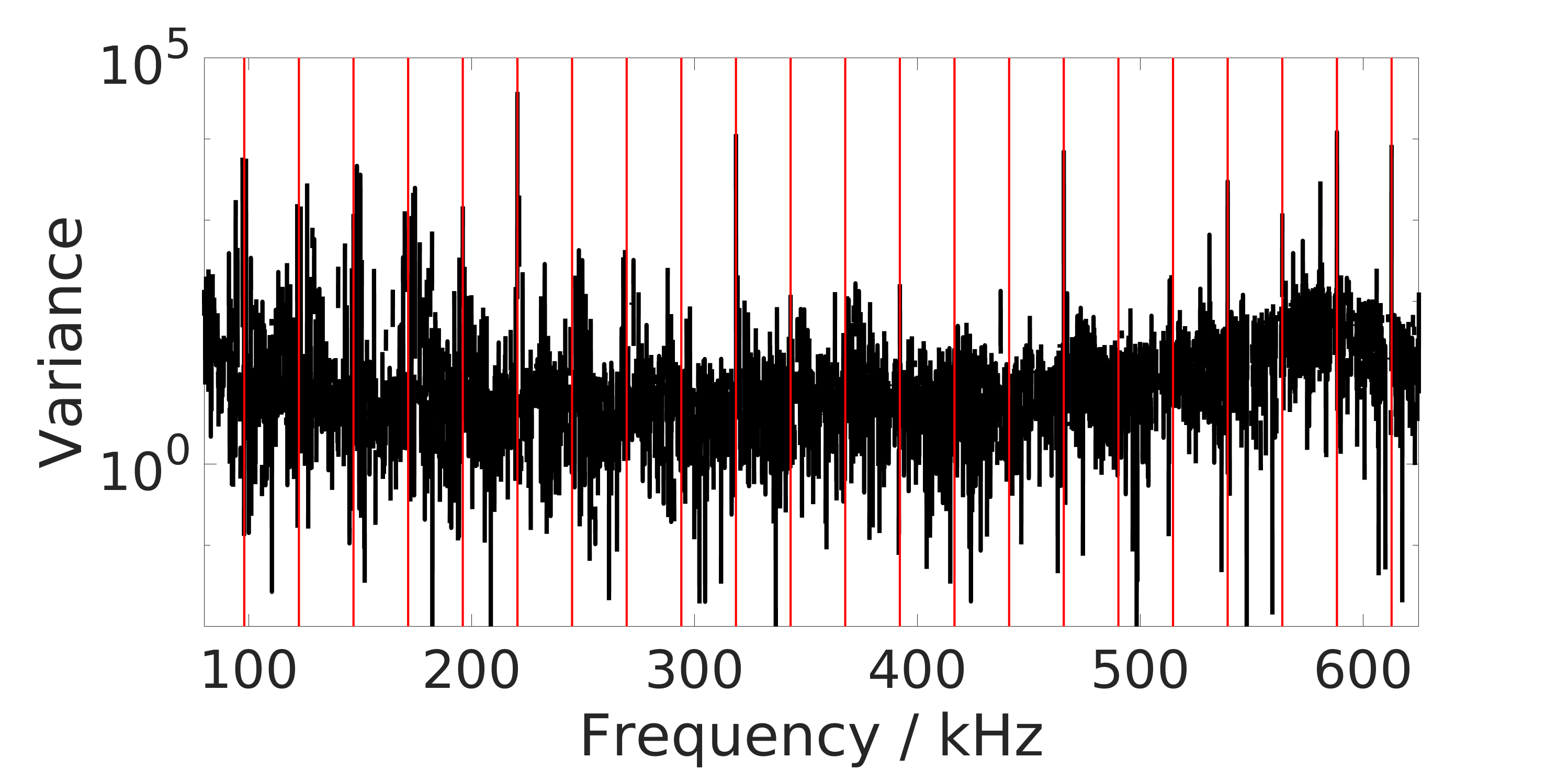} \\
\includegraphics[width=\textwidth]{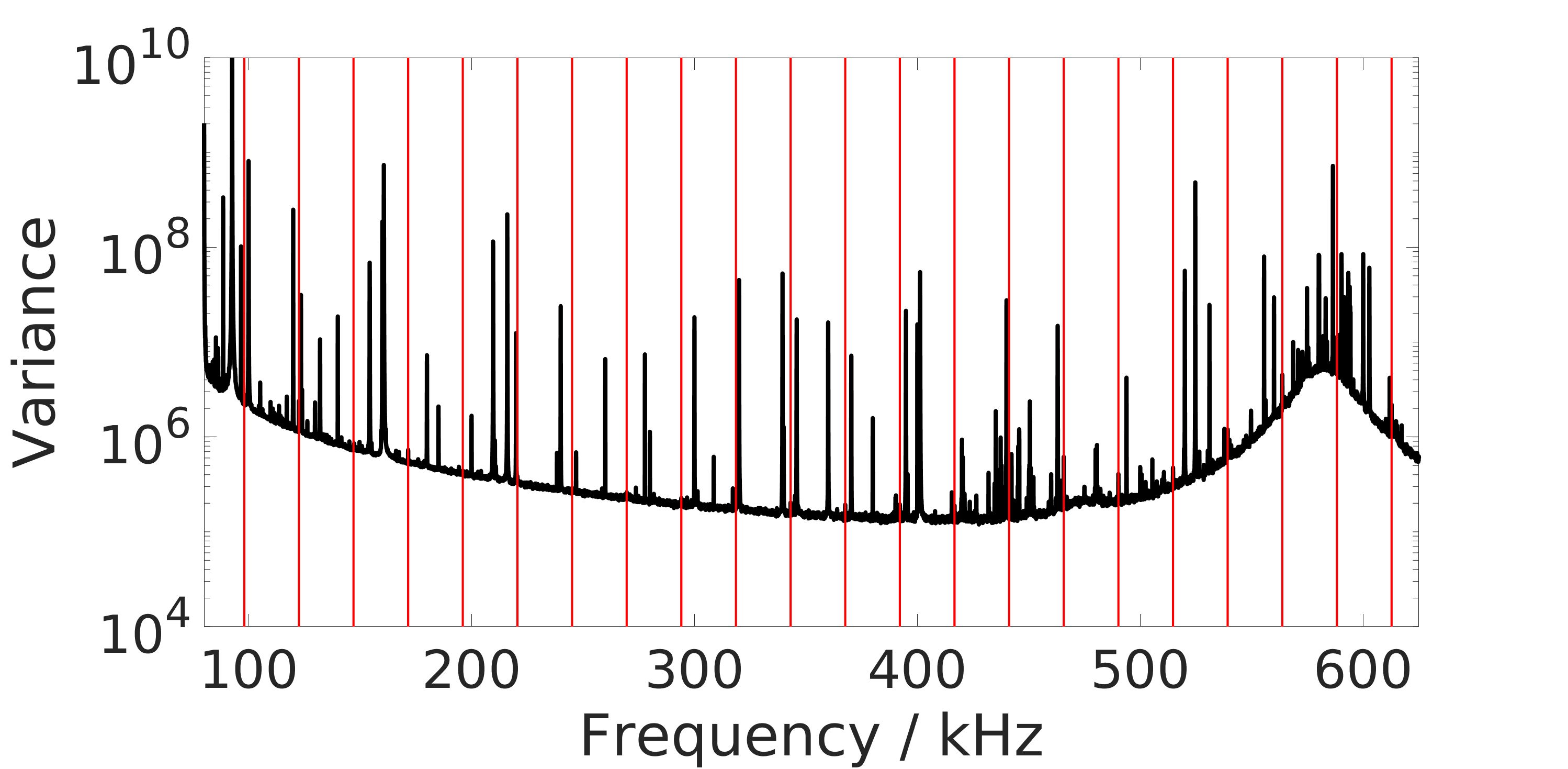} \\
 \end{minipage}
 \begin{minipage}{0.3\textwidth}
\centering
\ \\
$f_y$ \\
\includegraphics[width=\textwidth]{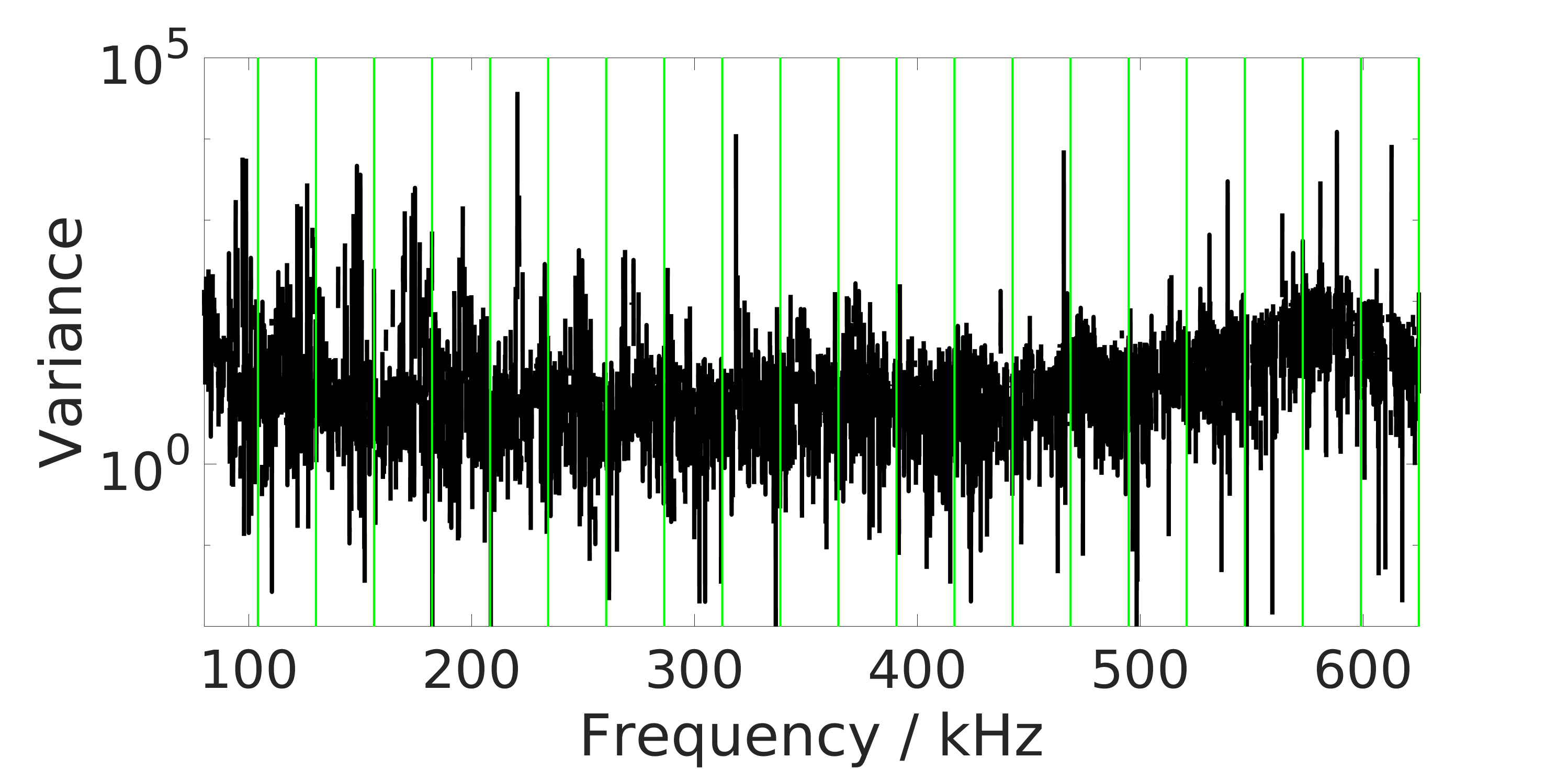} \\
\includegraphics[width=\textwidth]{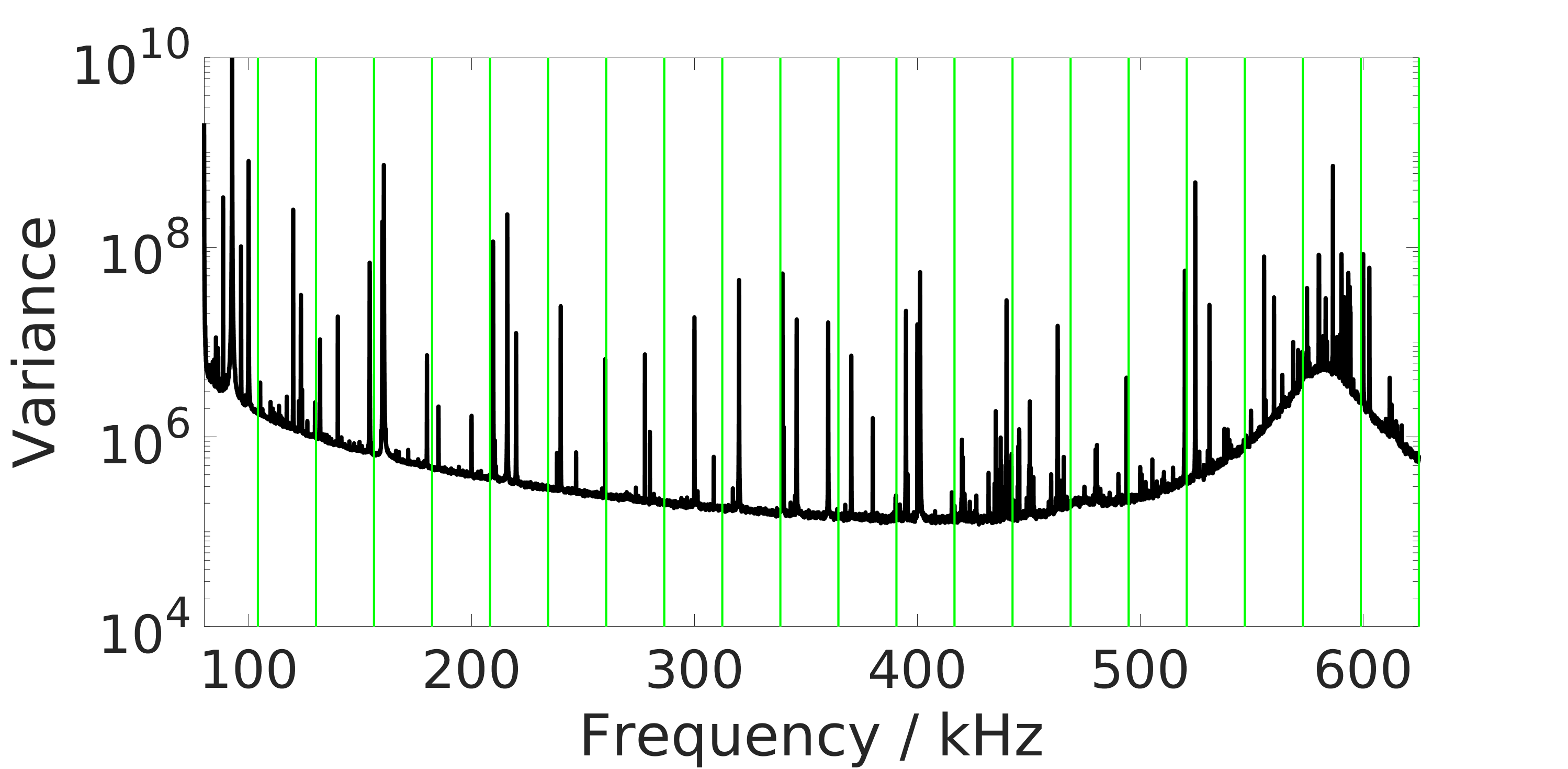} \\
 \end{minipage}
\begin{minipage}{0.3\textwidth}
\centering
\ \\
$f_z$ \\
\includegraphics[width=\textwidth]{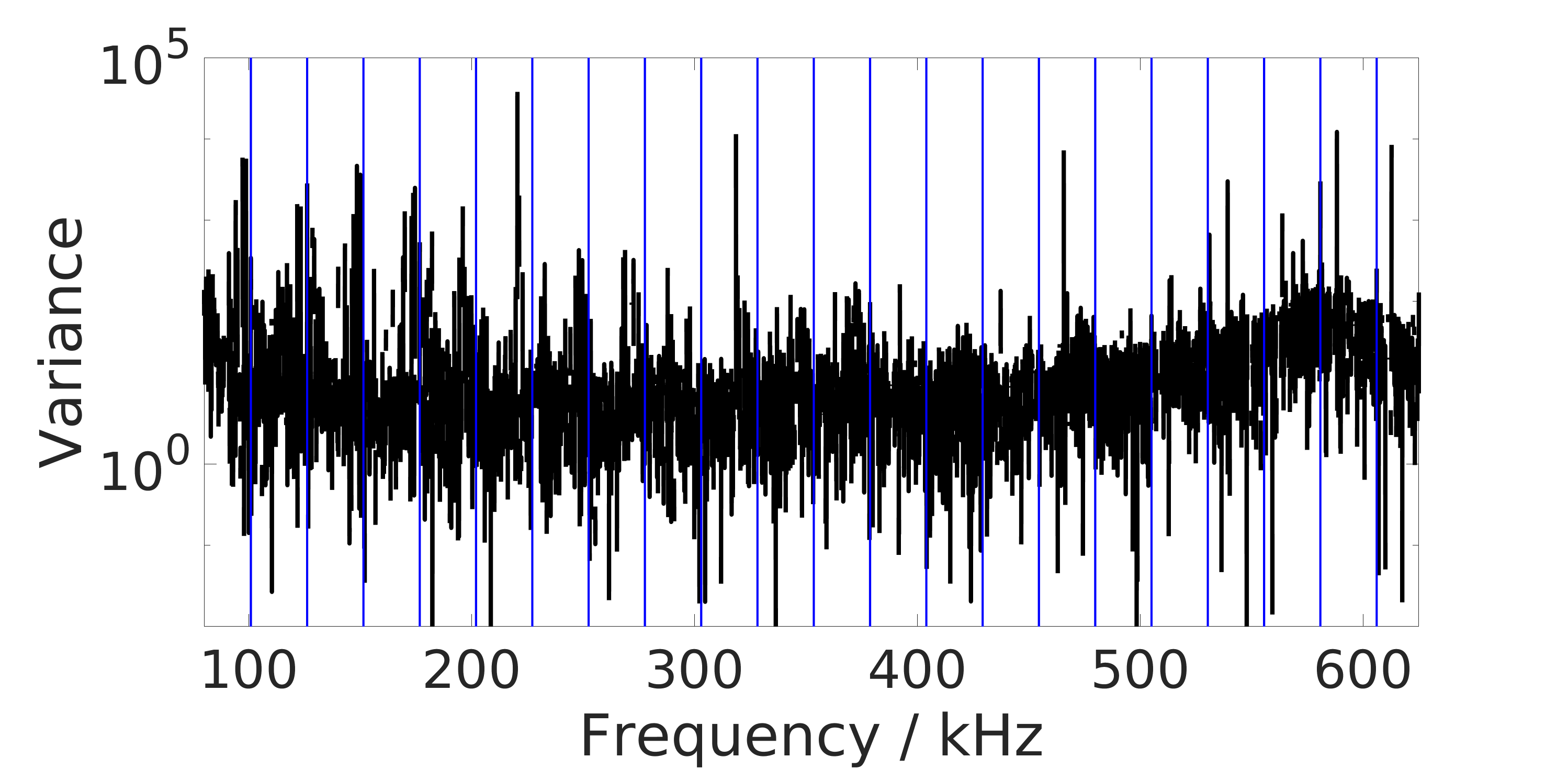} \\
\includegraphics[width=\textwidth]{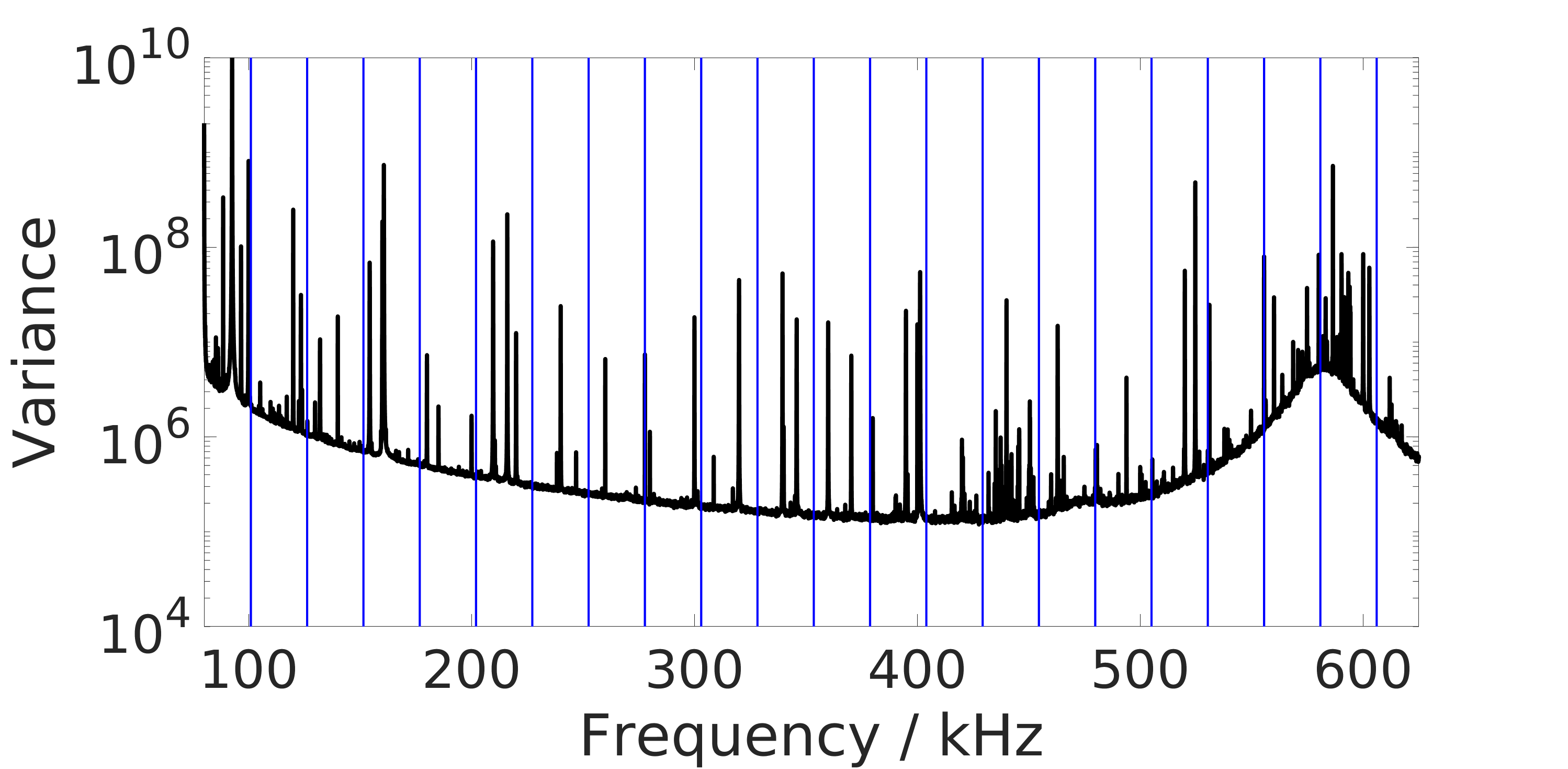} \\
 \end{minipage}
 }
 \vspace{0.8cm}
\scalebox{0.8}{
 \begin{minipage}{0.45\textwidth}
\centering
 All ($x$ receive coil) \\
Real part \\
\includegraphics[width=\textwidth]{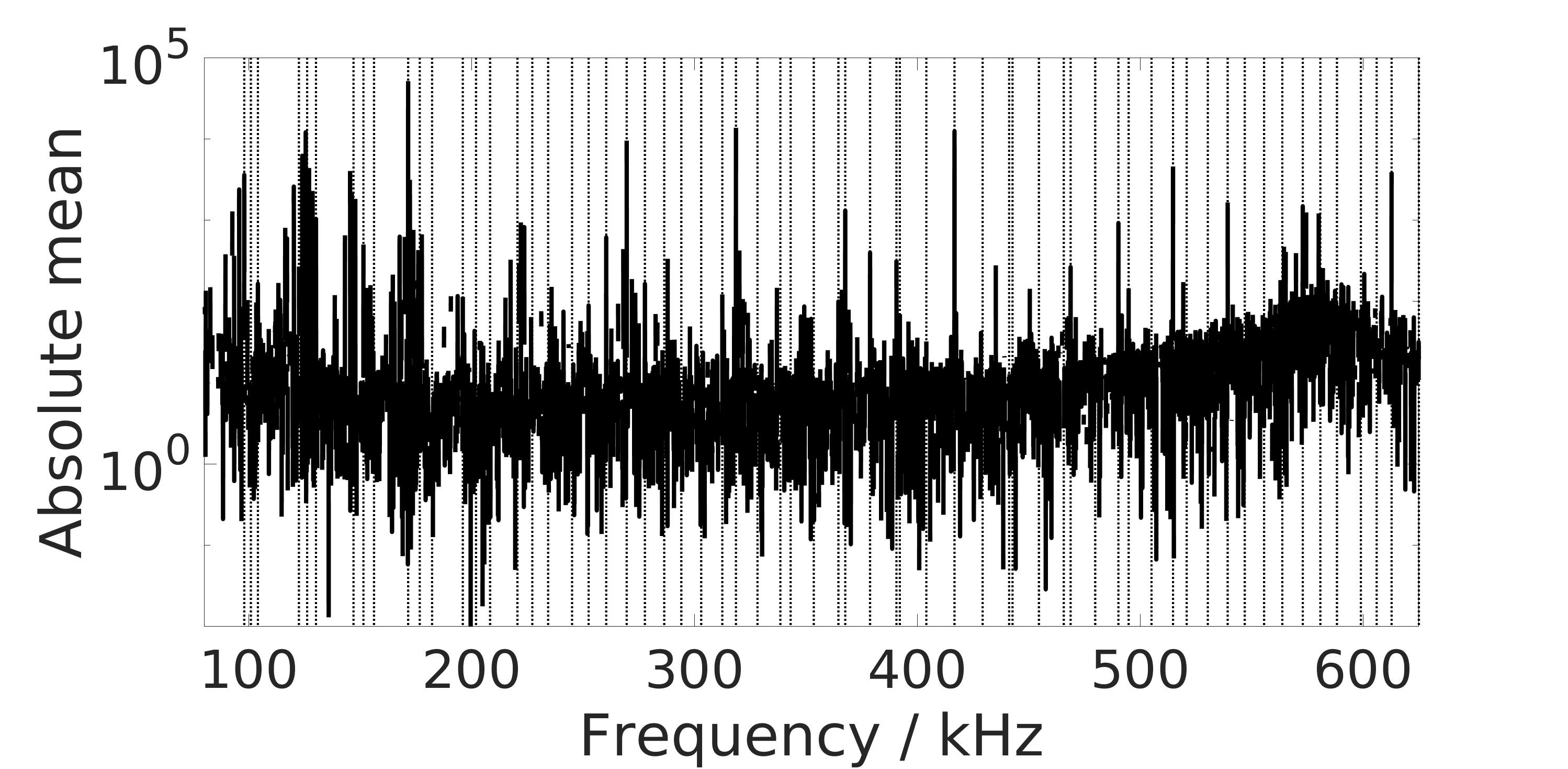} \\
\includegraphics[width=\textwidth]{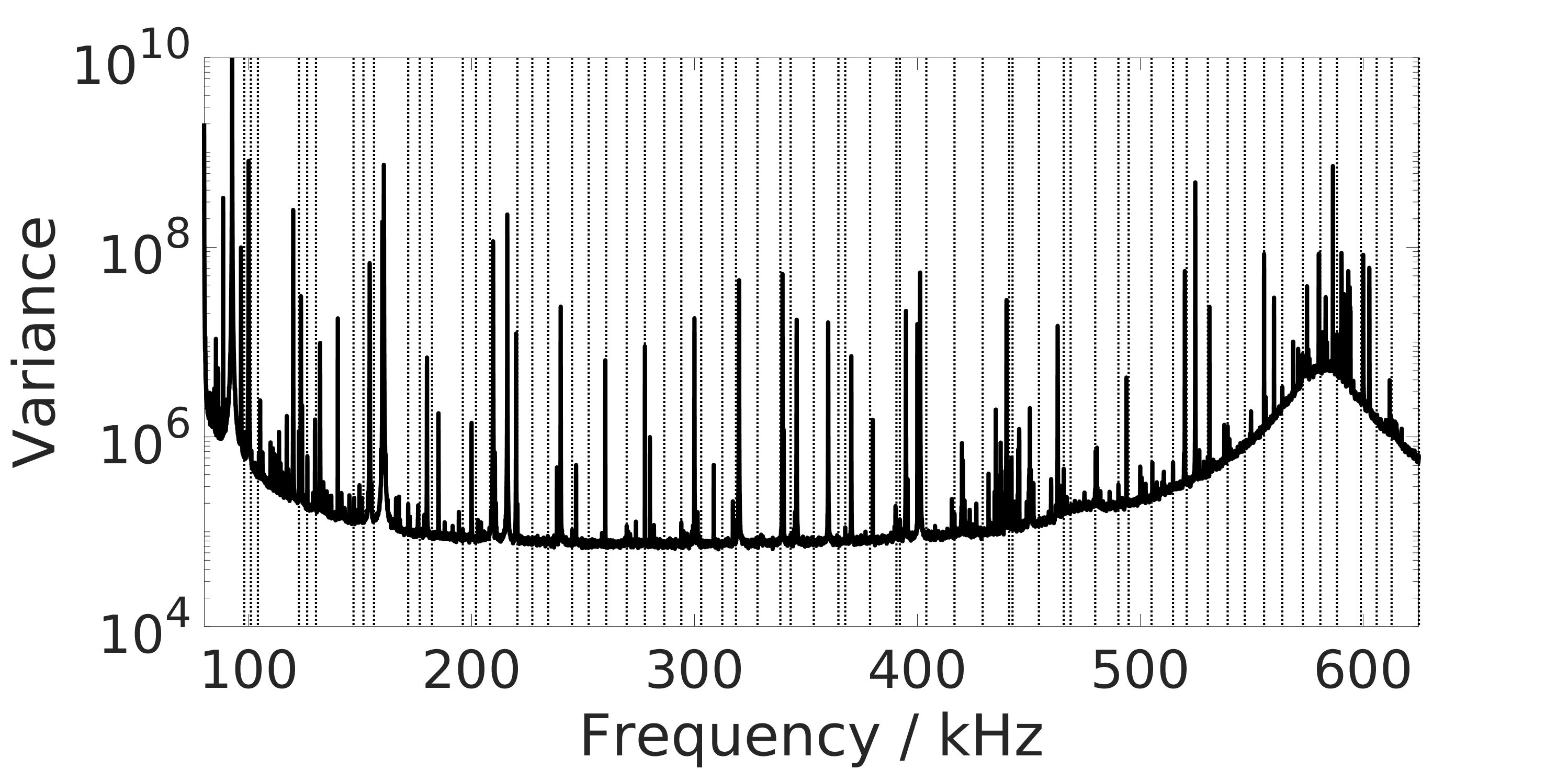} \\
 \end{minipage}
  \begin{minipage}{0.45\textwidth}
\centering
 \ \\
Imaginary \\
\includegraphics[width=\textwidth]{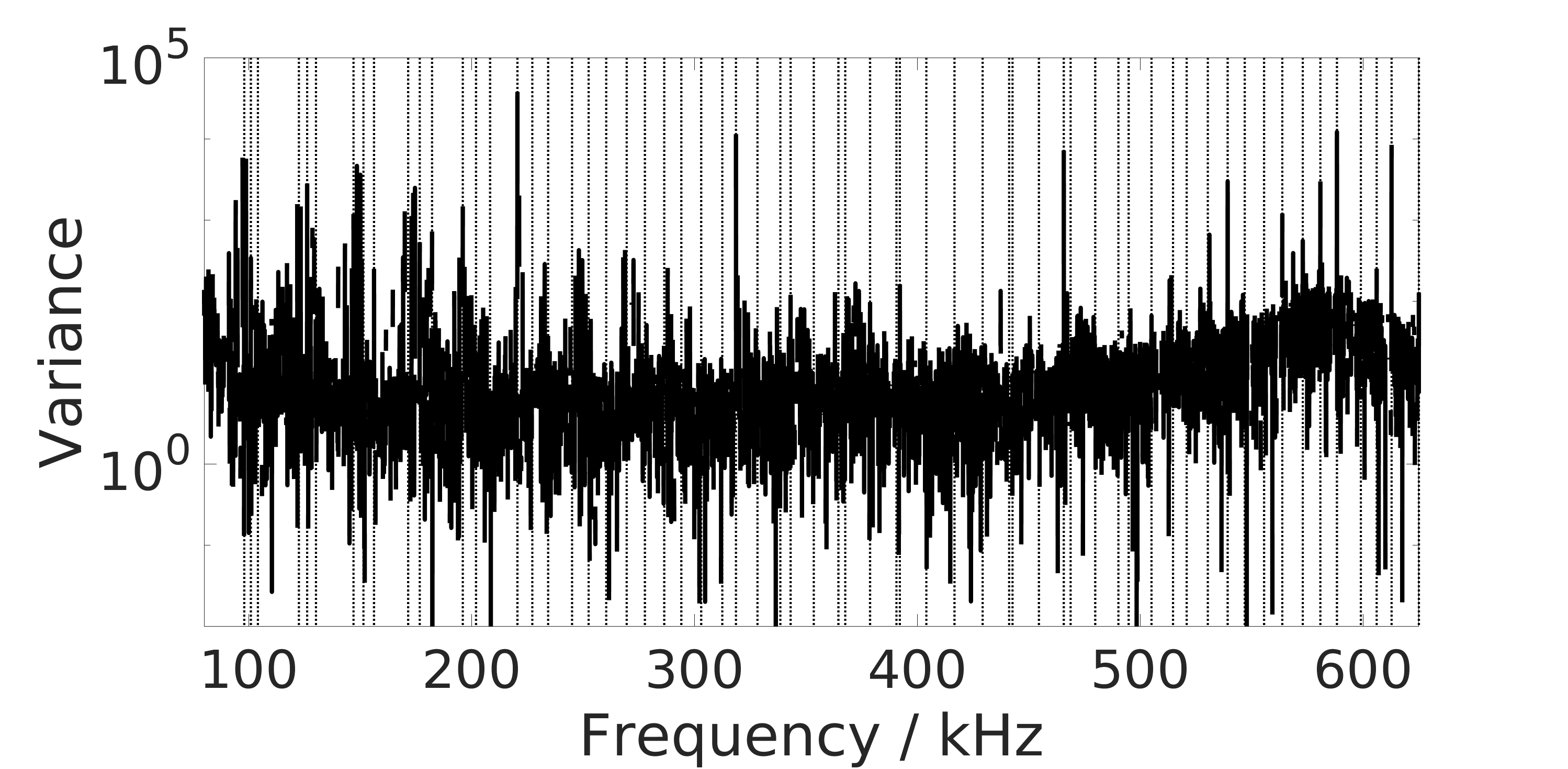} \\
\includegraphics[width=\textwidth]{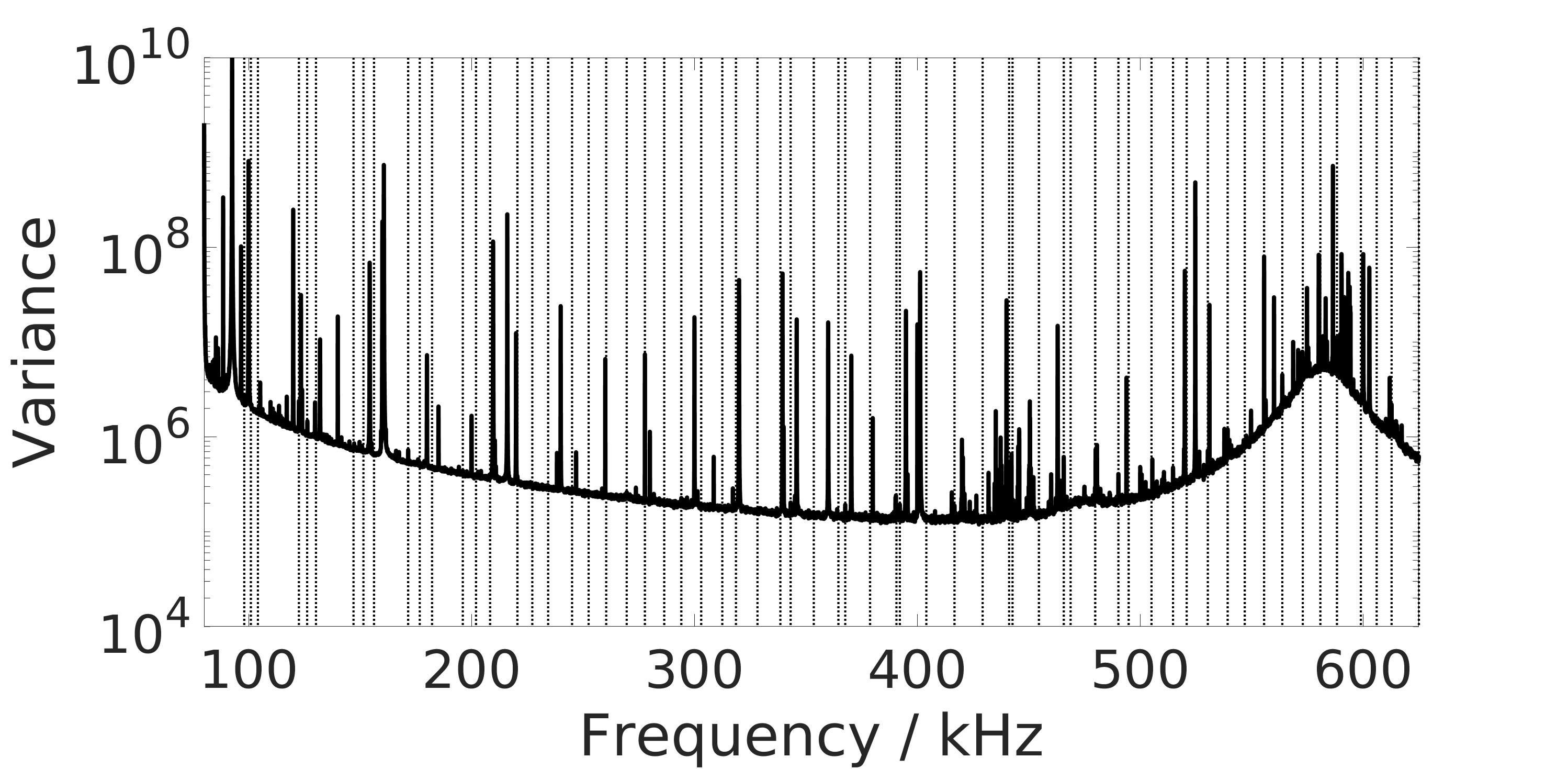} \\
 \end{minipage}

}
\caption{Mean and variance of the background measurement provided with the ``shape'' phantom from \texttt{Open MPI dataset},
computed from 1000 empty scans when using a 3D excitation in the preclinical Bruker MPI system. Visualized
for the $x$ receive coil with respect to the frequency; real part and
imaginary part. The higher harmonics of the excitations are highlighted by vertical lines for excitation frequencies $f_x$, $f_y$, $f_z$ (from left to right), and all at the bottom.
}\label{fig:cov_diag_harmonics_xcoil}
\end{figure}

\newpage
\section{Supplementary material: Standard approach iteration results - non-whitened}
\label{app:supplementsKaczmarz}

\begin{figure}[hbt!]%
\centering
\scalebox{1}{
\begin{tabular}{c|c}
PSNR & SSIM \\
 \hline
\multicolumn{2}{l}{$\tau=0$} \\
 \includegraphics[width=0.4\textwidth]{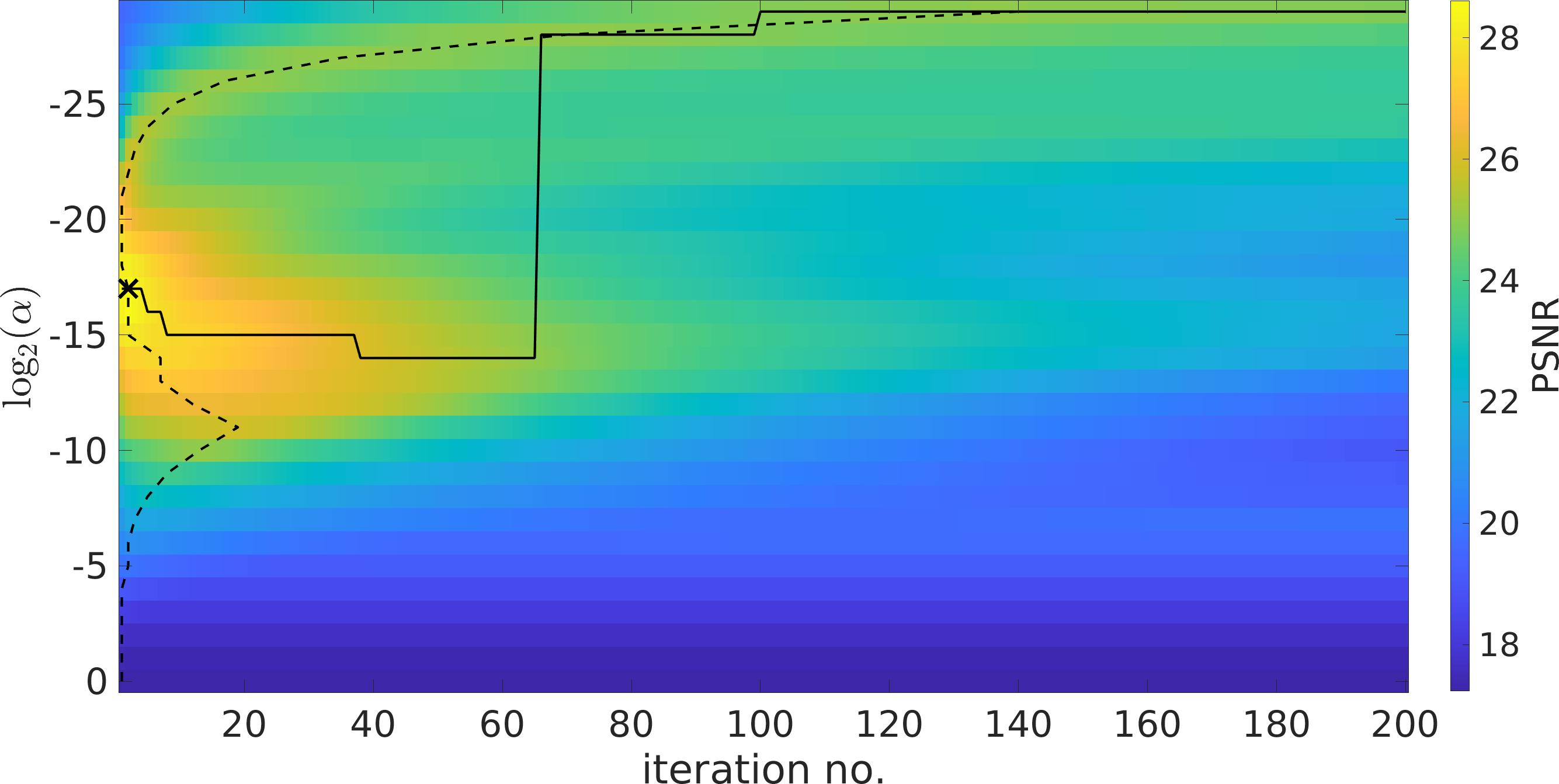} &
 \includegraphics[width=0.4\textwidth]{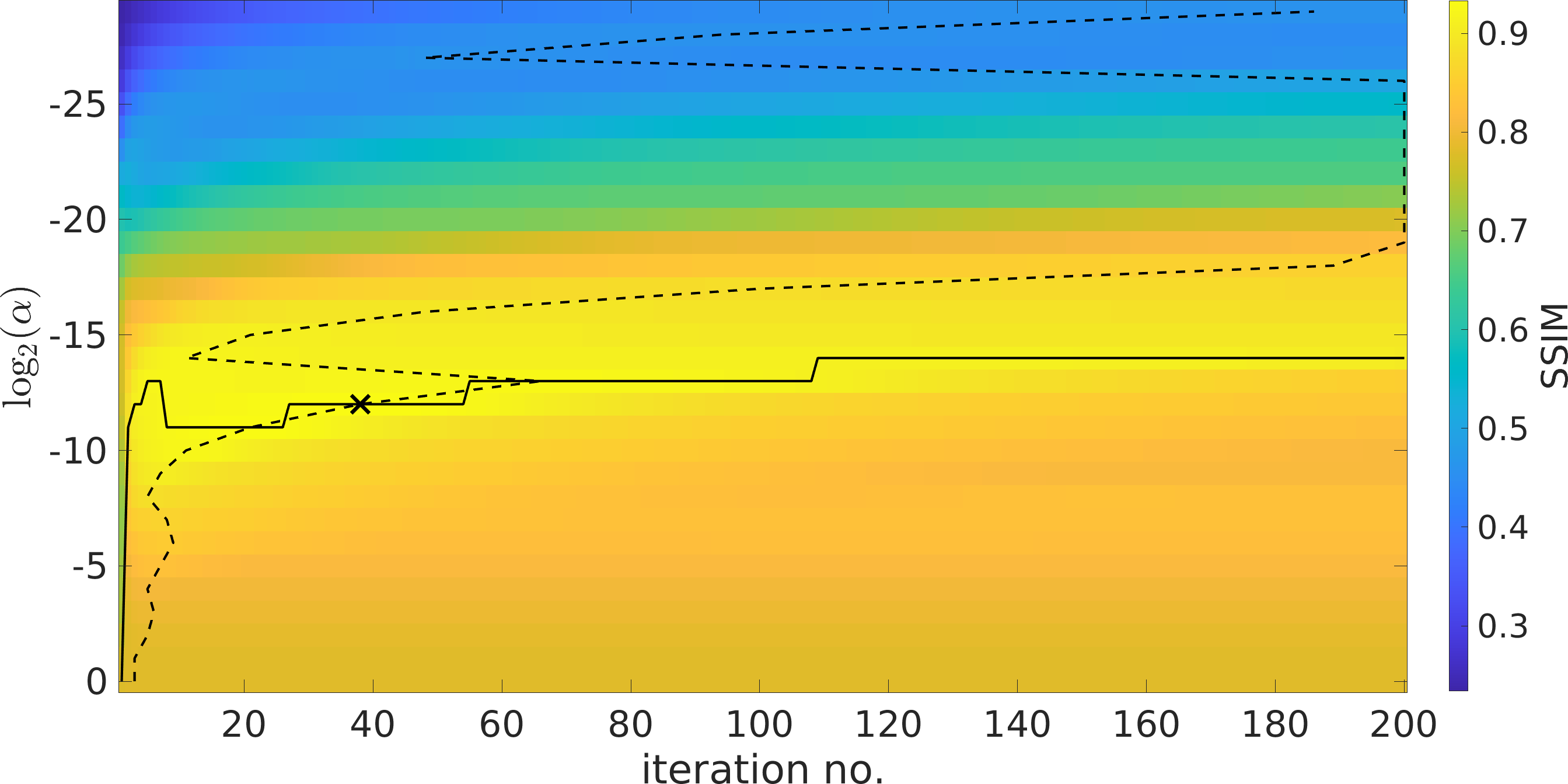} \\
 \hline
\multicolumn{2}{l}{$\tau=1$} \\
 \includegraphics[width=0.4\textwidth]{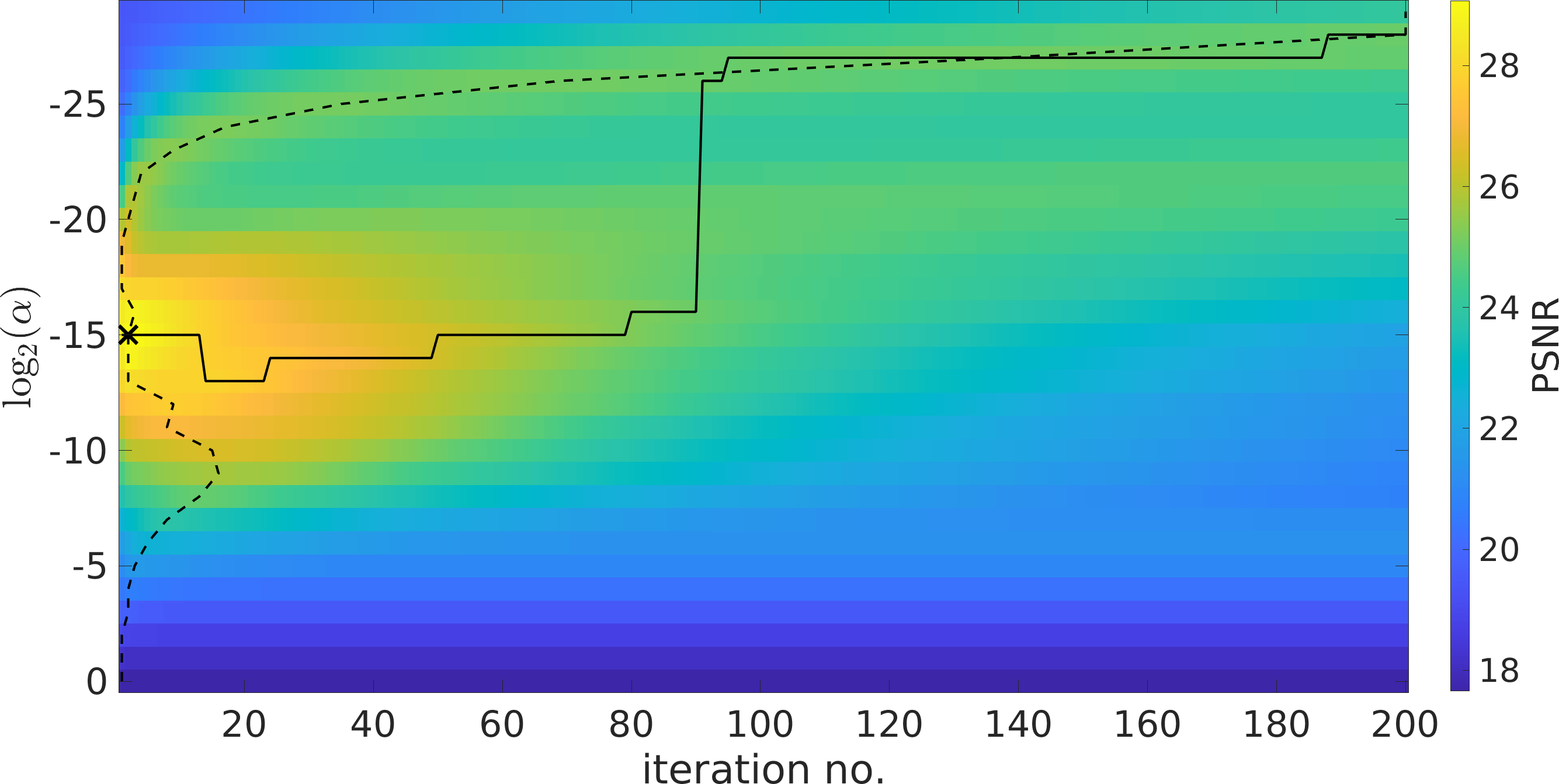} &
 \includegraphics[width=0.4\textwidth]{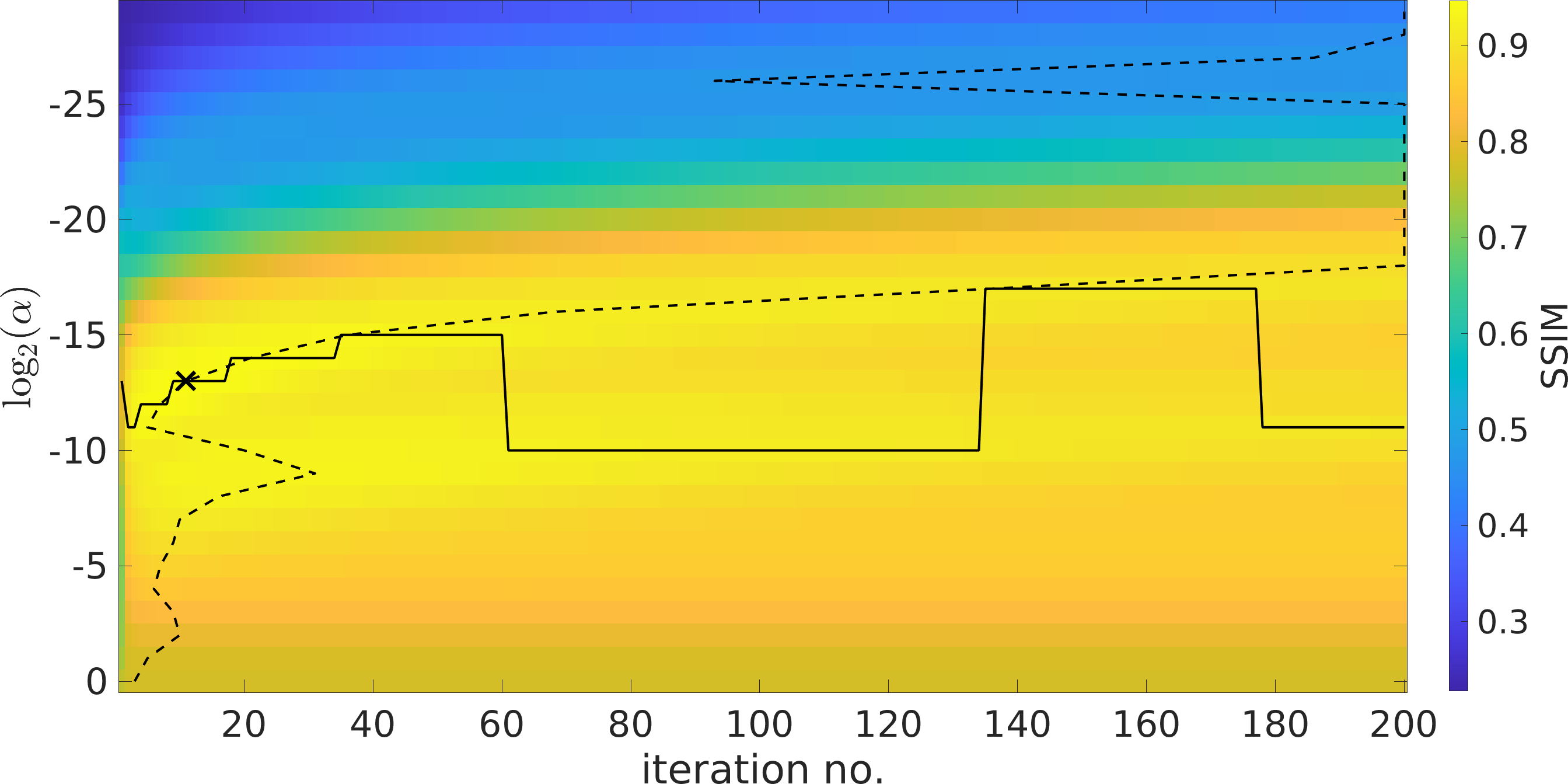} \\
  \hline
\multicolumn{2}{l}{$\tau=3$} \\
  \includegraphics[width=0.4\textwidth]{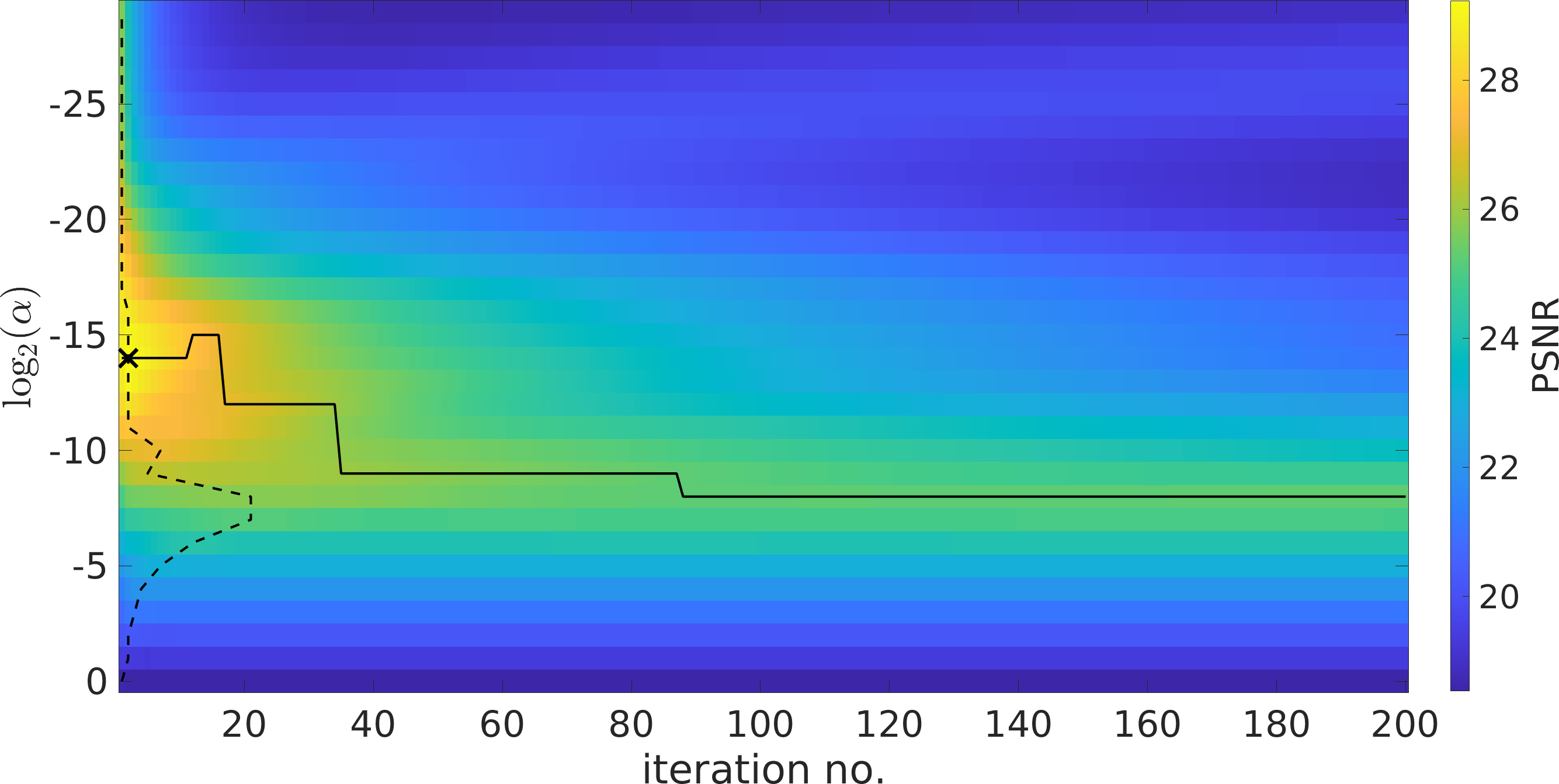} &
 \includegraphics[width=0.4\textwidth]{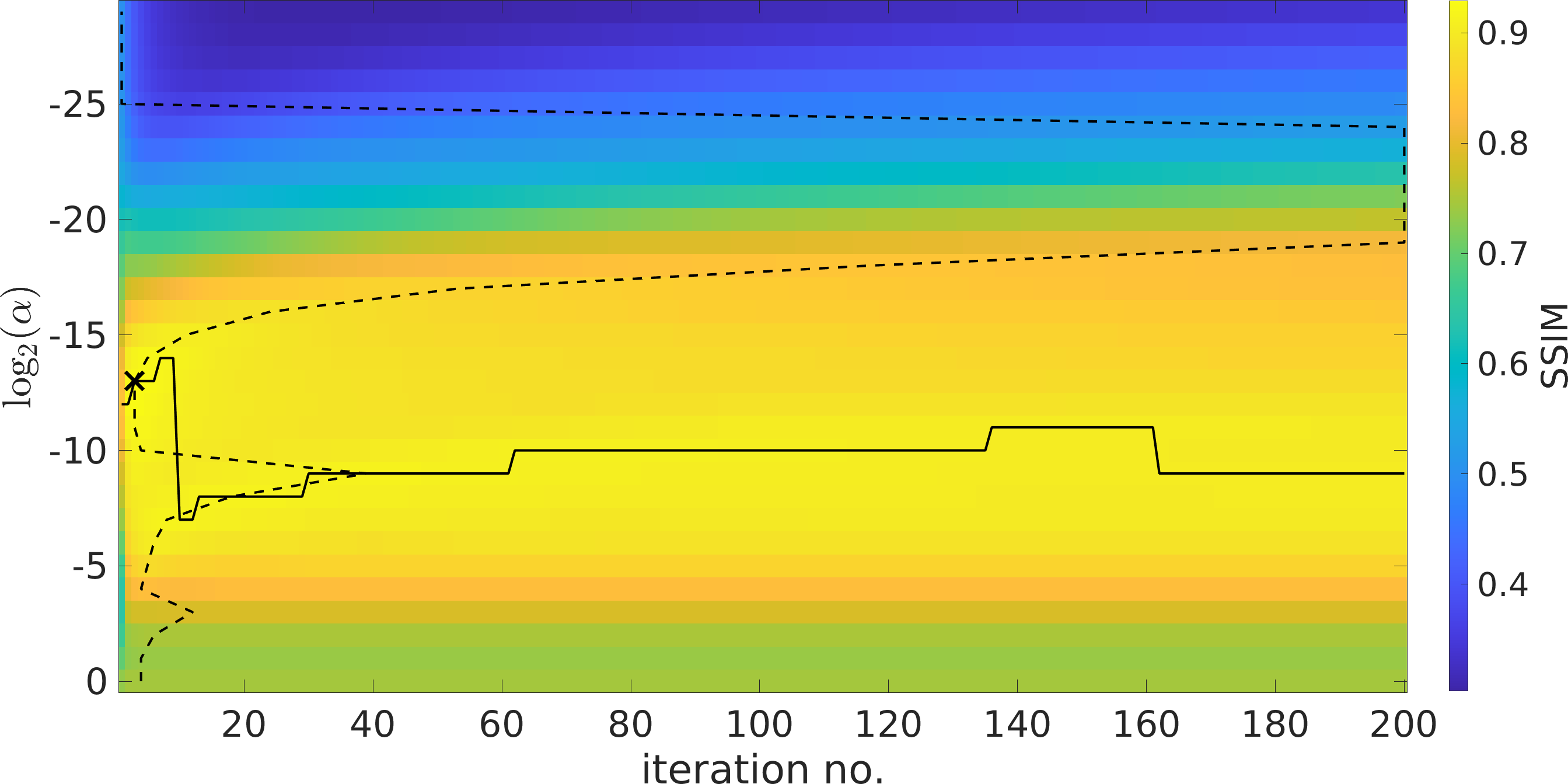} \\
  \hline
\multicolumn{2}{l}{$\tau=5$} \\
 \includegraphics[width=0.4\textwidth]{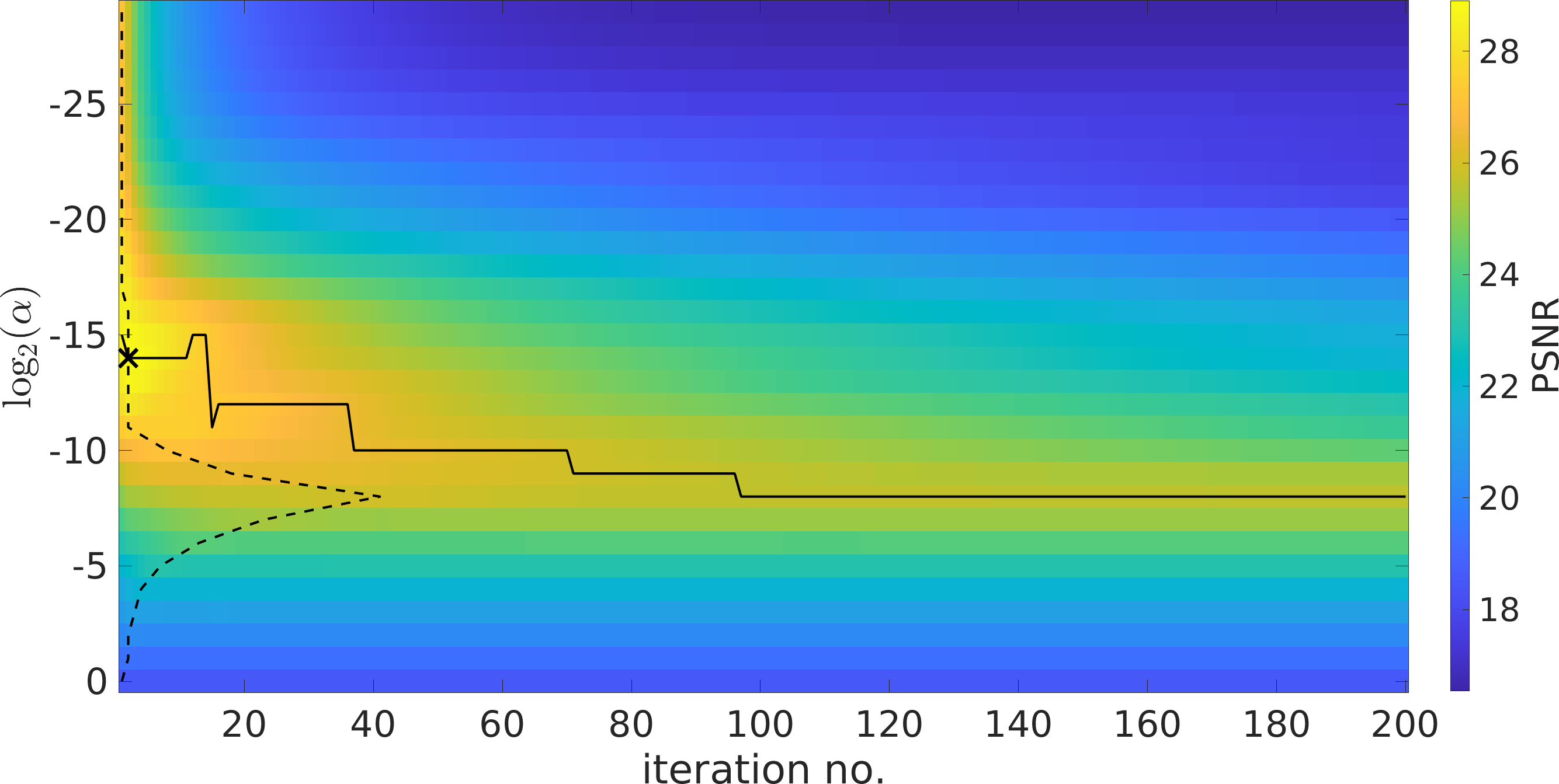} &
 \includegraphics[width=0.4\textwidth]{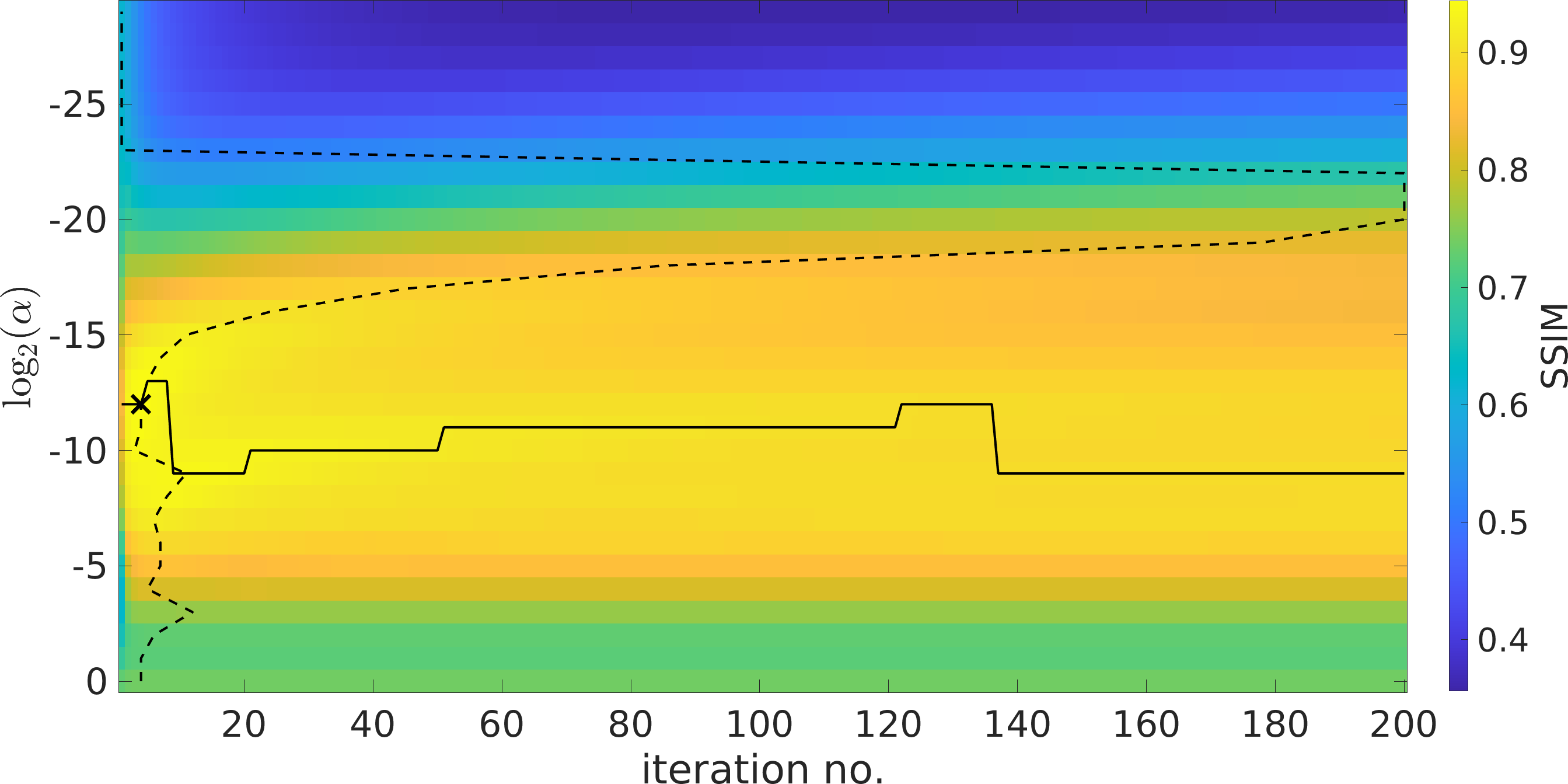} \\
\end{tabular}
}

\caption{Image quality measures with respect to $\alpha$ values and iteration number of l2-K  visualized for the ``shape'' phantom in the non-whitened case. The ``x'' marks the optimum. The solid line highlights the contour of the maximum image quality over $\alpha$ values for fixed iteration number $N$ (highlights the maximum of each column). The dashed line highlights the contour of the maximum image quality over iteration numbers $N$ for fixed $\alpha$ values (highlights the maximum of each row).}
\label{fig:shape_Kaczmarz_image_quality_nonwhitened}
\end{figure}

\begin{figure}[hbt!]%
\centering
\scalebox{1}{
\begin{tabular}{c|c}
PSNR & SSIM \\
 \hline
\multicolumn{2}{l}{$\tau=0$} \\
 \includegraphics[width=0.4\textwidth]{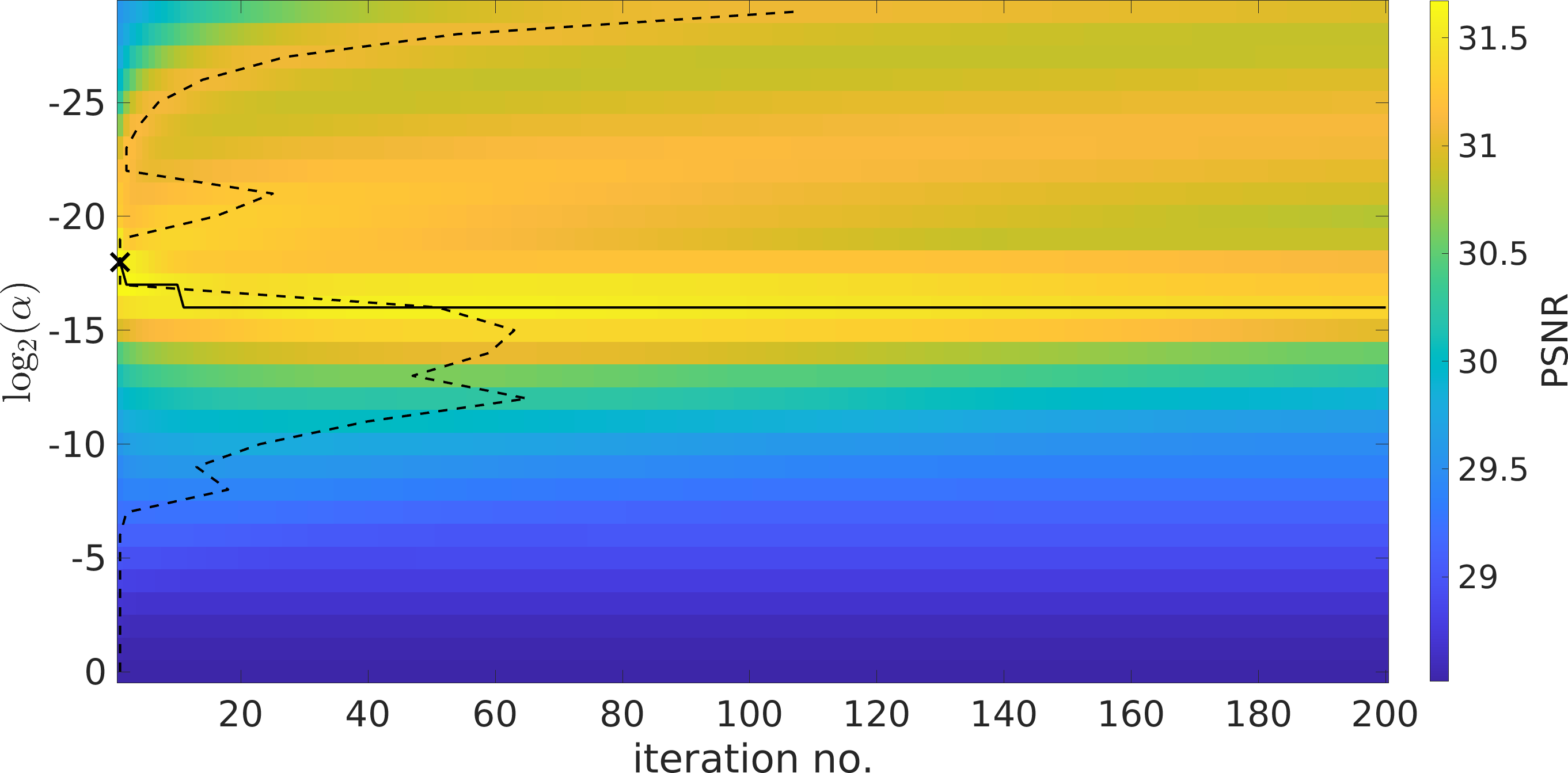} &
 \includegraphics[width=0.4\textwidth]{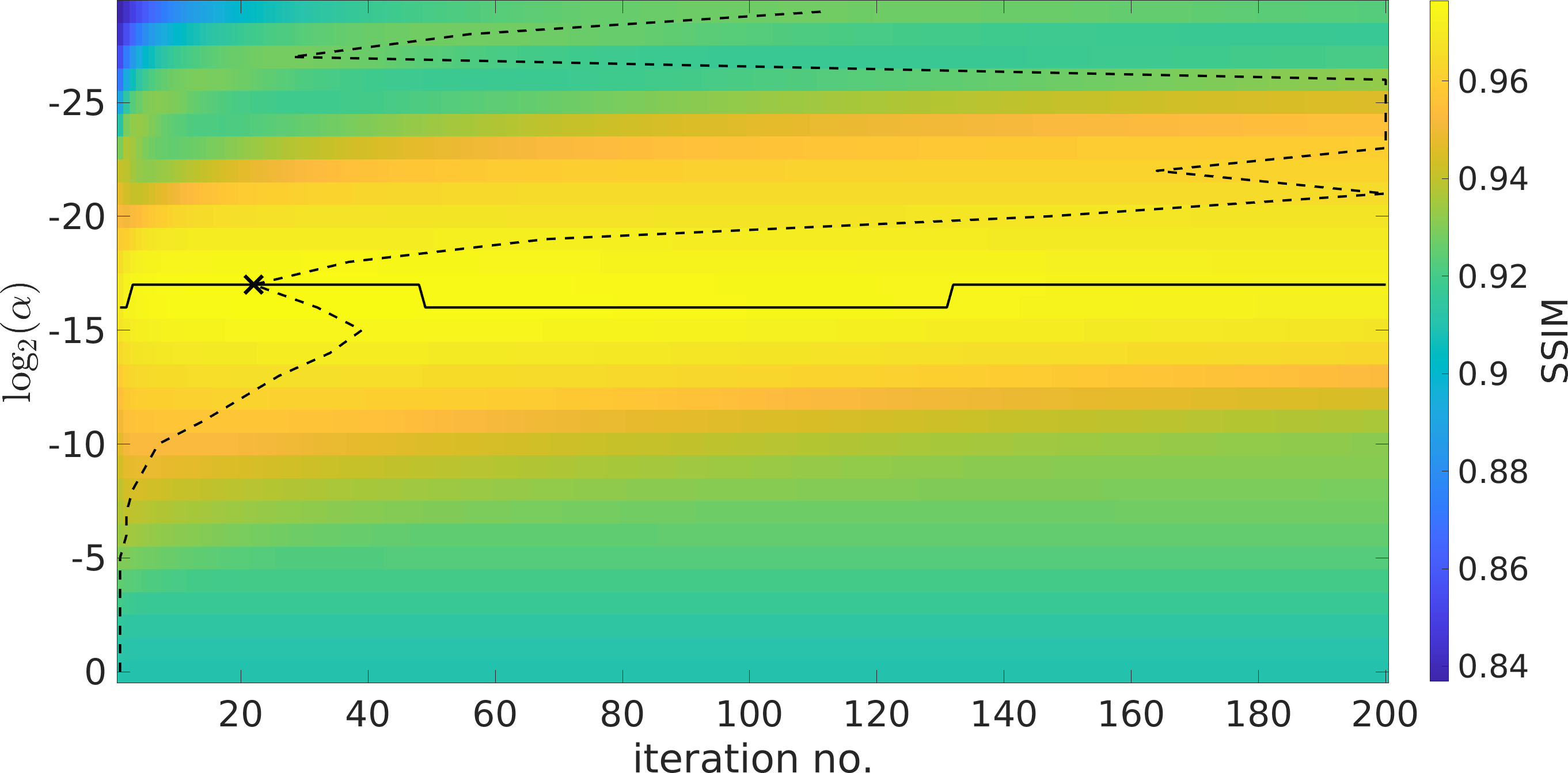} \\
 \hline
\multicolumn{2}{l}{$\tau=1$} \\
 \includegraphics[width=0.4\textwidth]{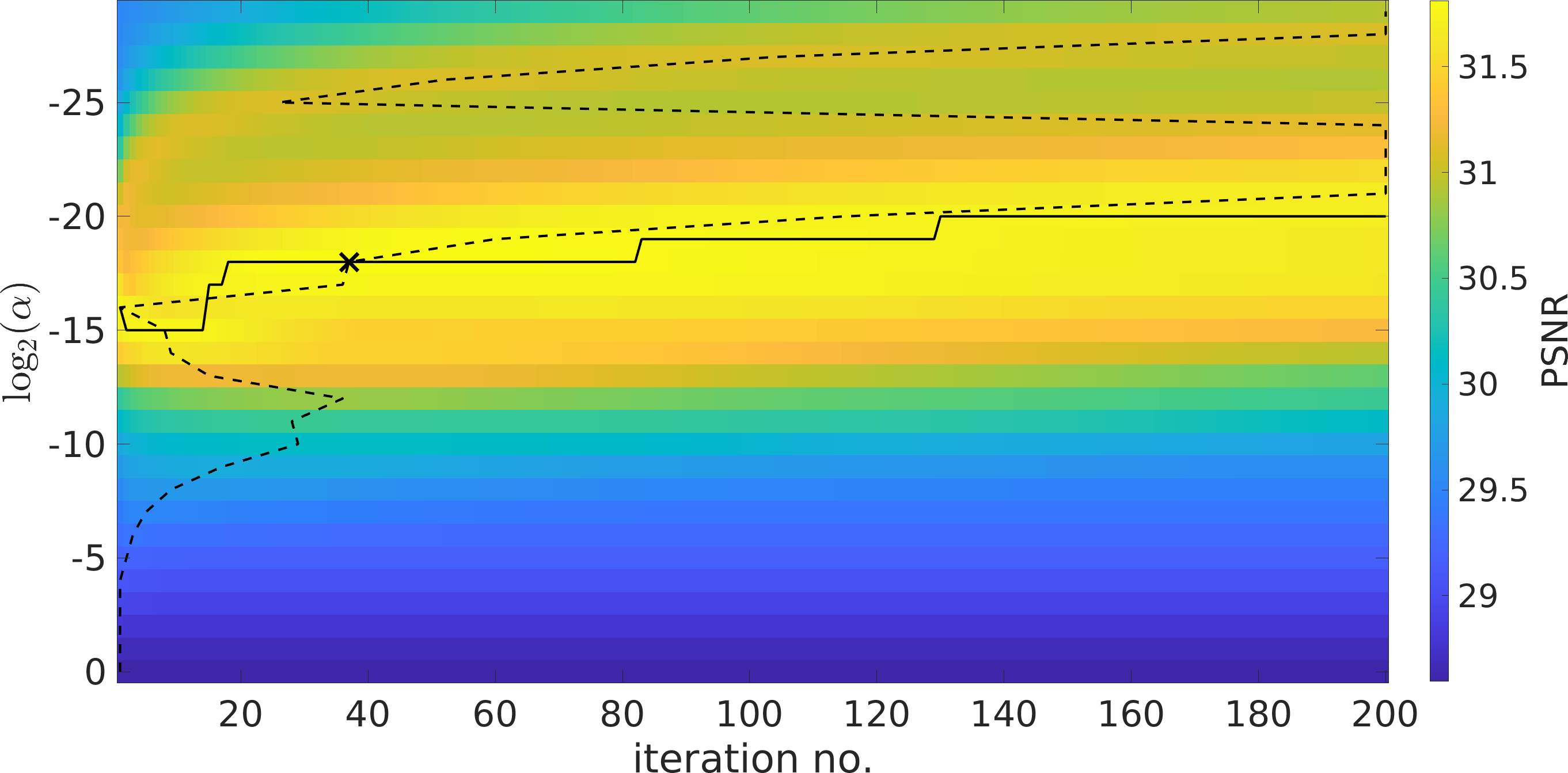} &
 \includegraphics[width=0.4\textwidth]{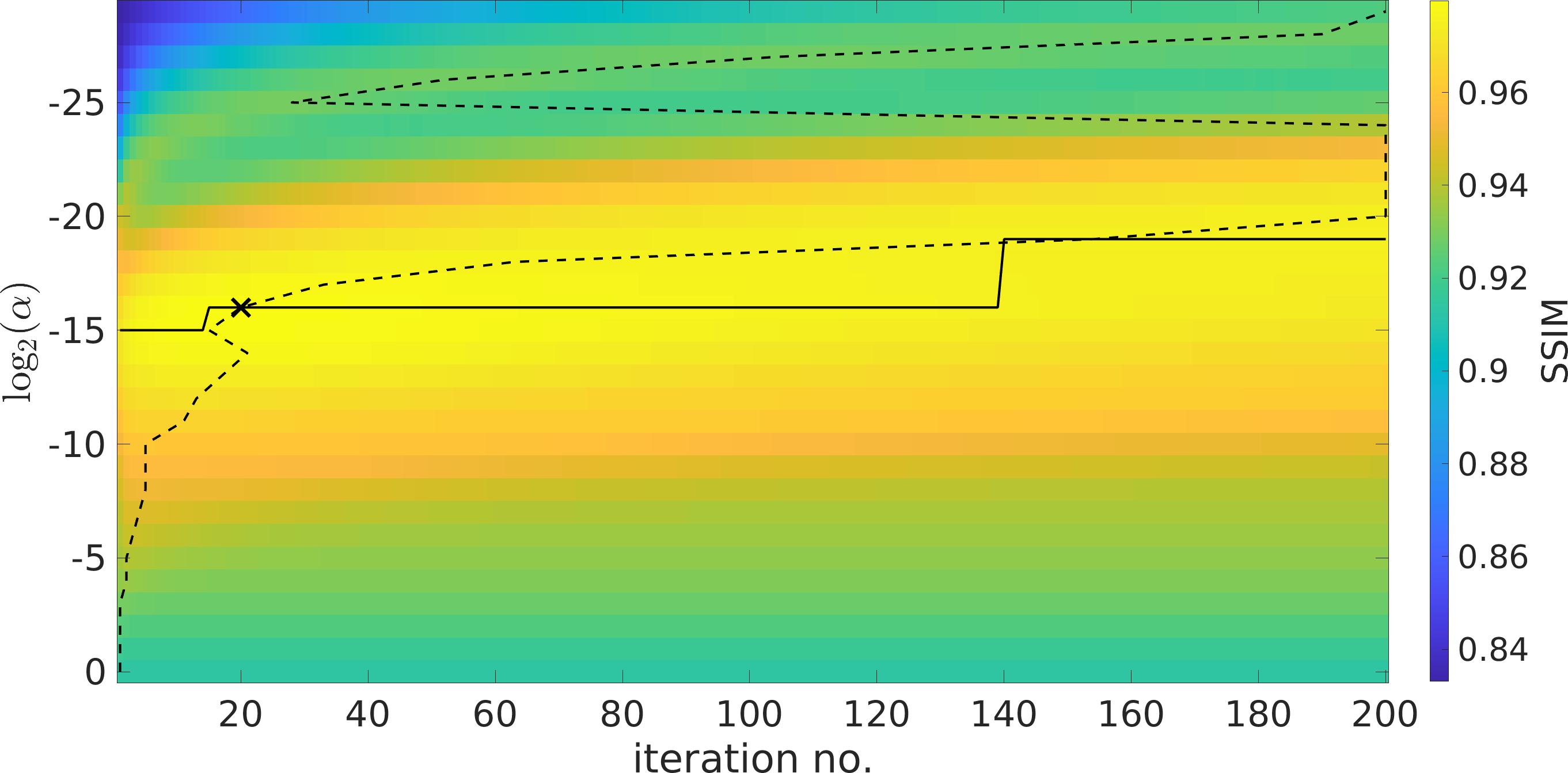} \\
  \hline
\multicolumn{2}{l}{$\tau=3$} \\
  \includegraphics[width=0.4\textwidth]{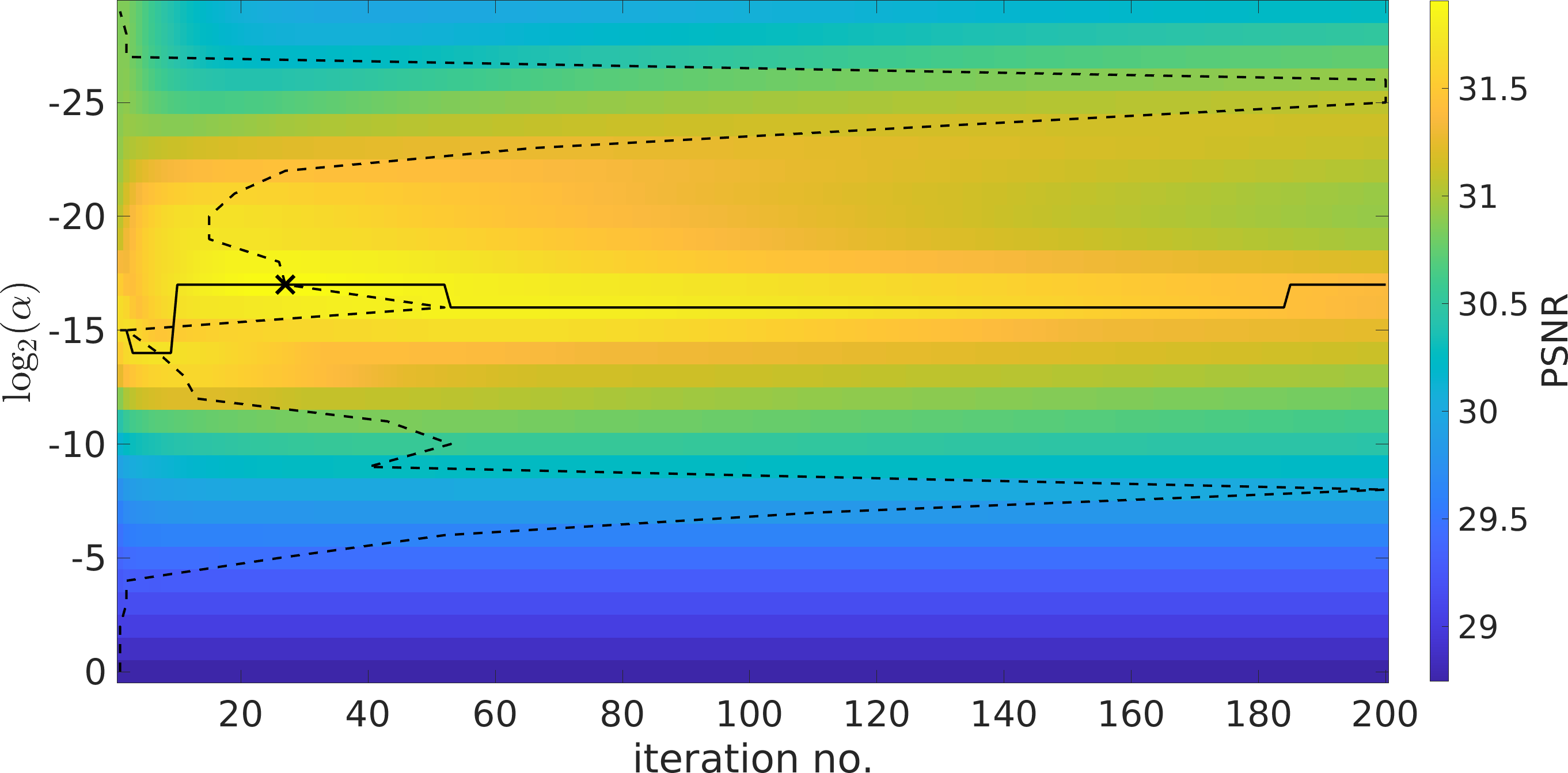} &
 \includegraphics[width=0.4\textwidth]{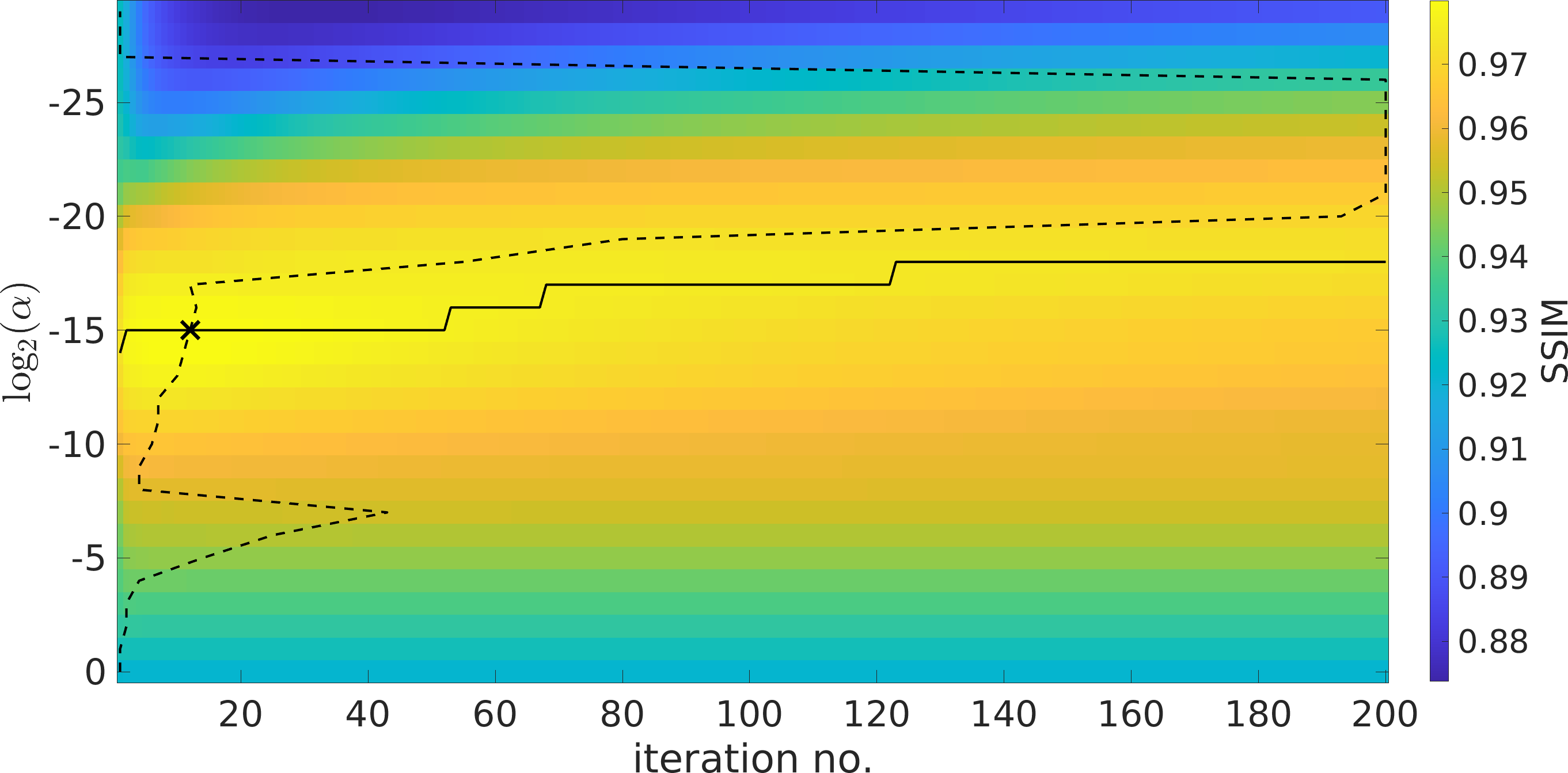} \\
  \hline
\multicolumn{2}{l}{$\tau=5$} \\
 \includegraphics[width=0.4\textwidth]{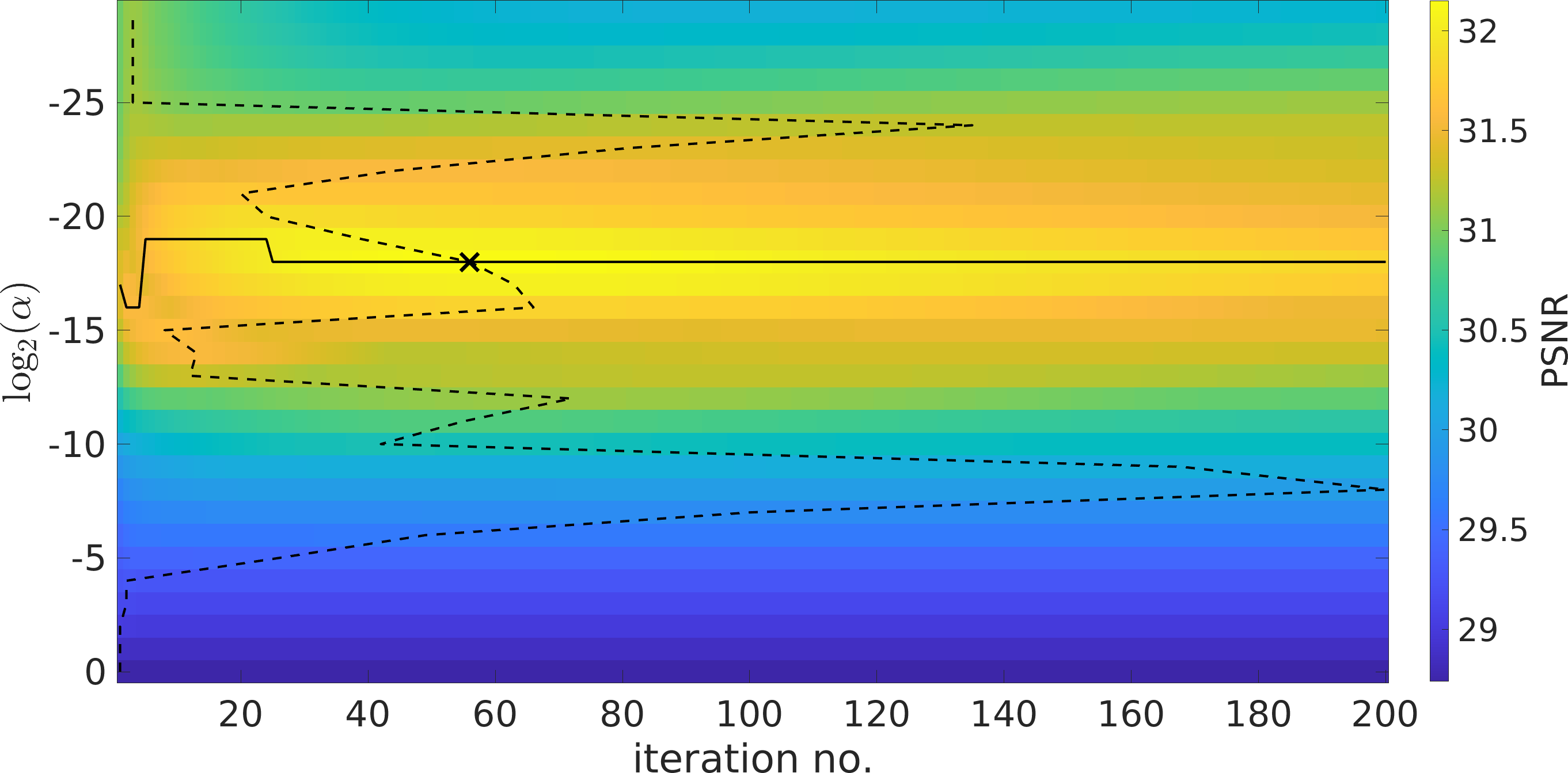} &
 \includegraphics[width=0.4\textwidth]{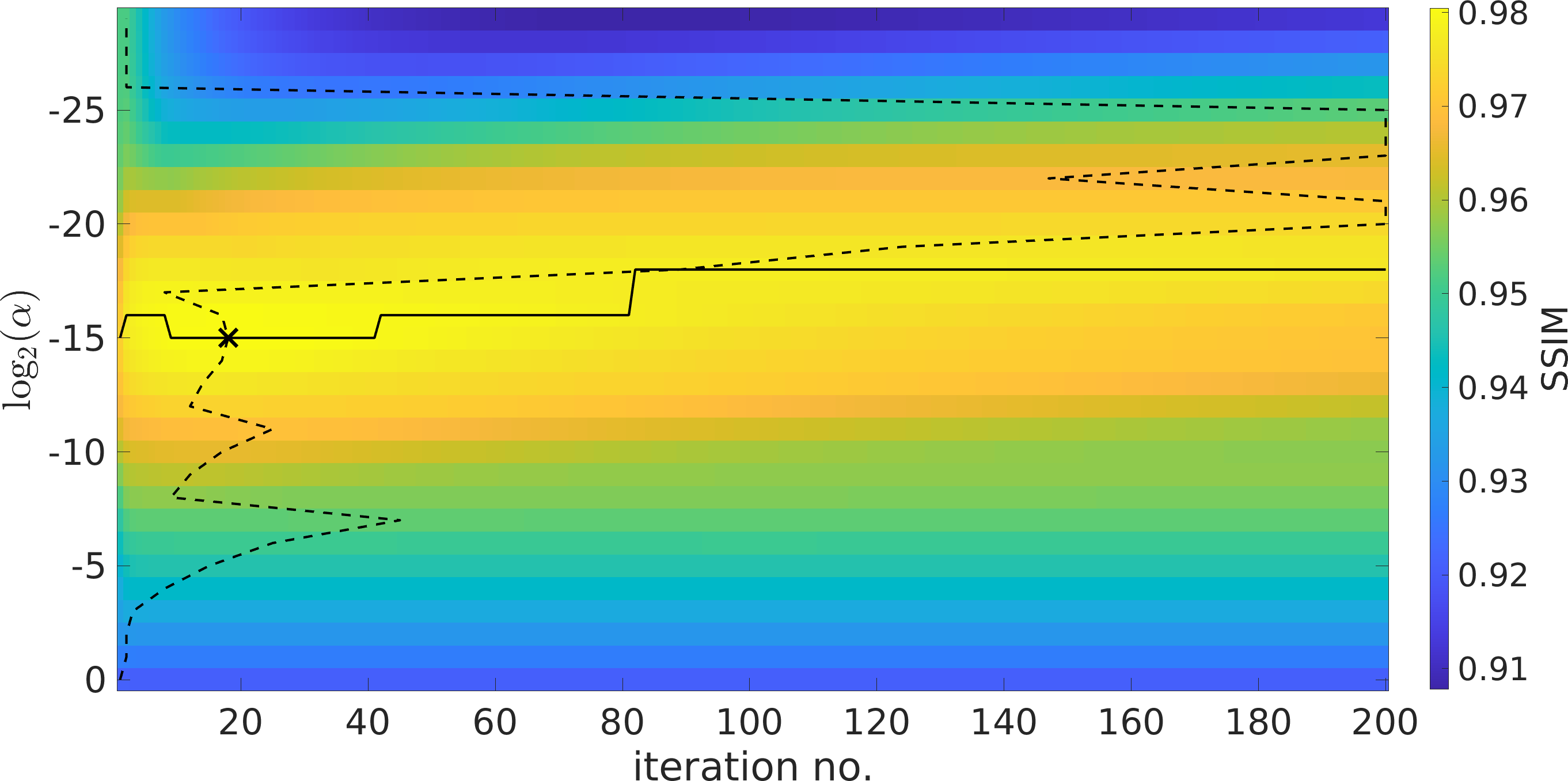} \\
\end{tabular}
}

\caption{Image quality measures with respect to $\alpha$ values and iteration number of l2-K  visualized for the ``resolution'' phantom in the non-whitened case. The ``x'' marks the optimum. The solid line highlights the contour of the maximum image quality over $\alpha$ values for fixed iteration number $N$ (highlights the maximum of each column). The dashed line highlights the contour of the maximum image quality over iteration numbers $N$ for fixed $\alpha$ values (highlights the maximum of each row).}
\label{fig:resolution_Kaczmarz_image_quality_nonwhitened}
\end{figure}

\newpage
\section{Supplementary material: Method comparison - inverted colormap}
\label{app:supplements_inverted_colormap}

\begin{figure}[hbt!]
\centering
\scalebox{0.85}{
\begin{tabular}{ccc|ccc}
\multicolumn{3}{c|}{non-whitened} & \multicolumn{3}{c}{whitened} \\
\hline
l1-L & l2-L & l2-K & l1-L & l2-L & l2-K \\
\hline
\multicolumn{6}{l}{$\tau=0$} \\
 \includegraphics[height=3.4cm]{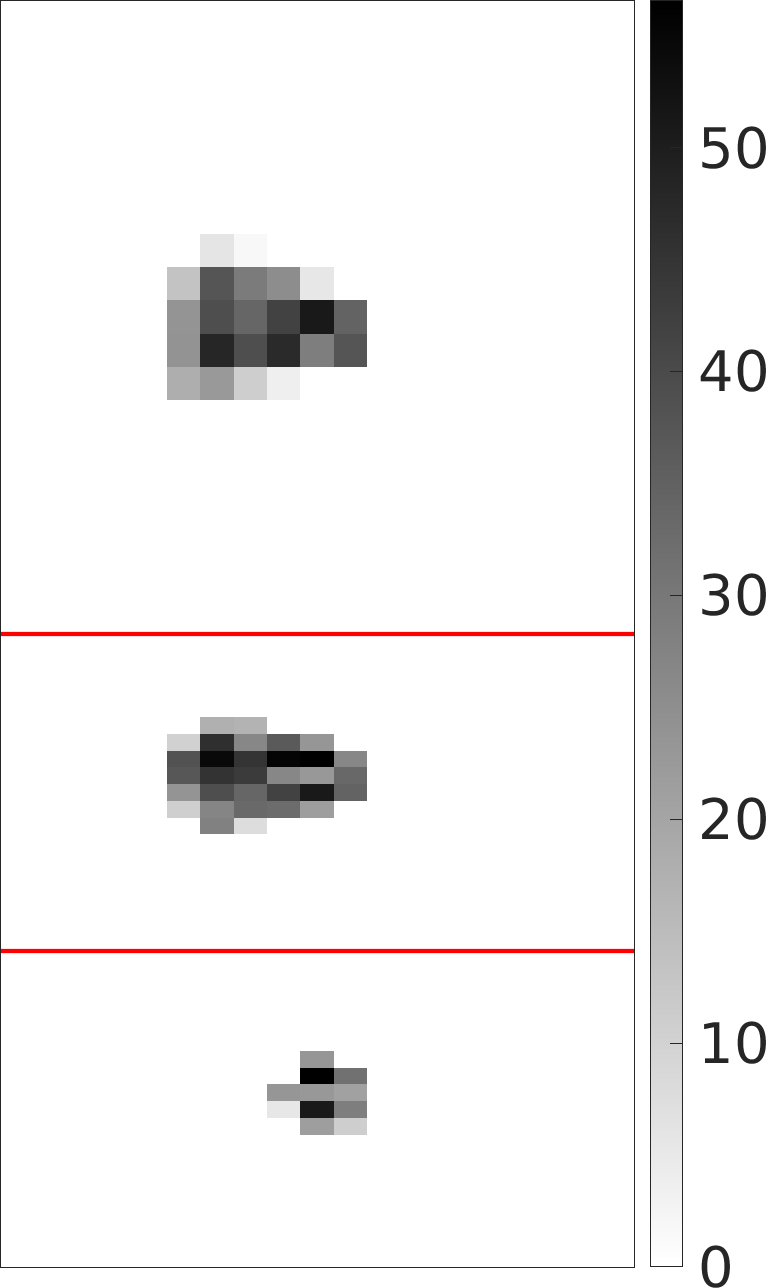}&
 \includegraphics[height=3.4cm]{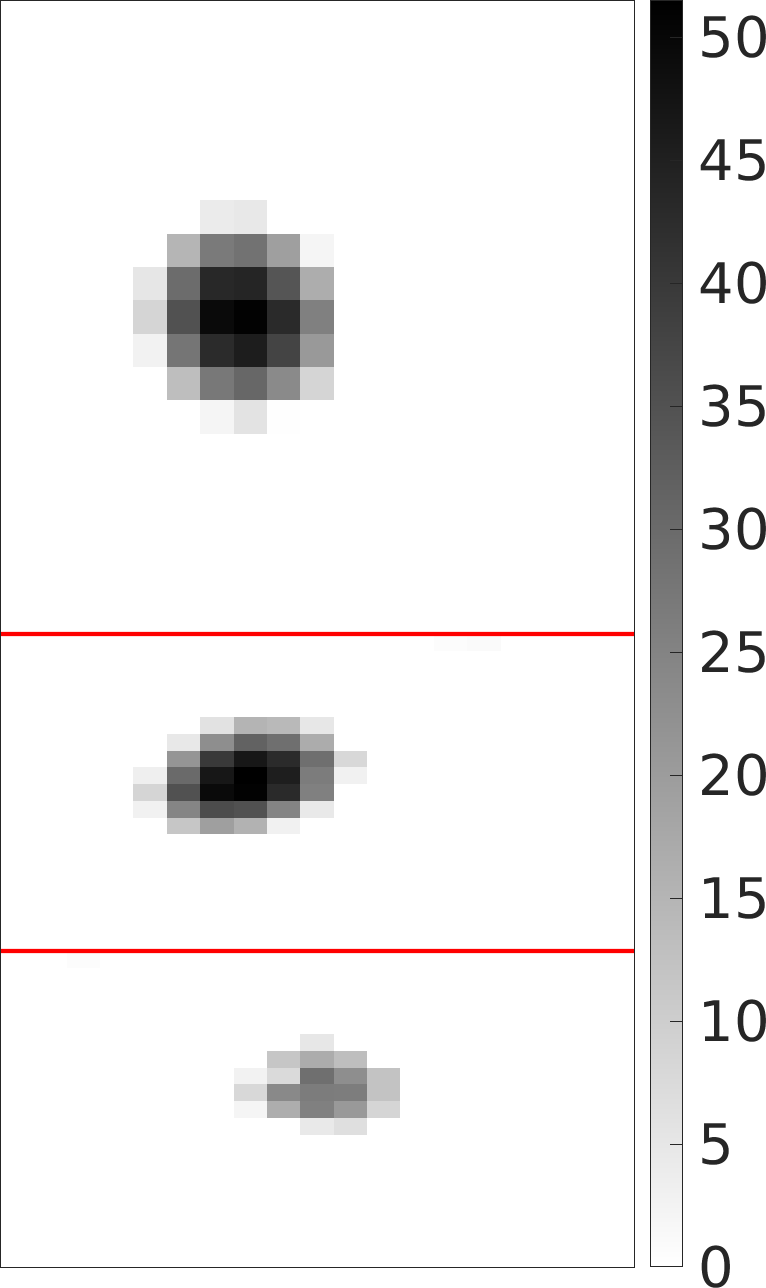}&
 \includegraphics[height=3.4cm]{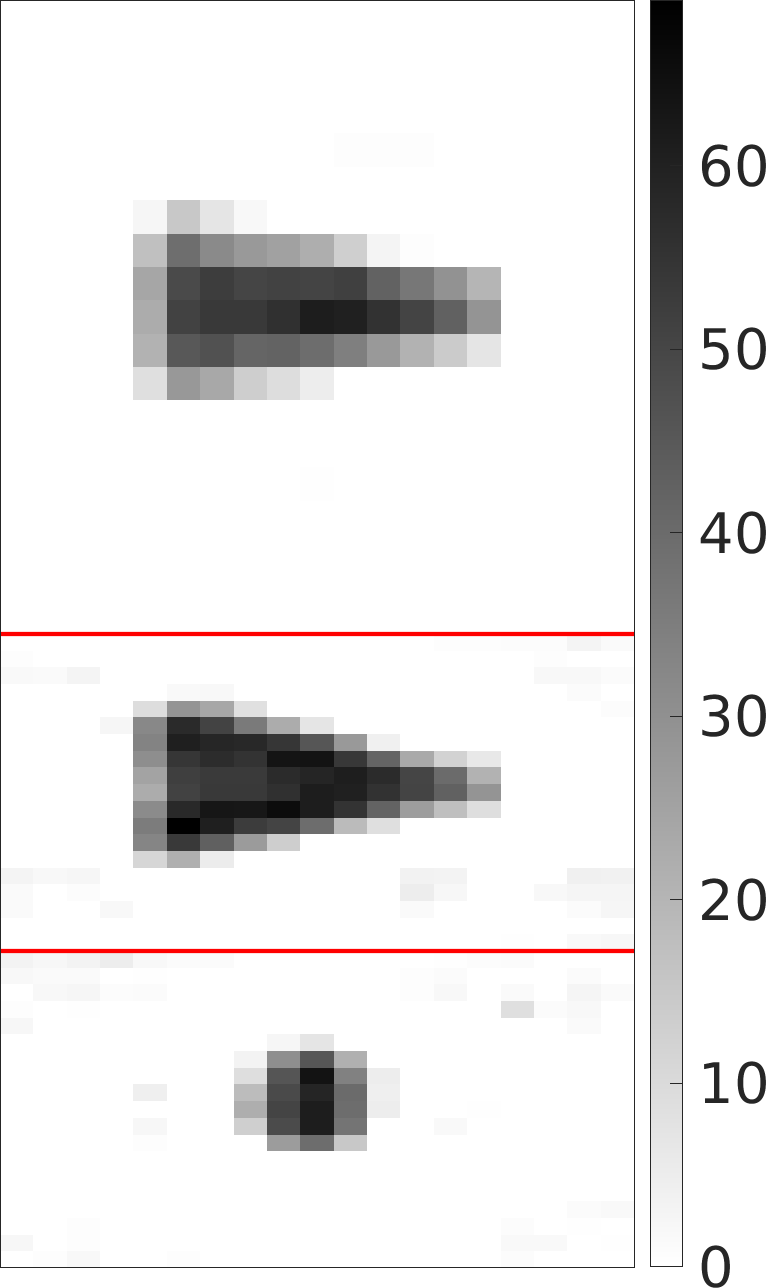}&
 \includegraphics[height=3.4cm]{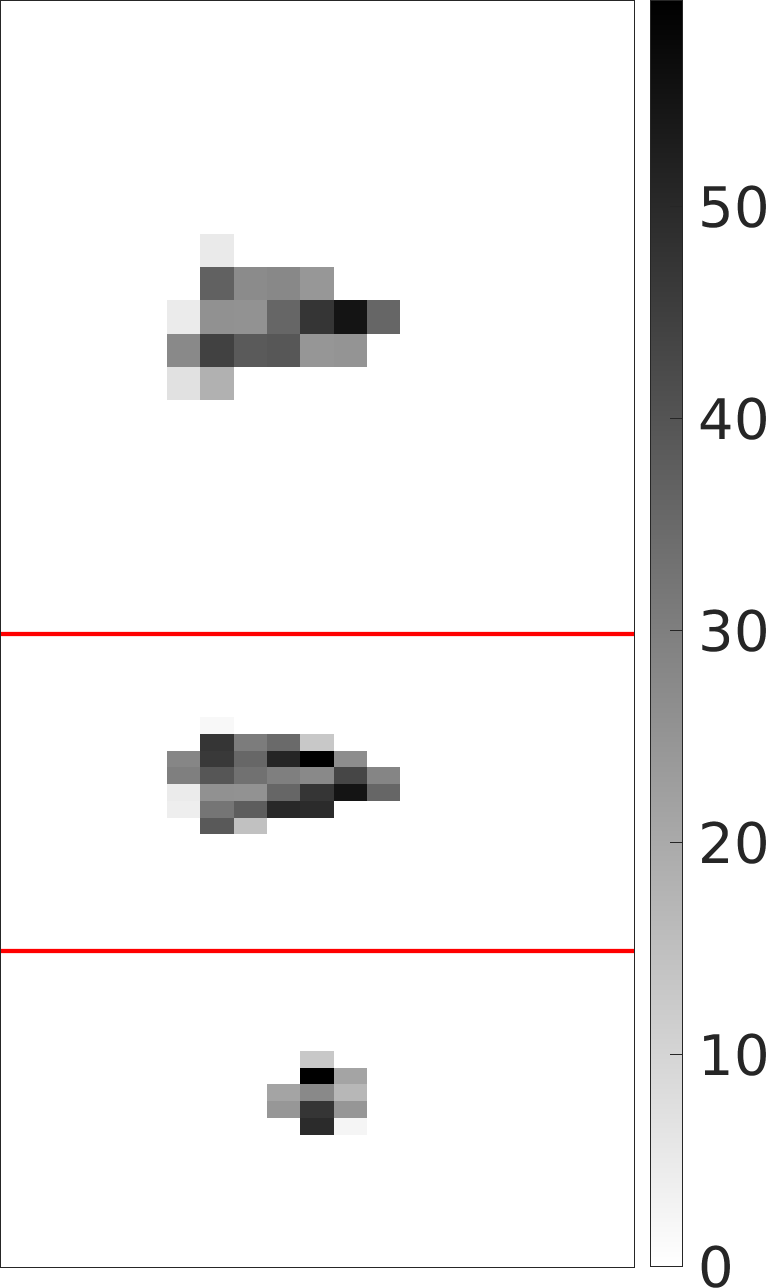}&
 \includegraphics[height=3.4cm]{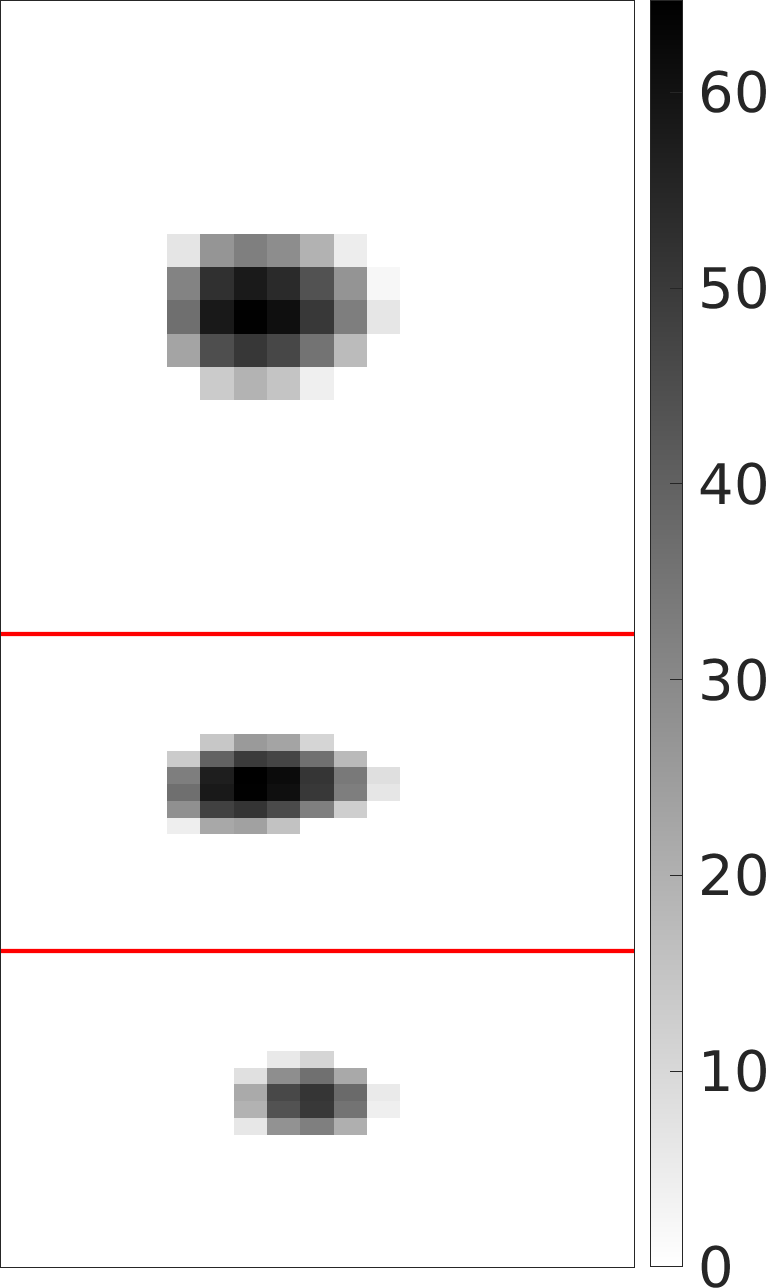}&
 \includegraphics[height=3.4cm]{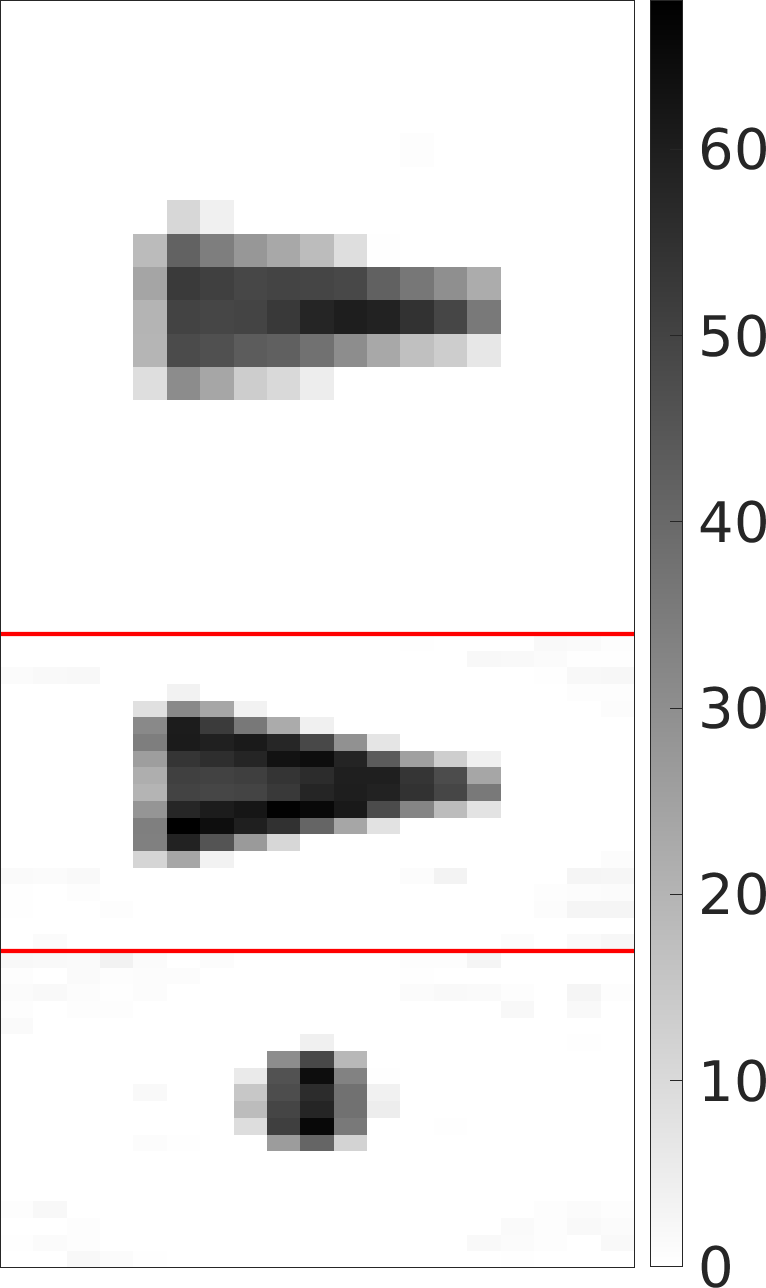}\\
\hline
\multicolumn{6}{l}{$\tau=1$} \\
 \includegraphics[height=3.4cm]{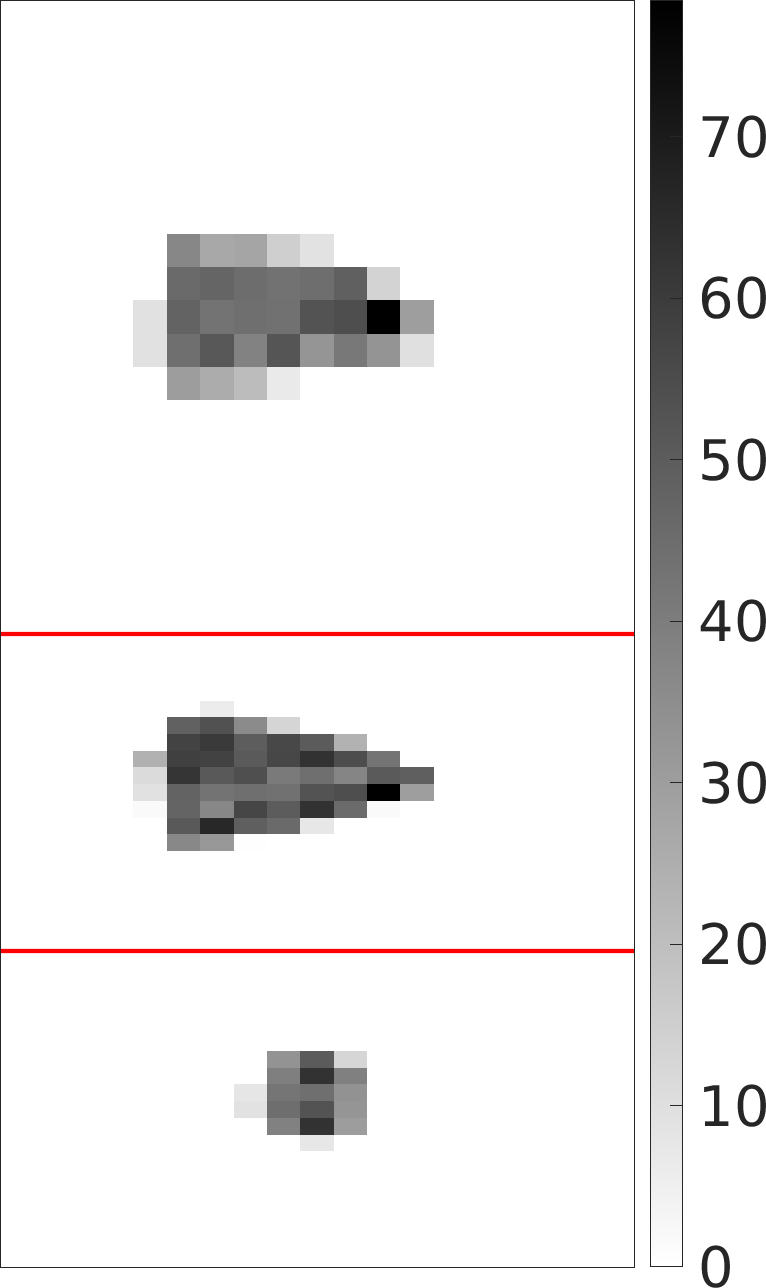}&
 \includegraphics[height=3.4cm]{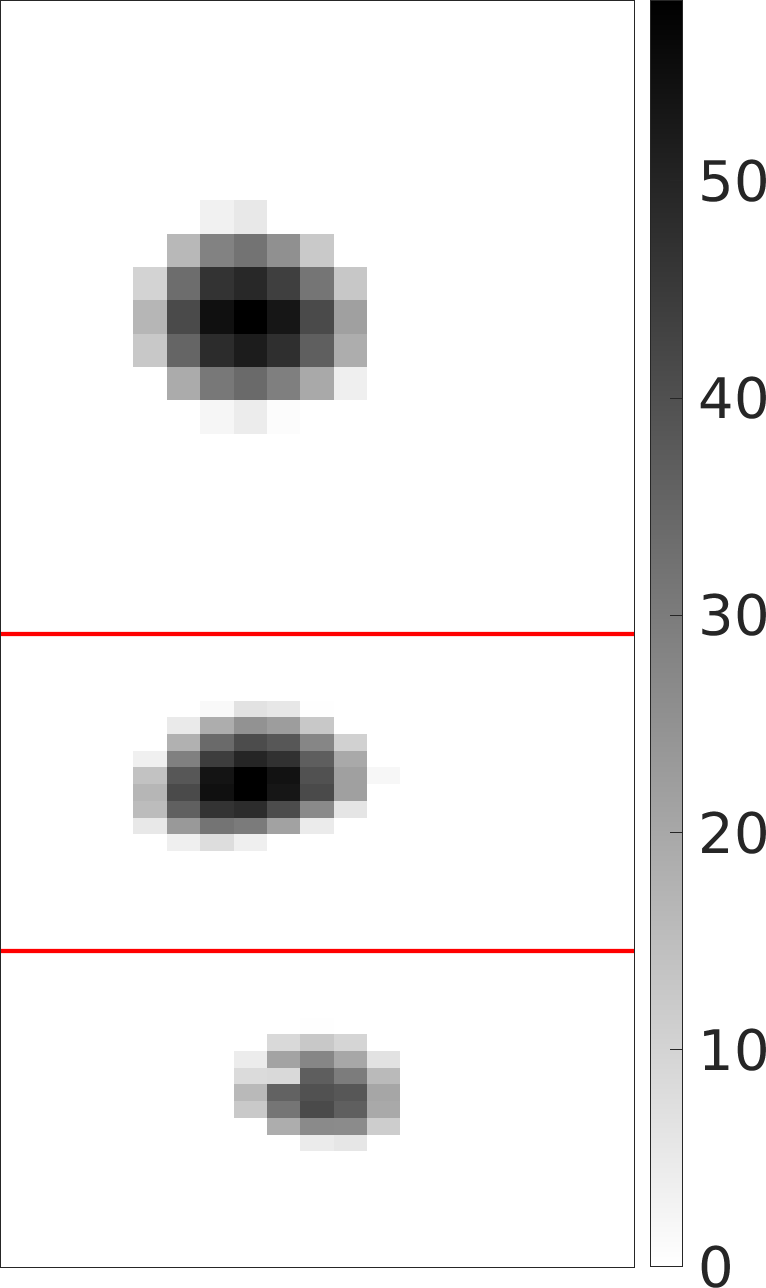}&
 \includegraphics[height=3.4cm]{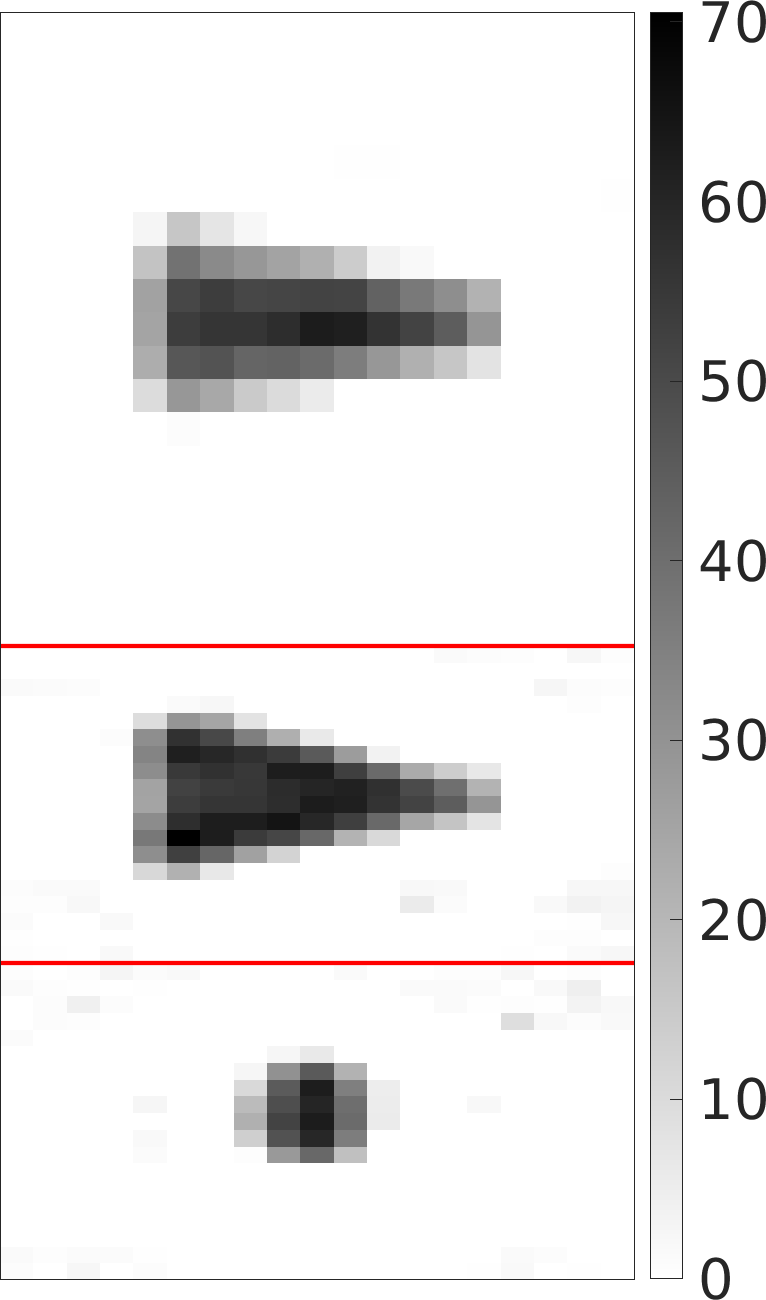}&
 \includegraphics[height=3.4cm]{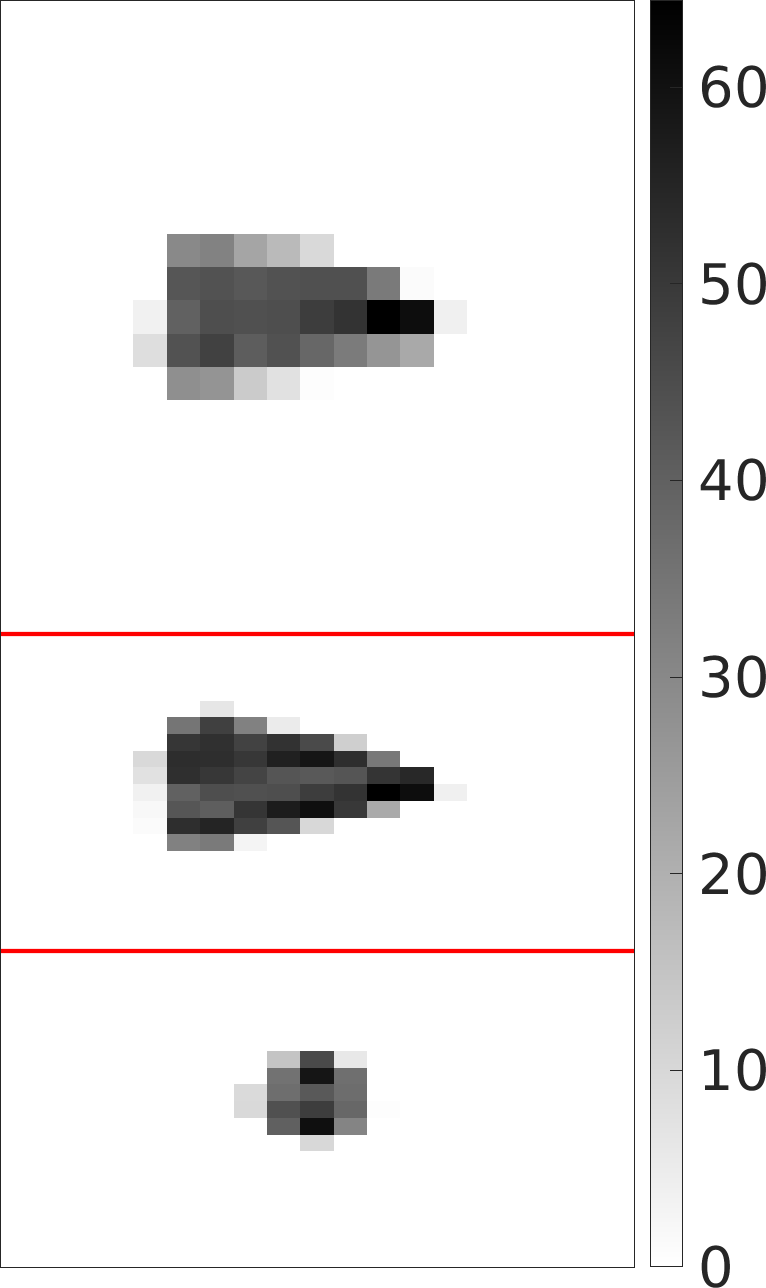}&
 \includegraphics[height=3.4cm]{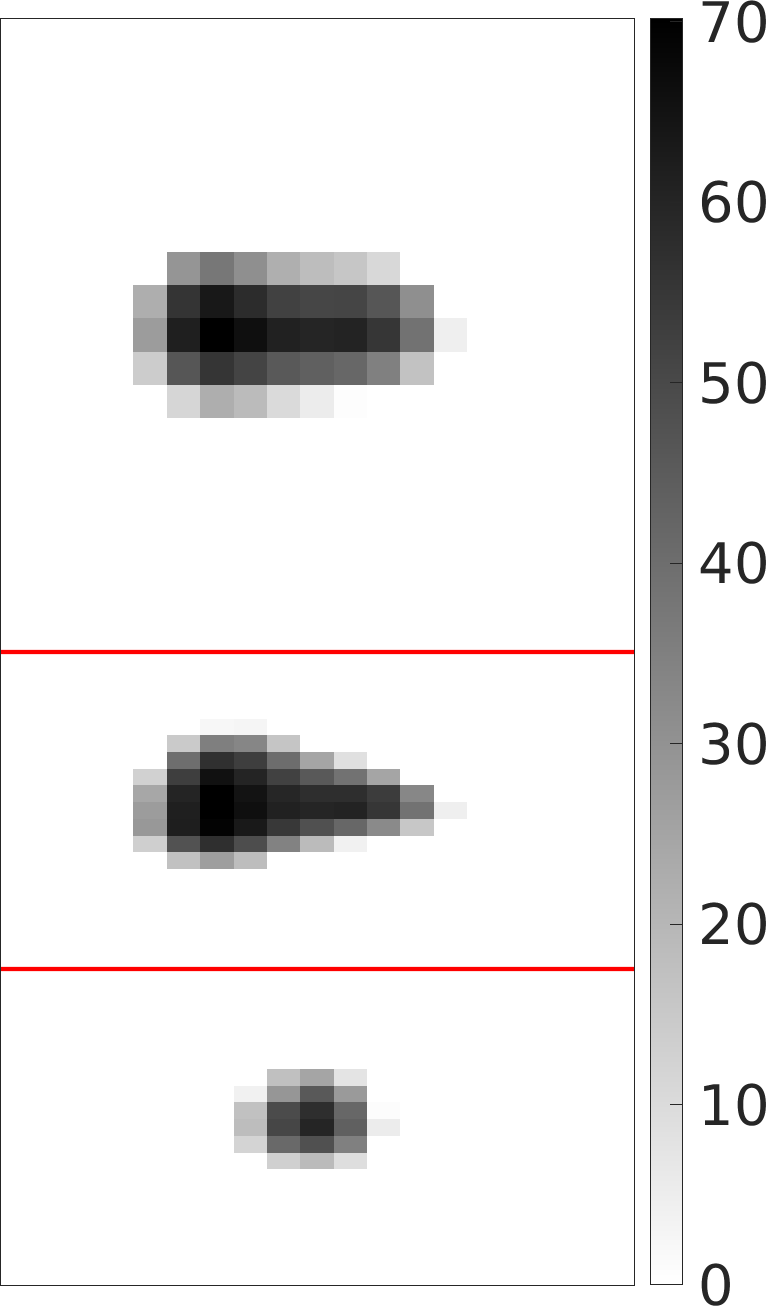}&
 \includegraphics[height=3.4cm]{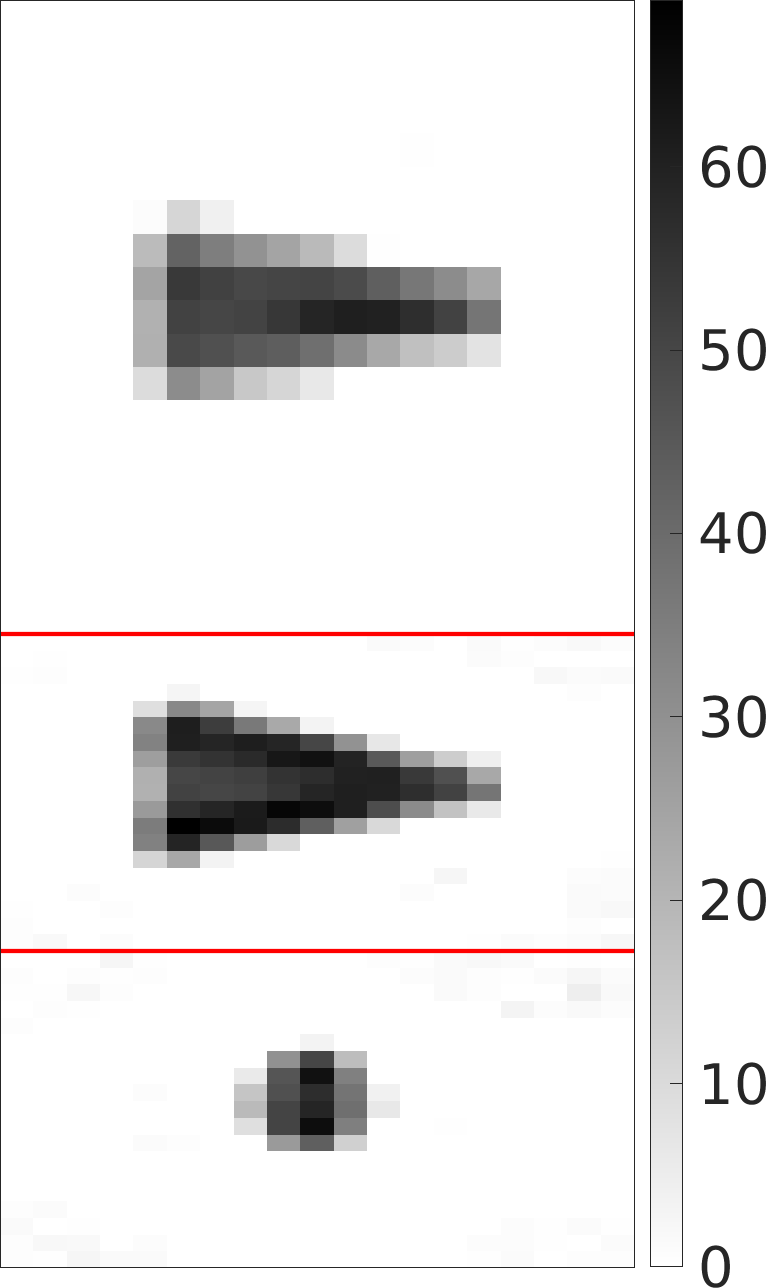}\\
\hline
\multicolumn{6}{l}{$\tau=3$} \\
 \includegraphics[height=3.4cm]{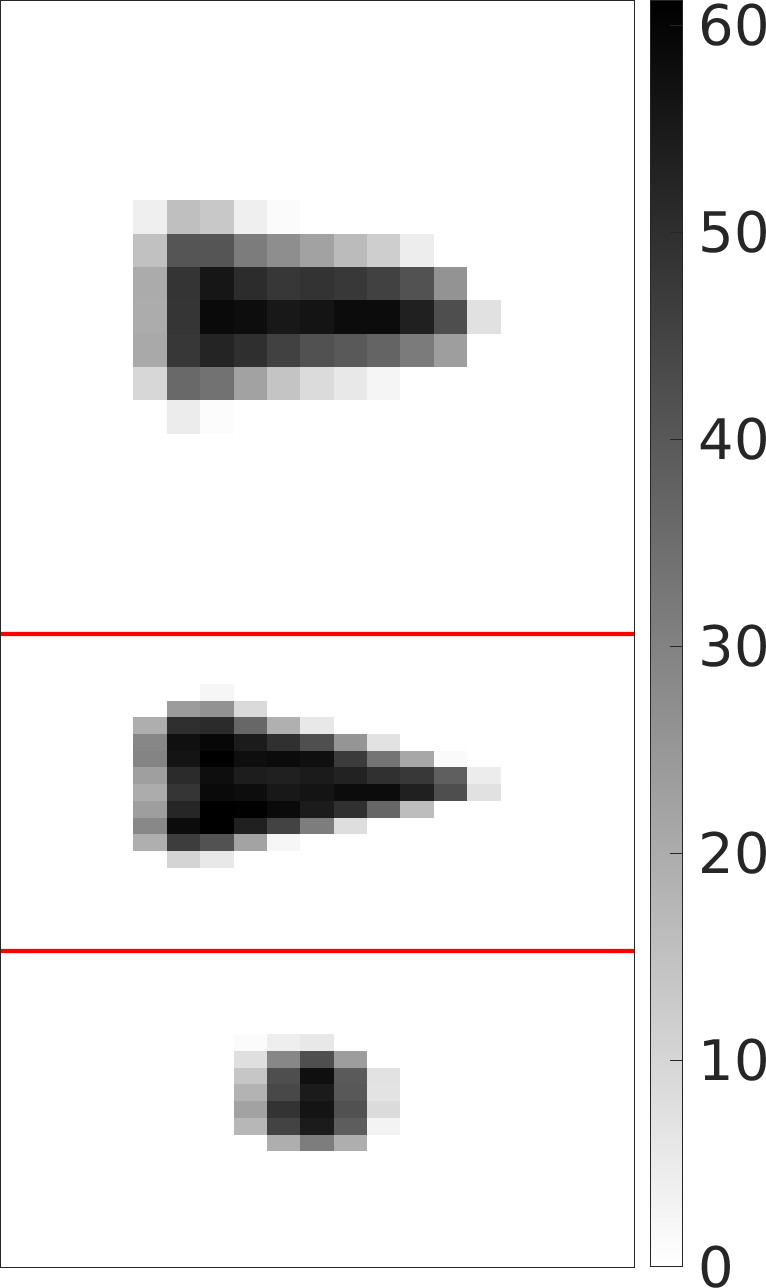}&
 \includegraphics[height=3.4cm]{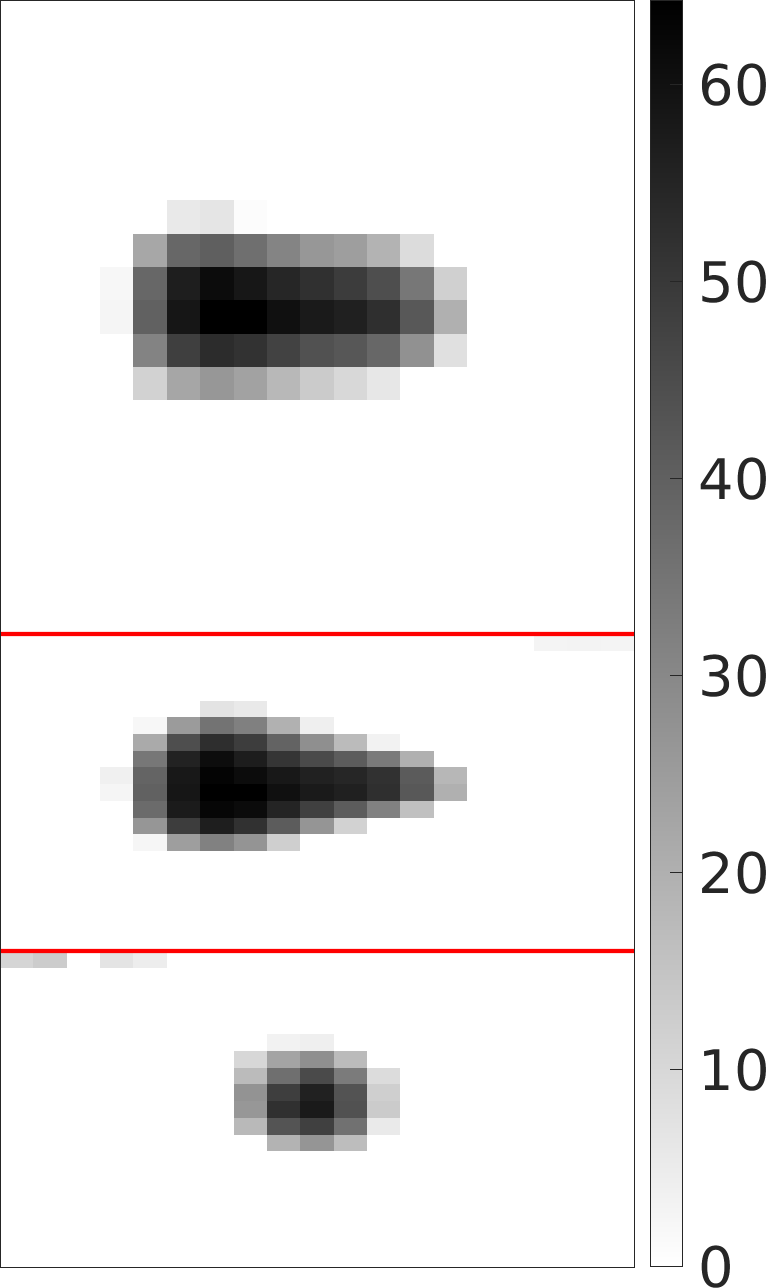}&
 \includegraphics[height=3.4cm]{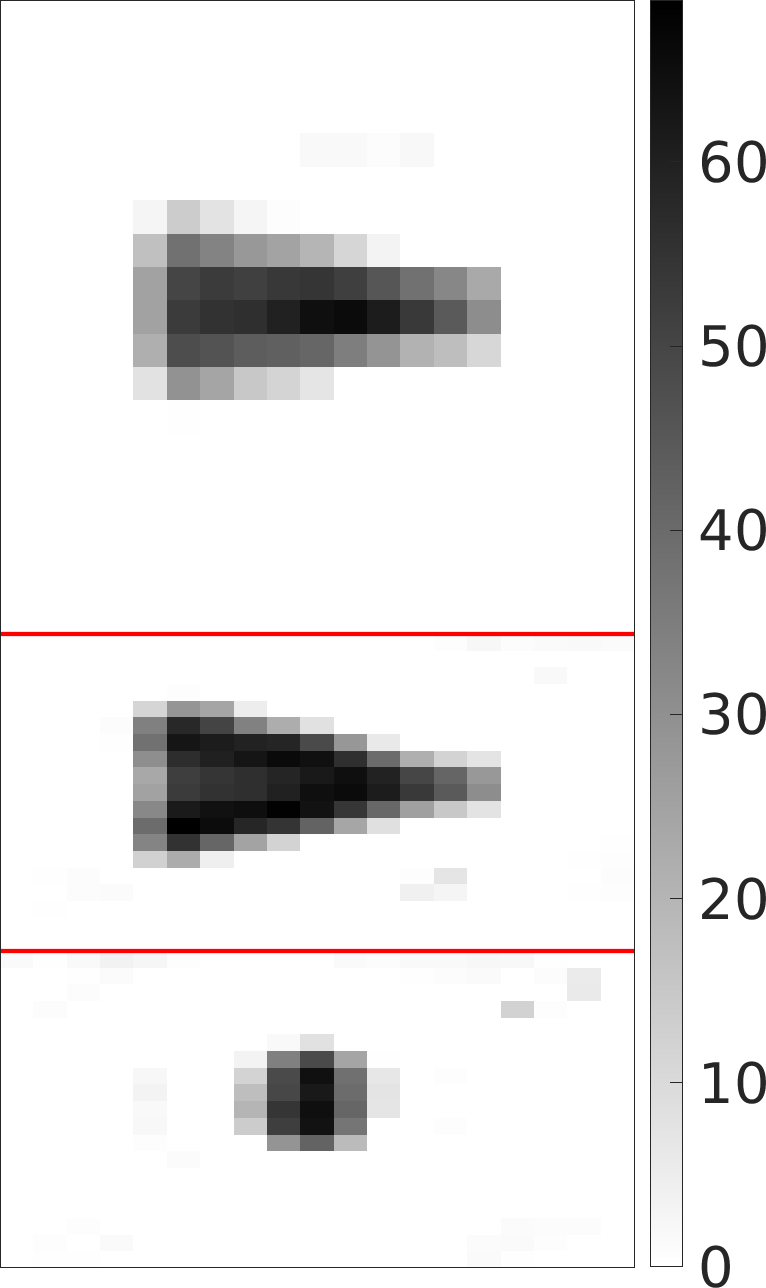}&
 \includegraphics[height=3.4cm]{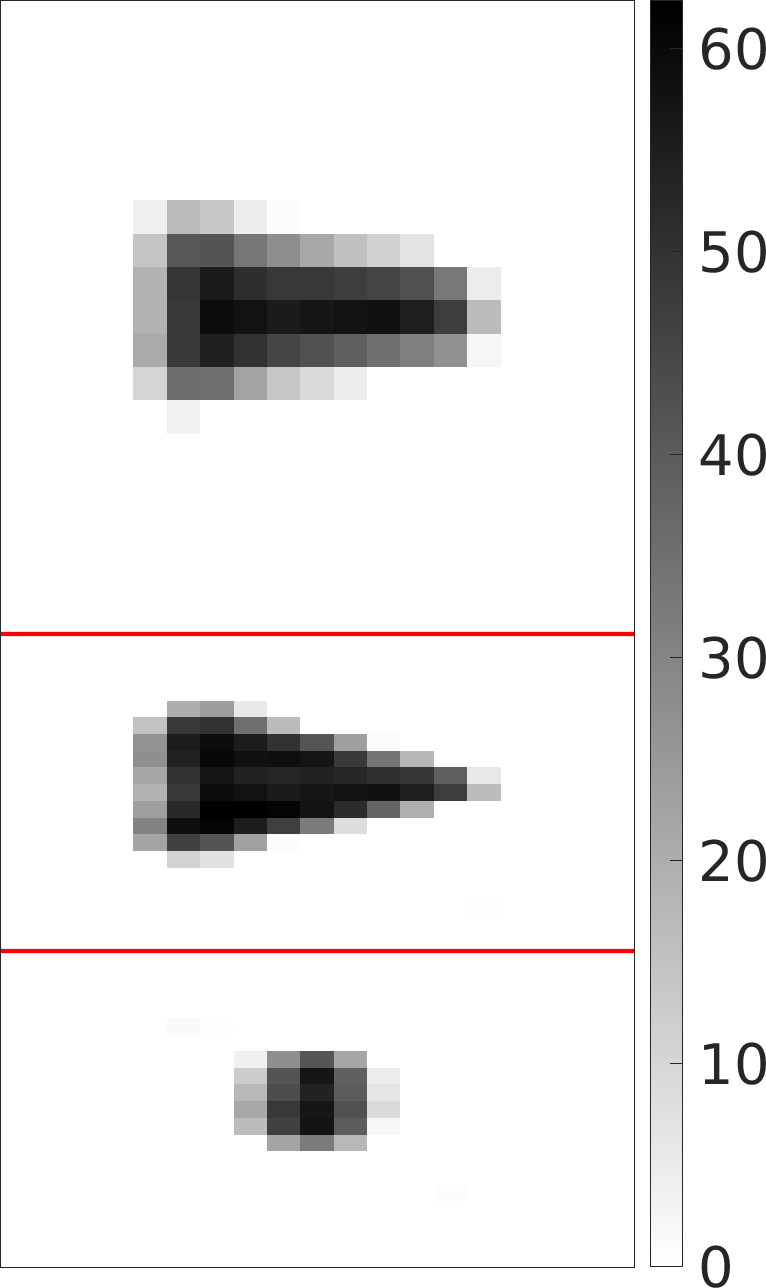}&
 \includegraphics[height=3.4cm]{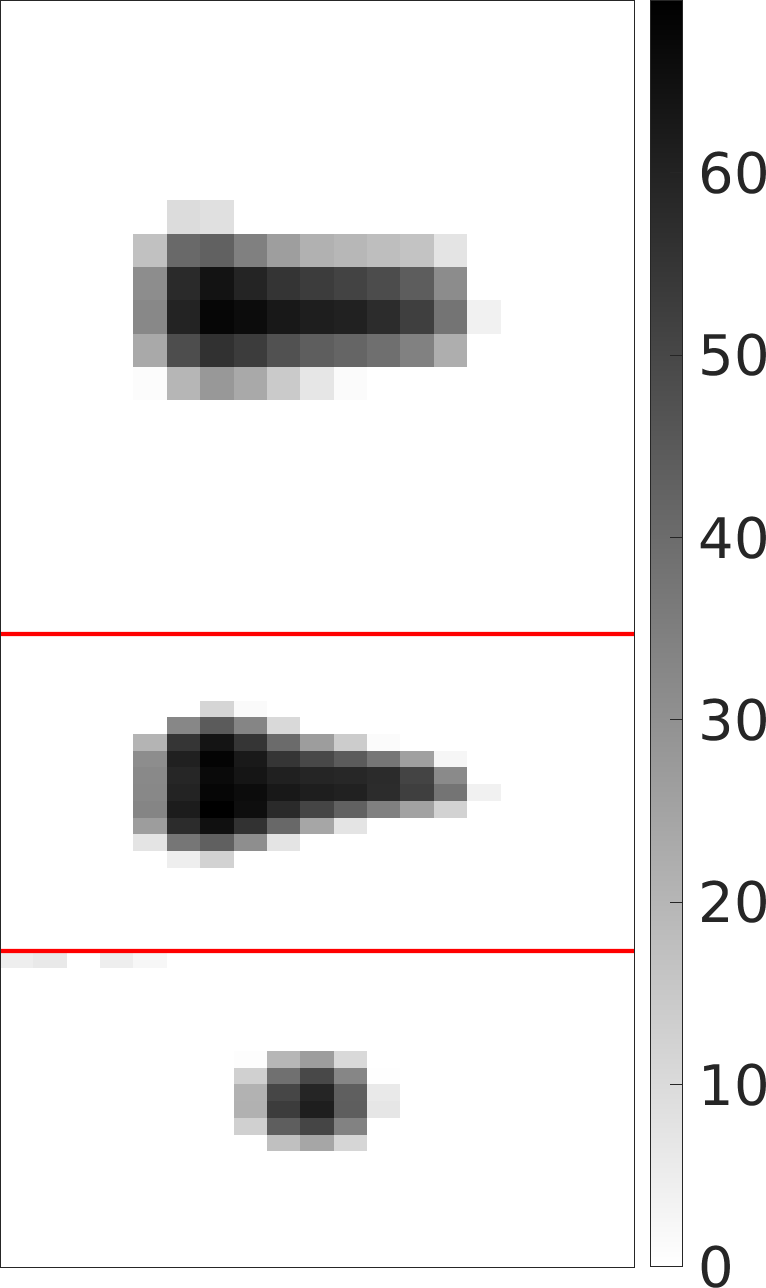}&
 \includegraphics[height=3.4cm]{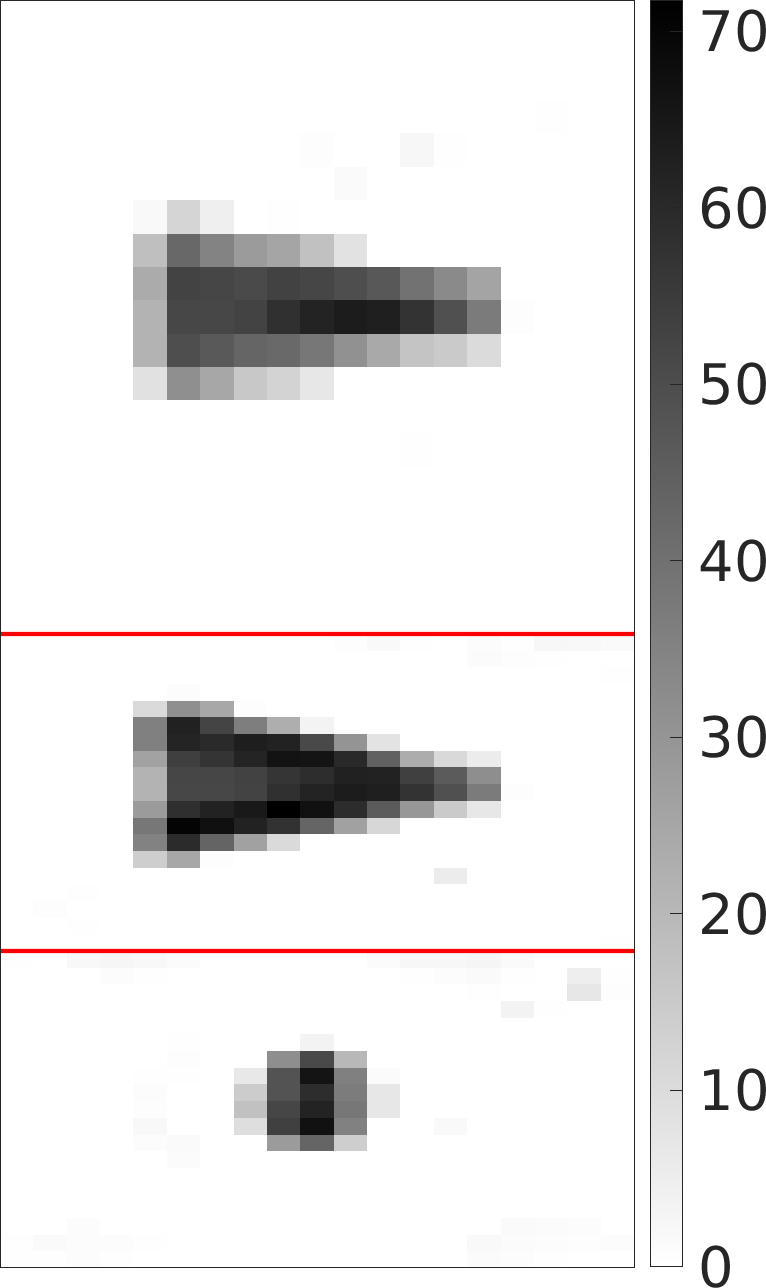}\\
 \hline
\multicolumn{6}{l}{$\tau=5$} \\
 \includegraphics[height=3.4cm]{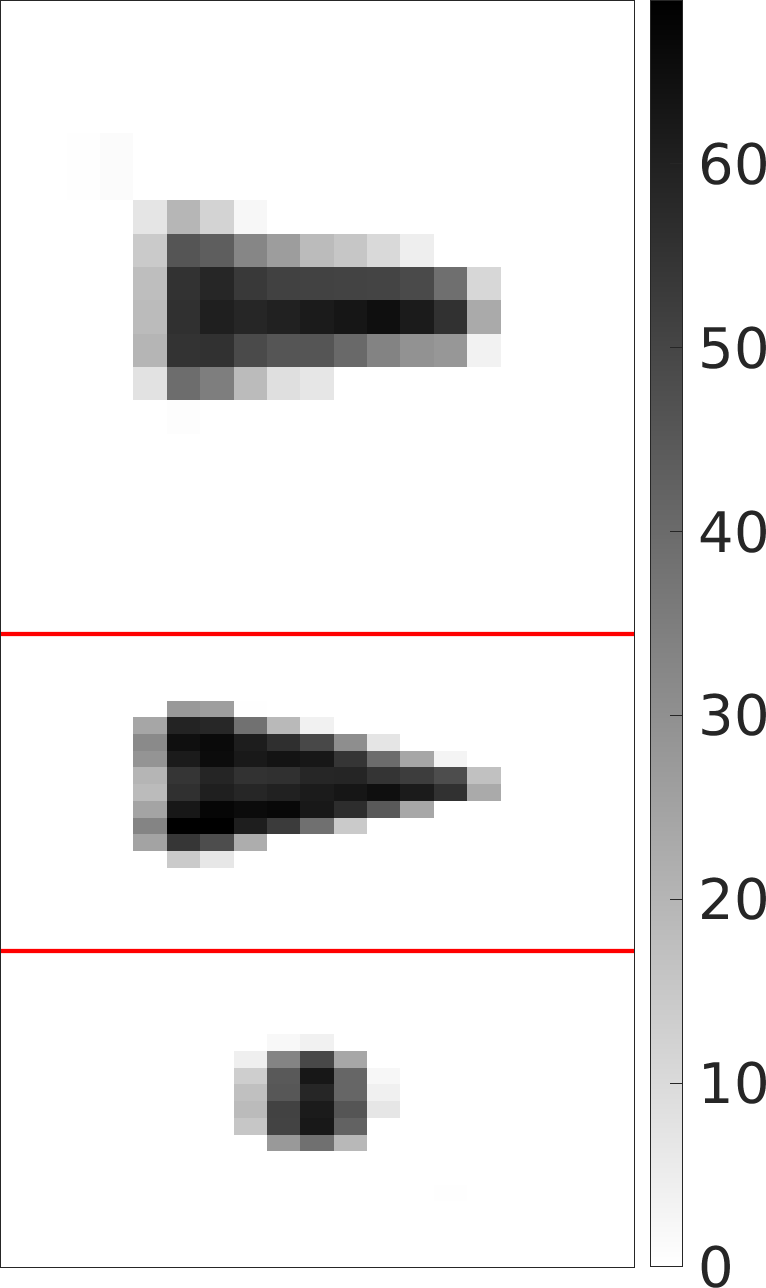}&
 \includegraphics[height=3.4cm]{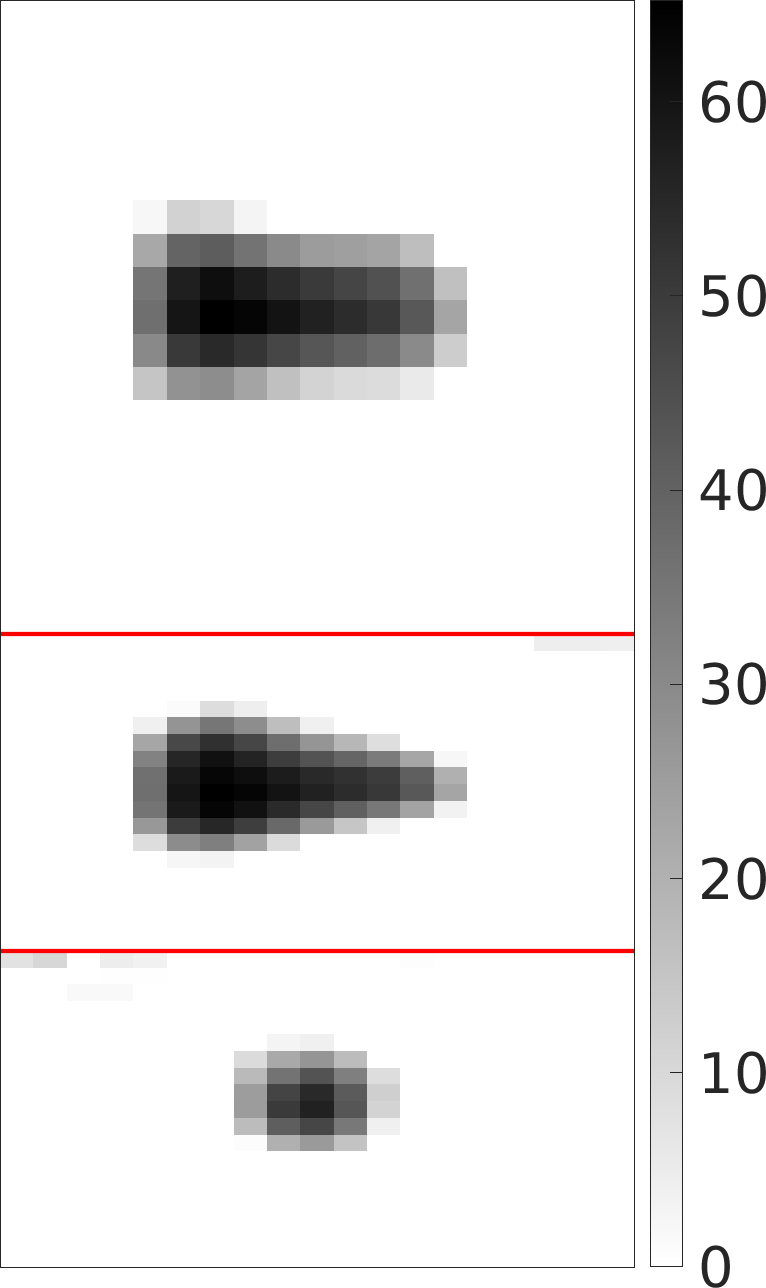}&
 \includegraphics[height=3.4cm]{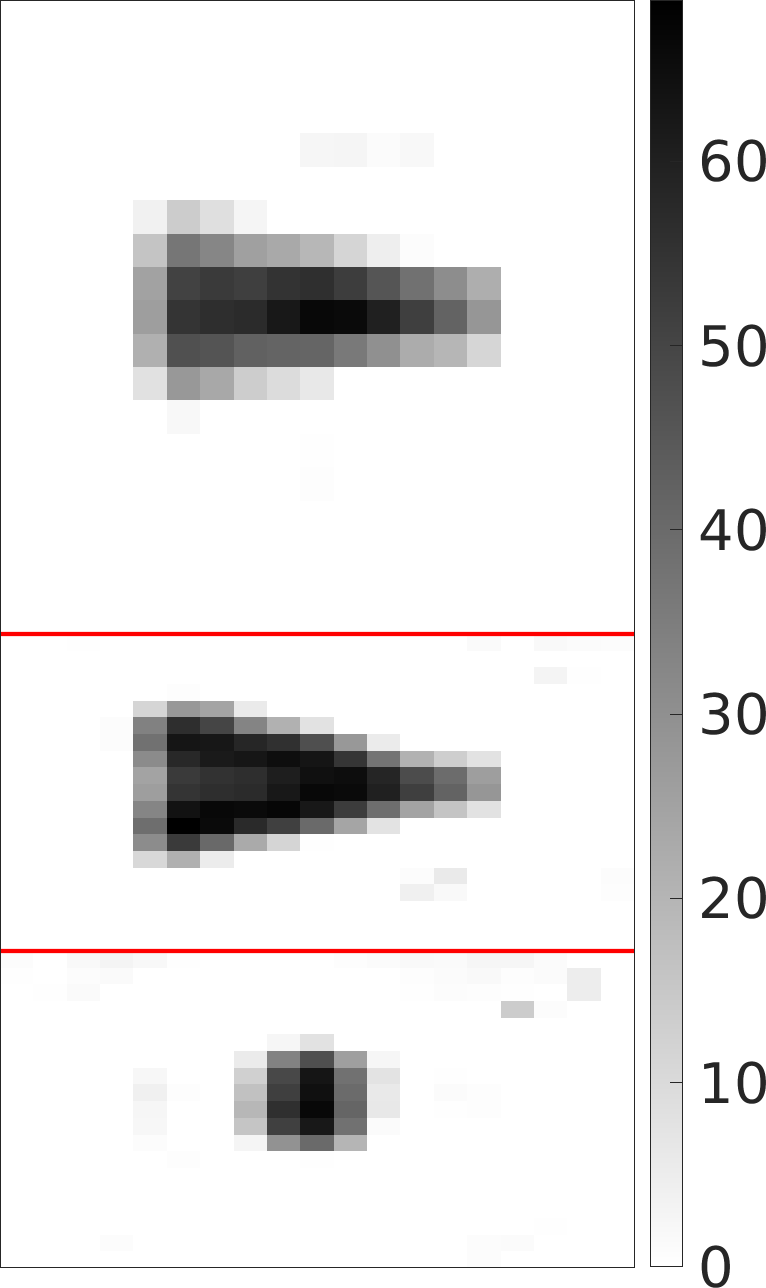}&
 \includegraphics[height=3.4cm]{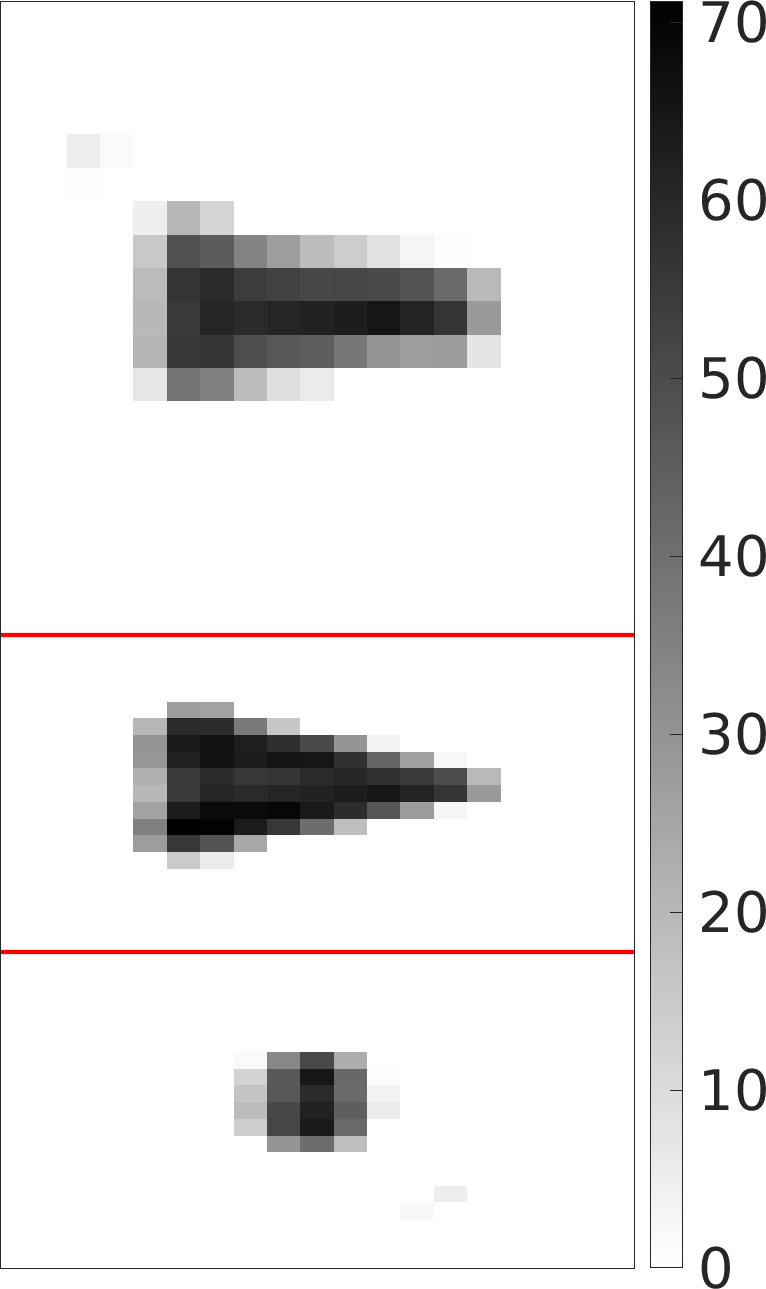}&
 \includegraphics[height=3.4cm]{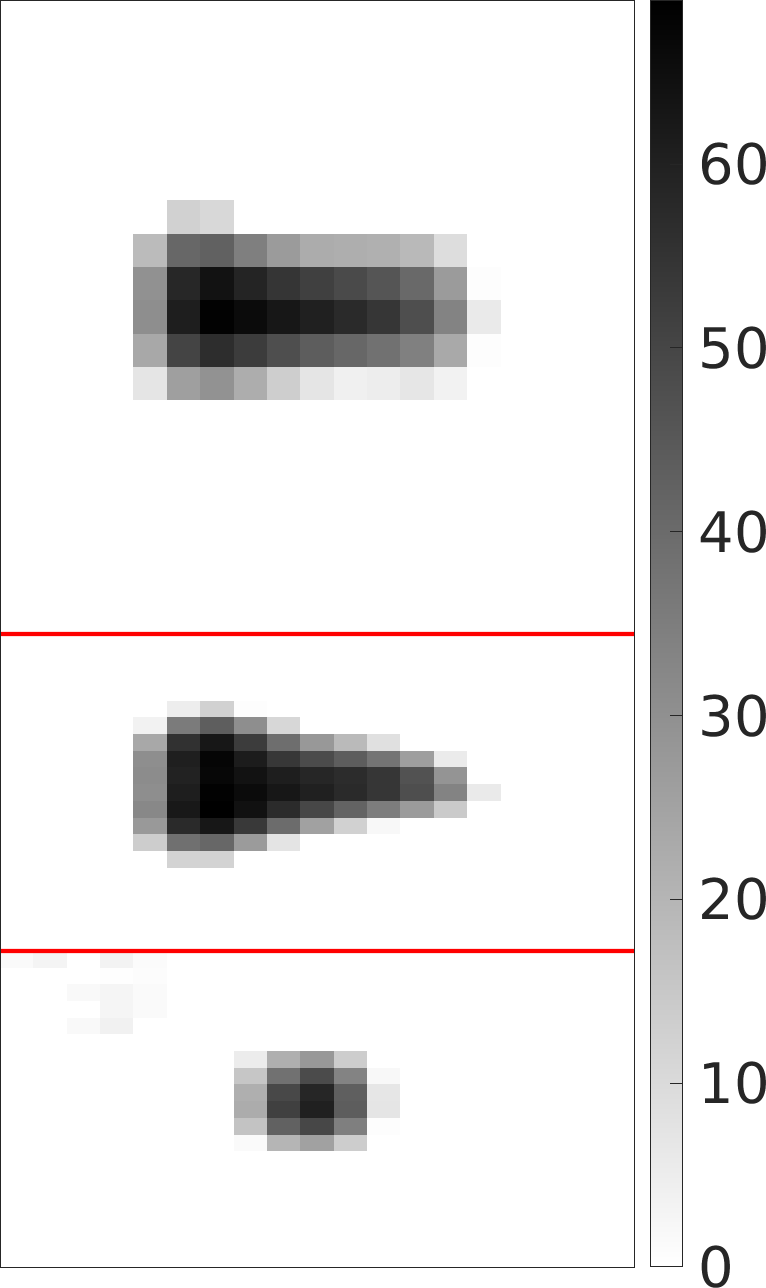}&
 \includegraphics[height=3.4cm]{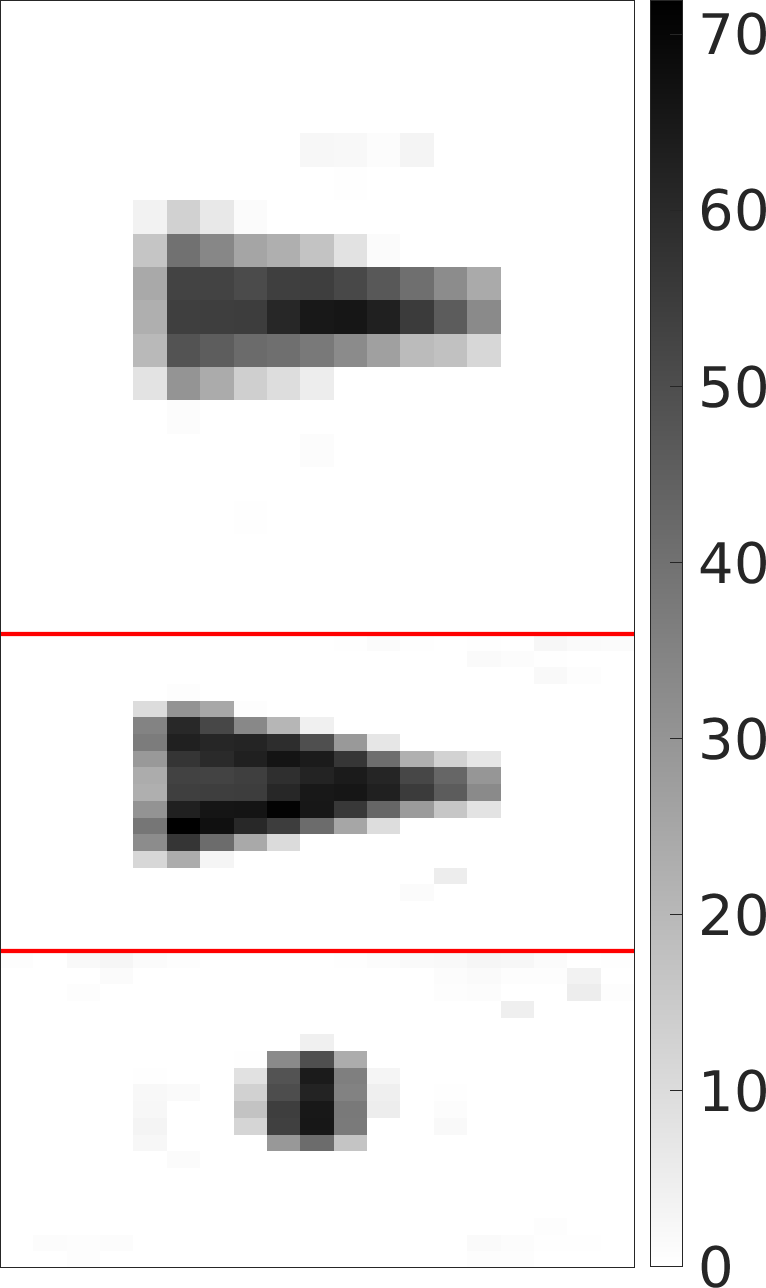}\\
\end{tabular}
}
\caption{Fig. \ref{fig:methods_nonwhitened_vs_whitened_shape_psnr} with inverted colormap: ``Shape'' phantom reconstructions, PSNR-optimized $\alpha$ and iteration number $N$ (for l2-K only) according to Table \ref{tab:psnr_nonwhitened_vs_whitened}.}
\label{fig:methods_nonwhitened_vs_whitened_shape_psnr_inverted_colormap}
\end{figure}

\begin{figure}[hbt!]
\centering
\scalebox{0.85}{
\begin{tabular}{ccc|ccc}
\multicolumn{3}{c|}{non-whitened} & \multicolumn{3}{c}{whitened} \\
\hline
l1-L & l2-L & l2-K & l1-L & l2-L & l2-K \\
\hline
\multicolumn{6}{l}{$\tau=0$} \\
 \includegraphics[height=3.4cm]{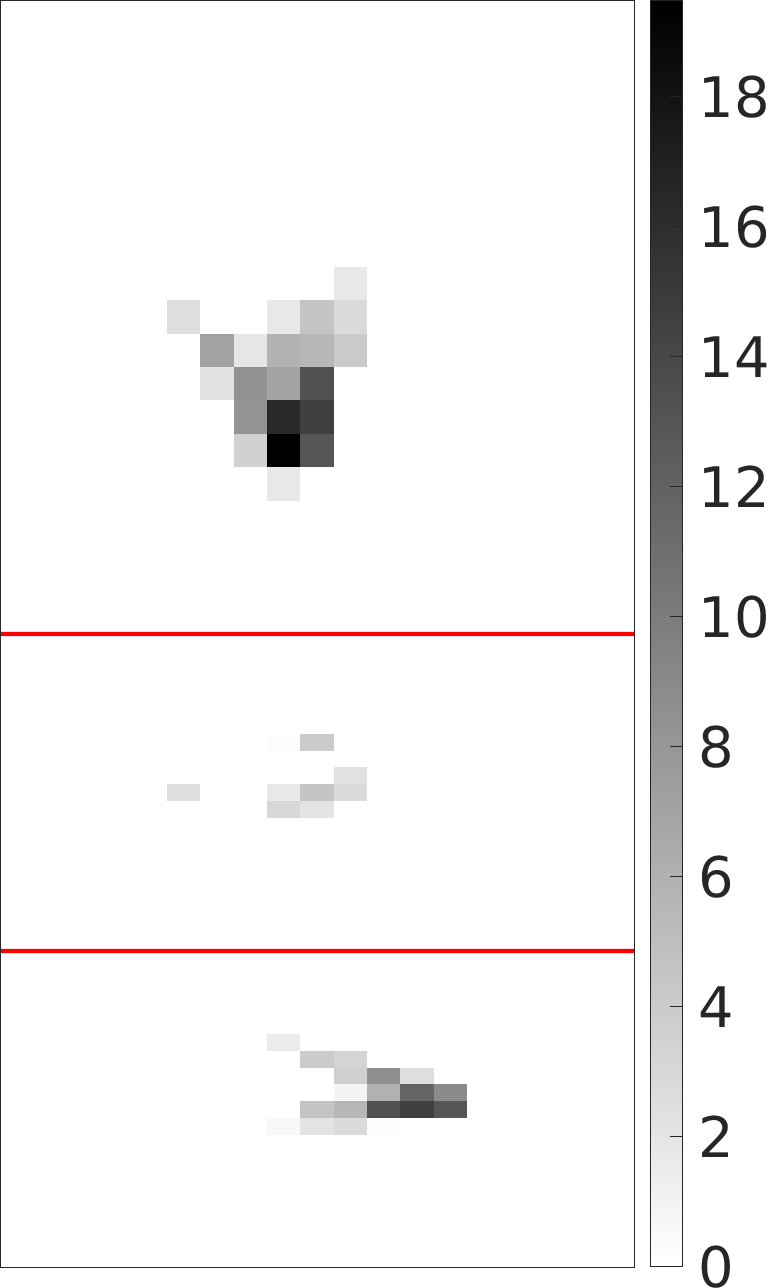}&
 \includegraphics[height=3.4cm]{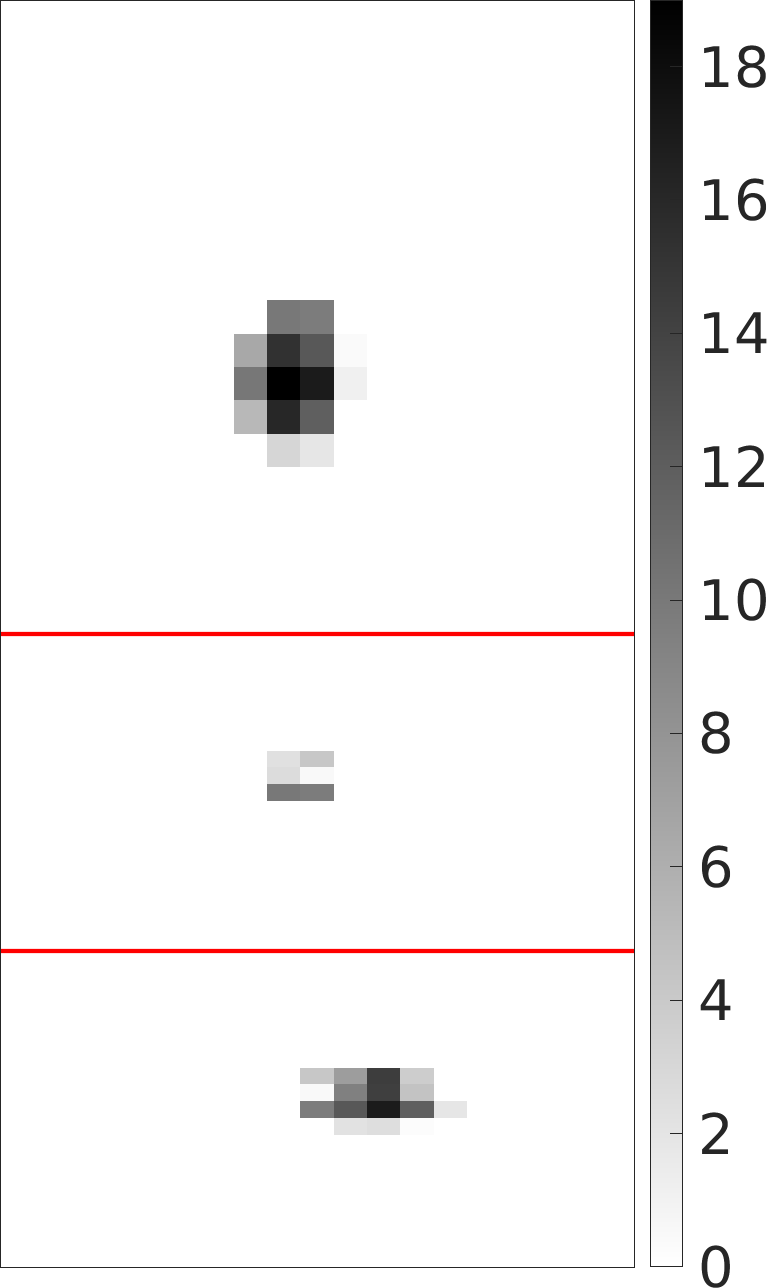}&
 \includegraphics[height=3.4cm]{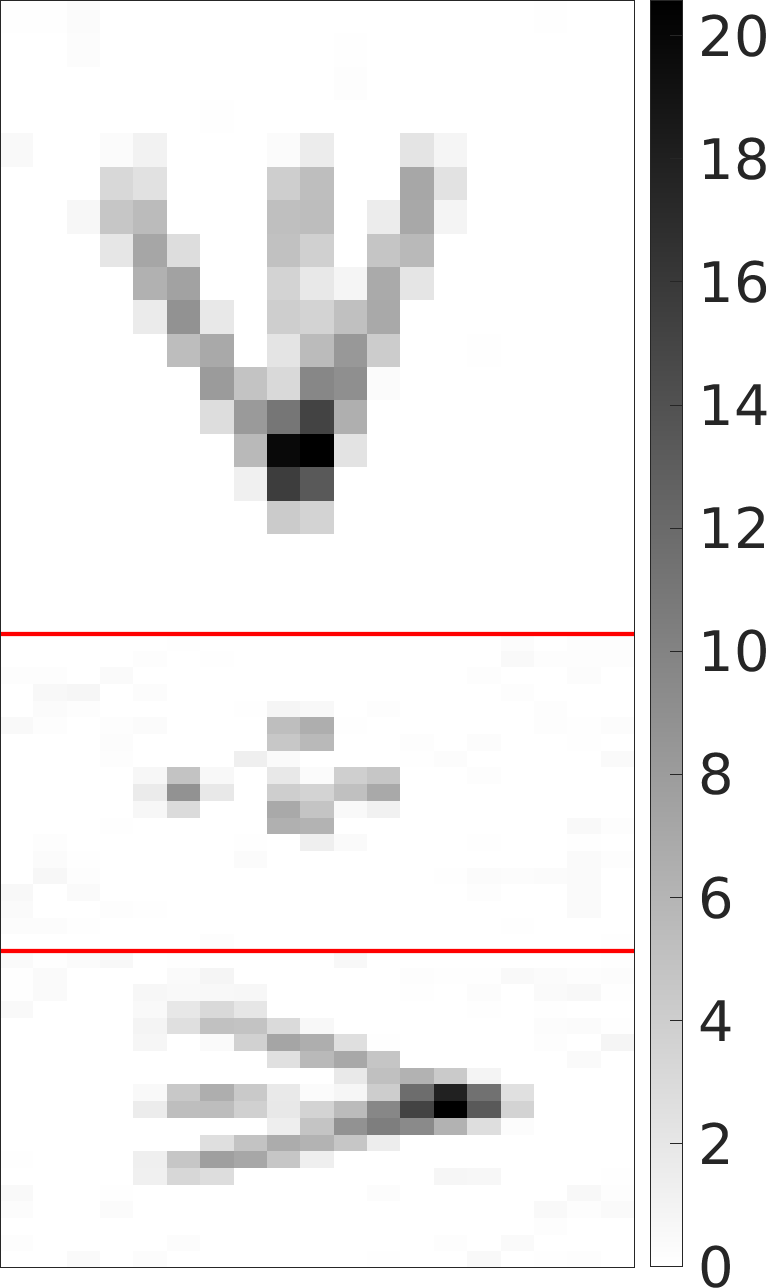}&
 \includegraphics[height=3.4cm]{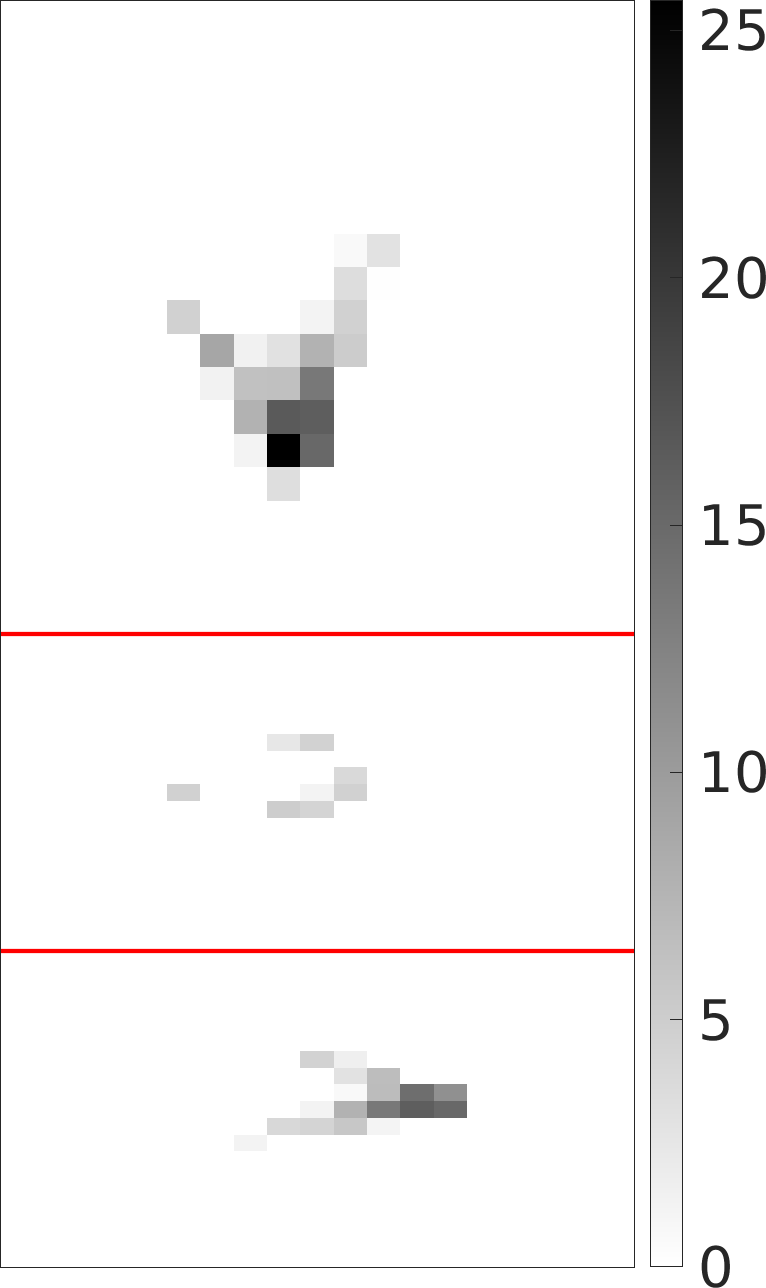}&
 \includegraphics[height=3.4cm]{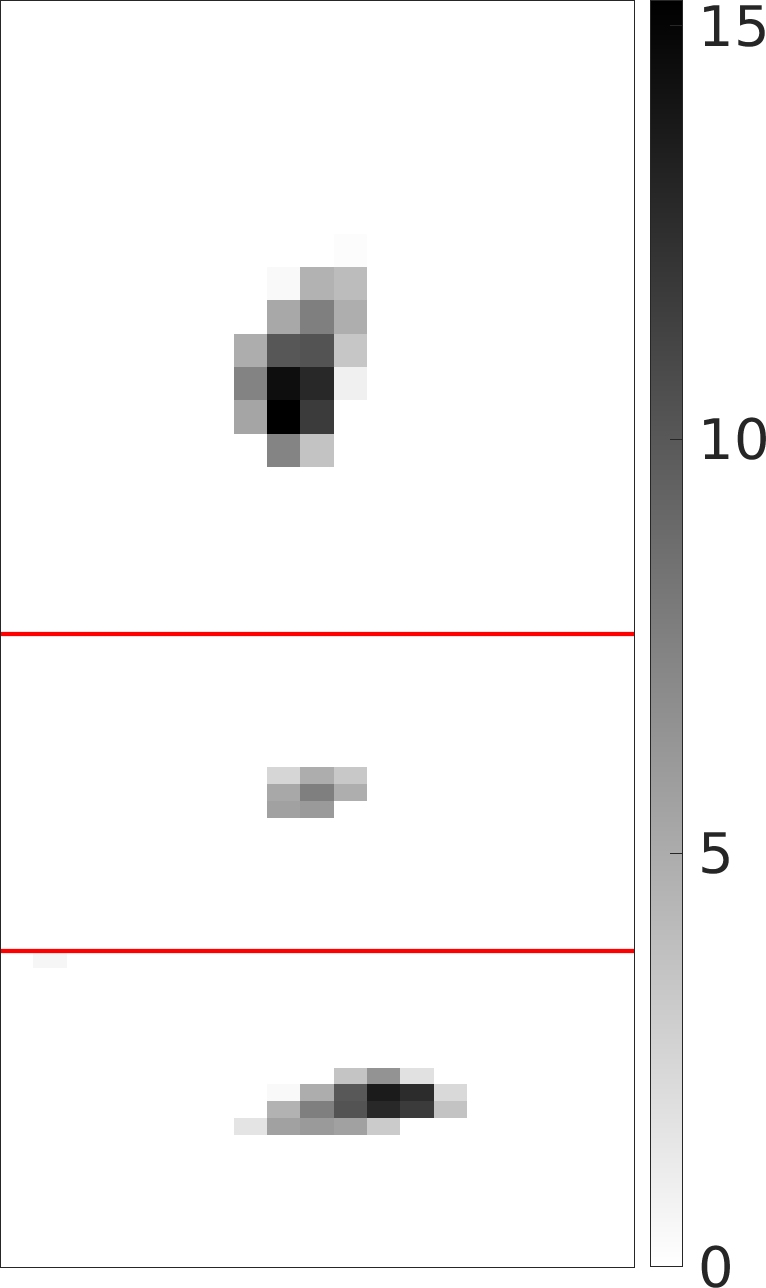}&
 \includegraphics[height=3.4cm]{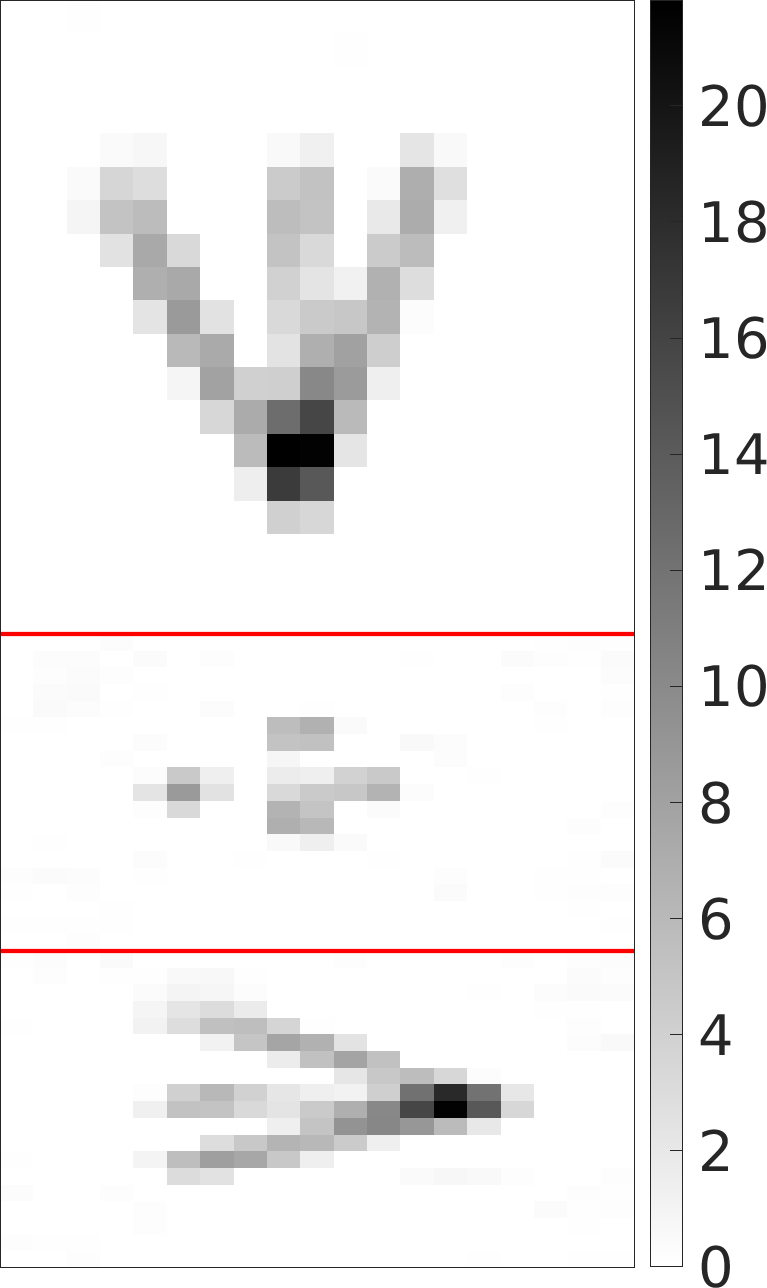}\\
\hline
\multicolumn{6}{l}{$\tau=1$} \\
 \includegraphics[height=3.4cm]{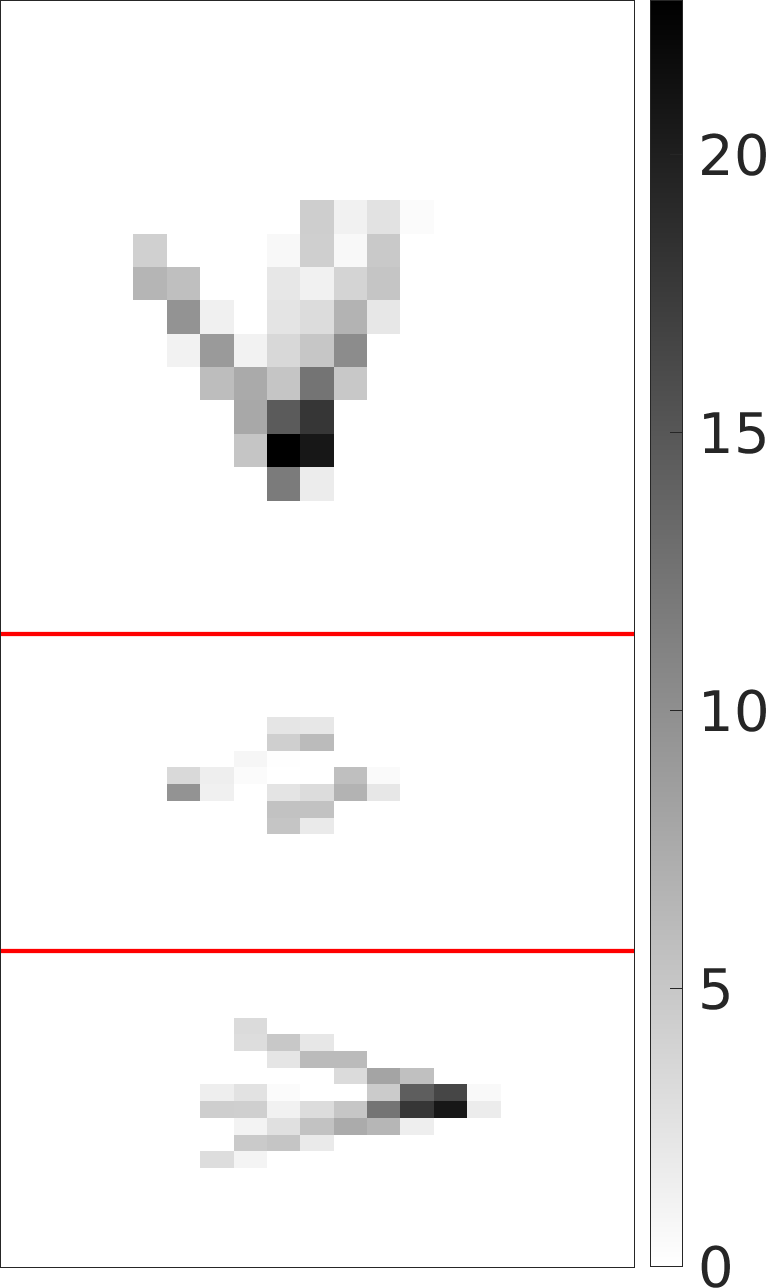}&
 \includegraphics[height=3.4cm]{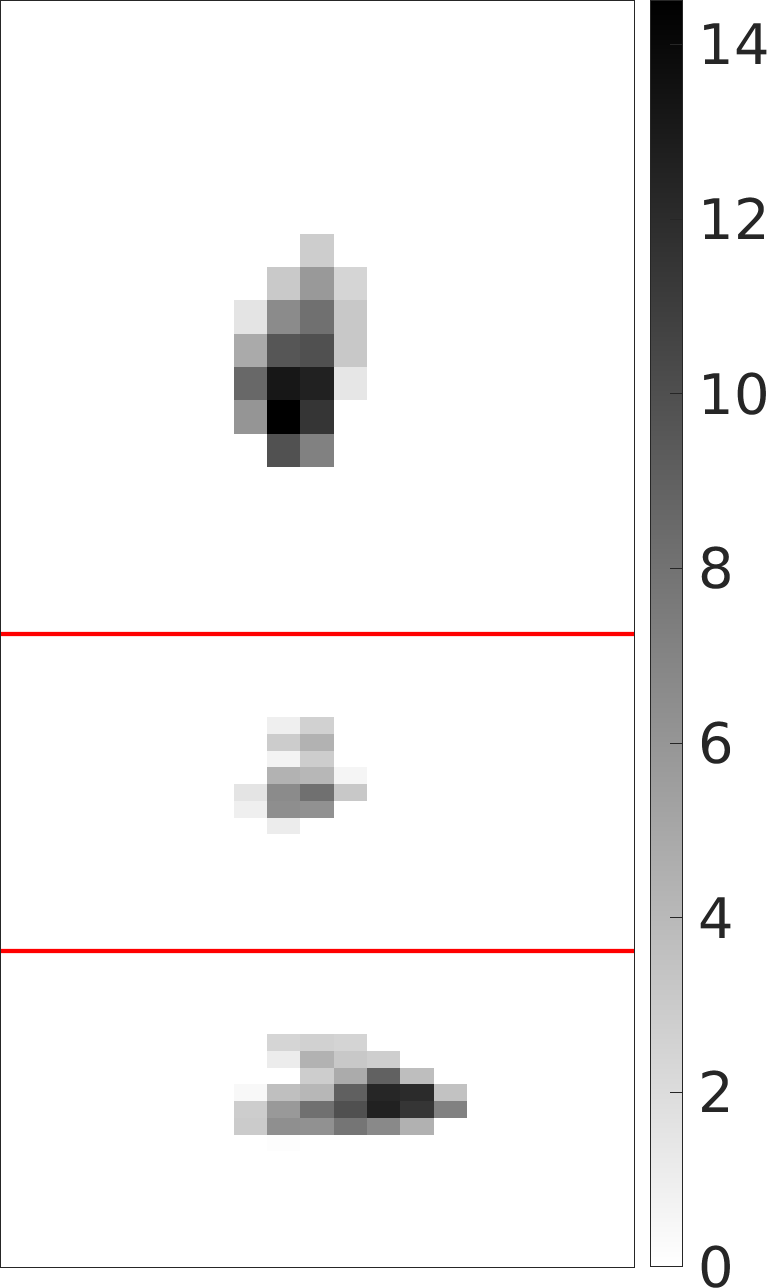}&
 \includegraphics[height=3.4cm]{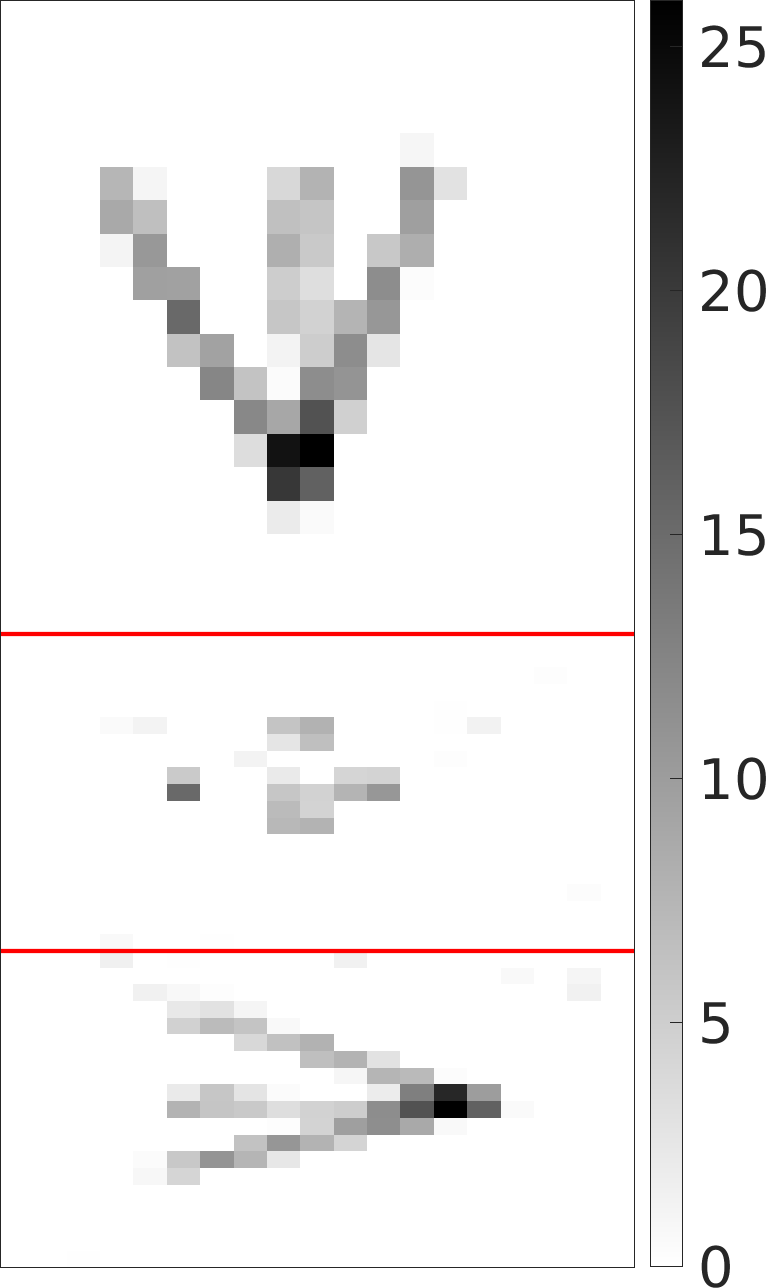}&
 \includegraphics[height=3.4cm]{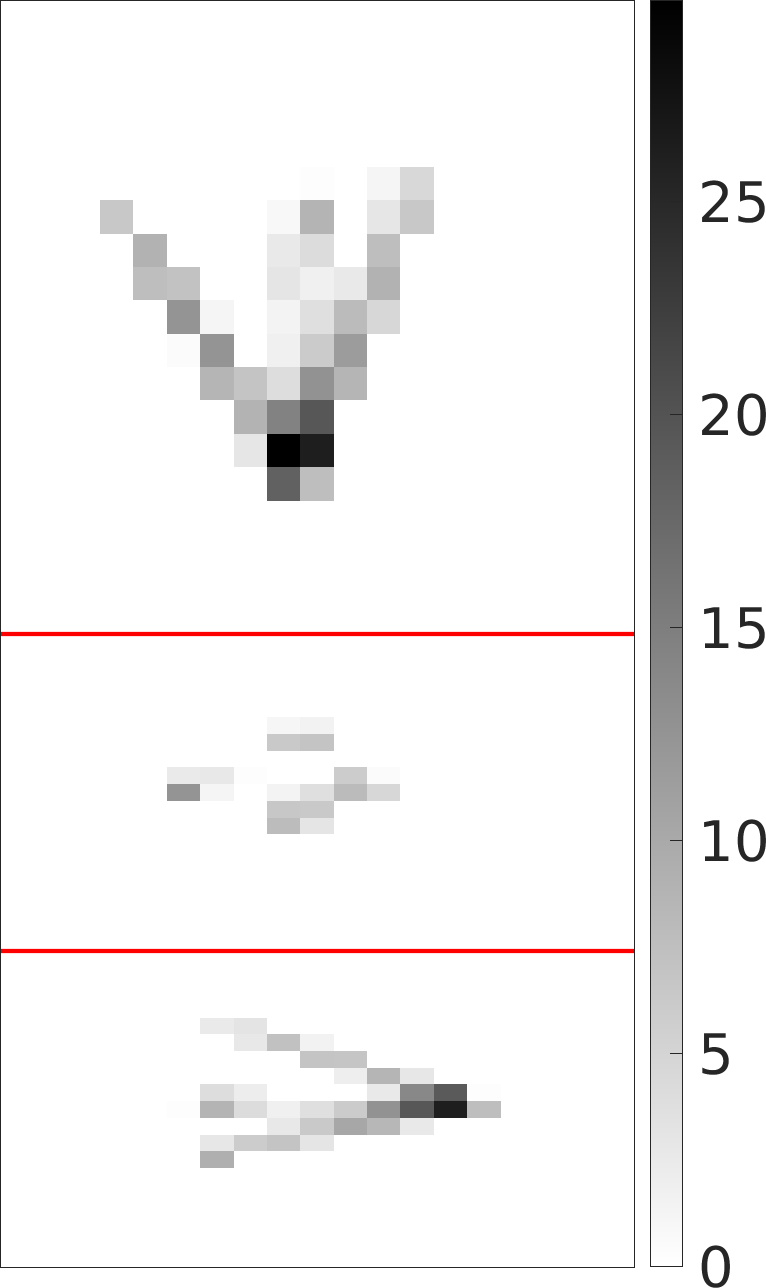}&
 \includegraphics[height=3.4cm]{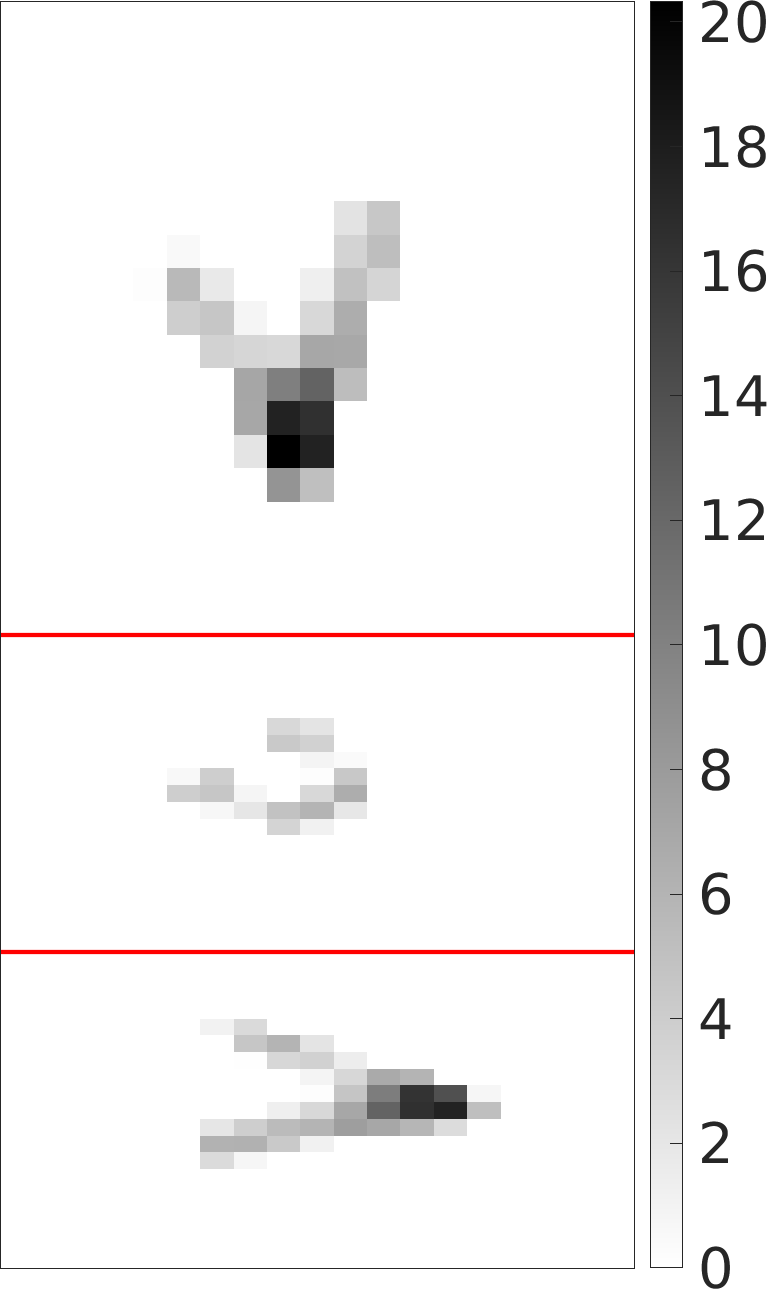}&
 \includegraphics[height=3.4cm]{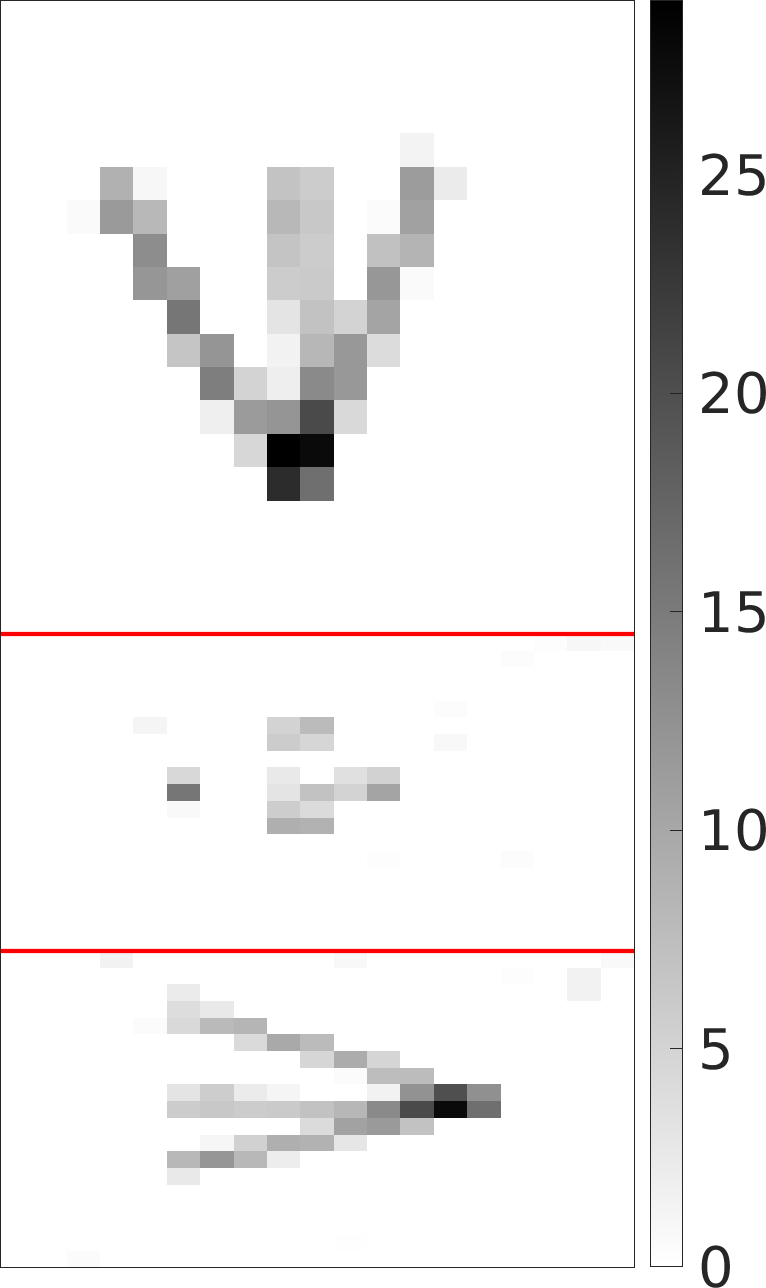}\\
\hline
\multicolumn{6}{l}{$\tau=3$} \\
 \includegraphics[height=3.4cm]{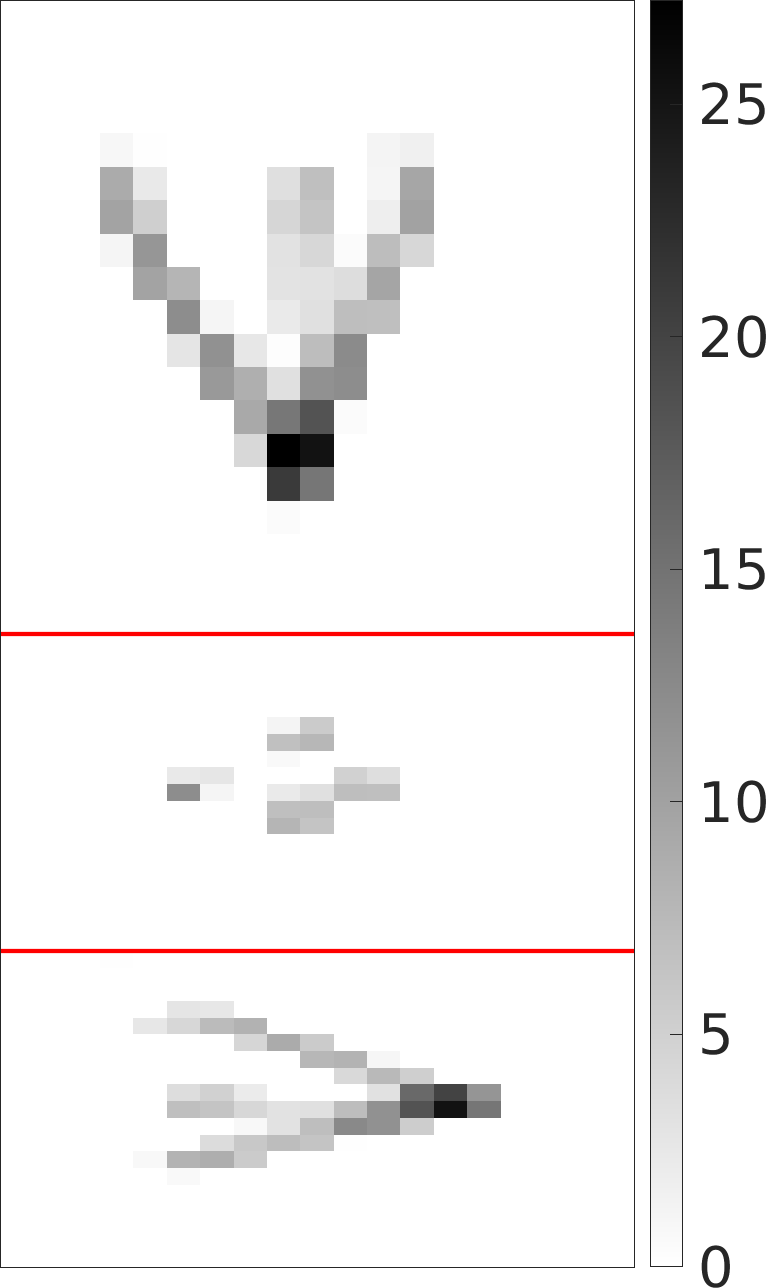}&
 \includegraphics[height=3.4cm]{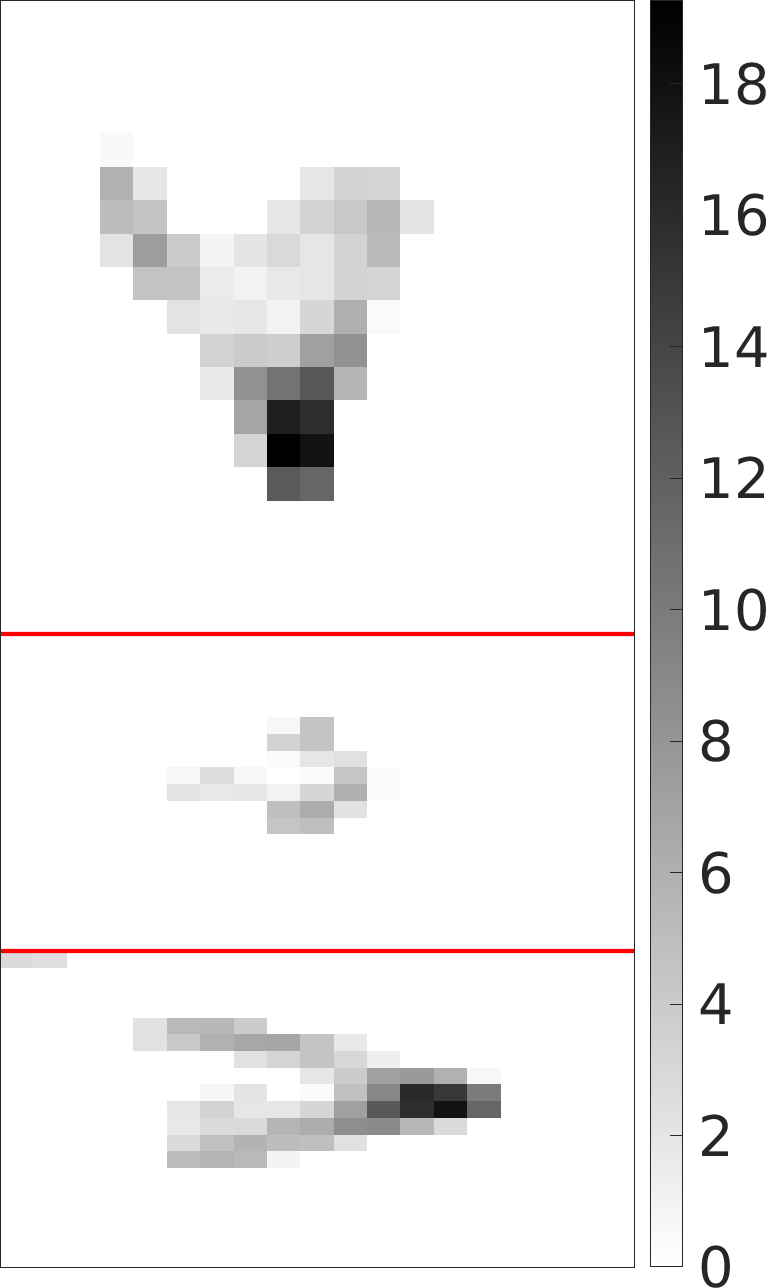}&
 \includegraphics[height=3.4cm]{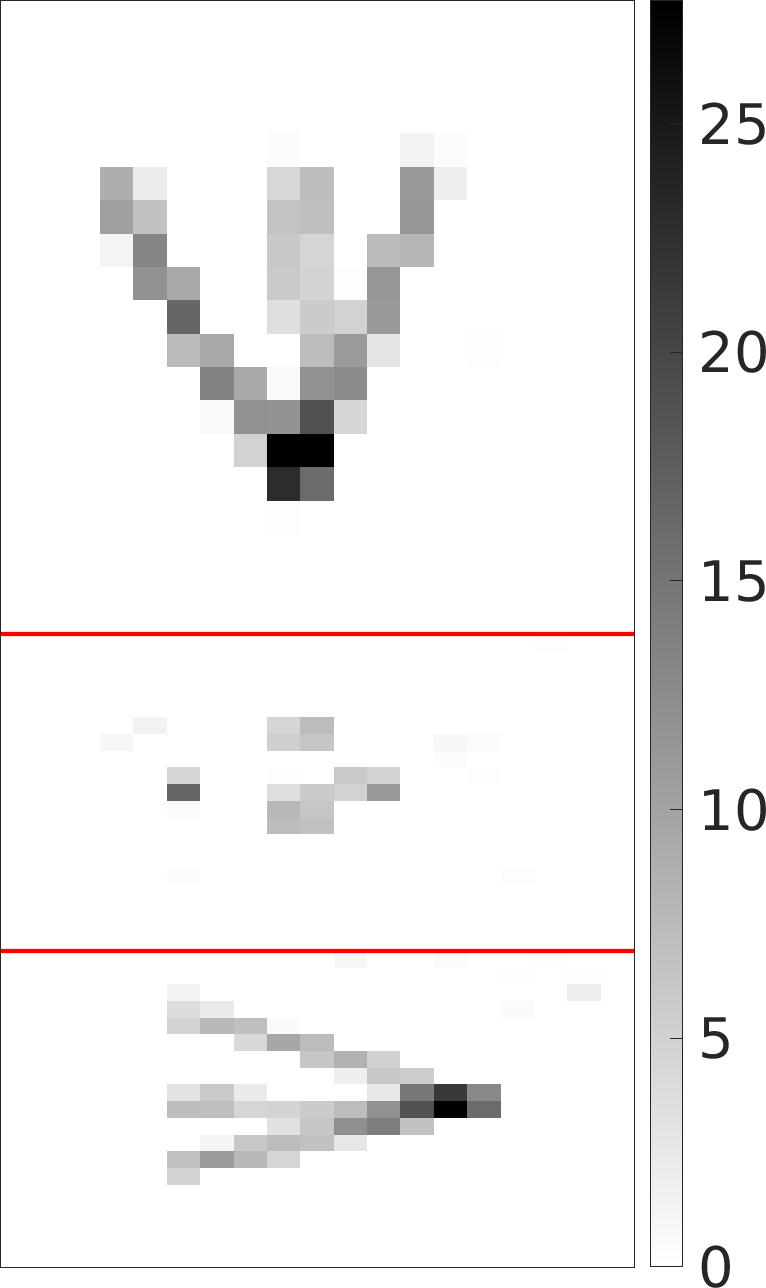}&
 \includegraphics[height=3.4cm]{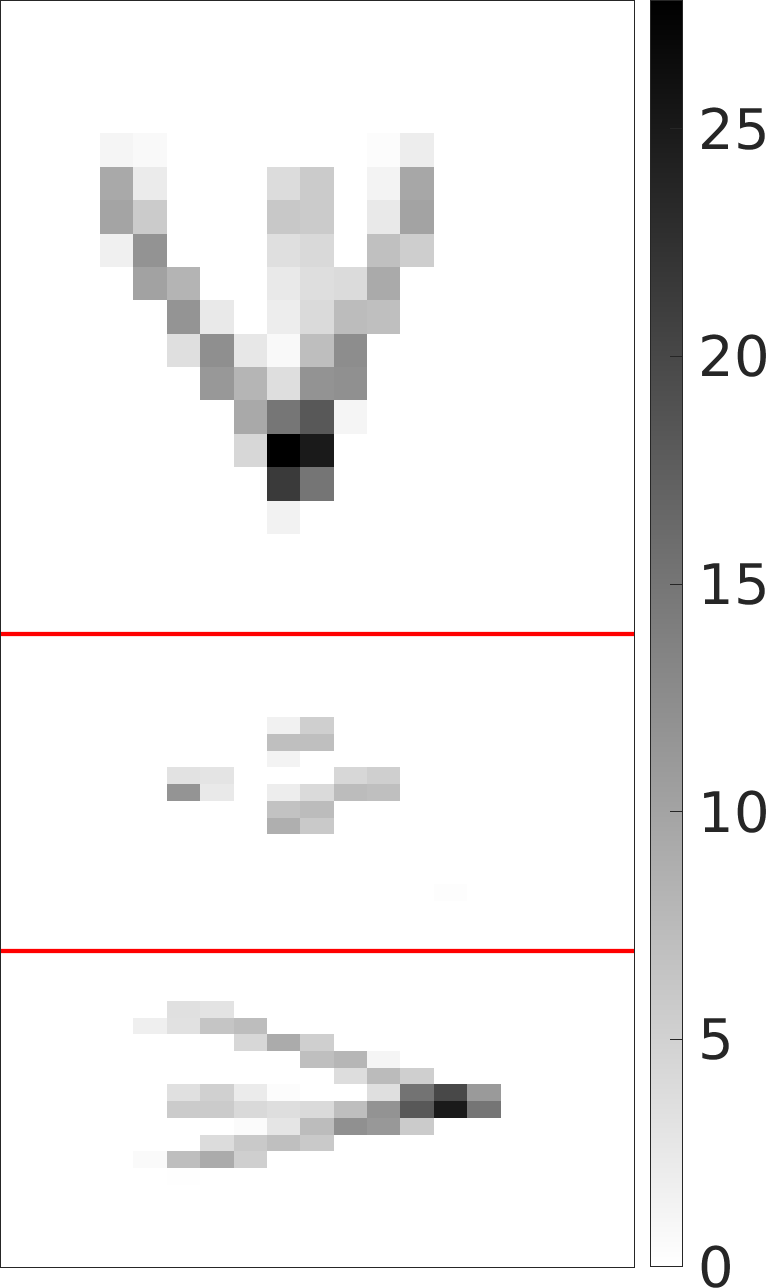}&
 \includegraphics[height=3.4cm]{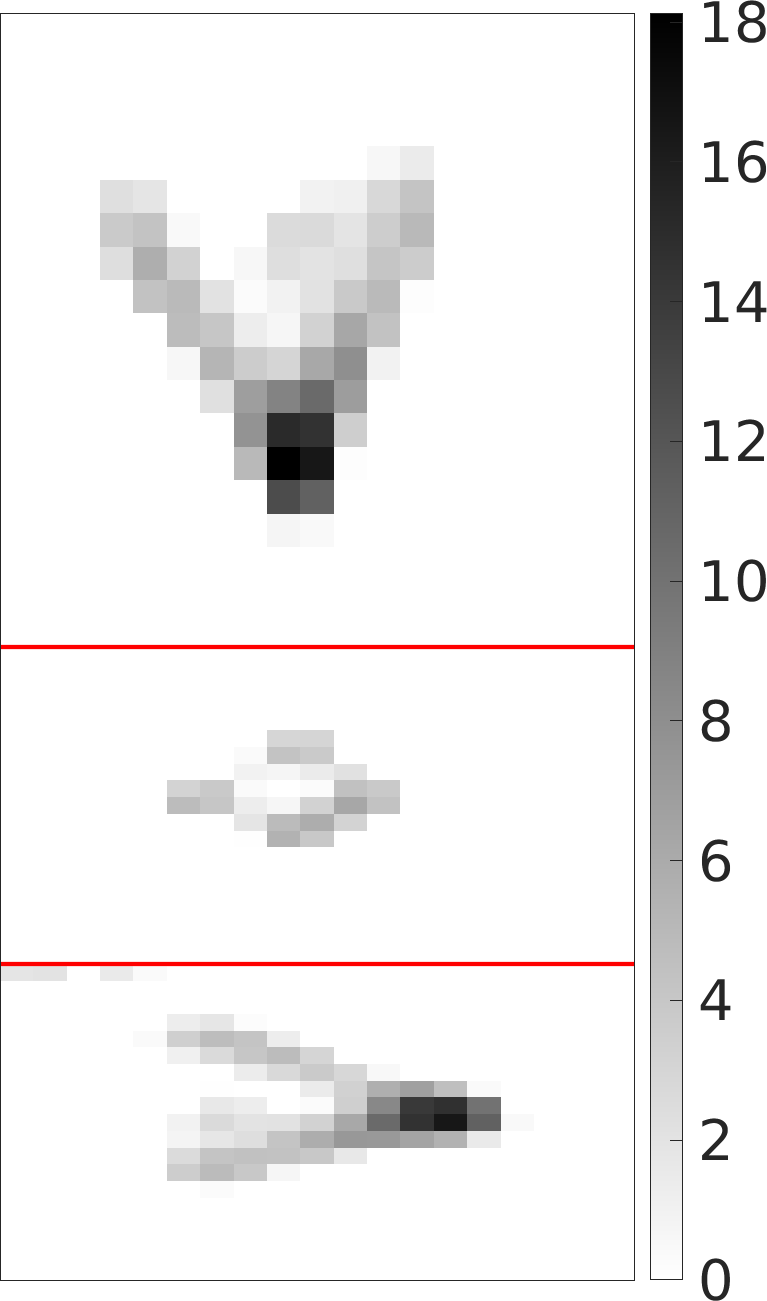}&
 \includegraphics[height=3.4cm]{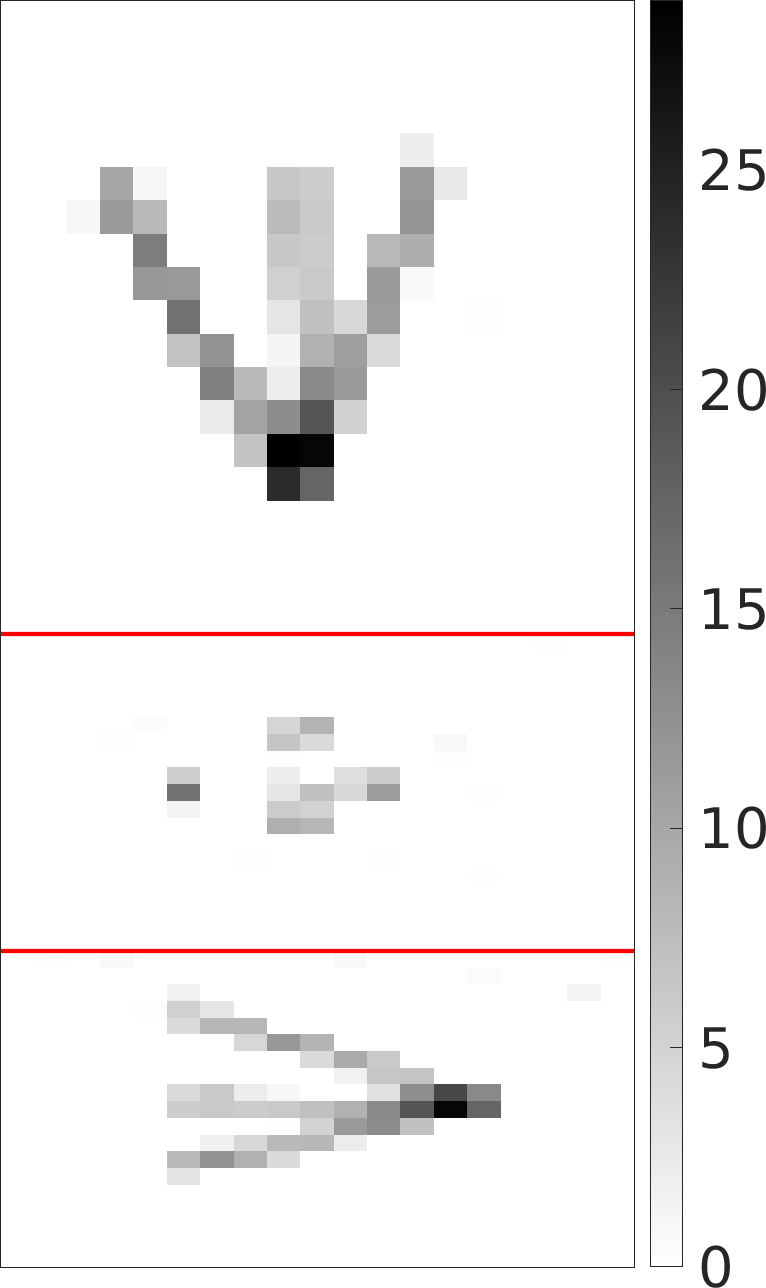}\\
 \hline
\multicolumn{6}{l}{$\tau=5$} \\
 \includegraphics[height=3.4cm]{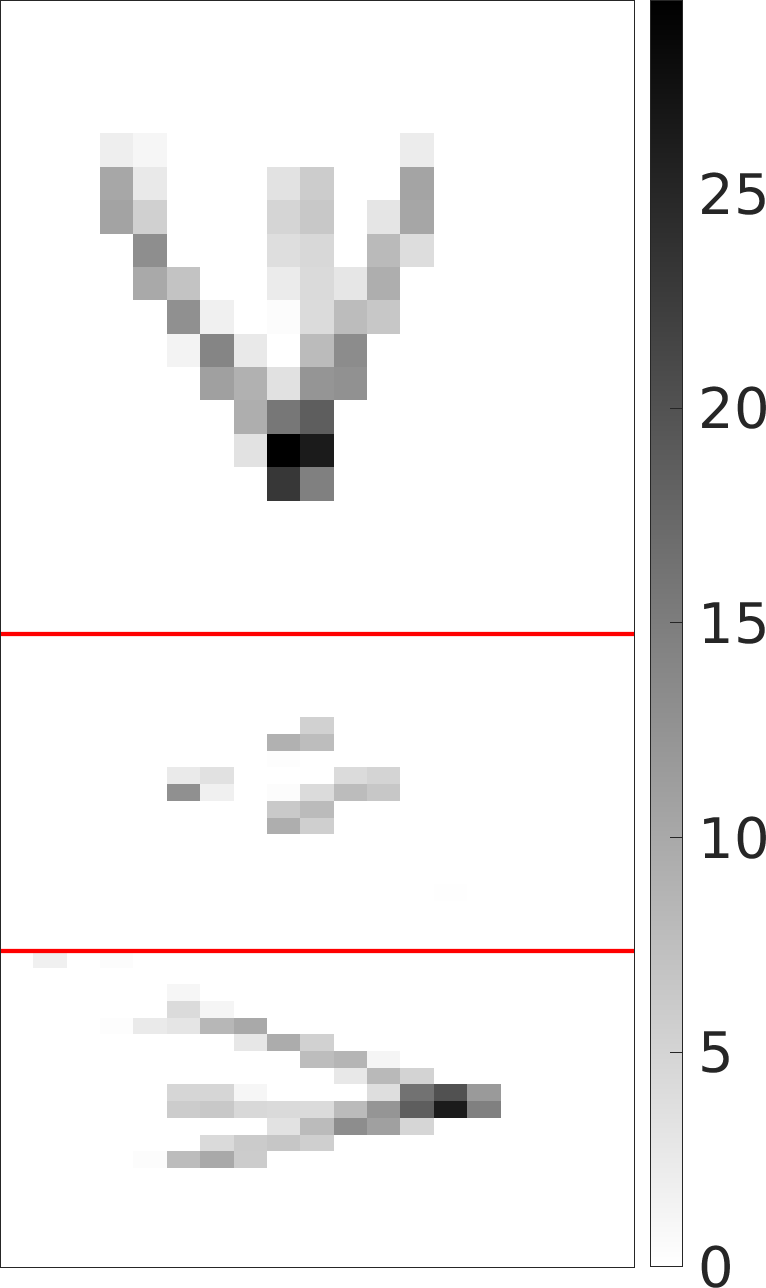}&
 \includegraphics[height=3.4cm]{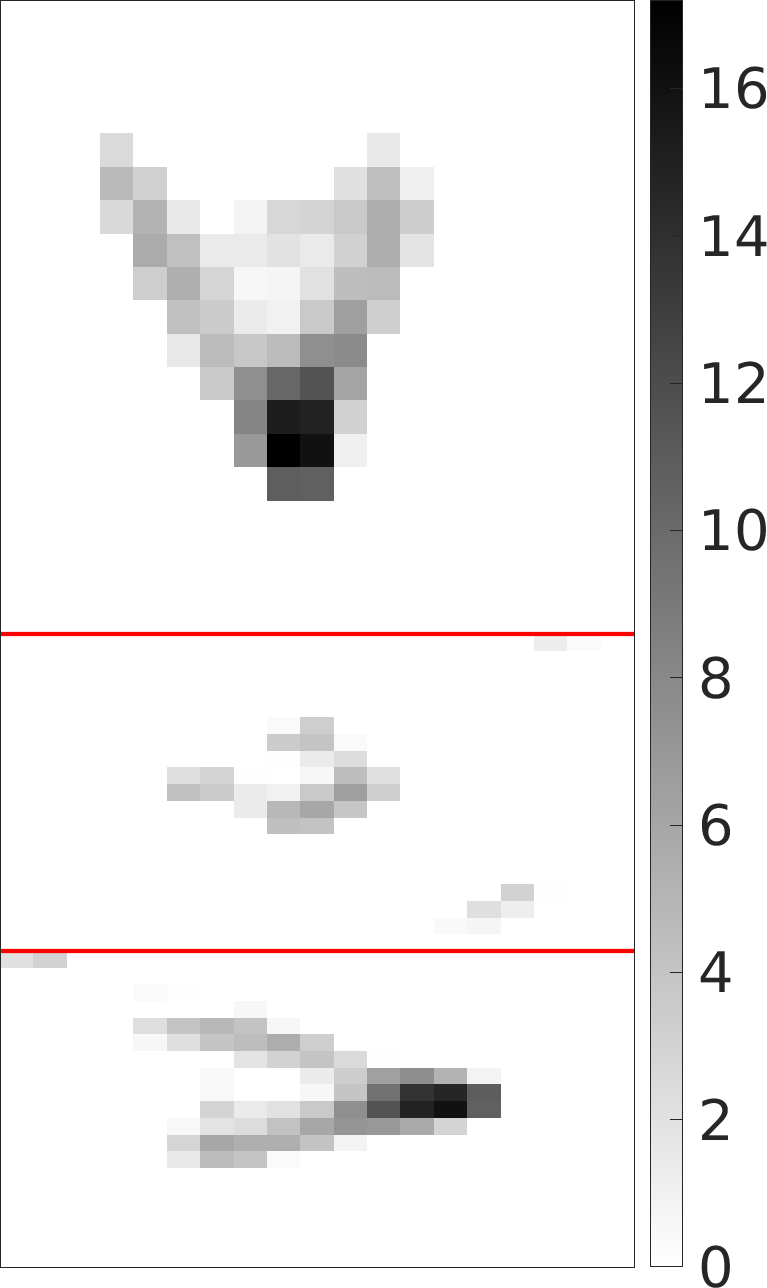}&
 \includegraphics[height=3.4cm]{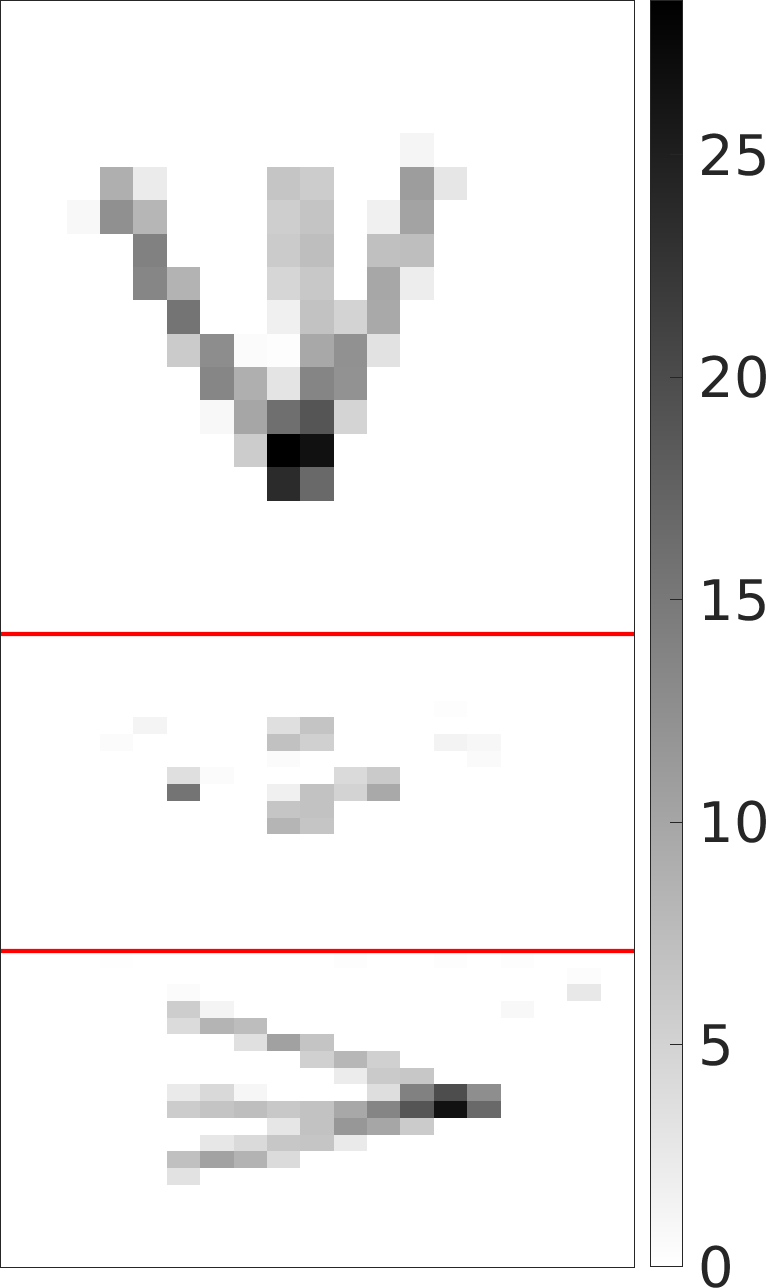}&
 \includegraphics[height=3.4cm]{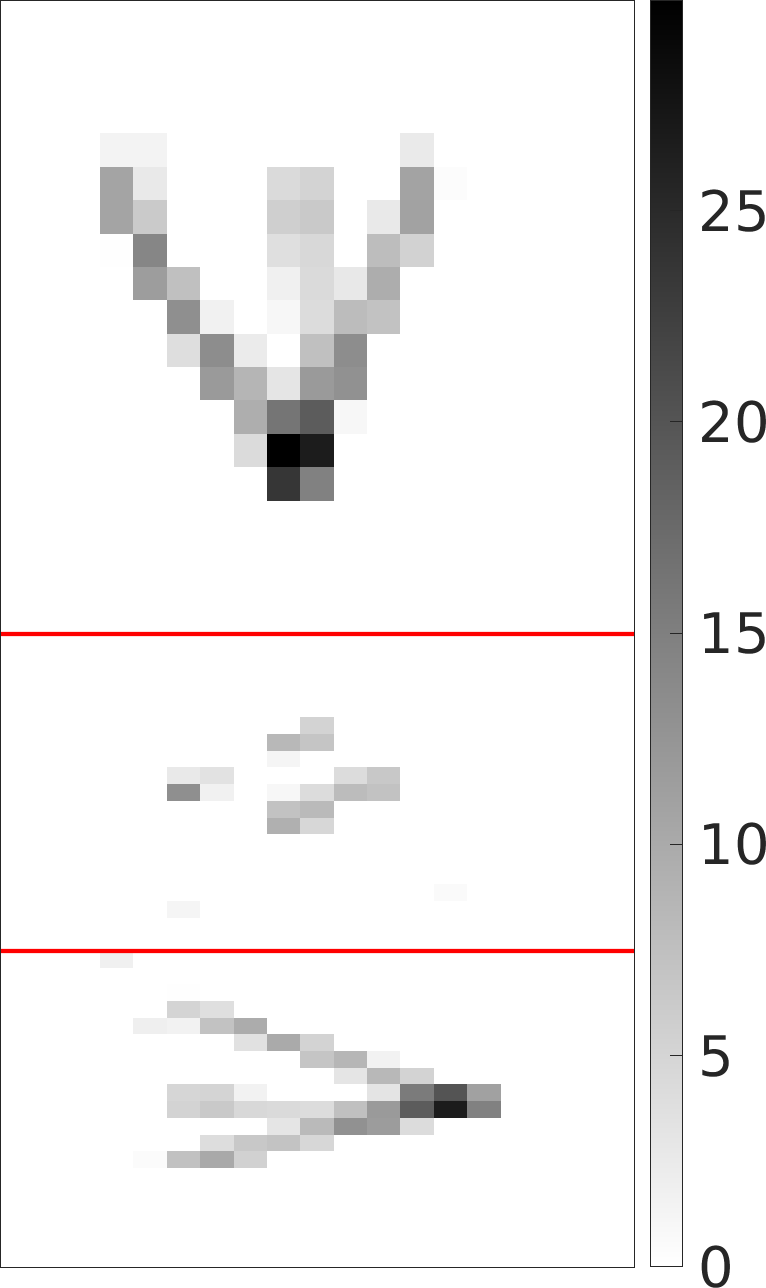}&
 \includegraphics[height=3.4cm]{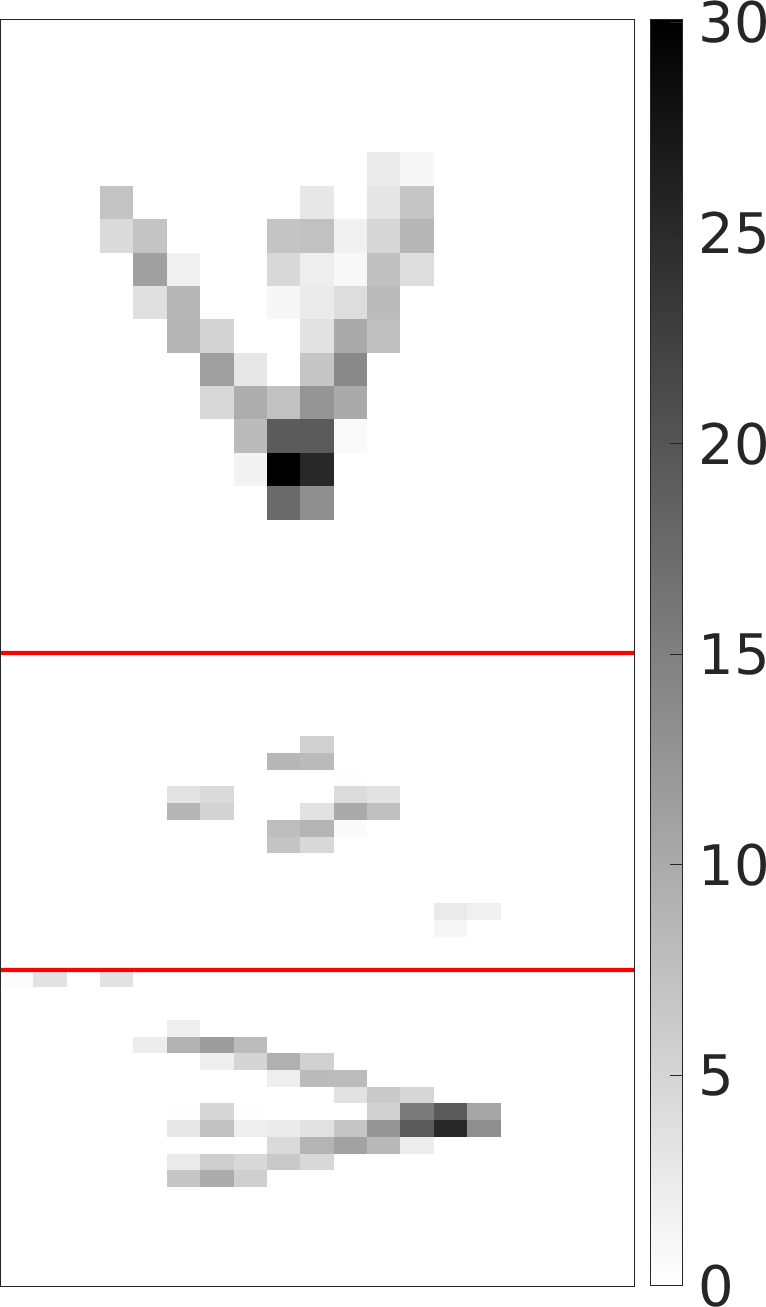}&
 \includegraphics[height=3.4cm]{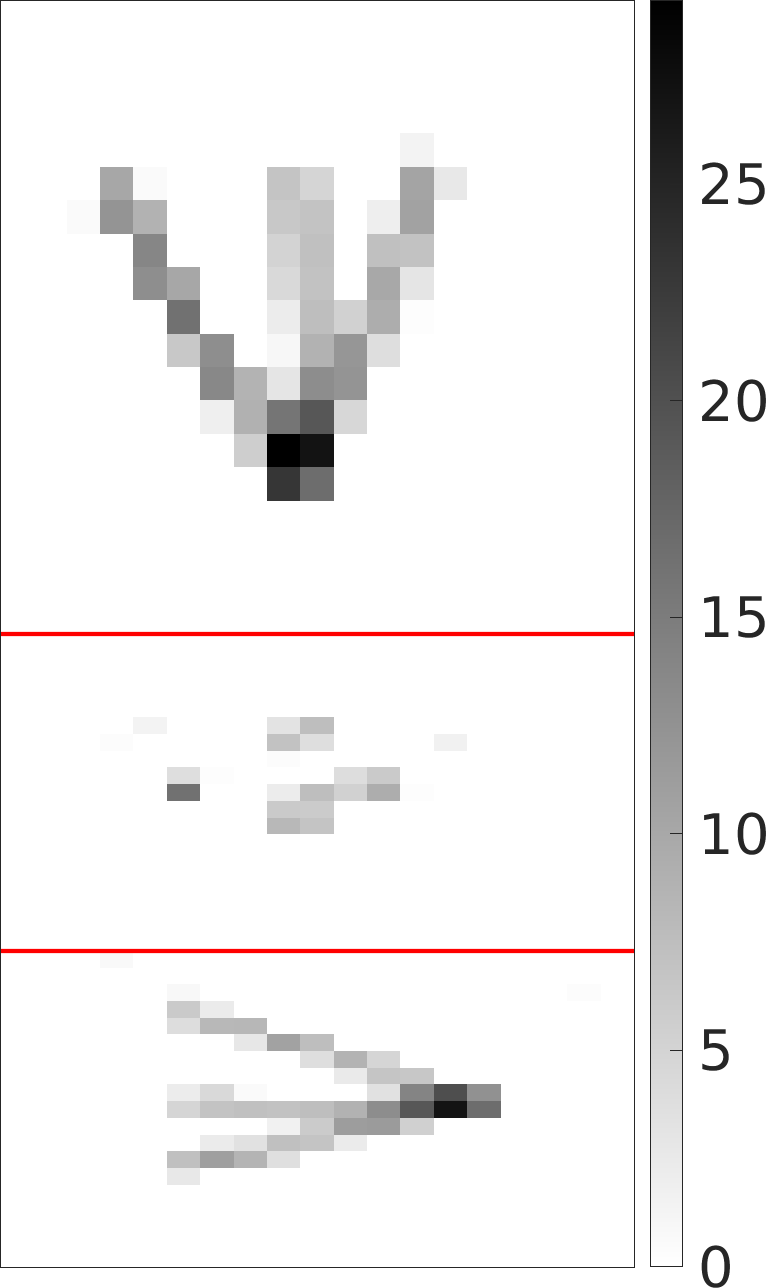}\\
\end{tabular}
}
\caption{Fig. \ref{fig:methods_nonwhitened_vs_whitened_resolution_psnr} with inverted colormap: ``Resolution'' phantom reconstructions, PSNR-optimized $\alpha$ and iteration number $N$ (for l2-K only) according to Table \ref{tab:psnr_nonwhitened_vs_whitened}.}
\label{fig:methods_nonwhitened_vs_whitened_resolution_psnr_inverted_colormap}
\end{figure}

\begin{figure}[hbt!]
\centering
\scalebox{0.85}{
\begin{tabular}{ccc|ccc}
\multicolumn{3}{c|}{non-whitened} & \multicolumn{3}{c}{whitened} \\
\hline
l1-L & l2-L & l2-K & l1-L & l2-L & l2-K \\
\hline
\multicolumn{6}{l}{$\tau=0$} \\
 \includegraphics[height=3.4cm]{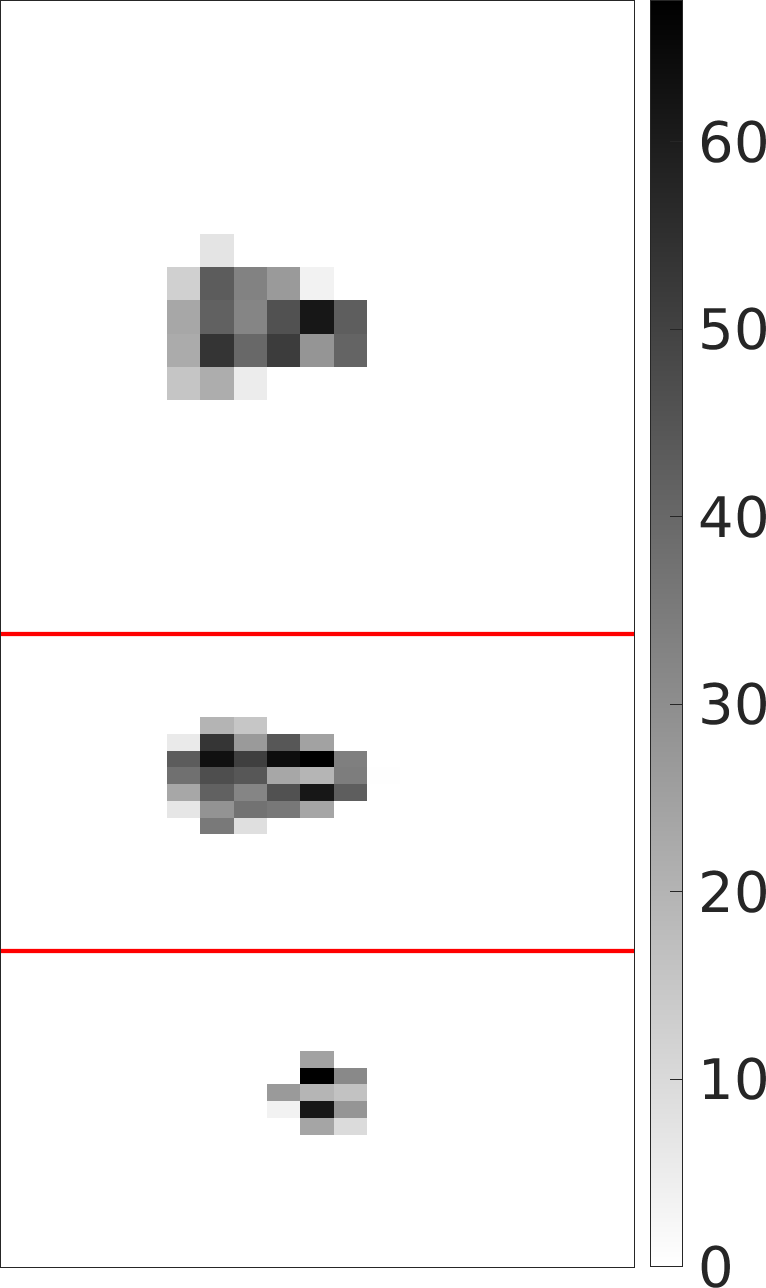}&
 \includegraphics[height=3.4cm]{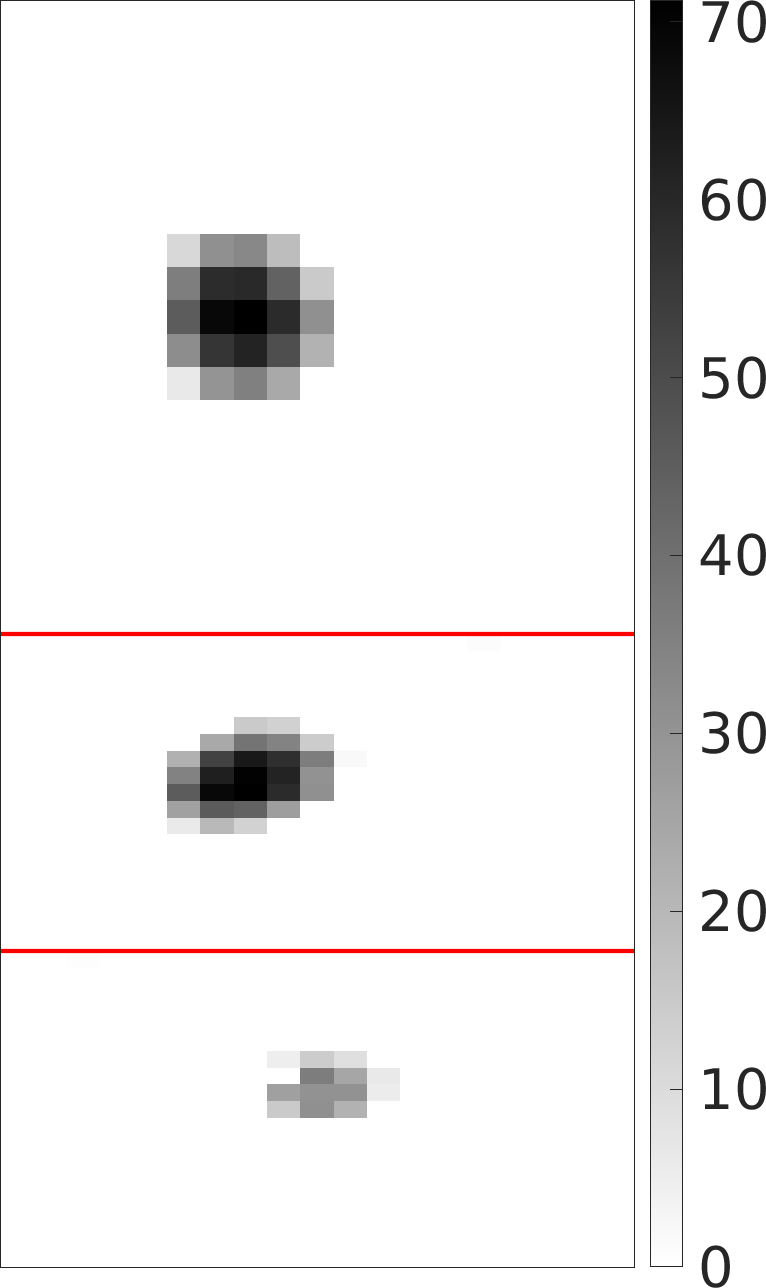}&
 \includegraphics[height=3.4cm]{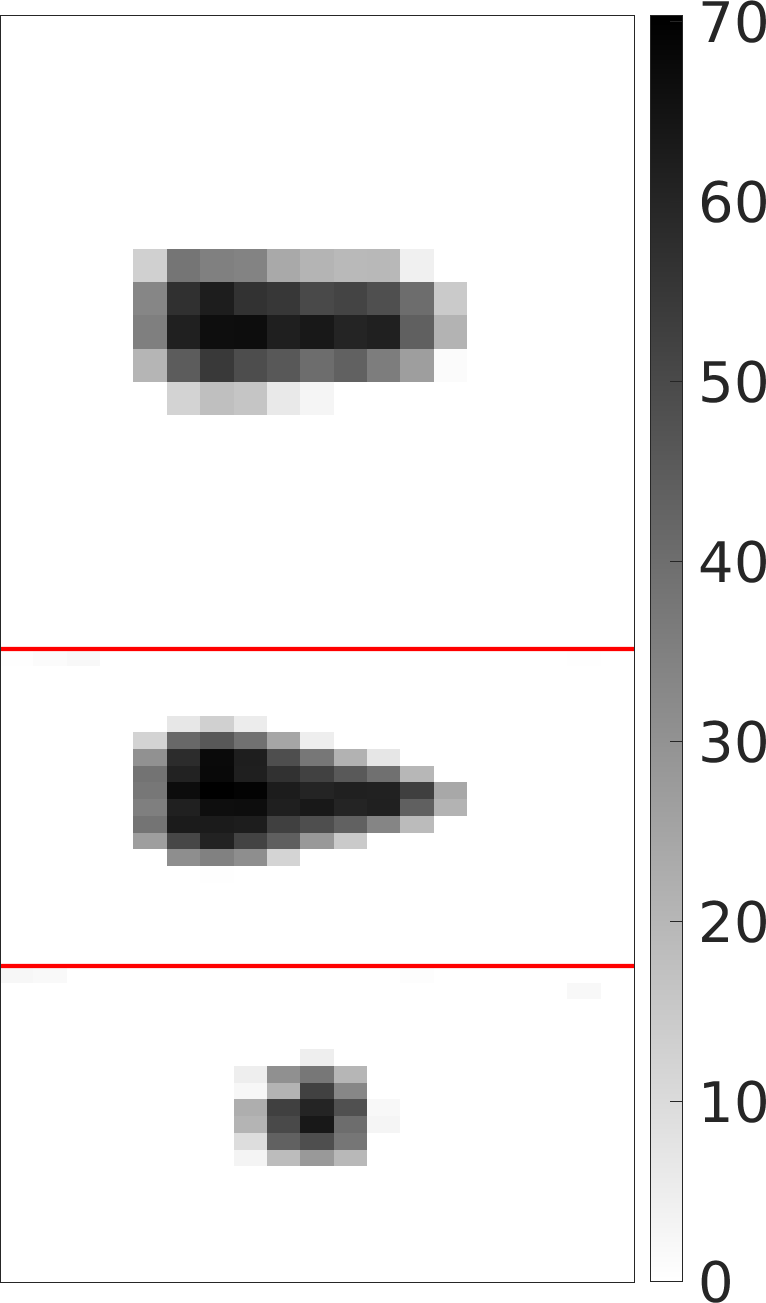}&
 \includegraphics[height=3.4cm]{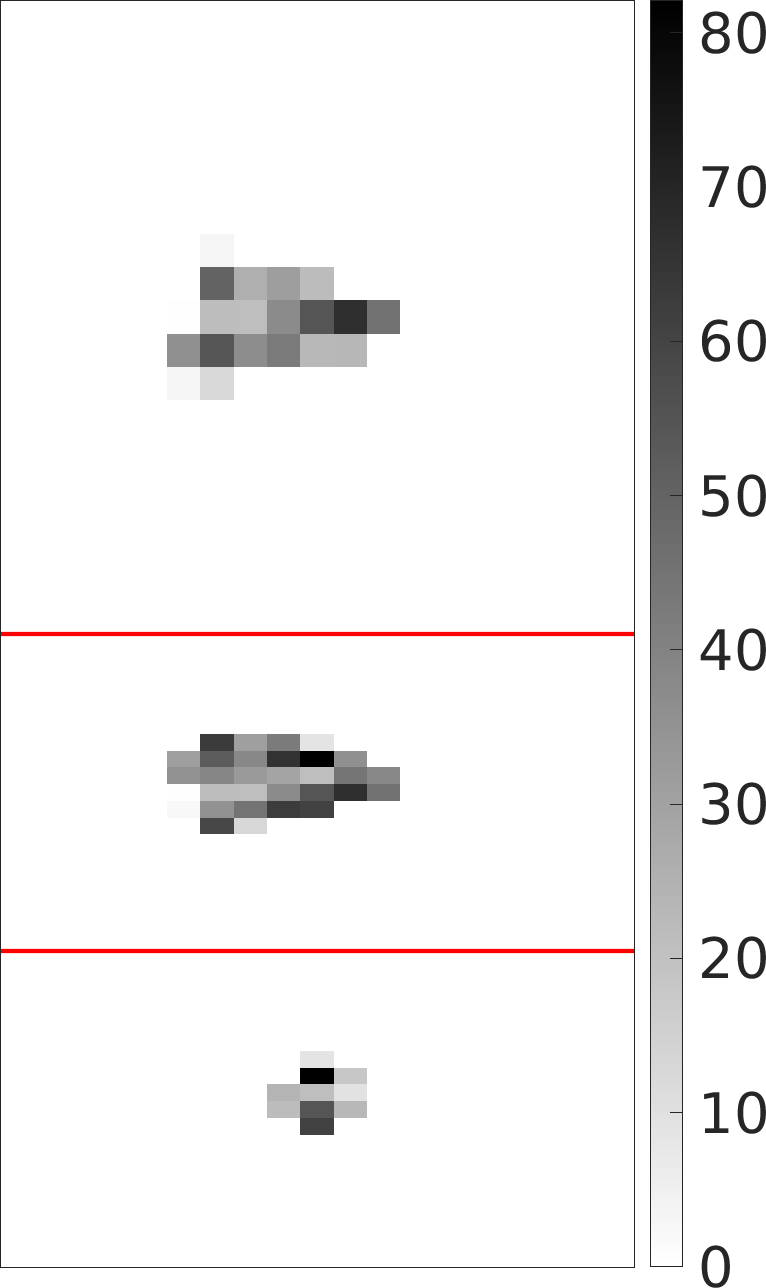}&
 \includegraphics[height=3.4cm]{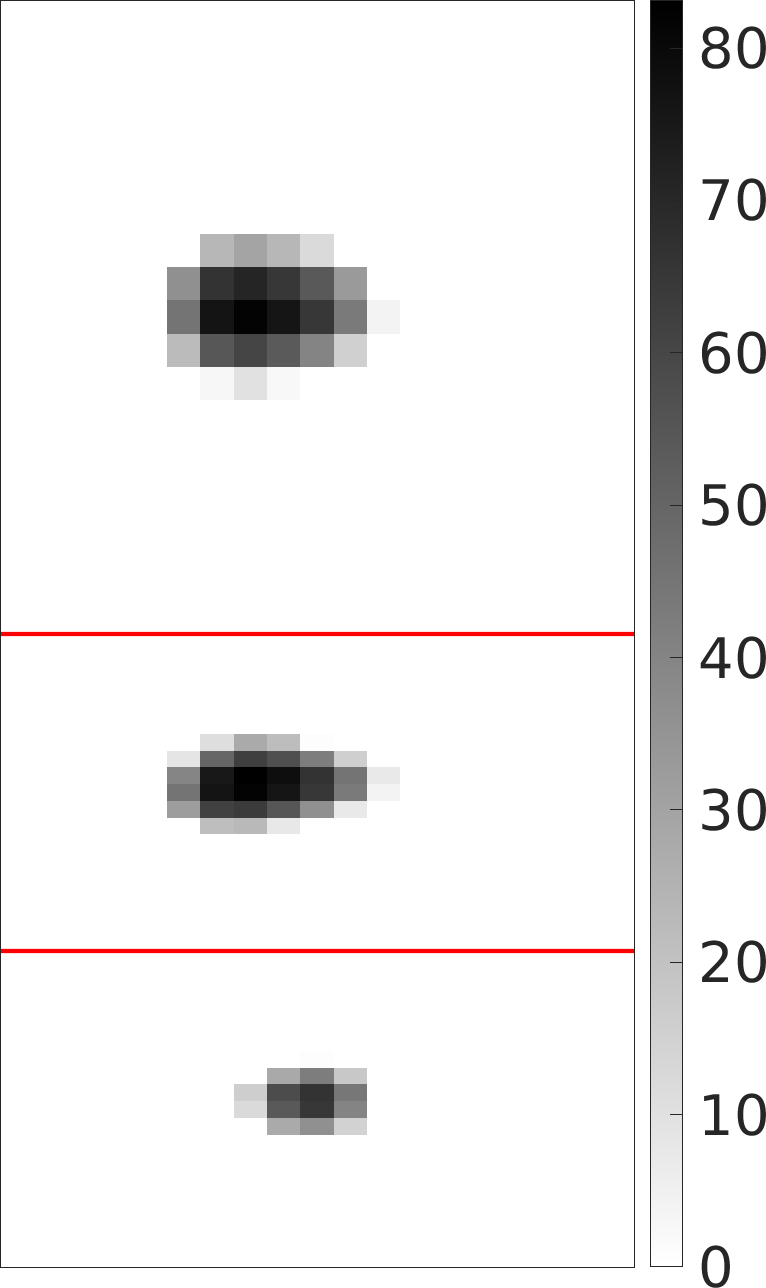}&
 \includegraphics[height=3.4cm]{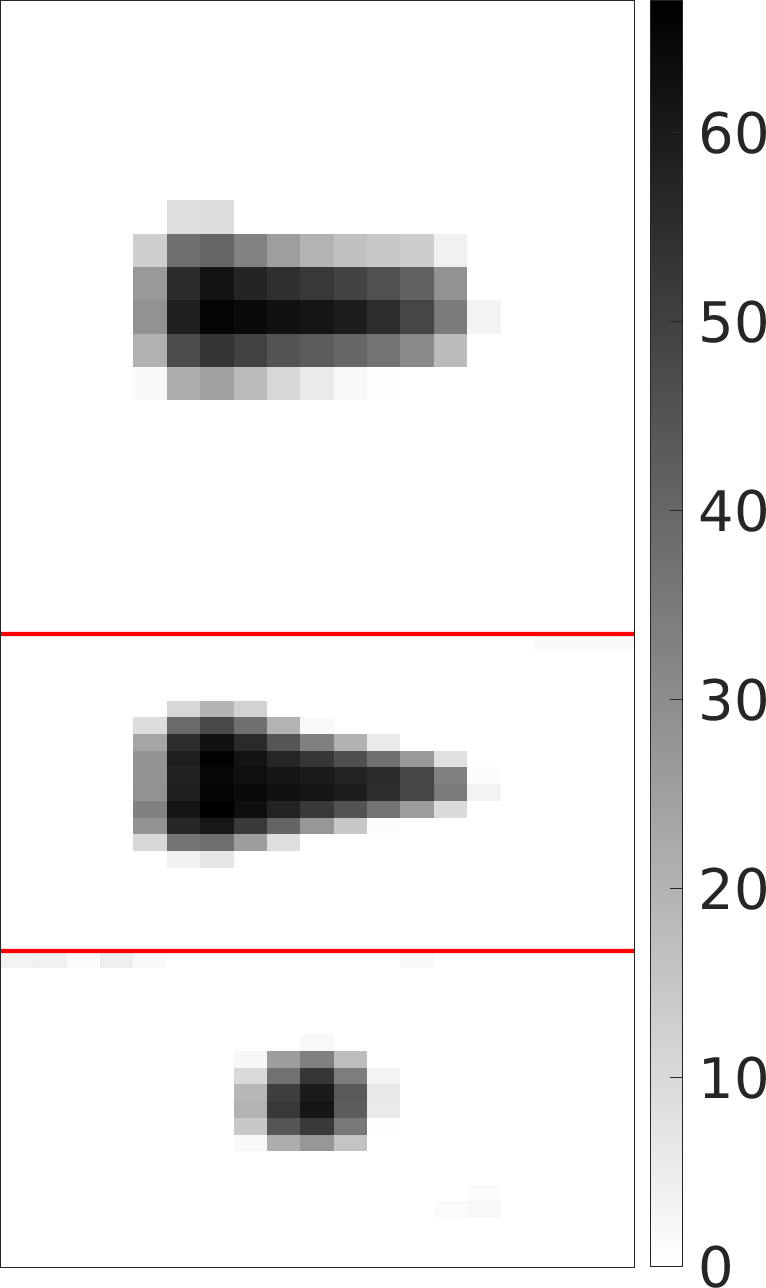}\\
\hline
\multicolumn{6}{l}{$\tau=1$} \\
 \includegraphics[height=3.4cm]{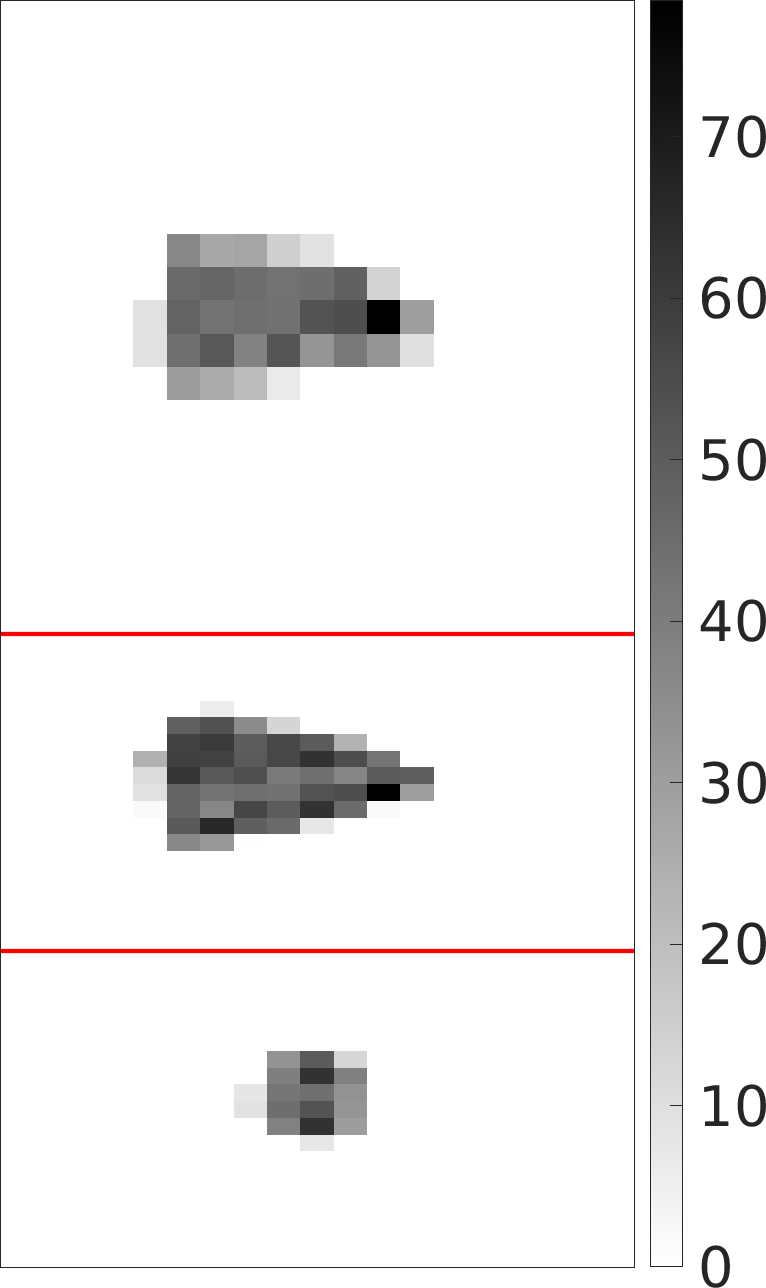}&
 \includegraphics[height=3.4cm]{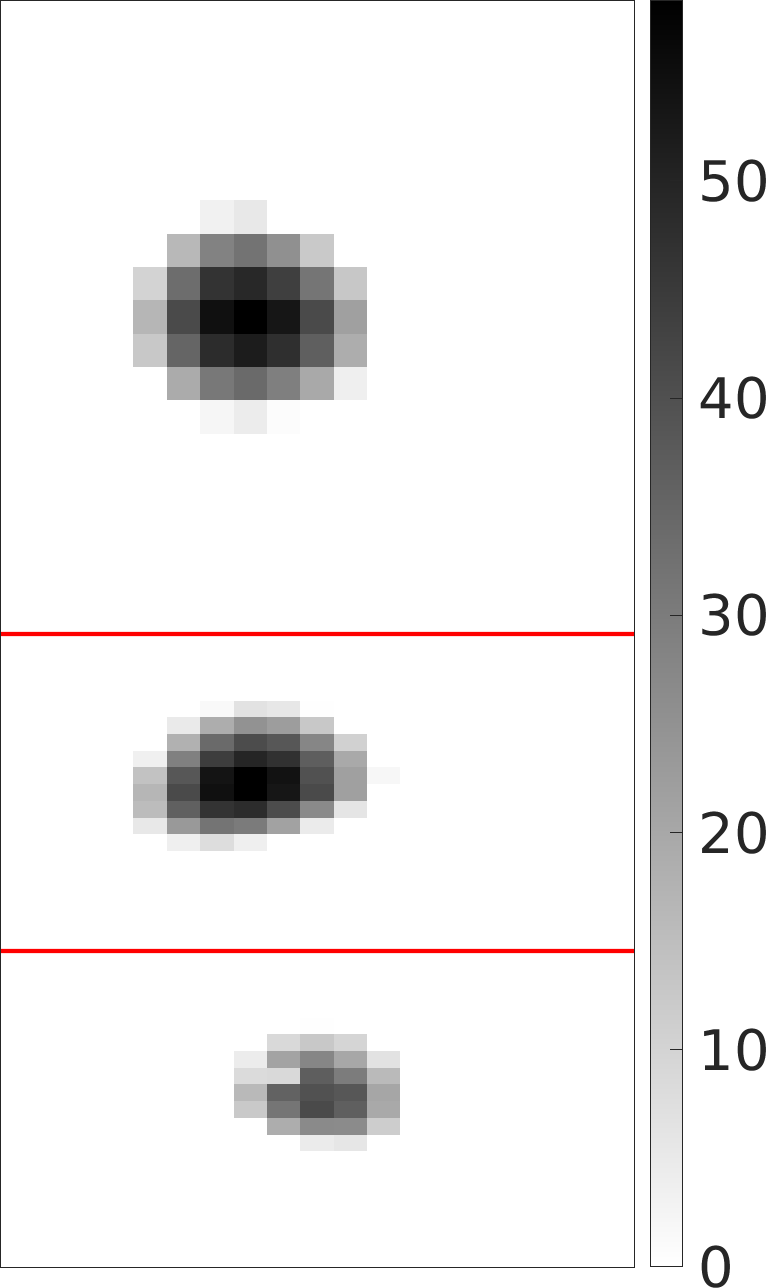}&
 \includegraphics[height=3.4cm]{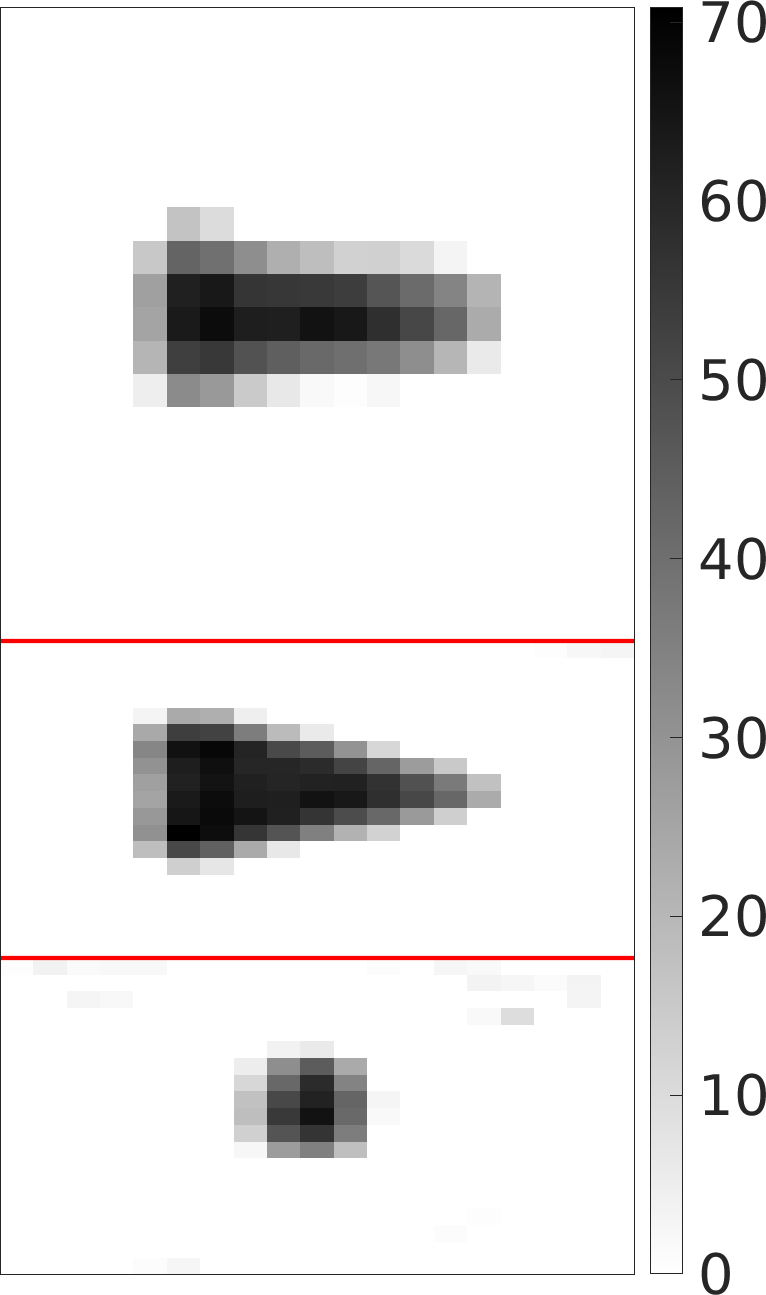}&
 \includegraphics[height=3.4cm]{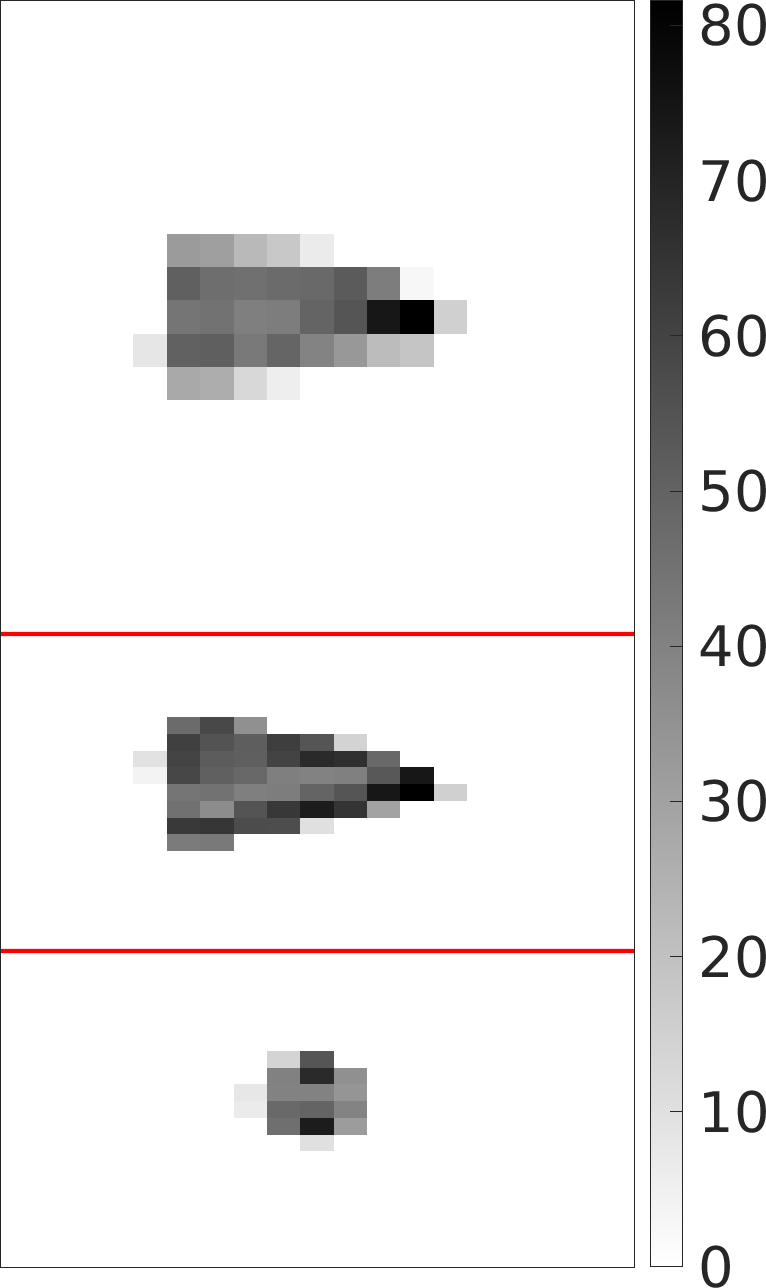}&
 \includegraphics[height=3.4cm]{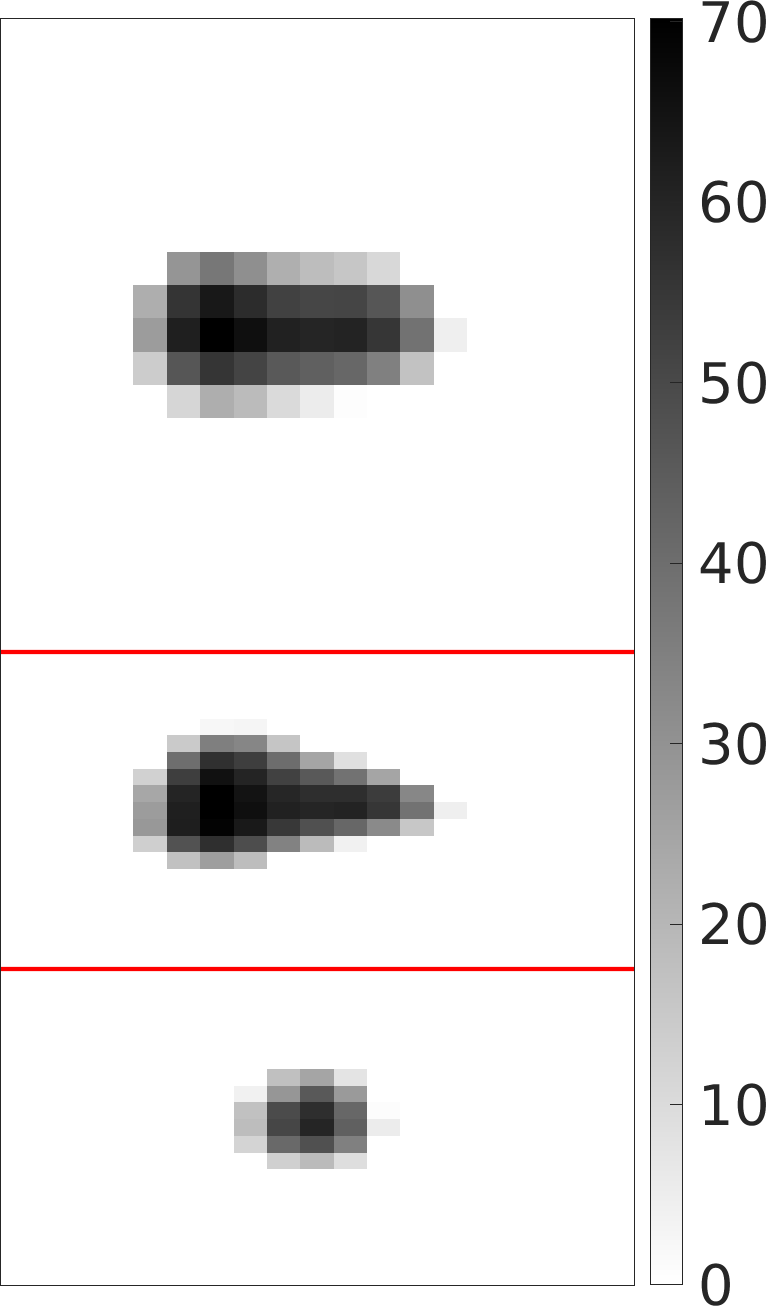}&
 \includegraphics[height=3.4cm]{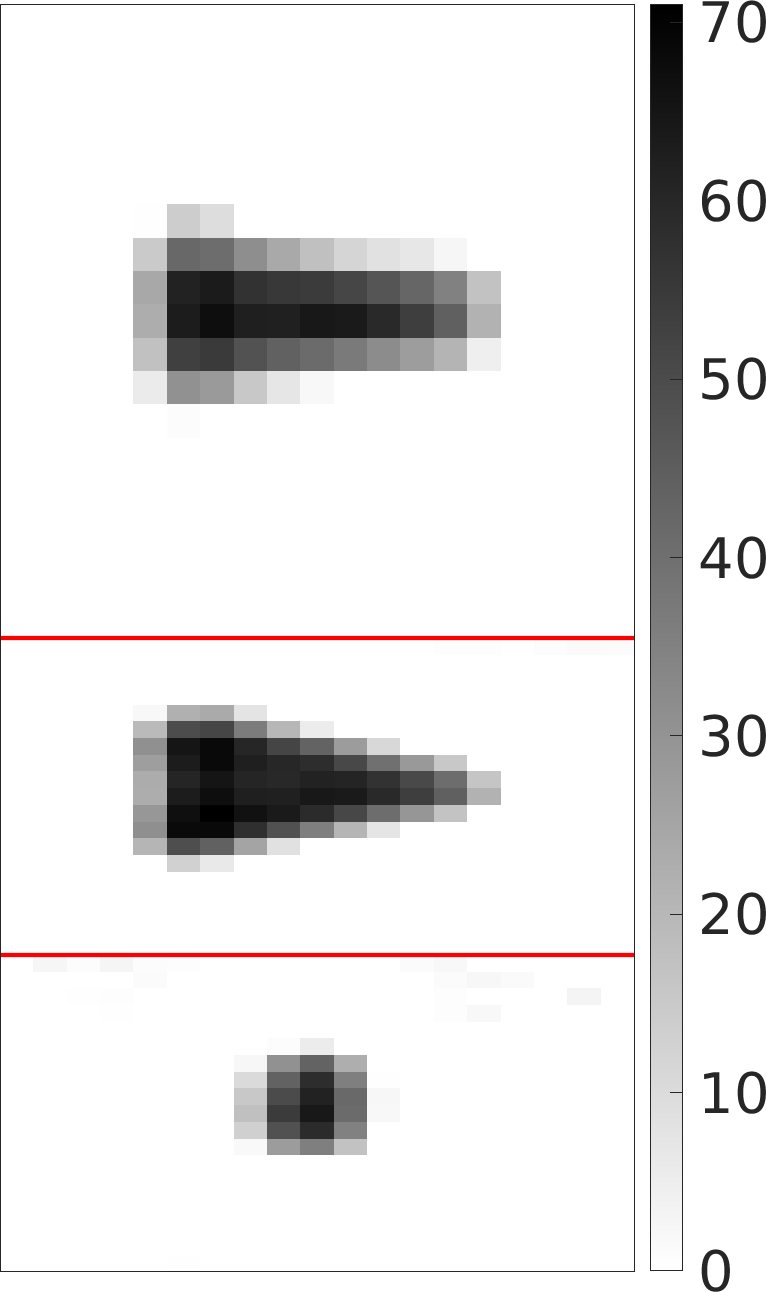}\\
\hline
\multicolumn{6}{l}{$\tau=3$} \\
 \includegraphics[height=3.4cm]{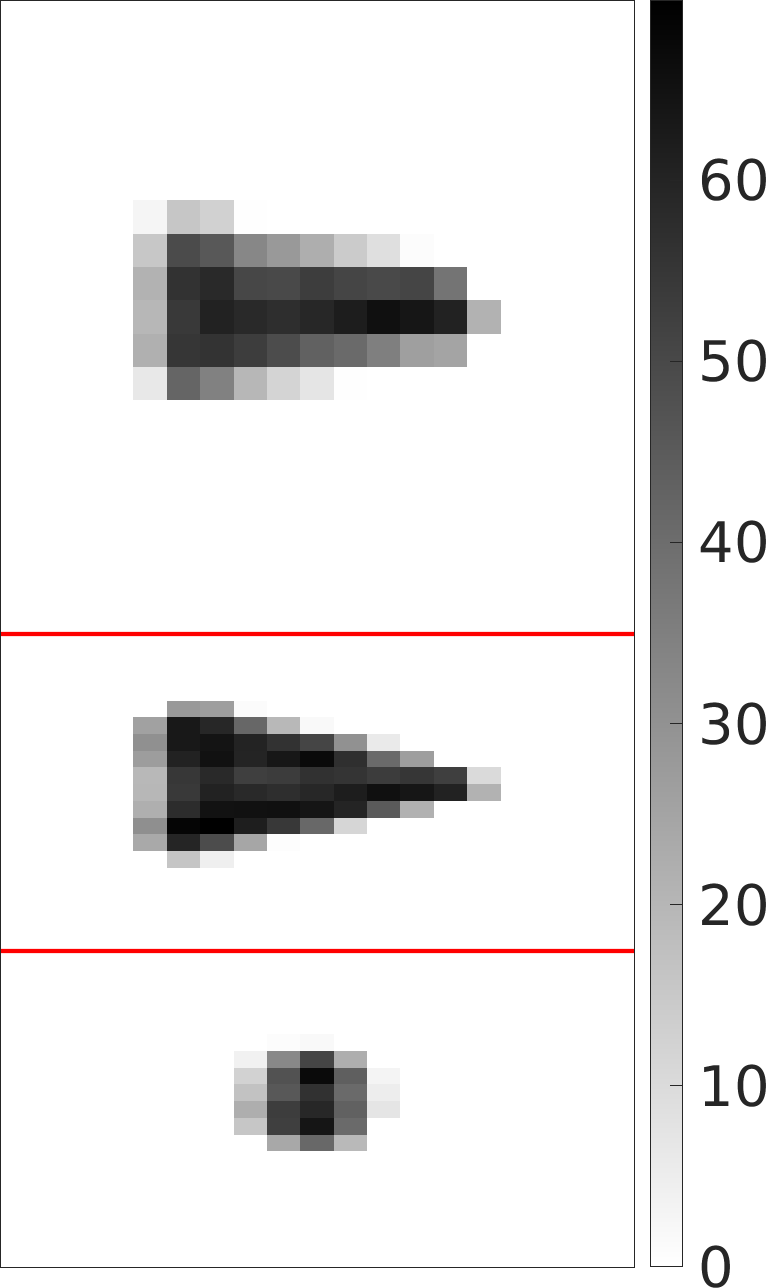}&
 \includegraphics[height=3.4cm]{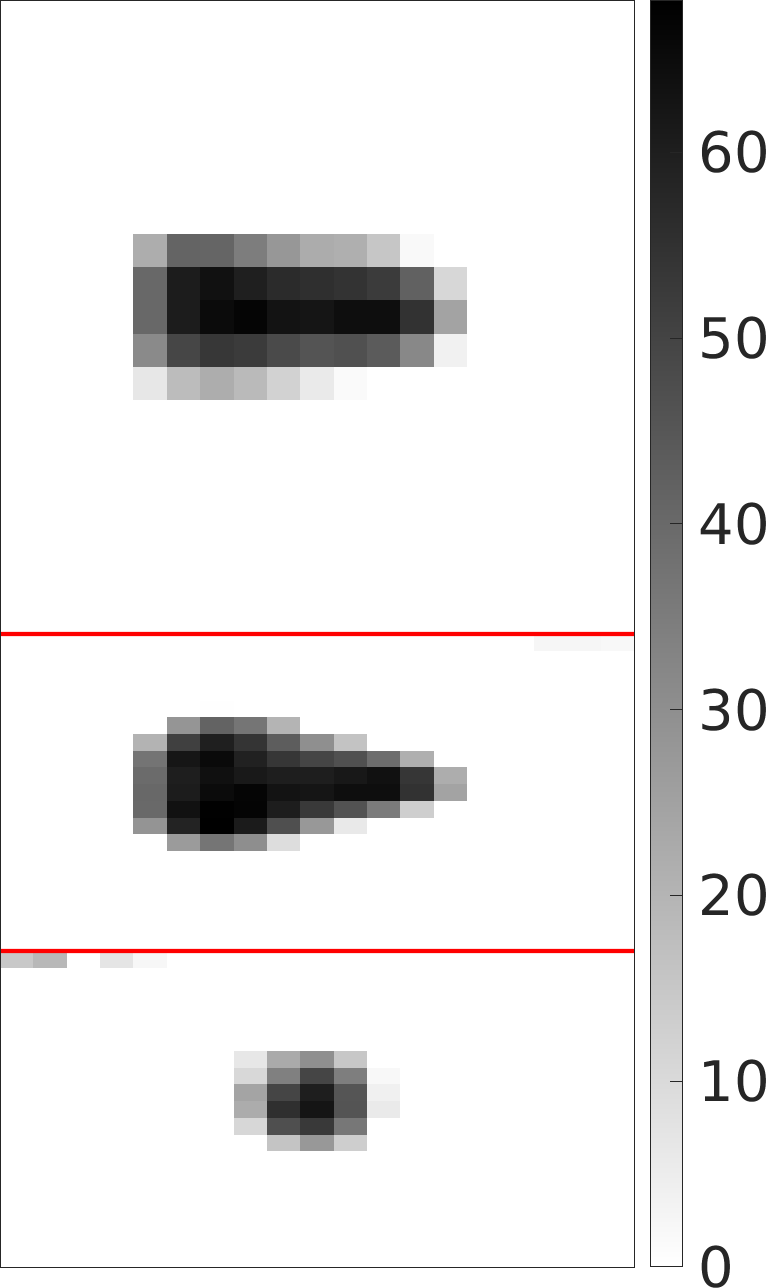}&
 \includegraphics[height=3.4cm]{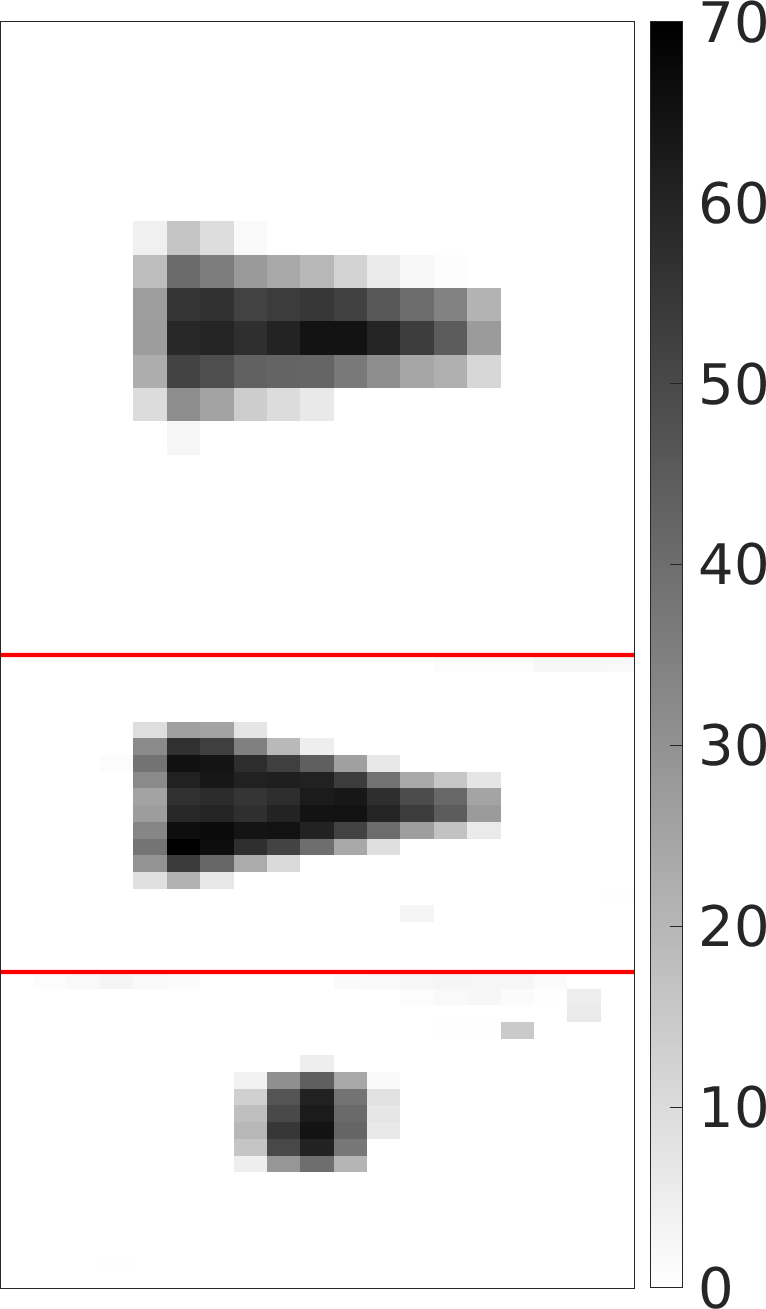}&
 \includegraphics[height=3.4cm]{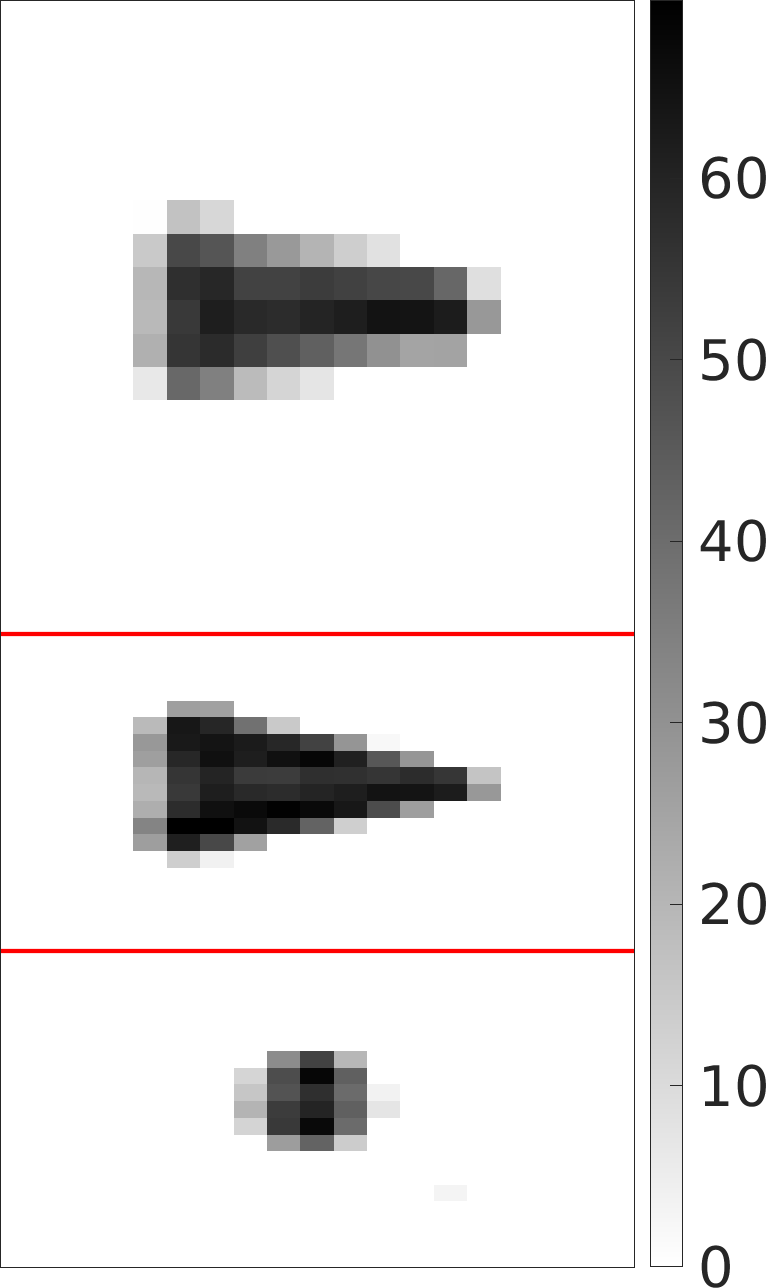}&
 \includegraphics[height=3.4cm]{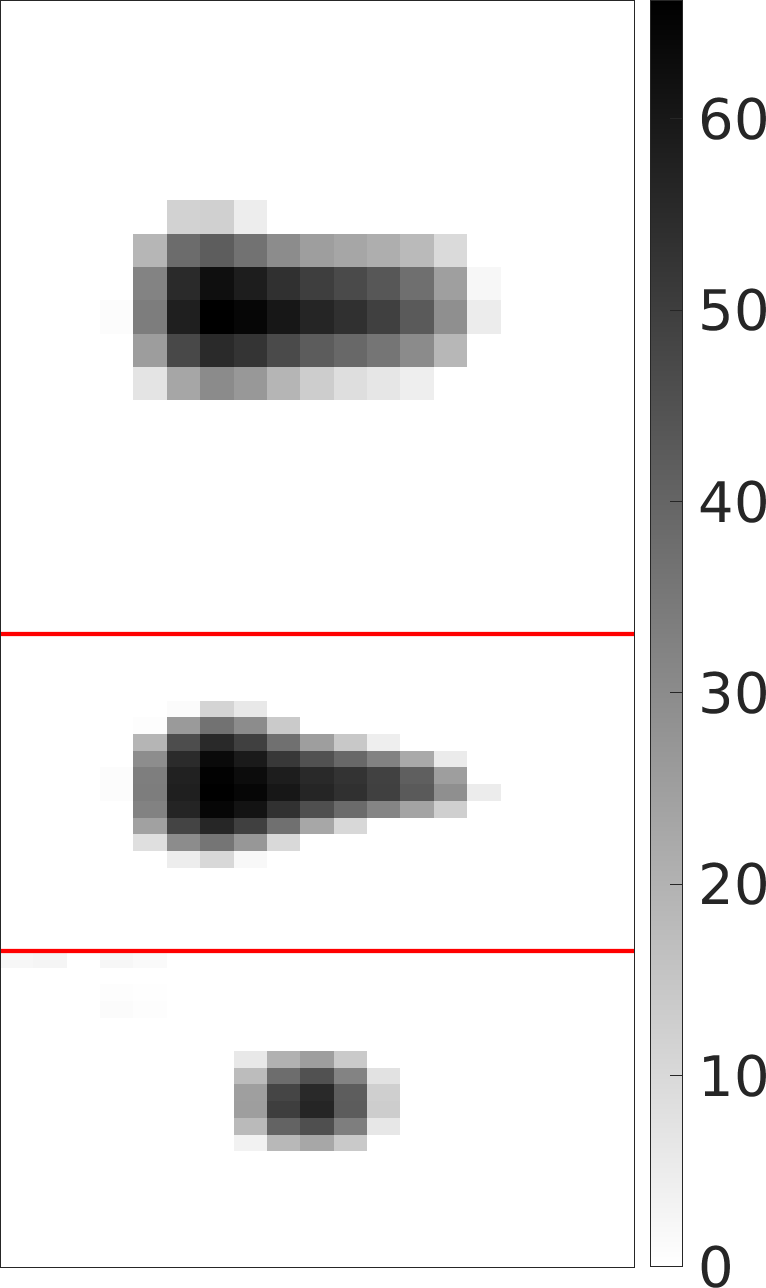}&
 \includegraphics[height=3.4cm]{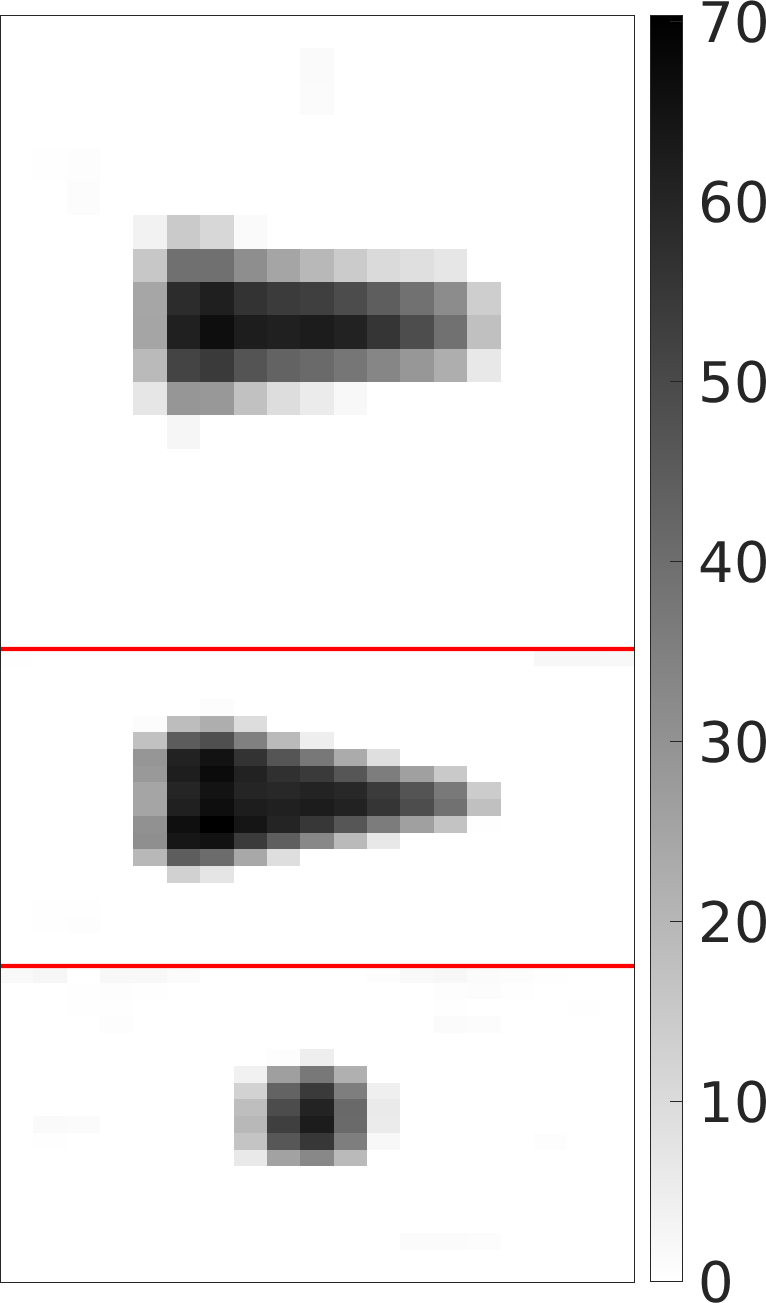}\\
 \hline
\multicolumn{6}{l}{$\tau=5$} \\
 \includegraphics[height=3.4cm]{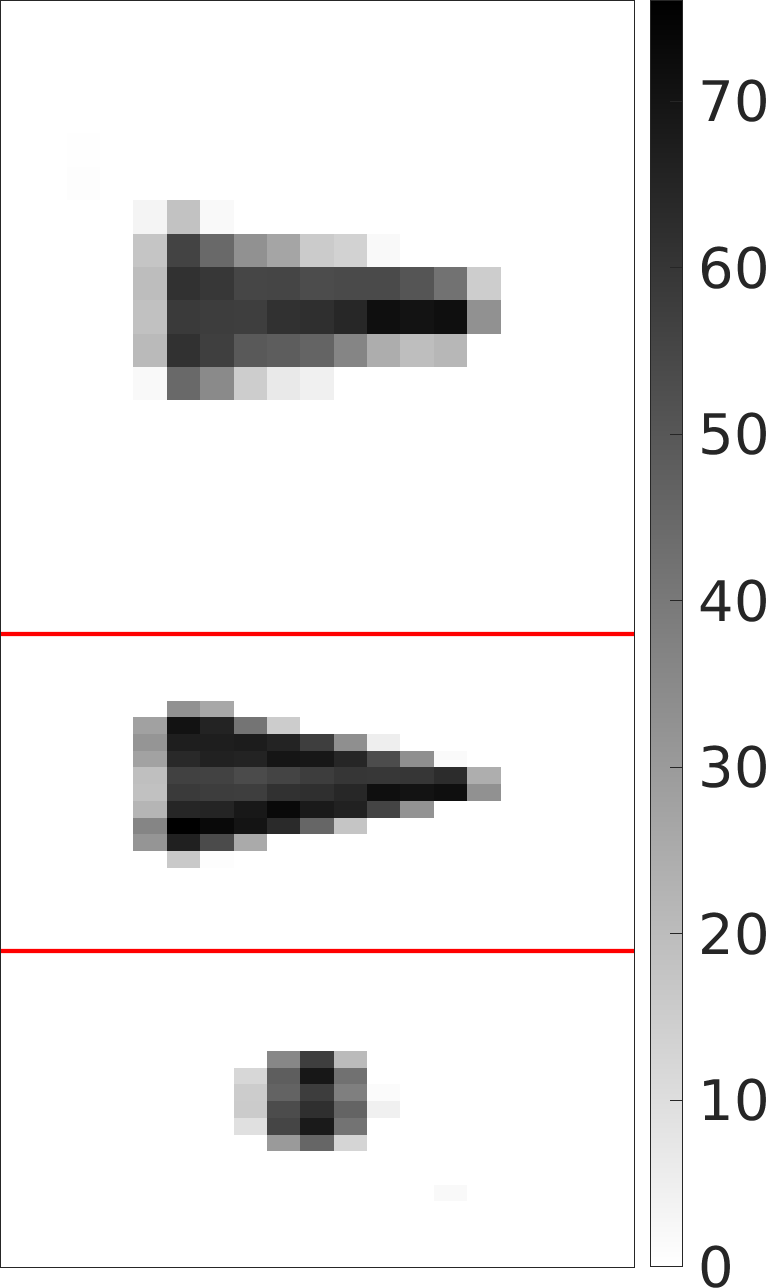}&
 \includegraphics[height=3.4cm]{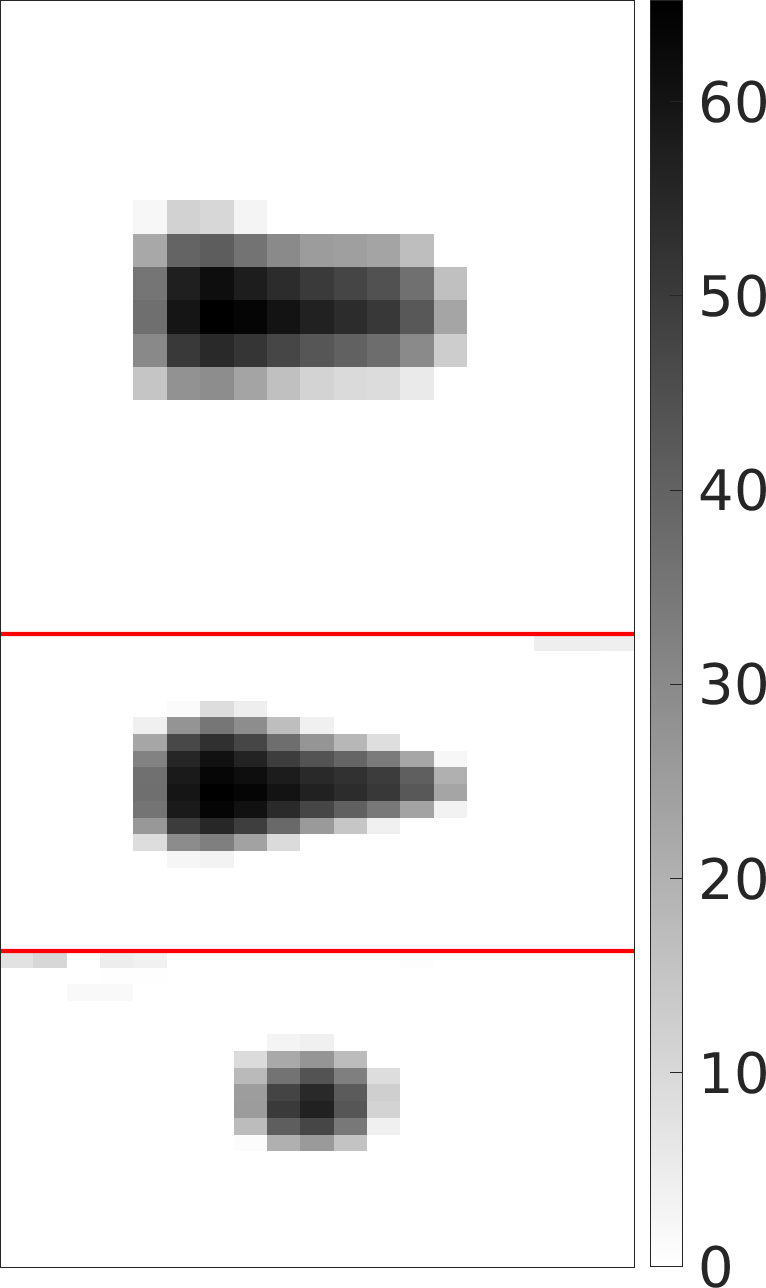}&
 \includegraphics[height=3.4cm]{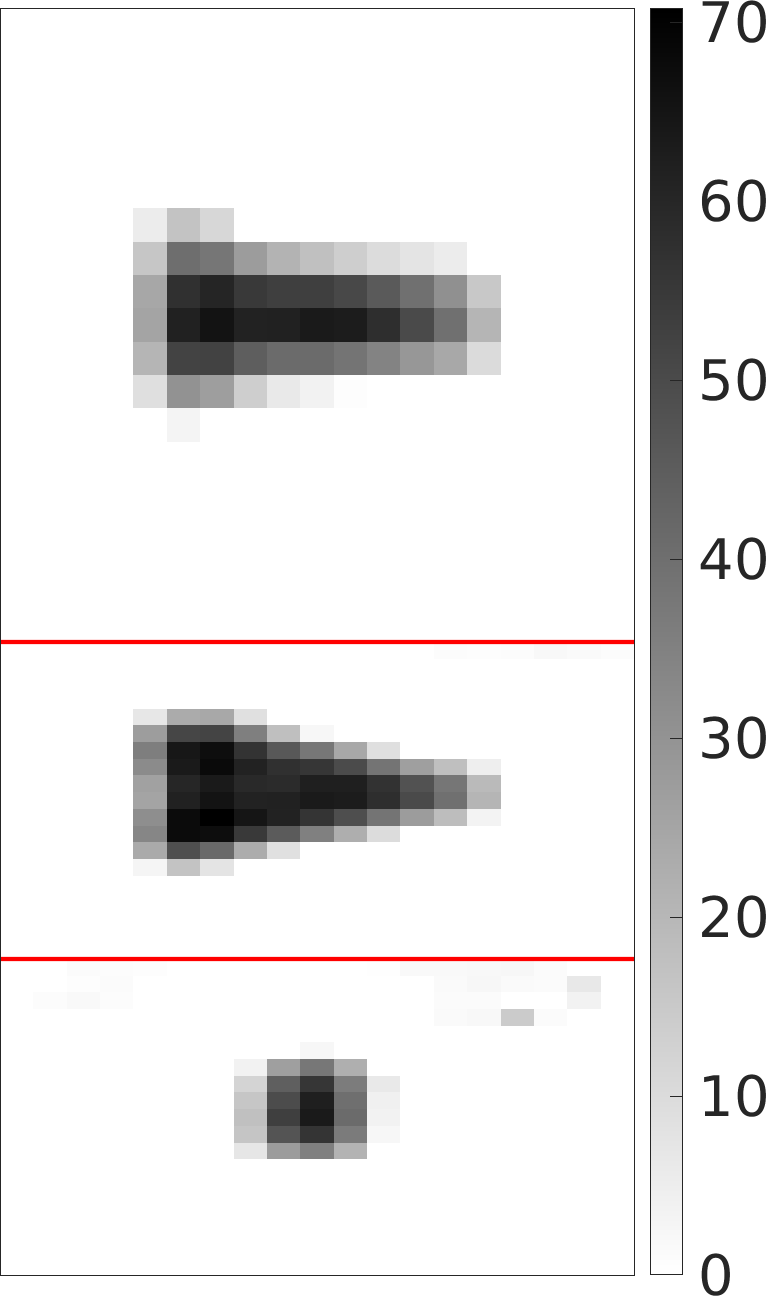}&
 \includegraphics[height=3.4cm]{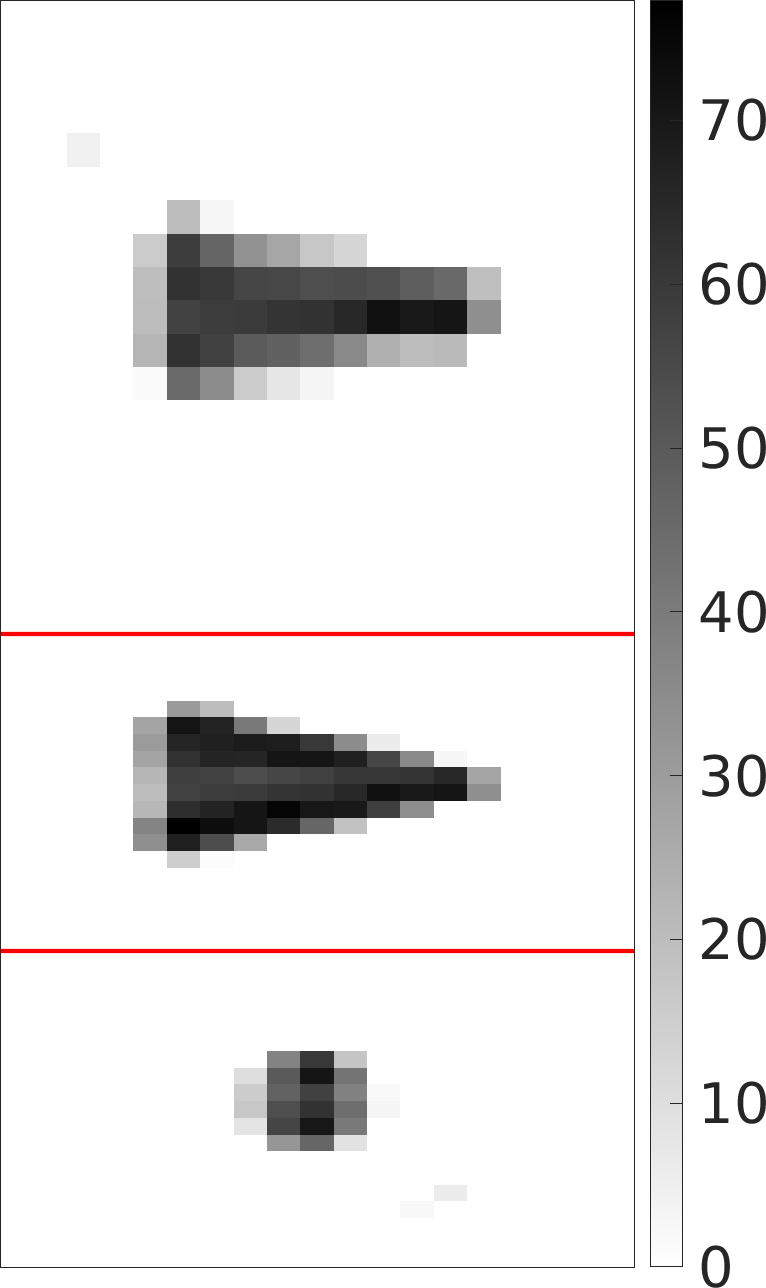}&
 \includegraphics[height=3.4cm]{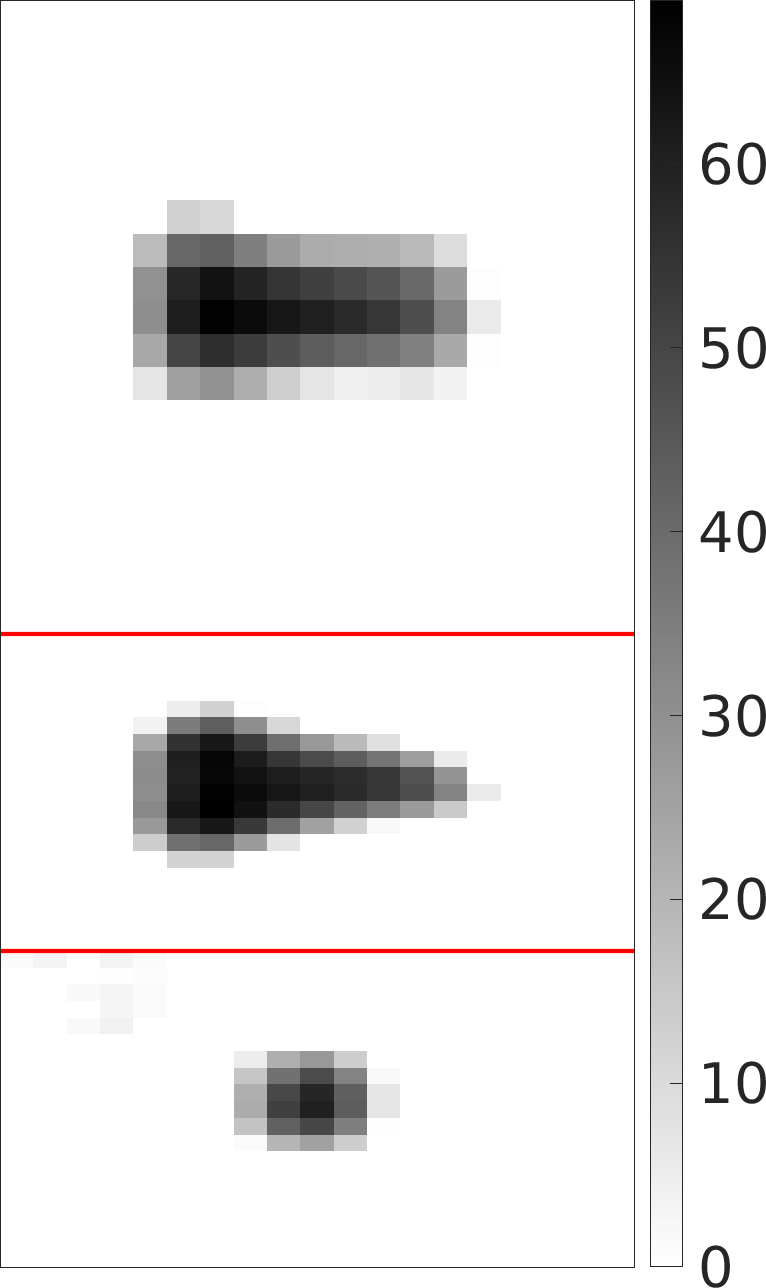}&
 \includegraphics[height=3.4cm]{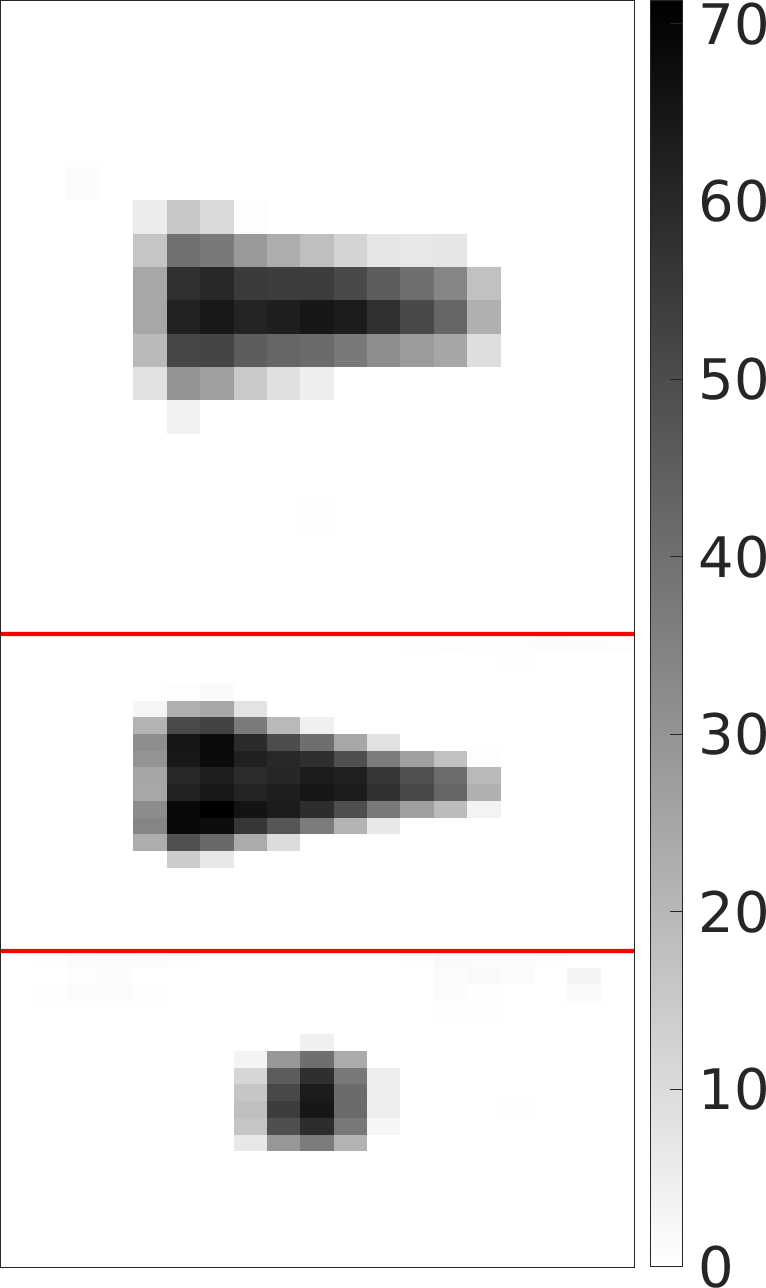}\\
\end{tabular}
}
\caption{Fig. \ref{fig:methods_nonwhitened_vs_whitened_shape_ssim} with inverted colormap: ``Shape'' phantom reconstructions, SSIM-optimized $\alpha$ and iteration number $N$ (for l2-K only) according to Table \ref{tab:ssim_nonwhitened_vs_whitened}. }
\label{fig:methods_nonwhitened_vs_whitened_shape_ssim_inverted_colormap}
\end{figure}

\begin{figure}[hbt!]%
\centering
\scalebox{0.85}{
\begin{tabular}{ccc|ccc}
\multicolumn{3}{c|}{non-whitened} & \multicolumn{3}{c}{whitened} \\
\hline
l1-L & l2-L & l2-K & l1-L & l2-L & l2-K \\
\hline
\multicolumn{6}{l}{$\tau=0$} \\
 \includegraphics[height=3.4cm]{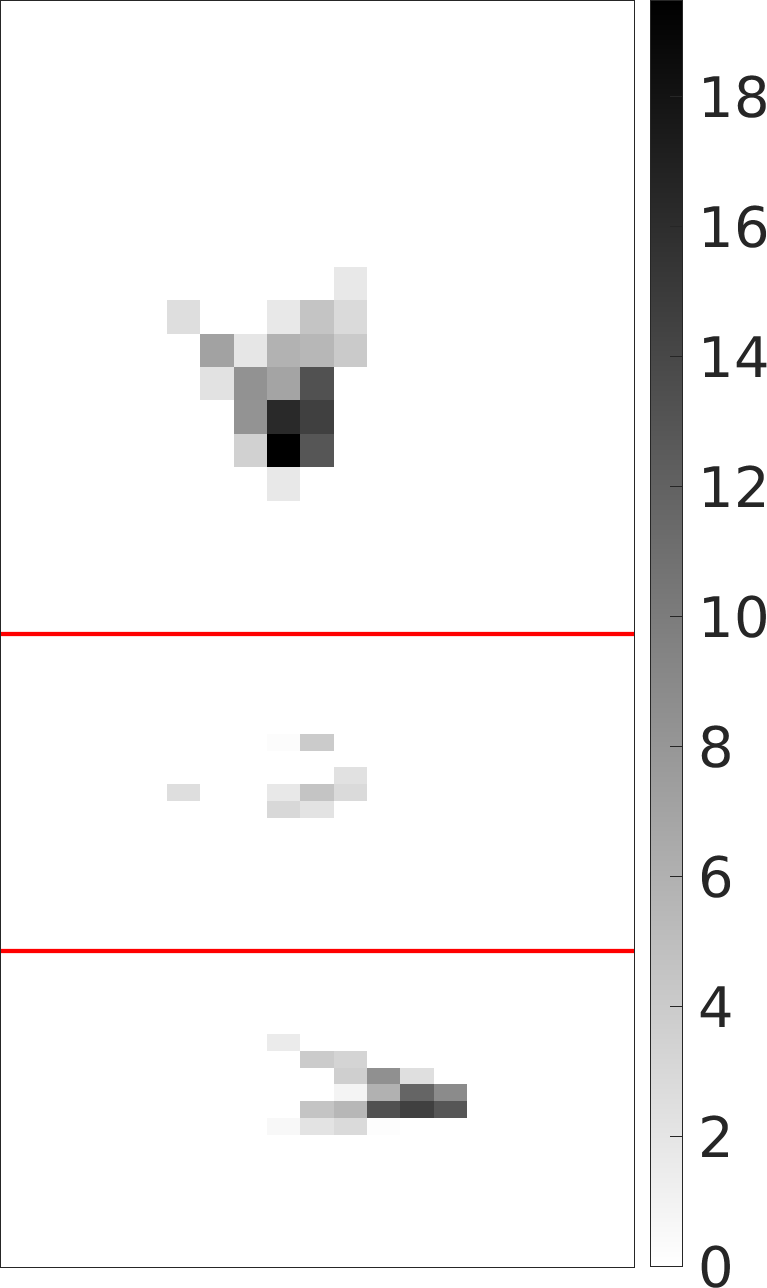}&
 \includegraphics[height=3.4cm]{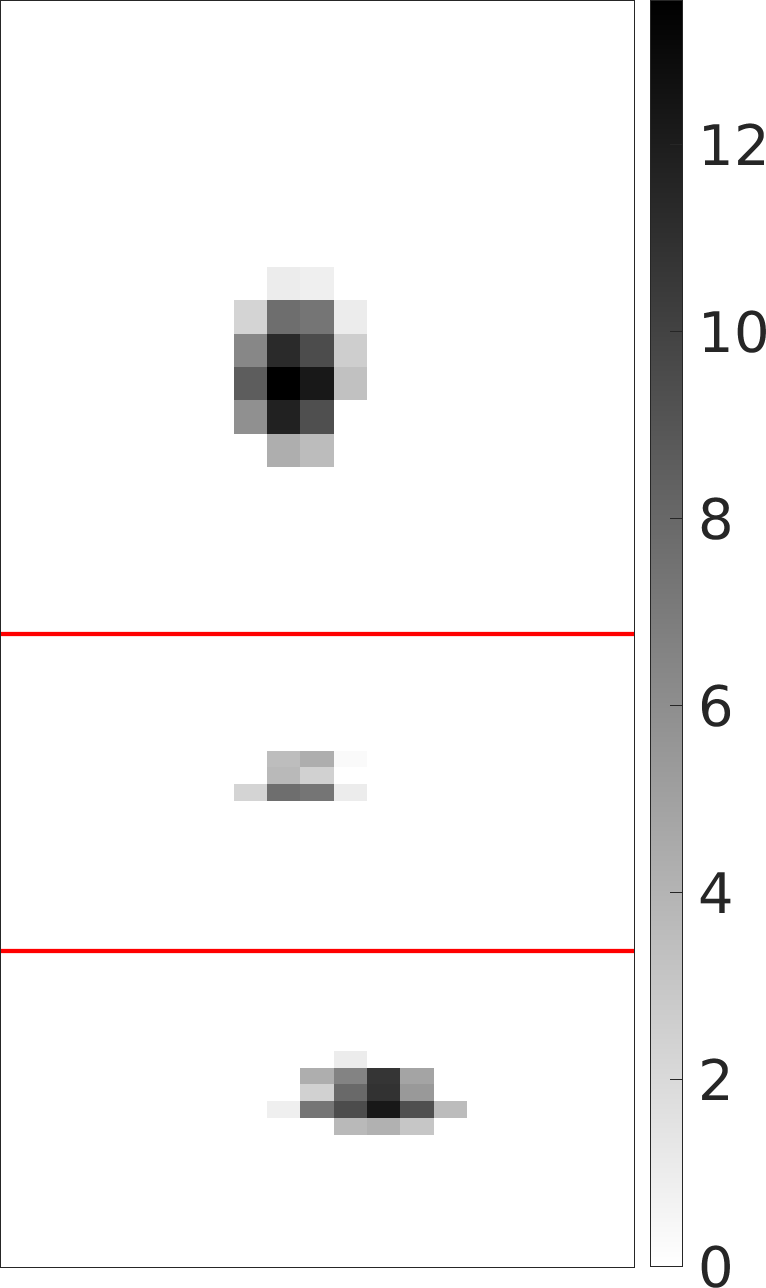}&
 \includegraphics[height=3.4cm]{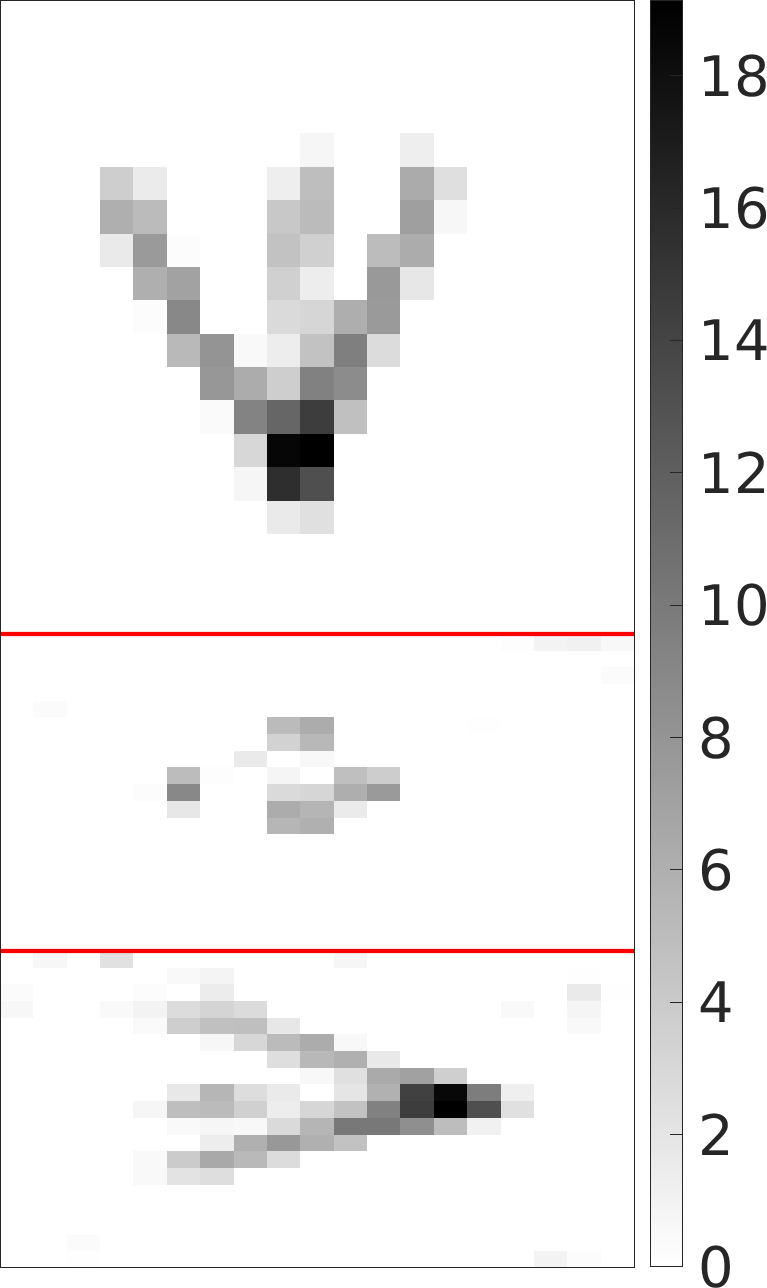}&
 \includegraphics[height=3.4cm]{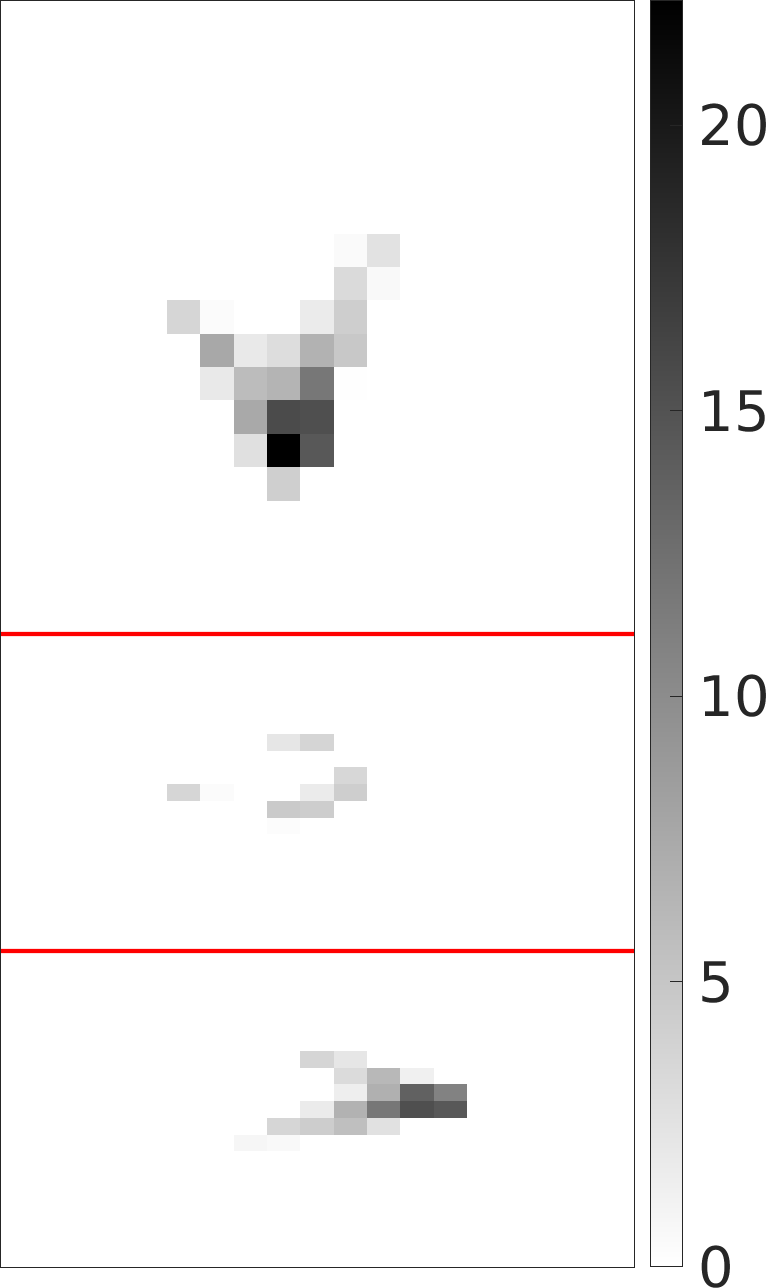}&
 \includegraphics[height=3.4cm]{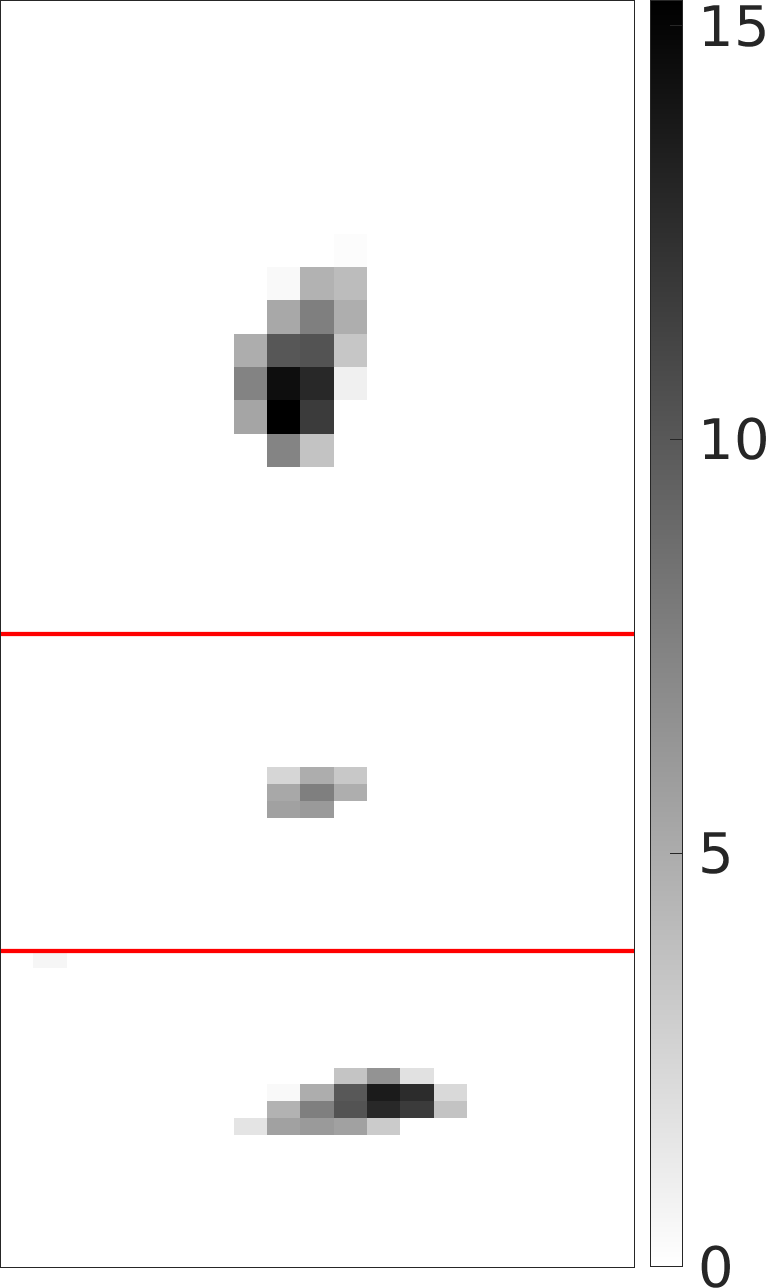}&
 \includegraphics[height=3.4cm]{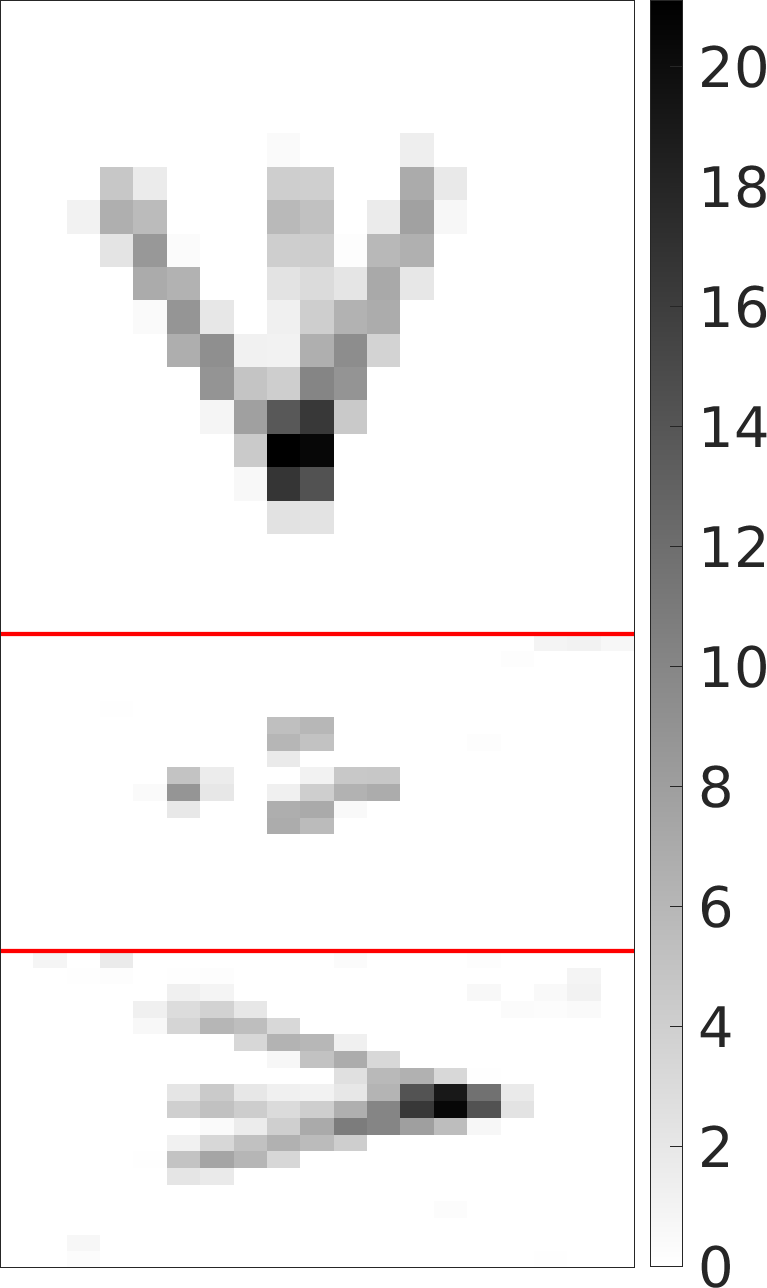}\\
\hline
\multicolumn{6}{l}{$\tau=1$} \\
 \includegraphics[height=3.4cm]{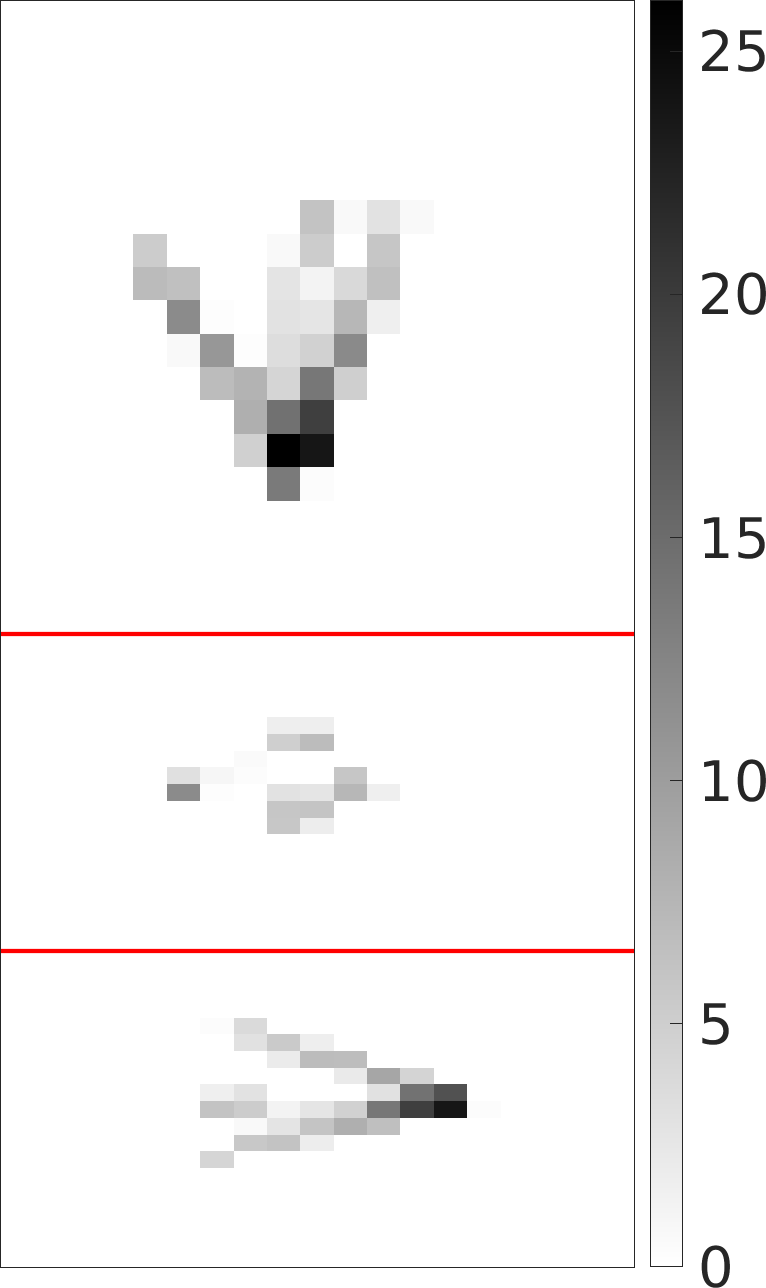}&
 \includegraphics[height=3.4cm]{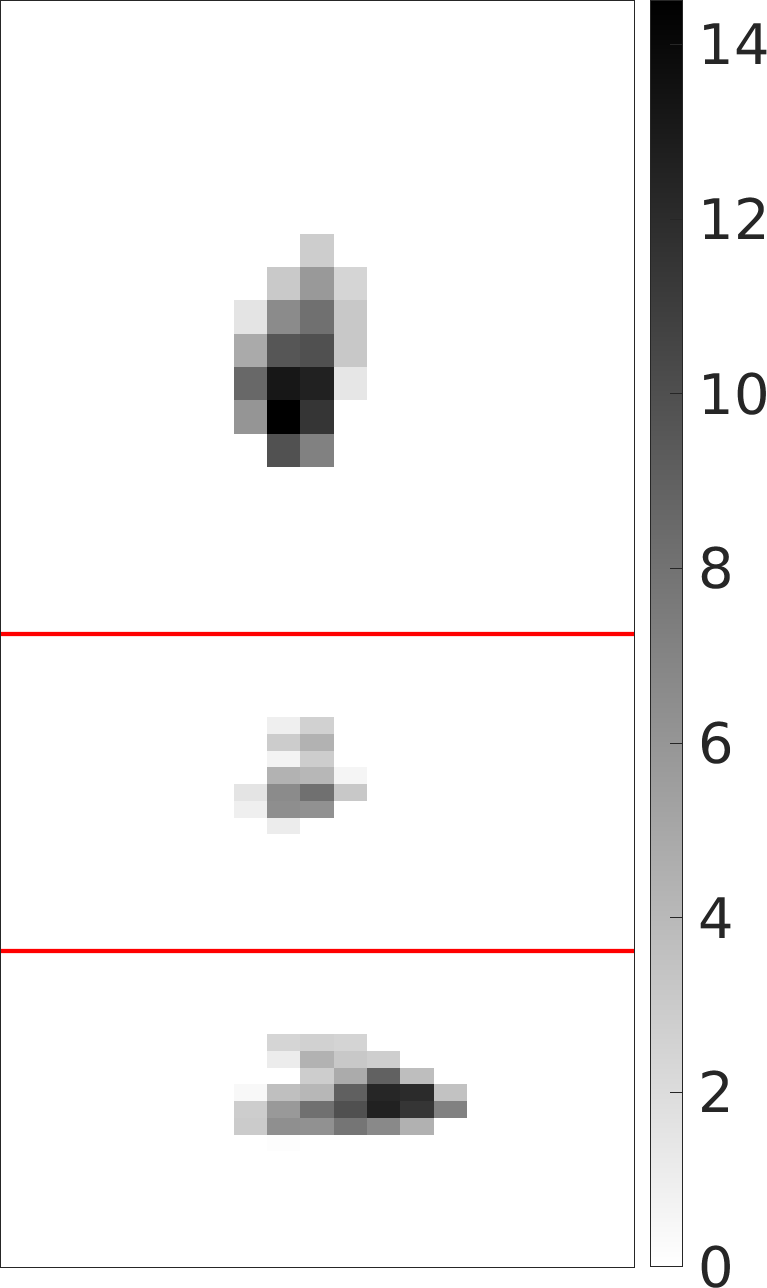}&
 \includegraphics[height=3.4cm]{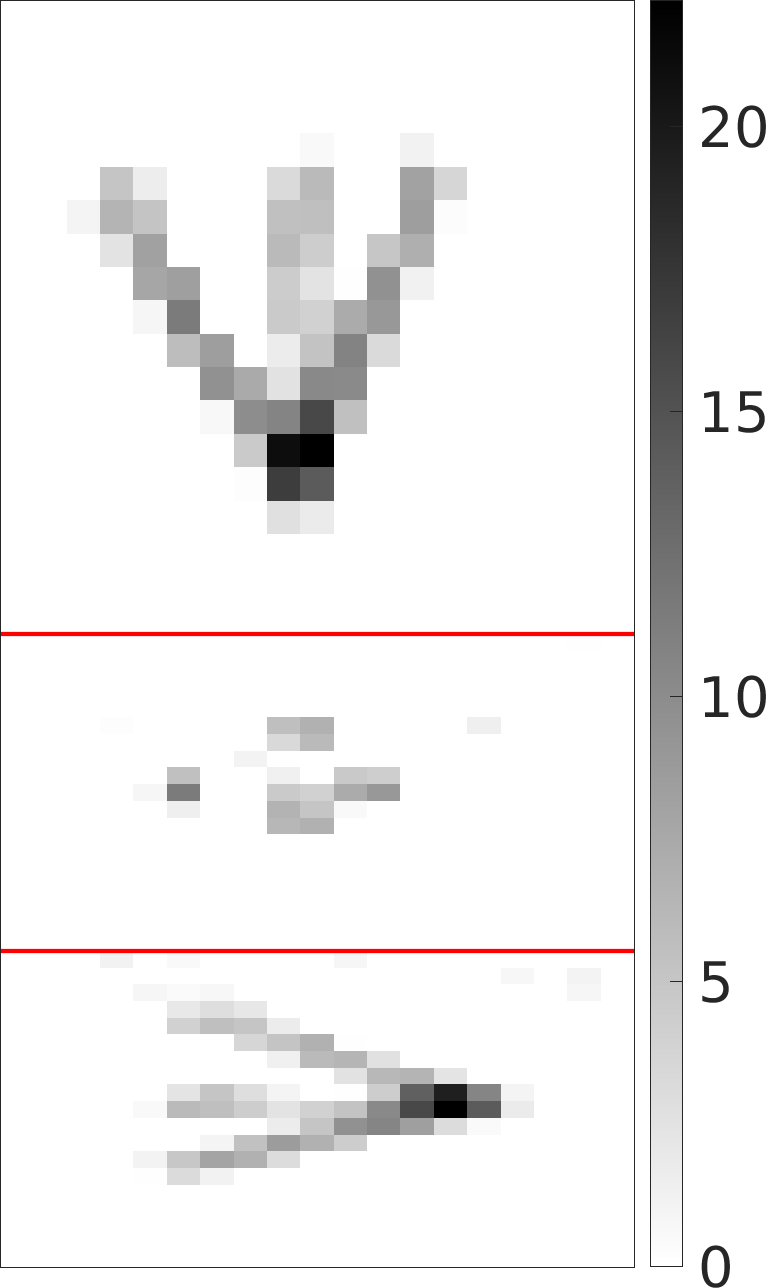}&
 \includegraphics[height=3.4cm]{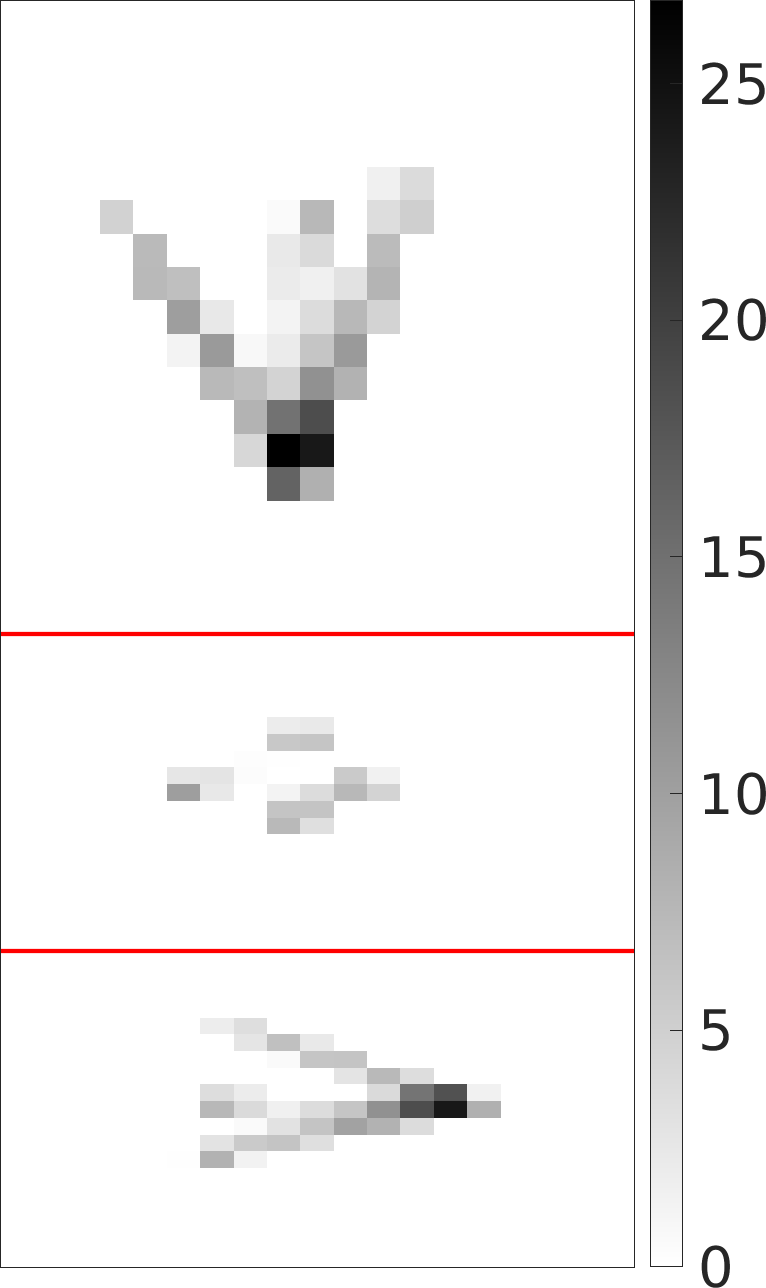}&
 \includegraphics[height=3.4cm]{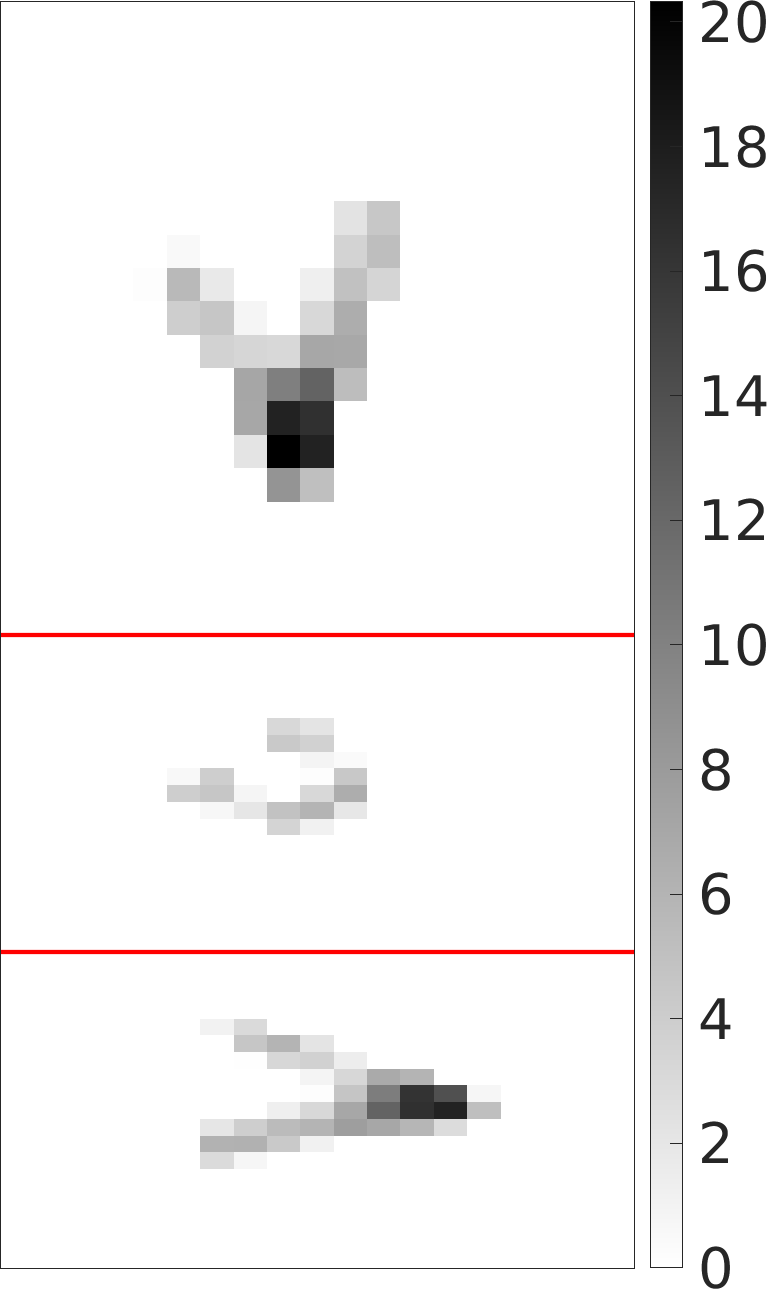}&
 \includegraphics[height=3.4cm]{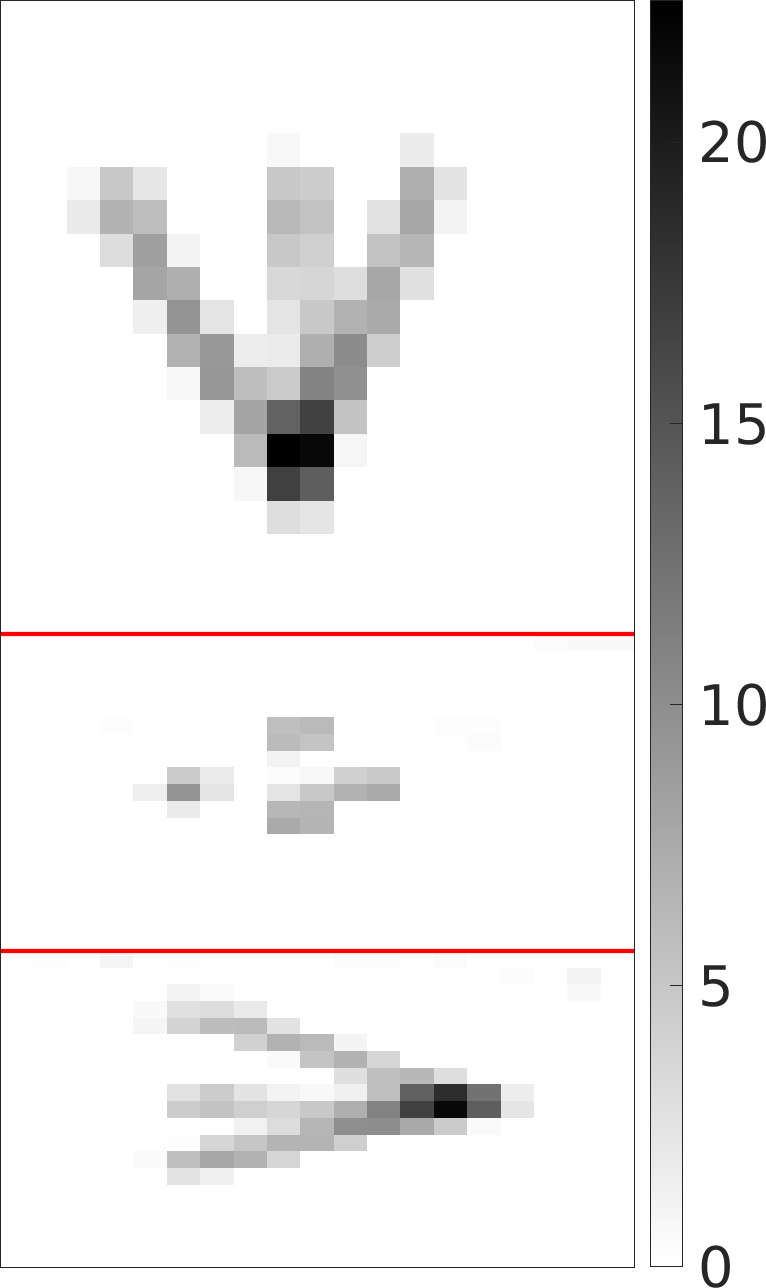}\\
\hline
\multicolumn{6}{l}{$\tau=3$} \\
 \includegraphics[height=3.4cm]{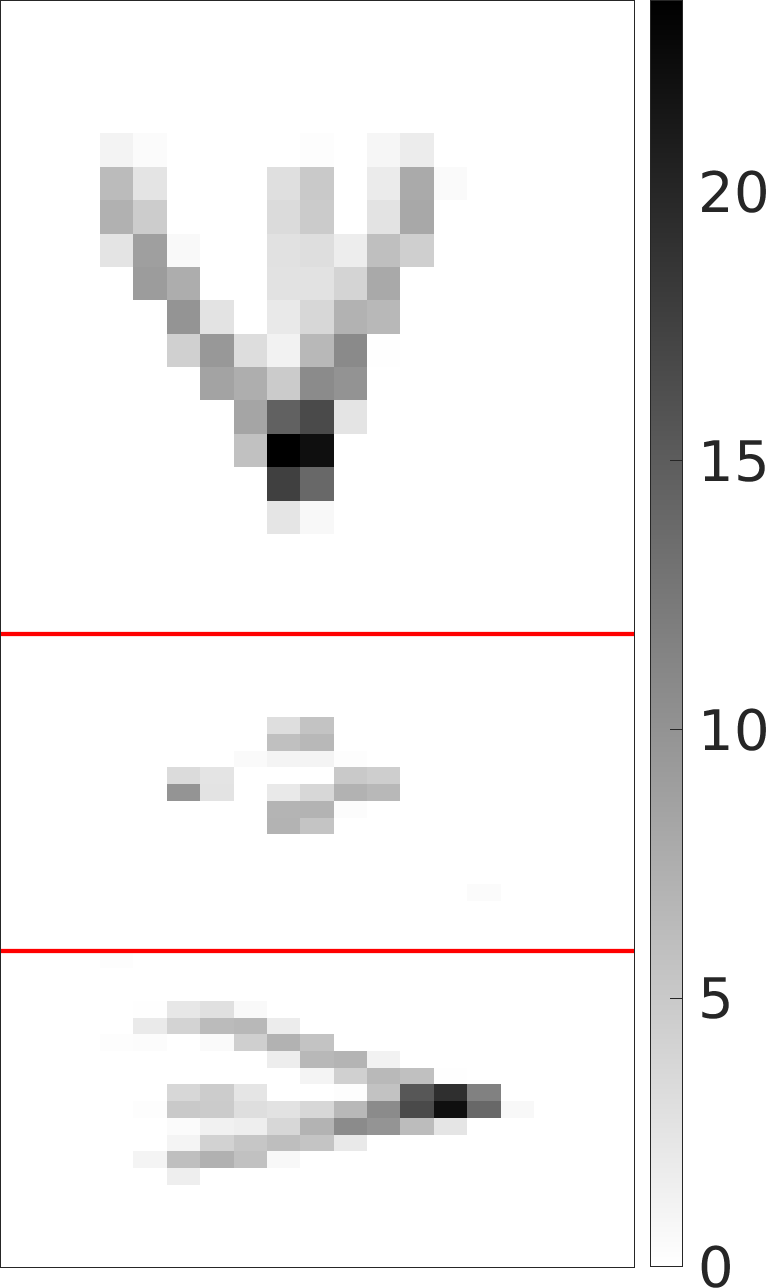}&
 \includegraphics[height=3.4cm]{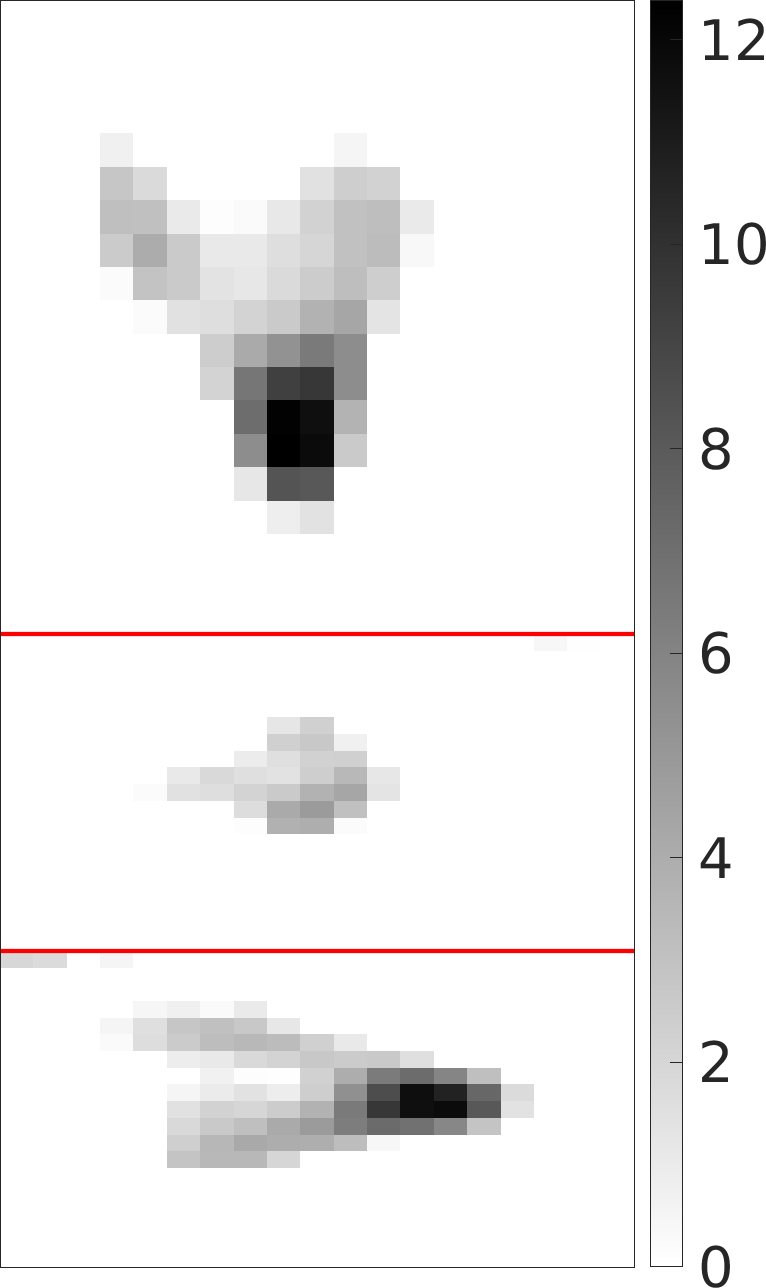}&
 \includegraphics[height=3.4cm]{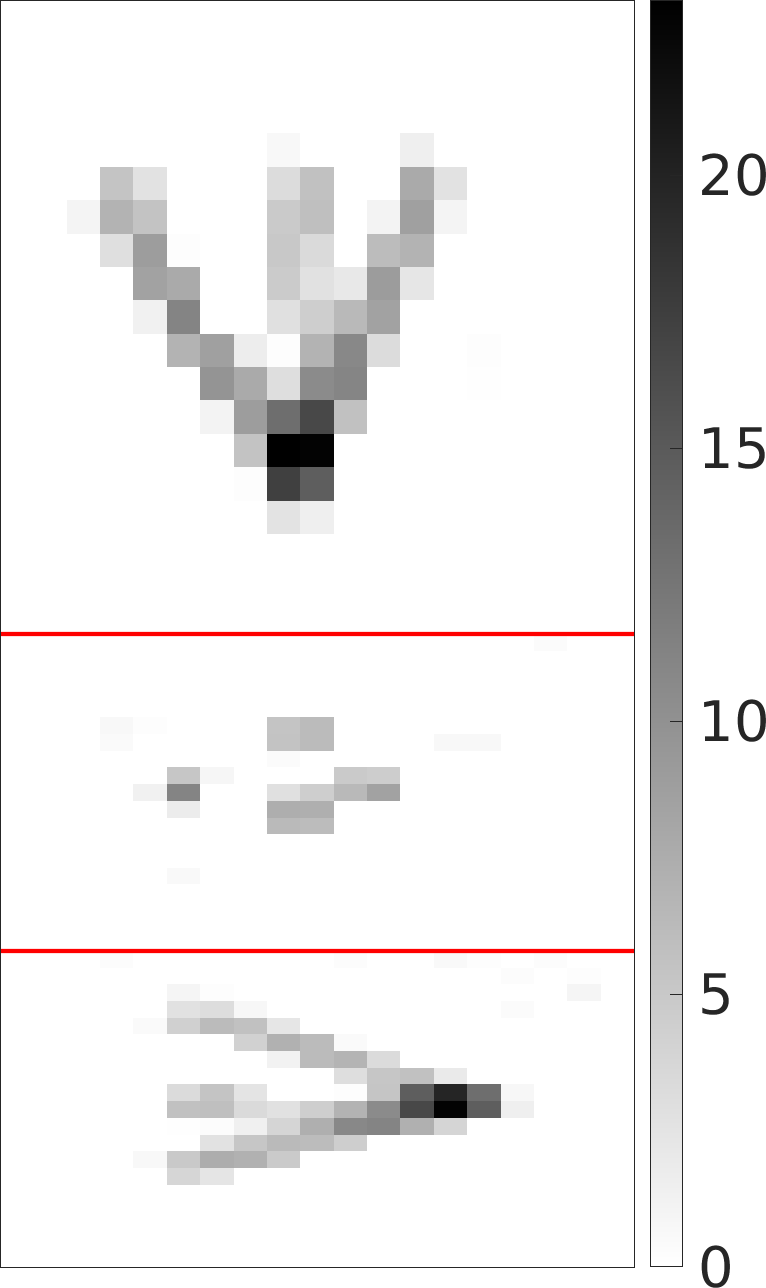}&
 \includegraphics[height=3.4cm]{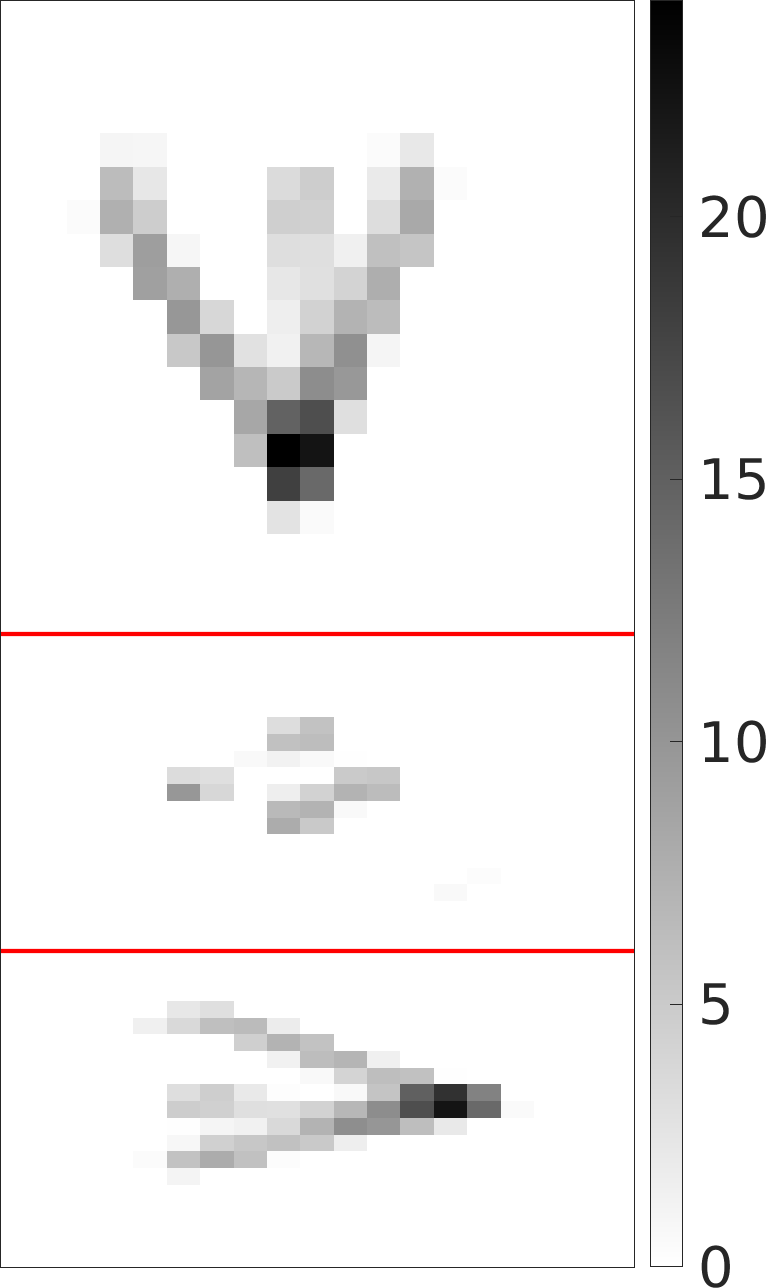}&
 \includegraphics[height=3.4cm]{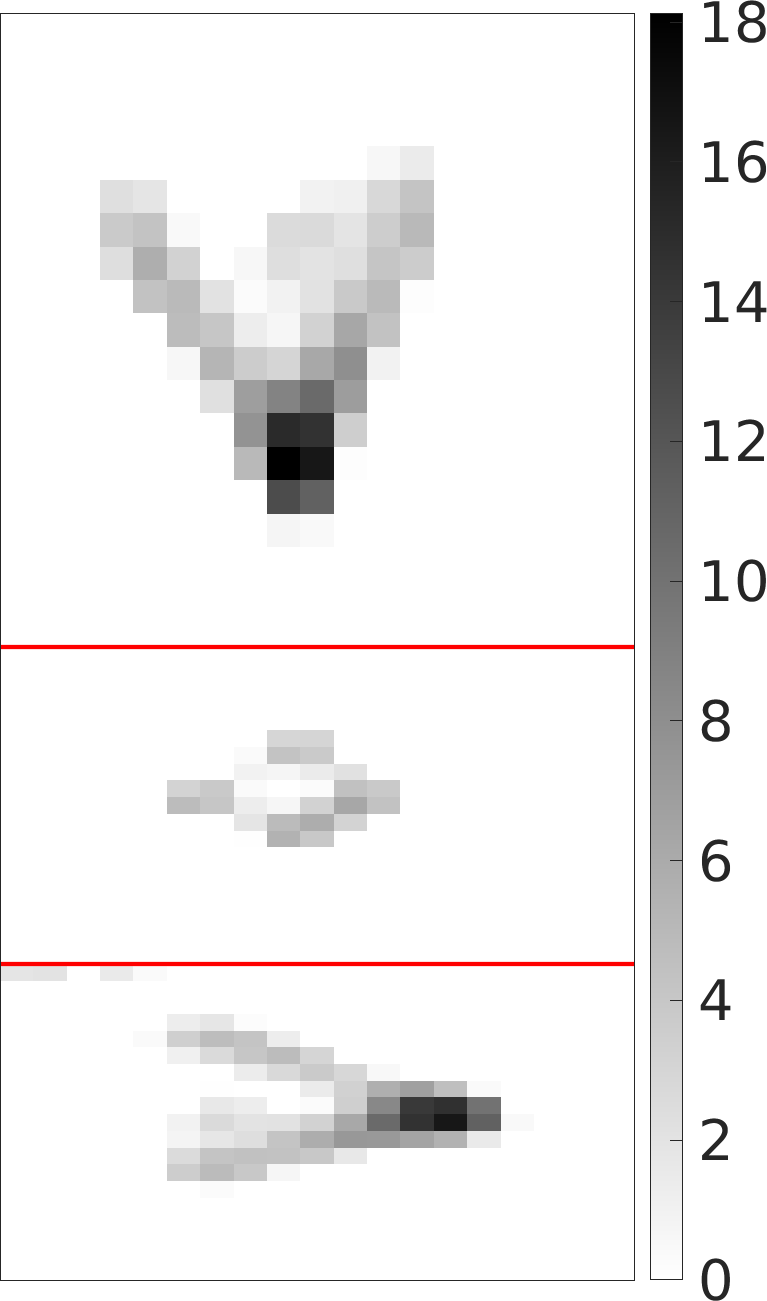}&
 \includegraphics[height=3.4cm]{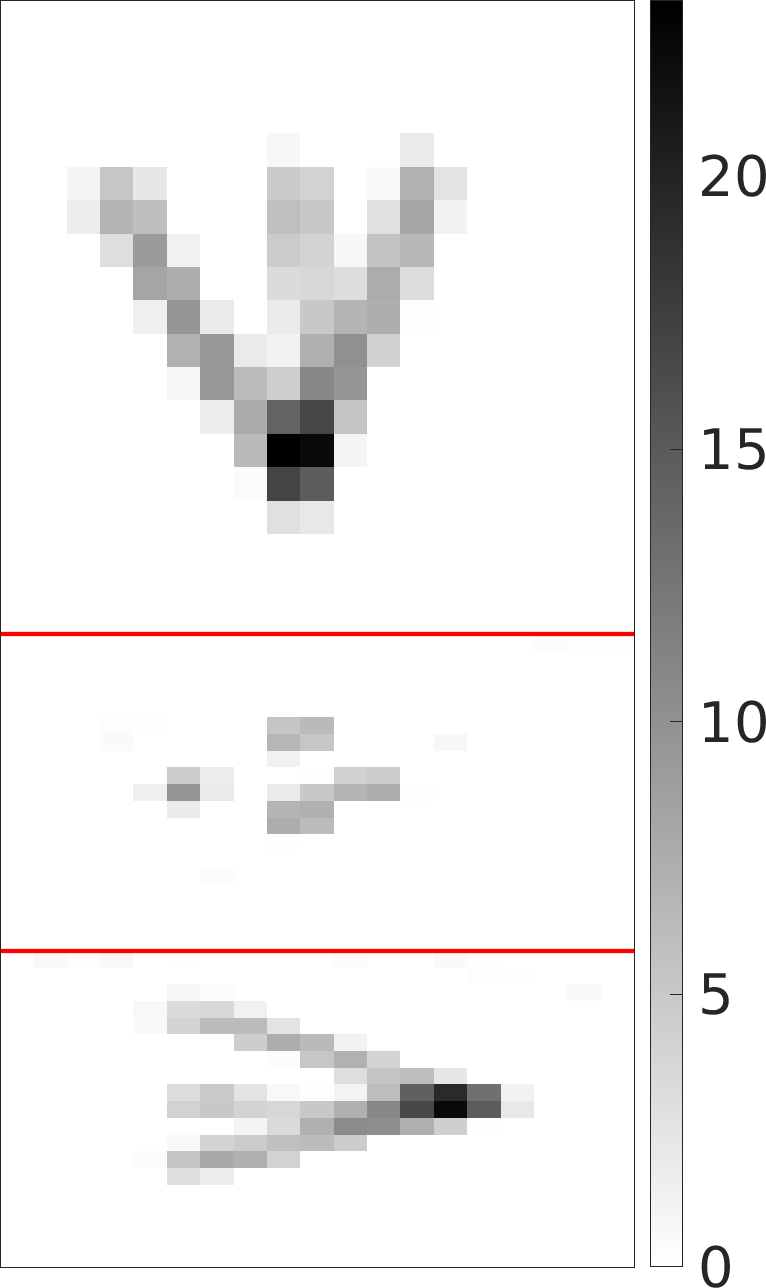}\\
 \hline
\multicolumn{6}{l}{$\tau=5$} \\
 \includegraphics[height=3.4cm]{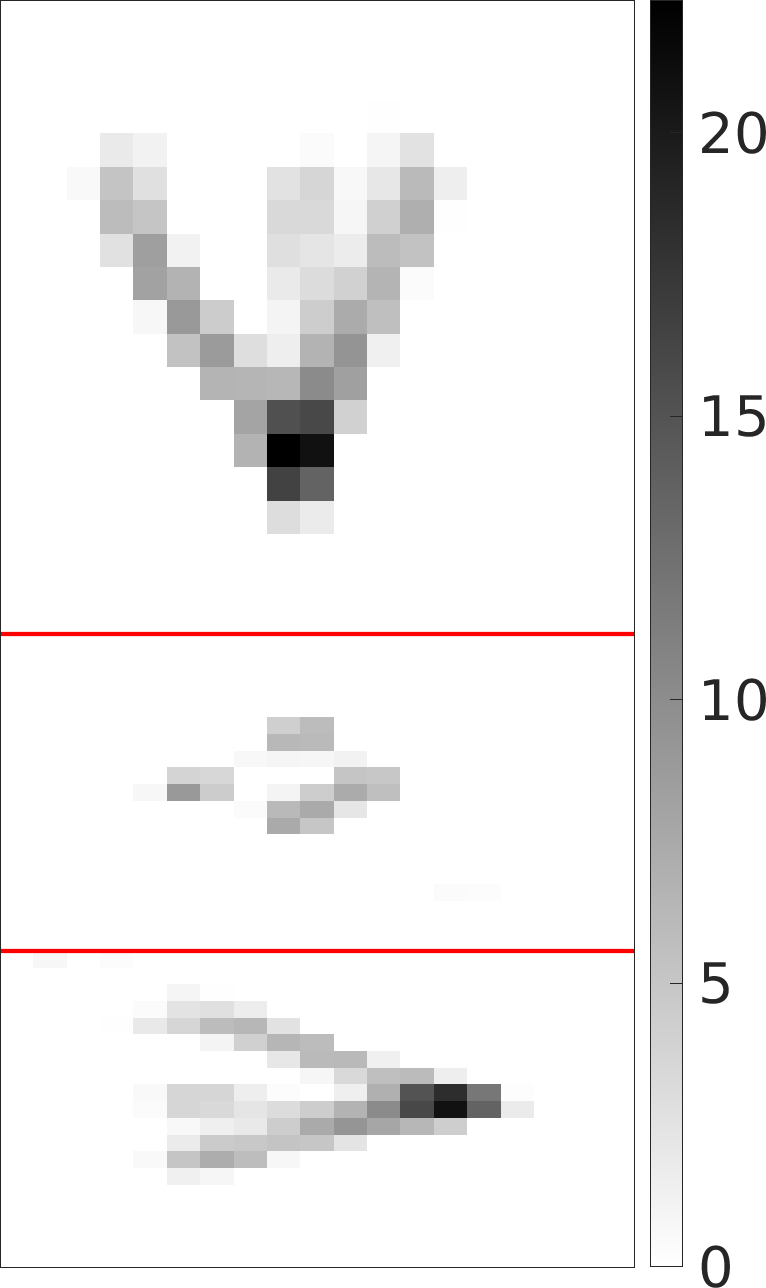}&
 \includegraphics[height=3.4cm]{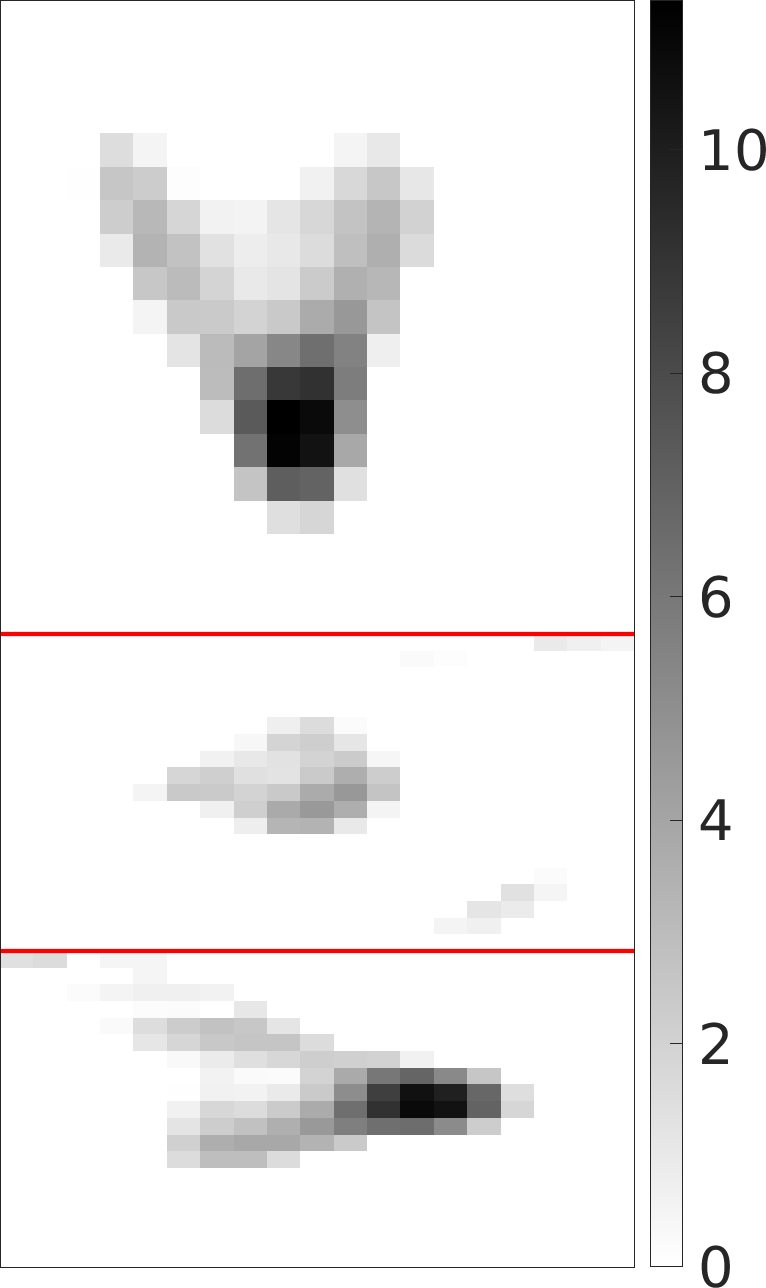}&
 \includegraphics[height=3.4cm]{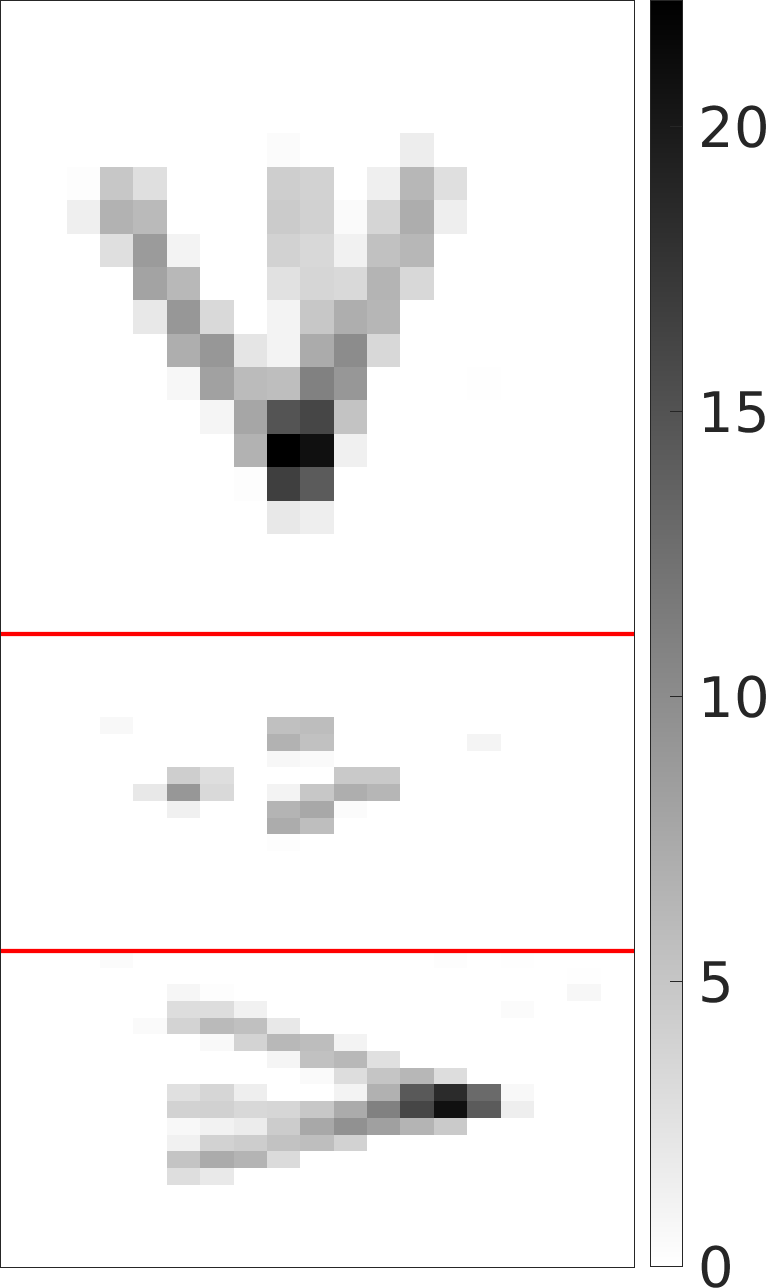}&
 \includegraphics[height=3.4cm]{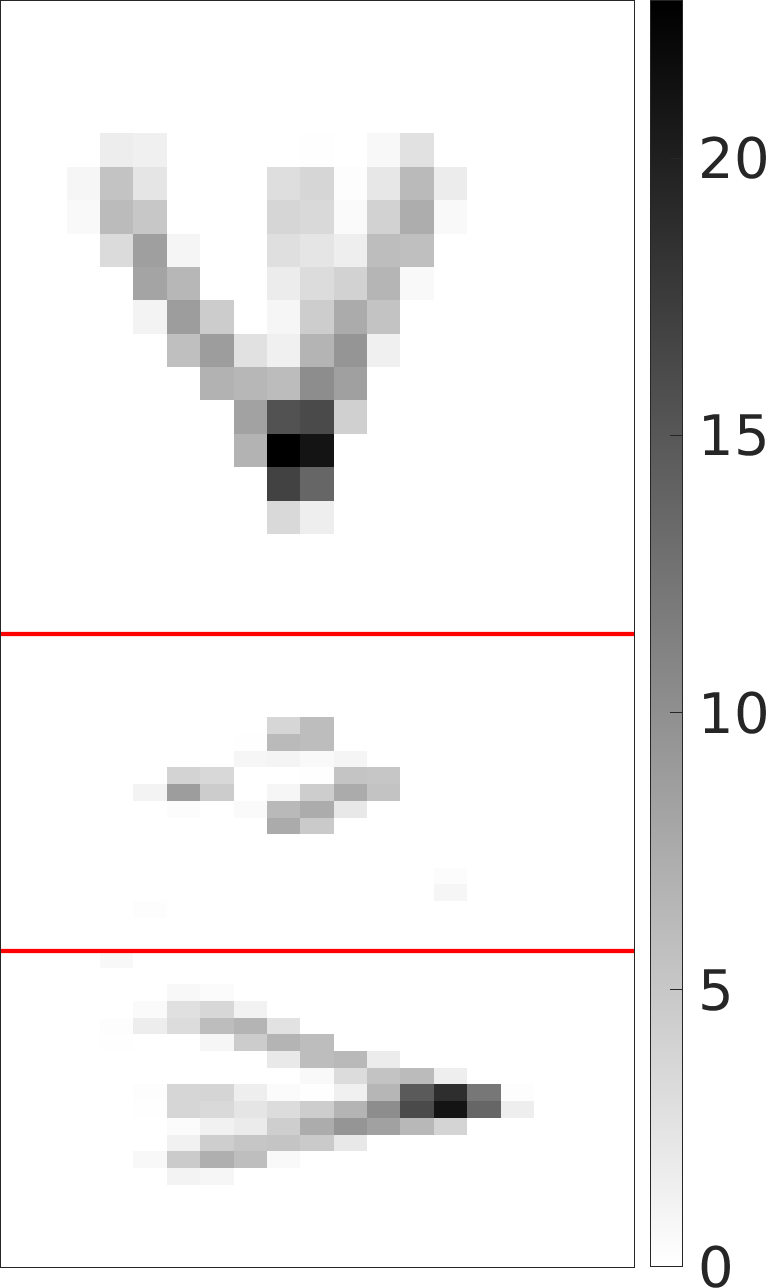}&
 \includegraphics[height=3.4cm]{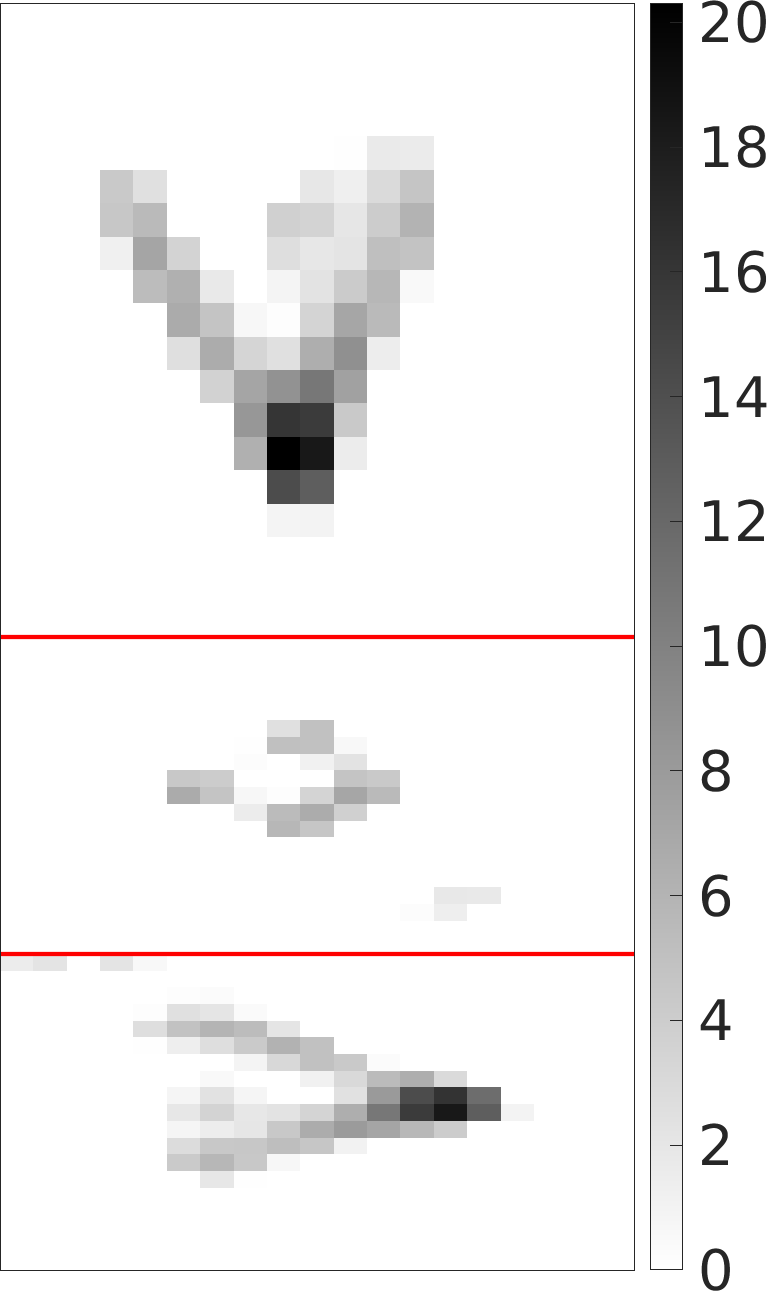}&
 \includegraphics[height=3.4cm]{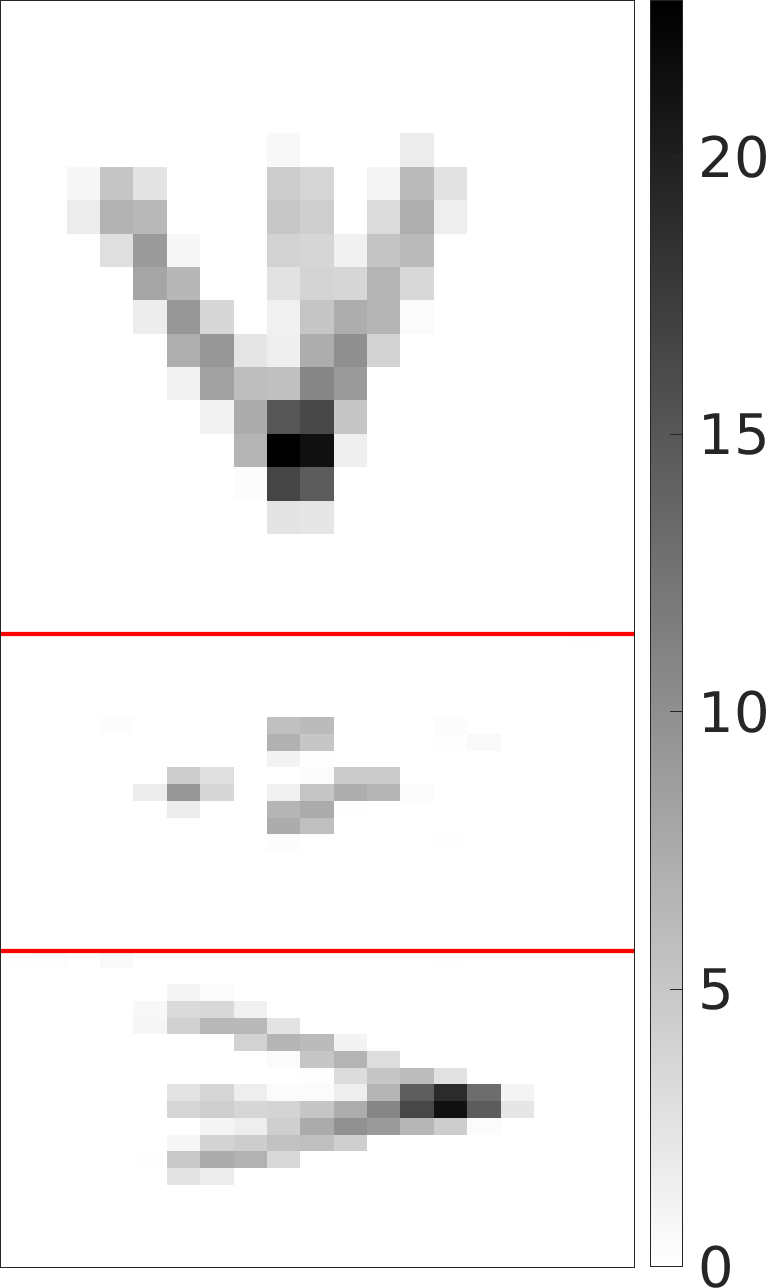}\\
\end{tabular}
}
\caption{Fig. \ref{fig:methods_nonwhitened_vs_whitened_resolution_ssim} with inverted colormap: ``Resolution'' phantom reconstructions, SSIM-optimized $\alpha$ and iteration number $N$ (for l2-K only) according to Table \ref{tab:ssim_nonwhitened_vs_whitened}.}
\label{fig:methods_nonwhitened_vs_whitened_resolution_ssim_inverted_colormap}
\end{figure}

\end{document}